\documentclass{article}

\usepackage{graphicx} % Required for inserting images
\usepackage{longtable}
\usepackage{lscape}
\usepackage{a4wide}
\usepackage{xcolor}
\usepackage{here}
\usepackage{amsmath}
\usepackage{amssymb}
\usepackage{natbib}
\usepackage{color}
\definecolor{darkblue}{rgb}{0.,0.,0.4}
\definecolor{darkred}{rgb}{0.5,0.,0.}
\usepackage[pdftex,colorlinks=true,linkcolor=darkblue,citecolor=darkred,urlcolor=blue]{hyperref}
\usepackage{authblk}
\usepackage{subcaption}
\usepackage{tikz,tikz-3dplot}
\usetikzlibrary{arrows,decorations.pathmorphing,backgrounds,fit,positioning,shapes.symbols,chains,mindmap,trees,calc}
\usepackage{listings}

\definecolor{blue}{RGB}{16,97,169} % bleu onera
\definecolor{grey}{RGB}{88, 102, 110} % gris onera 
\definecolor{green}{RGB}{65,209,204} % turquoise
\definecolor{orange}{RGB}{224, 131, 0} % corail
\definecolor{pink}{RGB}{255,105,180} % rose

\providecommand{\white}[1]{{\color{white}{#1}}} 
\providecommand{\blue}[1]{{\color{blue}{#1}}}

\providecommand{\orange}[1]{{\color{orange}{#1}}} \providecommand{\pink}[1]{{\color{pink}{#1}}}

\def\IR{{\mathbb R}}
\def\IC{{\mathbb C}}

\def\IL{{\mathbb L}}

\def\IN{{\mathbb N}}

\newcommand{\bC}{{\mathbf C}}
\newcommand{\bD}{{\mathbf D}}

\newcommand{\bG}{{\mathbf G}}

\newcommand{\bN}{{\mathbf N}}

\newcommand{\bH}{{\mathbf H}}

\newcommand{\bW}{{\mathbf W}}

\newcommand{\bV}{{\mathbf V}}

\newcommand{\bc}{{\mathbf c}}

\newcommand{\bn}{{\mathbf n}}
\newcommand{\bd}{{\mathbf d}}

\newcommand{\bv}{{\mathbf v}}
\newcommand{\bw}{{\mathbf w}}

\newcommand{\cN}{{\cal N}}

\newcommand{\bPhi}{ \boldsymbol{\Phi} }
\newcommand{\bphi}{ \boldsymbol{\phi} }

\newcommand{\bSigma}{\boldsymbol{\Sigma}}

\newcommand{\bone}{{\mathbf 1}}

\newcommand{\bvec}{{\mathbf{vec~}}}
\newcommand{\vargen}[2]{{#1}_{#2}}
\newcommand{\var}[1]{\vargen{x}{#1}}
\newcommand{\bvar}[1]{\mathbf{x}_{#1}}

\newcommand{\lan}[1]{\col{\vargen{\lambda}{#1}({j_{#1}})}}
\newcommand{\mun}[1]{\row{\vargen{\mu}{#1}({i_{#1}})}}
\newcommand{\lani}[2]{\col{\vargen{\lambda}{#1}({{#2}})}}
\newcommand{\muni}[2]{\row{\vargen{\mu}{#1}({{#2}})}}

\newcommand{\tableau}[1]{\mathcal{T}_{#1}^\otimes}
\newcommand{\ord}[0]{{n}}
\newcommand{\flop}[0]{{{\rm \texttt{flop}}}}
\newcommand{\bytes}[0]{{{\rm \texttt{Bytes}}}}

\newcommand{\pare}[1]{\left( #1\right)}

\providecommand{\col}[1]{{{#1}}} 
\providecommand{\row}[1]{{{#1}}} 

\providecommand{\COD}[0]{\textbf{C-o-D}} 

\newtheorem{theorem}{Theorem}
\newtheorem{remark}{Remark}
\newtheorem{corr}{Corrolary}

\def\CAS{50}

\title{Tensor-based multivariate function approximation: \\ methods benchmarking and comparison}

\author[1]{C. Poussot-Vassal\thanks{Corresponding author: \texttt{charles.poussot-vassal@onera.fr}}}
\author[2]{I-V. Gosea}
\author[1]{P. Vuillemin}
\author[3]{A.C. Antoulas}
\affil[1]{ONERA, DTIS, Universit\'e de Toulouse, France (\texttt{charles.poussot-vassal@onera.fr} and \texttt{pierre.vuillemin@onera.fr})}
\affil[2]{Max Planck Institute, CSC Group, Magdeburg, Germany (\texttt{gosea@mpi-magdeburg.mpg.de})}
\affil[4]{Department of Electrical and Computer Engineering, Rice University, Houston (\texttt{aca@rice.edu})}

\date{\today}

\begin{document}

\maketitle

\begin{abstract}
We evaluate some methods designed for tensor- (or data-) based multivariate model construction (approximation and compression). To this aim, a collection of multivariate functions and an evaluation methodology are suggested. First, these functions, with varying complexity (e.g., number and degree of the variables) and nature (e.g., rational, irrational, differentiable or not, symmetric, etc.) are used to build $\ord$-dimensional tensors, each of different dimension and memory size. Second, grounded on this tensor, we evaluate the performances of different methods and implementations leading to different types of surrogate models (e.g., rational functions, networks). The accuracy, the computational time, the parameter tuning impact, etc. are monitored and reported. One objective is to evaluate the different available strategies to guide users on the prospects, advantages, and limits of the various tools. The contributions are twofold: (i) to suggest a comprehensive benchmark collection together with a methodology for tensor approximation with a surrogate model and, in addition, (ii) to provide a digest and additional details of the multivariate Loewner Framework (\textbf{mLF}) approach \cite{Antoulas:2025}, as well as detailed examples and code. 
\end{abstract}

%(e.g. the approximation, the numerical analysis, the scientific computing, the dynamical systems and the neural network fields)

%%%%%%%%%%%%%%%%
\addtocontents{toc}{\setcounter{tocdepth}{2}} 
\tableofcontents
\newpage 

%%%%%%%%%%%%%%%%
\section{Introduction}
\label{sec:intro}
\subsection{Forewords}

Tensor compression and multivariate function approximation are the main topics of this report. From the lens of model approximation and surrogate construction, these two tasks are intrinsically connected as the common objective is to replace the multi-dimensional tensor or multivariate function, with a surrogate model such as a rational one, a Multi Layer Perceptron (MLP), a Kolmogorov Arnold Network  (KAN), etc.  This document is motivated by practical considerations: providing a benchmark collection and a methodology to evaluate different methodologies for performing compression and approximation. Different software are considered. Obviously, as the literature in surrogate, tensor approximation, function approximation is large and numerous, we do not claim an exhaustive evaluation. However, we provide practitioners an insight and evaluation on some recent methods. Here, attention is given to practical considerations, while for technical we refer to references. This article also aims at completing and providing a digest of the recently published contribution \cite{Antoulas:2025}, in terms of additional details and explanation, omitted out due to space limitations.

\subsection{Problem description}

\subsubsection{Multivariate functions}

A continuous $\ord$-variable function $\bH$ is defined as\footnote{Notice that $\bH$ may describe any mathematical expression, function, experimental setup, software or process.} 
\begin{equation}
\label{eq:H}
\begin{array}{ccl}
({\cal X}_1,{\cal X}_2,\cdots,{\cal X}_\ord) &\longrightarrow& \cal Y \\
(\var{1},\var{2},\cdots,\var{\ord}) &\longmapsto& y:=\bH(\var{1},\var{2},\cdots,\var{\ord}) 
\end{array},
\end{equation}
where $\var{l}\in{\cal X}_l$ ($l=1,\cdots,\ord$) is the $l$-th input variable of $\bH$, and $y\in{\cal Y}$ is the output variable. In a general (continuous) setting, these sets  denote either the set $\IR$ of real numbers or the set $\IC$ of complex numbers.

\begin{remark}[Domain restriction] \label{rem:domain}
We restrict our evaluation to the real and bounded domains, thus to real-valued multivariate functions, i.e., ${\cal X}_l:=\left[ \underline{x_l},\overline{x_l}\right]\subseteq \IR \text{~~and~~} {\cal Y}:=\left[ \underline{y},\overline{y}\right]\subseteq\IR$. These restrictions may be removed for some configurations and methods, but are necessary to compare as fairly as possible the approaches considered here.  
\end{remark}

\subsubsection{From multivariate functions to tensors (with grid structure)} 

\begin{figure}[H]
$$
\left.
\begin{array}{rcl}
{\bvar{1}} &=& \left[\begin{array}{cccc}\var{1}(1)&\var{1}(2)&\cdots&\var{1}(N_1)\end{array}\right]\\
{\bvar{2}} &=& \left[\begin{array}{cccc}\var{2}(1)&\var{2}(2)&\cdots&\var{2}(N_2)\end{array}\right]\\
&\vdots& \\
{\bvar{\ord}} &=& \left[\begin{array}{cccc}\var{\ord}(1)&\var{\ord}(2)&\cdots&\var{\ord}(N_\ord)\end{array}\right]\\
\end{array}
\right\}
\xrightarrow{~~~~\bH~~~~} \tableau{\ord}
\begin{array}{c}
%\tableau{\ord} \\
\scalebox{.35}{\input{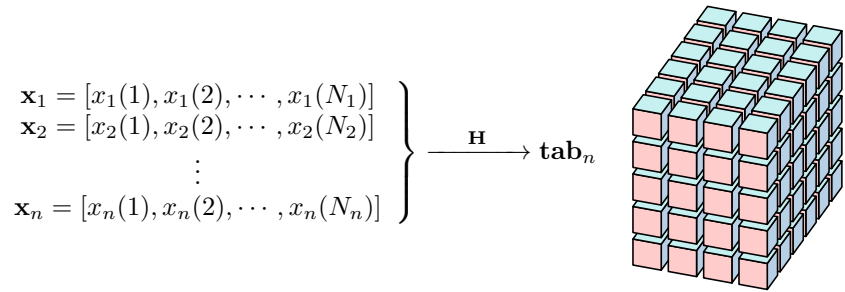}}
\end{array}
$$
\caption{Illustration of the data to tensor construction  (via $\bH$). On the left-hand side are the discrete data along each variables, while on the right-hand side is the tensor with grid structure  (here, graphical representation limited to $\ord=6$).}
\label{fig:tab}
\end{figure}

Evaluating \eqref{eq:H} over a finite discretization grid along each variable, each with finite dimension $\{N_1,N_2,\dots,N_\ord\}\in\IN$, leads to
\begin{equation}
\label{eq:tab}
\begin{array}{ccl}
({\cal X}^{N_1}_1,{\cal X}^{N_1}_2,\cdots,{\cal X}^{N_\ord}_\ord) &\longrightarrow& {\cal Y}^{N_1\times N_2\times \cdots \times N_\ord} \\
({\bvar{1}} , {\bvar{2}},\cdots,{\bvar{\ord}}) &\longmapsto& \tableau{\ord}:=\bH({\bvar{1}},{\bvar{2}},\cdots,{\bvar{\ord}}) 
\end{array},
\end{equation}
where  ${ \bvar{l}} \in{\cal X}_l^{N_l}$ ($l=1,\cdots,\ord$) is the discretized vector of the $l$-th variable, with dimension $N_l$ (i.e., discrete set of variable $\var{l}$ within the considered bounds). The resulting $\ord$-array tableau, denoted  $\tableau{\ord}$, is a \textbf{$\ord$-dimensional tensor evaluated on a dense grid}, as illustrated in Figure Figure \ref{fig:tab}, where $\var{l}(j_l)$ denotes the $j_l$-th element of the $l$-th variable ($j_l=1,\cdots,N_l$).

\subsubsection{Tensor-based model approximation: context and motivation}

Tensor approximation aims at constructing (exact, simplified or reduced-order) surrogate models (i.e., as a function, network, or realization format, etc.) that accurately capture the behavior of a potentially large-scale multi-dimensional tensor data-set, constructed from simulations or experiments, evaluated along the $\ord$ variables. Ultimately, one may expect to \textbf{discover the true underlying function that generated the tensor, its complexity and properties}. In general settings, these data may result from any measurements obtained from a parametrized experiment. Here, the simulator or the experiment is materialized by the function $\bH$, considered as unknown. The evaluations of $\bH$ along the grid points $\bvar{1}, \bvar{2},\dots,\bvar{\ord}$ generate the outputs (tensor) to be approximated, or function to be discovered. More specifically, being given \eqref{eq:tab}, we seek $\bG$ described as
\begin{equation}
\label{eq:G}
\begin{array}{ccl}
({\cal X}_1,{\cal X}_2,\cdots,{\cal X}_\ord) &\longrightarrow& {\cal \hat Y} \\
(\var{1},\var{2},\cdots,\var{\ord}) &\longmapsto& \hat y:=\bG(\var{1},\var{2},\cdots,\var{\ord}) 
\end{array},
\end{equation}
where $\hat y\in{\cal \hat Y}$ is the approximated function output and space, respectively. Obviously, we seek $\bG$ such that $\hat y \approx y$ and eventually $\bG \approx \bH$, i.e., recover or approximate the original model or system properties (either dynamic or static). Approximation is also connected to compression as the original model is also simplified.

\begin{remark}[The case of dynamical systems]
In the special context of dynamical systems governed by differential or algebraic equations, the multivariate nature comes from the parametric dependency of the underlying dynamical system:  the (first) variable being the dynamic $s$-(Laplace) or $z$-variable. This first variable accounts for the dynamical nature (frequency or time-dependency) while the rest of the parameters or variables account for physical characteristics such as mass, length, or material properties (in mechanical systems), flow speed, temperature (in fluid cases), chemical properties (in biological systems), age, weight, pressure (in clinical systems), etc. In many applications, the parameters are embedded within the model as tuning variables for the output of interest. One specific aspect is the physical meaning of this first dynamic variable, often complex, which deserves a specific treatment. In  \cite{Antoulas:2025},  this point is also considered through the construction of a multivariate Lagrange realization associated to $\bG$. This setting is out of the scope of the paper. Here we limit exposition to static multivariate functions.
\end{remark}

\subsubsection{Tensors and the curse of dimensionality (\COD)}

According to Richard E. Bellman, the \emph{"curse of dimensionality"} (\COD) refers to the diverse phenomena occurring when analyzing or ordering data in large dimensional spaces, that are not present in lower cases \cite{bellman1966dynamic}\footnote{See also the Wikipedia dedicated page \url{https://en.wikipedia.org/wiki/Curse_of_dimensionality}.}. In this note, and following  \cite{Antoulas:2025}, we are using the \COD \,term to refer to both the \textbf{computational} (floating point arithmetic, \flop) and to the \textbf{storage} (size on the disk, \texttt{Bytes}) limitations encountered when constructing multivariate model approximation from large multi-dimensional data sets as defined in \eqref{eq:tab}; as a side effect, we also claim that taming the \COD \, will also notably \textbf{improve the accuracy}\footnote{This claim may be discussed or amended, while notably observed in the conducted experiments.}. 

Accordingly, one important element presented is the impact of the dimension of the tensor in the ability of each method to be successful. In other terms, we evaluate the accuracy, the computational time and burden through complex examples (i.e., corresponding to tensors with elevated size and dimensions). \textbf{We believe that the scalability is an important feature in order for the  methods to achieve their full potential in real-life and industrial applications}.

\subsection{Contribution and structure} 

We present a benchmark collection and methodology to evaluate tensor-based multivariate function approximation methods and software. The evaluation considers the accuracy, the scalability, etc. The article is organized as follows. The current section has introduced the big picture and the problem setting. Section \ref{sec:methods} presents the benchmark problems under consideration (\CAS~ in total), as well as the evaluation methodology procedure. Then, Section \ref{sec:mlf} provides an overview of the method proposed by the authors, namely the \textbf{multivariate Loewner Framework (mLF)} \cite{Antoulas:2025}; this brief summary is accompanied by a  \texttt{MATLAB} code example\footnote{Code is based on the \texttt{MATLAB}  \texttt{+mLF} package available at \url{https://github.com/cpoussot/mLF}}. Then, Section \ref{sec:ex_overview} provides an overview of the results in term of accuracy, computational time, model complexity, etc., obtained by the \CAS~ examples considered. Comments including benefits and limitations of each method are discussed. Section \ref{sec:ex_details} details the statistics and best parametrization for the different methods\footnote{In addition, when not too long, in a  detailed analysis is given for the method exposed in \cite[Alg. 1]{Antoulas:2025}. We believe this exhaustive collection provides insightful details for researchers and practitioners, as well as a comprehensive view of the method presented in \cite{Antoulas:2025}.}. Conclusions and outlook are discussed in Section \ref{sec:conclusion}. The following contributions are made:
\begin{itemize}
\item[(i)] We gather a collection of test functions of different nature and complexity, and provide an evaluation procedure, including the accuracy, velocity and user-experience  (Section \ref{sec:methods});
\item[(ii)] We provide a comprehensive summary of the recently published work \cite{Antoulas:2025}, by means of additional unpublished practical details and explanations omitted due to space limitations. A brief tutorial of the implementation software tools (based on a library made available) is also provided (Section \ref{sec:mlf}).
\item[(iii)] We promote experimenting with different classes of methods with an emphasis on applications and practical features (Section \ref{sec:ex_overview} \& \ref{sec:ex_details}).
\end{itemize}

%\begin{remark}[Report evolutions]
%Authors insist that the present report is aimed at being updated along with time, with updated codes, methods and additional examples. In this philosophy, feedbacks from readers are welcome and will be carefully used for improvements of future versions.
%\end{remark}

\begin{remark}[Acknowledgements \& third party software]
In what follows, we investigate different methods aiming at constructing models on the basis of tensors. While the first three methods (later denoted by M1, M2 and M3) have been implemented by the authors, the remaining ones (later denoted by M4, M5, M6 and M7) are constructed by third parties. We want to give them credit in making the code available.%, second, we want to point out that as non-experts and not main authors, we may have badly parametrized them. Therefore, no conclusion regarding neither the nature nor the quality of these algorithms is intended. Last, authors want to point out that using these codes was actually relatively easy.
\end{remark}

\subsection{Reproducibility}

In order to ensure reproducibility, a regularly updated \texttt{MATLAB} code is provided at:
\begin{center}
\url{https://github.com/cpoussot/benchmark_tensor}
\end{center}
It is aimed at reproducing the results and figures. User should previously download the methods mentioned from the provided sources. See the above page for step by step details and description.

%%%%%%%%%%%%%%%%
%\newpage
\section{Benchmarks and evaluation methodology}
\label{sec:methods}
In this section, we first present the various methods and software tools under considerations (Section \ref{ssec:methods}), then we provide a brief description of the considered functions (Section \ref{ssec:functions}), and finally, we present the suggested evaluation methodology (Section \ref{ssec:methods_eval}).
\subsection{Overview of the evaluated methods}\label{ssec:methods}

We compare different tensor-driven  multivariate approximation methods (or software tools). Each method has its own tuning parameters. In what follows a subset of possible parametric configuration combinations is evaluated. %These configurations are detailed in what follows.

\paragraph{M1 - Method 1 \cite[Alg. 1]{Antoulas:2025}}  

\texttt{MATLAB} implementation of the \textbf{direct mLF} multivariate rational model approximation  \cite[Alg. 1]{Antoulas:2025}, with the following tunable parameters: %\texttt{\{tol\_ord / null\_method\}}.
\begin{itemize}
\item \texttt{tol\_ord}: $[1/2,10^{-1}, 10^{-2}, 10^{-3}, 10^{-4}, 10^{-6}, 10^{-9}, 10^{-10}, 10^{-11}, 10^{-12}, 10^{-13}, 10^{-14}]$, is the normalized singular values tolerance threshold used in the univariate Loewner step for the order selection;
\item \texttt{null\_method}: $[1,2,3]$ is the null space computation method used, being either:
\begin{enumerate} 
\item SVD decomposition (using last right singular vector);
\item QR decomposition (using the last right orthogonal factor vector);
\item linear resolution \texttt{\textbackslash} of the Loewner matrix first $k-1$ columns with the last $k$-th one;
\end{enumerate}
\end{itemize}
Additional details and code are available here: \url{https://github.com/cpoussot/mLF}. We also refer to Section \ref{sec:mlf} for details and examples.

\paragraph{M2 - Method 2 \cite[Alg. 2]{Antoulas:2025}} 

\texttt{MATLAB} implementation of the (AAA-like) \textbf{adaptive mLF} multivariate rational model approximation \cite[Alg. 2]{Antoulas:2025}, with the following tunable parameters: %\texttt{\{tol / null\_method\}}.
\begin{itemize}
\item \texttt{tol}: $[10^{-15}]$, is the maximal mismatch error tolerance (weighted by the maximal value);
\item \texttt{null\_method}: similar to M1.
\end{itemize}
Additional details and code are available here: \url{https://github.com/cpoussot/mLF}.

\paragraph{M3 - Method 3 \cite{mdspack}} 

\href{https://mordigitalsystems.fr/static/mdspack_html/MDSpack-guide.html}{\texttt{MDSPACK}} implementation of the \textbf{direct mLF} multivariate rational model approximation. It is an adaptation of M1, \cite[Algorithm 1]{Antoulas:2025}. It is a compiled \texttt{FORTRAN} code interfaced with \texttt{MATLAB}, \texttt{Python} and \texttt{Command Line} interfaces, developed at  \href{https://mordigitalsystems.fr}{\texttt{MOR Digital Systems}}, with the following tunable parameters: %\texttt{\{tol\_ord / tol\_k\}}.
\begin{itemize}
\item \texttt{tol\_ord}: similar as in M1 \& 0, meaning that no order detection is applied;
\item \texttt{tol\_k}: $[10^{-2}, 10^{-4}, 10^{-6}, 10^{-8}, 10^{-10}, 10^{-12}, 10^{-14}, 10^{-15}, -1]$, is the relative null space vector magnitude tolerance for each elements; if below, entry is removed ($-1$ means no deletion);
\item \texttt{'method'}: \texttt{'R'}, meaning that the recursive method is used. Notice that full \texttt{'F'} and full extended \texttt{'E'} can be applied (this is out of the scope of this paper).
\end{itemize}
Additional details are available at \url{https://mordigitalsystems.fr/static/mdspack_html/MDSpack-guide.html}. This algorithm is under continuous development.

\paragraph{M4 - Method 4 \cite{Poluektov:2025}} 

\texttt{MATLAB} implementation of a multivariate \textbf{ Kolmogorov Arnold Network (KAN)}-based method, with the following tunable parameters: %\texttt{\{method / alpha / Nrun /  lambda /  n /  q / p\}}.
\begin{itemize}
\item \texttt{method}: $[1,2,3,4]$, defines different basic functions in the KAN graph, being either:
\begin{enumerate} 
\item cubic splines, identification method - Gauss-Newton;
\item cubic splines, identification method - Newton-Kaczmarz, standard;
\item cubic splines, identification method - Newton-Kaczmarz, accelerated;
\item piecewise-linear, identification method - Newton-Kaczmarz, standard;
\end{enumerate}
\item \texttt{alpha}: $[0.95,1]$, is the damping factor (learning rate) for iterative parameter update;
\item \texttt{Nrun}: $[50]$, is the number of iterations;
\item \texttt{lambda}: $[0.01]$,  is the Tikhonov regularization parameter for Gauss-Newton method;
\item \texttt{n}: $[4,6, 10]$,  is the number of bottom nodes (input layer);
\item \texttt{q}: $[4,6, 12]$, is the number of top nodes (output layer);
\item \texttt{p}: $[2n+1]$, is the number of intermediate nodes (hidden layer); this choice is the optimal according to the Kolmogorov-Arnold theorem).
\end{itemize}
Additional details are available here: \url{https://github.com/andrewpolar}. We also refer to \cite{Poluektov:2025} for the exact results, proofs and notation. 

\paragraph{M5 - Method 5 \cite{Carracedo:2023}} 

\texttt{MATLAB} implementation of the \textbf{parametric Adaptive Antoulas Anderson (pAAA)}  rational model approximation, with the following tunable parameters: %\texttt{\{tol\}}.
\begin{itemize}
\item \texttt{tol}: $[10^{-3}, 10^{-6},10^{-9}]$, is the maximal mismatch error tolerance (weighted by the maximal value).
\end{itemize}
Additional details are available here: \url{https://github.com/lbalicki/parametric-AAA}. We also refer to \cite{Carracedo:2023} for the exact results, proofs and notation. We also use the \textbf{Tensor Toolbox} available at: \url{https://gitlab.com/tensors/tensor_toolbox} (downloaded on December 22nd, 2025).

\paragraph{M6 - Method 6 \cite{Balicki:2025}} 

\texttt{MATLAB} implementation of the \textbf{Low Rank parametric Adaptive Antoulas Anderson (LR-pAAA)}  rational model approximation approach, with the following tunable parameters: %\texttt{\{tol / rank\}}.
\begin{itemize}
\item \texttt{tol}: similar as in M5;
\item \texttt{rank}: $[2, 3, 4, 5]$, is a constraint for the number of terms included in the canonical polyadic (CP) decomposition used to represent the barycentric coefficients.
\end{itemize}
The authors of M6 introduce barycentric forms that are represented in the terms of separable functions, which lead to a so-called low-rank \textbf{p-AAA} algorithm. This leverages low-rank tensor decompositions in the setting of barycentric rational approximations. 
%A common challenge in multivariate approximation methods is that multivariate problems with a large number of variables often pose significant memory and computational demands. To tackle this hurdle in the setting of p-AAA, we first introduce barycentric forms that are represented in the terms of separable functions. This then leads to the low-rank p-AAA algorithm which leverages low-rank tensor decompositions in the setting of barycentric rational approximations. 
Additional details are available here: \url{https://github.com/lbalicki/parametric-AAA}. We also refer to \cite{Balicki:2025} for the exact results, proofs and notation. We also use the \textbf{Tensor Toolbox} available at: \url{https://gitlab.com/tensors/tensor_toolbox} (downloaded on December 22nd, 2025).

\paragraph{M7 - Method 7 \cite{TensorFlow}} 

\texttt{Python} implementation of the \textbf{Multi Layer Perceptron (MLP)} neural network approximation approach. This method has many tunable parameters. In this first version we fixed all of them and we limit the application of this methods to tensors with dimension $\ord\leq 4$.  Future investigation will consider different tunings and larger tensors. 
\begin{itemize}
\item \texttt{layer}: one single neuron layer is used with dense connections (according to the Universal Approximation Theorem, this should be enough for approximation).
\item \texttt{neurons}: 64, being the number of neurons.
\item \texttt{activation}: \texttt{'relu'}, being the activation function.
\item \texttt{optimizer}: \texttt{'adam'}, begin the optimization strategy.
\item \texttt{loss}: \texttt{'mse'}, being the cost function.
\item \texttt{epochs}: 500, being the number of allowed iterations.
\end{itemize}
Additional details are available here: \url{https://www.tensorflow.org/}. %Notice that examples \#29, 34, 35 are not tested since implementation was not straightforward. Then, \#9 failed but would benefit from more efforts in the future.

\paragraph{Comments on relevant similarities and differences of the methods}

It is complicated to compare algorithms that are grounded on different mathematical background or designed for different objectives. We point out  in Table \ref{tab:method-property} some important common points and differences between the evaluated approaches. 

\begin{table}%[H]
\centering\small
\begin{tabular}{ | l || c | c | c | c | c | c | c |}
\hline    & M1 &   M2  & M3 & M4 & M5 & M6 & M7  \\ \hline
$\bG$ model structure  & rational & rational & rational  & KAN & rational  & rational & MLP \\ 
Deal with complete tensor? & yes & yes & yes  & yes & yes   & yes & yes \\
Deal with incomplete tensor? & no & no & no & yes & no   & no & yes  \\
Deal with not-gridded tensor? & no & no & no & yes & no   & no & yes  \\
Deal with real variables? & yes & yes & yes & yes & yes   & yes & yes \\
Deal with complex variables? & yes & yes & yes & no & yes   & yes & no \\ 
Objective & interp. & interp. & interp. & MSE  & interp. & interp. & MSE \\ 
 &  &  &  &  & \& MSE & \& MSE &  \\ \hline
\end{tabular}
\normalsize
\caption{Some properties and features for each method. "interp.": interpolation; "MSE": mean square error.}
\label{tab:method-property}
\end{table} 

\subsection{List of examples and assumptions made}\label{ssec:functions}

We detail the considered benchmarks and assumptions. Table \ref{tab:examples} lists all the considered examples, with additional information such as dimensions and reference. The classification, e.g. rational, polynomial and irrational is also highlighted. Each example is accessible via the \texttt{mLF} package with the \texttt{[H,info]=mlf.examples(num)}, where \texttt{num} is an integer between 1 and \CAS.%, refering to the example number.% (see Table \ref{tab:method-property}), \texttt{H} is the handle function and \texttt{info} gathers informations such as bounds on the number of variables, tensor size, references, etc.

\newpage
%mlf.make_latex_examples(1:50)
    \begin{center} \begin{longtable}{|c|l|l|l|} \hline Case & Ref. & Information & Function \\ \hline\hline \orange{\#1} & [none] &  $\ord= 2$, $\mathbf{12.5}$ \textbf{KB} & $\mathrm{ReLU}(\var{1})+\frac{1}{100}\var{2}$ \\\orange{\#2} & \cite{Liu:2025} &  $\ord= 2$, $\mathbf{12.5}$ \textbf{KB} & $\mathrm{exp}\left(\sin(\var{1}) + \var{2}^2\right)$ \\\blue{\#3} & \cite{Liu:2025} &  $\ord= 2$, $\mathbf{12.5}$ \textbf{KB} & $\var{1} \var{2}$ \\\orange{\#4} & \cite{Liu:2025} &  $\ord= 3$, $\mathbf{500}$ \textbf{KB} & $\frac{1}{3} \sum_{i=1}^3 \sin(\pi x_i/2)^2$ \\\orange{\#5} & \cite{Liu:2025} &  $\ord= 4$, $\mathbf{19.5}$ \textbf{MB} & $\mathrm{exp}\left(1/2 \left( \sin(\pi(\var{1}^2+\var{2}^2) + \sin(\pi(\var{3}^2+\var{4}^2) \right) \right)$ \\\orange{\#6} & \cite{Austin:2021} &  $\ord= 2$, $\mathbf{12.5}$ \textbf{KB} & $\frac{\mathrm{exp}\left(\var{1} \var{2}\right)}{(\var{1}^2-1.44)(\var{2}^2-1.44)}$ \\\orange{\#7} & \cite{Austin:2021} &  $\ord= 2$, $\mathbf{12.5}$ \textbf{KB} & $\mathrm{log}(2.25-\var{1}^2-\var{2}^2)$ \\\orange{\#8} & \cite{Austin:2021} &  $\ord= 2$, $\mathbf{42.8}$ \textbf{KB} & $\mathrm{tanh}(4(\var{1}-\var{2}))$ \\\orange{\#9} & \cite{Austin:2021} &  $\ord= 2$, $\mathbf{12.5}$ \textbf{KB} & $\mathrm{exp}(\frac{-(\var{1}^2+\var{2}^2)}{1000})$ \\\orange{\#10} & \cite{Austin:2021} &  $\ord= 2$, $\mathbf{52.5}$ \textbf{KB} & $|\var{1}-\var{2}|^3$ \\{\#11} & \cite{Austin:2021} &  $\ord= 2$, $\mathbf{12.5}$ \textbf{KB} & $\frac{\var{1}+\var{2}^3}{\var{1}\var{2}^2+2}$ \\{\#12} & \cite{Austin:2021} &  $\ord= 2$, $\mathbf{12.5}$ \textbf{KB} & $\frac{\var{1}^2+\var{2}^2+\var{1}-\var{2}-1}{(\var{1}-1.1)(\var{2}-1.1)}$ \\{\#13} & \cite{Austin:2021} &  $\ord= 2$, $\mathbf{12.5}$ \textbf{KB} & $\frac{\var{1}^4+\var{2}^4+\var{1}^2\var{2}^2+\var{1}\var{2}}{(\var{1}-1.1)(\var{2}-1.1)}$ \\{\#14} & \cite{Austin:2021} &  $\ord= 4$, $\mathbf{1.22}$ \textbf{MB} & $\frac{\var{1}^2+\var{2}^2+\var{1}-\var{2}+1}{(\var{3}-1.5)(\var{4}-1.5)}$ \\{\#15} & \cite{Austin:2021} &  $\ord= 2$, $\mathbf{12.5}$ \textbf{KB} & $\frac{\var{1}^2+\var{2}^2+\var{1}-\var{2}-1}{\var{1}^3+\var{2}^3+4}$ \\{\#16} & \cite{Austin:2021} &  $\ord= 2$, $\mathbf{12.5}$ \textbf{KB} & $\frac{\var{1}^3+\var{2}^3}{\var{1}^2+\var{2}^2+3}$ \\{\#17} & \cite{Austin:2021} &  $\ord= 2$, $\mathbf{12.5}$ \textbf{KB} & $\frac{\var{1}^4+\var{2}^4+\var{1}^2\var{2}^2+\var{1}\var{2}}{\var{1}^2\var{2}^2-2\var{1}^2-2\var{2}^2+4}$ \\{\#18} & \cite{Austin:2021} &  $\ord= 2$, $\mathbf{12.5}$ \textbf{KB} & $\frac{\var{1}^3+\var{2}^3}{\var{1}^2\var{2}^2-2\var{1}^2-2\var{2}^2+4}$ \\{\#19} & \cite{Austin:2021} &  $\ord= 2$, $\mathbf{12.5}$ \textbf{KB} & $\frac{\var{1}^4+\var{2}^4+\var{1}^2\var{2}^2+\var{1}\var{2}}{\var{1}^3+\var{2}^3+4}$ \\\orange{\#20} & \cite{Austin:2021} &  $\ord= 3$, $\mathbf{500}$ \textbf{KB} & Breit Wigner function \\\orange{\#21} & \cite{Austin:2021} &  $\ord= 4$, $\mathbf{1.22}$ \textbf{MB} & $\frac{\sum_{i=1}^4\mathrm{atan}(x_i)}{\var{1}^2\var{2}^2-\var{1}^2-\var{2}^2+1}$ \\\orange{\#22} & \cite{Austin:2021} &  $\ord= 4$, $\mathbf{1.22}$ \textbf{MB} & $\frac{\mathrm{exp}(\var{1}\var{2}\var{3}\var{4})}{\var{1}^2+\var{2}^2-\var{3}\var{4}+3}$ \\\orange{\#23} & \cite{Austin:2021} &  $\ord= 4$, $\mathbf{1.79}$ \textbf{MB} & $10\prod_{i=1}^4\mathrm{sinc}(x_i)$ \\\orange{\#24} & \cite{Austin:2021} &  $\ord= 2$, $\mathbf{13.8}$ \textbf{KB} & $10\mathrm{sinc}(\var{1})\mathrm{sinc}(\var{2})$ \\\blue{\#25} & \cite{Austin:2021} &  $\ord= 2$, $\mathbf{12.5}$ \textbf{KB} & $\var{1}^2+\var{2}^2+\var{1}\var{2}-\var{2}+1$ \\\orange{\#26} & \cite{Balicki:2025} &  $\ord= 3$, $\mathbf{1.65}$ \textbf{MB} & $\frac{\var{1}+\var{2}+\var{3}}{6+\cos(\var{1})+\cos(\var{2})+\cos(\var{3})}$ \\\orange{\#27} & \cite{Balicki:2025} &  $\ord= 5$, $\mathbf{90.6}$ \textbf{MB} & $\frac{\var{1}+\var{2}+\var{3}+\var{4}+\var{5}}{10+\cos(\var{1})+\cos(\var{2})+\cos(\var{3})+\cos(\var{4})+\cos(\var{5})}$ \\\orange{\#28} & [none] &  $\ord= 2$, $\mathbf{30}$ \textbf{KB} & $\left(\frac{\var{1}}{\var{1}+1}\right)^4 (1+\mathrm{exp}(-\var{2}^2)) \left(1+\var{2} \cos(\var{2}) \mathrm{exp}\frac{(-\var{1}\var{2})}{\var{1}+1}\right)$ \\\orange{\#29} & [none] &  $\ord= 2$, $\mathbf{12.5}$ \textbf{KB} & $\min(10|\var{1}|,1)\mathrm{sign}(\var{1}) + \frac{\var{1}\var{2}^3}{10}$ \\\orange{\#30} & \cite{Surjanovic} &  $\ord= 8$, $\mathbf{128}$ \textbf{MB} & Borehole function \\\blue{\#31} & [none] &  $\ord= 6$, $\mathbf{128}$ \textbf{MB} & $\var{1}^2 \var{2}^3 \var{3} \var{4} - \var{5}^2 + \var{6}$ \\\orange{\#32} & [none] &  $\ord= 2$, $\mathbf{12.5}$ \textbf{KB} & $\mathrm{atan}(\var{1}) + \var{2}^3$ \\\orange{\#33} & [none] &  $\ord= 2$, $\mathbf{28.1}$ \textbf{KB} & $\frac{\var{1}+\var{2}}{\cos(\var{1})^2+\cos(\var{2}) + 3}$ \\\orange{\#34} & [none] &  $\ord= 2$, $\mathbf{1.22}$ \textbf{MB} & Riemann $\zeta$ function (real part) \\\orange{\#35} & [none] &  $\ord= 2$, $\mathbf{1.22}$ \textbf{MB} & Riemann $\zeta$ function (imaginary part) \\{\#36} & [none] &  $\ord= 3$, $\mathbf{62.5}$ \textbf{KB} & $\frac{\var{2}}{3+1/3 \var{2}\var{1}-\var{3}^2}$ \\\orange{\#37} & [none] &  $\ord= 4$, $\mathbf{1.22}$ \textbf{MB} & $\var{1}\var{4}^3+\sin(2\var{2})\var{3}$ \\{\#38} & [none] &  $\ord= 3$, $\mathbf{1.65}$ \textbf{MB} & $\frac{\var{1}^9 \var{2}^7 + \var{1}^3 + 5 \var{3}^2}{5 \var{1}^4 + 4 \var{1}^2 + \var{3}\var{2}^3 + 1}$ \\{\#39} & [none] &  $\ord= 3$, $\mathbf{500}$ \textbf{KB} & $\frac{\var{3}+\var{1}^4}{\var{1}^3+\var{2}^2+1}$ \\{\#40} & [none] &  $\ord= 4$, $\mathbf{19.5}$ \textbf{MB} & $\frac{\var{3}\var{1}}{\var{1}^2+\var{2}+\var{3}^2+1}+\var{4}^3$ \\{\#41} & [none] &  $\ord= 5$, $\mathbf{781}$ \textbf{KB} & $\frac{\var{5}^3\var{3}\var{1}+\var{3}^2}{\var{1}^3+\var{2}\var{3}+\var{4}}$ \\{\#42} & [none] &  $\ord= 6$, $\mathbf{7.63}$ \textbf{MB} & $\frac{\var{1}+\var{3}-\sqrt{2}\var{6}^2}{\var{1}^4+\var{2}\var{3}+\var{4}^3+\var{5}^2+\var{6}}$ \\{\#43} & [none] &  $\ord= 7$, $\mathbf{76.3}$ \textbf{MB} & $\frac{\var{3}\var{2}^3+1}{\var{1}^4+\var{2}^2\var{3}+\var{4}^2+\var{5}+\var{6}^3+\var{7}}$ \\{\#44} & [none] &  $\ord= 8$, $\mathbf{763}$ \textbf{MB} & $\frac{1}{\var{1}^4+\var{2}^2\var{3}+\var{4}^2+\var{5}+\var{6}+\var{7}+\var{8}}$ \\{\#45} & [none] &  $\ord= 9$, $\mathbf{76.9}$ \textbf{MB} & $\frac{1}{\var{1}^2+\var{2}^2\var{3}+\var{4}^2+\var{5}+\var{6}+\var{7}+\var{8}+\var{9}}$ \\{\#46} & [none] &  $\ord= 10$, $\mathbf{461}$ \textbf{MB} & $\frac{1}{\var{1}+\var{1}^2\var{2}\var{3}+\var{4}+\var{5}+\var{6}+\var{7}\var{8}+\var{9}^2+\var{10}}$ \\\blue{\#47} & \cite{GHK:2025} &  $\ord= 5$, $\mathbf{1.9}$ \textbf{MB} & $\begin{array}{c}(1 + 2\var{1})(-2 + \var{2})(-\var{3})(3 + \var{4})(2- 3\var{5}) \\ + (-1 + \var{1})(2\var{2})(1 + 3\var{3})(-\var{4})(1 -\var{5})\end{array}$ \\\blue{\#48} & \cite{Polya:1925} &  $\ord= 3$, $\mathbf{13.5}$ \textbf{KB} & $\var{1}\var{2}+\var{1}\var{3}+\var{2}\var{3}$ \\\orange{\#49} & [none] &  $\ord= 2$, $\mathbf{50}$ \textbf{KB} & Hankel function $H_0$ (real part) \\\orange{\#50} & [none] &  $\ord= 2$, $\mathbf{50}$ \textbf{KB} & Hankel function $H_0$ (imaginary part) \\\hline \caption{List of examples. \blue{Polynomial: 5}, {rational: 19}, \orange{irrational: 26}.} \label{tab:examples} \end{longtable} \end{center}\normalsize

To conduct the benchmarking, let us consider the following assumptions:
\begin{itemize}
\item[A1] Each function $\bH$ depends on $\ord$ variables with one single measured output, of the form \eqref{eq:H};
\item[A2] The input variables and measured output are real-valued, as given in Remark \ref{rem:domain};
\item[A3] The common input argument  for each method is a $\ord$-dimensional tensor $\tableau{\ord}$ (with the corresponding evaluation points), as given in  \eqref{eq:tab} and illustrated in Figure \ref{fig:tab};
\item[A4] No noise is considered on the measured output. In this regard, we refer  reader to signal processing or system identification communities, where filtering or output averaging methods are deployed, together with statistical tools (see e.g. \cite{Pintelon:2012}).  %Obviously, this is an axis for future investigations, possibly including other contributors.
\end{itemize}

\subsection{Evaluation methodology and metrics}
\label{ssec:methods_eval}

%\paragraph{Evaluation procedure}

The evaluation procedure is common to all methods, and is detailed step by step in what follows:
\begin{enumerate}
\item[S1] Consider the $\ord$-variable function \eqref{eq:H} (and Remark \ref{rem:domain});
\item[S2] Consider a discretization of the input space and compute the tensor $\tableau{\ord}$ as in \eqref{eq:tab}; 
\item[S3] For all methods $m=\{\text{M1,\,M2,\,\dots}\}$, enumerate all combinations of the tuning parameters configurations $p$ and construct a surrogate model $\bG_{m,p}$;
\item[S4] For 500 random draws of $\{\var{1},\var{2},\cdots,\var{\ord}\}$ within the considered domain bound ${\cal X}_l$, evaluate both $\bH$ (\eqref{eq:H}) and $\bG_{m,p}$ (\eqref{eq:G}) and compute the root mean square error (RMSE) given as
\begin{equation}
\text{RMSE}_{m,p} = \sqrt{\frac{1}{500}\sum_{j=1}^{500}\big(\bG_{m,p}(\var{1},\var{2},\cdots,\var{\ord})-\bH(\var{1},\var{2},\cdots,\var{\ord})\big)^2}.
\label{eq:rmse}
\end{equation}
\item[S5] For each method $m$, keep the best model along the possible parameter set (i.e.  the one achieving lowest RMSE), now denoted $\bG_{m}$;
\item[S6] Report and plot of the best candidate $\bG_{m}$, for each method.
\end{enumerate}

\begin{remark}[Computational setup]
The computations are carried out on \texttt{MATLAB} 2023b, with a MacBook Air (with Apple M1 with 16 GB memory). We note that results may vary with different architecture.
\end{remark}

%\begin{remark}[About RMSE]
%We believe the RMSE is an interesting metric to monitor since it usually makes sense for most engineers, and none of the method actually are specifically targeted to minimize this metric.
%\end{remark}

\begin{remark}[Additional details in Section \ref{sec:ex_details}]
For each case, the original function $\bH$ \eqref{eq:H} (used to generate the tensor), the reference where it has been used (if any), domain and bounds ${\cal X}_l$ ($l=1,\cdots,\ord$) given in Remark \ref{rem:domain}, and the tensor $\tableau{\ord}$ size, are first listed. Then, for each evaluated method, the tuning parameter configuration set $p$  leading to the lowest RMSE \eqref{eq:rmse} (evaluated over 500 random input variables draw) is reported and surrogate size, mean error and computation time are given (\texttt{NaN}, when no model has been found). In addition, the method presented in \cite[Alg. 1]{Antoulas:2025} is more specifically detailed (when not too long). 
\end{remark}

\section{The multivariate Loewner Framework (\textbf{mLF}) at a glance}
\label{sec:mlf}
This method, appropriately named  \textbf{multivariate Loewner Framework (mLF)}, was originally introduced in \cite{Antoulas:2025}. At its core it relies on the Loewner Framework (LF). One important contribution in this work is to address the problem of \textbf{dimensionality}, occurring essentially when the number of variables and tensor size increase. This is achieved thanks to a  \textbf{variable decoupling}. In addition, we present connections between the LF for rational interpolation of multivariate functions and the \textbf{Kolmogorov Superposition Theorem (KST) restricted to rational functions}, resulting in the formulation of the \textbf{KST} for this special function case (see Figure \ref{fig:mlf}). As a byproduct, taming the curse of dimensionality (\COD) in computational complexity, storage and numerical accuracy, is achieved\footnote{Notice that this framework is not restricted to the real domain since all variables may belong to the complex domain.}. We first introduce the data notation, Loewner matrix and Lagrange form in Section \ref{ssec:data}. Then, the variable decoupling and its benefits are exposed in Section \ref{ssec:decoupling}. Finally, Section \ref{ssec:mlf} illustrates and details all these results using the  \href{https://github.com/cpoussot/mLF}{\texttt{mLF}} \texttt{MATLAB} package, using a very simple multivariate example.

\begin{figure}[h]
\begin{center}
%\textbf{The \pink{Loewner framework} aims at bridging approximation and control theory.}\\~\\
\scalebox{.45}{\newcommand*{\info}[4][16.3]{%
  \node [ annotation, #3, scale=0.65, text width = #1em,
          inner sep = 2mm ] at (#2) {%
  \list{$\bullet$}{\topsep=0pt\itemsep=0pt\parsep=0pt
    \parskip=0pt\labelwidth=8pt\leftmargin=8pt
    \itemindent=0pt\labelsep=2pt}%
    #4
  \endlist
  };
}

\tikzset{
    invisible/.style={opacity=0},
    visible on/.style={alt={#1{}{invisible}}},
    alt/.code args={<#1>#2#3}{%
      \alt<#1>{\pgfkeysalso{#2}}{\pgfkeysalso{#3}} % \pgfkeysalso doesn't change the path
    },
  }
  
\begin{tikzpicture}
  %%% App theory
  \path[mindmap,concept color=blue!75,text=black]
    node[concept] {\LARGE\textbf{\white{Approx. \\theory}}} [clockwise from=180]
    child[concept color=yellow!75!black] { 
    node (basis) [concept] {\large\textbf{Basis fun.}} 
    [clockwise from=300]
    child {node (rbf) [concept] {\large\textbf{RBF}}}
    child {node (otho) [concept] {\large\textbf{Orth. basis}}}
    child {node (splin) [concept] {\large\textbf{Splin}}}
    }
    child[concept color=red!60!white] { 
    node (tensor) [concept] {\large\textbf{Tensor}} 
    [clockwise from=240]
    child {node (td) [concept] {\large\textbf{Tucker dec.}}}
    child {node (ttd) [concept] {\large\textbf{Tens. Train dec.}}}
    child {node (htd) [concept] {\large\textbf{Hier. Tucker dec.}}}
    child {node (tensOther) [concept] {\large\textbf{Others}}}
    }
    child[concept color=orange!75] { 
    node (interp) [concept] {\large\textbf{Interp.}}
    [clockwise from=120]
    child {node (hermite) [concept] {\large\textbf{Hermite}}}
    child {node (rational) [concept] {\large\textbf{Rational}}}
    child {node (polynomial) [concept] {\large\textbf{Poly.}}}
    }
    child[concept color=green!75!black] {
      node (nn) [concept] {\large\textbf{NN}}
      [clockwise from=60]
      child {node (kan) [concept] {\large KAN}}
      child {node (mlp) [concept] {\large MLP}}
    };
    
    %%% Systems
    \begin{scope}[mindmap,concept color=blue!75,text=black]
    \node (system) [mindmap,concept] at (18,0) {\LARGE\textbf{\white{Systems \\ theory}}}
    [clockwise from=180]
    child[concept color=orange!75] {node (red) [concept] {\large\textbf{Reduction \\ Approx.}} [clockwise from=210]}
    child[concept color=purple!75] {node (real) [concept] {\large\textbf{Realization}} [clockwise from=210]}
    child[concept color=green!75!black] {node (ident) [concept] {\large\textbf{Identif.}} [clockwise from=210]}
    child[concept color=yellow!75!black] {node (rest) [concept] {\large\textbf{Control\\ stability \\ ...}} [clockwise from=210]};
    \end{scope}	
    
    %\path (tensor) to[circle connection bar switch color=from (red!60!white) to (orange!75)] (interp);
    %\path (real) to[circle connection bar switch color=from (purple!75) to (orange!75)] (red);
    %\path (ident) to[circle connection bar switch color=from (green!75!black) to (purple!75)] (real);
    
    %%% Loewner
    \begin{scope}[mindmap,concept color=pink!75,text=black]
    \node (loewner) [mindmap,concept] at (8,12) {\LARGE\textbf{Loewner Framework}};
    \end{scope}	
    
    \path (loewner) to[circle connection bar switch color=from (pink!75) to (red!60!white)] (tensOther);
    \path (loewner) to[circle connection bar switch color=from (pink!75) to (orange!75)] (hermite);
    \path (loewner) to[circle connection bar switch color=from (pink!75) to (orange!75)] (rational);
    \path (loewner) to[circle connection bar switch color=from (pink!75) to (orange!75)] (polynomial);
    \path (loewner) to[circle connection bar switch color=from (pink!75) to (green!75!black)] (kan);
    \path (loewner) to[circle connection bar switch color=from (pink!75) to (purple!75)] (real);
    \path (loewner) to[circle connection bar switch color=from (pink!75) to (orange!75)] (red);
    
\end{tikzpicture}}\\
\end{center}
\caption{The \pink{Loewner framework} aims at bridging approximation and control theory.}
\label{fig:mlf}
\end{figure}

%%%%%%%%%%%%%%%%%%%%%%%%%%%%%%%%%%%%%%%
\subsection{$\ord$-D tensor data, Loewner matrix, and Lagrangian form}
\label{ssec:data}

%Let us first set the different elements used for the interpolation in the Loewner framework.

\subsubsection{Tensor data and columns/rows interpolation points} 

The data-set (tensor and evaluation variables, as presented in Section \ref{sec:intro}) is the primary ingredient of the $\ord$-variable \textbf{data-driven} rational interpolation and approximation. %This set is obtained by evaluating the $\ord$-variable function  $\bH(\var{1},\var{2},\cdots,\var{\ord})$, either through a computer simulation or directly an experimental set-up along the discretized grid $\bvar{1},\cdots,\bvar{\ord}$. 
In the context of the \textbf{LF}, these points are split into columns and rows \textbf{interpolation, or support points}\footnote{In what follows, $\var{l}(i)$ denotes the $l$-th variable evaluated at the $i$-th element. We also denote $j_l=1,\cdots,k_l$ and $i_l=1,\cdots,q_l$ with $l=1,\cdots,\ord$. $k_l$ and $q_l$ are the available data along the $l$-th variable.}.

Similarly to the classical univariate 1-D, bivariate 2-D, or 3-D  scenario \cite{Antoulas:2012,Ionita:2014}, evaluating the function $\bH(\var{1},\var{2},\cdots,\var{\ord})$ along with the combinations of the support points $\{\lan{1},\lan{2},\cdots,\lan{\ord}\}\in\IR$ and $\{\mun{1},\mun{2},\cdots,\mun{\ord}\}\in\IR$, thus forms a $\ord$-dimensional tensor, denoted $\tableau{\ord}$ ($j_l=1,\cdots,k_l$, $i_l=1,\cdots,q_l$ and $l=1,\cdots,\ord$). Here, $\lani{l}{j_l}$ and $\muni{l}{i_l}$ are a separation of $x_l$ with  $k_l+q_l=N_l$. To connect notations, we define $\bvar{1}=[\lan{1}, \mun{1}]$, $\bvar{2}=[\lan{2}, \mun{2}]$, \dots, $\bvar{\ord}=[\lan{\ord}, \mun{\ord}]$. %The $\ord$-variable measurement tensor, called $\tableau{\ord}$ \eqref{eq:tab}, is schematically recalled in Figure \ref{fig:tab}. 
%\end{minipage}
%\begin{minipage}{.45\textwidth}
%$$
%\left.
%\begin{array}{c}
%\bvar{1}=[\lan{1}, \mun{1}] \\
%\bvar{2}=[\lan{2}, \mun{2}] \\
%\vdots \\
%\bvar{\ord}=[\lan{\ord}, \mun{\ord}] 
%\end{array}
%\right\}
%\xrightarrow{~~\bH~~}\tableau{\ord}
%\begin{array}{c}
%\scalebox{.4}{\input{figs/tensor}}
%\end{array}
%$$
%\end{minipage}

Following the Loewner philosophy detailed e.g. in \cite{Mayo:2007,Antoulas:2012,Ionita:2014} and \cite[Section 3]{Antoulas:2025}, let us define  \col{$P_c^{(\ord)}$}, the \col{column} data, and \row{$P_r^{(\ord)}$}, the \row{row} data, being two subsets of the original $\ord$-D tensor $\tableau{\ord}$, leading to $\bw_{j_1,j_2,\cdots,j_\ord}$ and $\bv_{i_1,i_2,\cdots,i_\ord}$. More specifically, these subsets are given as follows:
\begin{equation}
\label{eq:data_n}
\left\{
\begin{array}{l}
P_c^{(\ord)}:=\left\{(\lan{1},\lan{2},\cdots,\lan{\ord};\bw_{j_1,j_2,\cdots,j_\ord}), ~j_l=1,\ldots,k_l, ~l=1,\ldots,\ord\right\}\\
P_r^{(\ord)}:=\left\{(\mun{1},\mun{2},\cdots,\mun{\ord};\bv_{i_1,i_2,\cdots,i_\ord}),~i_l=1,\ldots,q_l, ~l=1,\ldots,\ord \right\}
\end{array}
\right. .
\end{equation}

Then, we also denote by $\bW^{\otimes}_{k_1,\cdots,k_\ord}$ the $k_1\times \cdots \times k_\ord$ tensor and by $\bW^{\otimes}_{k_1,\cdots,k_\ord}(j_1,\cdots,j_\ord)$, its $(j_1,\cdots,j_\ord)$-th element, or simply $\bw_{j_1,j_2,\cdots,j_\ord}$.

\subsubsection{$\ord$-D Loewner matrix} 

The \textbf{tensor data with a grid structure} \eqref{eq:tab}, re-written in \eqref{eq:data_n} now serves for constructing the $\ord$-D Loewner matrix $\IL_\ord$, which may be viewed as an operator mapping the interpolation points and the $\ord$-D tensor onto a $Q\times K$ matrix, with \row{$Q=q_1q_2\dots q_\ord$ (rows)} and \col{$K=k_1k_2\dots k_\ord$ (columns)}, i.e.
\begin{equation}
\begin{array}{ccl}
   \left( \IC^{k_1} \times\IC^{q_1} \times \ldots \times \IC^{k_\ord}\times \IC^{q_\ord} \times \IC^{(k_1+q_1)\times \cdots \times (k_\ord+q_\ord)}\right) & \longrightarrow & \IC^{Q\times K} \\
    (\underbrace{\lan{1},\mun{1}}_{\bvar{1}},(\underbrace{\lan{2},\mun{2}}_{\bvar{2}},\ldots,\underbrace{\lan{\ord},\mun{\ord}}_{\bvar{\ord}},\tableau{\ord}) &  \longmapsto & \IL_\ord 
\end{array},
\end{equation}
where each entry of the $\IL_\ord$ matrix is given by
\begin{equation}
\ell_{j_1,j_2,\cdots,j_{\ord}}^{i_1,i_2,\cdots,i_{\ord}} = \dfrac{\bv_{i_1,i_2,\cdots,i_\ord}-\bw_{j_1,j_2,\cdots,j_\ord}}{\pare{\mun{1}-\lan{1}}\cdots \pare{\mun{\ord}-\lan{\ord}}}.
\end{equation}
%Notice that the \emph{interpolation (or support) points} are concentrated at the denominator and where the tensor data are located at the numerator. One essential step

\subsubsection{Lagrangian (barycentric) rational model}

By considering an appropriate number of column interpolation points $k_l$ ($l=1,\cdots,\ord$), one can compute $\IL_\ord \bc_\ord = 0$, the right null space of $\IL_\ord$, which contains the so-called \textbf{barycentric coefficients}, 
\begin{equation}
\bc_\ord^\top=
\left[
\begin{array}{ccc|c|ccc}\label{eq:null}
    c_{1,\cdots,1}& \cdots & c_{1,\cdots,k_\ord} &  \cdots & c_{k_1,\cdots,1} & \cdots & c_{k_1,\cdots,k_\ord}   
\end{array}
\right]\in \IC^{K}.
\end{equation}
Then, the multivariate Lagrangian (barycentric) form 
\begin{equation}
\begin{array}{rcl}
\bG_{\mathrm{lag}}(\var{1},\cdots,\var{\ord})
&=& \dfrac{\sum_{j_1=1}^{k_1}\cdots\sum_{j_\ord=1}^{k_\ord} \frac{c_{j_1,\cdots,j_\ord}\bw_{j_1,\cdots,j_\ord}}{\pare{\var{1}-\lan{1}}\cdots\pare{\var{\ord}-\lan{\ord}}}}{\sum_{j_1=1}^{k_1}\cdots\sum_{j_\ord=1}^{k_\ord} \frac{c_{j_1,\cdots,j_\ord}}{\pare{\var{1}-\lan{1}}\cdots\pare{\var{\ord}-\lan{\ord}}}}\\
&=& \dfrac{\sum_{j_1=1}^{k_1}\cdots\sum_{j_\ord=1}^{k_\ord} \frac{\bC^{\otimes}_{k_1,\cdots,k_\ord}(j_1,\cdots,j_\ord)\bW^{\otimes}_{k_1,\cdots,k_\ord}(j_1,\cdots,j_\ord)}{\pare{\var{1}-\lan{1}}\cdots\pare{\var{\ord}-\lan{\ord}}}}{\sum_{j_1=1}^{k_1}\cdots\sum_{j_\ord=1}^{k_\ord} \frac{\bC^{\otimes}_{k_1,\cdots,k_\ord}(j_1,\cdots,j_\ord)}{\pare{\var{1}-\lan{1}}\cdots\pare{\var{\ord}-\lan{\ord}}}},
\end{array}
\end{equation}
interpolates the $\ord$-D data tensor and eventually reveals the original underlying function $\bH$ (if rational). Reducing either $k_l$, or directly $K$ (the null space entries), reduces the complexity of $\bG_{\mathrm{lag}}$ and leads to tensor approximation. Notice that $\bC^{\otimes}_{k_1,\cdots,k_\ord}(j_1,\cdots,j_\ord)$ and $\bW^{\otimes}_{k_1,\cdots,k_\ord}(j_1,\cdots,j_\ord)$ respectively denote the $(j_1,\cdots,j_\ord)$-th element of the tensorized form of the vectors $\bc_\ord$ and $\bw$.

%%%%%%%%%%%%%%%%%%%%%%%%%%%%%%%%%%%%%%%
\subsection{Decoupling and taming the curse of dimensionality}
\label{ssec:decoupling}

\subsubsection{Decoupling of variables} 

Following \cite[Theorem 5.8]{Antoulas:2025},  Theorem \ref{thm:cod} describes how the $\ord$-D Loewner null space $\bc_\ord$ can be expressed as a linear combination of a 1-D Loewner matrix null space and $k_1$, $(\ord-1)$-D Loewner null spaces. 

\begin{theorem}\label{thm:cod}
Given the data $P_c^{(\ord)}$ and $P_r^{(\ord)}$  in response of the $\ord$-variable $\bH(\var{1},\cdots,\var{\ord})$ function, the null space $\bc_\ord$  of the corresponding $\ord$-D Loenwer matrix $\IL_\ord$, is spanned by
$$
    \cN(\IL_\ord)=
    \bvec \left[
    {\bc_{\ord-1}^{\lani{1}{1}}}\cdot\left[ \bc_1^{(\lani{2}{k_2},\lani{3}{k_3},\cdots,\lani{\ord}{k_\ord})}\right]_1,
    \cdots,
    {\bc_{\ord-1}^{\lani{1}{k_1}}}\cdot\left[\bc_1^{(\lani{2}{k_2},\lani{3}{k_3},\cdots,\lani{\ord}{k_\ord})}\right]_{k_1}
    \right],
$$
where (i) $\bc_1^{(\lani{2}{k_2},\lani{3}{k_3},\cdots,\lani{\ord}{k_\ord})}$ spans $\cN(\IL_1^{(\lani{2}{k_2},\lani{3}{k_3},\cdots,\lani{\ord}{k_\ord})})$, i.e. the null space of the {1-D Loewner matrix} for frozen $\{\lani{2}{k_2},\lani{3}{k_3},\cdots,\lani{\ord}{k_\ord}\}$, and (ii) {$c_{\ord-1}^{\lani{1}{j_1}}$} spans  $\cN(\IL_{\ord-1}^{\lani{1}{j_1}})$, i.e. the $j_1$-th null space of the {$(\ord-1)$-D Loewner matrix} for frozen $\var{1}=\{ \lani{1}{1},\cdots,\lani{1}{k_1}\}$.
\end{theorem}

A direct consequence of Theorem \ref{thm:cod} is summarized in the following decoupling Theorem \ref{thm:decoupling} \cite[Theorem 5.9]{Antoulas:2025}.
\begin{theorem} \label{thm:decoupling}
Given the data $P_c^{(\ord)}$ and $P_r^{(\ord)}$and Theorem \ref{thm:cod}, the latter achieves variable decoupling, and the null space $\cN(\IL_\ord)$ can be equivalently written/spanned by $\bc_\ord$ as:
    \begin{equation} \label{eq:decoupling}
        \bc_\ord = \underbrace{\bc^{\var{\ord}}}_{{\bf Bary}({\var{\ord})}} \odot 
         \underbrace{(\bc^{\var{\ord-1}} \otimes {\bone}_{k_\ord} )}_{{\bf Bary}({\var{\ord-1}})} \odot
        \underbrace{(\bc^{\var{\ord-2}} \otimes {\bone}_{k_\ord k_{\ord-1}} )}_{{\bf Bary}({\var{\ord-2}})} \odot 
        \cdots \odot 
        \underbrace{(\bc^{\var{1}} \otimes {\bone}_{k_{\ord}\dots k_2} )}_{{\bf Bary}({\var{1}})},
    \end{equation}
where $\bc^{\var{l}}$ denotes the vectorized barycentric coefficients related to the $l$-th variable.
\end{theorem}

As an illustration, in Theorem \ref{thm:decoupling}, $\bc^{\var{1}}=\bc_1^{(\lani{2}{k_2},\lani{3}{k_3},\cdots,\lani{\ord}{k_\ord})}$ while $\bc^{\var{2}}$ is the vectorized collection of $k_1$ vectors $\bc_1^{(\lani{1}{1},\lani{3}{k_3},\cdots,\lani{\ord}{k_\ord})}, \cdots, \bc_1^{(\lani{1}{k_1},\lani{3}{k_3},\cdots,\lani{\ord}{k_\ord})}$ and so on. Theorem \ref{thm:cod} and Theorem \ref{thm:decoupling} then suggest a recursive (or cascaded) scheme, where a collection of univariate null space computations (of small size) is needed, instead of one multivariate, large-scale computation. Next, we assess how much this contributes to taming the {\bf C-o-D}, both in terms of \flop \ and memory savings. For details and examples, we refer the reader to \cite[Section 5]{Antoulas:2025} and Section \ref{ssec:mlf}.

\subsubsection{Computational effort} 

Computing the null space vector in \eqref{eq:null} for $\IL_n\in\IC^{Q\times K}$ may be practically performed by means of an \texttt{SVD}. By noticing that $\row{Q}\times \col{K}$ matrix \texttt{SVD}  \flop \, estimation is $QK^2$ (if $Q>K$) or, in the most favorable case $N^3$ (if $Q=K=N$), the complexity curve of $\mathcal O(N^3)$ will limit the utilization of the method. By recursively applying the result in Theorem \ref{thm:cod} (or equivalently Theorem \ref{thm:decoupling}), it follows that the null space corresponding to a $\ord$-D Loewner matrix can be obtained only by means of 1-D Loewner matrices null space computations. See \cite[Theorem 5.10 \& Corollary 5.11]{Antoulas:2025} for details, proofs and didactic examples. The direct consequence in terms of \flop\, complexity is stated as follows.

\begin{theorem}\label{thm:complexity}
    The  \flop \ number for the recursive approach given in  Theorem \ref{thm:cod}, is:
    $$
    \text{\flop}_{1} = \displaystyle\sum_{l=1}^\ord \left( k_l^3 \prod_{j=1}^{l} k_{j-1}\right)\text{ where } k_0= 1.
    $$
\end{theorem}

\begin{figure}[h]
\centering
\includegraphics[width=.8\linewidth]{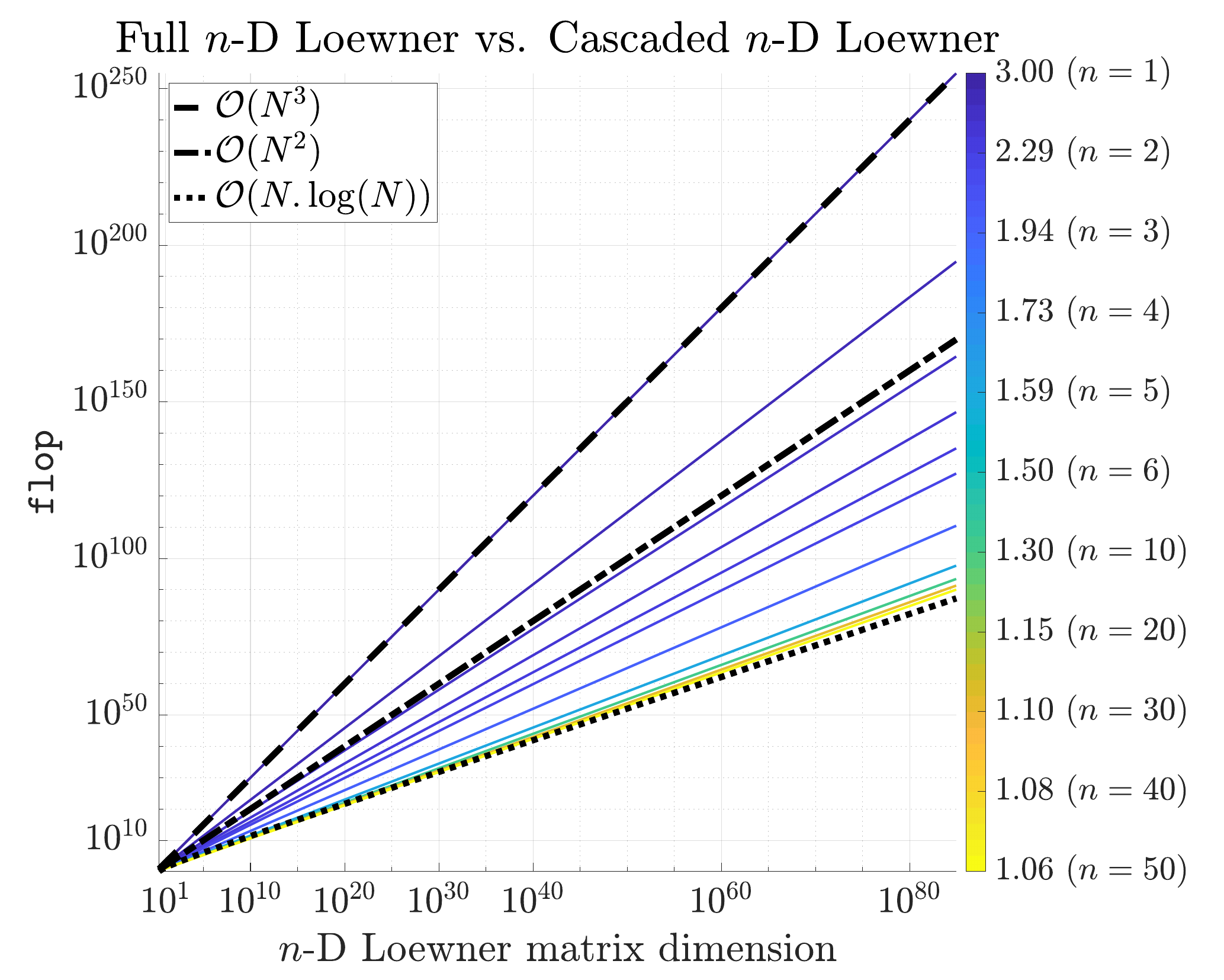}
 \caption{\flop \ comparison: cascaded/recursive $\ord$-D Loewner worst-case upper bounds $\overline{\text{\flop}_{1}}$ \eqref{eq:flop_worst} for varying number of variables $\ord$, while the full $\ord$-D Loewner is ${\cal O}(N^3)$ (black dashed); comparison with ${\cal O}(N^2)$ and ${\cal O}(N \log(N))$ references are shown in dash-dotted and dotted black lines.}
 \label{fig:flop}
\end{figure}

\begin{corr}
The most computationally demanding configuration occurs when all variables $\var{l}$ are of the same order $d_l=k_l-1=k-1$ ($l=1,\ldots,\ord$), thus requiring $k$ interpolation points each. In this configuration, the worst case \texttt{flop} is (note that $N=k^\ord$)
    \begin{equation}\label{eq:flop_worst}
    \overline{\text{\flop}_{1}}=\overbrace{k^3+k^4+\cdots+k^{n+2}}^{\text{$n$ terms}} 
    = k^3\dfrac{1-k^{n}}{1-k}=k^3\dfrac{1-N}{1-k},
    \end{equation}
    which is a ($n$ finite) geometric series of ratio $k$. 
\end{corr}
Consequently, an upper bound of \eqref{eq:flop_worst} can be estimated by assuming that $k>1$ and for a different number of variables $\ord$. As an example, the complexity is upper bounded by ${\cal O}(N^{3})$ for $\ord=1$, ${\cal O}(N^{2.29})$ for $\ord=2$, ${\cal O}(N^{1.94})$ for $\ord=3$, ${\cal O}(N^{1.73})$ for $\ord=4$ and already ${\cal O}(N^{1.5})$ for $\ord=6$. One can clearly observe that when the number of variables $\ord>1$, the \flop \, complexity drops, and this decreases as $\ord$ increases; e.g. for $\ord=50$ one gets ${\cal O}(N^{1.06})$. This is illustrated in Figure \ref{fig:flop} which shows the worst-case $\overline{\flop_1}$ \ as a function of $\ord$, and compares it to classical complexity references (notice that standard Loewner approaches is ${\cal O}(N^3)$).

%\noindent\begin{minipage}{.45\textwidth}
%\begin{result}
    %Following \ref{eq:flop_worst} the \flop \, complexity, is upper bounded as $\mathcal{O}(N^3)$ for $\ord=1$, $\mathcal{O}(N^{2.29})$ for $\ord=2$, $\mathcal{O}(N^{1.94})$ for $\ord=3$, \dots, $\mathcal{O}(N^{1.5})$  for $\ord=6$, \dots up to $\mathcal{O}(N^{1.06})$ for $\ord=50$. %Right frame provides \flop \ comparison: Cascaded $\ord$-D Loewner worst-case upper bounds for varying number of variables $\ord$, while the full $\ord$-D Loewner is ${\cal O}(N^3)$; comparison with ${\cal O}(N^2)$ and ${\cal O}(N \log(N))$ references.
%\end{result}
%\end{minipage}
%\hspace*{2mm}\begin{minipage}{.45\textwidth}
%    \includegraphics[width=\linewidth]{figures/charles/complexityGridWorst.pdf}
%\end{minipage}

\subsubsection{Storage effort} 

Second, and equally important as the computational burden, the storage needed for the $\IL_\ord$ matrix of dimension $\row{Q}\times \col{K}$ is $\dfrac{8}{2^{20}} QK$ \texttt{MB}, or simply $\dfrac{8}{2^{20}} N^2$ \texttt{MB} (if equal number of columns and rows are considered). Then, the following holds \cite[Theorem 5.13]{Antoulas:2025}.

\begin{theorem} \label{thm:memory}
Following Theorem \ref{thm:cod} \& \ref{thm:decoupling}, one only needs to sequentially construct single 1-D Loewner matrices, each of dimension $\IL_1 \in \IC^{k_l \times  k_l}$.  The largest stored matrix is $\IL_1 \in \IC^{\overline k\times \overline k}$, where $\overline k=\max_{l} k_l$ ($l=1,\cdots,\ord$). In complex and double precision, the maximum disk storage is $\dfrac{8}{2^{20}} \overline k^2$ \texttt{MB}.
\end{theorem}

%\newpage 
\subsubsection{Summary: full vs. cascaded  (with an example)}

The data  storage for $\tableau{\ord}$ (in complex arithmetic and double precision) is given by
$\dfrac{8}{2^{20}} \prod_{l=1}^\ord (q_l+k_l) ~~\texttt{MB}$. For example $\tableau{6}$ of dimension $2 \cdot [20,6,4,6,8,2]=2 \cdot [k_1,k_2,k_3,k_4,k_5,k_6]$ needs 45 \texttt{MB} (we assume $q_l=k_l$). Then, according to the full or the cascaded null space computation version, the following values hold.\\~~\\

\noindent\begin{minipage}[t]{.45\textwidth}
\begin{center}
    \textbf{Full $\ord$-D Loewner}\\
\end{center}
\begin{itemize}
    \item The construction of the $\IL_\ord \in \IC^{N\times N}$ matrix, where $N=K=Q$, needs
    $$
    \dfrac{8}{2^{20}} N^2 =  \dfrac{8}{2^{20}} 46,080^2~~\texttt{MB} 
    $$
    being $31.64$ \texttt{GB}.
    \item The required \flop\, is $N^3$, being $9.78\cdot 10^{13}$ \flop\, in our example.
\end{itemize}
\end{minipage}
\hspace*{2mm}\begin{minipage}[t]{.45\textwidth}
\begin{center}
    \textbf{Cascaded $\ord$-D Loewner}\\
\end{center}
\begin{itemize}
    \item The construction of the $\IL_1 \in \IC^{\overline k\times \overline k}$ matrix, where $\overline k=\max_{l} k_l$, needs
    $$
    \dfrac{8}{2^{20}} \overline k^2=\dfrac{8}{2^{20}}  20^2 ~~\texttt{MB} 
    $$
    being $6.25$ \texttt{KB}.
    \item The required \flop \, is $\flop_1$ as in Theorem \ref{thm:complexity}, being $1.78 \cdot 10^{6}$ \flop \, in our example\footnote{Note that \flop \, may be decreased to $8.13 \cdot 10^{5}$ \flop\, when variables optimally reordered, see \cite{Antoulas:2025}.}.
    \end{itemize}
\end{minipage}

%%%%%%%%%%%%%%%%%%%%%%%%%%%%%%%%%%%%%%%%%%%%%
\subsection{Simple detailed example with \texttt{MATLAB} codes}
\label{ssec:mlf}

Let us now illustrate with a very simple example, how to deploy the recursive null space computation scheme presented in this section and detailed in \cite[Alg. 1]{Antoulas:2025}. This is exemplified with the \texttt{MATLAB} package \texttt{+mLF} (multivariate Loewner Framework), available as a \href{https://github.com/cpoussot/mLF}{\texttt{GitHub}}  repository at the following link: \url{https://github.com/cpoussot/mLF}. 

\paragraph{Install \texttt{mLF}.}  

First clone the \texttt{GitHub} repository in the directory of your choice as follows (open a command line interface):
\definecolor{backcolour}{rgb}{0.95,0.95,0.92}
\begin{lstlisting}[basicstyle=\small]
cd "directory_for_mLF"
git clone git@github.com:cpoussot/mLF.git
\end{lstlisting}

Now in the \texttt{MATLAB} software, add the path of the cloned repository as follows, and  start using the available features:
\begin{lstlisting}[basicstyle=\small]%[style=Matlab-editor,basicstyle=\small]
addpath("directory_for_mLF") 
\end{lstlisting}

\paragraph{Define the problem and construct the tensor.}  We start the illustration with a two-variable function 
$$
\bH(\var{1},\var{2}) = \var{1}\,{\var{2}}^3+2\,\var{1}\,\var{2}-1

$$ 
i.e. $\ord=2$, and $l=\{1,2\}$, together with its bounds ${\cal X}_l:=\left[\begin{array}{cc} -1 & 1 \end{array}\right]
$ and linear discretization mesh $N_l=10
$, leading to $\bvar{l}= \left[\begin{array}{cccccccccc} -1 & -\frac{7}{9} & -\frac{5}{9} & -\frac{1}{3} & -\frac{1}{9} & \frac{1}{9} & \frac{1}{3} & \frac{5}{9} & \frac{7}{9} & 1 \end{array}\right]
 \in \mathcal{X}_{l}^{10
}$.
%\lstinputlisting[style=Matlab-editor,basicstyle=\small,backgroundcolor=\color{backcolour},firstline=9,lastline=18]{figures/code2vars_01/demo_01_2.m}
\begin{lstlisting}[basicstyle=\small,firstline=9,lastline=18]
%%% Define a multivariate problem 
syms x1 x2; 
Hsym    = x1*x2^3+2*x1*x2-1; 
Hf      = matlabFunction(Hsym);
H       = @(x) Hf(x(:,1),x(:,2));
n       = 2; % number of variables
xbnd    = [-1 1]; % bounds of both variables
Nip     = 10; % number of interpolation points
X1      = linspace(xbnd(1),xbnd(2),Nip);
X2      = linspace(xbnd(1),xbnd(2),Nip);
\end{lstlisting}

From the above setup, we select as columns $\lan{l}$ (denoted \texttt{p\_c\{1\}} and  \texttt{p\_c\{2\}}) and rows $\mun{l}$  (denoted \texttt{p\_r\{1\}} and  \texttt{p\_r\{2\}}) interpolatory set points, a subset of $\mathcal{X}_{1}^{10
}$ and $\mathcal{X}_{2}^{10
}$; in what follows, we chose equal dimensions (i.e. $k_l=5
$ and $q_l=5
$), where distribution alternates between columns and rows (notice that theoretically any arrangement is possible): 
%\lstinputlisting[style=Matlab-editor,basicstyle=\small,backgroundcolor=\color{backcolour},firstline=20,lastline=24]{figures/code2vars_01/demo_01_2.m}
\begin{lstlisting}[basicstyle=\small,firstline=20,lastline=24]
%%% Interpolation points 
p_c{1}  = X1(2:2:end);
p_r{1}  = X1(1:2:end);
p_c{2}  = X2(2:2:end);
p_r{2}  = X2(1:2:end);
\end{lstlisting}

After checking that  the interpolation points are disjoints, one may construct the 2-D tensor $\tableau{2}\in\IR^{(q_1+q_2)\times (k_1+k_2)}$.
%\lstinputlisting[style=Matlab-editor,basicstyle=\small,backgroundcolor=\color{backcolour},firstline=26,lastline=31]{figures/code2vars_01/demo_01_2.m}
\begin{lstlisting}[basicstyle=\small,firstline=26,lastline=31]
%%% Construct tensor
% Check that column/row IP are disjoint 
ok          = mlf.check(p_c,p_r);
% Construct tab_n 
[y,x,dim]   = mlf.make_tab_vec(H,p_c,p_r);
tab         = mlf.vec2mat(y,dim); 
\end{lstlisting}

\paragraph{Compute orders \& barycentric coefficients.} 

\begin{figure}[h]
\centering
\includegraphics[width=.9\textwidth]{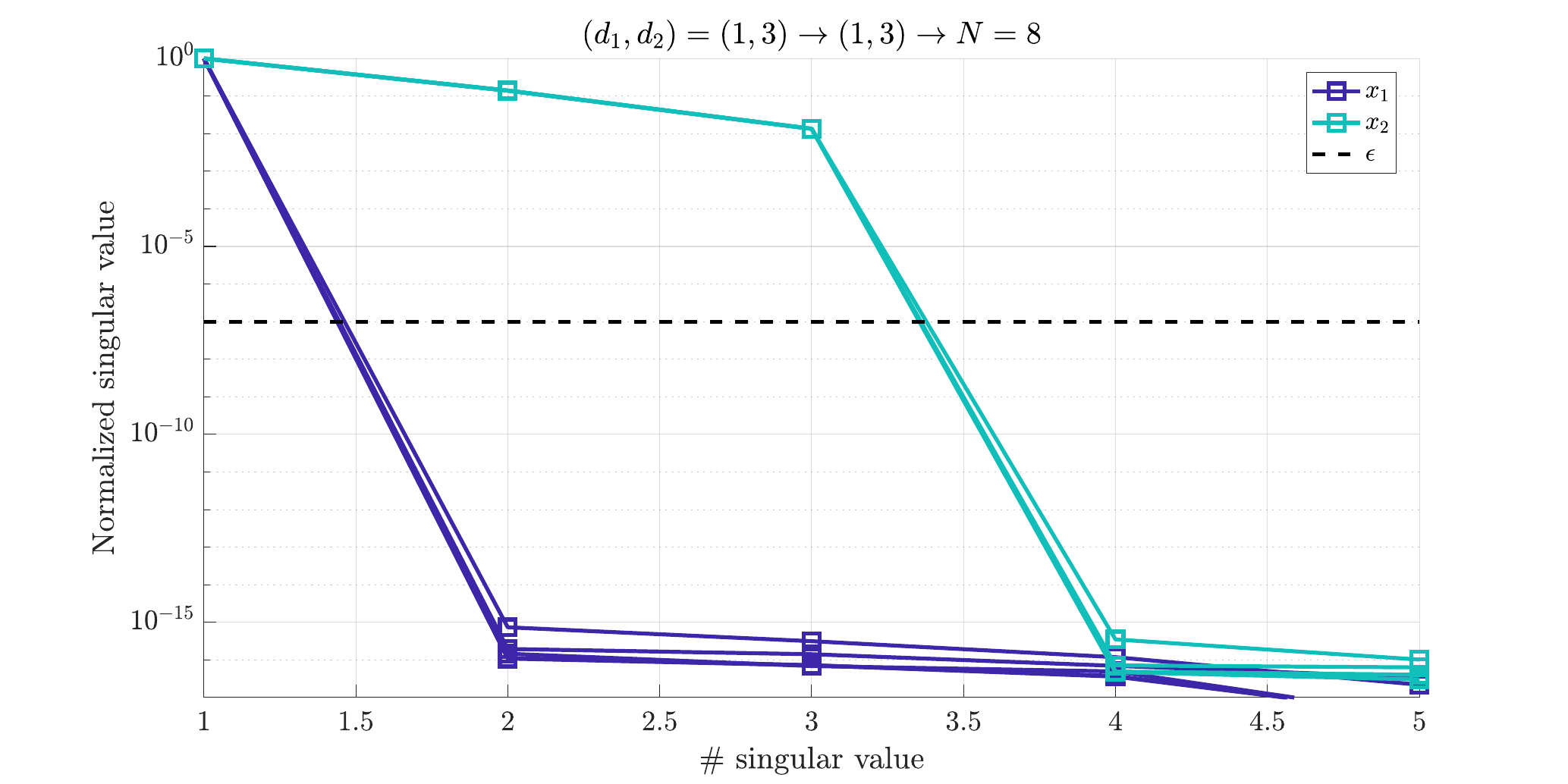}
\caption{Single variable Loewner matrices normalized singular values (order detection).}
\label{fig:matlab1}
\end{figure}

Following \cite[Alg. 1]{Antoulas:2025}, it is possible to estimate the order along each variables, then leading to Figure \ref{fig:matlab1}. 
%\lstinputlisting[style=Matlab-editor,basicstyle=\small,backgroundcolor=\color{backcolour},firstline=33,lastline=37]{figures/code2vars_01/demo_01_2.m}
\begin{lstlisting}[basicstyle=\small,firstline=33,lastline=37]
%%% Estimate orders along each variables and select a subset of IP
tol              = 1e-7; 
ord              = mlf.compute_order(p_c,p_r,tab,tol,[],5,true); 
[pc,pr,W,V,tabr] = mlf.points_selection(p_c,p_r,tab,ord,true);
w                = mlf.mat2vec(W);
\end{lstlisting}

Figure \ref{fig:matlab1} shows the normalized singular values of the single-variable Loewner matrices. The interpretation of  Figure \ref{fig:matlab1} concerns the dimension estimation along each variables. In this very simple case, $\var{1}$ is of degree $d_1=1
$ while  $\var{2}$ is of degree $d_2=3
$, thus one needs at least $(k_1,k_2)=(2
,4
)$ column interpolation points. This implies a barycentric vector of dimension $2
\times 4
=8
=N$. Adding more interpolation points would lead to \textbf{overfitting}. The selected  interpolation points are $\lan{1}=\left[\begin{array}{cc} -\frac{7}{9} & 1 \end{array}\right]
$ (denoted \texttt{pc\{1\}}) and $\lan{2}=\left[\begin{array}{cccc} -\frac{7}{9} & -\frac{1}{3} & \frac{1}{9} & 1 \end{array}\right]
$ (denoted \texttt{pc\{2\}}). The row interpolation points are  $\mun{1}=\left[\begin{array}{cc} -1 & \frac{7}{9} \end{array}\right]
$ (denoted \texttt{pr\{1\}}) and $\mun{2}=\left[\begin{array}{cccc} -1 & -\frac{5}{9} & -\frac{1}{9} & \frac{7}{9} \end{array}\right]
$ (denoted \texttt{pr\{2\}})\footnote{Notice that one may also chose as row interpolation point, all the remaining interpolation points (i.e. all points but the $\lan{l}$). This is the idea of \texttt{AAA} \cite{nakatsukasa2020algorithm} and \texttt{pAAA} \cite{Carracedo:2023,Balicki:2025}.}. Evaluating $\bH(\var{1},\var{2})$ at these combinations, then leads to the interpolation values $\bw$ (denoted \texttt{w}), and tensorized version $\bW^{\otimes}_{2
,4
}$ (denoted \texttt{W}). 
$$
\bw=\bH\left(\left[\begin{array}{cc} -\frac{7}{9} & -\frac{7}{9}\\[1mm] -\frac{7}{9} & -\frac{1}{3}\\[1mm] -\frac{7}{9} & \frac{1}{9}\\[1mm] -\frac{7}{9} & 1\\[1mm] 1 & -\frac{7}{9}\\[1mm] 1 & -\frac{1}{3}\\[1mm] 1 & \frac{1}{9}\\[1mm] 1 & 1 \end{array}\right]\right)
=\left[\begin{array}{c} \frac{3778}{6561}\\[1mm] -\frac{110}{243}\\[1mm] -\frac{7702}{6561}\\[1mm] -\frac{10}{3}\\[1mm] -\frac{2206}{729}\\[1mm] -\frac{46}{27}\\[1mm] -\frac{566}{729}\\[1mm] 2 \end{array}\right]
  \text{ and } \bW^{\otimes}_{2
,4
}=\left[\begin{array}{cccc} \frac{3778}{6561} & -\frac{110}{243} & -\frac{7702}{6561} & -\frac{10}{3}\\[1mm] -\frac{2206}{729} & -\frac{46}{27} & -\frac{566}{729} & 2 \end{array}\right]
.
$$

Now, following \cite[Theorem. 5.8 and Theorem. 5.9]{Antoulas:2025} (or Theorem \ref{thm:cod} and Theorem \ref{thm:decoupling}), one may compute the barycentric coefficients $\bc_2$ (denoted \texttt{c}) using the recursive scheme as follows:
%\lstinputlisting[style=Matlab-editor,basicstyle=\small,backgroundcolor=\color{backcolour},firstline=40,lastline=41]{figures/code2vars_01/demo_01_2.m}
\begin{lstlisting}[basicstyle=\small,firstline=40,lastline=41]
%%% Compute rec. LL null
[c,lag] = mlf.loewner_null_rec(pc,pr,tabr,'svd0');
\end{lstlisting}

This leads to the null space vector $\bc_2= \left[\begin{array}{cccccccc} 3 & -8 & 6 & -1 & -3 & 8 & -6 & 1 \end{array}\right]
^\top$ of the 2-D Loewner matrix, obtained without  constructing it. Then, following \cite[Theorem. 5.10 and Theorem. 5.13]{Antoulas:2025} (or Theorem \ref{thm:complexity} and Theorem \ref{thm:memory}), one may evaluate both the computational and storage efforts. Applying these  theorems lead to 136
 \flop\, and  16
 \bytes\, of maximal matrix storage (for the recursive approach, gathered in the \texttt{lag} variable) while one would need  512
 \flop\, and  64
 \bytes\, of maximal matrix storage (for the full approach). These numbers clearly show how much, even in such a very simple configuration, the complexity and computational load is reduced. This is the basis for our \textbf{claim of taming the curse of dimensionality}. 

\paragraph{Lagrangian basis.}

Then, based on the interpolation points $\{\lan{1},\lan{2}\}$, associated values $\bw$, and computed  null space vector (or barycentric, Lagrangian weight) $\bc_2$, the multivariate function in the Lagrangian basis $\bG_{\mathrm{lag}}$ (denoted \texttt{glag}) is obtained as follows:
%\lstinputlisting[style=Matlab-editor,basicstyle=\small,backgroundcolor=\color{backcolour},firstline=45,lastline=46]{figures/code2vars_01/demo_01_2.m}
\begin{lstlisting}[basicstyle=\small,firstline=45,lastline=46]
%%% Transfer function in Lagrangian form
[glag,ilag] = mlf.tf_lagrangian(pc,w,c,false);
\end{lstlisting}
where \texttt{ilag} gathers the basis, numerator and denominator coefficients:
$$
\left(\begin{array}{c||c|c}\mathcal{B}_\textrm{lag}(\var{1},\var{2}) & \bN_\textrm{lag} &\bD_\textrm{lag}\end{array}\right)=\left(\begin{array}{c||c|c} \left(\var{1}+\frac{7}{9}\right)\,\left(\var{2}+\frac{7}{9}\right) & \frac{3778}{2187} & 3\\[1mm] \left(\var{2}+\frac{1}{3}\right)\,\left(\var{1}+\frac{7}{9}\right) & \frac{880}{243} & -8\\[1mm] \left(\var{2}-\frac{1}{9}\right)\,\left(\var{1}+\frac{7}{9}\right) & -\frac{15404}{2187} & 6\\[1mm] \left(\var{2}-1\right)\,\left(\var{1}+\frac{7}{9}\right) & \frac{10}{3} & -1\\[1mm] \left(\var{1}-1\right)\,\left(\var{2}+\frac{7}{9}\right) & \frac{2206}{243} & -3\\[1mm] \left(\var{1}-1\right)\,\left(\var{2}+\frac{1}{3}\right) & -\frac{368}{27} & 8\\[1mm] \left(\var{1}-1\right)\,\left(\var{2}-\frac{1}{9}\right) & \frac{1132}{243} & -6\\[1mm] \left(\var{1}-1\right)\,\left(\var{2}-1\right) & 2 & 1 \end{array}\right)
.
$$

The above code block also provides the basis $\mathcal{B}_\textrm{lag}(\var{1},\var{2})$, the numerator  $ \bN_\textrm{lag}$ and denominator $ \bD_\textrm{lag}$ values, providing the result
$$
\begin{array}{rcl}
\bG_{\mathrm{lag}}(\var{1},\var{2})
=\dfrac{\bn_{\mathrm{lag}}(\var{1},\var{2})}{\bd_{\mathrm{lag}}(\var{1},\var{2})}
&=&\dfrac{\sum_{\textrm{row}} \bN_\textrm{lag} \odot\mathcal{B}^{-1}_\textrm{lag}(\var{1},\var{2})}{\sum_{\textrm{row}} \bD_\textrm{lag} \odot\mathcal{B}^{-1}_\textrm{lag}(\var{1},\var{2})} \\
&&\\
&=&\dfrac{\sum_{j_1=1}^{2
}\sum_{j_\ord=1}^{4
} \frac{\bC^{\otimes}_{2
,4
}(j_1,j_2)\bW^{\otimes}_{2
,4
}(j_1,j_2)}{\pare{\var{1}-\lan{1}}\pare{\var{2}-\lan{2}}}}{\sum_{j_1=1}^{2
}\sum_{j_2=1}^{4
} \frac{\bC^{\otimes}_{2
,4
}(j_1,j_2)}{\pare{\var{1}-\lan{1}}\pare{\var{2}-\lan{2}}}},
\end{array}
$$
where, $\lan{1}=\left[\begin{array}{cc} -\frac{7}{9} & 1 \end{array}\right]
$ and $\lan{2}=\left[\begin{array}{cccc} -\frac{7}{9} & -\frac{1}{3} & \frac{1}{9} & 1 \end{array}\right]
$  (with $j_1=1,\cdots,2
$ and $j_2=1,\cdots,4
$) and where the tensorized information read
$$
\bC^{\otimes}_{2
,4
} =\left[\begin{array}{cccc} 3 & -8 & 6 & -1\\[1mm] -3 & 8 & -6 & 1 \end{array}\right]
 \text{ and } \bW^{\otimes}_{2
,4
}=\left[\begin{array}{cccc} \frac{3778}{6561} & -\frac{110}{243} & -\frac{7702}{6561} & -\frac{10}{3}\\[1mm] -\frac{2206}{729} & -\frac{46}{27} & -\frac{566}{729} & 2 \end{array}\right]
.
$$
%$$
%\mathcal B^{\otimes}_{{\mathrm{lag}}_{\input{figures/code2vars_01/k1},\input{figures/code2vars_01/k2}}}=\input{figures/code2vars_01/Blag}.
%$$

%\input{figures/code2vars_01/NNnum}
At this stage, one can also view connections of the multivariate Lagrangian basis with the KAN (Kolmogorov Neural Network), see also \cite{Poluektov:2025,Liu:2025} (where activation function is in a Lagrangian form).
\begin{figure}[h]\begin{center} \scalebox{.7}{\begin{tikzpicture}[line width=0.4mm]\tikzstyle{place}=[circle, draw=black, minimum size = 8mm]\tikzstyle{placeInOut}=[circle, draw=orange, minimum size = 8mm]\node at (0,-2) [placeInOut] (first_1){$\var{1}$};\node at (0,-4) [placeInOut] (first_2){$\var{2}$};\node at (5,-2) [place] (secondL1_1){$\frac{1}{\var{1}-\lani{1}{1}}$};\node at (5,-4) [place] (secondL1_2){$\frac{1}{\var{1}-\lani{1}{2}}$};\node at (5,-6) [place] (secondL2_1){$\frac{1}{\var{2}-\lani{2}{1}}$};\node at (5,-8) [place] (secondL2_2){$\frac{1}{\var{2}-\lani{2}{2}}$};\node at (5,-10) [place] (secondL2_3){$\frac{1}{\var{2}-\lani{2}{3}}$};\node at (5,-12) [place] (secondL2_4){$\frac{1}{\var{2}-\lani{2}{4}}$};\node at (10,-2) [place] (third_1){$\prod$};\node at (10,-4) [place] (third_2){$\prod$};\node at (10,-6) [place] (third_3){$\prod$};\node at (10,-8) [place] (third_4){$\prod$};\node at (10,-10) [place] (third_5){$\prod$};\node at (10,-12) [place] (third_6){$\prod$};\node at (10,-14) [place] (third_7){$\prod$};\node at (10,-16) [place] (third_8){$\prod$};\node at (15,-9) [placeInOut] (output){$\bSigma$};\draw[->] (first_1)--(secondL1_1) node[above,sloped,pos=0.75] { };\draw[->] (first_1)--(secondL1_2) node[above,sloped,pos=0.75] { };\draw[->] (first_2)--(secondL2_1) node[above,sloped,pos=0.75] { };\draw[->] (first_2)--(secondL2_2) node[above,sloped,pos=0.75] { };\draw[->] (first_2)--(secondL2_3) node[above,sloped,pos=0.75] { };\draw[->] (first_2)--(secondL2_4) node[above,sloped,pos=0.75] { };\draw[->] (secondL1_1)--(third_1) node[above,sloped,pos=0.25] {};\draw[->] (secondL1_1)--(third_2) node[above,sloped,pos=0.25] {};\draw[->] (secondL1_1)--(third_3) node[above,sloped,pos=0.25] {};\draw[->] (secondL1_1)--(third_4) node[above,sloped,pos=0.25] {};\draw[->] (secondL1_2)--(third_5) node[above,sloped,pos=0.25] {};\draw[->] (secondL1_2)--(third_6) node[above,sloped,pos=0.25] {};\draw[->] (secondL1_2)--(third_7) node[above,sloped,pos=0.25] {};\draw[->] (secondL1_2)--(third_8) node[above,sloped,pos=0.25] {};\draw[->] (secondL2_1)--(third_1) node[above,sloped,pos=0.25] {};\draw[->] (secondL2_2)--(third_2) node[above,sloped,pos=0.25] {};\draw[->] (secondL2_3)--(third_3) node[above,sloped,pos=0.25] {};\draw[->] (secondL2_4)--(third_4) node[above,sloped,pos=0.25] {};\draw[->] (secondL2_1)--(third_5) node[above,sloped,pos=0.25] {};\draw[->] (secondL2_2)--(third_6) node[above,sloped,pos=0.25] {};\draw[->] (secondL2_3)--(third_7) node[above,sloped,pos=0.25] {};\draw[->] (secondL2_4)--(third_8) node[above,sloped,pos=0.25] {};\draw[->] (third_1)--(output) node[above,sloped,pos=0.25] {3};\draw[->] (third_2)--(output) node[above,sloped,pos=0.25] {-8};\draw[->] (third_3)--(output) node[above,sloped,pos=0.25] {6};\draw[->] (third_4)--(output) node[above,sloped,pos=0.25] {-1};\draw[->] (third_5)--(output) node[above,sloped,pos=0.25] {-3};\draw[->] (third_6)--(output) node[above,sloped,pos=0.25] {8};\draw[->] (third_7)--(output) node[above,sloped,pos=0.25] {-6};\draw[->] (third_8)--(output) node[above,sloped,pos=0.25] {1};\end{tikzpicture}} \caption{Equivalent NN representation of the denominator $\bd_{\textrm{lag}}$.}\end{center}\end{figure}

\paragraph{Monomial basis}

Introducing the projection Vandermonde matrices $\bV=\bV_{\var{1}} \otimes \bV_{\var{2}}$, where
$$
\bV_{\var{1}} = \left[\begin{array}{cc} 1 & 1\\[1mm] -1 & \frac{7}{9} \end{array}\right]
 \text{ and }
\bV_{\var{2}} = \left[\begin{array}{cccc} 1 & 1 & 1 & 1\\[1mm] -\frac{7}{9} & -\frac{1}{3} & \frac{1}{9} & 1\\[1mm] -\frac{7}{27} & -\frac{61}{81} & -\frac{23}{27} & \frac{11}{81}\\[1mm] \frac{1}{27} & \frac{7}{81} & -\frac{7}{27} & -\frac{7}{243} \end{array}\right]
  
$$
one gets the monomial coefficients for the numerator as $\bV (\bw\odot \bc)$ and denominator as $\bV \bc$. Computational details on the projectors can be found in \cite[Section 2.1.a]{Gander:2005}. The multivariate model in the monomial basis $\bG_{\mathrm{mon}}$ (denoted \texttt{gmon}) is obtained as follows:
%\lstinputlisting[style=Matlab-editor,basicstyle=\small,backgroundcolor=\color{backcolour},firstline=48,lastline=49]{figures/code2vars_01/demo_01_2.m}
\begin{lstlisting}[basicstyle=\small,firstline=48,lastline=49]
%%% Transfer function in Monomial form
[gmon,imon] = mlf.tf_monomial(pc,w,c,false);
\end{lstlisting}
where \texttt{imon} gathers the basis, numerator and denominator coefficients:
$$
\left(\begin{array}{c||c|c}\mathcal{B}_\textrm{mon}(\var{1},\var{2}) & \bN_\textrm{mon} &\bD_\textrm{mon}\end{array}\right)=
\left(\begin{array}{c||c|c} \var{1}\,{\var{2}}^3 & 1 & 0\\[1mm] \var{1}\,{\var{2}}^2 & 0 & 0\\[1mm] \var{1}\,\var{2} & 2 & 0\\[1mm] \var{1} & 0 & 0\\[1mm] {\var{2}}^3 & 0 & 0\\[1mm] {\var{2}}^2 & 0 & 0\\[1mm] \var{2} & 0 & 0\\[1mm] 1 & -1 & 1 \end{array}\right)
.
$$

Similarly to the Lagrangian case, the above code block also provides the basis $\mathcal{B}_\textrm{mon}(\var{1},\var{2})$, the numerator  $ \bN_\textrm{mon}$ and denominator $ \bD_\textrm{mon}$ values, providing the result
$$
\begin{array}{rcl}
\bG_{\mathrm{mon}}(\var{1},\var{2})
=\dfrac{\bn_{\mathrm{mon}}(\var{1},\var{2})}{\bd_{\mathrm{mon}}(\var{1},\var{2})}
&=&\dfrac{\sum_{\textrm{row}} \bN_\textrm{mon} \odot \mathcal{B}_\textrm{mon}(\var{1},\var{2})}{\sum_{\textrm{row}} \bD_\textrm{mon} \odot\mathcal{B}_\textrm{mon}(\var{1},\var{2})} \\
&&\\
&=&\dfrac{\sum_{j_1=1}^{2
}\sum_{j_\ord=1}^{4
} \bN^{\otimes}_{2
,4
}(j_1,j_2)\pare{\var{1}^{(j_1-1)}\var{2}^{(j_2-1)}}}{\sum_{j_1=1}^{2
}\sum_{j_\ord=1}^{4
} \bD^{\otimes}_{2
,4
}(j_1,j_2)\pare{\var{1}^{(j_1-1)}\var{2}^{(j_2-1)}}},
\end{array}
$$
where the tensorized information reads
$$
\bN^{\otimes}_{2
,4
} =\left[\begin{array}{cccc} 1 & 0 & 2 & 0\\[1mm] 0 & 0 & 0 & -1 \end{array}\right]
 \text{ and } \bD^{\otimes}_{2
,4
}=\left[\begin{array}{cccc} 0 & 0 & 0 & 0\\[1mm] 0 & 0 & 0 & 1 \end{array}\right]
. %\text{ and } \mathcal B^{\otimes}_{{\mathrm{mon}}_{\input{figures/code2vars_01/k1},\input{figures/code2vars_01/k2}}}=\input{figures/code2vars_01/Bmon}.
$$

\paragraph{Kolmogorov Supperposition Theorem equivalent form.}

An other important contribution presented in \cite[Section 6]{Antoulas:2025} concerns the connection with the Kolmogorov Superposition Theorem (KST). From the recursive null space construction, the KST elements can be obtained as:
%\lstinputlisting[style=Matlab-editor,basicstyle=\small,backgroundcolor=\color{backcolour},firstline=51,lastline=52]{figures/code2vars_01/demo_01_2.m}
\begin{lstlisting}[basicstyle=\small,firstline=51,lastline=52]
%%% Decoupling toward KST
[Bary,Lag,Cx]  = mlf.decoupling(pc,lag);
\end{lstlisting}

Where \texttt{[Bary\{1\} Bary\{2\}]} (resp. \texttt{[Lag\{1\}, Lag\{2\}]}) gathers the barycentric coefficients (resp.  basis) along $\var{1}$ and $\var{2}$, with the following values:
$$
\left[\begin{array}{cc} \textbf{Bary}(\var{1}) & \textbf{Bary}(\var{2}) \end{array}\right] = 
\left[\begin{array}{cc} -1 & -3\\[1mm] -1 & 8\\[1mm] -1 & -6\\[1mm] -1 & 1\\[1mm] 1 & -3\\[1mm] 1 & 8\\[1mm] 1 & -6\\[1mm] 1 & 1 \end{array}\right]

\mbox{ , }
\left[\begin{array}{cc} \textbf{Lag}(\var{1}) & \textbf{Lag}(\var{2}) \end{array}\right] = 
\left[\begin{array}{cc} \frac{1}{\var{1}+\frac{7}{9}} & \frac{1}{\var{2}+\frac{7}{9}}\\[1mm] \frac{1}{\var{1}+\frac{7}{9}} & \frac{1}{\var{2}+\frac{1}{3}}\\[1mm] \frac{1}{\var{1}+\frac{7}{9}} & \frac{1}{\var{2}-\frac{1}{9}}\\[1mm] \frac{1}{\var{1}+\frac{7}{9}} & \frac{1}{\var{2}-1}\\[1mm] \frac{1}{\var{1}-1} & \frac{1}{\var{2}+\frac{7}{9}}\\[1mm] \frac{1}{\var{1}-1} & \frac{1}{\var{2}+\frac{1}{3}}\\[1mm] \frac{1}{\var{1}-1} & \frac{1}{\var{2}-\frac{1}{9}}\\[1mm] \frac{1}{\var{1}-1} & \frac{1}{\var{2}-1} \end{array}\right]
,
$$
which are obtained from Theorem \ref{thm:decoupling} as $\textbf{Bary}(\var{1}) = \bc^{\var{1}} \otimes \bone_{4
}$ and $\textbf{Bary}(\var{2}) = \mathbf{vec}(\bc^{\var{2}})$, where
$$
\bc^{\var{1}} = \left[\begin{array}{c} -1\\[1mm] 1 \end{array}\right]
 \mbox{ and }
\bc^{\var{2}} = \left[\begin{array}{cc} -3 & -3\\[1mm] 8 & 8\\[1mm] -6 & -6\\[1mm] 1 & 1 \end{array}\right]
.
$$
Then, with the above notation, one defines the following univariate vector functions 
$$
\left\{
\begin{array}{rcl}
\bPhi_1(\var{1}) &=& \textbf{Bary}(\var{1}) \odot \textbf{Lag}(\var{1})\\
\bPhi_2(\var{2}) &=& \textbf{Bary}(\var{2}) \odot \textbf{Lag}(\var{2})
\end{array}
\right. ,
$$
and the resulting KST equivalent rational interpolant is now obtained as
$$
\left.
\begin{array}{rcl}
\mathbf{n}_{\mathrm{kst}}(\var{1},\var{2}) &=& \sum_{\text{rows}} \bw \odot  \bPhi_1(\var{1}) \odot \bPhi_2(\var{2}) \\
\mathbf{d}_{\mathrm{kst}}(\var{1},\var{2}) &=& \sum_{\text{rows}} \bPhi_1(\var{1}) \odot \bPhi_2(\var{2})
\end{array}
\right\}
\Rightarrow \bG_{\mathrm{kst}}(\var{1},\var{2}) = \dfrac{\mathbf{n}_{\mathrm{kst}}(\var{1},\var{2})}{\mathbf{d}_{\mathrm{kst}}(\var{1},\var{2}) }.
$$

\begin{figure}[h]
\centering
\includegraphics[width=\textwidth]{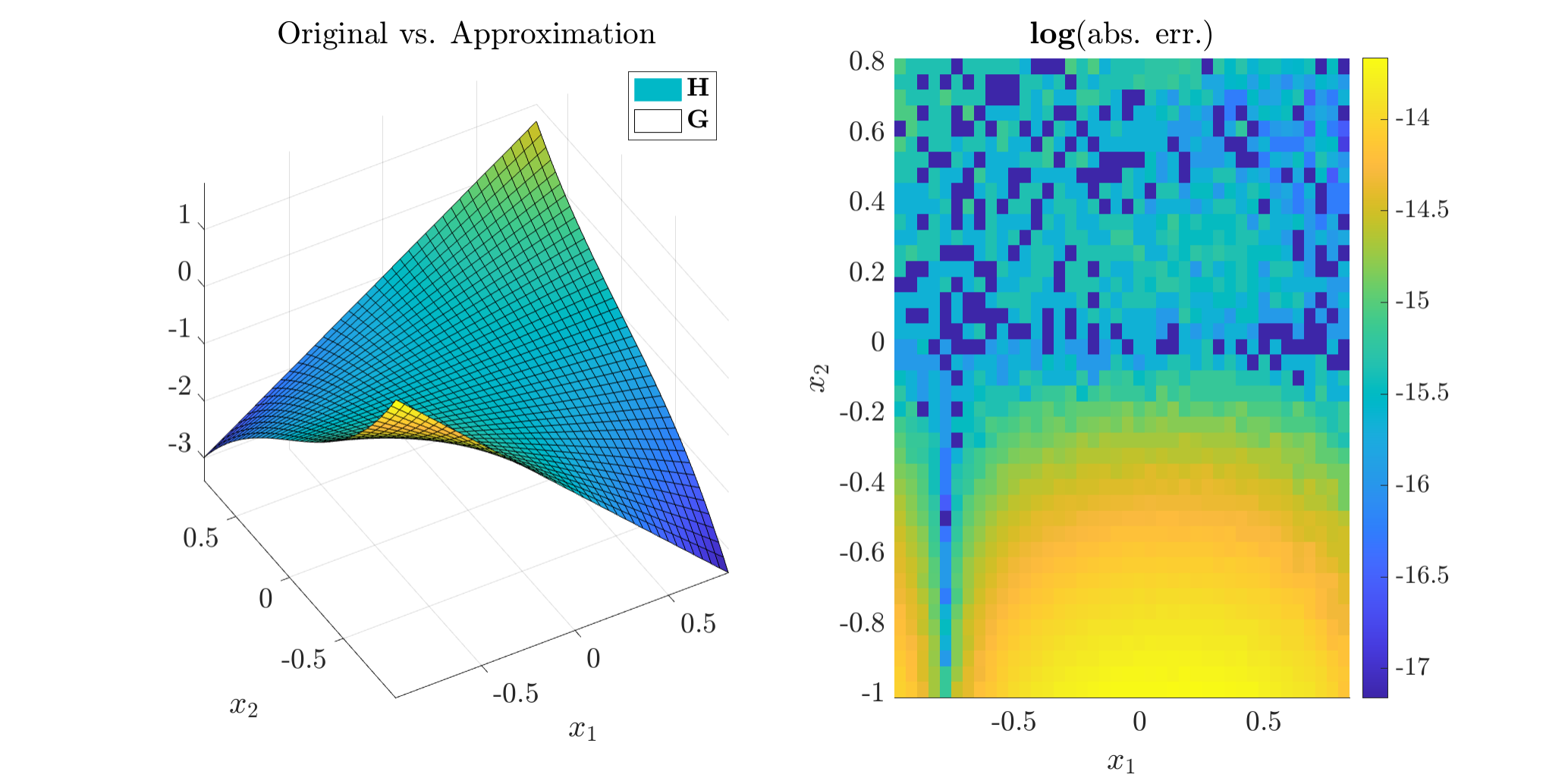}
\caption{Left: original $\bH$ (grid) and approximated $\bG_{\mathrm{lag}}$ (equivalently $\bG_{\mathrm{mon}}$ and $\bG_{\mathrm{kst}}$) (colored surface) models. Right: absolute mismatch error (in log-scale). Computation is performed in double precision.}
\label{fig:matlab2}
\end{figure}

The evaluation  of both $\bH$ and $\bG_{\mathrm{lag}}$ (equivalently $\bG_{\mathrm{mon}}$ and $\bG_{\mathrm{kst}}$), left frame of Figure \ref{fig:matlab2}, as well as the mismatch absolute error, in a log-scale  (right frame of Figure \ref{fig:matlab2}) are displayed. It is notable to mention that this error is at level of machine precision. In this case, the function is recovered exactly.

\begin{remark}[Additional functions] 
Note that the \texttt{+mLF} package also provides higher level functions, namely (i) \texttt{mlf.alg1}, an implementation of \cite[Alg. 1]{Antoulas:2025} and \texttt{mlf.alg2}, an implementation of \cite[Alg. 2]{Antoulas:2025}. More details are available on \url{https://github.com/cpoussot/mLF}.
\end{remark}

%%%%%%%%%%%%%%%%%%%%%%%%%%%%%%%%%%%%%%%%%%%%%
\subsection{Summary}

In \cite{Antoulas:2025}, thanks to the \textbf{cascaded (or recursive) null space construction}, we achieve \textbf{variable decoupling} and thus provide an equivalent alternative to the standard brute force null space computation of multivariate Loewner matrices. Through this recursive null space construction, we avoid the costly intermediate large $\ord$-D Loewner matrix construction, thus saving disk access time and memory and thus taming the \COD. In addition, the $\ord$-D rational function construction problem is recast as a collection of 1-D problems, \textbf{simpler to store and solve in practice}, and leading to overall more accurate results. This statement also shows how much \textbf{the variable decoupling is intrinsically achieved by this process}. By this, in addition to the data-driven multivariate rational approximation feature, we believe we also provide a viable solution to tensor approximation (with grid data structure)  via rational functions. Finally, by connecting this result to the \textbf{Kolmogorov Superposition Theorem}, we bridge  the gap between \textbf{Kolmogorov Arnold Networks (KANs)} and \textbf{rational approximation}. In the rest of the document, we focus on the accuracy and scalability of this method and compare it to other methods, stressing the potential of the method.

%%%%%%%%%%%%%%%% 
%\newpage
\section{Overview of the results}
\label{sec:ex_overview}
We now follow the process presented in Section \ref{sec:methods} and evaluate the methods using the examples listed in Table \ref{tab:examples}. We conclude after listing the obtained performances.

\subsection{Approximation performance statistics}

\begin{figure}[H]
\centering
\includegraphics[width=\textwidth]{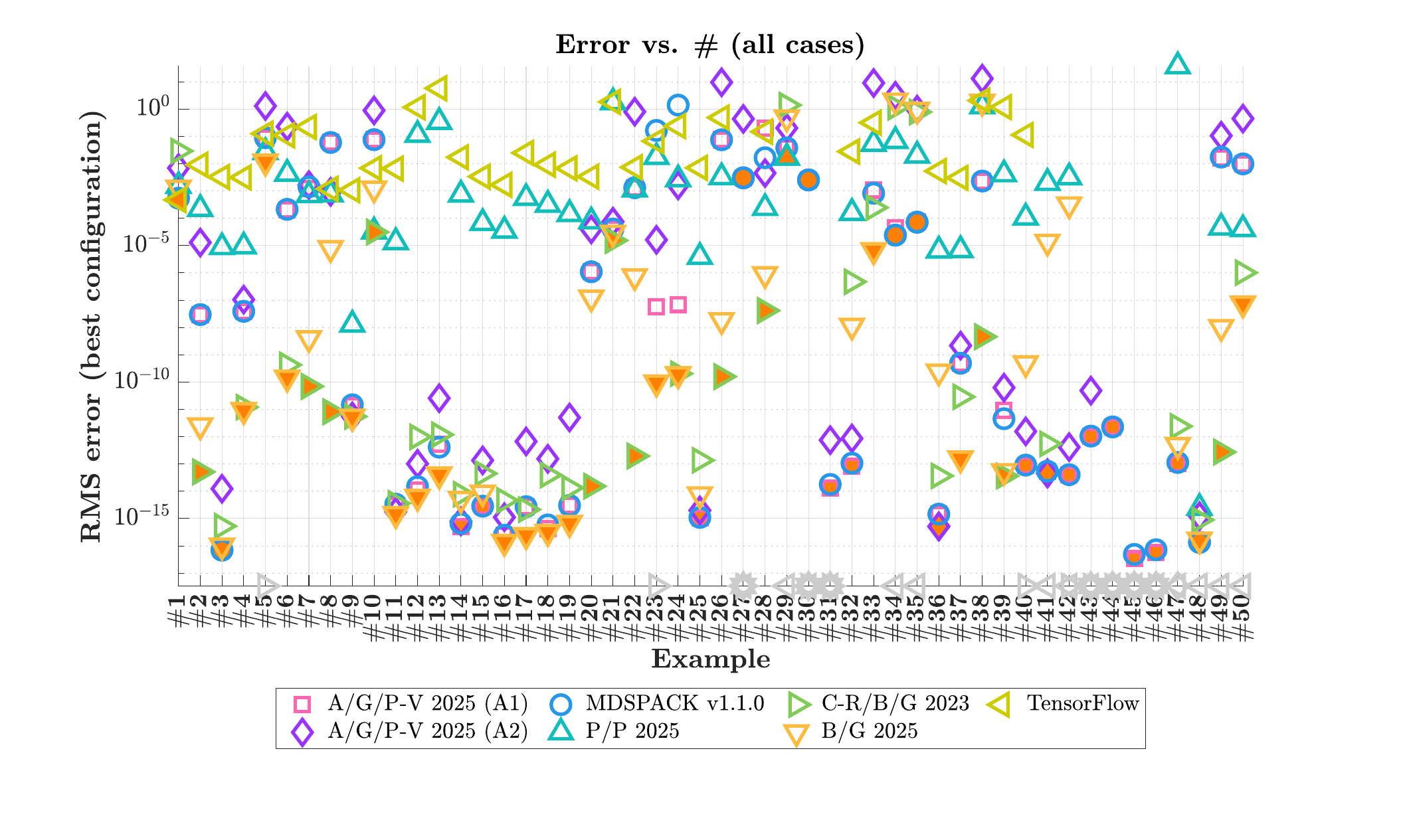}
\caption{Average RMS error of the best candidate $\bG_m$ as a function of the case number. Grey symbols indicate that the method has not converged. Red filled are used to mark the best candidate.} %Recall methods ordering: M1 \cite[Alg. 1]{Antoulas:2025}, M2 \cite[Alg. 2]{Antoulas:2025}, M3 \cite{mdspack}, M4 \cite[KAN]{Poluektov:2025}, M5 \cite[pAAA]{Balicki:2025} and M6 \cite[LR-pAAA]{Balicki:2025}, and M7 \cite[MLP by Tensor Flow]{TensorFlow}.}
\label{fig:err}
\end{figure}

\begin{figure}[H]
\centering
\includegraphics[width=\textwidth]{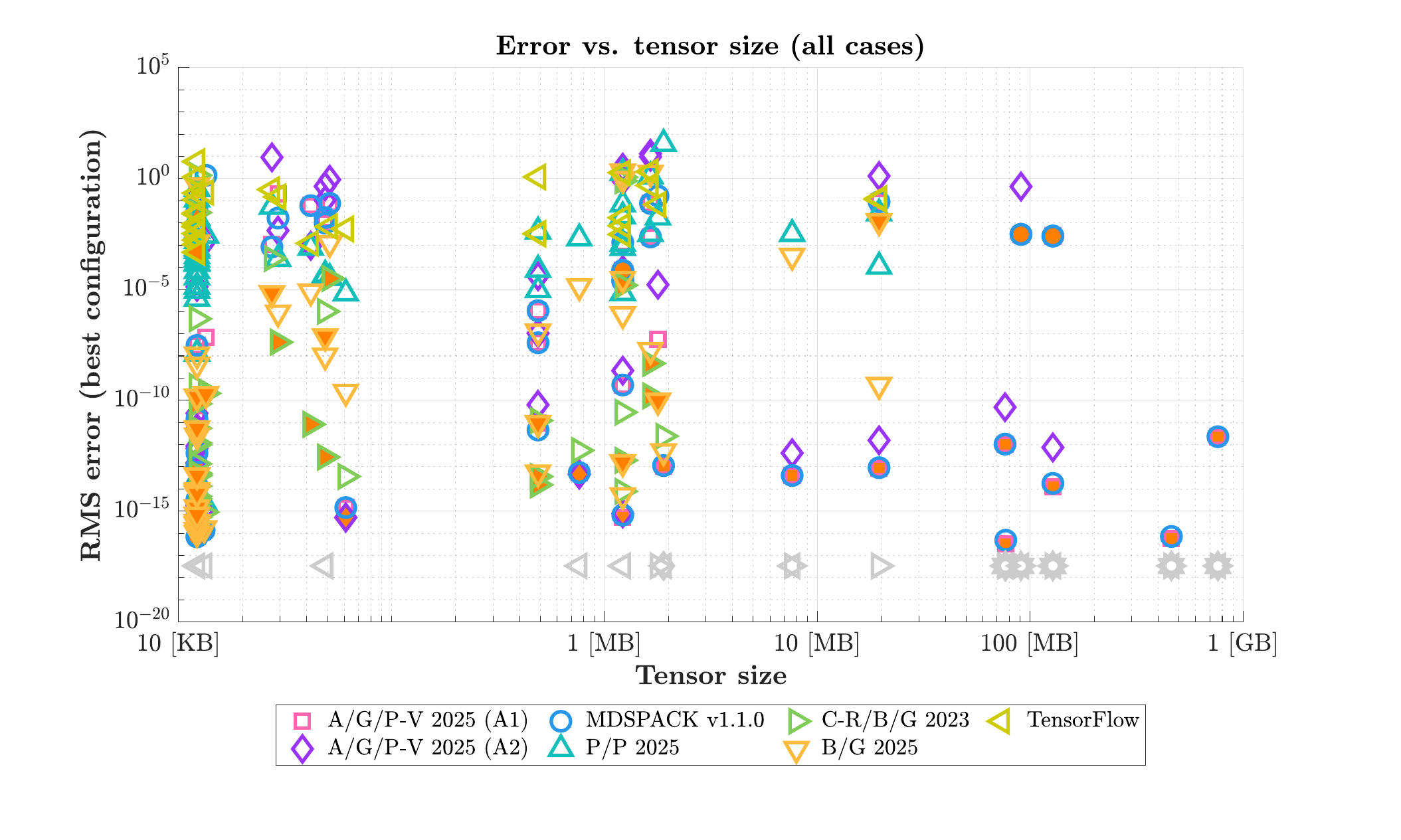}
\caption{Average RMS error of the best candidate $\bG_m$ as a function of the tensor size. Grey symbols indicate that the method has not converged. Red filled are used to mark the best candidate.} %Recall methods ordering: M1 \cite[Alg. 1]{Antoulas:2025}, M2 \cite[Alg. 2]{Antoulas:2025}, M3 \cite{mdspack}, M4 \cite[KAN]{Poluektov:2025}, M5 \cite[pAAA]{Balicki:2025} and M6 \cite[LR-pAAA]{Balicki:2025}, and M7 \cite[MLP by Tensor Flow]{TensorFlow}.}
\label{fig:err_size}
\end{figure}

For each of the \CAS\, considered examples, we evaluate the mismatch between the original and obtained surrogate $\bG_{m,p}$ model using 500 random experiments (i.e. input variable draw). Then, selecting the best average parametrization,  we obtain $\bG_m$ (i.e. the best surrogate candidate). Figure \ref{fig:err} \& \ref{fig:err_size} show the best average RMS error obtained for all cases, and as a function of the tensor size. This figure retains the best parametrization obtained with each method, only. Not-converged cases are marked with grey symbols and the best method is filled with red color\footnote{For some examples, especially the ones involving high dimensional tensors, some methods either  do not converge or were stopped as they were stuck or crashed because of memory limitations. Therefore we conclude that scalability limits were reached and are marked as  "not converged" instead.}.

\begin{figure}[H]
\centering
\includegraphics[width=.8\textwidth]{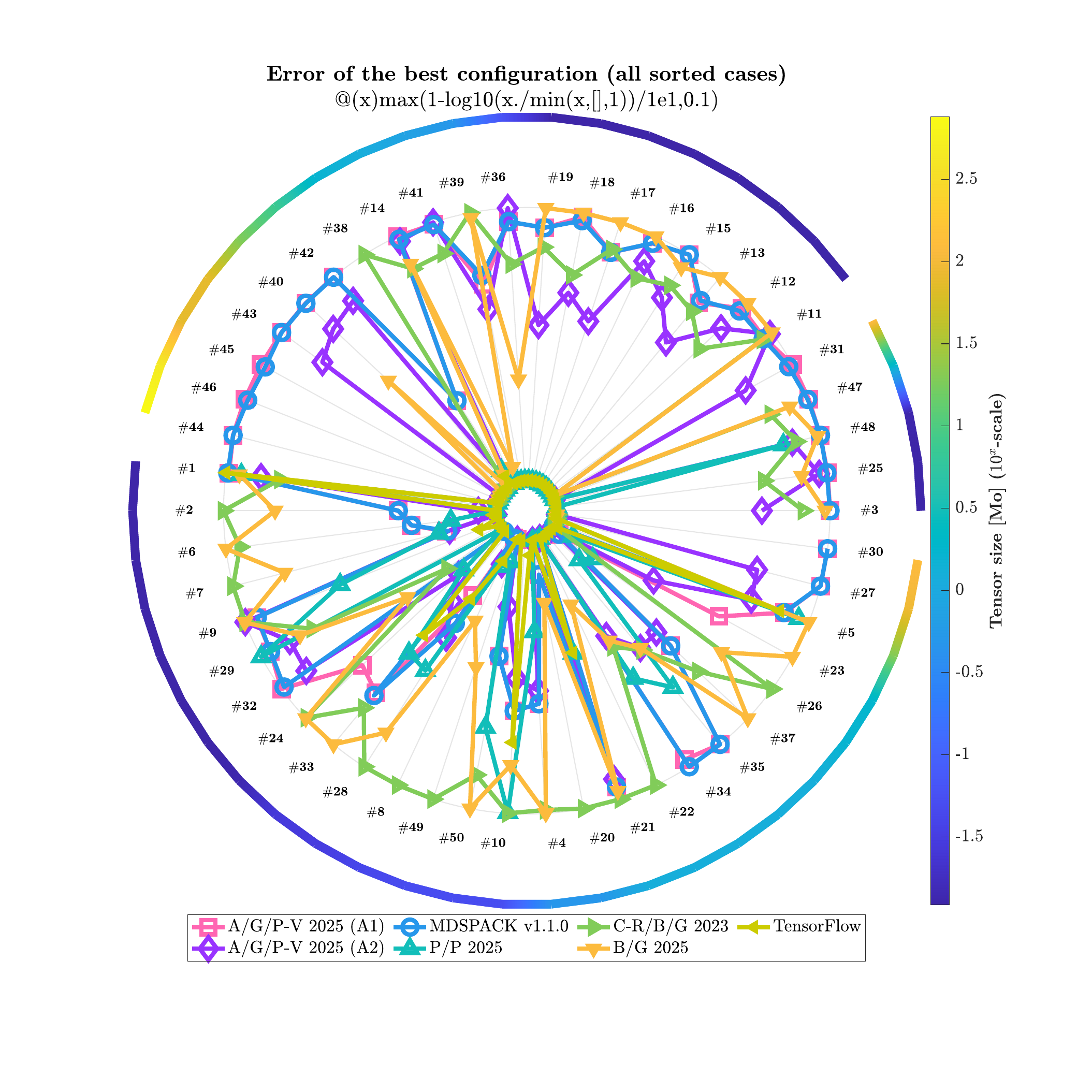}
\caption{For all function type (polynomial, rational and irrational), radar plot of the RMS error, normalized with the minimal one. Closer to unit circle are the best cases, un-converged and $10^9$ times the higher errors are on the circe with radius 0.1.}
\label{fig:radar-err-AS}
\end{figure}

Figure \ref{fig:radar-err-AS} then provides a radar representation of the RMSE errors for all methods. The circle is subdivided in three parts (outside circle color bars); starting from degree zero and moving anti-clockwise, we first group polynomial, then rational and finally irrational cases. The color gradient indicates the size of the original tensor. The symbols are located following the formula in the title. As a consequence, the closer to the unit circle the symbol is, the better the method is (this is a relative representation). Notice that this plot should be read with Figure \ref{fig:err} as in some cases, multiple methods provide very good results hard to dissociate since they are  close to machine precision

Here, we point out that only M1, M2 and M3 are able to scale for large tensors and provide an accurate solution for polynomial and rational cases. Conclusions for irrational cases are more intricate. Indeed, for low order tensors, M5 and M6 seem to be a good candidates. However they fail for high complexity.

\begin{figure}[H]
\includegraphics[width=\textwidth]{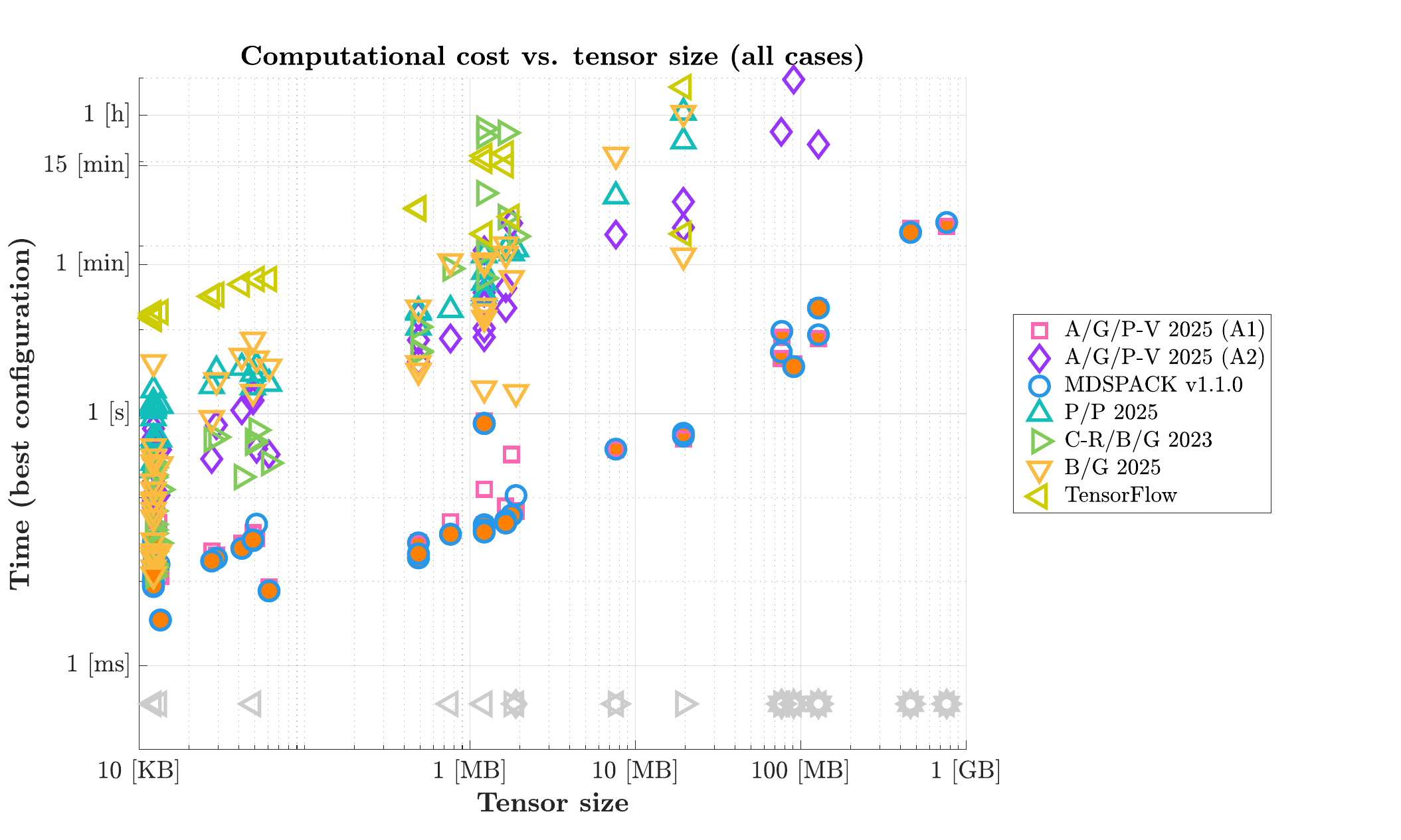}
\caption{Best candidate $\bG_m$ model construction time as a function of the tensor size. Grey symbols mean that the method has not converged. Red filled are used to mark the best method.} %Recall methods ordering: M1 \cite[Alg. 1]{Antoulas:2025}, M2 \cite[Alg. 2]{Antoulas:2025}, M3 \cite{mdspack}, M4 \cite[KAN]{Poluektov:2025}, M5 \cite[pAAA]{Balicki:2025} and M6 \cite[LR-pAAA]{Balicki:2025}, and M7 \cite[MLP by Tensor Flow]{TensorFlow}.}
%\caption{Best configuration $\bG_m$ model construction computation time vs. $\tableau{\ord}$ tensor size. Recall methods ordering: M1 \cite[Alg. 1]{Antoulas:2025}, M2 \cite[Alg. 2]{Antoulas:2025}, M3 \cite{mdspack}, M4 \cite[KAN]{Poluektov:2025}, M5 \cite[pAAA]{Balicki:2025} and M6 \cite[LR-pAAA]{Balicki:2025}, and M7 \cite[MLP by Tensor Flow]{TensorFlow}. }
\label{fig:stat-data_vs_time}
\end{figure}

Next, Figure \ref{fig:stat-data_vs_time} reports the computation time (for the best candidate only) for each method as a function of the original tensor $\tableau{\ord}$ size. It illustrates the scalability feature of the methods. Indeed, M1, M2 and M3, all based on the decoupling proposition of \cite{Antoulas:2025} scale well with tensors of dimension close to 1 \texttt{GBytes}, the others suffer of the curse of dimensionality. At this stage, notice that M2 is not that appropriate to deal with very large tensors due to the greedy iterations. %This is an axis for investigations and improvements.

\begin{figure}[H]
\centering
\includegraphics[width=.8\textwidth]{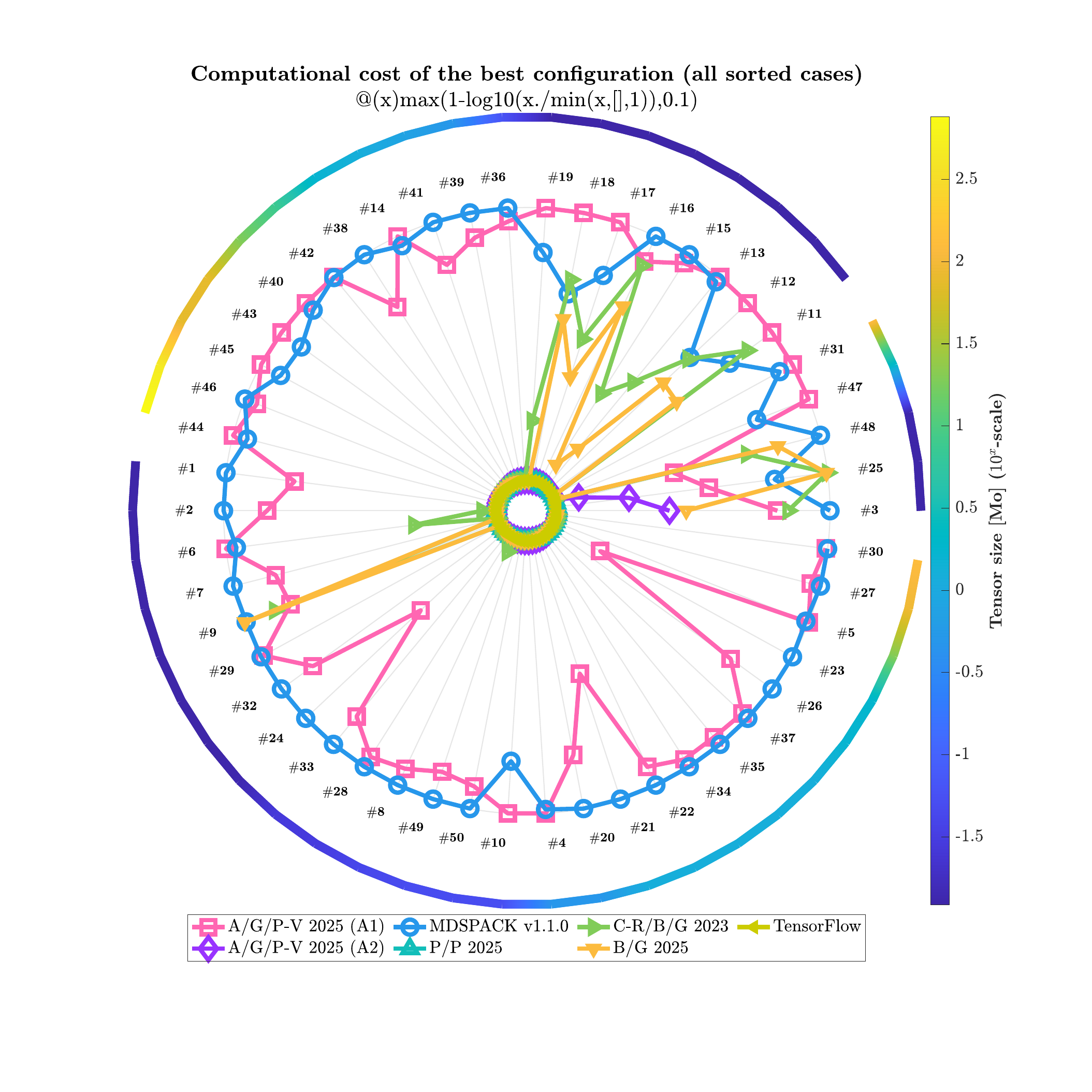}
\caption{For all function type (polynomial, rational and irrational), radar plot of the computational time, normalized with the minimal one. Closer to unit circle are the best cases, un-converged and $\approx 8$ times the slowest are on the circe with radius 0.1.}
\label{fig:radar-time-AS}
\end{figure}

Similarly, Figure \ref{fig:radar-time-AS} shows a radar representation of the computation time for all methods. The circle is again subdivided in three parts and symbols are located following the formula in the title. The closer to the unit circle the symbol is, the faster the method is. Regarding this matter, it is clear that the decoupling allows solving tensor approximation problem much faster than the other methods.

\begin{figure}[H]
\centering 
\includegraphics[width=.8\textwidth]{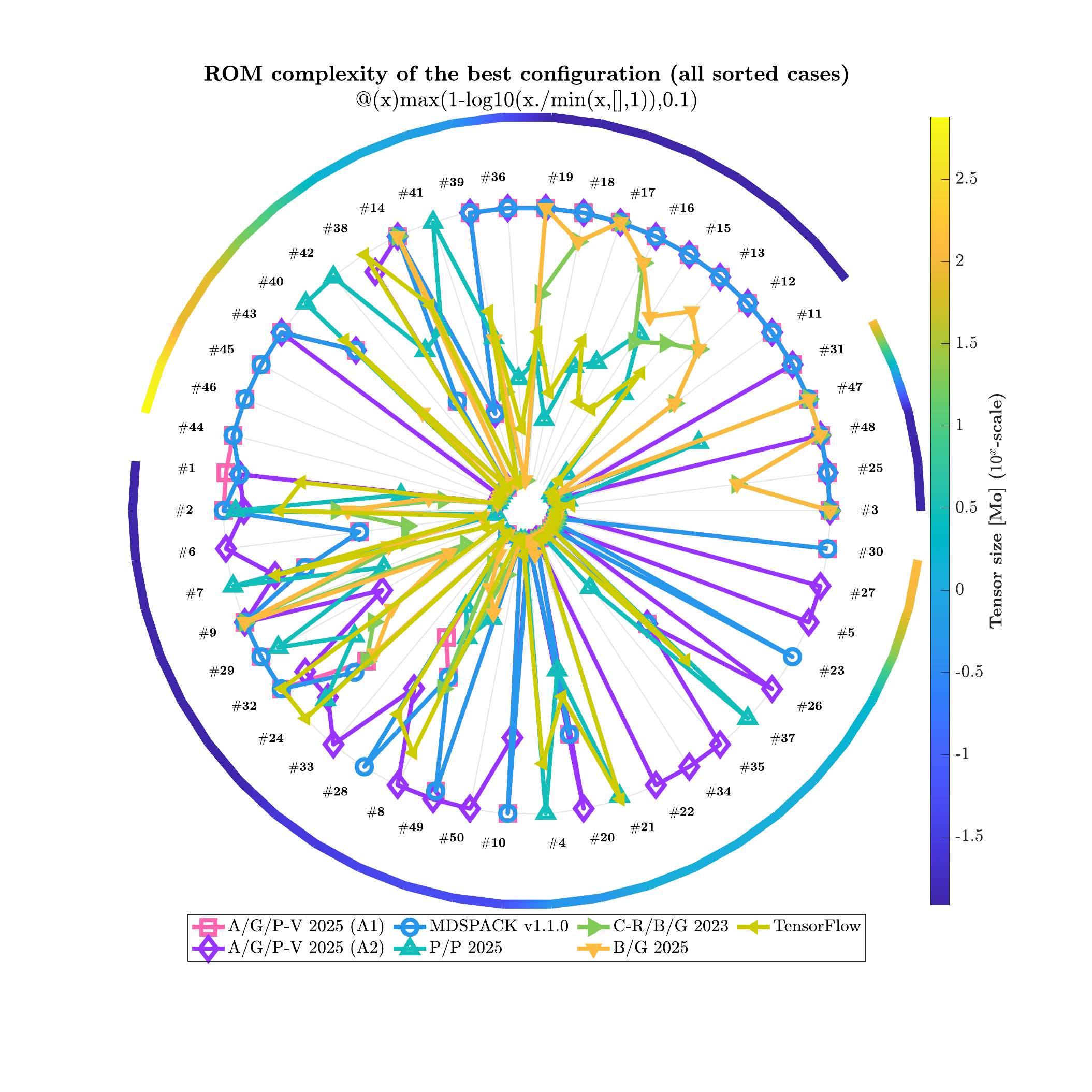}
\caption{For all function type (polynomial, rational and irrational), radar plot of the surrogate model complexity, normalized with the minimal one. Closer to unit circle are the best cases, un-converged and $\approx 8$ times the slowest are stuck to the circe with radius 0.1.}
\label{fig:radar-data-AS}
\end{figure}

Finally, Figure \ref{fig:radar-data-AS} shows a radar representation of the complexity of each surrogate model for all methods. The circle is again subdivided in three parts and symbols are located following the formula in the title. The closer to the unit circle the symbol is, the simpler the model is. Notice that we do not iterate to obtain the simplest model, but rather the more accurate. Therefore, this information can be considered as a side effect. %Of specific interest are the polynomial and rational areas (form 0  to about 180 degrees) where 

%\begin{remark}[Feedback to authors]
%This report is aimed at being regularly updated and improved. Please send any feedback or suggestions to \texttt{charles.poussot-vassal@onera.fr}.
%\end{remark}

%\foreach \n in {1,...,\CASBLK}{
%\begin{figure}[H]
%\begin{subfigure}{\textwidth}
%\centering
%\includegraphics[width=.8\textwidth]{../../../Benchmarks/mLF_evaluation/v2/tex_pdf/figures/stat/err_\n}
%\caption{Mismatch absolute error $\lvert\bG_m-\bH\lvert$ of the best configuration only (log-scale). }
%\label{fig:stat-error-\n}
%\end{subfigure}
%\begin{subfigure}{\textwidth}
%\centering
%\includegraphics[width=.8\textwidth]{../../../Benchmarks/mLF_evaluation/v2/tex_pdf/figures/stat/time_\n}
%\caption{Model $\bG_m$ computation time of the best configuration only (log-scale).}
%\label{fig:stat-time-\n}
%\end{subfigure}
%\begin{subfigure}{\textwidth}
%\centering
%\includegraphics[width=.8\textwidth]{../../../Benchmarks/mLF_evaluation/v2/tex_pdf/figures/stat/data_\n}
%\caption{Model $\bG_m$ complexity of the best configuration only (log-scale).}
%\label{fig:stat-data-\n}
%\end{subfigure}
%\caption{Examples set \n.}
%\end{figure}
%}

%\begin{figure}[h]
%\includegraphics[width=\textwidth]{figures/stat_data_time}
%\caption{Average $\bG_m$ model construction computation time vs. $\tableau{\ord}$ tensor size. Recall methods ordering: M1 \cite[Alg. 1]{Antoulas:2025}, M2 \cite[Alg. 2]{Antoulas:2025}, M3 \cite{mdspack}, M4 \cite[KAN]{Poluektov:2025}, M5 \cite[pAAA]{Balicki:2025} and M6 \cite[LR-pAAA]{Balicki:2025}, and M7 \cite[MLP by Tensor Flow]{TensorFlow}. }
%%\label{fig:stat-data_vs_time}
%\end{figure}

\subsection{Preliminary remarks}

\textbf{Speed and scalability.}  By inspecting Figure \ref{fig:stat-data_vs_time} \& \ref{fig:radar-time-AS} it is clear that M1, M2 and M3 provide a very fast solution to the tensor approximation problem, with computational times way faster than the other methods. In addition with reference to Figure \ref{fig:stat-data_vs_time} we demonstrate that these methods are able to address high dimensional tensors, where others fail or lead to prohibitive computational time. This is an essential feature for taming the \COD. We claim that \textbf{the recursive null space computation proposed in \cite{Antoulas:2025} is a  cutting edge solution to both the scalability issue and the user experience improvement}. In  \cite{Antoulas:2025}, an illustration of a function with 20 variables is also shown. \\%Beside this aspect, fast algorithm also permit to iterate on the tuning parameters.  \\

\noindent \textbf{Tuning parameters.}  Remembering that each method has multiple tunable parameters, finding the adequate/optimal parameter set is not a trivial task. In Section \ref{sec:ex_details} we show how the "optimal" parameters vary from an example to an other. In this light, there is no clear optimal tuning parameter set for all functions. While in the single variable case (stopping) criteria are quite well understood (e.g. Loewner rank, singularities/eigenvalues, root mean square error, complexity, etc.), in the multivariate case, such criteria still need to be defined, analyzed and adjusted. We believe this is an open research question, not only from an engineering point of view, but also from a theoretical one (this is especially true for non-rational cases). %Moreover, it is still not so clear which parameters play an important role in the success or failure of some methods (especially for non-rational cases). 
This collection of examples gives some insight of issues encountered and potential solutions, but a rather generic solution needs to be discovered. This will be the purpose of future investigations. One interesting argument for M1, M2 and M3, is that the fast computational process allows for multiple (greedy) iterations on the tuning parameter configurations.\\

\noindent \textbf{Accuracy.} By inspecting Figure \ref{fig:err}, \ref{fig:err_size} \& \ref{fig:radar-err-AS}  we notice that some models obtained with M1, M2, M3, M5 and M6 lead to RMSE close or even below machine precision. This generally means that the underlying model generating the tensor has been recovered (or discovered) from data only. This property is (always) fulfilled when the generating system $\bH$ is a polynomial or rational model. Indeed the barycentric form of the surrogate model fits the rational form. When the generating function $\bH$ is irrational, machine precision mismatch error may not be obtained, or at the price of a complex surrogate. Even if machine precision is not reached, we would like to mention that most methods perform well overall. % which is important to point out, giving also credit to third parties softwares. 
However, when the tensor size or the number of variables $\ord$ increases, M4, M5, M6 and M7, which do not benefit from the recursive and decoupling scheme, fail due to memory issues or prohibitive computational time. This highlights the scalability issue, proposed in \cite{Antoulas:2025}.\\

\noindent \textbf{Complexity.}  By inspecting Figure \ref{fig:radar-data-AS} we may first notice an overall homogeneity in the complexity of the achieved models. However a precise inspection reveals that M4, M5 and M6 (and sometimes M2) tend to overfit by constructing models with too many variables. This is especially true when the original function is a rational function, where complexity and minimality can be evaluated exactly. Regarding M7, no conclusion can be made since one single configuration has been tested.\\

\noindent \textbf{Comments on KAN and MLP based methods.}  M4 is an algorithm implementing KAN-based models, while M7 is an algorithm implementing MLP-based models. In both cases, the topology of the network influence the optimizer, that is chosen by the user. %As non experts in network modeling, we humbly chose some configurations that are amendable. Authors are open for feedbacks from more expert users, leading to improved  parameterizations and results.

%%%%%%%%%%%%%%%%
%\newpage
\section{Detailed examples exposition and results}
\label{sec:ex_details}

For the considered examples, we now list for all the cases the  best configuration and provide some details on the results.

%\foreach \n in {1,...,\CAS}{
%\input{figures/case_\n/text_main}
%}
\newpage \subsection{Function \#1 (${\ord=2}$ variables, tensor size: 12.5 \textbf{KB})} $$\mathrm{ReLU}(\var{1})+\frac{1}{100}\var{2}$$ \subsubsection{Setup and results overview}\begin{itemize}\item Reference: Personal communication, [none]\item Domain: $\mathbb{R}$\item Tensor size: 12.5 \textbf{KB} ($40^{2}$ points)\item Bounds: $ \left(\begin{array}{cc} -1 & 1 \end{array}\right) \times \left(\begin{array}{cc} -1 & -\frac{1}{10000000000} \end{array}\right)$ \end{itemize} \begin{table}[H] \centering \begin{tabular}{llllll}
$\#$ & Alg. & Parameters & Dim. & CPU [s] & RMSE \\ 
\hline 
$\mathbf{\#1}$ & A/G/P-V 2025 (A1) & $1 \cdot 10^{-11},3$ & $\mathbf{1.4 \cdot 10^{02}}$ & $0.017$ & $0.0006$ \\ 
 & A/G/P-V 2025 (A2) & $1 \cdot 10^{-15},1$ & $1.6 \cdot 10^{02}$ & $0.66$ & $0.0069$ \\ 
 & MDSPACK v1.1.0 & $-1,9$ & $1.6 \cdot 10^{02}$ & $\mathbf{0.0099}$ & $0.00055$ \\ 
 & P/P 2025 & $1,0.95,50,0.01,10,6,21$ & $5.5 \cdot 10^{02}$ & $1.1$ & $0.0015$ \\ 
 & C-R/B/G 2023 & $1 \cdot 10^{-06},20$ & $7.6 \cdot 10^{02}$ & $0.26$ & $0.029$ \\ 
 & B/G 2025 & $1 \cdot 10^{-06},20,4$ & $6.8 \cdot 10^{02}$ & $0.26$ & $0.0013$ \\ 
 & TensorFlow & $$ & $2.6 \cdot 10^{02}$ & $16$ & $\mathbf{0.00047}$ \\ 
\hline 
\end{tabular} \caption{Function \#1: best model configuration and performances per methods.} \end{table}\begin{figure}[H] \centering  \includegraphics[width=\textwidth]{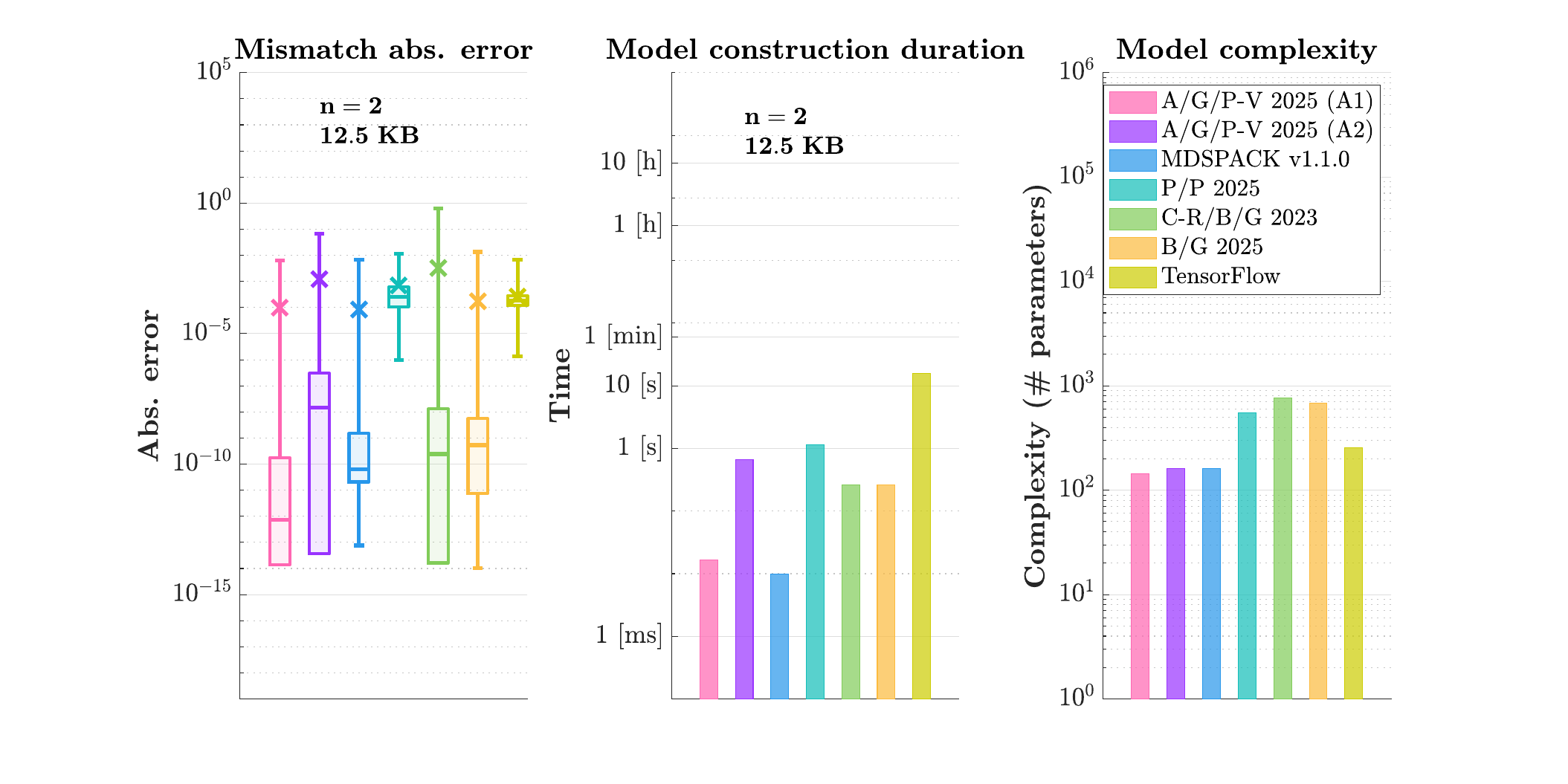} \caption{Function \#1: graphical view of the best model performances.} \end{figure}\begin{figure}[H] \centering  \includegraphics[width=\textwidth]{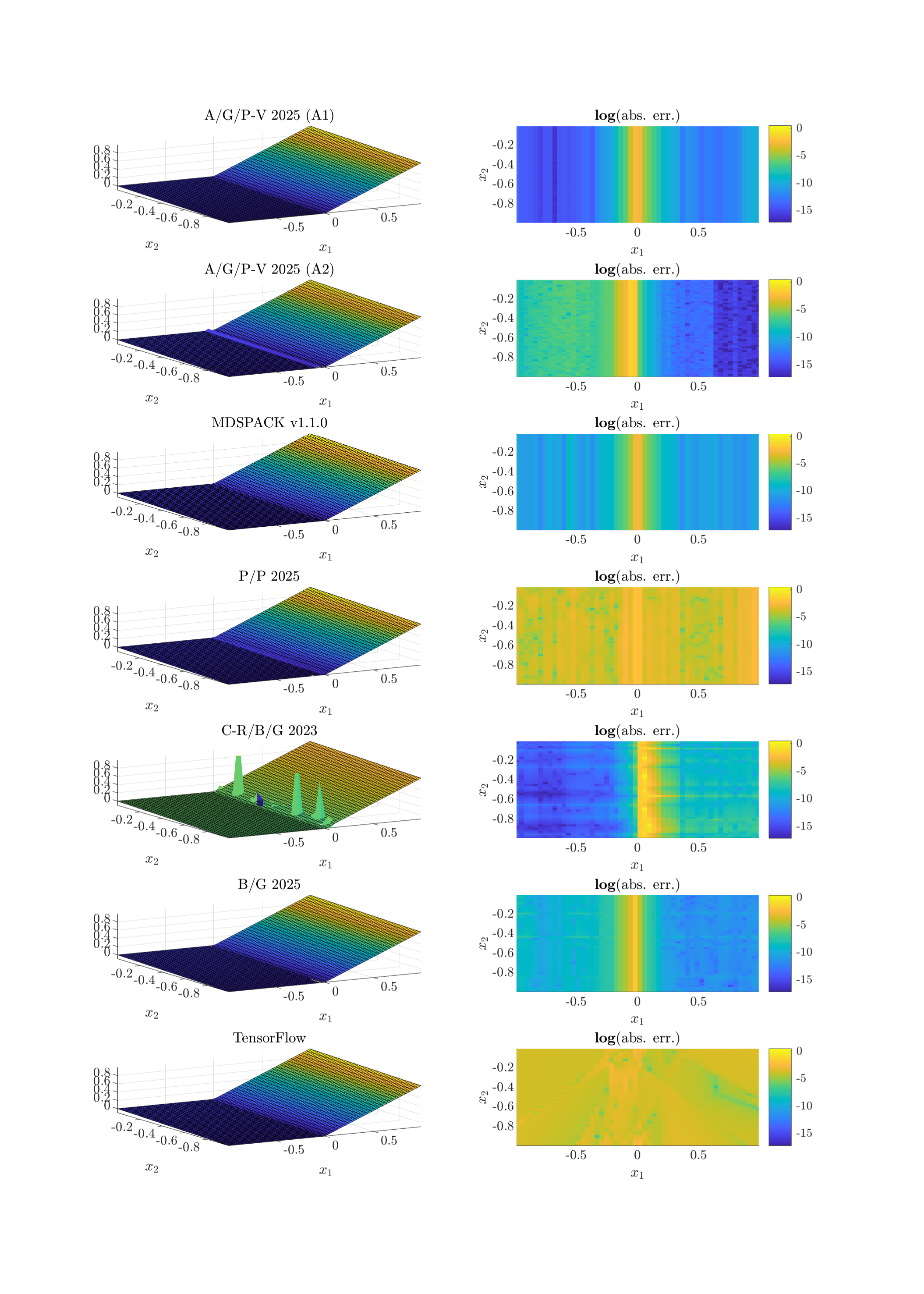} \caption{Function \#1: left side, evaluation of the original (mesh) vs. approximated (coloured surface) and right side, absolute errors (in log-scale).} \end{figure}\subsubsection{mLF detailed informations (M1)} \noindent \textbf{Right interpolation points}: $k_l=\left(\begin{array}{cc} 18 & 2 \end{array}\right)$, where $l=1,\cdots,\ord$.$$ \begin{array}{rcl}\lan{1} &\in& \IC^{18} \text{ , linearly spaced between bounds}\\\lan{2} &\in& \IC^{2} \text{ , linearly spaced between bounds}\\\end{array} $$\noindent \textbf{$\ord$-D Loewner matrix, barycentric weights and Lagrangian basis}:$$ \begin{array}{rcl}\IL & \in & \IC^{36 \times 36}\\\bc & \in & \IC^{36}\\\bw & \in & \IC^{36}\\\bc\odot \bw & \in & \IC^{36}\\\mathbf{Lag}(\var{1},\var{2}) & \in & \IC^{36}\\\end{array} $$

\newpage \subsection{Function \#2 (${\ord=2}$ variables, tensor size: 12.5 \textbf{KB})} $$\mathrm{exp}\left(\sin(\var{1}) + \var{2}^2\right)$$ \subsubsection{Setup and results overview}\begin{itemize}\item Reference: L/al. 2024, \cite{Liu:2025}\item Domain: $\mathbb{R}$\item Tensor size: 12.5 \textbf{KB} ($40^{2}$ points)\item Bounds: $ \left(\begin{array}{cc} -1 & 1 \end{array}\right) \times \left(\begin{array}{cc} -1 & 1 \end{array}\right)$ \end{itemize} \begin{table}[H] \centering \begin{tabular}{llllll}
$\#$ & Alg. & Parameters & Dim. & CPU [s] & RMSE \\ 
\hline 
$\mathbf{\#2}$ & A/G/P-V 2025 (A1) & $1 \cdot 10^{-06},1$ & $\mathbf{1.7 \cdot 10^{02}}$ & $0.014$ & $3 \cdot 10^{-08}$ \\ 
 & A/G/P-V 2025 (A2) & $1 \cdot 10^{-15},3$ & $2 \cdot 10^{02}$ & $0.29$ & $1.3 \cdot 10^{-05}$ \\ 
 & MDSPACK v1.1.0 & $1 \cdot 10^{-06},3$ & $1.7 \cdot 10^{02}$ & $\mathbf{0.0097}$ & $3 \cdot 10^{-08}$ \\ 
 & P/P 2025 & $1,1,50,0.01,4,12,9$ & $1.8 \cdot 10^{02}$ & $0.51$ & $0.00022$ \\ 
 & C-R/B/G 2023 & $1 \cdot 10^{-09},20$ & $4 \cdot 10^{02}$ & $0.069$ & $\mathbf{5.1 \cdot 10^{-14}}$ \\ 
 & B/G 2025 & $1 \cdot 10^{-09},20,3$ & $4.3 \cdot 10^{02}$ & $0.16$ & $2.6 \cdot 10^{-12}$ \\ 
 & TensorFlow & $$ & $2.6 \cdot 10^{02}$ & $15$ & $0.0091$ \\ 
\hline 
\end{tabular} \caption{Function \#2: best model configuration and performances per methods.} \end{table}\begin{figure}[H] \centering  \includegraphics[width=\textwidth]{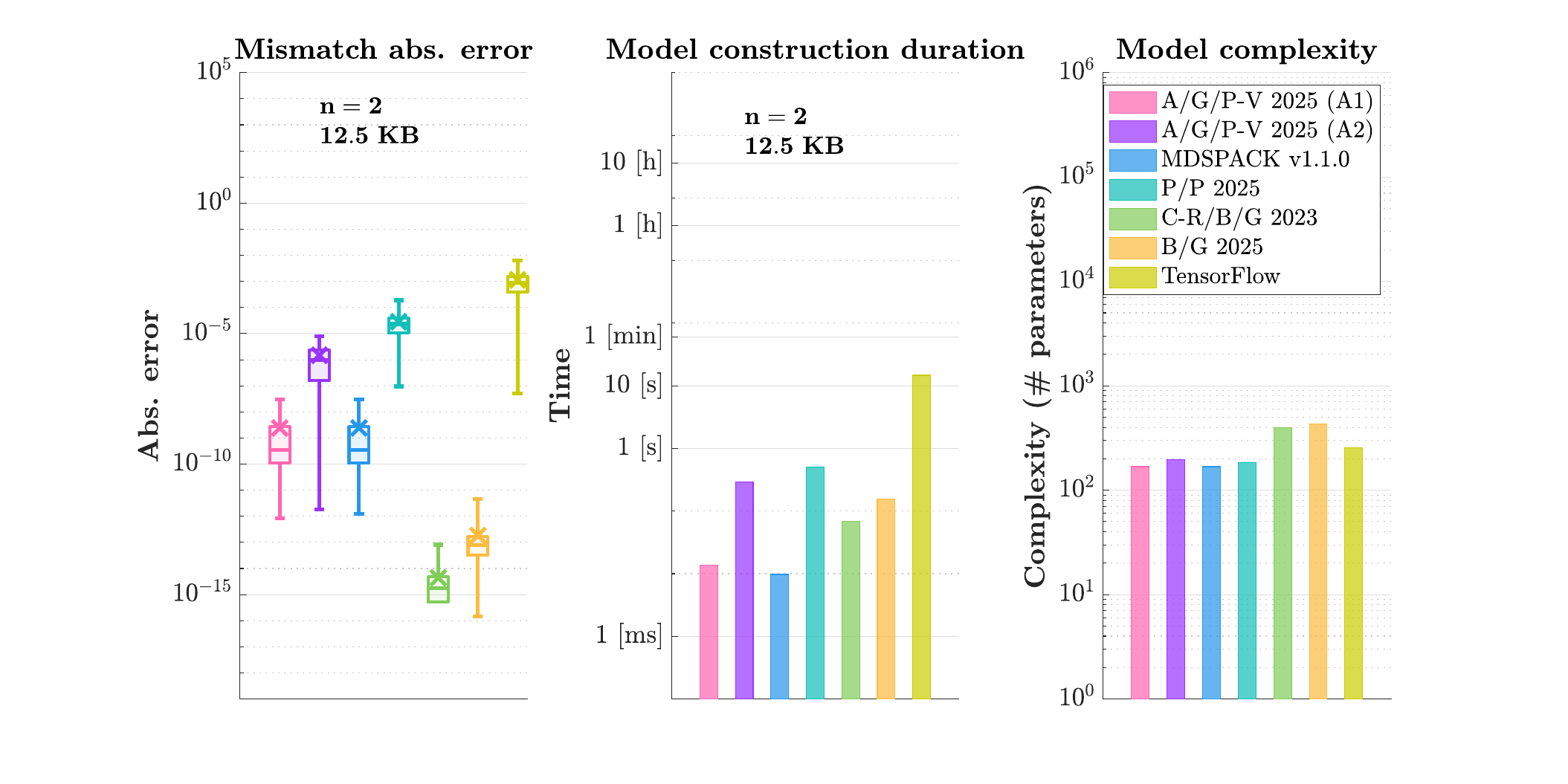} \caption{Function \#2: graphical view of the best model performances.} \end{figure}\begin{figure}[H] \centering  \includegraphics[width=\textwidth]{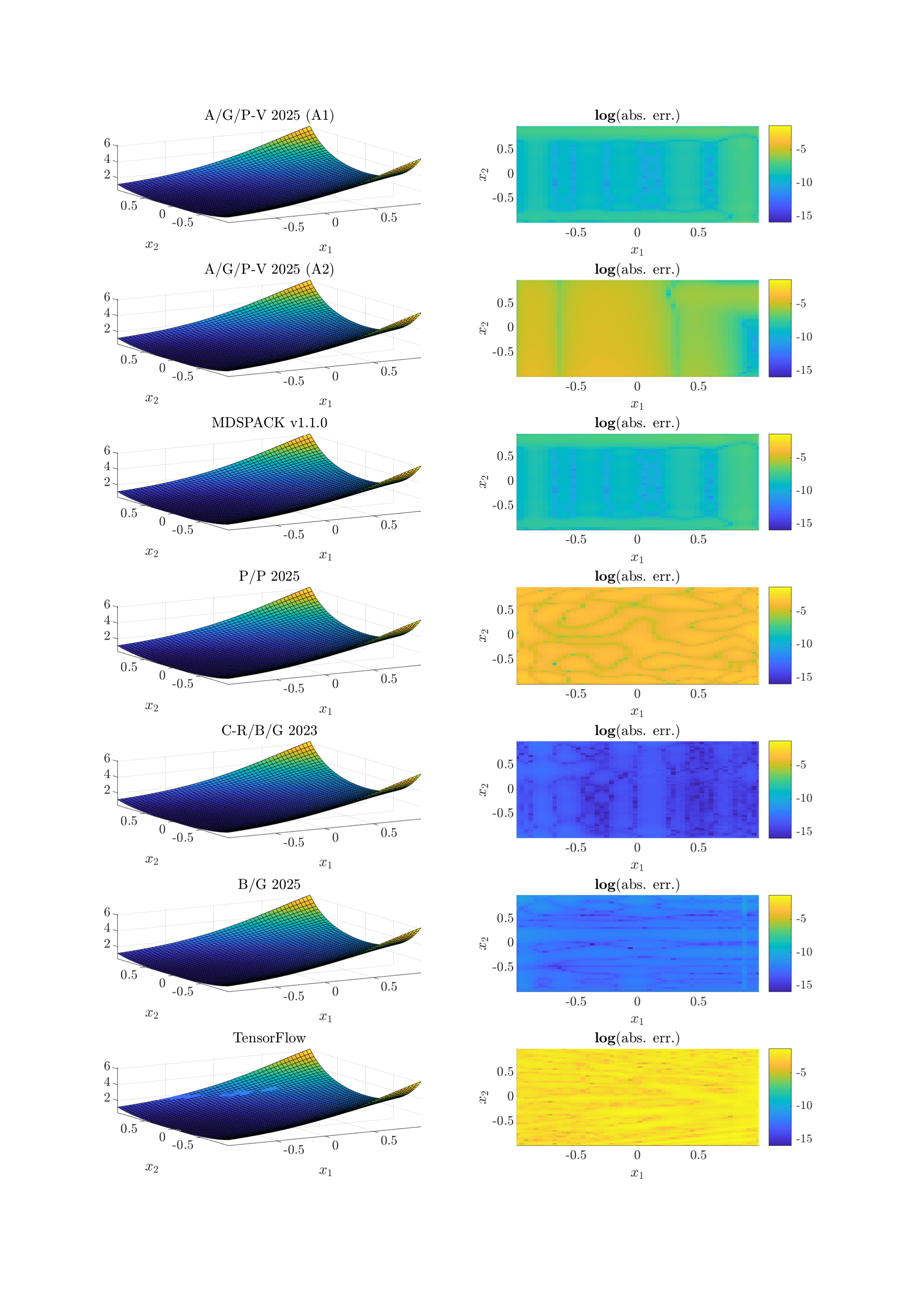} \caption{Function \#2: left side, evaluation of the original (mesh) vs. approximated (coloured surface) and right side, absolute errors (in log-scale).} \end{figure}\subsubsection{mLF detailed informations (M1)} \noindent \textbf{Right interpolation points}: $k_l=\left(\begin{array}{cc} 6 & 7 \end{array}\right)$, where $l=1,\cdots,\ord$.$$ \begin{array}{rcl}\lan{1} &\in& \IC^{6} \text{ , linearly spaced between bounds}\\\lan{2} &\in& \IC^{7} \text{ , linearly spaced between bounds}\\\end{array} $$\noindent \textbf{$\ord$-D Loewner matrix, barycentric weights and Lagrangian basis}:$$ \begin{array}{rcl}\IL & \in & \IC^{42 \times 42}\\\bc & \in & \IC^{42}\\\bw & \in & \IC^{42}\\\bc\odot \bw & \in & \IC^{42}\\\mathbf{Lag}(\var{1},\var{2}) & \in & \IC^{42}\\\end{array} $$

\newpage \subsection{Function \#3 (${\ord=2}$ variables, tensor size: 12.5 \textbf{KB})} $$\var{1} \var{2}$$ \subsubsection{Setup and results overview}\begin{itemize}\item Reference: L/al. 2024, \cite{Liu:2025}\item Domain: $\mathbb{R}$\item Tensor size: 12.5 \textbf{KB} ($40^{2}$ points)\item Bounds: $ \left(\begin{array}{cc} -1 & 1 \end{array}\right) \times \left(\begin{array}{cc} -1 & 1 \end{array}\right)$ \end{itemize} \begin{table}[H] \centering \begin{tabular}{llllll}
$\#$ & Alg. & Parameters & Dim. & CPU [s] & RMSE \\ 
\hline 
$\mathbf{\#3}$ & A/G/P-V 2025 (A1) & $0.5,1$ & $\mathbf{16}$ & $0.045$ & $\mathbf{7 \cdot 10^{-17}}$ \\ 
 & A/G/P-V 2025 (A2) & $1 \cdot 10^{-15},1$ & $16$ & $0.1$ & $1.2 \cdot 10^{-14}$ \\ 
 & MDSPACK v1.1.0 & $0.01,1$ & $16$ & $\mathbf{0.03}$ & $7 \cdot 10^{-17}$ \\ 
 & P/P 2025 & $1,1,50,0.01,4,6,9$ & $1.3 \cdot 10^{02}$ & $0.26$ & $8.8 \cdot 10^{-06}$ \\ 
 & C-R/B/G 2023 & $0.001,20$ & $16$ & $0.041$ & $5.3 \cdot 10^{-16}$ \\ 
 & B/G 2025 & $0.001,20,2$ & $16$ & $0.09$ & $1 \cdot 10^{-16}$ \\ 
 & TensorFlow & $$ & $2.6 \cdot 10^{02}$ & $14$ & $0.0033$ \\ 
\hline 
\end{tabular} \caption{Function \#3: best model configuration and performances per methods.} \end{table}\begin{figure}[H] \centering  \includegraphics[width=\textwidth]{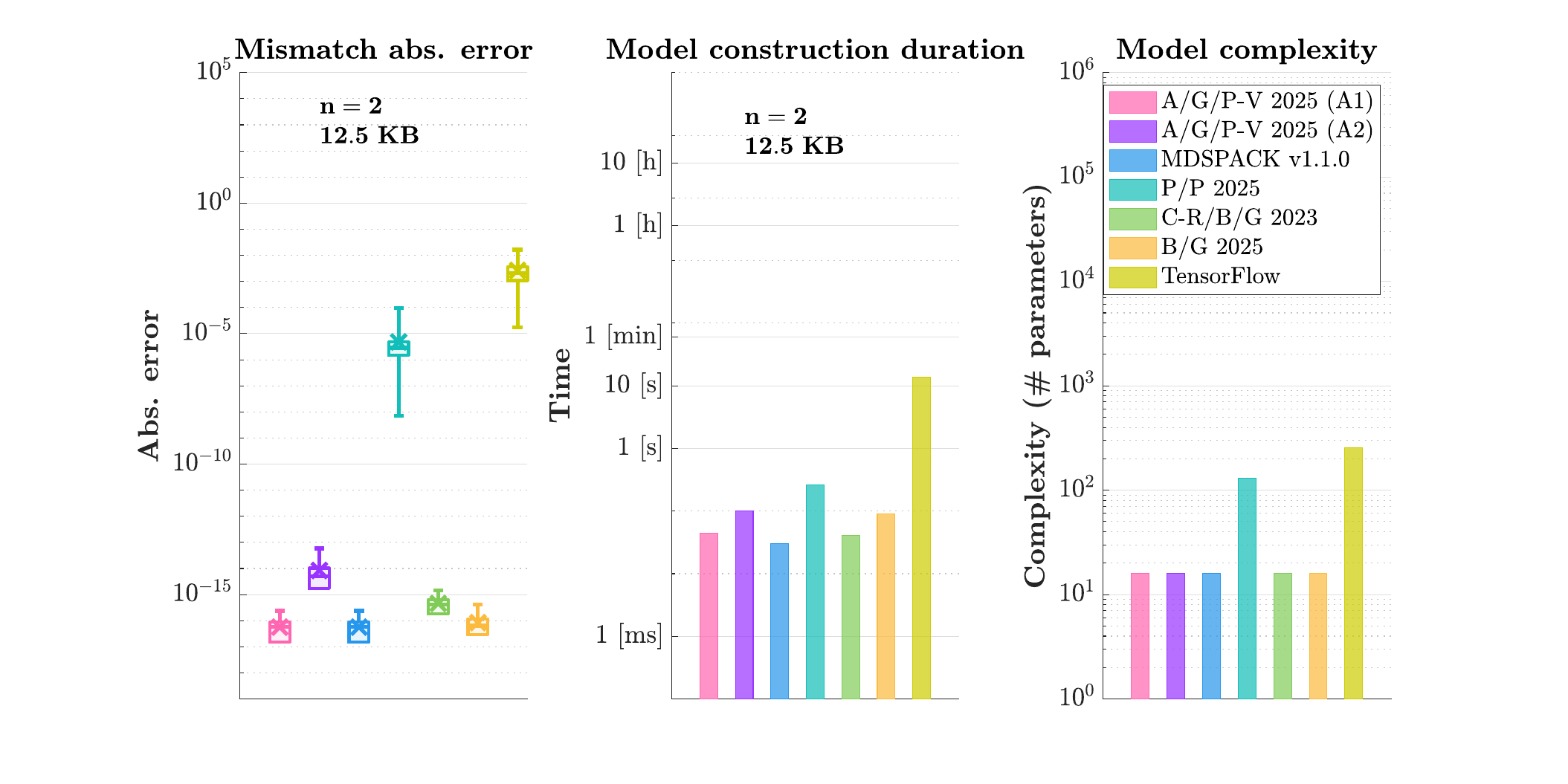} \caption{Function \#3: graphical view of the best model performances.} \end{figure}\begin{figure}[H] \centering  \includegraphics[width=\textwidth]{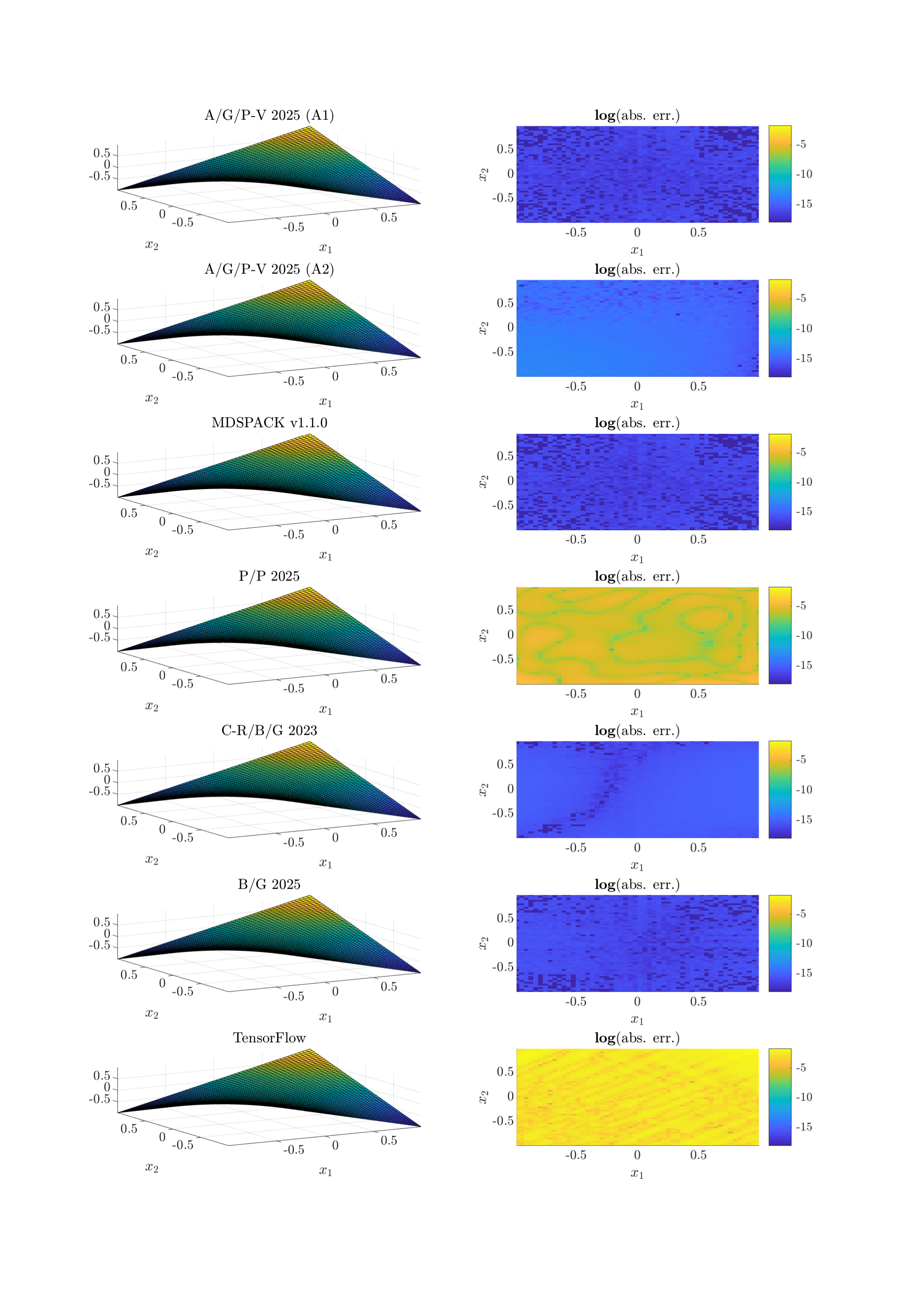} \caption{Function \#3: left side, evaluation of the original (mesh) vs. approximated (coloured surface) and right side, absolute errors (in log-scale).} \end{figure}\subsubsection{mLF detailed informations (M1)} \noindent \textbf{Right interpolation points} ($k_l=\left(\begin{array}{cc} 2 & 2 \end{array}\right)$, where $l=1,\cdots,\ord$):$$ \begin{array}{rcl}\lan{1} &=& \left(\begin{array}{cc} -1 & 1 \end{array}\right)\\\lan{2} &=& \left(\begin{array}{cc} -1 & 1 \end{array}\right)\\\end{array} $$\noindent \textbf{Lagrangian weights}: $$\left(\begin{array}{ccc} \bc & \bw & \bc\odot\bw\\ 1.0 & 1.0 & 1.0\\ -1.0 & -1.0 & 1.0\\ -1.0 & -1.0 & 1.0\\ 1.0 & 1.0 & 1.0 \end{array}\right)$$\noindent \textbf{Lagrangian form} (basis, numerator and denominator coefficients):$$\left(\begin{array}{ccc}\mathcal{B}_\textrm{lag}(\var{1},\var{2}) & \bN_\textrm{lag} &\bD_\textrm{lag}\end{array}\right) =$$ $$\left(\begin{array}{ccc} \left(\var{1}+1.0\right)\,\left(\var{2}+1.0\right) & 1.0 & 1.0\\ \left(\var{1}+1.0\right)\,\left(\var{2}-1.0\right) & 1.0 & -1.0\\ \left(\var{1}-1.0\right)\,\left(\var{2}+1.0\right) & 1.0 & -1.0\\ \left(\var{1}-1.0\right)\,\left(\var{2}-1.0\right) & 1.0 & 1.0 \end{array}\right).$$\noindent The corresponding function is:$$\begin{array}{rcl}\bG_{\textrm{lag}}(\var{1},\var{2}) &=& \dfrac{\bn_{\textrm{lag}}(\var{1},\var{2})}{\bd_{\textrm{lag}}(\var{1},\var{2})}\\ && \\&=& \dfrac{\sum_{\textrm{row}} \bN_\textrm{lag} \odot\mathcal{B}^{-1}_\textrm{lag}(\var{1},\var{2})}{\sum_{\textrm{row}} \bD_\textrm{lag} \odot\mathcal{B}^{-1}_\textrm{lag}(\var{1},\var{2})}, \end{array}$$\noindent where,\\$\bn_{\textrm{lag}}(\var{1},\var{2}) = \var{1}\,\var{2}$ \\~~\\$\bd_{\textrm{lag}}(\var{1},\var{2}) = 1.0$ \\~~\\\noindent \textbf{Monomial form} (basis, numerator and denominator coefficients - entries $<10^{-12}$ removed):$$\left(\begin{array}{ccc}\mathcal{B}_\textrm{mon}(\var{1},\var{2}) & \bN_\textrm{mon} &\bD_\textrm{mon}\end{array}\right) =$$ $$\left(\begin{array}{ccc} \var{1}\,\var{2} & 1.0 & 0\\ \var{1} & 0 & 0\\ \var{2} & 0 & 0\\ 1.0 & 0 & 1.0 \end{array}\right)$$\noindent The corresponding function is:$$\begin{array}{rcl}\bG_{\textrm{mon}}(\var{1},\var{2}) &=& \dfrac{\bn_{\textrm{mon}}(\var{1},\var{2})}{\bd_{\textrm{mon}}(\var{1},\var{2})}\\ && \\&=& \dfrac{\sum_{\textrm{row}} \bN_\textrm{mon} \odot \mathcal{B}_\textrm{mon}(\var{1},\var{2})}{\sum_{\textrm{row}} \bD_\textrm{mon} \odot\mathcal{B}_\textrm{mon}(\var{1},\var{2})},  \end{array}$$\noindent where,\\$\bn_{\textrm{mon}}(\var{1},\var{2}) = \var{1}\,\var{2}$ \\~~\\$\bd_{\textrm{mon}}(\var{1},\var{2}) = 1.0$ \\~~\\\noindent \textbf{KST equivalent decoupling pattern} (Barycentric weights $\bc^{\var{l}}$): $$\begin{array}{rclll}\var{2}&: & \left(\begin{array}{cc} -1.0 & -1.0\\ 1.0 & 1.0 \end{array}\right)& \textrm{vec}(.) & := \textbf{Bary}(\var{2}) \\\var{1}&: & \left(\begin{array}{c} -1.0\\ 1.0 \end{array}\right)& \textrm{vec}(.) \otimes \bone_{k_{2}} & := \textbf{Bary}(\var{1}) \\\end{array}$$~\\ Then, with the above notations, one defines the following univariate vector functions:  $$ \left\{ \begin{array}{rcl}\bPhi_{1}(\var{1}) &:=& \textbf{Bary}(\var{1}) \odot \mathbf{Lag}(\var{1}) \\\bPhi_{2}(\var{2}) &:=& \textbf{Bary}(\var{2}) \odot \mathbf{Lag}(\var{2}) \\\end{array} \right. $$\noindent The corresponding function is:$$\begin{array}{rcl}\bG_{\textrm{kst}}(\var{1},\var{2}) &=& \dfrac{\bn_{\textrm{kst}}(\var{1},\var{2})}{\bd_{\textrm{kst}}(\var{1},\var{2})}\\ && \\ &=& \dfrac{\sum_{\text{rows}} \bw \odot \bPhi_{1}(\var{1}) \odot \cdots \odot\bPhi_{2}(\var{2})}{\sum_{\text{rows}} \bPhi_{1}(\var{1}) \odot \cdots \odot\bPhi_{2}(\var{2})} . \end{array}$$~\\ \noindent \textbf{KST-like univariate functions} (equivalent scaled univariate functions $\bphi_{1,\cdots,2}$): $$\left\{\begin{array}{rcrcl}z_{1} &=&\bphi_{1}(\var{1}) &=& \var{1}\\z_{2} &=&\bphi_{2}(\var{2}) &=& -1.0\,\var{2}\\\end{array} \right. .$$\noindent \textbf{Connection with Neural Networks} (equivalent numerator $\bn_{\textrm{lag}}$ representation):\\ \begin{figure}[H]\begin{center} \scalebox{.7}{\begin{tikzpicture}[line width=0.4mm]\tikzstyle{place}=[circle, draw=black, minimum size = 8mm]\tikzstyle{placeInOut}=[circle, draw=orange, minimum size = 8mm]\node at (0,-2) [placeInOut] (first_1){$\var{1}$};\node at (0,-4) [placeInOut] (first_2){$\var{2}$};\node at (5,-2) [place] (secondL1_1){$\frac{1}{\var{1}-\lani{1}{1}}$};\node at (5,-4) [place] (secondL1_2){$\frac{1}{\var{1}-\lani{1}{2}}$};\node at (5,-6) [place] (secondL2_1){$\frac{1}{\var{2}-\lani{2}{1}}$};\node at (5,-8) [place] (secondL2_2){$\frac{1}{\var{2}-\lani{2}{2}}$};\node at (10,-2) [place] (third_1){$\prod$};\node at (10,-4) [place] (third_2){$\prod$};\node at (10,-6) [place] (third_3){$\prod$};\node at (10,-8) [place] (third_4){$\prod$};\node at (15,-5) [placeInOut] (output){$\bSigma$};\draw[->] (first_1)--(secondL1_1) node[above,sloped,pos=0.75] { };\draw[->] (first_1)--(secondL1_2) node[above,sloped,pos=0.75] { };\draw[->] (first_2)--(secondL2_1) node[above,sloped,pos=0.75] { };\draw[->] (first_2)--(secondL2_2) node[above,sloped,pos=0.75] { };\draw[->] (secondL1_1)--(third_1) node[above,sloped,pos=0.25] {};\draw[->] (secondL1_1)--(third_2) node[above,sloped,pos=0.25] {};\draw[->] (secondL1_2)--(third_3) node[above,sloped,pos=0.25] {};\draw[->] (secondL1_2)--(third_4) node[above,sloped,pos=0.25] {};\draw[->] (secondL2_1)--(third_1) node[above,sloped,pos=0.25] {};\draw[->] (secondL2_2)--(third_2) node[above,sloped,pos=0.25] {};\draw[->] (secondL2_1)--(third_3) node[above,sloped,pos=0.25] {};\draw[->] (secondL2_2)--(third_4) node[above,sloped,pos=0.25] {};\draw[->] (third_1)--(output) node[above,sloped,pos=0.25] {1};\draw[->] (third_2)--(output) node[above,sloped,pos=0.25] {1};\draw[->] (third_3)--(output) node[above,sloped,pos=0.25] {1};\draw[->] (third_4)--(output) node[above,sloped,pos=0.25] {1};\end{tikzpicture}} \caption{Equivalent NN representation of the numerator $\bn_{\textrm{lag}}$.}\end{center}\end{figure}

\newpage \subsection{Function \#4 (${\ord=3}$ variables, tensor size: 500 \textbf{KB})} $$\frac{1}{3} \sum_{i=1}^3 \sin(\pi x_i/2)^2$$ \subsubsection{Setup and results overview}\begin{itemize}\item Reference: L/al. 2024, \cite{Liu:2025}\item Domain: $\mathbb{R}$\item Tensor size: 500 \textbf{KB} ($40^{3}$ points)\item Bounds: $ \left(\begin{array}{cc} -1 & 1 \end{array}\right) \times \left(\begin{array}{cc} -1 & 1 \end{array}\right) \times \left(\begin{array}{cc} -1 & 1 \end{array}\right)$ \end{itemize} \begin{table}[H] \centering \begin{tabular}{llllll}
$\#$ & Alg. & Parameters & Dim. & CPU [s] & RMSE \\ 
\hline 
$\mathbf{\#4}$ & A/G/P-V 2025 (A1) & $1 \cdot 10^{-06},3$ & $1.7 \cdot 10^{03}$ & $\mathbf{0.028}$ & $3.9 \cdot 10^{-08}$ \\ 
 & A/G/P-V 2025 (A2) & $1 \cdot 10^{-15},1$ & $1.7 \cdot 10^{03}$ & $15$ & $1 \cdot 10^{-07}$ \\ 
 & MDSPACK v1.1.0 & $0.0001,2$ & $1.7 \cdot 10^{03}$ & $0.029$ & $3.9 \cdot 10^{-08}$ \\ 
 & P/P 2025 & $1,1,50,0.01,4,12,9$ & $\mathbf{2.2 \cdot 10^{02}}$ & $11$ & $9.3 \cdot 10^{-06}$ \\ 
 & C-R/B/G 2023 & $0.001,20$ & $3.6 \cdot 10^{03}$ & $11$ & $1.2 \cdot 10^{-11}$ \\ 
 & B/G 2025 & $1 \cdot 10^{-09},20,2$ & $6.4 \cdot 10^{03}$ & $3.2$ & $\mathbf{9.3 \cdot 10^{-12}}$ \\ 
 & TensorFlow & $$ & $3.2 \cdot 10^{02}$ & $2.8 \cdot 10^{02}$ & $0.0031$ \\ 
\hline 
\end{tabular} \caption{Function \#4: best model configuration and performances per methods.} \end{table}\begin{figure}[H] \centering  \includegraphics[width=\textwidth]{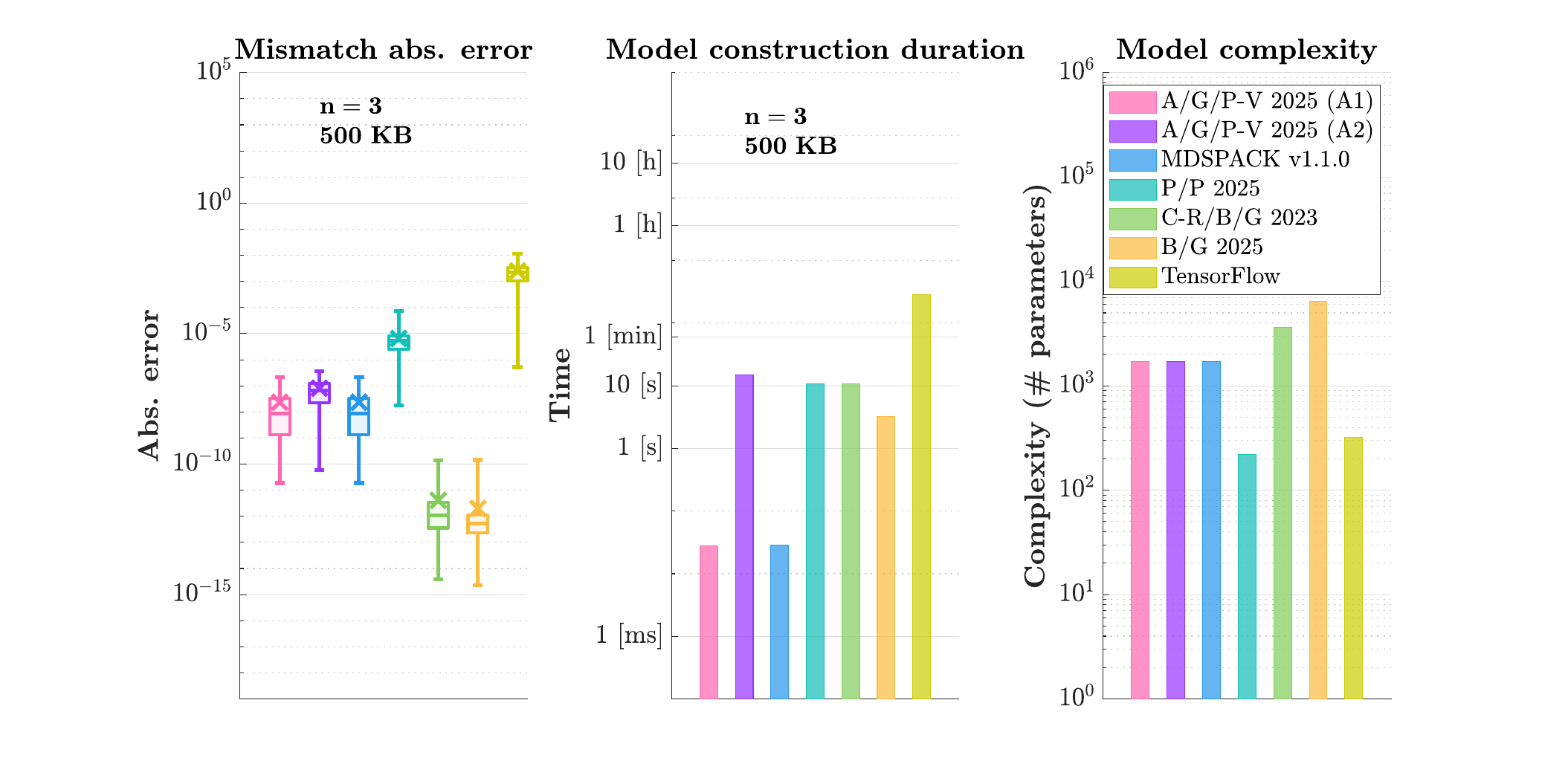} \caption{Function \#4: graphical view of the best model performances.} \end{figure}\begin{figure}[H] \centering  \includegraphics[width=\textwidth]{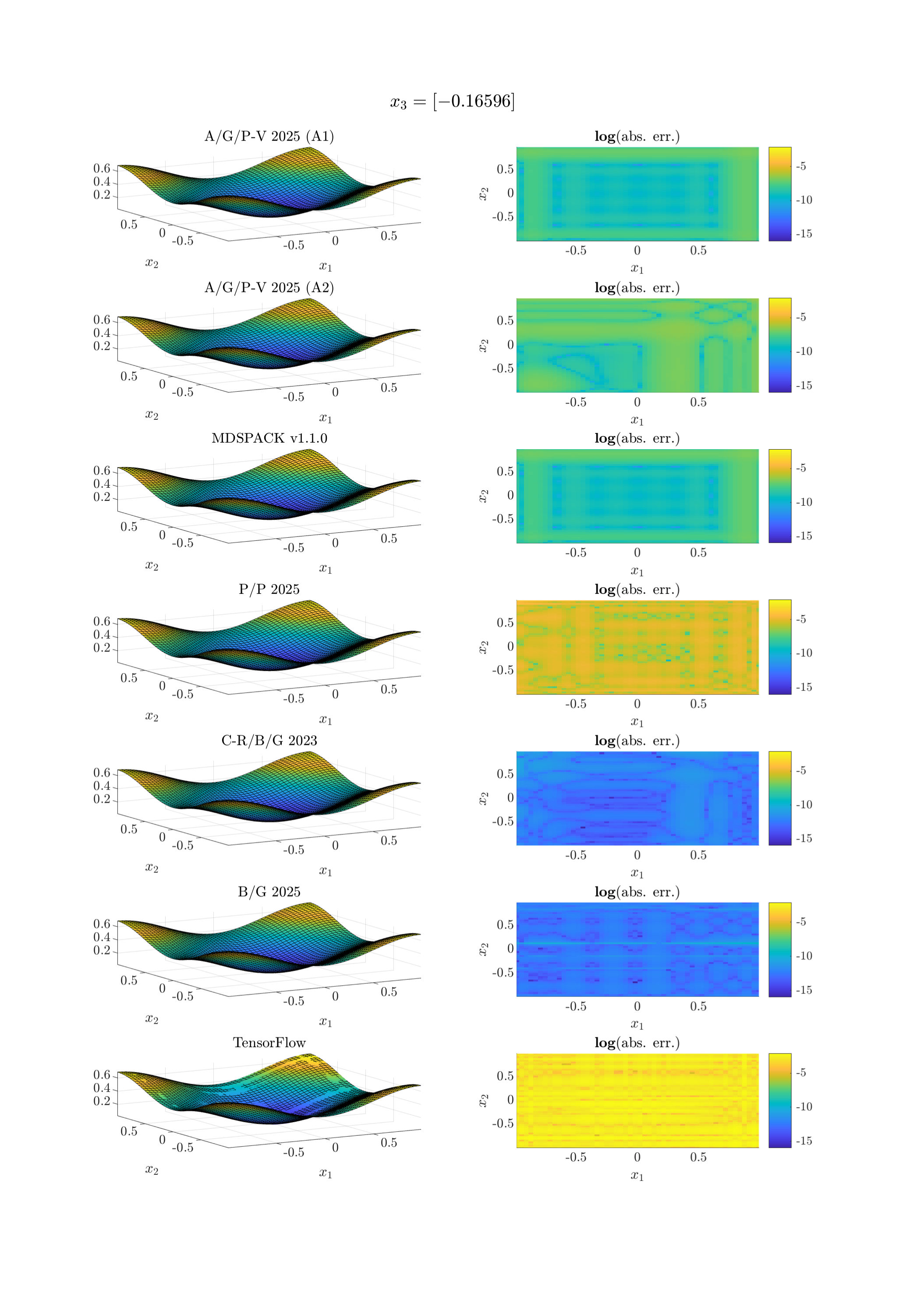} \caption{Function \#4: left side, evaluation of the original (mesh) vs. approximated (coloured surface) and right side, absolute errors (in log-scale).} \end{figure}\subsubsection{mLF detailed informations (M1)} \noindent \textbf{Right interpolation points}: $k_l=\left(\begin{array}{ccc} 7 & 7 & 7 \end{array}\right)$, where $l=1,\cdots,\ord$.$$ \begin{array}{rcl}\lan{1} &\in& \IC^{7} \text{ , linearly spaced between bounds}\\\lan{2} &\in& \IC^{7} \text{ , linearly spaced between bounds}\\\lan{3} &\in& \IC^{7} \text{ , linearly spaced between bounds}\\\end{array} $$\noindent \textbf{$\ord$-D Loewner matrix, barycentric weights and Lagrangian basis}:$$ \begin{array}{rcl}\IL & \in & \IC^{343 \times 343}\\\bc & \in & \IC^{343}\\\bw & \in & \IC^{343}\\\bc\odot \bw & \in & \IC^{343}\\\mathbf{Lag}(\var{1},\var{2},\var{3}) & \in & \IC^{343}\\\end{array} $$

\newpage \subsection{Function \#5 (${\ord=4}$ variables, tensor size: 19.5 \textbf{MB})} $$\mathrm{exp}\left(1/2 \left( \sin(\pi(\var{1}^2+\var{2}^2) + \sin(\pi(\var{3}^2+\var{4}^2) \right) \right)$$ \subsubsection{Setup and results overview}\begin{itemize}\item Reference: L/al. 2024, \cite{Liu:2025}\item Domain: $\mathbb{R}$\item Tensor size: 19.5 \textbf{MB} ($40^{4}$ points)\item Bounds: $ \left(\begin{array}{cc} -1 & 1 \end{array}\right) \times \left(\begin{array}{cc} -1 & 1 \end{array}\right) \times \left(\begin{array}{cc} -1 & 1 \end{array}\right) \times \left(\begin{array}{cc} -1 & 1 \end{array}\right)$ \end{itemize} \begin{table}[H] \centering \begin{tabular}{llllll}
$\#$ & Alg. & Parameters & Dim. & CPU [s] & RMSE \\ 
\hline 
$\mathbf{\#5}$ & A/G/P-V 2025 (A1) & $0.5,2$ & $3.8 \cdot 10^{03}$ & $\mathbf{0.57}$ & $0.086$ \\ 
 & A/G/P-V 2025 (A2) & $1 \cdot 10^{-15},1$ & $\mathbf{6}$ & $3.4 \cdot 10^{02}$ & $1.3$ \\ 
 & MDSPACK v1.1.0 & $0.01,1$ & $3.8 \cdot 10^{03}$ & $0.58$ & $0.09$ \\ 
 & P/P 2025 & $1,0.95,50,0.01,6,12,13$ & $4.7 \cdot 10^{02}$ & $3.9 \cdot 10^{03}$ & $0.027$ \\ 
 & C-R/B/G 2023 & $NaN$ & $NaN$ & $NaN$ & $NaN$ \\ 
 & B/G 2025 & $0.001,20,4$ & $3.7 \cdot 10^{05}$ & $3.8 \cdot 10^{03}$ & $\mathbf{0.012}$ \\ 
 & TensorFlow & $$ & $3.8 \cdot 10^{02}$ & $7.8 \cdot 10^{03}$ & $0.13$ \\ 
\hline 
\end{tabular} \caption{Function \#5: best model configuration and performances per methods.} \end{table}\begin{figure}[H] \centering  \includegraphics[width=\textwidth]{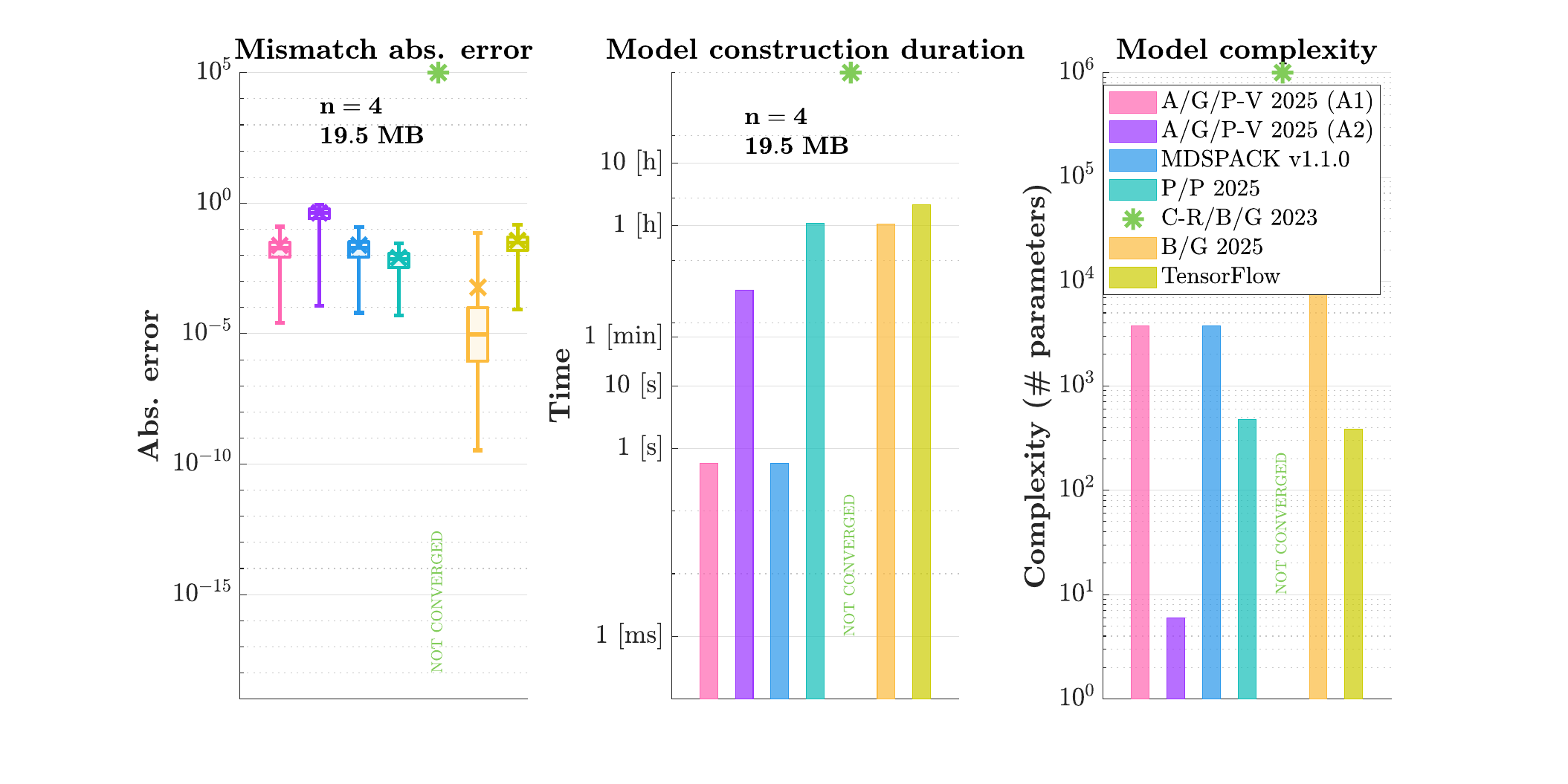} \caption{Function \#5: graphical view of the best model performances.} \end{figure}\begin{figure}[H] \centering  \includegraphics[width=\textwidth]{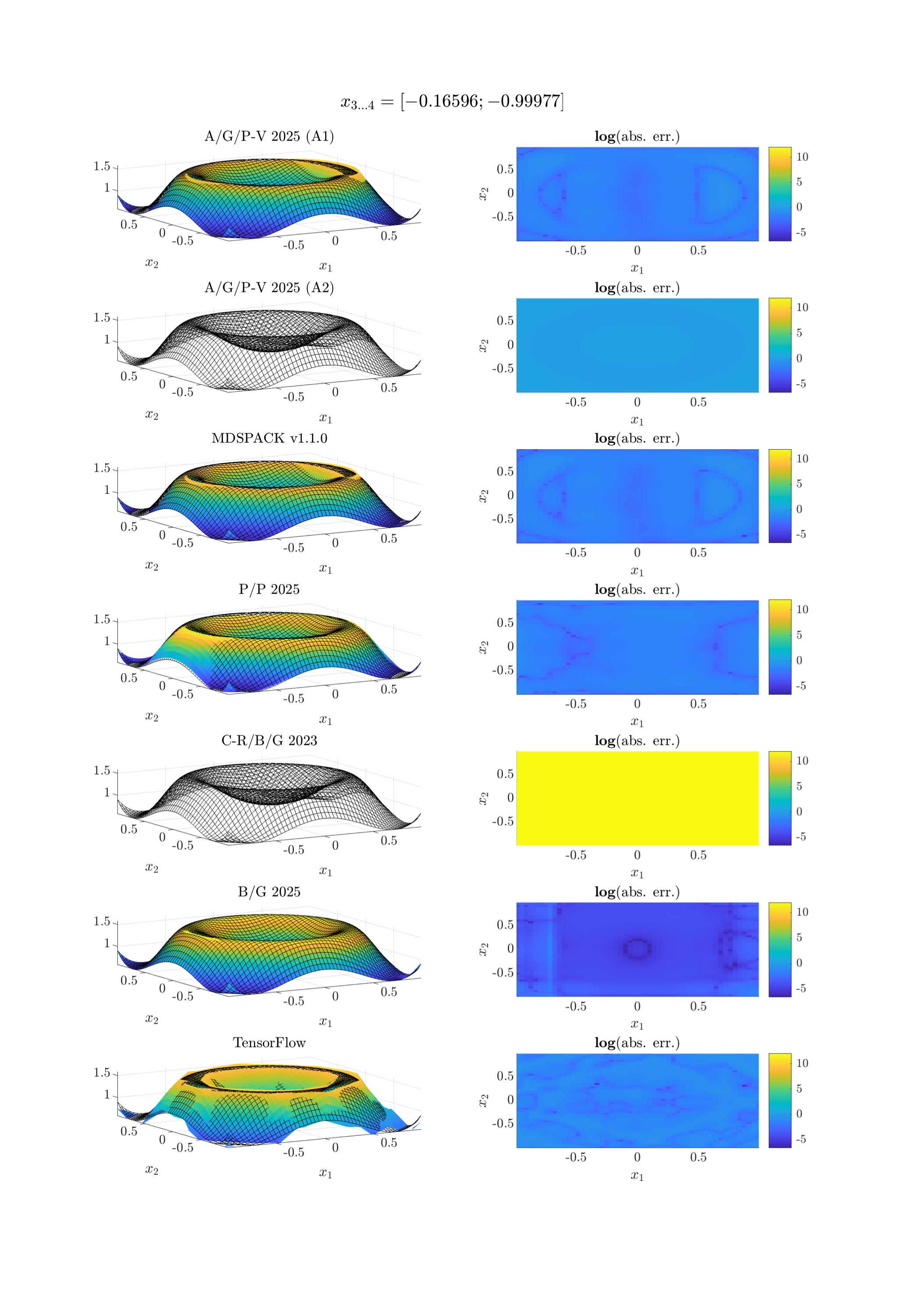} \caption{Function \#5: left side, evaluation of the original (mesh) vs. approximated (coloured surface) and right side, absolute errors (in log-scale).} \end{figure}\subsubsection{mLF detailed informations (M1)} \noindent \textbf{Right interpolation points}: $k_l=\left(\begin{array}{cccc} 5 & 5 & 5 & 5 \end{array}\right)$, where $l=1,\cdots,\ord$.$$ \begin{array}{rcl}\lan{1} &\in& \IC^{5} \text{ , linearly spaced between bounds}\\\lan{2} &\in& \IC^{5} \text{ , linearly spaced between bounds}\\\lan{3} &\in& \IC^{5} \text{ , linearly spaced between bounds}\\\lan{4} &\in& \IC^{5} \text{ , linearly spaced between bounds}\\\end{array} $$\noindent \textbf{$\ord$-D Loewner matrix, barycentric weights and Lagrangian basis}:$$ \begin{array}{rcl}\IL & \in & \IC^{625 \times 625}\\\bc & \in & \IC^{625}\\\bw & \in & \IC^{625}\\\bc\odot \bw & \in & \IC^{625}\\\mathbf{Lag}(\var{1},\var{2},\var{3},\var{4}) & \in & \IC^{625}\\\end{array} $$

\newpage \subsection{Function \#6 (${\ord=2}$ variables, tensor size: 12.5 \textbf{KB})} $$\frac{\mathrm{exp}\left(\var{1} \var{2}\right)}{(\var{1}^2-1.44)(\var{2}^2-1.44)}$$ \subsubsection{Setup and results overview}\begin{itemize}\item Reference: A/al. 2021 (A.5.1), \cite{Austin:2021}\item Domain: $\mathbb{R}$\item Tensor size: 12.5 \textbf{KB} ($40^{2}$ points)\item Bounds: $ \left(\begin{array}{cc} -1 & 1 \end{array}\right) \times \left(\begin{array}{cc} -1 & 1 \end{array}\right)$ \end{itemize} \begin{table}[H] \centering \begin{tabular}{llllll}
$\#$ & Alg. & Parameters & Dim. & CPU [s] & RMSE \\ 
\hline 
$\mathbf{\#6}$ & A/G/P-V 2025 (A1) & $1 \cdot 10^{-06},3$ & $1 \cdot 10^{02}$ & $\mathbf{0.0089}$ & $0.00021$ \\ 
 & A/G/P-V 2025 (A2) & $1 \cdot 10^{-15},1$ & $\mathbf{36}$ & $0.39$ & $0.23$ \\ 
 & MDSPACK v1.1.0 & $1 \cdot 10^{-10},5$ & $1 \cdot 10^{02}$ & $0.0096$ & $0.00021$ \\ 
 & P/P 2025 & $1,0.95,50,0.01,6,12,13$ & $3.2 \cdot 10^{02}$ & $1.1$ & $0.0043$ \\ 
 & C-R/B/G 2023 & $1 \cdot 10^{-09},20$ & $1.4 \cdot 10^{02}$ & $0.037$ & $4.2 \cdot 10^{-10}$ \\ 
 & B/G 2025 & $1 \cdot 10^{-09},20,4$ & $2.6 \cdot 10^{02}$ & $0.41$ & $\mathbf{1.4 \cdot 10^{-10}}$ \\ 
 & TensorFlow & $$ & $2.6 \cdot 10^{02}$ & $14$ & $0.11$ \\ 
\hline 
\end{tabular} \caption{Function \#6: best model configuration and performances per methods.} \end{table}\begin{figure}[H] \centering  \includegraphics[width=\textwidth]{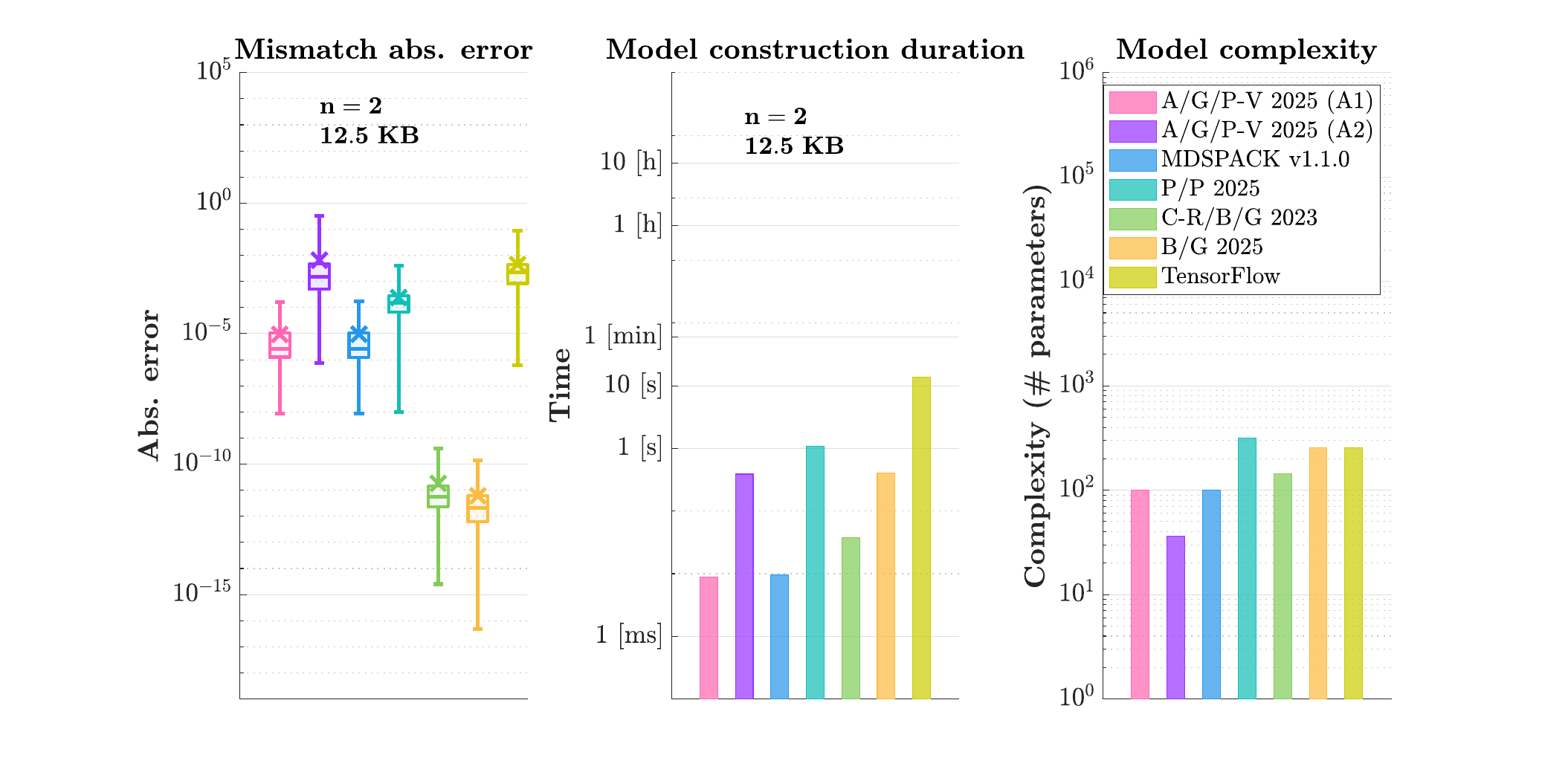} \caption{Function \#6: graphical view of the best model performances.} \end{figure}\begin{figure}[H] \centering  \includegraphics[width=\textwidth]{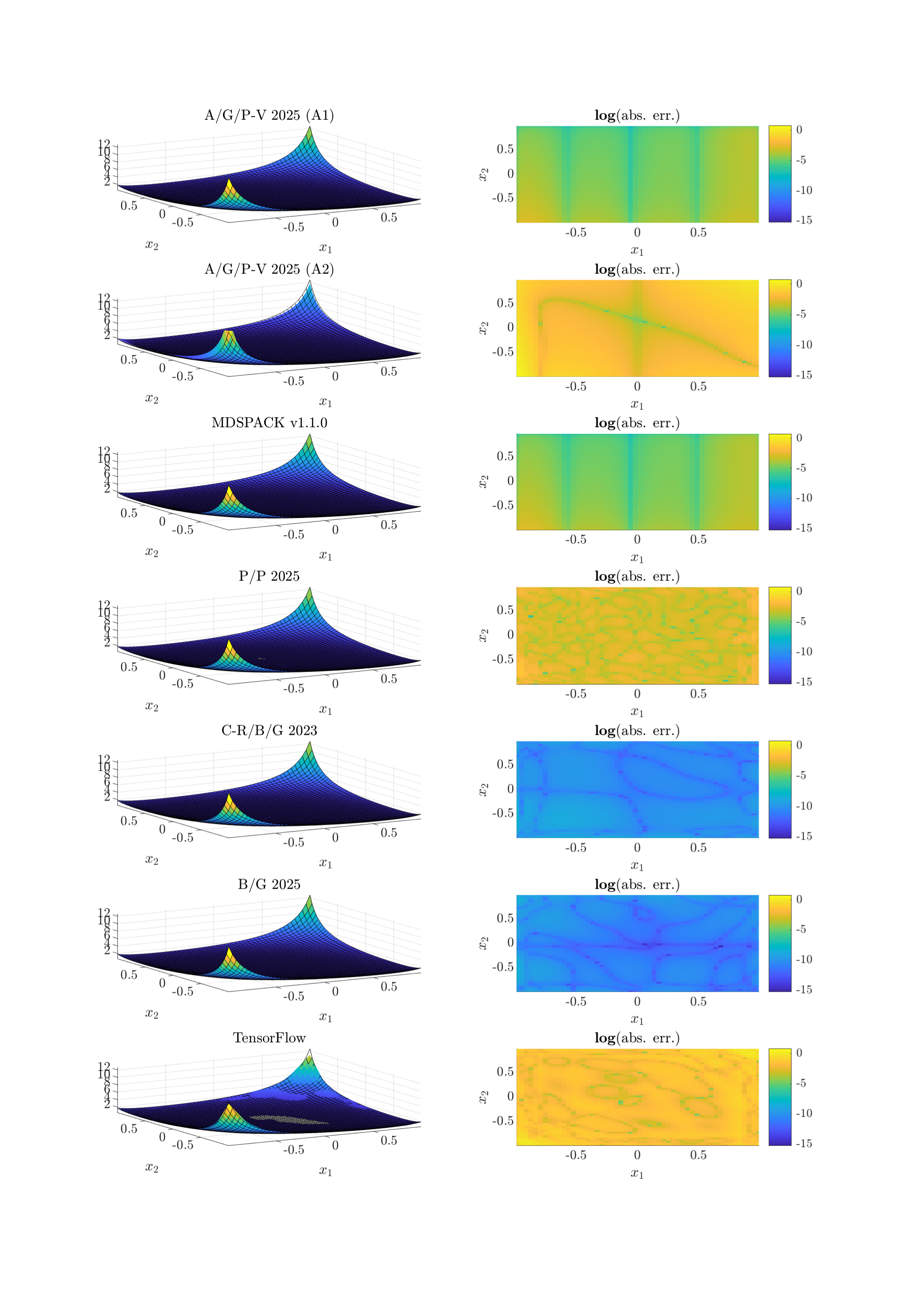} \caption{Function \#6: left side, evaluation of the original (mesh) vs. approximated (coloured surface) and right side, absolute errors (in log-scale).} \end{figure}\subsubsection{mLF detailed informations (M1)} \noindent \textbf{Right interpolation points} ($k_l=\left(\begin{array}{cc} 5 & 5 \end{array}\right)$, where $l=1,\cdots,\ord$):$$ \begin{array}{rcl}\lan{1} &=& \left(\begin{array}{ccccc} -1 & -\frac{11}{19} & -\frac{1}{19} & \frac{9}{19} & 1 \end{array}\right)\\\lan{2} &=& \left(\begin{array}{ccccc} -1 & -\frac{11}{19} & -\frac{1}{19} & \frac{9}{19} & 1 \end{array}\right)\\\end{array} $$\noindent \textbf{Lagrangian weights}: $$\left(\begin{array}{ccc} \bc & \bw & \bc\odot\bw\\ 0.07376 & 14.04 & 1.036\\ -0.6906 & 3.67 & -2.535\\ 1.437 & 1.667 & 2.396\\ -0.9262 & 1.164 & -1.078\\ 0.0965 & 1.9 & 0.1834\\ -0.6898 & 3.67 & -2.532\\ 6.027 & 1.145 & 6.903\\ -11.6 & 0.6493 & -7.53\\ 6.968 & 0.566 & 3.944\\ -0.6818 & 1.153 & -0.7861\\ 1.437 & 1.667 & 2.395\\ -11.6 & 0.6493 & -7.532\\ 20.31 & 0.4855 & 9.86\\ -11.15 & 0.5583 & -6.222\\ 1.0 & 1.5 & 1.5\\ -0.927 & 1.164 & -1.079\\ 6.974 & 0.566 & 3.947\\ -11.15 & 0.5583 & -6.225\\ 5.574 & 0.8469 & 4.721\\ -0.4542 & 3.002 & -1.364\\ 0.0965 & 1.9 & 0.1834\\ -0.6818 & 1.153 & -0.7861\\ 1.0 & 1.5 & 1.5\\ -0.4542 & 3.002 & -1.364\\ 0.03329 & 14.04 & 0.4674 \end{array}\right)$$\noindent \textbf{Lagrangian form} (basis, numerator and denominator coefficients):$$\left(\begin{array}{ccc}\mathcal{B}_\textrm{lag}(\var{1},\var{2}) & \bN_\textrm{lag} &\bD_\textrm{lag}\end{array}\right) =$$ $$\left(\begin{array}{ccc} \left(\var{1}+1.0\right)\,\left(\var{2}+1.0\right) & 1.036 & 0.07376\\ \left(\var{1}+1.0\right)\,\left(\var{2}+0.5789\right) & -2.535 & -0.6906\\ \left(\var{1}+1.0\right)\,\left(\var{2}+0.05263\right) & 2.396 & 1.437\\ \left(\var{1}+1.0\right)\,\left(\var{2}-0.4737\right) & -1.078 & -0.9262\\ \left(\var{1}+1.0\right)\,\left(\var{2}-1.0\right) & 0.1834 & 0.0965\\ \left(\var{2}+1.0\right)\,\left(\var{1}+0.5789\right) & -2.532 & -0.6898\\ \left(\var{1}+0.5789\right)\,\left(\var{2}+0.5789\right) & 6.903 & 6.027\\ \left(\var{1}+0.5789\right)\,\left(\var{2}+0.05263\right) & -7.53 & -11.6\\ \left(\var{1}+0.5789\right)\,\left(\var{2}-0.4737\right) & 3.944 & 6.968\\ \left(\var{2}-1.0\right)\,\left(\var{1}+0.5789\right) & -0.7861 & -0.6818\\ \left(\var{2}+1.0\right)\,\left(\var{1}+0.05263\right) & 2.395 & 1.437\\ \left(\var{2}+0.5789\right)\,\left(\var{1}+0.05263\right) & -7.532 & -11.6\\ \left(\var{1}+0.05263\right)\,\left(\var{2}+0.05263\right) & 9.86 & 20.31\\ \left(\var{1}+0.05263\right)\,\left(\var{2}-0.4737\right) & -6.222 & -11.15\\ \left(\var{2}-1.0\right)\,\left(\var{1}+0.05263\right) & 1.5 & 1.0\\ \left(\var{2}+1.0\right)\,\left(\var{1}-0.4737\right) & -1.079 & -0.927\\ \left(\var{2}+0.5789\right)\,\left(\var{1}-0.4737\right) & 3.947 & 6.974\\ \left(\var{2}+0.05263\right)\,\left(\var{1}-0.4737\right) & -6.225 & -11.15\\ \left(\var{1}-0.4737\right)\,\left(\var{2}-0.4737\right) & 4.721 & 5.574\\ \left(\var{2}-1.0\right)\,\left(\var{1}-0.4737\right) & -1.364 & -0.4542\\ \left(\var{1}-1.0\right)\,\left(\var{2}+1.0\right) & 0.1834 & 0.0965\\ \left(\var{1}-1.0\right)\,\left(\var{2}+0.5789\right) & -0.7861 & -0.6818\\ \left(\var{1}-1.0\right)\,\left(\var{2}+0.05263\right) & 1.5 & 1.0\\ \left(\var{1}-1.0\right)\,\left(\var{2}-0.4737\right) & -1.364 & -0.4542\\ \left(\var{1}-1.0\right)\,\left(\var{2}-1.0\right) & 0.4674 & 0.03329 \end{array}\right).$$\noindent The corresponding function is:$$\begin{array}{rcl}\bG_{\textrm{lag}}(\var{1},\var{2}) &=& \dfrac{\bn_{\textrm{lag}}(\var{1},\var{2})}{\bd_{\textrm{lag}}(\var{1},\var{2})}\\ && \\&=& \dfrac{\sum_{\textrm{row}} \bN_\textrm{lag} \odot\mathcal{B}^{-1}_\textrm{lag}(\var{1},\var{2})}{\sum_{\textrm{row}} \bD_\textrm{lag} \odot\mathcal{B}^{-1}_\textrm{lag}(\var{1},\var{2})}, \end{array}$$\noindent where,\\$\bn_{\textrm{lag}}(\var{1},\var{2}) = 0.001349\,{\var{1}}^4\,{\var{2}}^4-4.372 \cdot 10^{-5}\,{\var{1}}^4\,{\var{2}}^3+0.0004675\,{\var{1}}^4\,{\var{2}}^2-0.0003926\,{\var{1}}^4\,\var{2}-3.637 \cdot 10^{-5}\,{\var{1}}^4-3.703 \cdot 10^{-6}\,{\var{1}}^3\,{\var{2}}^4+0.01623\,{\var{1}}^3\,{\var{2}}^3-9.082 \cdot 10^{-5}\,{\var{1}}^3\,{\var{2}}^2+0.0006661\,{\var{1}}^3\,\var{2}-0.0006219\,{\var{1}}^3-2.66 \cdot 10^{-6}\,{\var{1}}^2\,{\var{2}}^4+2.761 \cdot 10^{-5}\,{\var{1}}^2\,{\var{2}}^3+0.09643\,{\var{1}}^2\,{\var{2}}^2+0.0002922\,{\var{1}}^2\,\var{2}+9.731 \cdot 10^{-5}\,{\var{1}}^2+1.055 \cdot 10^{-6}\,\var{1}\,{\var{2}}^4-3.302 \cdot 10^{-5}\,\var{1}\,{\var{2}}^3+3.476 \cdot 10^{-5}\,\var{1}\,{\var{2}}^2+0.3215\,\var{1}\,\var{2}+0.0005703\,\var{1}+6.223 \cdot 10^{-8}\,{\var{2}}^4-1.787 \cdot 10^{-6}\,{\var{2}}^3+1.608 \cdot 10^{-6}\,{\var{2}}^2+6.719 \cdot 10^{-6}\,\var{2}+0.4823$ \\~~\\$\bd_{\textrm{lag}}(\var{1},\var{2}) = 0.0161\,{\var{1}}^4\,{\var{2}}^4-0.0004017\,{\var{1}}^4\,{\var{2}}^3-0.02298\,{\var{1}}^4\,{\var{2}}^2+0.0005784\,{\var{1}}^4\,\var{2}-0.0002892\,{\var{1}}^4+9.806 \cdot 10^{-6}\,{\var{1}}^3\,{\var{2}}^4-0.1609\,{\var{1}}^3\,{\var{2}}^3+0.0006639\,{\var{1}}^3\,{\var{2}}^2+0.2317\,{\var{1}}^3\,\var{2}-0.0009761\,{\var{1}}^3-0.02309\,{\var{1}}^2\,{\var{2}}^4+0.0004373\,{\var{1}}^2\,{\var{2}}^3+0.5153\,{\var{1}}^2\,{\var{2}}^2-0.0006299\,{\var{1}}^2\,\var{2}-0.6941\,{\var{1}}^2-1.496 \cdot 10^{-5}\,\var{1}\,{\var{2}}^4+0.2313\,\var{1}\,{\var{2}}^3-0.0006336\,\var{1}\,{\var{2}}^2-0.3331\,\var{1}\,\var{2}+0.0009434\,\var{1}-1.073 \cdot 10^{-6}\,{\var{2}}^4-1.727 \cdot 10^{-5}\,{\var{2}}^3-0.6944\,{\var{2}}^2+2.489 \cdot 10^{-5}\,\var{2}+1.0$ \\~~\\\noindent \textbf{Monomial form} (basis, numerator and denominator coefficients - entries $<10^{-12}$ removed):$$\left(\begin{array}{ccc}\mathcal{B}_\textrm{mon}(\var{1},\var{2}) & \bN_\textrm{mon} &\bD_\textrm{mon}\end{array}\right) =$$ $$\left(\begin{array}{ccc} {\var{1}}^4\,{\var{2}}^4 & 0.001349 & 0.0161\\ {\var{1}}^4\,{\var{2}}^3 & -4.372 \cdot 10^{-5} & -0.0004017\\ {\var{1}}^4\,{\var{2}}^2 & 0.0004675 & -0.02298\\ {\var{1}}^4\,\var{2} & -0.0003926 & 0.0005784\\ {\var{1}}^4 & -3.637 \cdot 10^{-5} & -0.0002892\\ {\var{1}}^3\,{\var{2}}^4 & -3.703 \cdot 10^{-6} & 9.806 \cdot 10^{-6}\\ {\var{1}}^3\,{\var{2}}^3 & 0.01623 & -0.1609\\ {\var{1}}^3\,{\var{2}}^2 & -9.082 \cdot 10^{-5} & 0.0006639\\ {\var{1}}^3\,\var{2} & 0.0006661 & 0.2317\\ {\var{1}}^3 & -0.0006219 & -0.0009761\\ {\var{1}}^2\,{\var{2}}^4 & -2.66 \cdot 10^{-6} & -0.02309\\ {\var{1}}^2\,{\var{2}}^3 & 2.761 \cdot 10^{-5} & 0.0004373\\ {\var{1}}^2\,{\var{2}}^2 & 0.09643 & 0.5153\\ {\var{1}}^2\,\var{2} & 0.0002922 & -0.0006299\\ {\var{1}}^2 & 9.731 \cdot 10^{-5} & -0.6941\\ \var{1}\,{\var{2}}^4 & 1.055 \cdot 10^{-6} & -1.496 \cdot 10^{-5}\\ \var{1}\,{\var{2}}^3 & -3.302 \cdot 10^{-5} & 0.2313\\ \var{1}\,{\var{2}}^2 & 3.476 \cdot 10^{-5} & -0.0006336\\ \var{1}\,\var{2} & 0.3215 & -0.3331\\ \var{1} & 0.0005703 & 0.0009434\\ {\var{2}}^4 & 6.223 \cdot 10^{-8} & -1.073 \cdot 10^{-6}\\ {\var{2}}^3 & -1.787 \cdot 10^{-6} & -1.727 \cdot 10^{-5}\\ {\var{2}}^2 & 1.608 \cdot 10^{-6} & -0.6944\\ \var{2} & 6.719 \cdot 10^{-6} & 2.489 \cdot 10^{-5}\\ 1.0 & 0.4823 & 1.0 \end{array}\right)$$\noindent The corresponding function is:$$\begin{array}{rcl}\bG_{\textrm{mon}}(\var{1},\var{2}) &=& \dfrac{\bn_{\textrm{mon}}(\var{1},\var{2})}{\bd_{\textrm{mon}}(\var{1},\var{2})}\\ && \\&=& \dfrac{\sum_{\textrm{row}} \bN_\textrm{mon} \odot \mathcal{B}_\textrm{mon}(\var{1},\var{2})}{\sum_{\textrm{row}} \bD_\textrm{mon} \odot\mathcal{B}_\textrm{mon}(\var{1},\var{2})},  \end{array}$$\noindent where,\\$\bn_{\textrm{mon}}(\var{1},\var{2}) = 0.001349\,{\var{1}}^4\,{\var{2}}^4-4.372 \cdot 10^{-5}\,{\var{1}}^4\,{\var{2}}^3+0.0004675\,{\var{1}}^4\,{\var{2}}^2-0.0003926\,{\var{1}}^4\,\var{2}-3.637 \cdot 10^{-5}\,{\var{1}}^4-3.703 \cdot 10^{-6}\,{\var{1}}^3\,{\var{2}}^4+0.01623\,{\var{1}}^3\,{\var{2}}^3-9.082 \cdot 10^{-5}\,{\var{1}}^3\,{\var{2}}^2+0.0006661\,{\var{1}}^3\,\var{2}-0.0006219\,{\var{1}}^3-2.66 \cdot 10^{-6}\,{\var{1}}^2\,{\var{2}}^4+2.761 \cdot 10^{-5}\,{\var{1}}^2\,{\var{2}}^3+0.09643\,{\var{1}}^2\,{\var{2}}^2+0.0002922\,{\var{1}}^2\,\var{2}+9.731 \cdot 10^{-5}\,{\var{1}}^2+1.055 \cdot 10^{-6}\,\var{1}\,{\var{2}}^4-3.302 \cdot 10^{-5}\,\var{1}\,{\var{2}}^3+3.476 \cdot 10^{-5}\,\var{1}\,{\var{2}}^2+0.3215\,\var{1}\,\var{2}+0.0005703\,\var{1}+6.223 \cdot 10^{-8}\,{\var{2}}^4-1.787 \cdot 10^{-6}\,{\var{2}}^3+1.608 \cdot 10^{-6}\,{\var{2}}^2+6.719 \cdot 10^{-6}\,\var{2}+0.4823$ \\~~\\$\bd_{\textrm{mon}}(\var{1},\var{2}) = 0.0161\,{\var{1}}^4\,{\var{2}}^4-0.0004017\,{\var{1}}^4\,{\var{2}}^3-0.02298\,{\var{1}}^4\,{\var{2}}^2+0.0005784\,{\var{1}}^4\,\var{2}-0.0002892\,{\var{1}}^4+9.806 \cdot 10^{-6}\,{\var{1}}^3\,{\var{2}}^4-0.1609\,{\var{1}}^3\,{\var{2}}^3+0.0006639\,{\var{1}}^3\,{\var{2}}^2+0.2317\,{\var{1}}^3\,\var{2}-0.0009761\,{\var{1}}^3-0.02309\,{\var{1}}^2\,{\var{2}}^4+0.0004373\,{\var{1}}^2\,{\var{2}}^3+0.5153\,{\var{1}}^2\,{\var{2}}^2-0.0006299\,{\var{1}}^2\,\var{2}-0.6941\,{\var{1}}^2-1.496 \cdot 10^{-5}\,\var{1}\,{\var{2}}^4+0.2313\,\var{1}\,{\var{2}}^3-0.0006336\,\var{1}\,{\var{2}}^2-0.3331\,\var{1}\,\var{2}+0.0009434\,\var{1}-1.073 \cdot 10^{-6}\,{\var{2}}^4-1.727 \cdot 10^{-5}\,{\var{2}}^3-0.6944\,{\var{2}}^2+2.489 \cdot 10^{-5}\,\var{2}+1.0$ \\~~\\\noindent \textbf{KST equivalent decoupling pattern} (Barycentric weights $\bc^{\var{l}}$): $$\begin{array}{rclll}\var{2}&: & \left(\begin{array}{ccccc} 0.7643 & 1.012 & 1.437 & 2.041 & 2.899\\ -7.156 & -8.839 & -11.6 & -15.35 & -20.48\\ 14.89 & 17.01 & 20.31 & 24.55 & 30.04\\ -9.598 & -10.22 & -11.15 & -12.27 & -13.64\\ 1.0 & 1.0 & 1.0 & 1.0 & 1.0 \end{array}\right)& \textrm{vec}(.) & := \textbf{Bary}(\var{2}) \\\var{1}&: & \left(\begin{array}{c} 0.0965\\ -0.6818\\ 1.0\\ -0.4542\\ 0.03329 \end{array}\right)& \textrm{vec}(.) \otimes \bone_{k_{2}} & := \textbf{Bary}(\var{1}) \\\end{array}$$~\\ Then, with the above notations, one defines the following univariate vector functions:  $$ \left\{ \begin{array}{rcl}\bPhi_{1}(\var{1}) &:=& \textbf{Bary}(\var{1}) \odot \mathbf{Lag}(\var{1}) \\\bPhi_{2}(\var{2}) &:=& \textbf{Bary}(\var{2}) \odot \mathbf{Lag}(\var{2}) \\\end{array} \right. $$\noindent The corresponding function is:$$\begin{array}{rcl}\bG_{\textrm{kst}}(\var{1},\var{2}) &=& \dfrac{\bn_{\textrm{kst}}(\var{1},\var{2})}{\bd_{\textrm{kst}}(\var{1},\var{2})}\\ && \\ &=& \dfrac{\sum_{\text{rows}} \bw \odot \bPhi_{1}(\var{1}) \odot \cdots \odot\bPhi_{2}(\var{2})}{\sum_{\text{rows}} \bPhi_{1}(\var{1}) \odot \cdots \odot\bPhi_{2}(\var{2})} . \end{array}$$~\\ \noindent \textbf{KST-like univariate functions} (equivalent scaled univariate functions $\bphi_{1,\cdots,2}$): $$\left\{\begin{array}{rcrcl}z_{1} &=&\bphi_{1}(\var{1}) &=& \cfrac{\bn_{1}}{\bd_{1}} \\z_{2} &=&\bphi_{2}(\var{2}) &=& \cfrac{\bn_{2}}{\bd_{2}} \\\end{array} \right. .$$\noindent where, \\ \noindent $\bn_{1}=0.004399\,{\var{1}}^4+0.05295\,{\var{1}}^3+0.3169\,{\var{1}}^2+1.054\,\var{1}+1.578$ and \\ \noindent $\bd_{1}=-0.02289\,{\var{1}}^4+0.2307\,{\var{1}}^3-0.6615\,{\var{1}}^2-0.3321\,\var{1}+1.0$, \\ \noindent $\bn_{2}=0.004415\,{\var{2}}^4-0.05305\,{\var{2}}^3+0.3172\,{\var{2}}^2-1.054\,\var{2}+1.578$ and \\ \noindent $\bd_{2}=-0.02285\,{\var{2}}^4-0.2305\,{\var{2}}^3-0.6615\,{\var{2}}^2+0.3319\,\var{2}+1.0$, \\

\newpage \subsection{Function \#7 (${\ord=2}$ variables, tensor size: 12.5 \textbf{KB})} $$\mathrm{log}(2.25-\var{1}^2-\var{2}^2)$$ \subsubsection{Setup and results overview}\begin{itemize}\item Reference: A/al. 2021 (A.5.2), \cite{Austin:2021}\item Domain: $\mathbb{R}$\item Tensor size: 12.5 \textbf{KB} ($40^{2}$ points)\item Bounds: $ \left(\begin{array}{cc} -1 & 1 \end{array}\right) \times \left(\begin{array}{cc} -1 & 1 \end{array}\right)$ \end{itemize} \begin{table}[H] \centering \begin{tabular}{llllll}
$\#$ & Alg. & Parameters & Dim. & CPU [s] & RMSE \\ 
\hline 
$\mathbf{\#7}$ & A/G/P-V 2025 (A1) & $1 \cdot 10^{-06},2$ & $3.2 \cdot 10^{02}$ & $0.012$ & $0.0013$ \\ 
 & A/G/P-V 2025 (A2) & $1 \cdot 10^{-15},3$ & $2.6 \cdot 10^{02}$ & $0.35$ & $0.0016$ \\ 
 & MDSPACK v1.1.0 & $1 \cdot 10^{-10},5$ & $3.2 \cdot 10^{02}$ & $\mathbf{0.0087}$ & $0.0014$ \\ 
 & P/P 2025 & $1,1,50,0.01,4,12,9$ & $\mathbf{1.8 \cdot 10^{02}}$ & $0.46$ & $0.00071$ \\ 
 & C-R/B/G 2023 & $1 \cdot 10^{-09},20$ & $7.2 \cdot 10^{02}$ & $0.23$ & $\mathbf{7 \cdot 10^{-11}}$ \\ 
 & B/G 2025 & $1 \cdot 10^{-09},20,4$ & $6.2 \cdot 10^{02}$ & $4.1$ & $4 \cdot 10^{-09}$ \\ 
 & TensorFlow & $$ & $2.6 \cdot 10^{02}$ & $14$ & $0.22$ \\ 
\hline 
\end{tabular} \caption{Function \#7: best model configuration and performances per methods.} \end{table}\begin{figure}[H] \centering  \includegraphics[width=\textwidth]{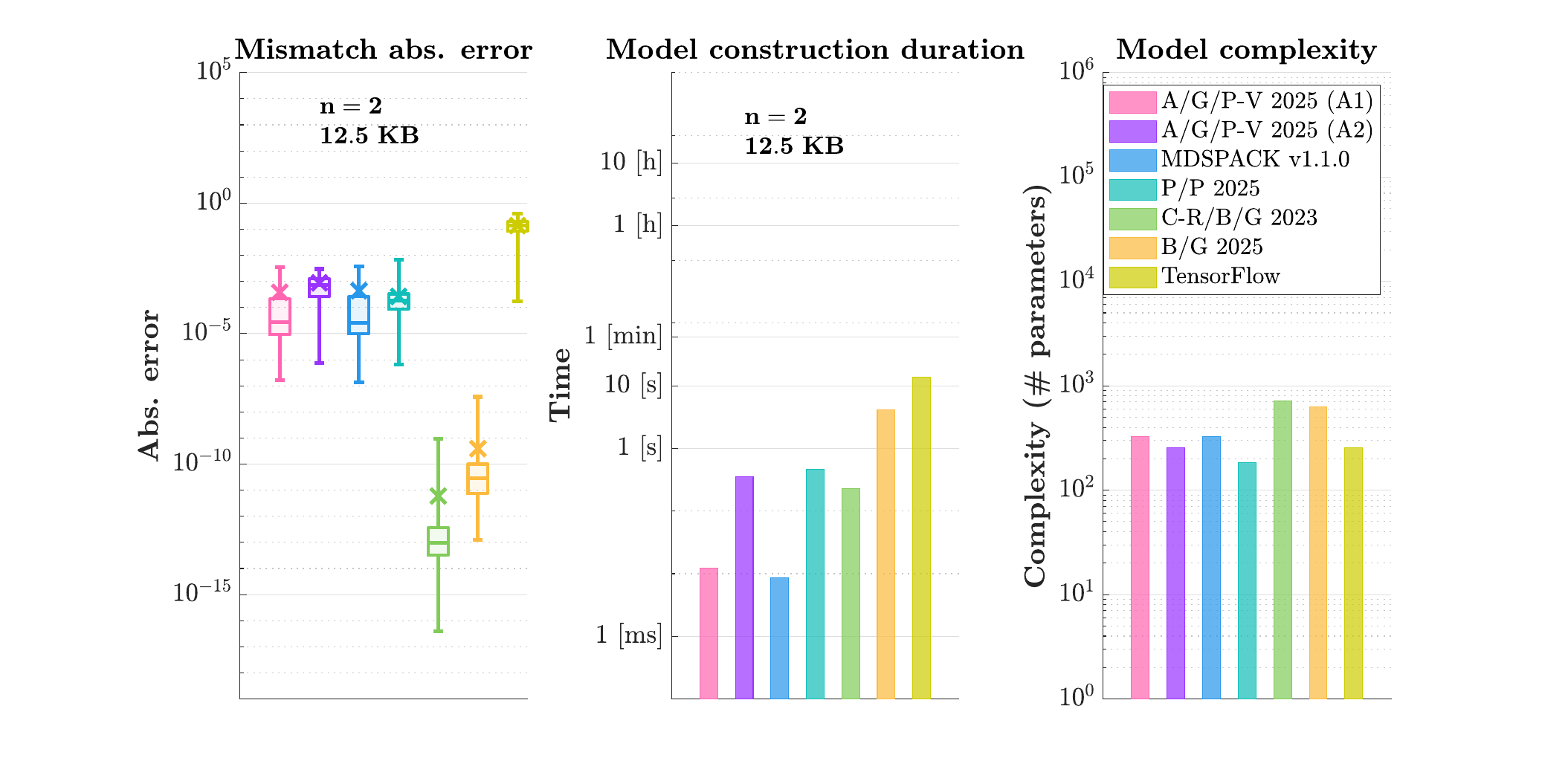} \caption{Function \#7: graphical view of the best model performances.} \end{figure}\begin{figure}[H] \centering  \includegraphics[width=\textwidth]{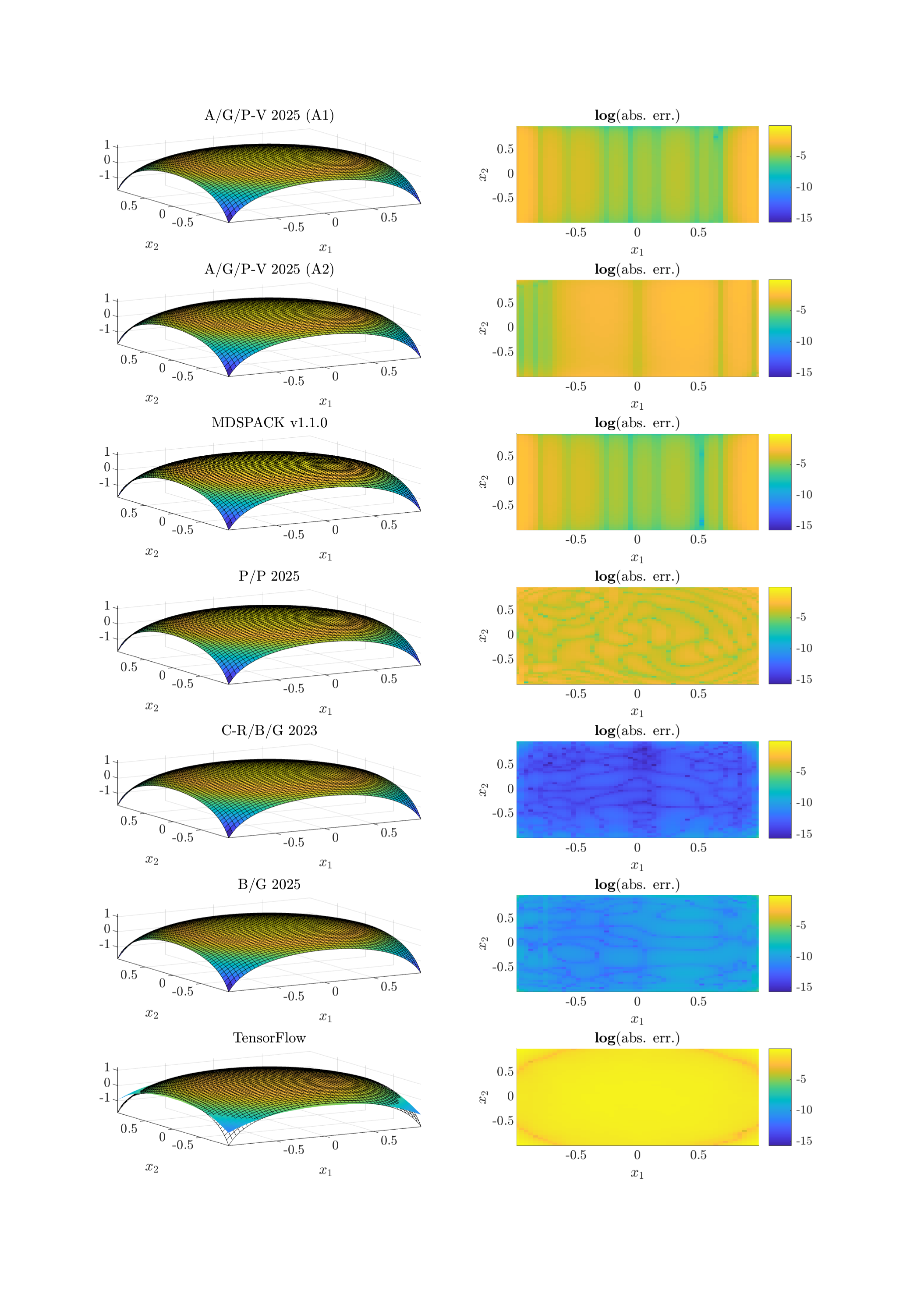} \caption{Function \#7: left side, evaluation of the original (mesh) vs. approximated (coloured surface) and right side, absolute errors (in log-scale).} \end{figure}\subsubsection{mLF detailed informations (M1)} \noindent \textbf{Right interpolation points}: $k_l=\left(\begin{array}{cc} 9 & 9 \end{array}\right)$, where $l=1,\cdots,\ord$.$$ \begin{array}{rcl}\lan{1} &\in& \IC^{9} \text{ , linearly spaced between bounds}\\\lan{2} &\in& \IC^{9} \text{ , linearly spaced between bounds}\\\end{array} $$\noindent \textbf{$\ord$-D Loewner matrix, barycentric weights and Lagrangian basis}:$$ \begin{array}{rcl}\IL & \in & \IC^{81 \times 81}\\\bc & \in & \IC^{81}\\\bw & \in & \IC^{81}\\\bc\odot \bw & \in & \IC^{81}\\\mathbf{Lag}(\var{1},\var{2}) & \in & \IC^{81}\\\end{array} $$

\newpage \subsection{Function \#8 (${\ord=2}$ variables, tensor size: 42.8 \textbf{KB})} $$\mathrm{tanh}(4(\var{1}-\var{2}))$$ \subsubsection{Setup and results overview}\begin{itemize}\item Reference: A/al. 2021 (A.5.3), \cite{Austin:2021}\item Domain: $\mathbb{R}$\item Tensor size: 42.8 \textbf{KB} ($74^{2}$ points)\item Bounds: $ \left(\begin{array}{cc} -1 & 1 \end{array}\right) \times \left(\begin{array}{cc} -1 & 1 \end{array}\right)$ \end{itemize} \begin{table}[H] \centering \begin{tabular}{llllll}
$\#$ & Alg. & Parameters & Dim. & CPU [s] & RMSE \\ 
\hline 
$\mathbf{\#8}$ & A/G/P-V 2025 (A1) & $1 \cdot 10^{-12},1$ & $4.8 \cdot 10^{02}$ & $0.028$ & $0.061$ \\ 
 & A/G/P-V 2025 (A2) & $1 \cdot 10^{-15},1$ & $\mathbf{2 \cdot 10^{02}}$ & $1.1$ & $0.00093$ \\ 
 & MDSPACK v1.1.0 & $1 \cdot 10^{-14},7$ & $4.8 \cdot 10^{02}$ & $\mathbf{0.025}$ & $0.06$ \\ 
 & P/P 2025 & $1,1,50,0.01,10,12,21$ & $6.8 \cdot 10^{02}$ & $3.5$ & $0.00076$ \\ 
 & C-R/B/G 2023 & $1 \cdot 10^{-09},20$ & $4.4 \cdot 10^{02}$ & $0.18$ & $\mathbf{8.2 \cdot 10^{-12}}$ \\ 
 & B/G 2025 & $1 \cdot 10^{-09},20,4$ & $1 \cdot 10^{03}$ & $4.9$ & $8 \cdot 10^{-06}$ \\ 
 & TensorFlow & $$ & $2.6 \cdot 10^{02}$ & $35$ & $0.0012$ \\ 
\hline 
\end{tabular} \caption{Function \#8: best model configuration and performances per methods.} \end{table}\begin{figure}[H] \centering  \includegraphics[width=\textwidth]{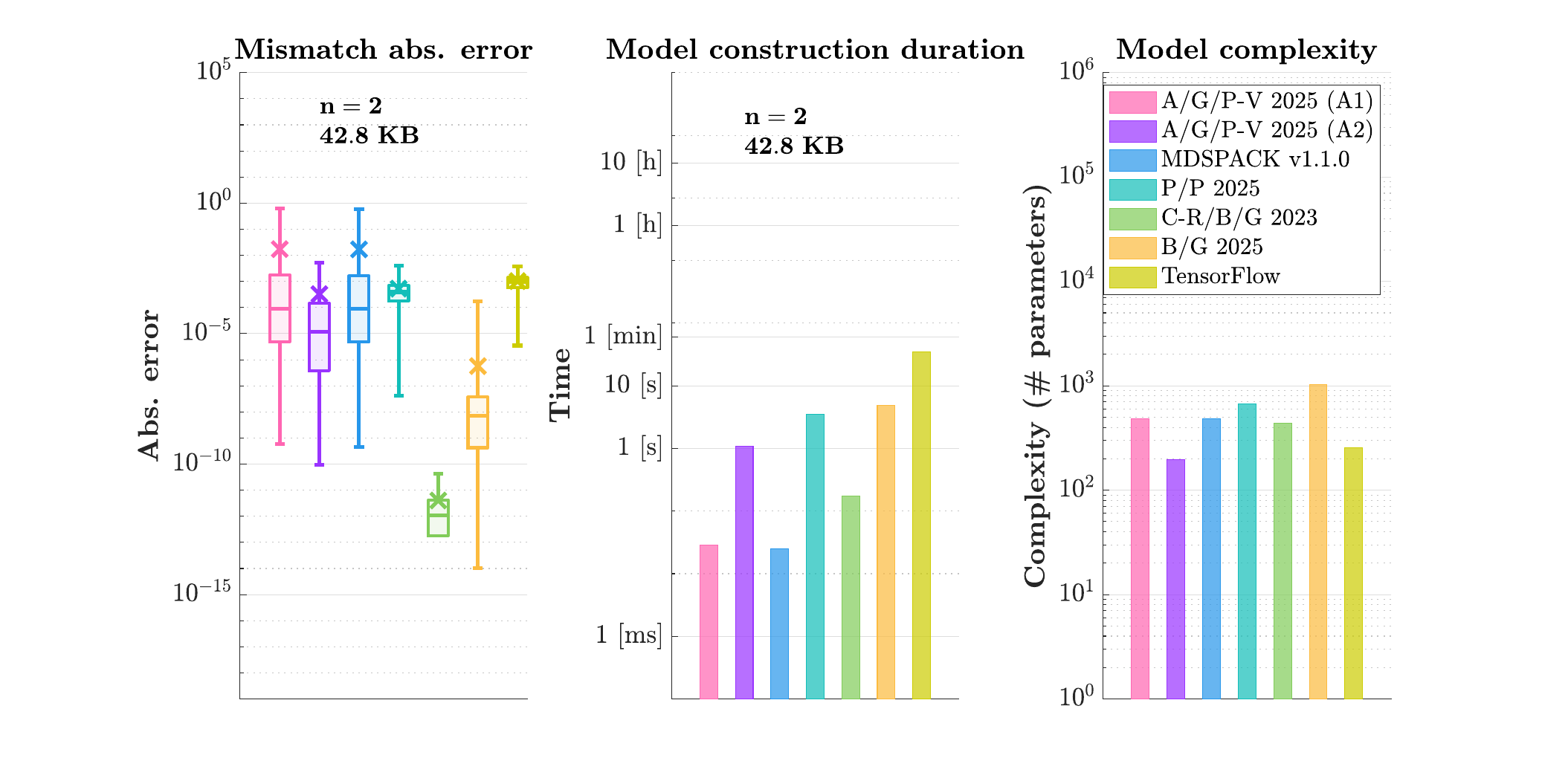} \caption{Function \#8: graphical view of the best model performances.} \end{figure}\begin{figure}[H] \centering  \includegraphics[width=\textwidth]{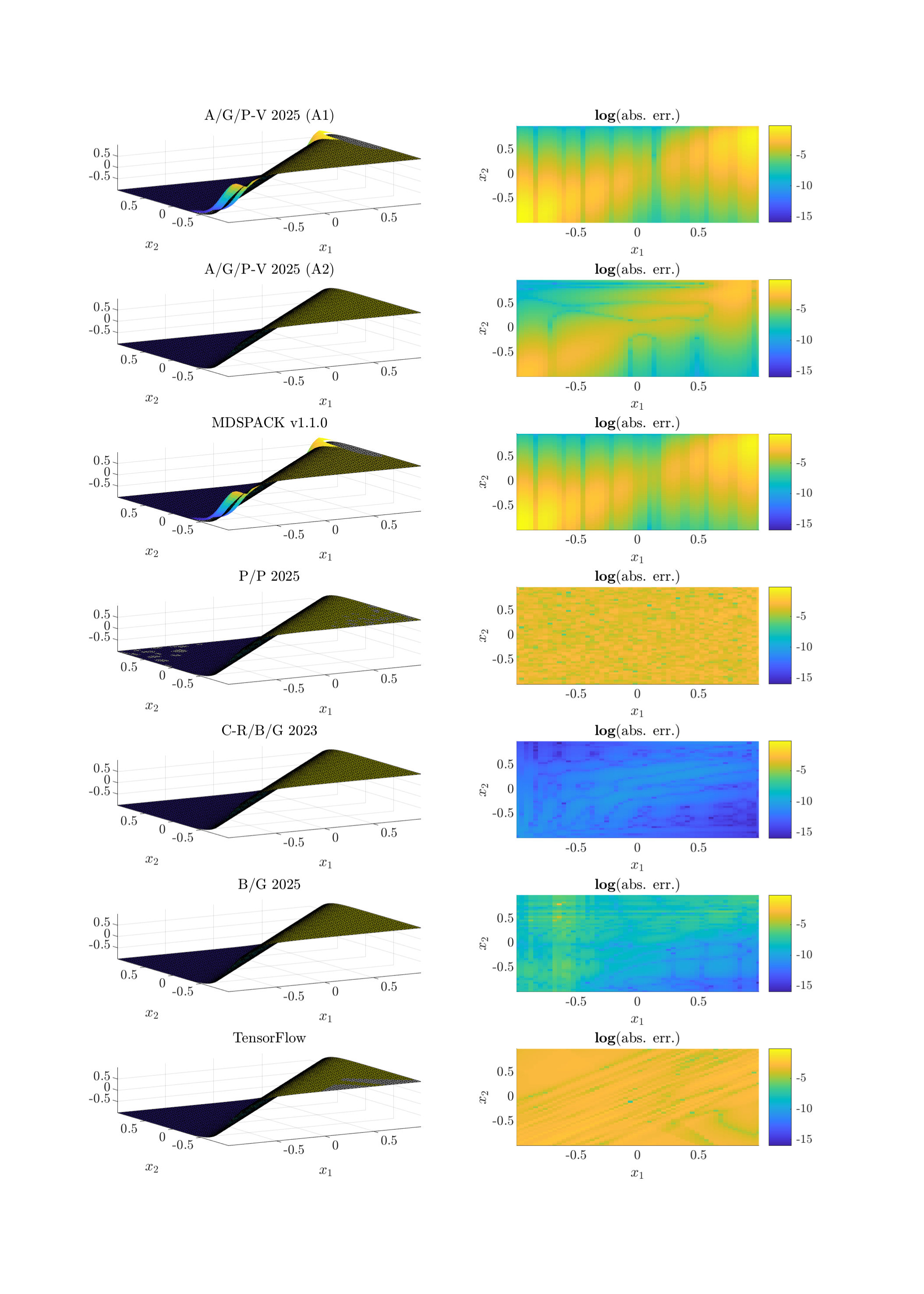} \caption{Function \#8: left side, evaluation of the original (mesh) vs. approximated (coloured surface) and right side, absolute errors (in log-scale).} \end{figure}\subsubsection{mLF detailed informations (M1)} \noindent \textbf{Right interpolation points}: $k_l=\left(\begin{array}{cc} 11 & 11 \end{array}\right)$, where $l=1,\cdots,\ord$.$$ \begin{array}{rcl}\lan{1} &\in& \IC^{11} \text{ , linearly spaced between bounds}\\\lan{2} &\in& \IC^{11} \text{ , linearly spaced between bounds}\\\end{array} $$\noindent \textbf{$\ord$-D Loewner matrix, barycentric weights and Lagrangian basis}:$$ \begin{array}{rcl}\IL & \in & \IC^{121 \times 121}\\\bc & \in & \IC^{121}\\\bw & \in & \IC^{121}\\\bc\odot \bw & \in & \IC^{121}\\\mathbf{Lag}(\var{1},\var{2}) & \in & \IC^{121}\\\end{array} $$

\newpage \subsection{Function \#9 (${\ord=2}$ variables, tensor size: 12.5 \textbf{KB})} $$\mathrm{exp}(\frac{-(\var{1}^2+\var{2}^2)}{1000})$$ \subsubsection{Setup and results overview}\begin{itemize}\item Reference: A/al. 2021 (A.5.4), \cite{Austin:2021}\item Domain: $\mathbb{R}$\item Tensor size: 12.5 \textbf{KB} ($40^{2}$ points)\item Bounds: $ \left(\begin{array}{cc} -1 & 1 \end{array}\right) \times \left(\begin{array}{cc} -1 & 1 \end{array}\right)$ \end{itemize} \begin{table}[H] \centering \begin{tabular}{llllll}
$\#$ & Alg. & Parameters & Dim. & CPU [s] & RMSE \\ 
\hline 
$\mathbf{\#9}$ & A/G/P-V 2025 (A1) & $0.5,2$ & $\mathbf{36}$ & $0.018$ & $1.4 \cdot 10^{-11}$ \\ 
 & A/G/P-V 2025 (A2) & $1 \cdot 10^{-15},1$ & $36$ & $0.15$ & $5.5 \cdot 10^{-12}$ \\ 
 & MDSPACK v1.1.0 & $0.0001,2$ & $36$ & $0.012$ & $1.5 \cdot 10^{-11}$ \\ 
 & P/P 2025 & $1,1,50,0.01,4,4,9$ & $1.1 \cdot 10^{02}$ & $0.22$ & $1.3 \cdot 10^{-08}$ \\ 
 & C-R/B/G 2023 & $0.001,20$ & $36$ & $0.015$ & $5.5 \cdot 10^{-12}$ \\ 
 & B/G 2025 & $1 \cdot 10^{-09},20,4$ & $36$ & $\mathbf{0.012}$ & $\mathbf{5.3 \cdot 10^{-12}}$ \\ 
 & TensorFlow & $$ & $2.6 \cdot 10^{02}$ & $14$ & $0.0011$ \\ 
\hline 
\end{tabular} \caption{Function \#9: best model configuration and performances per methods.} \end{table}\begin{figure}[H] \centering  \includegraphics[width=\textwidth]{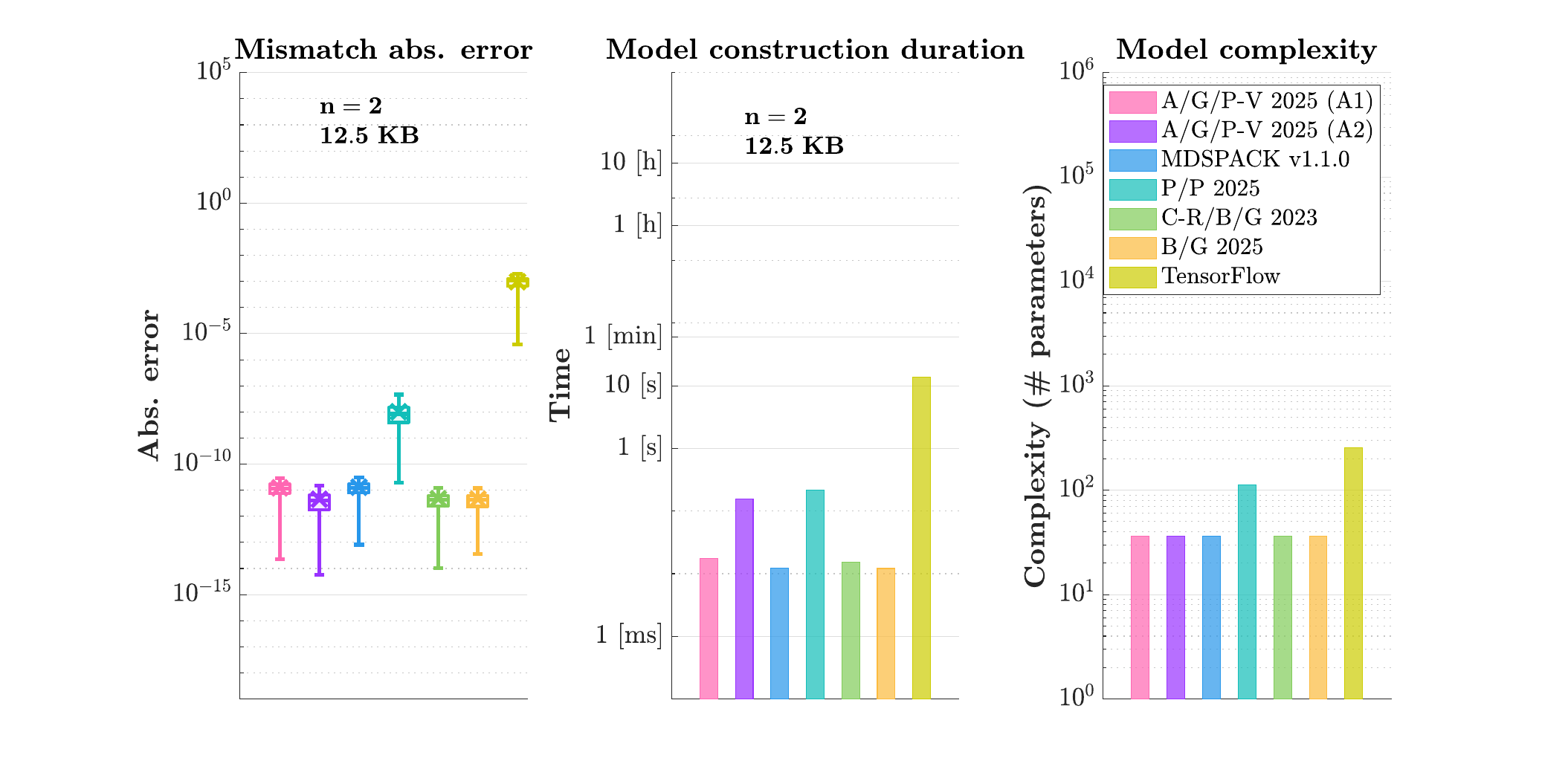} \caption{Function \#9: graphical view of the best model performances.} \end{figure}\begin{figure}[H] \centering  \includegraphics[width=\textwidth]{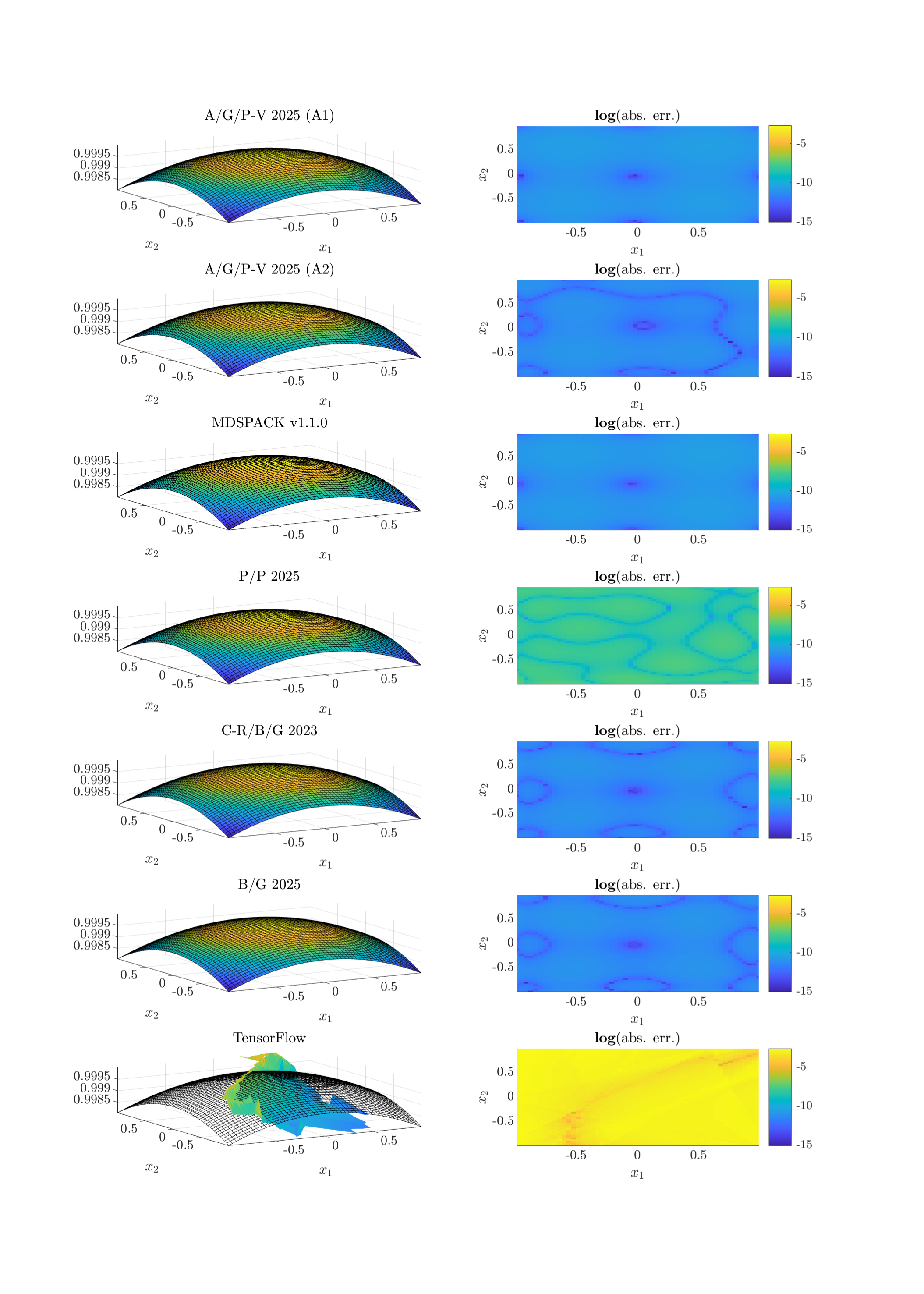} \caption{Function \#9: left side, evaluation of the original (mesh) vs. approximated (coloured surface) and right side, absolute errors (in log-scale).} \end{figure}\subsubsection{mLF detailed informations (M1)} \noindent \textbf{Right interpolation points} ($k_l=\left(\begin{array}{cc} 3 & 3 \end{array}\right)$, where $l=1,\cdots,\ord$):$$ \begin{array}{rcl}\lan{1} &=& \left(\begin{array}{ccc} -1 & -\frac{1}{19} & 1 \end{array}\right)\\\lan{2} &=& \left(\begin{array}{ccc} -1 & -\frac{1}{19} & 1 \end{array}\right)\\\end{array} $$\noindent \textbf{Lagrangian weights}: $$\left(\begin{array}{ccc} \bc & \bw & \bc\odot\bw\\ -0.5851 & 0.998 & -0.5839\\ 1.111 & 0.999 & 1.11\\ -0.5266 & 0.998 & -0.5255\\ 1.111 & 0.999 & 1.11\\ -2.11 & 1.0 & -2.11\\ 1.0 & 0.999 & 0.999\\ -0.5266 & 0.998 & -0.5255\\ 1.0 & 0.999 & 0.999\\ -0.4739 & 0.998 & -0.473 \end{array}\right)$$\noindent \textbf{Lagrangian form} (basis, numerator and denominator coefficients):$$\left(\begin{array}{ccc}\mathcal{B}_\textrm{lag}(\var{1},\var{2}) & \bN_\textrm{lag} &\bD_\textrm{lag}\end{array}\right) =$$ $$\left(\begin{array}{ccc} \left(\var{1}+1.0\right)\,\left(\var{2}+1.0\right) & -0.5839 & -0.5851\\ \left(\var{1}+1.0\right)\,\left(\var{2}+0.05263\right) & 1.11 & 1.111\\ \left(\var{1}+1.0\right)\,\left(\var{2}-1.0\right) & -0.5255 & -0.5266\\ \left(\var{2}+1.0\right)\,\left(\var{1}+0.05263\right) & 1.11 & 1.111\\ \left(\var{1}+0.05263\right)\,\left(\var{2}+0.05263\right) & -2.11 & -2.11\\ \left(\var{2}-1.0\right)\,\left(\var{1}+0.05263\right) & 0.999 & 1.0\\ \left(\var{1}-1.0\right)\,\left(\var{2}+1.0\right) & -0.5255 & -0.5266\\ \left(\var{1}-1.0\right)\,\left(\var{2}+0.05263\right) & 0.999 & 1.0\\ \left(\var{1}-1.0\right)\,\left(\var{2}-1.0\right) & -0.473 & -0.4739 \end{array}\right).$$\noindent The corresponding function is:$$\begin{array}{rcl}\bG_{\textrm{lag}}(\var{1},\var{2}) &=& \dfrac{\bn_{\textrm{lag}}(\var{1},\var{2})}{\bd_{\textrm{lag}}(\var{1},\var{2})}\\ && \\&=& \dfrac{\sum_{\textrm{row}} \bN_\textrm{lag} \odot\mathcal{B}^{-1}_\textrm{lag}(\var{1},\var{2})}{\sum_{\textrm{row}} \bD_\textrm{lag} \odot\mathcal{B}^{-1}_\textrm{lag}(\var{1},\var{2})}, \end{array}$$\noindent where,\\$\bn_{\textrm{lag}}(\var{1},\var{2}) = 2.498 \cdot 10^{-7}\,{\var{1}}^2\,{\var{2}}^2-4.247 \cdot 10^{-12}\,{\var{1}}^2\,\var{2}-0.0004998\,{\var{1}}^2-2.076 \cdot 10^{-12}\,\var{1}\,{\var{2}}^2-5.687 \cdot 10^{-17}\,\var{1}\,\var{2}+4.153 \cdot 10^{-9}\,\var{1}-0.0004998\,{\var{2}}^2+7.311 \cdot 10^{-9}\,\var{2}+1.0$ \\~~\\$\bd_{\textrm{lag}}(\var{1},\var{2}) = 2.502 \cdot 10^{-7}\,{\var{1}}^2\,{\var{2}}^2+3.066 \cdot 10^{-12}\,{\var{1}}^2\,\var{2}+0.0005002\,{\var{1}}^2+2.079 \cdot 10^{-12}\,\var{1}\,{\var{2}}^2-7.621 \cdot 10^{-18}\,\var{1}\,\var{2}+4.157 \cdot 10^{-9}\,\var{1}+0.0005002\,{\var{2}}^2+7.318 \cdot 10^{-9}\,\var{2}+1.0$ \\~~\\\noindent \textbf{Monomial form} (basis, numerator and denominator coefficients - entries $<10^{-12}$ removed):$$\left(\begin{array}{ccc}\mathcal{B}_\textrm{mon}(\var{1},\var{2}) & \bN_\textrm{mon} &\bD_\textrm{mon}\end{array}\right) =$$ $$\left(\begin{array}{ccc} {\var{1}}^2\,{\var{2}}^2 & -2.498 \cdot 10^{-7} & -2.502 \cdot 10^{-7}\\ {\var{1}}^2\,\var{2} & 4.247 \cdot 10^{-12} & -3.066 \cdot 10^{-12}\\ {\var{1}}^2 & 0.0004998 & -0.0005002\\ \var{1}\,{\var{2}}^2 & 2.076 \cdot 10^{-12} & -2.079 \cdot 10^{-12}\\ \var{1}\,\var{2} & 0 & 0\\ \var{1} & -4.153 \cdot 10^{-9} & -4.157 \cdot 10^{-9}\\ {\var{2}}^2 & 0.0004998 & -0.0005002\\ \var{2} & -7.311 \cdot 10^{-9} & -7.318 \cdot 10^{-9}\\ 1.0 & -1.0 & -1.0 \end{array}\right)$$\noindent The corresponding function is:$$\begin{array}{rcl}\bG_{\textrm{mon}}(\var{1},\var{2}) &=& \dfrac{\bn_{\textrm{mon}}(\var{1},\var{2})}{\bd_{\textrm{mon}}(\var{1},\var{2})}\\ && \\&=& \dfrac{\sum_{\textrm{row}} \bN_\textrm{mon} \odot \mathcal{B}_\textrm{mon}(\var{1},\var{2})}{\sum_{\textrm{row}} \bD_\textrm{mon} \odot\mathcal{B}_\textrm{mon}(\var{1},\var{2})},  \end{array}$$\noindent where,\\$\bn_{\textrm{mon}}(\var{1},\var{2}) = 2.498 \cdot 10^{-7}\,{\var{1}}^2\,{\var{2}}^2-4.247 \cdot 10^{-12}\,{\var{1}}^2\,\var{2}-0.0004998\,{\var{1}}^2-2.076 \cdot 10^{-12}\,\var{1}\,{\var{2}}^2+4.153 \cdot 10^{-9}\,\var{1}-0.0004998\,{\var{2}}^2+7.311 \cdot 10^{-9}\,\var{2}+1.0$ \\~~\\$\bd_{\textrm{mon}}(\var{1},\var{2}) = 2.502 \cdot 10^{-7}\,{\var{1}}^2\,{\var{2}}^2+3.066 \cdot 10^{-12}\,{\var{1}}^2\,\var{2}+0.0005002\,{\var{1}}^2+2.079 \cdot 10^{-12}\,\var{1}\,{\var{2}}^2+4.157 \cdot 10^{-9}\,\var{1}+0.0005002\,{\var{2}}^2+7.318 \cdot 10^{-9}\,\var{2}+1.0$ \\~~\\\noindent \textbf{KST equivalent decoupling pattern} (Barycentric weights $\bc^{\var{l}}$): $$\begin{array}{rclll}\var{2}&: & \left(\begin{array}{ccc} 1.111 & 1.111 & 1.111\\ -2.11 & -2.11 & -2.11\\ 1.0 & 1.0 & 1.0 \end{array}\right)& \textrm{vec}(.) & := \textbf{Bary}(\var{2}) \\\var{1}&: & \left(\begin{array}{c} -0.5266\\ 1.0\\ -0.4739 \end{array}\right)& \textrm{vec}(.) \otimes \bone_{k_{2}} & := \textbf{Bary}(\var{1}) \\\end{array}$$~\\ Then, with the above notations, one defines the following univariate vector functions:  $$ \left\{ \begin{array}{rcl}\bPhi_{1}(\var{1}) &:=& \textbf{Bary}(\var{1}) \odot \mathbf{Lag}(\var{1}) \\\bPhi_{2}(\var{2}) &:=& \textbf{Bary}(\var{2}) \odot \mathbf{Lag}(\var{2}) \\\end{array} \right. $$\noindent The corresponding function is:$$\begin{array}{rcl}\bG_{\textrm{kst}}(\var{1},\var{2}) &=& \dfrac{\bn_{\textrm{kst}}(\var{1},\var{2})}{\bd_{\textrm{kst}}(\var{1},\var{2})}\\ && \\ &=& \dfrac{\sum_{\text{rows}} \bw \odot \bPhi_{1}(\var{1}) \odot \cdots \odot\bPhi_{2}(\var{2})}{\sum_{\text{rows}} \bPhi_{1}(\var{1}) \odot \cdots \odot\bPhi_{2}(\var{2})} . \end{array}$$~\\ \noindent \textbf{KST-like univariate functions} (equivalent scaled univariate functions $\bphi_{1,\cdots,2}$): $$\left\{\begin{array}{rcrcl}z_{1} &=&\bphi_{1}(\var{1}) &=& \cfrac{\bn_{1}}{\bd_{1}} \\z_{2} &=&\bphi_{2}(\var{2}) &=& \cfrac{\bn_{2}}{\bd_{2}} \\\end{array} \right. .$$\noindent where, \\ \noindent $\bn_{1}=-0.0004993\,{\var{1}}^2+4.149 \cdot 10^{-9}\,\var{1}+0.999$ and \\ \noindent $\bd_{1}=0.0005002\,{\var{1}}^2+4.157 \cdot 10^{-9}\,\var{1}+1.0$, \\ \noindent $\bn_{2}=-0.0004993\,{\var{2}}^2+7.303 \cdot 10^{-9}\,\var{2}+0.999$ and \\ \noindent $\bd_{2}=0.0005002\,{\var{2}}^2+7.318 \cdot 10^{-9}\,\var{2}+1.0$, \\ \noindent \textbf{Connection with Neural Networks} (equivalent numerator $\bn_{\textrm{lag}}$ representation):\\ \begin{figure}[H]\begin{center} \scalebox{.7}{\begin{tikzpicture}[line width=0.4mm]\tikzstyle{place}=[circle, draw=black, minimum size = 8mm]\tikzstyle{placeInOut}=[circle, draw=orange, minimum size = 8mm]\node at (0,-2) [placeInOut] (first_1){$\var{1}$};\node at (0,-4) [placeInOut] (first_2){$\var{2}$};\node at (5,-2) [place] (secondL1_1){$\frac{1}{\var{1}-\lani{1}{1}}$};\node at (5,-4) [place] (secondL1_2){$\frac{1}{\var{1}-\lani{1}{2}}$};\node at (5,-6) [place] (secondL1_3){$\frac{1}{\var{1}-\lani{1}{3}}$};\node at (5,-8) [place] (secondL2_1){$\frac{1}{\var{2}-\lani{2}{1}}$};\node at (5,-10) [place] (secondL2_2){$\frac{1}{\var{2}-\lani{2}{2}}$};\node at (5,-12) [place] (secondL2_3){$\frac{1}{\var{2}-\lani{2}{3}}$};\node at (10,-2) [place] (third_1){$\prod$};\node at (10,-4) [place] (third_2){$\prod$};\node at (10,-6) [place] (third_3){$\prod$};\node at (10,-8) [place] (third_4){$\prod$};\node at (10,-10) [place] (third_5){$\prod$};\node at (10,-12) [place] (third_6){$\prod$};\node at (10,-14) [place] (third_7){$\prod$};\node at (10,-16) [place] (third_8){$\prod$};\node at (10,-18) [place] (third_9){$\prod$};\node at (15,-10) [placeInOut] (output){$\bSigma$};\draw[->] (first_1)--(secondL1_1) node[above,sloped,pos=0.75] { };\draw[->] (first_1)--(secondL1_2) node[above,sloped,pos=0.75] { };\draw[->] (first_1)--(secondL1_3) node[above,sloped,pos=0.75] { };\draw[->] (first_2)--(secondL2_1) node[above,sloped,pos=0.75] { };\draw[->] (first_2)--(secondL2_2) node[above,sloped,pos=0.75] { };\draw[->] (first_2)--(secondL2_3) node[above,sloped,pos=0.75] { };\draw[->] (secondL1_1)--(third_1) node[above,sloped,pos=0.25] {};\draw[->] (secondL1_1)--(third_2) node[above,sloped,pos=0.25] {};\draw[->] (secondL1_1)--(third_3) node[above,sloped,pos=0.25] {};\draw[->] (secondL1_2)--(third_4) node[above,sloped,pos=0.25] {};\draw[->] (secondL1_2)--(third_5) node[above,sloped,pos=0.25] {};\draw[->] (secondL1_2)--(third_6) node[above,sloped,pos=0.25] {};\draw[->] (secondL1_3)--(third_7) node[above,sloped,pos=0.25] {};\draw[->] (secondL1_3)--(third_8) node[above,sloped,pos=0.25] {};\draw[->] (secondL1_3)--(third_9) node[above,sloped,pos=0.25] {};\draw[->] (secondL2_1)--(third_1) node[above,sloped,pos=0.25] {};\draw[->] (secondL2_2)--(third_2) node[above,sloped,pos=0.25] {};\draw[->] (secondL2_3)--(third_3) node[above,sloped,pos=0.25] {};\draw[->] (secondL2_1)--(third_4) node[above,sloped,pos=0.25] {};\draw[->] (secondL2_2)--(third_5) node[above,sloped,pos=0.25] {};\draw[->] (secondL2_3)--(third_6) node[above,sloped,pos=0.25] {};\draw[->] (secondL2_1)--(third_7) node[above,sloped,pos=0.25] {};\draw[->] (secondL2_2)--(third_8) node[above,sloped,pos=0.25] {};\draw[->] (secondL2_3)--(third_9) node[above,sloped,pos=0.25] {};\draw[->] (third_1)--(output) node[above,sloped,pos=0.25] {-0.58392};\draw[->] (third_2)--(output) node[above,sloped,pos=0.25] {1.11};\draw[->] (third_3)--(output) node[above,sloped,pos=0.25] {-0.52553};\draw[->] (third_4)--(output) node[above,sloped,pos=0.25] {1.11};\draw[->] (third_5)--(output) node[above,sloped,pos=0.25] {-2.11};\draw[->] (third_6)--(output) node[above,sloped,pos=0.25] {0.999};\draw[->] (third_7)--(output) node[above,sloped,pos=0.25] {-0.52553};\draw[->] (third_8)--(output) node[above,sloped,pos=0.25] {0.999};\draw[->] (third_9)--(output) node[above,sloped,pos=0.25] {-0.47297};\end{tikzpicture}} \caption{Equivalent NN representation of the numerator $\bn_{\textrm{lag}}$.}\end{center}\end{figure}

\newpage \subsection{Function \#10 (${\ord=2}$ variables, tensor size: 52.5 \textbf{KB})} $$|\var{1}-\var{2}|^3$$ \subsubsection{Setup and results overview}\begin{itemize}\item Reference: A/al. 2021 (A.5.5), \cite{Austin:2021}\item Domain: $\mathbb{R}$\item Tensor size: 52.5 \textbf{KB} ($82^{2}$ points)\item Bounds: $ \left(\begin{array}{cc} -1 & 1 \end{array}\right) \times \left(\begin{array}{cc} -1 & 1 \end{array}\right)$ \end{itemize} \begin{table}[H] \centering \begin{tabular}{llllll}
$\#$ & Alg. & Parameters & Dim. & CPU [s] & RMSE \\ 
\hline 
$\mathbf{\#10}$ & A/G/P-V 2025 (A1) & $0.5,3$ & $\mathbf{36}$ & $\mathbf{0.032}$ & $0.076$ \\ 
 & A/G/P-V 2025 (A2) & $1 \cdot 10^{-15},1$ & $64$ & $0.38$ & $0.89$ \\ 
 & MDSPACK v1.1.0 & $0.01,1$ & $36$ & $0.048$ & $0.076$ \\ 
 & P/P 2025 & $1,1,50,0.01,10,6,21$ & $5.5 \cdot 10^{02}$ & $3.7$ & $3.3 \cdot 10^{-05}$ \\ 
 & C-R/B/G 2023 & $1 \cdot 10^{-06},20$ & $1 \cdot 10^{03}$ & $0.64$ & $\mathbf{3.1 \cdot 10^{-05}}$ \\ 
 & B/G 2025 & $1 \cdot 10^{-06},20,4$ & $1.2 \cdot 10^{03}$ & $4.5$ & $0.0012$ \\ 
 & TensorFlow & $$ & $2.6 \cdot 10^{02}$ & $40$ & $0.0069$ \\ 
\hline 
\end{tabular} \caption{Function \#10: best model configuration and performances per methods.} \end{table}\begin{figure}[H] \centering  \includegraphics[width=\textwidth]{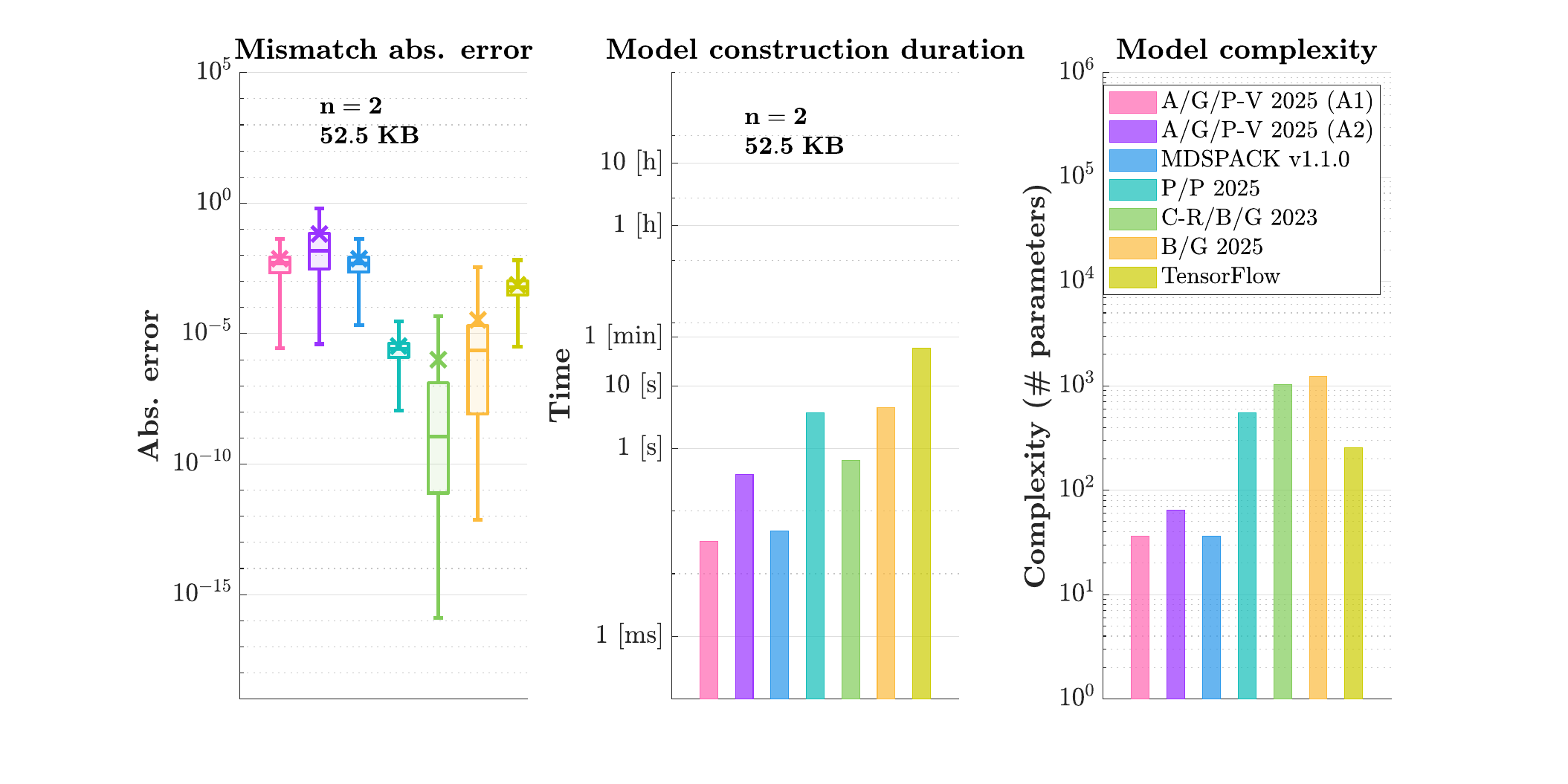} \caption{Function \#10: graphical view of the best model performances.} \end{figure}\begin{figure}[H] \centering  \includegraphics[width=\textwidth]{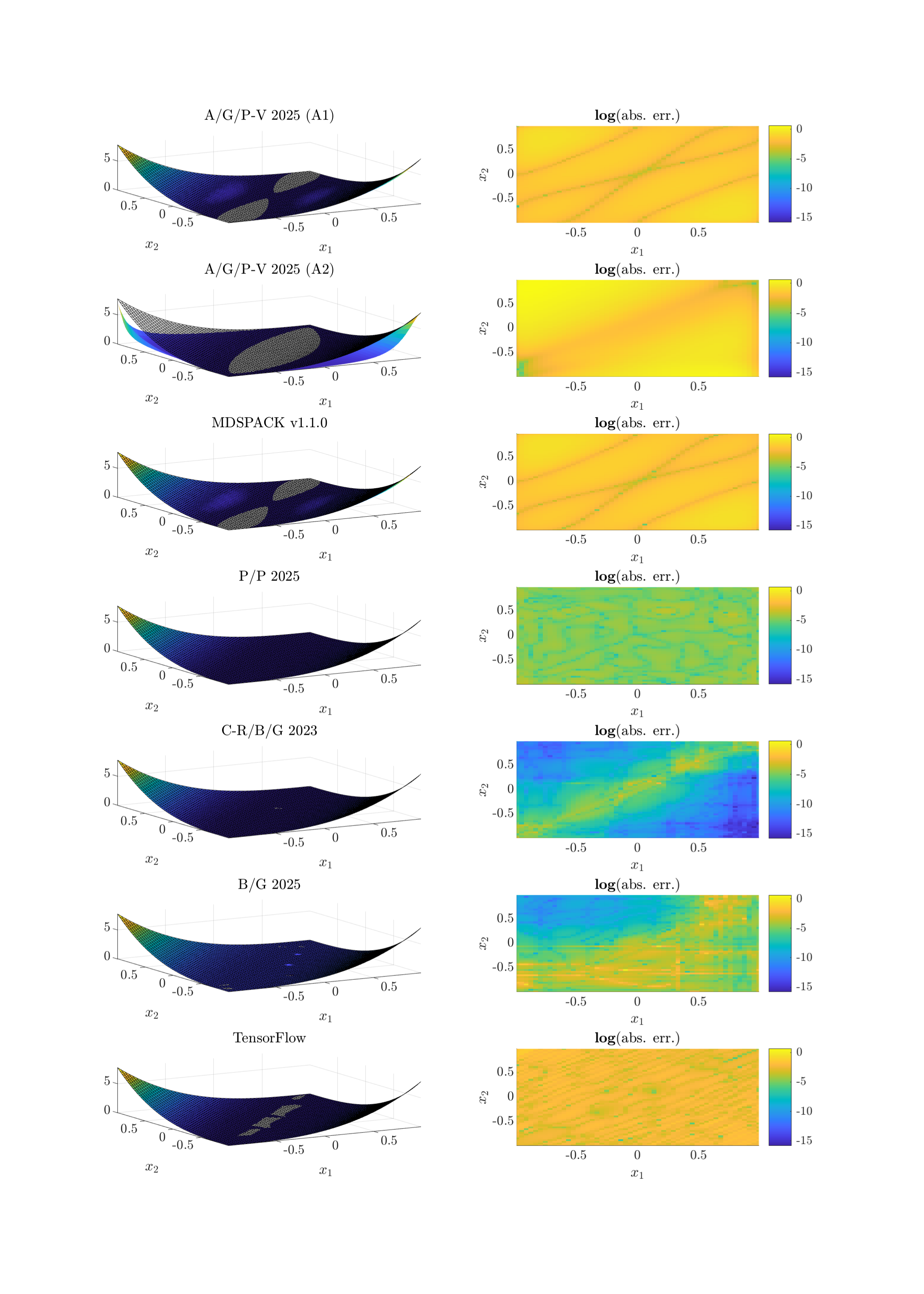} \caption{Function \#10: left side, evaluation of the original (mesh) vs. approximated (coloured surface) and right side, absolute errors (in log-scale).} \end{figure}\subsubsection{mLF detailed informations (M1)} \noindent \textbf{Right interpolation points} ($k_l=\left(\begin{array}{cc} 3 & 3 \end{array}\right)$, where $l=1,\cdots,\ord$):$$ \begin{array}{rcl}\lan{1} &=& \left(\begin{array}{ccc} -1 & 0 & 1 \end{array}\right)\\\lan{2} &=& \left(\begin{array}{ccc} -1 & 0 & 1 \end{array}\right)\\\end{array} $$\noindent \textbf{Lagrangian weights}: $$\left(\begin{array}{ccc} \bc & \bw & \bc\odot\bw\\ -0.8972 & 0 & 0\\ 0.9808 & 1.0 & 0.9808\\ -0.2752 & 8.0 & -2.202\\ 0.9884 & 1.0 & 0.9884\\ -2.982 & 0 & 0\\ 1.0 & 1.0 & 1.0\\ -0.2799 & 8.0 & -2.239\\ 1.009 & 1.0 & 1.009\\ -0.9374 & 0 & 0 \end{array}\right)$$\noindent \textbf{Lagrangian form} (basis, numerator and denominator coefficients):$$\left(\begin{array}{ccc}\mathcal{B}_\textrm{lag}(\var{1},\var{2}) & \bN_\textrm{lag} &\bD_\textrm{lag}\end{array}\right) =$$ $$\left(\begin{array}{ccc} \left(\var{1}+1.0\right)\,\left(\var{2}+1.0\right) & 0 & -0.8972\\ \var{2}\,\left(\var{1}+1.0\right) & 0.9808 & 0.9808\\ \left(\var{1}+1.0\right)\,\left(\var{2}-1.0\right) & -2.202 & -0.2752\\ \var{1}\,\left(\var{2}+1.0\right) & 0.9884 & 0.9884\\ \var{1}\,\var{2} & 0 & -2.982\\ \var{1}\,\left(\var{2}-1.0\right) & 1.0 & 1.0\\ \left(\var{1}-1.0\right)\,\left(\var{2}+1.0\right) & -2.239 & -0.2799\\ \var{2}\,\left(\var{1}-1.0\right) & 1.009 & 1.009\\ \left(\var{1}-1.0\right)\,\left(\var{2}-1.0\right) & 0 & -0.9374 \end{array}\right).$$\noindent The corresponding function is:$$\begin{array}{rcl}\bG_{\textrm{lag}}(\var{1},\var{2}) &=& \dfrac{\bn_{\textrm{lag}}(\var{1},\var{2})}{\bd_{\textrm{lag}}(\var{1},\var{2})}\\ && \\&=& \dfrac{\sum_{\textrm{row}} \bN_\textrm{lag} \odot\mathcal{B}^{-1}_\textrm{lag}(\var{1},\var{2})}{\sum_{\textrm{row}} \bD_\textrm{lag} \odot\mathcal{B}^{-1}_\textrm{lag}(\var{1},\var{2})}, \end{array}$$\noindent where,\\$\bn_{\textrm{lag}}(\var{1},\var{2}) = 0.1551\,{\var{1}}^2\,{\var{2}}^2-0.0165\,{\var{1}}^2\,\var{2}+0.6674\,{\var{1}}^2+0.003143\,\var{1}\,{\var{2}}^2-1.489\,\var{1}\,\var{2}+0.00947\,\var{1}+0.6669\,{\var{2}}^2+0.003889\,\var{2}$ \\~~\\$\bd_{\textrm{lag}}(\var{1},\var{2}) = 0.4672\,{\var{1}}^2\,{\var{2}}^2+0.008018\,{\var{1}}^2\,\var{2}-0.3326\,{\var{1}}^2+0.005589\,\var{1}\,{\var{2}}^2+0.4291\,\var{1}\,\var{2}+0.00947\,\var{1}-0.3331\,{\var{2}}^2+0.003889\,\var{2}+1.0$ \\~~\\\noindent \textbf{Monomial form} (basis, numerator and denominator coefficients - entries $<10^{-12}$ removed):$$\left(\begin{array}{ccc}\mathcal{B}_\textrm{mon}(\var{1},\var{2}) & \bN_\textrm{mon} &\bD_\textrm{mon}\end{array}\right) =$$ $$\left(\begin{array}{ccc} {\var{1}}^2\,{\var{2}}^2 & -0.1551 & -0.4672\\ {\var{1}}^2\,\var{2} & 0.0165 & -0.008018\\ {\var{1}}^2 & -0.6674 & 0.3326\\ \var{1}\,{\var{2}}^2 & -0.003143 & -0.005589\\ \var{1}\,\var{2} & 1.489 & -0.4291\\ \var{1} & -0.00947 & -0.00947\\ {\var{2}}^2 & -0.6669 & 0.3331\\ \var{2} & -0.003889 & -0.003889\\ 1.0 & 0 & -1.0 \end{array}\right)$$\noindent The corresponding function is:$$\begin{array}{rcl}\bG_{\textrm{mon}}(\var{1},\var{2}) &=& \dfrac{\bn_{\textrm{mon}}(\var{1},\var{2})}{\bd_{\textrm{mon}}(\var{1},\var{2})}\\ && \\&=& \dfrac{\sum_{\textrm{row}} \bN_\textrm{mon} \odot \mathcal{B}_\textrm{mon}(\var{1},\var{2})}{\sum_{\textrm{row}} \bD_\textrm{mon} \odot\mathcal{B}_\textrm{mon}(\var{1},\var{2})},  \end{array}$$\noindent where,\\$\bn_{\textrm{mon}}(\var{1},\var{2}) = 0.1551\,{\var{1}}^2\,{\var{2}}^2-0.0165\,{\var{1}}^2\,\var{2}+0.6674\,{\var{1}}^2+0.003143\,\var{1}\,{\var{2}}^2-1.489\,\var{1}\,\var{2}+0.00947\,\var{1}+0.6669\,{\var{2}}^2+0.003889\,\var{2}$ \\~~\\$\bd_{\textrm{mon}}(\var{1},\var{2}) = 0.4672\,{\var{1}}^2\,{\var{2}}^2+0.008018\,{\var{1}}^2\,\var{2}-0.3326\,{\var{1}}^2+0.005589\,\var{1}\,{\var{2}}^2+0.4291\,\var{1}\,\var{2}+0.00947\,\var{1}-0.3331\,{\var{2}}^2+0.003889\,\var{2}+1.0$ \\~~\\\noindent \textbf{KST equivalent decoupling pattern} (Barycentric weights $\bc^{\var{l}}$): $$\begin{array}{rclll}\var{2}&: & \left(\begin{array}{ccc} 3.26 & 0.9884 & 0.2986\\ -3.564 & -2.982 & -1.076\\ 1.0 & 1.0 & 1.0 \end{array}\right)& \textrm{vec}(.) & := \textbf{Bary}(\var{2}) \\\var{1}&: & \left(\begin{array}{c} -0.2752\\ 1.0\\ -0.9374 \end{array}\right)& \textrm{vec}(.) \otimes \bone_{k_{2}} & := \textbf{Bary}(\var{1}) \\\end{array}$$~\\ Then, with the above notations, one defines the following univariate vector functions:  $$ \left\{ \begin{array}{rcl}\bPhi_{1}(\var{1}) &:=& \textbf{Bary}(\var{1}) \odot \mathbf{Lag}(\var{1}) \\\bPhi_{2}(\var{2}) &:=& \textbf{Bary}(\var{2}) \odot \mathbf{Lag}(\var{2}) \\\end{array} \right. $$\noindent The corresponding function is:$$\begin{array}{rcl}\bG_{\textrm{kst}}(\var{1},\var{2}) &=& \dfrac{\bn_{\textrm{kst}}(\var{1},\var{2})}{\bd_{\textrm{kst}}(\var{1},\var{2})}\\ && \\ &=& \dfrac{\sum_{\text{rows}} \bw \odot \bPhi_{1}(\var{1}) \odot \cdots \odot\bPhi_{2}(\var{2})}{\sum_{\text{rows}} \bPhi_{1}(\var{1}) \odot \cdots \odot\bPhi_{2}(\var{2})} . \end{array}$$~\\ \noindent \textbf{KST-like univariate functions} (equivalent scaled univariate functions $\bphi_{1,\cdots,2}$): $$\left\{\begin{array}{rcrcl}z_{1} &=&\bphi_{1}(\var{1}) &=& \cfrac{\bn_{1}}{\bd_{1}} \\z_{2} &=&\bphi_{2}(\var{2}) &=& \cfrac{\bn_{2}}{\bd_{2}} \\\end{array} \right. .$$\noindent where, \\ \noindent $\bn_{1}=1.202\,{\var{1}}^2-2.202\,\var{1}+1.0$ and \\ \noindent $\bd_{1}=0.2126\,{\var{1}}^2+0.6622\,\var{1}+1.0$, \\ \noindent $\bn_{2}=1.245\,{\var{2}}^2+2.245\,\var{2}+1.0$ and \\ \noindent $\bd_{2}=0.1953\,{\var{2}}^2-0.6342\,\var{2}+1.0$, \\ \noindent \textbf{Connection with Neural Networks} (equivalent numerator $\bn_{\textrm{lag}}$ representation):\\ \begin{figure}[H]\begin{center} \scalebox{.7}{\begin{tikzpicture}[line width=0.4mm]\tikzstyle{place}=[circle, draw=black, minimum size = 8mm]\tikzstyle{placeInOut}=[circle, draw=orange, minimum size = 8mm]\node at (0,-2) [placeInOut] (first_1){$\var{1}$};\node at (0,-4) [placeInOut] (first_2){$\var{2}$};\node at (5,-2) [place] (secondL1_1){$\frac{1}{\var{1}-\lani{1}{1}}$};\node at (5,-4) [place] (secondL1_2){$\frac{1}{\var{1}-\lani{1}{2}}$};\node at (5,-6) [place] (secondL1_3){$\frac{1}{\var{1}-\lani{1}{3}}$};\node at (5,-8) [place] (secondL2_1){$\frac{1}{\var{2}-\lani{2}{1}}$};\node at (5,-10) [place] (secondL2_2){$\frac{1}{\var{2}-\lani{2}{2}}$};\node at (5,-12) [place] (secondL2_3){$\frac{1}{\var{2}-\lani{2}{3}}$};\node at (10,-2) [place] (third_1){$\prod$};\node at (10,-4) [place] (third_2){$\prod$};\node at (10,-6) [place] (third_3){$\prod$};\node at (10,-8) [place] (third_4){$\prod$};\node at (10,-10) [place] (third_5){$\prod$};\node at (10,-12) [place] (third_6){$\prod$};\node at (10,-14) [place] (third_7){$\prod$};\node at (10,-16) [place] (third_8){$\prod$};\node at (10,-18) [place] (third_9){$\prod$};\node at (15,-10) [placeInOut] (output){$\bSigma$};\draw[->] (first_1)--(secondL1_1) node[above,sloped,pos=0.75] { };\draw[->] (first_1)--(secondL1_2) node[above,sloped,pos=0.75] { };\draw[->] (first_1)--(secondL1_3) node[above,sloped,pos=0.75] { };\draw[->] (first_2)--(secondL2_1) node[above,sloped,pos=0.75] { };\draw[->] (first_2)--(secondL2_2) node[above,sloped,pos=0.75] { };\draw[->] (first_2)--(secondL2_3) node[above,sloped,pos=0.75] { };\draw[->] (secondL1_1)--(third_1) node[above,sloped,pos=0.25] {};\draw[->] (secondL1_1)--(third_2) node[above,sloped,pos=0.25] {};\draw[->] (secondL1_1)--(third_3) node[above,sloped,pos=0.25] {};\draw[->] (secondL1_2)--(third_4) node[above,sloped,pos=0.25] {};\draw[->] (secondL1_2)--(third_5) node[above,sloped,pos=0.25] {};\draw[->] (secondL1_2)--(third_6) node[above,sloped,pos=0.25] {};\draw[->] (secondL1_3)--(third_7) node[above,sloped,pos=0.25] {};\draw[->] (secondL1_3)--(third_8) node[above,sloped,pos=0.25] {};\draw[->] (secondL1_3)--(third_9) node[above,sloped,pos=0.25] {};\draw[->] (secondL2_1)--(third_1) node[above,sloped,pos=0.25] {};\draw[->] (secondL2_2)--(third_2) node[above,sloped,pos=0.25] {};\draw[->] (secondL2_3)--(third_3) node[above,sloped,pos=0.25] {};\draw[->] (secondL2_1)--(third_4) node[above,sloped,pos=0.25] {};\draw[->] (secondL2_2)--(third_5) node[above,sloped,pos=0.25] {};\draw[->] (secondL2_3)--(third_6) node[above,sloped,pos=0.25] {};\draw[->] (secondL2_1)--(third_7) node[above,sloped,pos=0.25] {};\draw[->] (secondL2_2)--(third_8) node[above,sloped,pos=0.25] {};\draw[->] (secondL2_3)--(third_9) node[above,sloped,pos=0.25] {};\draw[->] (third_1)--(output) node[above,sloped,pos=0.25] {0};\draw[->] (third_2)--(output) node[above,sloped,pos=0.25] {0.98079};\draw[->] (third_3)--(output) node[above,sloped,pos=0.25] {-2.2015};\draw[->] (third_4)--(output) node[above,sloped,pos=0.25] {0.98841};\draw[->] (third_5)--(output) node[above,sloped,pos=0.25] {0};\draw[->] (third_6)--(output) node[above,sloped,pos=0.25] {1};\draw[->] (third_7)--(output) node[above,sloped,pos=0.25] {-2.2391};\draw[->] (third_8)--(output) node[above,sloped,pos=0.25] {1.009};\draw[->] (third_9)--(output) node[above,sloped,pos=0.25] {0};\end{tikzpicture}} \caption{Equivalent NN representation of the numerator $\bn_{\textrm{lag}}$.}\end{center}\end{figure}

\newpage \subsection{Function \#11 (${\ord=2}$ variables, tensor size: 12.5 \textbf{KB})} $$\frac{\var{1}+\var{2}^3}{\var{1}\var{2}^2+2}$$ \subsubsection{Setup and results overview}\begin{itemize}\item Reference: A/al. 2021 (A.5.6), \cite{Austin:2021}\item Domain: $\mathbb{R}$\item Tensor size: 12.5 \textbf{KB} ($40^{2}$ points)\item Bounds: $ \left(\begin{array}{cc} \frac{1}{10000000000} & 1 \end{array}\right) \times \left(\begin{array}{cc} \frac{1}{10000000000} & 1 \end{array}\right)$ \end{itemize} \begin{table}[H] \centering \begin{tabular}{llllll}
$\#$ & Alg. & Parameters & Dim. & CPU [s] & RMSE \\ 
\hline 
$\mathbf{\#11}$ & A/G/P-V 2025 (A1) & $0.01,2$ & $\mathbf{32}$ & $\mathbf{0.0093}$ & $3 \cdot 10^{-15}$ \\ 
 & A/G/P-V 2025 (A2) & $1 \cdot 10^{-15},2$ & $32$ & $0.085$ & $1.8 \cdot 10^{-15}$ \\ 
 & MDSPACK v1.1.0 & $0.01,1$ & $32$ & $0.014$ & $3.4 \cdot 10^{-15}$ \\ 
 & P/P 2025 & $1,0.95,50,0.01,6,12,13$ & $3.2 \cdot 10^{02}$ & $1.2$ & $1.4 \cdot 10^{-05}$ \\ 
 & C-R/B/G 2023 & $0.001,20$ & $80$ & $0.012$ & $3.6 \cdot 10^{-15}$ \\ 
 & B/G 2025 & $1 \cdot 10^{-09},20,4$ & $80$ & $0.023$ & $\mathbf{1.5 \cdot 10^{-15}}$ \\ 
 & TensorFlow & $$ & $2.6 \cdot 10^{02}$ & $14$ & $0.0066$ \\ 
\hline 
\end{tabular} \caption{Function \#11: best model configuration and performances per methods.} \end{table}\begin{figure}[H] \centering  \includegraphics[width=\textwidth]{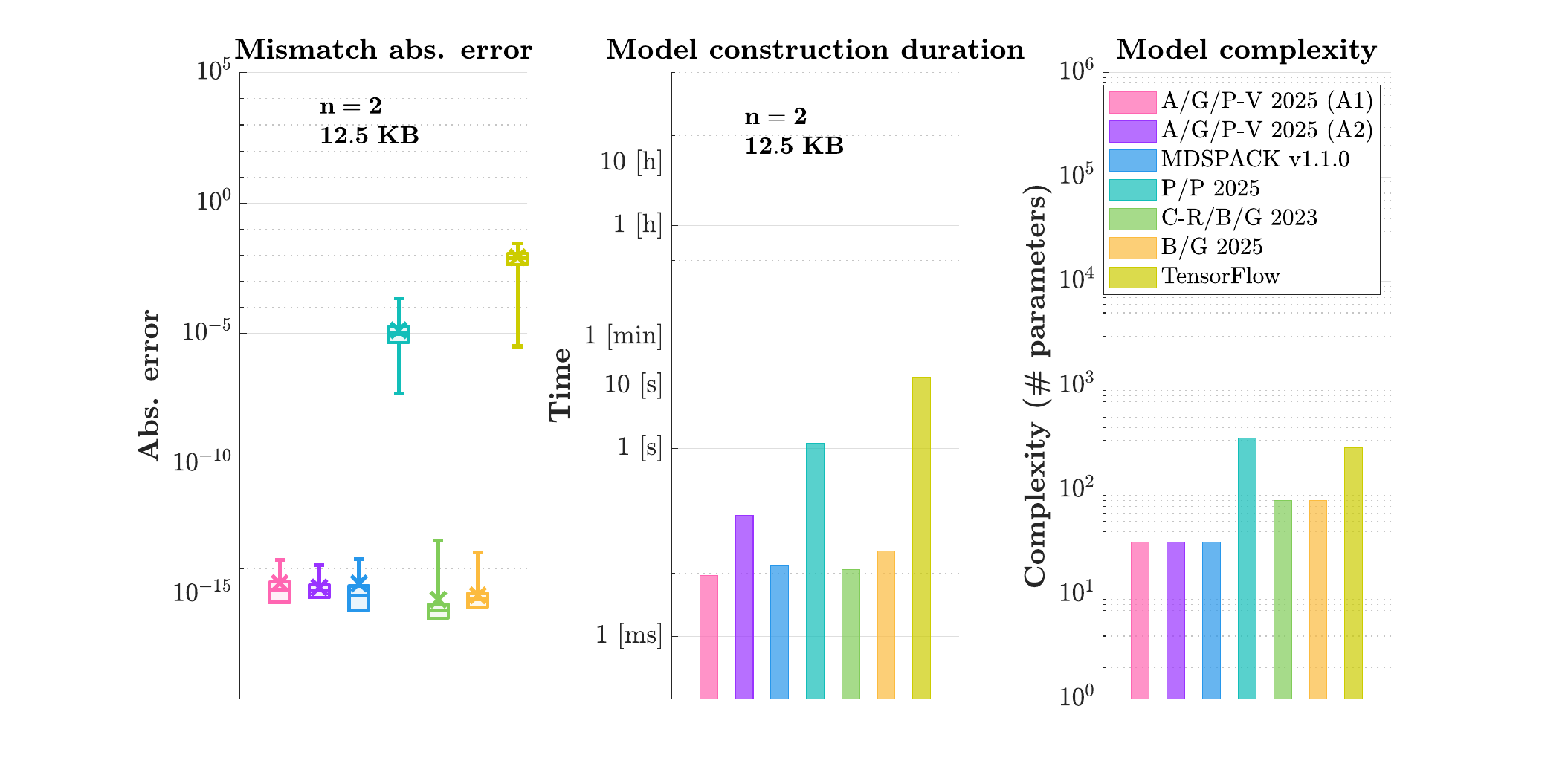} \caption{Function \#11: graphical view of the best model performances.} \end{figure}\begin{figure}[H] \centering  \includegraphics[width=\textwidth]{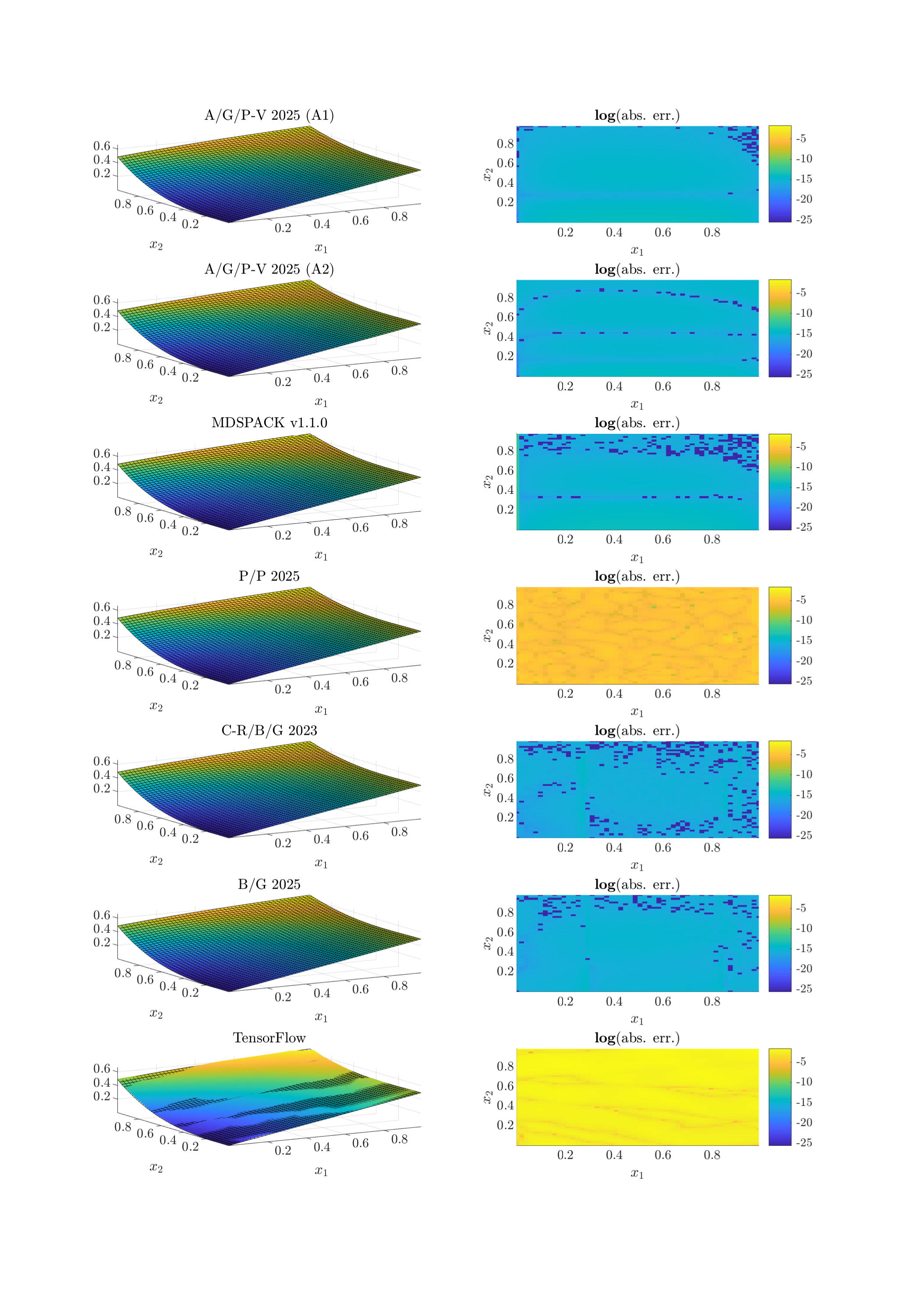} \caption{Function \#11: left side, evaluation of the original (mesh) vs. approximated (coloured surface) and right side, absolute errors (in log-scale).} \end{figure}\subsubsection{mLF detailed informations (M1)} \noindent \textbf{Right interpolation points} ($k_l=\left(\begin{array}{cc} 2 & 4 \end{array}\right)$, where $l=1,\cdots,\ord$):$$ \begin{array}{rcl}\lan{1} &=& \left(\begin{array}{cc} \frac{1}{10000000000} & 1 \end{array}\right)\\\lan{2} &=& \left(\begin{array}{cccc} \frac{1}{10000000000} & \frac{5688757425279507}{18014398509481984} & \frac{5688757424378787}{9007199254740992} & 1 \end{array}\right)\\\end{array} $$\noindent \textbf{Lagrangian weights}: $$\left(\begin{array}{ccc} \bc & \bw & \bc\odot\bw\\ 0.8426 & 5.0 \cdot 10^{-11} & 4.213 \cdot 10^{-11}\\ -2.463 & 0.01575 & -0.03878\\ 2.287 & 0.126 & 0.2881\\ -0.6667 & 0.5 & -0.3333\\ -0.8426 & 0.5 & -0.4213\\ 2.586 & 0.4913 & 1.27\\ -2.743 & 0.5219 & -1.432\\ 1.0 & 0.6667 & 0.6667 \end{array}\right)$$\noindent \textbf{Lagrangian form} (basis, numerator and denominator coefficients):$$\left(\begin{array}{ccc}\mathcal{B}_\textrm{lag}(\var{1},\var{2}) & \bN_\textrm{lag} &\bD_\textrm{lag}\end{array}\right) =$$ $$\left(\begin{array}{ccc} \left(\var{1}-1.0 \cdot 10^{-10}\right)\,\left(\var{2}-1.0 \cdot 10^{-10}\right) & 4.213 \cdot 10^{-11} & 0.8426\\ \left(\var{2}-0.3158\right)\,\left(\var{1}-1.0 \cdot 10^{-10}\right) & -0.03878 & -2.463\\ \left(\var{2}-0.6316\right)\,\left(\var{1}-1.0 \cdot 10^{-10}\right) & 0.2881 & 2.287\\ \left(\var{2}-1.0\right)\,\left(\var{1}-1.0 \cdot 10^{-10}\right) & -0.3333 & -0.6667\\ \left(\var{1}-1.0\right)\,\left(\var{2}-1.0 \cdot 10^{-10}\right) & -0.4213 & -0.8426\\ \left(\var{1}-1.0\right)\,\left(\var{2}-0.3158\right) & 1.27 & 2.586\\ \left(\var{1}-1.0\right)\,\left(\var{2}-0.6316\right) & -1.432 & -2.743\\ \left(\var{1}-1.0\right)\,\left(\var{2}-1.0\right) & 0.6667 & 1.0 \end{array}\right).$$\noindent The corresponding function is:$$\begin{array}{rcl}\bG_{\textrm{lag}}(\var{1},\var{2}) &=& \dfrac{\bn_{\textrm{lag}}(\var{1},\var{2})}{\bd_{\textrm{lag}}(\var{1},\var{2})}\\ && \\&=& \dfrac{\sum_{\textrm{row}} \bN_\textrm{lag} \odot\mathcal{B}^{-1}_\textrm{lag}(\var{1},\var{2})}{\sum_{\textrm{row}} \bD_\textrm{lag} \odot\mathcal{B}^{-1}_\textrm{lag}(\var{1},\var{2})}, \end{array}$$\noindent where,\\$\bn_{\textrm{lag}}(\var{1},\var{2}) = 0.5\,\var{1}+2.846 \cdot 10^{-16}\,\var{2}+3.804 \cdot 10^{-13}\,\var{1}\,\var{2}-5.086 \cdot 10^{-13}\,\var{1}\,{\var{2}}^2+2.161 \cdot 10^{-13}\,\var{1}\,{\var{2}}^3-5.936 \cdot 10^{-15}\,{\var{2}}^2+0.5\,{\var{2}}^3+5.737 \cdot 10^{-24}$ \\~~\\$\bd_{\textrm{lag}}(\var{1},\var{2}) = 1.306 \cdot 10^{-13}\,\var{2}-1.153 \cdot 10^{-13}\,\var{1}+6.298 \cdot 10^{-13}\,\var{1}\,\var{2}+0.5\,\var{1}\,{\var{2}}^2+3.984 \cdot 10^{-13}\,\var{1}\,{\var{2}}^3-1.206 \cdot 10^{-13}\,{\var{2}}^2+4.096 \cdot 10^{-14}\,{\var{2}}^3+1.0$ \\~~\\\noindent \textbf{Monomial form} (basis, numerator and denominator coefficients - entries $<10^{-12}$ removed):$$\left(\begin{array}{ccc}\mathcal{B}_\textrm{mon}(\var{1},\var{2}) & \bN_\textrm{mon} &\bD_\textrm{mon}\end{array}\right) =$$ $$\left(\begin{array}{ccc} \var{1}\,{\var{2}}^3 & 0 & 0\\ \var{1}\,{\var{2}}^2 & 0 & 0.5\\ \var{1}\,\var{2} & 0 & 0\\ \var{1} & 0.5 & 0\\ {\var{2}}^3 & 0.5 & 0\\ {\var{2}}^2 & 0 & 0\\ \var{2} & 0 & 0\\ 1.0 & 0 & 1.0 \end{array}\right)$$\noindent The corresponding function is:$$\begin{array}{rcl}\bG_{\textrm{mon}}(\var{1},\var{2}) &=& \dfrac{\bn_{\textrm{mon}}(\var{1},\var{2})}{\bd_{\textrm{mon}}(\var{1},\var{2})}\\ && \\&=& \dfrac{\sum_{\textrm{row}} \bN_\textrm{mon} \odot \mathcal{B}_\textrm{mon}(\var{1},\var{2})}{\sum_{\textrm{row}} \bD_\textrm{mon} \odot\mathcal{B}_\textrm{mon}(\var{1},\var{2})},  \end{array}$$\noindent where,\\$\bn_{\textrm{mon}}(\var{1},\var{2}) = 0.5\,{\var{2}}^3+0.5\,\var{1}$ \\~~\\$\bd_{\textrm{mon}}(\var{1},\var{2}) = 0.5\,\var{1}\,{\var{2}}^2+1.0$ \\~~\\\noindent \textbf{KST equivalent decoupling pattern} (Barycentric weights $\bc^{\var{l}}$): $$\begin{array}{rclll}\var{2}&: & \left(\begin{array}{cc} -1.264 & -0.8426\\ 3.694 & 2.586\\ -3.431 & -2.743\\ 1.0 & 1.0 \end{array}\right)& \textrm{vec}(.) & := \textbf{Bary}(\var{2}) \\\var{1}&: & \left(\begin{array}{c} -0.6667\\ 1.0 \end{array}\right)& \textrm{vec}(.) \otimes \bone_{k_{2}} & := \textbf{Bary}(\var{1}) \\\end{array}$$~\\ Then, with the above notations, one defines the following univariate vector functions:  $$ \left\{ \begin{array}{rcl}\bPhi_{1}(\var{1}) &:=& \textbf{Bary}(\var{1}) \odot \mathbf{Lag}(\var{1}) \\\bPhi_{2}(\var{2}) &:=& \textbf{Bary}(\var{2}) \odot \mathbf{Lag}(\var{2}) \\\end{array} \right. $$\noindent The corresponding function is:$$\begin{array}{rcl}\bG_{\textrm{kst}}(\var{1},\var{2}) &=& \dfrac{\bn_{\textrm{kst}}(\var{1},\var{2})}{\bd_{\textrm{kst}}(\var{1},\var{2})}\\ && \\ &=& \dfrac{\sum_{\text{rows}} \bw \odot \bPhi_{1}(\var{1}) \odot \cdots \odot\bPhi_{2}(\var{2})}{\sum_{\text{rows}} \bPhi_{1}(\var{1}) \odot \cdots \odot\bPhi_{2}(\var{2})} . \end{array}$$~\\ \noindent \textbf{KST-like univariate functions} (equivalent scaled univariate functions $\bphi_{1,\cdots,2}$): $$\left\{\begin{array}{rcrcl}z_{1} &=&\bphi_{1}(\var{1}) &=& \cfrac{1.48 \cdot 10^{-16}\,\left(2.476 \cdot 10^{37}\,\var{1}+2.476 \cdot 10^{37}\right)}{3.665 \cdot 10^{21}\,\var{1}+7.33 \cdot 10^{21}}\\z_{2} &=&\bphi_{2}(\var{2}) &=& \cfrac{\bn_{2}}{\bd_{2}} \\\end{array} \right. .$$\noindent where, \\ \noindent $\bn_{2}=0.5\,{\var{2}}^3-6.125 \cdot 10^{-15}\,{\var{2}}^2+3.284 \cdot 10^{-16}\,\var{2}+5.0 \cdot 10^{-11}$ and \\ \noindent $\bd_{2}=4.052 \cdot 10^{-14}\,{\var{2}}^3+4.988 \cdot 10^{-11}\,{\var{2}}^2+1.302 \cdot 10^{-13}\,\var{2}+1.0$, \\ \noindent \textbf{Connection with Neural Networks} (equivalent numerator $\bn_{\textrm{lag}}$ representation):\\ \begin{figure}[H]\begin{center} \scalebox{.7}{\begin{tikzpicture}[line width=0.4mm]\tikzstyle{place}=[circle, draw=black, minimum size = 8mm]\tikzstyle{placeInOut}=[circle, draw=orange, minimum size = 8mm]\node at (0,-2) [placeInOut] (first_1){$\var{1}$};\node at (0,-4) [placeInOut] (first_2){$\var{2}$};\node at (5,-2) [place] (secondL1_1){$\frac{1}{\var{1}-\lani{1}{1}}$};\node at (5,-4) [place] (secondL1_2){$\frac{1}{\var{1}-\lani{1}{2}}$};\node at (5,-6) [place] (secondL2_1){$\frac{1}{\var{2}-\lani{2}{1}}$};\node at (5,-8) [place] (secondL2_2){$\frac{1}{\var{2}-\lani{2}{2}}$};\node at (5,-10) [place] (secondL2_3){$\frac{1}{\var{2}-\lani{2}{3}}$};\node at (5,-12) [place] (secondL2_4){$\frac{1}{\var{2}-\lani{2}{4}}$};\node at (10,-2) [place] (third_1){$\prod$};\node at (10,-4) [place] (third_2){$\prod$};\node at (10,-6) [place] (third_3){$\prod$};\node at (10,-8) [place] (third_4){$\prod$};\node at (10,-10) [place] (third_5){$\prod$};\node at (10,-12) [place] (third_6){$\prod$};\node at (10,-14) [place] (third_7){$\prod$};\node at (10,-16) [place] (third_8){$\prod$};\node at (15,-9) [placeInOut] (output){$\bSigma$};\draw[->] (first_1)--(secondL1_1) node[above,sloped,pos=0.75] { };\draw[->] (first_1)--(secondL1_2) node[above,sloped,pos=0.75] { };\draw[->] (first_2)--(secondL2_1) node[above,sloped,pos=0.75] { };\draw[->] (first_2)--(secondL2_2) node[above,sloped,pos=0.75] { };\draw[->] (first_2)--(secondL2_3) node[above,sloped,pos=0.75] { };\draw[->] (first_2)--(secondL2_4) node[above,sloped,pos=0.75] { };\draw[->] (secondL1_1)--(third_1) node[above,sloped,pos=0.25] {};\draw[->] (secondL1_1)--(third_2) node[above,sloped,pos=0.25] {};\draw[->] (secondL1_1)--(third_3) node[above,sloped,pos=0.25] {};\draw[->] (secondL1_1)--(third_4) node[above,sloped,pos=0.25] {};\draw[->] (secondL1_2)--(third_5) node[above,sloped,pos=0.25] {};\draw[->] (secondL1_2)--(third_6) node[above,sloped,pos=0.25] {};\draw[->] (secondL1_2)--(third_7) node[above,sloped,pos=0.25] {};\draw[->] (secondL1_2)--(third_8) node[above,sloped,pos=0.25] {};\draw[->] (secondL2_1)--(third_1) node[above,sloped,pos=0.25] {};\draw[->] (secondL2_2)--(third_2) node[above,sloped,pos=0.25] {};\draw[->] (secondL2_3)--(third_3) node[above,sloped,pos=0.25] {};\draw[->] (secondL2_4)--(third_4) node[above,sloped,pos=0.25] {};\draw[->] (secondL2_1)--(third_5) node[above,sloped,pos=0.25] {};\draw[->] (secondL2_2)--(third_6) node[above,sloped,pos=0.25] {};\draw[->] (secondL2_3)--(third_7) node[above,sloped,pos=0.25] {};\draw[->] (secondL2_4)--(third_8) node[above,sloped,pos=0.25] {};\draw[->] (third_1)--(output) node[above,sloped,pos=0.25] {4.213e-11};\draw[->] (third_2)--(output) node[above,sloped,pos=0.25] {-0.038781};\draw[->] (third_3)--(output) node[above,sloped,pos=0.25] {0.28809};\draw[->] (third_4)--(output) node[above,sloped,pos=0.25] {-0.33333};\draw[->] (third_5)--(output) node[above,sloped,pos=0.25] {-0.4213};\draw[->] (third_6)--(output) node[above,sloped,pos=0.25] {1.2703};\draw[->] (third_7)--(output) node[above,sloped,pos=0.25] {-1.4316};\draw[->] (third_8)--(output) node[above,sloped,pos=0.25] {0.66667};\end{tikzpicture}} \caption{Equivalent NN representation of the numerator $\bn_{\textrm{lag}}$.}\end{center}\end{figure}

\newpage \subsection{Function \#12 (${\ord=2}$ variables, tensor size: 12.5 \textbf{KB})} $$\frac{\var{1}^2+\var{2}^2+\var{1}-\var{2}-1}{(\var{1}-1.1)(\var{2}-1.1)}$$ \subsubsection{Setup and results overview}\begin{itemize}\item Reference: A/al. 2021 (A.5.7), \cite{Austin:2021}\item Domain: $\mathbb{R}$\item Tensor size: 12.5 \textbf{KB} ($40^{2}$ points)\item Bounds: $ \left(\begin{array}{cc} -1 & 1 \end{array}\right) \times \left(\begin{array}{cc} -1 & 1 \end{array}\right)$ \end{itemize} \begin{table}[H] \centering \begin{tabular}{llllll}
$\#$ & Alg. & Parameters & Dim. & CPU [s] & RMSE \\ 
\hline 
$\mathbf{\#12}$ & A/G/P-V 2025 (A1) & $0.01,3$ & $\mathbf{36}$ & $\mathbf{0.0087}$ & $1.1 \cdot 10^{-14}$ \\ 
 & A/G/P-V 2025 (A2) & $1 \cdot 10^{-15},2$ & $36$ & $0.077$ & $1 \cdot 10^{-13}$ \\ 
 & MDSPACK v1.1.0 & $0.01,1$ & $36$ & $0.016$ & $1.5 \cdot 10^{-14}$ \\ 
 & P/P 2025 & $1,0.95,50,0.01,6,6,13$ & $2.4 \cdot 10^{02}$ & $1.3$ & $0.11$ \\ 
 & C-R/B/G 2023 & $0.001,20$ & $60$ & $0.016$ & $9.8 \cdot 10^{-13}$ \\ 
 & B/G 2025 & $0.001,20,4$ & $60$ & $0.021$ & $\mathbf{6.2 \cdot 10^{-15}}$ \\ 
 & TensorFlow & $$ & $2.6 \cdot 10^{02}$ & $14$ & $1.2$ \\ 
\hline 
\end{tabular} \caption{Function \#12: best model configuration and performances per methods.} \end{table}\begin{figure}[H] \centering  \includegraphics[width=\textwidth]{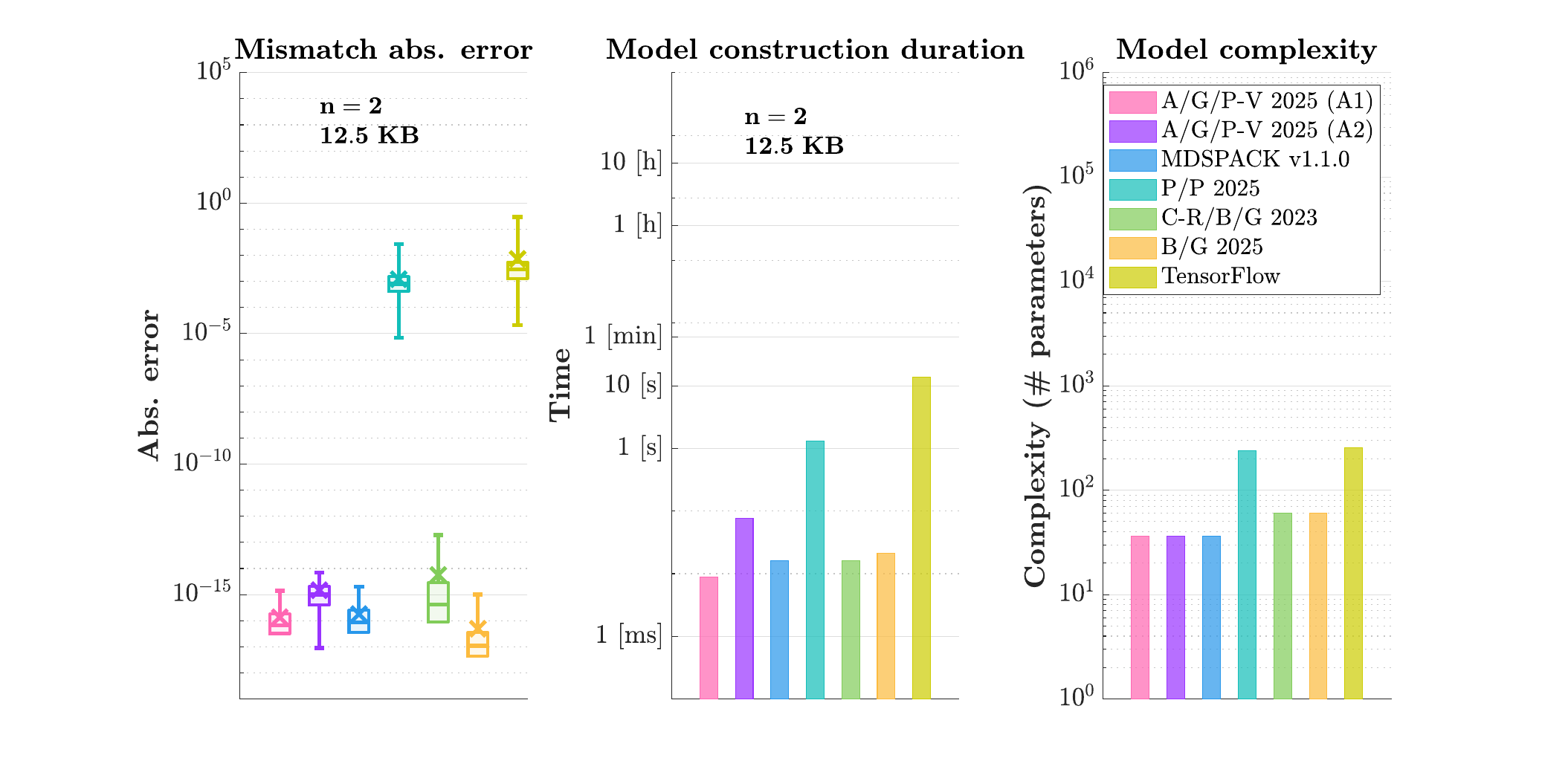} \caption{Function \#12: graphical view of the best model performances.} \end{figure}\begin{figure}[H] \centering  \includegraphics[width=\textwidth]{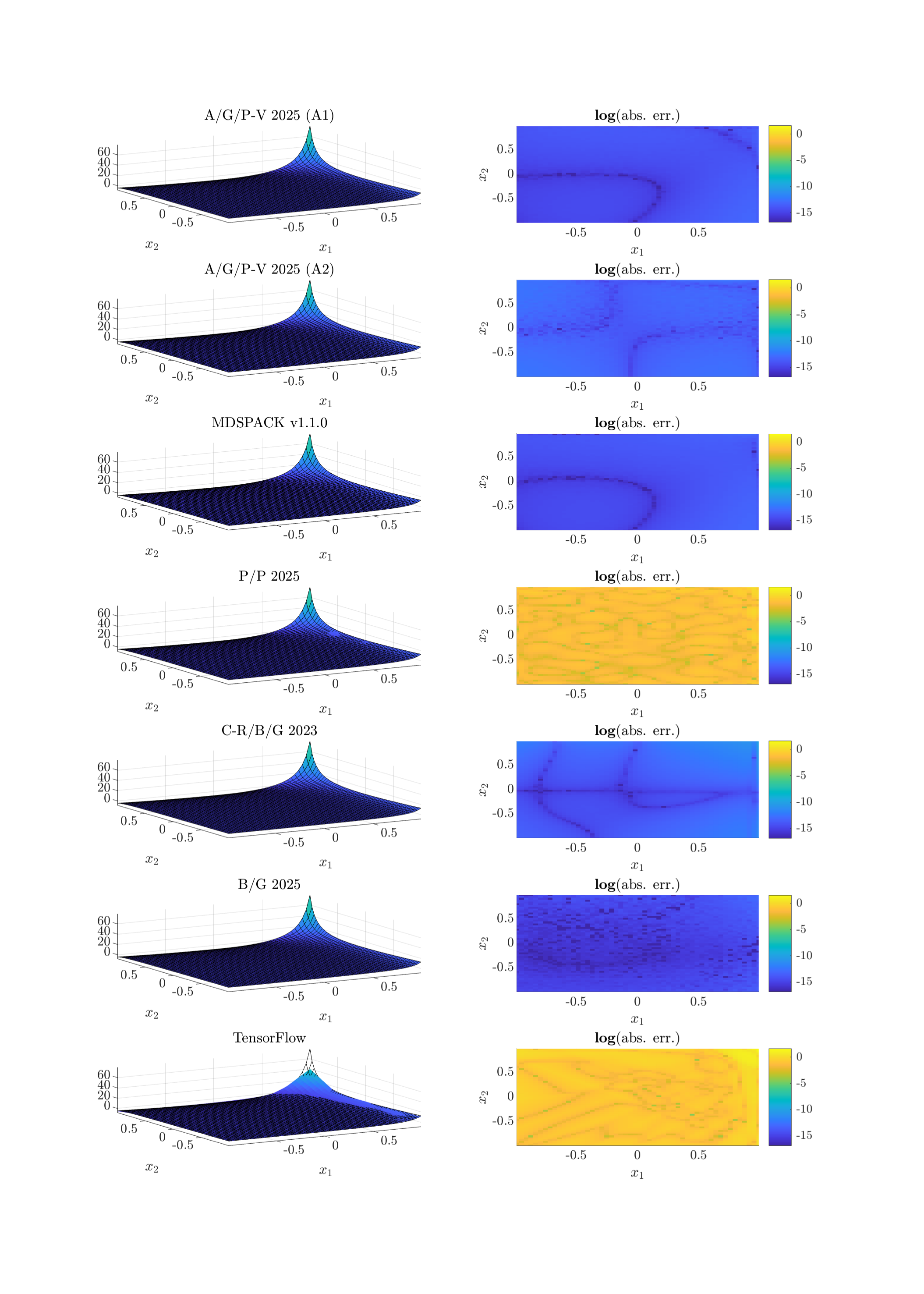} \caption{Function \#12: left side, evaluation of the original (mesh) vs. approximated (coloured surface) and right side, absolute errors (in log-scale).} \end{figure}\subsubsection{mLF detailed informations (M1)} \noindent \textbf{Right interpolation points} ($k_l=\left(\begin{array}{cc} 3 & 3 \end{array}\right)$, where $l=1,\cdots,\ord$):$$ \begin{array}{rcl}\lan{1} &=& \left(\begin{array}{ccc} -1 & -\frac{1}{19} & 1 \end{array}\right)\\\lan{2} &=& \left(\begin{array}{ccc} -1 & -\frac{1}{19} & 1 \end{array}\right)\\\end{array} $$\noindent \textbf{Lagrangian weights}: $$\left(\begin{array}{ccc} \bc & \bw & \bc\odot\bw\\ -22.37 & 0.2268 & -5.074\\ 23.33 & -0.3902 & -9.106\\ -0.9589 & -4.762 & 4.566\\ 23.33 & 0.3925 & 9.159\\ -24.33 & -0.7485 & 18.21\\ 1.0 & -9.108 & -9.108\\ -0.9589 & 14.29 & -13.7\\ 1.0 & 9.156 & 9.156\\ -0.0411 & 100.0 & -4.11 \end{array}\right)$$\noindent \textbf{Lagrangian form} (basis, numerator and denominator coefficients):$$\left(\begin{array}{ccc}\mathcal{B}_\textrm{lag}(\var{1},\var{2}) & \bN_\textrm{lag} &\bD_\textrm{lag}\end{array}\right) =$$ $$\left(\begin{array}{ccc} \left(\var{1}+1.0\right)\,\left(\var{2}+1.0\right) & -5.074 & -22.37\\ \left(\var{1}+1.0\right)\,\left(\var{2}+0.05263\right) & -9.106 & 23.33\\ \left(\var{1}+1.0\right)\,\left(\var{2}-1.0\right) & 4.566 & -0.9589\\ \left(\var{2}+1.0\right)\,\left(\var{1}+0.05263\right) & 9.159 & 23.33\\ \left(\var{1}+0.05263\right)\,\left(\var{2}+0.05263\right) & 18.21 & -24.33\\ \left(\var{2}-1.0\right)\,\left(\var{1}+0.05263\right) & -9.108 & 1.0\\ \left(\var{1}-1.0\right)\,\left(\var{2}+1.0\right) & -13.7 & -0.9589\\ \left(\var{1}-1.0\right)\,\left(\var{2}+0.05263\right) & 9.156 & 1.0\\ \left(\var{1}-1.0\right)\,\left(\var{2}-1.0\right) & -4.11 & -0.0411 \end{array}\right).$$\noindent The corresponding function is:$$\begin{array}{rcl}\bG_{\textrm{lag}}(\var{1},\var{2}) &=& \dfrac{\bn_{\textrm{lag}}(\var{1},\var{2})}{\bd_{\textrm{lag}}(\var{1},\var{2})}\\ && \\&=& \dfrac{\sum_{\textrm{row}} \bN_\textrm{lag} \odot\mathcal{B}^{-1}_\textrm{lag}(\var{1},\var{2})}{\sum_{\textrm{row}} \bD_\textrm{lag} \odot\mathcal{B}^{-1}_\textrm{lag}(\var{1},\var{2})}, \end{array}$$\noindent where,\\$\bn_{\textrm{lag}}(\var{1},\var{2}) = 4.232 \cdot 10^{-15}\,{\var{1}}^2\,{\var{2}}^2-5.462 \cdot 10^{-15}\,{\var{1}}^2\,\var{2}+0.8264\,{\var{1}}^2+1.453 \cdot 10^{-15}\,\var{1}\,{\var{2}}^2+0.8264\,\var{1}+0.8264\,{\var{2}}^2-0.8264\,\var{2}-0.8264$ \\~~\\$\bd_{\textrm{lag}}(\var{1},\var{2}) = 0.008264\,\left(10.0\,\var{1}-11.0\right)\,\left(10.0\,\var{2}-11.0\right)$ \\~~\\\noindent \textbf{Monomial form} (basis, numerator and denominator coefficients - entries $<10^{-12}$ removed):$$\left(\begin{array}{ccc}\mathcal{B}_\textrm{mon}(\var{1},\var{2}) & \bN_\textrm{mon} &\bD_\textrm{mon}\end{array}\right) =$$ $$\left(\begin{array}{ccc} {\var{1}}^2\,{\var{2}}^2 & 0 & 0\\ {\var{1}}^2\,\var{2} & 0 & 0\\ {\var{1}}^2 & -0.8264 & 0\\ \var{1}\,{\var{2}}^2 & 0 & 0\\ \var{1}\,\var{2} & 0 & -0.8264\\ \var{1} & -0.8264 & 0.9091\\ {\var{2}}^2 & -0.8264 & 0\\ \var{2} & 0.8264 & 0.9091\\ 1.0 & 0.8264 & -1.0 \end{array}\right)$$\noindent The corresponding function is:$$\begin{array}{rcl}\bG_{\textrm{mon}}(\var{1},\var{2}) &=& \dfrac{\bn_{\textrm{mon}}(\var{1},\var{2})}{\bd_{\textrm{mon}}(\var{1},\var{2})}\\ && \\&=& \dfrac{\sum_{\textrm{row}} \bN_\textrm{mon} \odot \mathcal{B}_\textrm{mon}(\var{1},\var{2})}{\sum_{\textrm{row}} \bD_\textrm{mon} \odot\mathcal{B}_\textrm{mon}(\var{1},\var{2})},  \end{array}$$\noindent where,\\$\bn_{\textrm{mon}}(\var{1},\var{2}) = 0.8264\,{\var{1}}^2+0.8264\,\var{1}+0.8264\,{\var{2}}^2-0.8264\,\var{2}-0.8264$ \\~~\\$\bd_{\textrm{mon}}(\var{1},\var{2}) = 0.008264\,\left(10.0\,\var{1}-11.0\right)\,\left(10.0\,\var{2}-11.0\right)$ \\~~\\\noindent \textbf{KST equivalent decoupling pattern} (Barycentric weights $\bc^{\var{l}}$): $$\begin{array}{rclll}\var{2}&: & \left(\begin{array}{ccc} 23.33 & 23.33 & 23.33\\ -24.33 & -24.33 & -24.33\\ 1.0 & 1.0 & 1.0 \end{array}\right)& \textrm{vec}(.) & := \textbf{Bary}(\var{2}) \\\var{1}&: & \left(\begin{array}{c} -0.9589\\ 1.0\\ -0.0411 \end{array}\right)& \textrm{vec}(.) \otimes \bone_{k_{2}} & := \textbf{Bary}(\var{1}) \\\end{array}$$~\\ Then, with the above notations, one defines the following univariate vector functions:  $$ \left\{ \begin{array}{rcl}\bPhi_{1}(\var{1}) &:=& \textbf{Bary}(\var{1}) \odot \mathbf{Lag}(\var{1}) \\\bPhi_{2}(\var{2}) &:=& \textbf{Bary}(\var{2}) \odot \mathbf{Lag}(\var{2}) \\\end{array} \right. $$\noindent The corresponding function is:$$\begin{array}{rcl}\bG_{\textrm{kst}}(\var{1},\var{2}) &=& \dfrac{\bn_{\textrm{kst}}(\var{1},\var{2})}{\bd_{\textrm{kst}}(\var{1},\var{2})}\\ && \\ &=& \dfrac{\sum_{\text{rows}} \bw \odot \bPhi_{1}(\var{1}) \odot \cdots \odot\bPhi_{2}(\var{2})}{\sum_{\text{rows}} \bPhi_{1}(\var{1}) \odot \cdots \odot\bPhi_{2}(\var{2})} . \end{array}$$~\\ \noindent \textbf{KST-like univariate functions} (equivalent scaled univariate functions $\bphi_{1,\cdots,2}$): $$\left\{\begin{array}{rcrcl}z_{1} &=&\bphi_{1}(\var{1}) &=& -\cfrac{1.0\,\left(100.0\,{\var{1}}^2+100.0\,\var{1}-100.0\right)}{10.0\,\var{1}-11.0}\\z_{2} &=&\bphi_{2}(\var{2}) &=& \cfrac{-100.0\,{\var{2}}^2+100.0\,\var{2}+100.0}{210.0\,\var{2}-231.0}\\\end{array} \right. .$$\noindent \textbf{Connection with Neural Networks} (equivalent numerator $\bn_{\textrm{lag}}$ representation):\\ \begin{figure}[H]\begin{center} \scalebox{.7}{\begin{tikzpicture}[line width=0.4mm]\tikzstyle{place}=[circle, draw=black, minimum size = 8mm]\tikzstyle{placeInOut}=[circle, draw=orange, minimum size = 8mm]\node at (0,-2) [placeInOut] (first_1){$\var{1}$};\node at (0,-4) [placeInOut] (first_2){$\var{2}$};\node at (5,-2) [place] (secondL1_1){$\frac{1}{\var{1}-\lani{1}{1}}$};\node at (5,-4) [place] (secondL1_2){$\frac{1}{\var{1}-\lani{1}{2}}$};\node at (5,-6) [place] (secondL1_3){$\frac{1}{\var{1}-\lani{1}{3}}$};\node at (5,-8) [place] (secondL2_1){$\frac{1}{\var{2}-\lani{2}{1}}$};\node at (5,-10) [place] (secondL2_2){$\frac{1}{\var{2}-\lani{2}{2}}$};\node at (5,-12) [place] (secondL2_3){$\frac{1}{\var{2}-\lani{2}{3}}$};\node at (10,-2) [place] (third_1){$\prod$};\node at (10,-4) [place] (third_2){$\prod$};\node at (10,-6) [place] (third_3){$\prod$};\node at (10,-8) [place] (third_4){$\prod$};\node at (10,-10) [place] (third_5){$\prod$};\node at (10,-12) [place] (third_6){$\prod$};\node at (10,-14) [place] (third_7){$\prod$};\node at (10,-16) [place] (third_8){$\prod$};\node at (10,-18) [place] (third_9){$\prod$};\node at (15,-10) [placeInOut] (output){$\bSigma$};\draw[->] (first_1)--(secondL1_1) node[above,sloped,pos=0.75] { };\draw[->] (first_1)--(secondL1_2) node[above,sloped,pos=0.75] { };\draw[->] (first_1)--(secondL1_3) node[above,sloped,pos=0.75] { };\draw[->] (first_2)--(secondL2_1) node[above,sloped,pos=0.75] { };\draw[->] (first_2)--(secondL2_2) node[above,sloped,pos=0.75] { };\draw[->] (first_2)--(secondL2_3) node[above,sloped,pos=0.75] { };\draw[->] (secondL1_1)--(third_1) node[above,sloped,pos=0.25] {};\draw[->] (secondL1_1)--(third_2) node[above,sloped,pos=0.25] {};\draw[->] (secondL1_1)--(third_3) node[above,sloped,pos=0.25] {};\draw[->] (secondL1_2)--(third_4) node[above,sloped,pos=0.25] {};\draw[->] (secondL1_2)--(third_5) node[above,sloped,pos=0.25] {};\draw[->] (secondL1_2)--(third_6) node[above,sloped,pos=0.25] {};\draw[->] (secondL1_3)--(third_7) node[above,sloped,pos=0.25] {};\draw[->] (secondL1_3)--(third_8) node[above,sloped,pos=0.25] {};\draw[->] (secondL1_3)--(third_9) node[above,sloped,pos=0.25] {};\draw[->] (secondL2_1)--(third_1) node[above,sloped,pos=0.25] {};\draw[->] (secondL2_2)--(third_2) node[above,sloped,pos=0.25] {};\draw[->] (secondL2_3)--(third_3) node[above,sloped,pos=0.25] {};\draw[->] (secondL2_1)--(third_4) node[above,sloped,pos=0.25] {};\draw[->] (secondL2_2)--(third_5) node[above,sloped,pos=0.25] {};\draw[->] (secondL2_3)--(third_6) node[above,sloped,pos=0.25] {};\draw[->] (secondL2_1)--(third_7) node[above,sloped,pos=0.25] {};\draw[->] (secondL2_2)--(third_8) node[above,sloped,pos=0.25] {};\draw[->] (secondL2_3)--(third_9) node[above,sloped,pos=0.25] {};\draw[->] (third_1)--(output) node[above,sloped,pos=0.25] {-5.0736};\draw[->] (third_2)--(output) node[above,sloped,pos=0.25] {-9.1057};\draw[->] (third_3)--(output) node[above,sloped,pos=0.25] {4.5662};\draw[->] (third_4)--(output) node[above,sloped,pos=0.25] {9.1591};\draw[->] (third_5)--(output) node[above,sloped,pos=0.25] {18.2141};\draw[->] (third_6)--(output) node[above,sloped,pos=0.25] {-9.1084};\draw[->] (third_7)--(output) node[above,sloped,pos=0.25] {-13.6986};\draw[->] (third_8)--(output) node[above,sloped,pos=0.25] {9.1565};\draw[->] (third_9)--(output) node[above,sloped,pos=0.25] {-4.1096};\end{tikzpicture}} \caption{Equivalent NN representation of the numerator $\bn_{\textrm{lag}}$.}\end{center}\end{figure}

\newpage \subsection{Function \#13 (${\ord=2}$ variables, tensor size: 12.5 \textbf{KB})} $$\frac{\var{1}^4+\var{2}^4+\var{1}^2\var{2}^2+\var{1}\var{2}}{(\var{1}-1.1)(\var{2}-1.1)}$$ \subsubsection{Setup and results overview}\begin{itemize}\item Reference: A/al. 2021 (A.5.8), \cite{Austin:2021}\item Domain: $\mathbb{R}$\item Tensor size: 12.5 \textbf{KB} ($40^{2}$ points)\item Bounds: $ \left(\begin{array}{cc} -1 & 1 \end{array}\right) \times \left(\begin{array}{cc} -1 & 1 \end{array}\right)$ \end{itemize} \begin{table}[H] \centering \begin{tabular}{llllll}
$\#$ & Alg. & Parameters & Dim. & CPU [s] & RMSE \\ 
\hline 
$\mathbf{\#13}$ & A/G/P-V 2025 (A1) & $0.001,2$ & $\mathbf{1 \cdot 10^{02}}$ & $\mathbf{0.011}$ & $5.2 \cdot 10^{-13}$ \\ 
 & A/G/P-V 2025 (A2) & $1 \cdot 10^{-15},2$ & $1 \cdot 10^{02}$ & $0.11$ & $2.6 \cdot 10^{-11}$ \\ 
 & MDSPACK v1.1.0 & $0.01,1$ & $1 \cdot 10^{02}$ & $0.011$ & $4.1 \cdot 10^{-13}$ \\ 
 & P/P 2025 & $1,1,50,0.01,6,12,13$ & $3.2 \cdot 10^{02}$ & $1.2$ & $0.34$ \\ 
 & C-R/B/G 2023 & $1 \cdot 10^{-06},20$ & $1.9 \cdot 10^{02}$ & $0.03$ & $1.2 \cdot 10^{-12}$ \\ 
 & B/G 2025 & $0.001,20,3$ & $1.4 \cdot 10^{02}$ & $0.058$ & $\mathbf{4.1 \cdot 10^{-14}}$ \\ 
 & TensorFlow & $$ & $2.6 \cdot 10^{02}$ & $14$ & $5.8$ \\ 
\hline 
\end{tabular} \caption{Function \#13: best model configuration and performances per methods.} \end{table}\begin{figure}[H] \centering  \includegraphics[width=\textwidth]{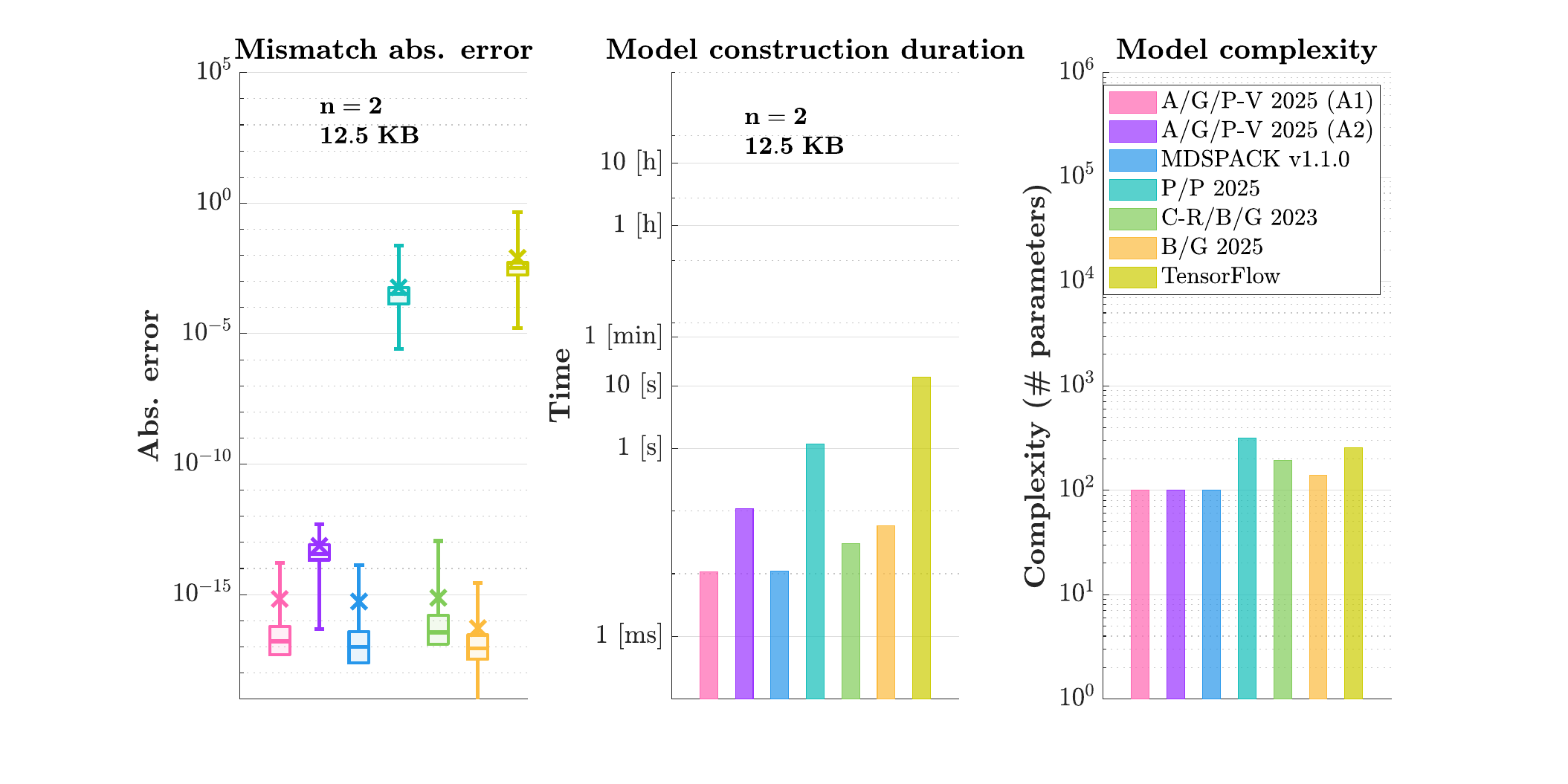} \caption{Function \#13: graphical view of the best model performances.} \end{figure}\begin{figure}[H] \centering  \includegraphics[width=\textwidth]{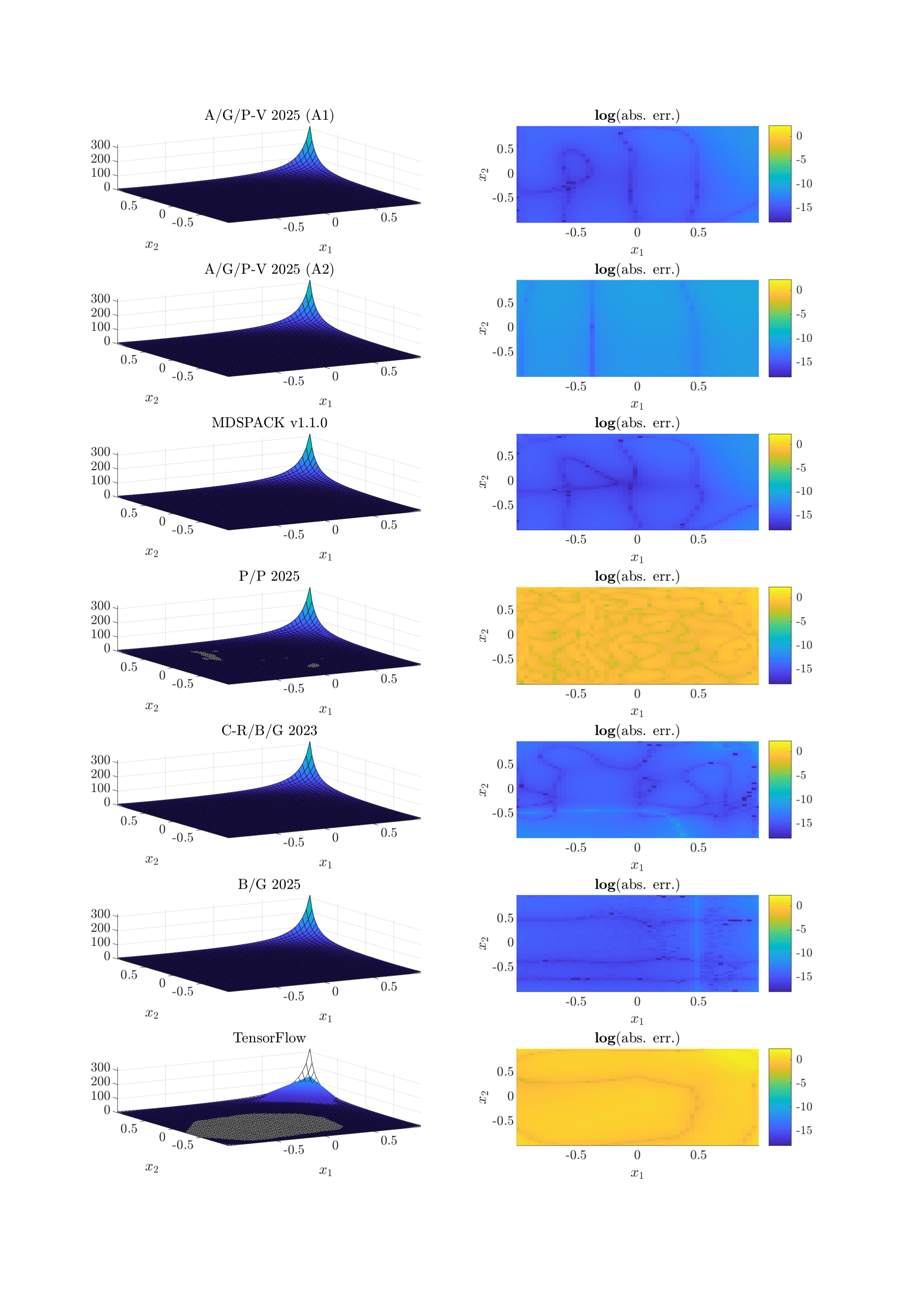} \caption{Function \#13: left side, evaluation of the original (mesh) vs. approximated (coloured surface) and right side, absolute errors (in log-scale).} \end{figure}\subsubsection{mLF detailed informations (M1)} \noindent \textbf{Right interpolation points} ($k_l=\left(\begin{array}{cc} 5 & 5 \end{array}\right)$, where $l=1,\cdots,\ord$):$$ \begin{array}{rcl}\lan{1} &=& \left(\begin{array}{ccccc} -1 & -\frac{11}{19} & -\frac{1}{19} & \frac{9}{19} & 1 \end{array}\right)\\\lan{2} &=& \left(\begin{array}{ccccc} -1 & -\frac{11}{19} & -\frac{1}{19} & \frac{9}{19} & 1 \end{array}\right)\\\end{array} $$\noindent \textbf{Lagrangian weights}: $$\left(\begin{array}{ccc} \bc & \bw & \bc\odot\bw\\ -12.25 & 0.907 & -11.11\\ 31.25 & 0.5748 & 17.96\\ -28.61 & 0.436 & -12.47\\ 9.992 & 0.609 & 6.086\\ -0.3918 & 9.524 & -3.732\\ 31.25 & 0.5748 & 17.96\\ -79.75 & 0.2385 & -19.02\\ 73.0 & 0.07428 & 5.423\\ -25.5 & -0.03456 & 0.8813\\ 1.0 & 5.173 & 5.173\\ -28.61 & 0.436 & -12.47\\ 73.0 & 0.07428 & 5.423\\ -66.82 & 0.002102 & -0.1405\\ 23.34 & 0.03608 & 0.8421\\ -0.9154 & 8.243 & -7.546\\ 9.992 & 0.609 & 6.086\\ -25.5 & -0.03456 & 0.8813\\ 23.34 & 0.03608 & 0.8421\\ -8.154 & 0.957 & -7.803\\ 0.3197 & 27.92 & 8.926\\ -0.3918 & 9.524 & -3.732\\ 1.0 & 5.173 & 5.173\\ -0.9154 & 8.243 & -7.546\\ 0.3197 & 27.92 & 8.926\\ -0.01254 & 400.0 & -5.016 \end{array}\right)$$\noindent \textbf{Lagrangian form} (basis, numerator and denominator coefficients):$$\left(\begin{array}{ccc}\mathcal{B}_\textrm{lag}(\var{1},\var{2}) & \bN_\textrm{lag} &\bD_\textrm{lag}\end{array}\right) =$$ $$\left(\begin{array}{ccc} \left(\var{1}+1.0\right)\,\left(\var{2}+1.0\right) & -11.11 & -12.25\\ \left(\var{1}+1.0\right)\,\left(\var{2}+0.5789\right) & 17.96 & 31.25\\ \left(\var{1}+1.0\right)\,\left(\var{2}+0.05263\right) & -12.47 & -28.61\\ \left(\var{1}+1.0\right)\,\left(\var{2}-0.4737\right) & 6.086 & 9.992\\ \left(\var{1}+1.0\right)\,\left(\var{2}-1.0\right) & -3.732 & -0.3918\\ \left(\var{2}+1.0\right)\,\left(\var{1}+0.5789\right) & 17.96 & 31.25\\ \left(\var{1}+0.5789\right)\,\left(\var{2}+0.5789\right) & -19.02 & -79.75\\ \left(\var{1}+0.5789\right)\,\left(\var{2}+0.05263\right) & 5.423 & 73.0\\ \left(\var{1}+0.5789\right)\,\left(\var{2}-0.4737\right) & 0.8813 & -25.5\\ \left(\var{2}-1.0\right)\,\left(\var{1}+0.5789\right) & 5.173 & 1.0\\ \left(\var{2}+1.0\right)\,\left(\var{1}+0.05263\right) & -12.47 & -28.61\\ \left(\var{2}+0.5789\right)\,\left(\var{1}+0.05263\right) & 5.423 & 73.0\\ \left(\var{1}+0.05263\right)\,\left(\var{2}+0.05263\right) & -0.1405 & -66.82\\ \left(\var{1}+0.05263\right)\,\left(\var{2}-0.4737\right) & 0.8421 & 23.34\\ \left(\var{2}-1.0\right)\,\left(\var{1}+0.05263\right) & -7.546 & -0.9154\\ \left(\var{2}+1.0\right)\,\left(\var{1}-0.4737\right) & 6.086 & 9.992\\ \left(\var{2}+0.5789\right)\,\left(\var{1}-0.4737\right) & 0.8813 & -25.5\\ \left(\var{2}+0.05263\right)\,\left(\var{1}-0.4737\right) & 0.8421 & 23.34\\ \left(\var{1}-0.4737\right)\,\left(\var{2}-0.4737\right) & -7.803 & -8.154\\ \left(\var{2}-1.0\right)\,\left(\var{1}-0.4737\right) & 8.926 & 0.3197\\ \left(\var{1}-1.0\right)\,\left(\var{2}+1.0\right) & -3.732 & -0.3918\\ \left(\var{1}-1.0\right)\,\left(\var{2}+0.5789\right) & 5.173 & 1.0\\ \left(\var{1}-1.0\right)\,\left(\var{2}+0.05263\right) & -7.546 & -0.9154\\ \left(\var{1}-1.0\right)\,\left(\var{2}-0.4737\right) & 8.926 & 0.3197\\ \left(\var{1}-1.0\right)\,\left(\var{2}-1.0\right) & -5.016 & -0.01254 \end{array}\right).$$\noindent The corresponding function is:$$\begin{array}{rcl}\bG_{\textrm{lag}}(\var{1},\var{2}) &=& \dfrac{\bn_{\textrm{lag}}(\var{1},\var{2})}{\bd_{\textrm{lag}}(\var{1},\var{2})}\\ && \\&=& \dfrac{\sum_{\textrm{row}} \bN_\textrm{lag} \odot\mathcal{B}^{-1}_\textrm{lag}(\var{1},\var{2})}{\sum_{\textrm{row}} \bD_\textrm{lag} \odot\mathcal{B}^{-1}_\textrm{lag}(\var{1},\var{2})}, \end{array}$$\noindent where,\\$\bn_{\textrm{lag}}(\var{1},\var{2}) = -1.309 \cdot 10^{-13}\,{\var{1}}^4\,{\var{2}}^4+2.835 \cdot 10^{-14}\,{\var{1}}^4\,{\var{2}}^3+6.973 \cdot 10^{-14}\,{\var{1}}^4\,{\var{2}}^2+3.124 \cdot 10^{-14}\,{\var{1}}^4\,\var{2}+0.8264\,{\var{1}}^4-1.333 \cdot 10^{-14}\,{\var{1}}^3\,{\var{2}}^4-1.719 \cdot 10^{-14}\,{\var{1}}^3\,{\var{2}}^3+1.08 \cdot 10^{-13}\,{\var{1}}^3\,{\var{2}}^2+1.731 \cdot 10^{-14}\,{\var{1}}^3\,\var{2}+7.291 \cdot 10^{-14}\,{\var{1}}^3+1.894 \cdot 10^{-13}\,{\var{1}}^2\,{\var{2}}^4-5.007 \cdot 10^{-14}\,{\var{1}}^2\,{\var{2}}^3+0.8264\,{\var{1}}^2\,{\var{2}}^2+4.32 \cdot 10^{-14}\,{\var{1}}^2\,\var{2}-1.803 \cdot 10^{-14}\,{\var{1}}^2+7.386 \cdot 10^{-14}\,\var{1}\,{\var{2}}^4-1.239 \cdot 10^{-14}\,\var{1}\,{\var{2}}^3-3.02 \cdot 10^{-14}\,\var{1}\,{\var{2}}^2+0.8264\,\var{1}\,\var{2}-1.316 \cdot 10^{-14}\,\var{1}+0.8264\,{\var{2}}^4+7.59 \cdot 10^{-15}\,{\var{2}}^3-1.075 \cdot 10^{-15}\,{\var{2}}^2-2.787 \cdot 10^{-16}\,\var{2}-6.271 \cdot 10^{-16}$ \\~~\\$\bd_{\textrm{lag}}(\var{1},\var{2}) = 3.384 \cdot 10^{-13}\,{\var{1}}^4\,{\var{2}}^4-1.002 \cdot 10^{-13}\,{\var{1}}^4\,{\var{2}}^3-5.232 \cdot 10^{-13}\,{\var{1}}^4\,{\var{2}}^2+1.496 \cdot 10^{-13}\,{\var{1}}^4\,\var{2}+1.318 \cdot 10^{-13}\,{\var{1}}^4+3.385 \cdot 10^{-16}\,{\var{1}}^3\,{\var{2}}^4-2.613 \cdot 10^{-14}\,{\var{1}}^3\,{\var{2}}^3+4.956 \cdot 10^{-14}\,{\var{1}}^3\,{\var{2}}^2+5.334 \cdot 10^{-14}\,{\var{1}}^3\,\var{2}-7.915 \cdot 10^{-14}\,{\var{1}}^3-4.264 \cdot 10^{-13}\,{\var{1}}^2\,{\var{2}}^4+1.212 \cdot 10^{-13}\,{\var{1}}^2\,{\var{2}}^3+5.695 \cdot 10^{-13}\,{\var{1}}^2\,{\var{2}}^2-1.319 \cdot 10^{-13}\,{\var{1}}^2\,\var{2}-1.284 \cdot 10^{-13}\,{\var{1}}^2-1.624 \cdot 10^{-14}\,\var{1}\,{\var{2}}^4+4.902 \cdot 10^{-14}\,\var{1}\,{\var{2}}^3+1.205 \cdot 10^{-14}\,\var{1}\,{\var{2}}^2+0.8264\,\var{1}\,\var{2}-0.9091\,\var{1}+1.018 \cdot 10^{-13}\,{\var{2}}^4-4.047 \cdot 10^{-14}\,{\var{2}}^3-1.049 \cdot 10^{-13}\,{\var{2}}^2-0.9091\,\var{2}+1.0$ \\~~\\\noindent \textbf{Monomial form} (basis, numerator and denominator coefficients - entries $<10^{-12}$ removed):$$\left(\begin{array}{ccc}\mathcal{B}_\textrm{mon}(\var{1},\var{2}) & \bN_\textrm{mon} &\bD_\textrm{mon}\end{array}\right) =$$ $$\left(\begin{array}{ccc} {\var{1}}^4\,{\var{2}}^4 & 0 & 0\\ {\var{1}}^4\,{\var{2}}^3 & 0 & 0\\ {\var{1}}^4\,{\var{2}}^2 & 0 & 0\\ {\var{1}}^4\,\var{2} & 0 & 0\\ {\var{1}}^4 & -0.8264 & 0\\ {\var{1}}^3\,{\var{2}}^4 & 0 & 0\\ {\var{1}}^3\,{\var{2}}^3 & 0 & 0\\ {\var{1}}^3\,{\var{2}}^2 & 0 & 0\\ {\var{1}}^3\,\var{2} & 0 & 0\\ {\var{1}}^3 & 0 & 0\\ {\var{1}}^2\,{\var{2}}^4 & 0 & 0\\ {\var{1}}^2\,{\var{2}}^3 & 0 & 0\\ {\var{1}}^2\,{\var{2}}^2 & -0.8264 & 0\\ {\var{1}}^2\,\var{2} & 0 & 0\\ {\var{1}}^2 & 0 & 0\\ \var{1}\,{\var{2}}^4 & 0 & 0\\ \var{1}\,{\var{2}}^3 & 0 & 0\\ \var{1}\,{\var{2}}^2 & 0 & 0\\ \var{1}\,\var{2} & -0.8264 & -0.8264\\ \var{1} & 0 & 0.9091\\ {\var{2}}^4 & -0.8264 & 0\\ {\var{2}}^3 & 0 & 0\\ {\var{2}}^2 & 0 & 0\\ \var{2} & 0 & 0.9091\\ 1.0 & 0 & -1.0 \end{array}\right)$$\noindent The corresponding function is:$$\begin{array}{rcl}\bG_{\textrm{mon}}(\var{1},\var{2}) &=& \dfrac{\bn_{\textrm{mon}}(\var{1},\var{2})}{\bd_{\textrm{mon}}(\var{1},\var{2})}\\ && \\&=& \dfrac{\sum_{\textrm{row}} \bN_\textrm{mon} \odot \mathcal{B}_\textrm{mon}(\var{1},\var{2})}{\sum_{\textrm{row}} \bD_\textrm{mon} \odot\mathcal{B}_\textrm{mon}(\var{1},\var{2})},  \end{array}$$\noindent where,\\$\bn_{\textrm{mon}}(\var{1},\var{2}) = 0.8264\,{\var{1}}^4+0.8264\,{\var{1}}^2\,{\var{2}}^2+0.8264\,\var{1}\,\var{2}+0.8264\,{\var{2}}^4$ \\~~\\$\bd_{\textrm{mon}}(\var{1},\var{2}) = 0.8264\,\var{1}\,\var{2}-0.9091\,\var{2}-0.9091\,\var{1}+1.0$ \\~~\\\noindent \textbf{KST equivalent decoupling pattern} (Barycentric weights $\bc^{\var{l}}$): $$\begin{array}{rclll}\var{2}&: & \left(\begin{array}{ccccc} 31.25 & 31.25 & 31.25 & 31.25 & 31.25\\ -79.75 & -79.75 & -79.75 & -79.75 & -79.75\\ 73.0 & 73.0 & 73.0 & 73.0 & 73.0\\ -25.5 & -25.5 & -25.5 & -25.5 & -25.5\\ 1.0 & 1.0 & 1.0 & 1.0 & 1.0 \end{array}\right)& \textrm{vec}(.) & := \textbf{Bary}(\var{2}) \\\var{1}&: & \left(\begin{array}{c} -0.3918\\ 1.0\\ -0.9154\\ 0.3197\\ -0.01254 \end{array}\right)& \textrm{vec}(.) \otimes \bone_{k_{2}} & := \textbf{Bary}(\var{1}) \\\end{array}$$~\\ Then, with the above notations, one defines the following univariate vector functions:  $$ \left\{ \begin{array}{rcl}\bPhi_{1}(\var{1}) &:=& \textbf{Bary}(\var{1}) \odot \mathbf{Lag}(\var{1}) \\\bPhi_{2}(\var{2}) &:=& \textbf{Bary}(\var{2}) \odot \mathbf{Lag}(\var{2}) \\\end{array} \right. $$\noindent The corresponding function is:$$\begin{array}{rcl}\bG_{\textrm{kst}}(\var{1},\var{2}) &=& \dfrac{\bn_{\textrm{kst}}(\var{1},\var{2})}{\bd_{\textrm{kst}}(\var{1},\var{2})}\\ && \\ &=& \dfrac{\sum_{\text{rows}} \bw \odot \bPhi_{1}(\var{1}) \odot \cdots \odot\bPhi_{2}(\var{2})}{\sum_{\text{rows}} \bPhi_{1}(\var{1}) \odot \cdots \odot\bPhi_{2}(\var{2})} . \end{array}$$~\\ \noindent \textbf{KST-like univariate functions} (equivalent scaled univariate functions $\bphi_{1,\cdots,2}$): $$\left\{\begin{array}{rcrcl}z_{1} &=&\bphi_{1}(\var{1}) &=& \cfrac{\bn_{1}}{\bd_{1}} \\z_{2} &=&\bphi_{2}(\var{2}) &=& \cfrac{\bn_{2}}{\bd_{2}} \\\end{array} \right. .$$\noindent where, \\ \noindent $\bn_{1}=9.091\,{\var{1}}^4+1.893 \cdot 10^{-12}\,{\var{1}}^3+9.091\,{\var{1}}^2+9.091\,\var{1}+9.091$ and \\ \noindent $\bd_{1}=-4.044 \cdot 10^{-14}\,{\var{1}}^4-2.248 \cdot 10^{-14}\,{\var{1}}^3+4.42 \cdot 10^{-14}\,{\var{1}}^2-0.9091\,\var{1}+1.0$, \\ \noindent $\bn_{2}=0.4329\,{\var{2}}^4+2.846 \cdot 10^{-14}\,{\var{2}}^3+0.4329\,{\var{2}}^2-0.4329\,\var{2}+0.4329$ and \\ \noindent $\bd_{2}=1.052 \cdot 10^{-13}\,{\var{2}}^4-6.218 \cdot 10^{-16}\,{\var{2}}^3-1.699 \cdot 10^{-13}\,{\var{2}}^2-0.9091\,\var{2}+1.0$, \\

\newpage \subsection{Function \#14 (${\ord=4}$ variables, tensor size: 1.22 \textbf{MB})} $$\frac{\var{1}^2+\var{2}^2+\var{1}-\var{2}+1}{(\var{3}-1.5)(\var{4}-1.5)}$$ \subsubsection{Setup and results overview}\begin{itemize}\item Reference: A/al. 2021 (A.5.9), \cite{Austin:2021}\item Domain: $\mathbb{R}$\item Tensor size: 1.22 \textbf{MB} ($20^{4}$ points)\item Bounds: $ \left(\begin{array}{cc} -1 & 1 \end{array}\right) \times \left(\begin{array}{cc} -1 & 1 \end{array}\right) \times \left(\begin{array}{cc} -1 & 1 \end{array}\right) \times \left(\begin{array}{cc} -1 & 1 \end{array}\right)$ \end{itemize} \begin{table}[H] \centering \begin{tabular}{llllll}
$\#$ & Alg. & Parameters & Dim. & CPU [s] & RMSE \\ 
\hline 
$\mathbf{\#14}$ & A/G/P-V 2025 (A1) & $0.1,3$ & $\mathbf{2.2 \cdot 10^{02}}$ & $\mathbf{0.044}$ & $\mathbf{5.2 \cdot 10^{-16}}$ \\ 
 & A/G/P-V 2025 (A2) & $1 \cdot 10^{-15},1$ & $2.2 \cdot 10^{02}$ & $8.2$ & $7.9 \cdot 10^{-16}$ \\ 
 & MDSPACK v1.1.0 & $0.01,1$ & $2.2 \cdot 10^{02}$ & $0.048$ & $6.6 \cdot 10^{-16}$ \\ 
 & P/P 2025 & $1,0.95,50,0.01,6,12,13$ & $4.7 \cdot 10^{02}$ & $49$ & $0.00074$ \\ 
 & C-R/B/G 2023 & $0.001,20$ & $2.2 \cdot 10^{02}$ & $41$ & $7.8 \cdot 10^{-15}$ \\ 
 & B/G 2025 & $0.001,20,3$ & $2.2 \cdot 10^{02}$ & $2$ & $5.2 \cdot 10^{-15}$ \\ 
 & TensorFlow & $$ & $3.8 \cdot 10^{02}$ & $1 \cdot 10^{03}$ & $0.017$ \\ 
\hline 
\end{tabular} \caption{Function \#14: best model configuration and performances per methods.} \end{table}\begin{figure}[H] \centering  \includegraphics[width=\textwidth]{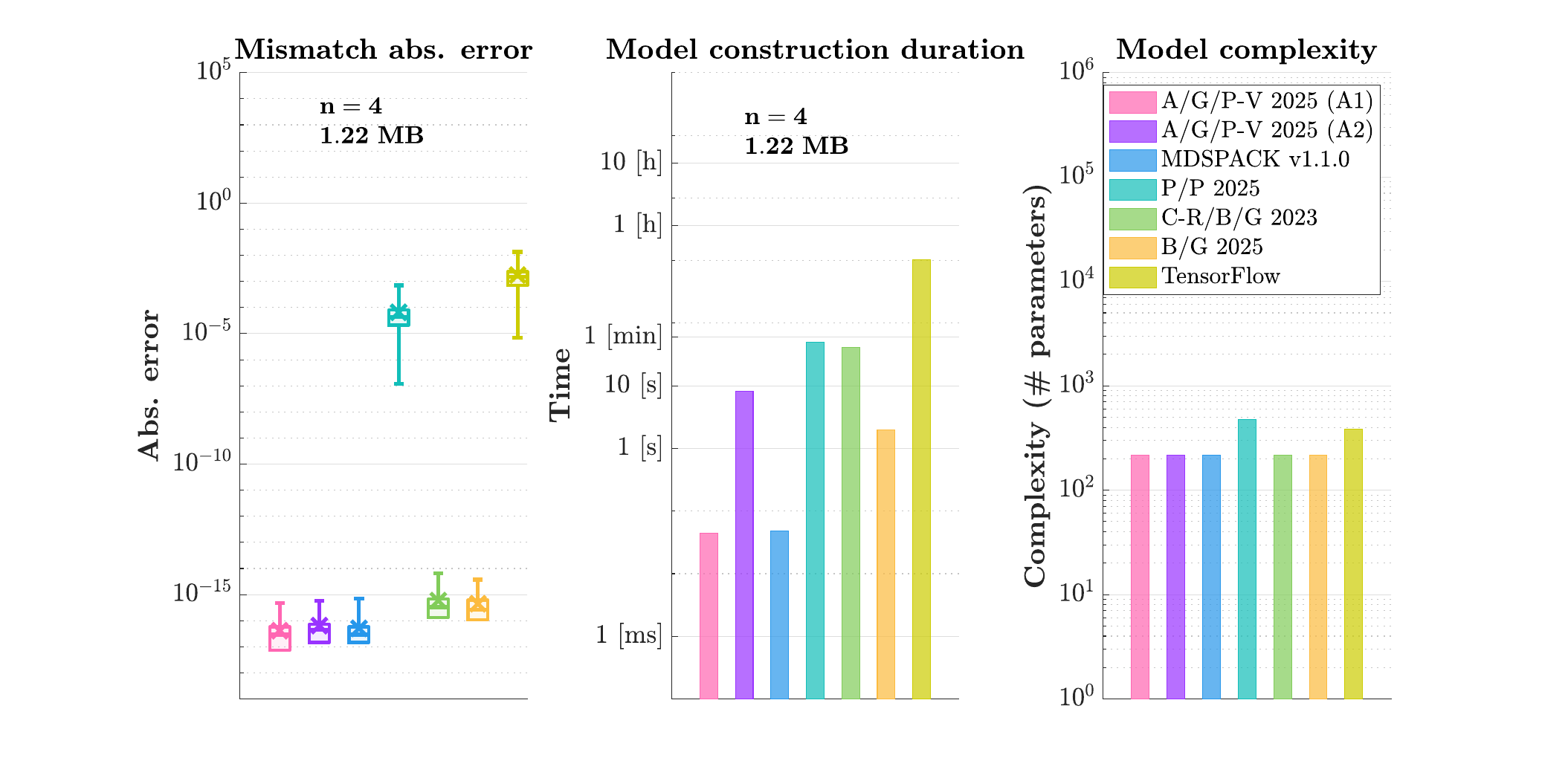} \caption{Function \#14: graphical view of the best model performances.} \end{figure}\begin{figure}[H] \centering  \includegraphics[width=\textwidth]{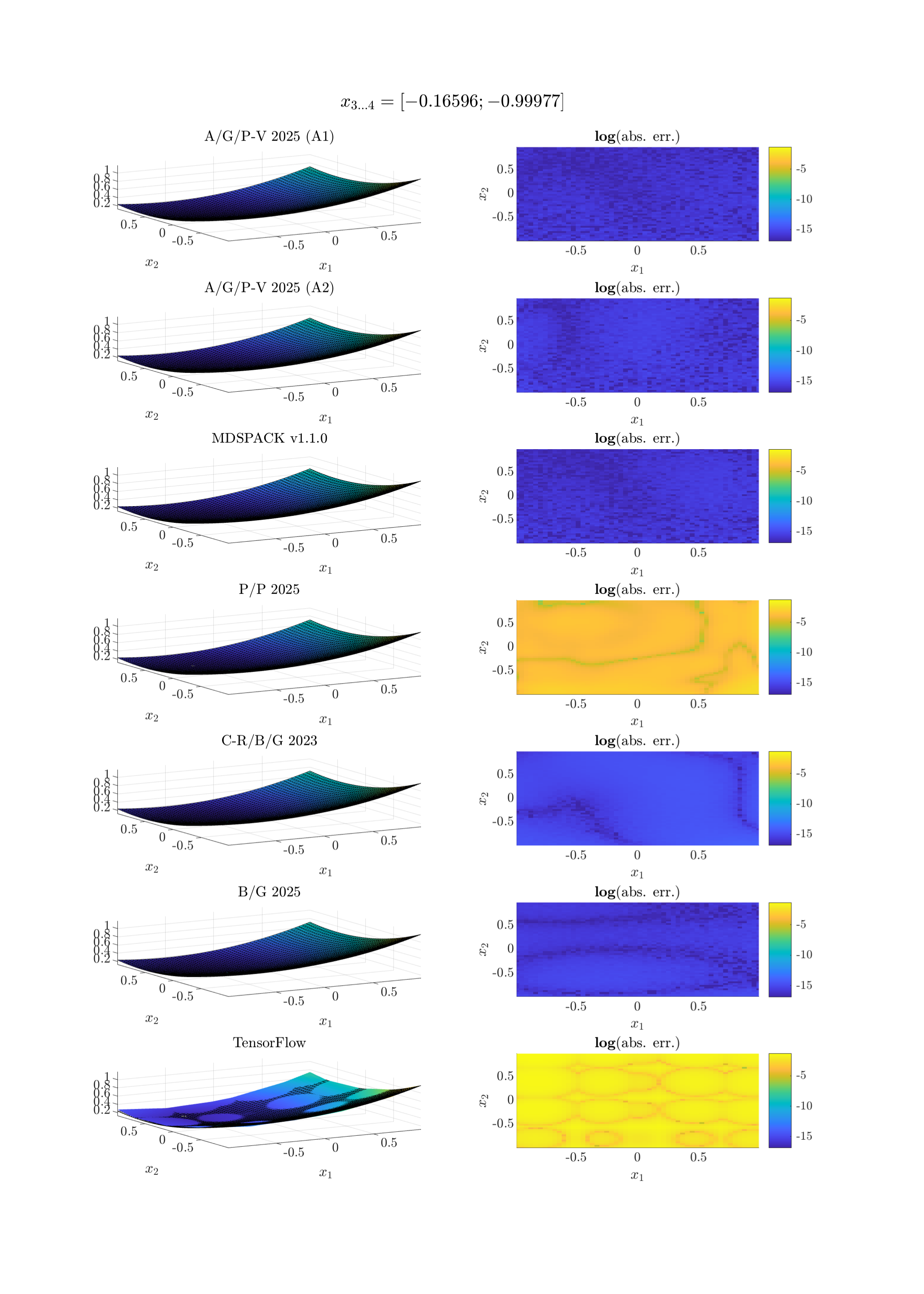} \caption{Function \#14: left side, evaluation of the original (mesh) vs. approximated (coloured surface) and right side, absolute errors (in log-scale).} \end{figure}\subsubsection{mLF detailed informations (M1)} \noindent \textbf{Right interpolation points}: $k_l=\left(\begin{array}{cccc} 3 & 3 & 2 & 2 \end{array}\right)$, where $l=1,\cdots,\ord$.$$ \begin{array}{rcl}\lan{1} &\in& \IC^{3} \text{ , linearly spaced between bounds}\\\lan{2} &\in& \IC^{3} \text{ , linearly spaced between bounds}\\\lan{3} &\in& \IC^{2} \text{ , linearly spaced between bounds}\\\lan{4} &\in& \IC^{2} \text{ , linearly spaced between bounds}\\\end{array} $$\noindent \textbf{$\ord$-D Loewner matrix, barycentric weights and Lagrangian basis}:$$ \begin{array}{rcl}\IL & \in & \IC^{36 \times 36}\\\bc & \in & \IC^{36}\\\bw & \in & \IC^{36}\\\bc\odot \bw & \in & \IC^{36}\\\mathbf{Lag}(\var{1},\var{2},\var{3},\var{4}) & \in & \IC^{36}\\\end{array} $$

\newpage \subsection{Function \#15 (${\ord=2}$ variables, tensor size: 12.5 \textbf{KB})} $$\frac{\var{1}^2+\var{2}^2+\var{1}-\var{2}-1}{\var{1}^3+\var{2}^3+4}$$ \subsubsection{Setup and results overview}\begin{itemize}\item Reference: A/al. 2021 (A.5.10), \cite{Austin:2021}\item Domain: $\mathbb{R}$\item Tensor size: 12.5 \textbf{KB} ($40^{2}$ points)\item Bounds: $ \left(\begin{array}{cc} -1 & 1 \end{array}\right) \times \left(\begin{array}{cc} -1 & 1 \end{array}\right)$ \end{itemize} \begin{table}[H] \centering \begin{tabular}{llllll}
$\#$ & Alg. & Parameters & Dim. & CPU [s] & RMSE \\ 
\hline 
$\mathbf{\#15}$ & A/G/P-V 2025 (A1) & $0.01,1$ & $\mathbf{64}$ & $0.012$ & $\mathbf{2.8 \cdot 10^{-15}}$ \\ 
 & A/G/P-V 2025 (A2) & $1 \cdot 10^{-15},3$ & $64$ & $0.091$ & $1.4 \cdot 10^{-13}$ \\ 
 & MDSPACK v1.1.0 & $0.0001,2$ & $64$ & $\mathbf{0.011}$ & $2.8 \cdot 10^{-15}$ \\ 
 & P/P 2025 & $1,0.95,50,0.01,4,6,9$ & $1.3 \cdot 10^{02}$ & $0.21$ & $6.7 \cdot 10^{-05}$ \\ 
 & C-R/B/G 2023 & $0.001,20$ & $1.4 \cdot 10^{02}$ & $0.038$ & $4.5 \cdot 10^{-14}$ \\ 
 & B/G 2025 & $1 \cdot 10^{-09},20,3$ & $1.1 \cdot 10^{02}$ & $0.071$ & $8.9 \cdot 10^{-15}$ \\ 
 & TensorFlow & $$ & $2.6 \cdot 10^{02}$ & $14$ & $0.0034$ \\ 
\hline 
\end{tabular} \caption{Function \#15: best model configuration and performances per methods.} \end{table}\begin{figure}[H] \centering  \includegraphics[width=\textwidth]{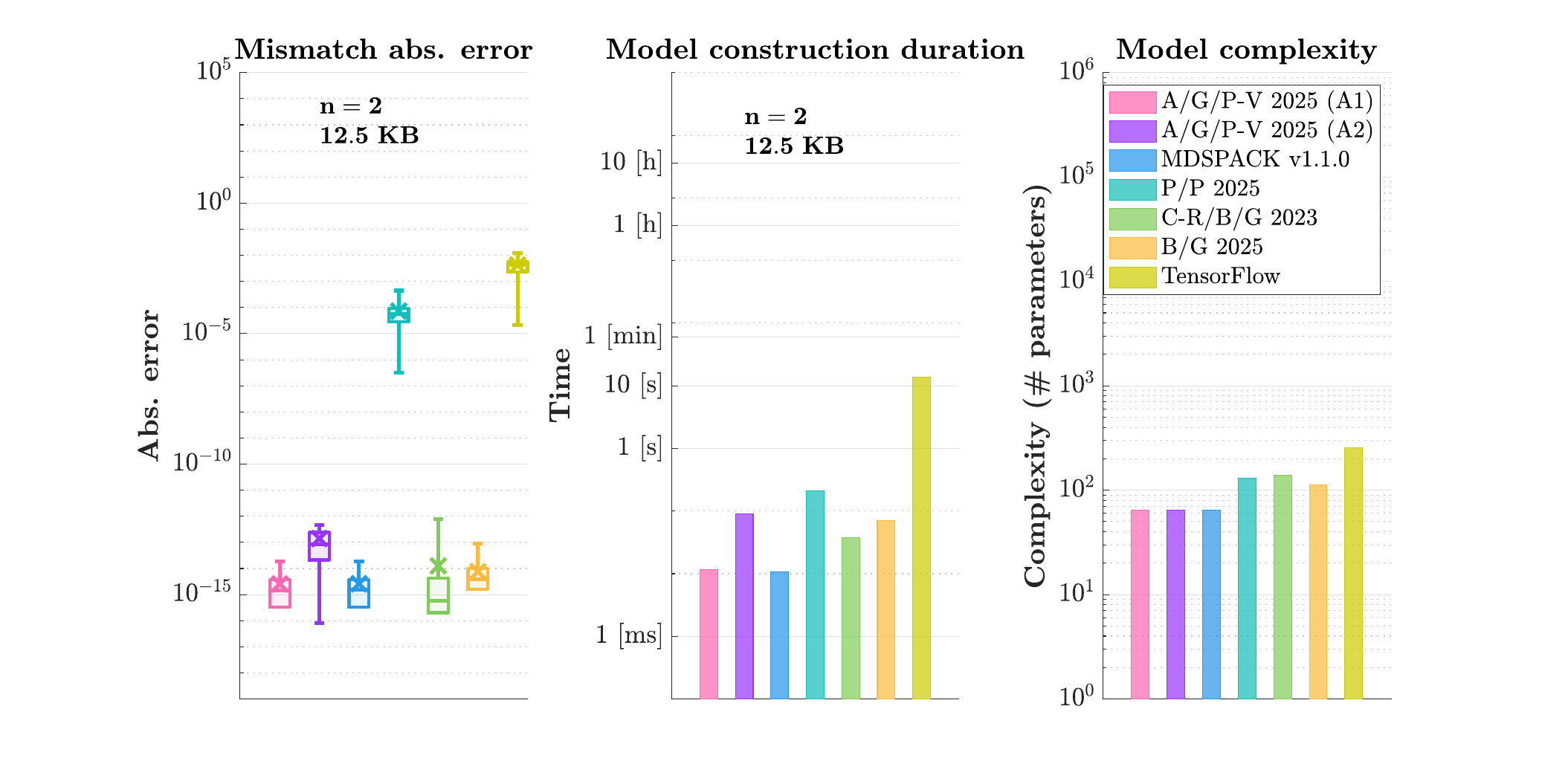} \caption{Function \#15: graphical view of the best model performances.} \end{figure}\begin{figure}[H] \centering  \includegraphics[width=\textwidth]{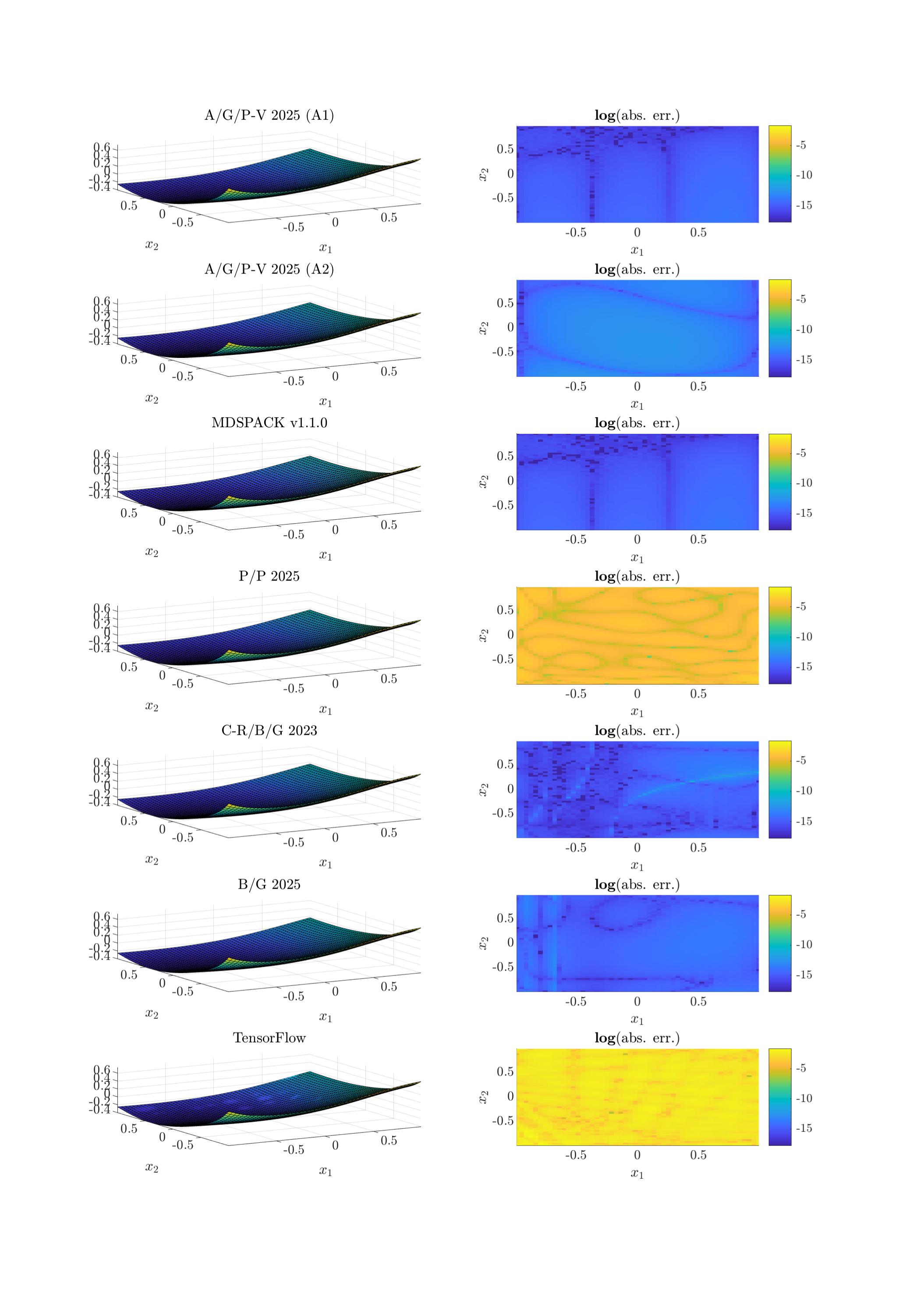} \caption{Function \#15: left side, evaluation of the original (mesh) vs. approximated (coloured surface) and right side, absolute errors (in log-scale).} \end{figure}\subsubsection{mLF detailed informations (M1)} \noindent \textbf{Right interpolation points} ($k_l=\left(\begin{array}{cc} 4 & 4 \end{array}\right)$, where $l=1,\cdots,\ord$):$$ \begin{array}{rcl}\lan{1} &=& \left(\begin{array}{cccc} -1 & -\frac{7}{19} & \frac{5}{19} & 1 \end{array}\right)\\\lan{2} &=& \left(\begin{array}{cccc} -1 & -\frac{7}{19} & \frac{5}{19} & 1 \end{array}\right)\\\end{array} $$\noindent \textbf{Lagrangian weights}: $$\left(\begin{array}{ccc} \bc & \bw & \bc\odot\bw\\ 0.1747 & 0.5 & 0.08735\\ -0.7532 & -0.1681 & 0.1266\\ 0.7156 & -0.3956 & -0.2831\\ -0.2764 & -0.25 & 0.06911\\ -0.7532 & 0.2601 & -0.1959\\ 2.911 & -0.1868 & -0.5437\\ -2.75 & -0.3595 & 0.9887\\ 1.0 & -0.249 & -0.249\\ 0.7156 & 0.4415 & 0.3159\\ -2.75 & -0.04119 & 0.1133\\ 2.598 & -0.2134 & -0.5544\\ -0.9414 & -0.133 & 0.1252\\ -0.2764 & 0.75 & -0.2073\\ 1.0 & 0.3039 & 0.3039\\ -0.9414 & 0.1606 & -0.1512\\ 0.3281 & 0.1667 & 0.05468 \end{array}\right)$$\noindent \textbf{Lagrangian form} (basis, numerator and denominator coefficients):$$\left(\begin{array}{ccc}\mathcal{B}_\textrm{lag}(\var{1},\var{2}) & \bN_\textrm{lag} &\bD_\textrm{lag}\end{array}\right) =$$ $$\left(\begin{array}{ccc} \left(\var{1}+1.0\right)\,\left(\var{2}+1.0\right) & 0.08735 & 0.1747\\ \left(\var{1}+1.0\right)\,\left(\var{2}+0.3684\right) & 0.1266 & -0.7532\\ \left(\var{1}+1.0\right)\,\left(\var{2}-0.2632\right) & -0.2831 & 0.7156\\ \left(\var{1}+1.0\right)\,\left(\var{2}-1.0\right) & 0.06911 & -0.2764\\ \left(\var{2}+1.0\right)\,\left(\var{1}+0.3684\right) & -0.1959 & -0.7532\\ \left(\var{1}+0.3684\right)\,\left(\var{2}+0.3684\right) & -0.5437 & 2.911\\ \left(\var{2}-0.2632\right)\,\left(\var{1}+0.3684\right) & 0.9887 & -2.75\\ \left(\var{2}-1.0\right)\,\left(\var{1}+0.3684\right) & -0.249 & 1.0\\ \left(\var{2}+1.0\right)\,\left(\var{1}-0.2632\right) & 0.3159 & 0.7156\\ \left(\var{1}-0.2632\right)\,\left(\var{2}+0.3684\right) & 0.1133 & -2.75\\ \left(\var{1}-0.2632\right)\,\left(\var{2}-0.2632\right) & -0.5544 & 2.598\\ \left(\var{2}-1.0\right)\,\left(\var{1}-0.2632\right) & 0.1252 & -0.9414\\ \left(\var{1}-1.0\right)\,\left(\var{2}+1.0\right) & -0.2073 & -0.2764\\ \left(\var{1}-1.0\right)\,\left(\var{2}+0.3684\right) & 0.3039 & 1.0\\ \left(\var{1}-1.0\right)\,\left(\var{2}-0.2632\right) & -0.1512 & -0.9414\\ \left(\var{1}-1.0\right)\,\left(\var{2}-1.0\right) & 0.05468 & 0.3281 \end{array}\right).$$\noindent The corresponding function is:$$\begin{array}{rcl}\bG_{\textrm{lag}}(\var{1},\var{2}) &=& \dfrac{\bn_{\textrm{lag}}(\var{1},\var{2})}{\bd_{\textrm{lag}}(\var{1},\var{2})}\\ && \\&=& \dfrac{\sum_{\textrm{row}} \bN_\textrm{lag} \odot\mathcal{B}^{-1}_\textrm{lag}(\var{1},\var{2})}{\sum_{\textrm{row}} \bD_\textrm{lag} \odot\mathcal{B}^{-1}_\textrm{lag}(\var{1},\var{2})}, \end{array}$$\noindent where,\\$\bn_{\textrm{lag}}(\var{1},\var{2}) = -2.397 \cdot 10^{-14}\,{\var{1}}^3\,{\var{2}}^3+2.866 \cdot 10^{-14}\,{\var{1}}^3\,{\var{2}}^2+3.027 \cdot 10^{-15}\,{\var{1}}^3\,\var{2}-1.348 \cdot 10^{-14}\,{\var{1}}^3-7.299 \cdot 10^{-15}\,{\var{1}}^2\,{\var{2}}^3+1.793 \cdot 10^{-14}\,{\var{1}}^2\,{\var{2}}^2-1.278 \cdot 10^{-14}\,{\var{1}}^2\,\var{2}+0.25\,{\var{1}}^2+1.667 \cdot 10^{-14}\,\var{1}\,{\var{2}}^3-3.508 \cdot 10^{-14}\,\var{1}\,{\var{2}}^2+8.576 \cdot 10^{-15}\,\var{1}\,\var{2}+0.25\,\var{1}+1.563 \cdot 10^{-15}\,{\var{2}}^3+0.25\,{\var{2}}^2-0.25\,\var{2}-0.25$ \\~~\\$\bd_{\textrm{lag}}(\var{1},\var{2}) = 1.543 \cdot 10^{-14}\,{\var{1}}^3\,{\var{2}}^3-1.885 \cdot 10^{-15}\,{\var{1}}^3\,{\var{2}}^2-7.321 \cdot 10^{-14}\,{\var{1}}^3\,\var{2}+0.25\,{\var{1}}^3+5.067 \cdot 10^{-15}\,{\var{1}}^2\,{\var{2}}^3+8.393 \cdot 10^{-16}\,{\var{1}}^2\,{\var{2}}^2-1.741 \cdot 10^{-14}\,{\var{1}}^2\,\var{2}+2.557 \cdot 10^{-14}\,{\var{1}}^2-3.217 \cdot 10^{-14}\,\var{1}\,{\var{2}}^3+5.336 \cdot 10^{-15}\,\var{1}\,{\var{2}}^2+5.951 \cdot 10^{-14}\,\var{1}\,\var{2}-9.161 \cdot 10^{-14}\,\var{1}+0.25\,{\var{2}}^3+8.504 \cdot 10^{-16}\,{\var{2}}^2+5.868 \cdot 10^{-15}\,\var{2}+1.0$ \\~~\\\noindent \textbf{Monomial form} (basis, numerator and denominator coefficients - entries $<10^{-12}$ removed):$$\left(\begin{array}{ccc}\mathcal{B}_\textrm{mon}(\var{1},\var{2}) & \bN_\textrm{mon} &\bD_\textrm{mon}\end{array}\right) =$$ $$\left(\begin{array}{ccc} {\var{1}}^3\,{\var{2}}^3 & 0 & 0\\ {\var{1}}^3\,{\var{2}}^2 & 0 & 0\\ {\var{1}}^3\,\var{2} & 0 & 0\\ {\var{1}}^3 & 0 & 0.25\\ {\var{1}}^2\,{\var{2}}^3 & 0 & 0\\ {\var{1}}^2\,{\var{2}}^2 & 0 & 0\\ {\var{1}}^2\,\var{2} & 0 & 0\\ {\var{1}}^2 & 0.25 & 0\\ \var{1}\,{\var{2}}^3 & 0 & 0\\ \var{1}\,{\var{2}}^2 & 0 & 0\\ \var{1}\,\var{2} & 0 & 0\\ \var{1} & 0.25 & 0\\ {\var{2}}^3 & 0 & 0.25\\ {\var{2}}^2 & 0.25 & 0\\ \var{2} & -0.25 & 0\\ 1.0 & -0.25 & 1.0 \end{array}\right)$$\noindent The corresponding function is:$$\begin{array}{rcl}\bG_{\textrm{mon}}(\var{1},\var{2}) &=& \dfrac{\bn_{\textrm{mon}}(\var{1},\var{2})}{\bd_{\textrm{mon}}(\var{1},\var{2})}\\ && \\&=& \dfrac{\sum_{\textrm{row}} \bN_\textrm{mon} \odot \mathcal{B}_\textrm{mon}(\var{1},\var{2})}{\sum_{\textrm{row}} \bD_\textrm{mon} \odot\mathcal{B}_\textrm{mon}(\var{1},\var{2})},  \end{array}$$\noindent where,\\$\bn_{\textrm{mon}}(\var{1},\var{2}) = 0.25\,{\var{1}}^2+0.25\,\var{1}+0.25\,{\var{2}}^2-0.25\,\var{2}-0.25$ \\~~\\$\bd_{\textrm{mon}}(\var{1},\var{2}) = 0.25\,{\var{1}}^3+0.25\,{\var{2}}^3+1.0$ \\~~\\\noindent \textbf{KST equivalent decoupling pattern} (Barycentric weights $\bc^{\var{l}}$): $$\begin{array}{rclll}\var{2}&: & \left(\begin{array}{cccc} -0.6319 & -0.7532 & -0.7602 & -0.8426\\ 2.725 & 2.911 & 2.921 & 3.048\\ -2.589 & -2.75 & -2.759 & -2.869\\ 1.0 & 1.0 & 1.0 & 1.0 \end{array}\right)& \textrm{vec}(.) & := \textbf{Bary}(\var{2}) \\\var{1}&: & \left(\begin{array}{c} -0.2764\\ 1.0\\ -0.9414\\ 0.3281 \end{array}\right)& \textrm{vec}(.) \otimes \bone_{k_{2}} & := \textbf{Bary}(\var{1}) \\\end{array}$$~\\ Then, with the above notations, one defines the following univariate vector functions:  $$ \left\{ \begin{array}{rcl}\bPhi_{1}(\var{1}) &:=& \textbf{Bary}(\var{1}) \odot \mathbf{Lag}(\var{1}) \\\bPhi_{2}(\var{2}) &:=& \textbf{Bary}(\var{2}) \odot \mathbf{Lag}(\var{2}) \\\end{array} \right. $$\noindent The corresponding function is:$$\begin{array}{rcl}\bG_{\textrm{kst}}(\var{1},\var{2}) &=& \dfrac{\bn_{\textrm{kst}}(\var{1},\var{2})}{\bd_{\textrm{kst}}(\var{1},\var{2})}\\ && \\ &=& \dfrac{\sum_{\text{rows}} \bw \odot \bPhi_{1}(\var{1}) \odot \cdots \odot\bPhi_{2}(\var{2})}{\sum_{\text{rows}} \bPhi_{1}(\var{1}) \odot \cdots \odot\bPhi_{2}(\var{2})} . \end{array}$$~\\ \noindent \textbf{KST-like univariate functions} (equivalent scaled univariate functions $\bphi_{1,\cdots,2}$): $$\left\{\begin{array}{rcrcl}z_{1} &=&\bphi_{1}(\var{1}) &=& \cfrac{\bn_{1}}{\bd_{1}} \\z_{2} &=&\bphi_{2}(\var{2}) &=& \cfrac{\bn_{2}}{\bd_{2}} \\\end{array} \right. .$$\noindent where, \\ \noindent $\bn_{1}=-4.608 \cdot 10^{-15}\,{\var{1}}^3+0.2\,{\var{1}}^2+0.2\,\var{1}-0.2$ and \\ \noindent $\bd_{1}=0.2\,{\var{1}}^3+1.125 \cdot 10^{-14}\,{\var{1}}^2-4.716 \cdot 10^{-14}\,\var{1}+1.0$, \\ \noindent $\bn_{2}=9.149 \cdot 10^{-16}\,{\var{2}}^3+0.3333\,{\var{2}}^2-0.3333\,\var{2}-0.3333$ and \\ \noindent $\bd_{2}=0.3333\,{\var{2}}^3-8.522 \cdot 10^{-16}\,{\var{2}}^2+2.313 \cdot 10^{-15}\,\var{2}+1.0$, \\ \noindent \textbf{Connection with Neural Networks} (equivalent numerator $\bn_{\textrm{lag}}$ representation):\\ \begin{figure}[H]\begin{center} \scalebox{.7}{\begin{tikzpicture}[line width=0.4mm]\tikzstyle{place}=[circle, draw=black, minimum size = 8mm]\tikzstyle{placeInOut}=[circle, draw=orange, minimum size = 8mm]\node at (0,-2) [placeInOut] (first_1){$\var{1}$};\node at (0,-4) [placeInOut] (first_2){$\var{2}$};\node at (5,-2) [place] (secondL1_1){$\frac{1}{\var{1}-\lani{1}{1}}$};\node at (5,-4) [place] (secondL1_2){$\frac{1}{\var{1}-\lani{1}{2}}$};\node at (5,-6) [place] (secondL1_3){$\frac{1}{\var{1}-\lani{1}{3}}$};\node at (5,-8) [place] (secondL1_4){$\frac{1}{\var{1}-\lani{1}{4}}$};\node at (5,-10) [place] (secondL2_1){$\frac{1}{\var{2}-\lani{2}{1}}$};\node at (5,-12) [place] (secondL2_2){$\frac{1}{\var{2}-\lani{2}{2}}$};\node at (5,-14) [place] (secondL2_3){$\frac{1}{\var{2}-\lani{2}{3}}$};\node at (5,-16) [place] (secondL2_4){$\frac{1}{\var{2}-\lani{2}{4}}$};\node at (10,-2) [place] (third_1){$\prod$};\node at (10,-4) [place] (third_2){$\prod$};\node at (10,-6) [place] (third_3){$\prod$};\node at (10,-8) [place] (third_4){$\prod$};\node at (10,-10) [place] (third_5){$\prod$};\node at (10,-12) [place] (third_6){$\prod$};\node at (10,-14) [place] (third_7){$\prod$};\node at (10,-16) [place] (third_8){$\prod$};\node at (10,-18) [place] (third_9){$\prod$};\node at (10,-20) [place] (third_10){$\prod$};\node at (10,-22) [place] (third_11){$\prod$};\node at (10,-24) [place] (third_12){$\prod$};\node at (10,-26) [place] (third_13){$\prod$};\node at (10,-28) [place] (third_14){$\prod$};\node at (10,-30) [place] (third_15){$\prod$};\node at (10,-32) [place] (third_16){$\prod$};\node at (15,-17) [placeInOut] (output){$\bSigma$};\draw[->] (first_1)--(secondL1_1) node[above,sloped,pos=0.75] { };\draw[->] (first_1)--(secondL1_2) node[above,sloped,pos=0.75] { };\draw[->] (first_1)--(secondL1_3) node[above,sloped,pos=0.75] { };\draw[->] (first_1)--(secondL1_4) node[above,sloped,pos=0.75] { };\draw[->] (first_2)--(secondL2_1) node[above,sloped,pos=0.75] { };\draw[->] (first_2)--(secondL2_2) node[above,sloped,pos=0.75] { };\draw[->] (first_2)--(secondL2_3) node[above,sloped,pos=0.75] { };\draw[->] (first_2)--(secondL2_4) node[above,sloped,pos=0.75] { };\draw[->] (secondL1_1)--(third_1) node[above,sloped,pos=0.25] {};\draw[->] (secondL1_1)--(third_2) node[above,sloped,pos=0.25] {};\draw[->] (secondL1_1)--(third_3) node[above,sloped,pos=0.25] {};\draw[->] (secondL1_1)--(third_4) node[above,sloped,pos=0.25] {};\draw[->] (secondL1_2)--(third_5) node[above,sloped,pos=0.25] {};\draw[->] (secondL1_2)--(third_6) node[above,sloped,pos=0.25] {};\draw[->] (secondL1_2)--(third_7) node[above,sloped,pos=0.25] {};\draw[->] (secondL1_2)--(third_8) node[above,sloped,pos=0.25] {};\draw[->] (secondL1_3)--(third_9) node[above,sloped,pos=0.25] {};\draw[->] (secondL1_3)--(third_10) node[above,sloped,pos=0.25] {};\draw[->] (secondL1_3)--(third_11) node[above,sloped,pos=0.25] {};\draw[->] (secondL1_3)--(third_12) node[above,sloped,pos=0.25] {};\draw[->] (secondL1_4)--(third_13) node[above,sloped,pos=0.25] {};\draw[->] (secondL1_4)--(third_14) node[above,sloped,pos=0.25] {};\draw[->] (secondL1_4)--(third_15) node[above,sloped,pos=0.25] {};\draw[->] (secondL1_4)--(third_16) node[above,sloped,pos=0.25] {};\draw[->] (secondL2_1)--(third_1) node[above,sloped,pos=0.25] {};\draw[->] (secondL2_2)--(third_2) node[above,sloped,pos=0.25] {};\draw[->] (secondL2_3)--(third_3) node[above,sloped,pos=0.25] {};\draw[->] (secondL2_4)--(third_4) node[above,sloped,pos=0.25] {};\draw[->] (secondL2_1)--(third_5) node[above,sloped,pos=0.25] {};\draw[->] (secondL2_2)--(third_6) node[above,sloped,pos=0.25] {};\draw[->] (secondL2_3)--(third_7) node[above,sloped,pos=0.25] {};\draw[->] (secondL2_4)--(third_8) node[above,sloped,pos=0.25] {};\draw[->] (secondL2_1)--(third_9) node[above,sloped,pos=0.25] {};\draw[->] (secondL2_2)--(third_10) node[above,sloped,pos=0.25] {};\draw[->] (secondL2_3)--(third_11) node[above,sloped,pos=0.25] {};\draw[->] (secondL2_4)--(third_12) node[above,sloped,pos=0.25] {};\draw[->] (secondL2_1)--(third_13) node[above,sloped,pos=0.25] {};\draw[->] (secondL2_2)--(third_14) node[above,sloped,pos=0.25] {};\draw[->] (secondL2_3)--(third_15) node[above,sloped,pos=0.25] {};\draw[->] (secondL2_4)--(third_16) node[above,sloped,pos=0.25] {};\draw[->] (third_1)--(output) node[above,sloped,pos=0.25] {0.08735};\draw[->] (third_2)--(output) node[above,sloped,pos=0.25] {0.1266};\draw[->] (third_3)--(output) node[above,sloped,pos=0.25] {-0.28307};\draw[->] (third_4)--(output) node[above,sloped,pos=0.25] {0.069112};\draw[->] (third_5)--(output) node[above,sloped,pos=0.25] {-0.19592};\draw[->] (third_6)--(output) node[above,sloped,pos=0.25] {-0.54374};\draw[->] (third_7)--(output) node[above,sloped,pos=0.25] {0.98869};\draw[->] (third_8)--(output) node[above,sloped,pos=0.25] {-0.24903};\draw[->] (third_9)--(output) node[above,sloped,pos=0.25] {0.31591};\draw[->] (third_10)--(output) node[above,sloped,pos=0.25] {0.11327};\draw[->] (third_11)--(output) node[above,sloped,pos=0.25] {-0.55441};\draw[->] (third_12)--(output) node[above,sloped,pos=0.25] {0.12523};\draw[->] (third_13)--(output) node[above,sloped,pos=0.25] {-0.20734};\draw[->] (third_14)--(output) node[above,sloped,pos=0.25] {0.30387};\draw[->] (third_15)--(output) node[above,sloped,pos=0.25] {-0.15122};\draw[->] (third_16)--(output) node[above,sloped,pos=0.25] {0.054682};\end{tikzpicture}} \caption{Equivalent NN representation of the numerator $\bn_{\textrm{lag}}$.}\end{center}\end{figure}

\newpage \subsection{Function \#16 (${\ord=2}$ variables, tensor size: 12.5 \textbf{KB})} $$\frac{\var{1}^3+\var{2}^3}{\var{1}^2+\var{2}^2+3}$$ \subsubsection{Setup and results overview}\begin{itemize}\item Reference: A/al. 2021 (A.5.11), \cite{Austin:2021}\item Domain: $\mathbb{R}$\item Tensor size: 12.5 \textbf{KB} ($40^{2}$ points)\item Bounds: $ \left(\begin{array}{cc} -1 & 1 \end{array}\right) \times \left(\begin{array}{cc} -1 & 1 \end{array}\right)$ \end{itemize} \begin{table}[H] \centering \begin{tabular}{llllll}
$\#$ & Alg. & Parameters & Dim. & CPU [s] & RMSE \\ 
\hline 
$\mathbf{\#16}$ & A/G/P-V 2025 (A1) & $0.1,1$ & $\mathbf{64}$ & $0.013$ & $2.5 \cdot 10^{-16}$ \\ 
 & A/G/P-V 2025 (A2) & $1 \cdot 10^{-15},3$ & $64$ & $0.086$ & $1.1 \cdot 10^{-15}$ \\ 
 & MDSPACK v1.1.0 & $0.01,1$ & $64$ & $\mathbf{0.01}$ & $2.5 \cdot 10^{-16}$ \\ 
 & P/P 2025 & $1,0.95,50,0.01,4,12,9$ & $1.8 \cdot 10^{02}$ & $0.88$ & $3.5 \cdot 10^{-05}$ \\ 
 & C-R/B/G 2023 & $0.001,20$ & $80$ & $0.013$ & $4.6 \cdot 10^{-15}$ \\ 
 & B/G 2025 & $1 \cdot 10^{-06},20,3$ & $80$ & $0.019$ & $\mathbf{1.4 \cdot 10^{-16}}$ \\ 
 & TensorFlow & $$ & $2.6 \cdot 10^{02}$ & $14$ & $0.0017$ \\ 
\hline 
\end{tabular} \caption{Function \#16: best model configuration and performances per methods.} \end{table}\begin{figure}[H] \centering  \includegraphics[width=\textwidth]{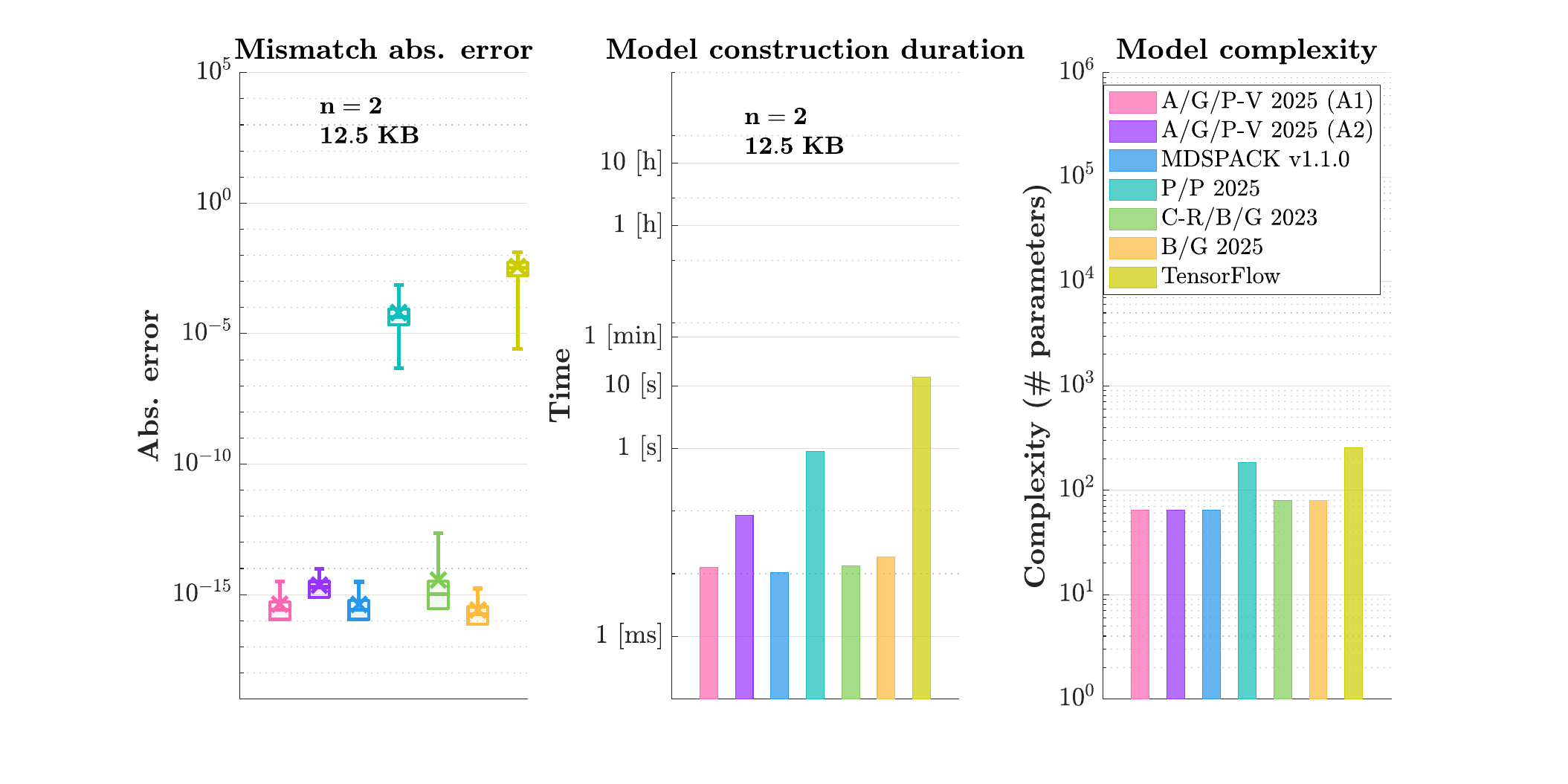} \caption{Function \#16: graphical view of the best model performances.} \end{figure}\begin{figure}[H] \centering  \includegraphics[width=\textwidth]{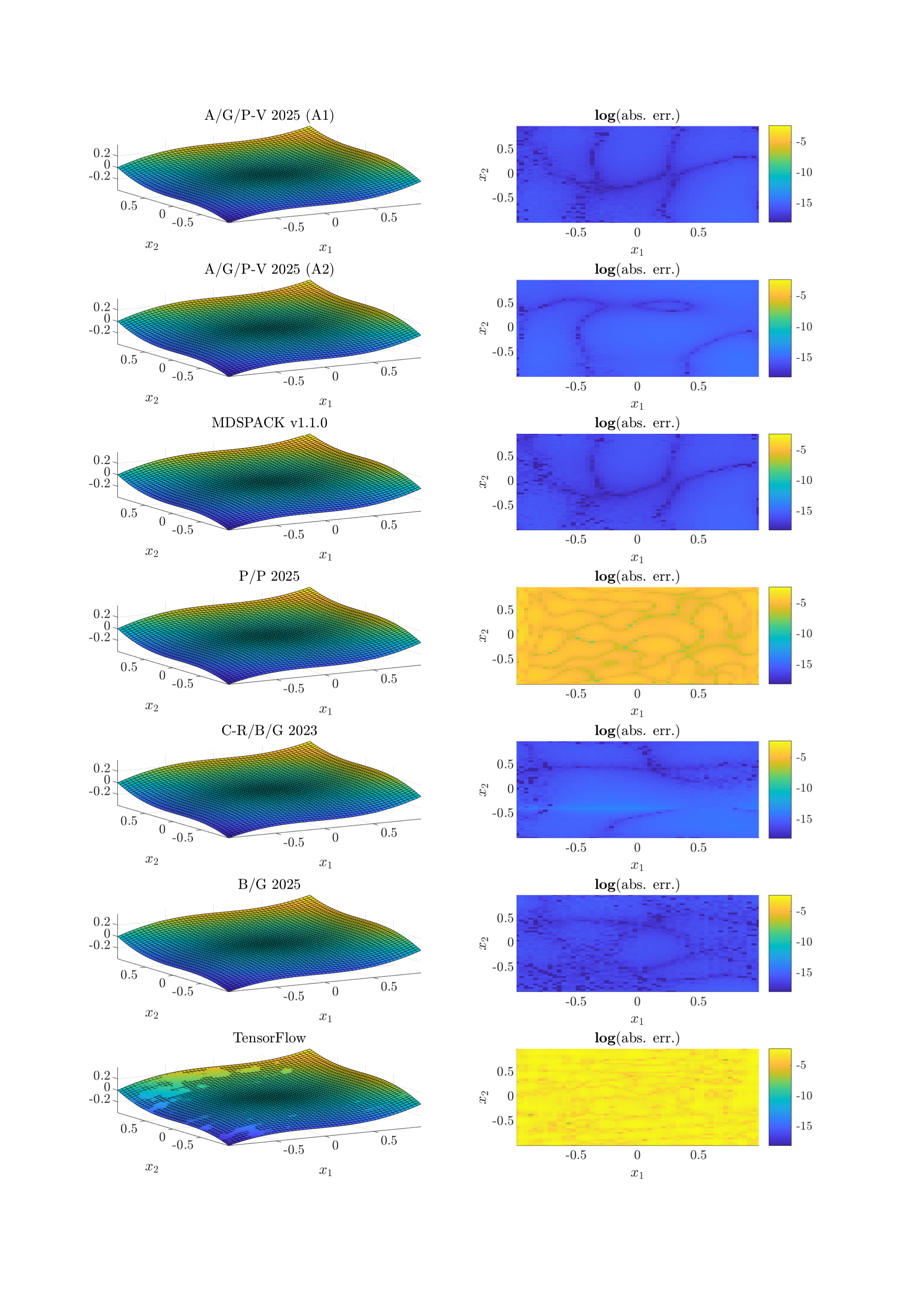} \caption{Function \#16: left side, evaluation of the original (mesh) vs. approximated (coloured surface) and right side, absolute errors (in log-scale).} \end{figure}\subsubsection{mLF detailed informations (M1)} \noindent \textbf{Right interpolation points} ($k_l=\left(\begin{array}{cc} 4 & 4 \end{array}\right)$, where $l=1,\cdots,\ord$):$$ \begin{array}{rcl}\lan{1} &=& \left(\begin{array}{cccc} -1 & -\frac{7}{19} & \frac{5}{19} & 1 \end{array}\right)\\\lan{2} &=& \left(\begin{array}{cccc} -1 & -\frac{7}{19} & \frac{5}{19} & 1 \end{array}\right)\\\end{array} $$\noindent \textbf{Lagrangian weights}: $$\left(\begin{array}{ccc} \bc & \bw & \bc\odot\bw\\ 0.5227 & -0.4 & -0.2091\\ -1.264 & -0.2539 & 0.3209\\ 1.155 & -0.2413 & -0.2786\\ -0.4136 & 0 & 0\\ -1.264 & -0.2539 & 0.3209\\ 2.922 & -0.03057 & -0.08934\\ -2.659 & -0.009917 & 0.02636\\ 1.0 & 0.2297 & 0.2297\\ 1.155 & -0.2413 & -0.2786\\ -2.659 & -0.009917 & 0.02636\\ 2.417 & 0.01161 & 0.02807\\ -0.9136 & 0.2502 & -0.2286\\ -0.4136 & 0 & 0\\ 1.0 & 0.2297 & 0.2297\\ -0.9136 & 0.2502 & -0.2286\\ 0.3272 & 0.4 & 0.1309 \end{array}\right)$$\noindent \textbf{Lagrangian form} (basis, numerator and denominator coefficients):$$\left(\begin{array}{ccc}\mathcal{B}_\textrm{lag}(\var{1},\var{2}) & \bN_\textrm{lag} &\bD_\textrm{lag}\end{array}\right) =$$ $$\left(\begin{array}{ccc} \left(\var{1}+1.0\right)\,\left(\var{2}+1.0\right) & -0.2091 & 0.5227\\ \left(\var{1}+1.0\right)\,\left(\var{2}+0.3684\right) & 0.3209 & -1.264\\ \left(\var{1}+1.0\right)\,\left(\var{2}-0.2632\right) & -0.2786 & 1.155\\ \left(\var{1}+1.0\right)\,\left(\var{2}-1.0\right) & 0 & -0.4136\\ \left(\var{2}+1.0\right)\,\left(\var{1}+0.3684\right) & 0.3209 & -1.264\\ \left(\var{1}+0.3684\right)\,\left(\var{2}+0.3684\right) & -0.08934 & 2.922\\ \left(\var{2}-0.2632\right)\,\left(\var{1}+0.3684\right) & 0.02636 & -2.659\\ \left(\var{2}-1.0\right)\,\left(\var{1}+0.3684\right) & 0.2297 & 1.0\\ \left(\var{2}+1.0\right)\,\left(\var{1}-0.2632\right) & -0.2786 & 1.155\\ \left(\var{1}-0.2632\right)\,\left(\var{2}+0.3684\right) & 0.02636 & -2.659\\ \left(\var{1}-0.2632\right)\,\left(\var{2}-0.2632\right) & 0.02807 & 2.417\\ \left(\var{2}-1.0\right)\,\left(\var{1}-0.2632\right) & -0.2286 & -0.9136\\ \left(\var{1}-1.0\right)\,\left(\var{2}+1.0\right) & 0 & -0.4136\\ \left(\var{1}-1.0\right)\,\left(\var{2}+0.3684\right) & 0.2297 & 1.0\\ \left(\var{1}-1.0\right)\,\left(\var{2}-0.2632\right) & -0.2286 & -0.9136\\ \left(\var{1}-1.0\right)\,\left(\var{2}-1.0\right) & 0.1309 & 0.3272 \end{array}\right).$$\noindent The corresponding function is:$$\begin{array}{rcl}\bG_{\textrm{lag}}(\var{1},\var{2}) &=& \dfrac{\bn_{\textrm{lag}}(\var{1},\var{2})}{\bd_{\textrm{lag}}(\var{1},\var{2})}\\ && \\&=& \dfrac{\sum_{\textrm{row}} \bN_\textrm{lag} \odot\mathcal{B}^{-1}_\textrm{lag}(\var{1},\var{2})}{\sum_{\textrm{row}} \bD_\textrm{lag} \odot\mathcal{B}^{-1}_\textrm{lag}(\var{1},\var{2})}, \end{array}$$\noindent where,\\$\bn_{\textrm{lag}}(\var{1},\var{2}) = 8.506 \cdot 10^{-15}\,{\var{1}}^3\,{\var{2}}^3+1.224 \cdot 10^{-15}\,{\var{1}}^3\,{\var{2}}^2-2.329 \cdot 10^{-15}\,{\var{1}}^3\,\var{2}+0.3333\,{\var{1}}^3+6.906 \cdot 10^{-15}\,{\var{1}}^2\,{\var{2}}^3+1.574 \cdot 10^{-15}\,{\var{1}}^2\,{\var{2}}^2-4.698 \cdot 10^{-15}\,{\var{1}}^2\,\var{2}-1.371 \cdot 10^{-15}\,{\var{1}}^2-5.858 \cdot 10^{-15}\,\var{1}\,{\var{2}}^3-4.555 \cdot 10^{-16}\,\var{1}\,{\var{2}}^2-4.036 \cdot 10^{-17}\,\var{1}\,\var{2}-1.944 \cdot 10^{-16}\,\var{1}+0.3333\,{\var{2}}^3+4.36 \cdot 10^{-17}\,{\var{2}}^2+5.885 \cdot 10^{-16}\,\var{2}+3.233 \cdot 10^{-17}$ \\~~\\$\bd_{\textrm{lag}}(\var{1},\var{2}) = -2.49 \cdot 10^{-14}\,{\var{1}}^3\,{\var{2}}^3-2.29 \cdot 10^{-15}\,{\var{1}}^3\,{\var{2}}^2+1.082 \cdot 10^{-14}\,{\var{1}}^3\,\var{2}+1.637 \cdot 10^{-14}\,{\var{1}}^3-5.1 \cdot 10^{-15}\,{\var{1}}^2\,{\var{2}}^3-5.091 \cdot 10^{-15}\,{\var{1}}^2\,{\var{2}}^2+1.07 \cdot 10^{-15}\,{\var{1}}^2\,\var{2}+0.3333\,{\var{1}}^2+2.613 \cdot 10^{-14}\,\var{1}\,{\var{2}}^3+3.899 \cdot 10^{-15}\,\var{1}\,{\var{2}}^2-1.277 \cdot 10^{-14}\,\var{1}\,\var{2}-1.725 \cdot 10^{-14}\,\var{1}+3.877 \cdot 10^{-15}\,{\var{2}}^3+0.3333\,{\var{2}}^2+8.849 \cdot 10^{-16}\,\var{2}+1.0$ \\~~\\\noindent \textbf{Monomial form} (basis, numerator and denominator coefficients - entries $<10^{-12}$ removed):$$\left(\begin{array}{ccc}\mathcal{B}_\textrm{mon}(\var{1},\var{2}) & \bN_\textrm{mon} &\bD_\textrm{mon}\end{array}\right) =$$ $$\left(\begin{array}{ccc} {\var{1}}^3\,{\var{2}}^3 & 0 & 0\\ {\var{1}}^3\,{\var{2}}^2 & 0 & 0\\ {\var{1}}^3\,\var{2} & 0 & 0\\ {\var{1}}^3 & 0.3333 & 0\\ {\var{1}}^2\,{\var{2}}^3 & 0 & 0\\ {\var{1}}^2\,{\var{2}}^2 & 0 & 0\\ {\var{1}}^2\,\var{2} & 0 & 0\\ {\var{1}}^2 & 0 & 0.3333\\ \var{1}\,{\var{2}}^3 & 0 & 0\\ \var{1}\,{\var{2}}^2 & 0 & 0\\ \var{1}\,\var{2} & 0 & 0\\ \var{1} & 0 & 0\\ {\var{2}}^3 & 0.3333 & 0\\ {\var{2}}^2 & 0 & 0.3333\\ \var{2} & 0 & 0\\ 1.0 & 0 & 1.0 \end{array}\right)$$\noindent The corresponding function is:$$\begin{array}{rcl}\bG_{\textrm{mon}}(\var{1},\var{2}) &=& \dfrac{\bn_{\textrm{mon}}(\var{1},\var{2})}{\bd_{\textrm{mon}}(\var{1},\var{2})}\\ && \\&=& \dfrac{\sum_{\textrm{row}} \bN_\textrm{mon} \odot \mathcal{B}_\textrm{mon}(\var{1},\var{2})}{\sum_{\textrm{row}} \bD_\textrm{mon} \odot\mathcal{B}_\textrm{mon}(\var{1},\var{2})},  \end{array}$$\noindent where,\\$\bn_{\textrm{mon}}(\var{1},\var{2}) = 0.3333\,{\var{1}}^3+0.3333\,{\var{2}}^3$ \\~~\\$\bd_{\textrm{mon}}(\var{1},\var{2}) = 0.3333\,{\var{1}}^2+0.3333\,{\var{2}}^2+1.0$ \\~~\\\noindent \textbf{KST equivalent decoupling pattern} (Barycentric weights $\bc^{\var{l}}$): $$\begin{array}{rclll}\var{2}&: & \left(\begin{array}{cccc} -1.264 & -1.264 & -1.264 & -1.264\\ 3.056 & 2.922 & 2.91 & 3.056\\ -2.792 & -2.659 & -2.646 & -2.792\\ 1.0 & 1.0 & 1.0 & 1.0 \end{array}\right)& \textrm{vec}(.) & := \textbf{Bary}(\var{2}) \\\var{1}&: & \left(\begin{array}{c} -0.4136\\ 1.0\\ -0.9136\\ 0.3272 \end{array}\right)& \textrm{vec}(.) \otimes \bone_{k_{2}} & := \textbf{Bary}(\var{1}) \\\end{array}$$~\\ Then, with the above notations, one defines the following univariate vector functions:  $$ \left\{ \begin{array}{rcl}\bPhi_{1}(\var{1}) &:=& \textbf{Bary}(\var{1}) \odot \mathbf{Lag}(\var{1}) \\\bPhi_{2}(\var{2}) &:=& \textbf{Bary}(\var{2}) \odot \mathbf{Lag}(\var{2}) \\\end{array} \right. $$\noindent The corresponding function is:$$\begin{array}{rcl}\bG_{\textrm{kst}}(\var{1},\var{2}) &=& \dfrac{\bn_{\textrm{kst}}(\var{1},\var{2})}{\bd_{\textrm{kst}}(\var{1},\var{2})}\\ && \\ &=& \dfrac{\sum_{\text{rows}} \bw \odot \bPhi_{1}(\var{1}) \odot \cdots \odot\bPhi_{2}(\var{2})}{\sum_{\text{rows}} \bPhi_{1}(\var{1}) \odot \cdots \odot\bPhi_{2}(\var{2})} . \end{array}$$~\\ \noindent \textbf{KST-like univariate functions} (equivalent scaled univariate functions $\bphi_{1,\cdots,2}$): $$\left\{\begin{array}{rcrcl}z_{1} &=&\bphi_{1}(\var{1}) &=& \cfrac{{\var{1}}^3+1.0}{{\var{1}}^2+4.0}\\z_{2} &=&\bphi_{2}(\var{2}) &=& \cfrac{{\var{2}}^3-1.0}{{\var{2}}^2+4.0}\\\end{array} \right. .$$\noindent \textbf{Connection with Neural Networks} (equivalent numerator $\bn_{\textrm{lag}}$ representation):\\ \begin{figure}[H]\begin{center} \scalebox{.7}{\begin{tikzpicture}[line width=0.4mm]\tikzstyle{place}=[circle, draw=black, minimum size = 8mm]\tikzstyle{placeInOut}=[circle, draw=orange, minimum size = 8mm]\node at (0,-2) [placeInOut] (first_1){$\var{1}$};\node at (0,-4) [placeInOut] (first_2){$\var{2}$};\node at (5,-2) [place] (secondL1_1){$\frac{1}{\var{1}-\lani{1}{1}}$};\node at (5,-4) [place] (secondL1_2){$\frac{1}{\var{1}-\lani{1}{2}}$};\node at (5,-6) [place] (secondL1_3){$\frac{1}{\var{1}-\lani{1}{3}}$};\node at (5,-8) [place] (secondL1_4){$\frac{1}{\var{1}-\lani{1}{4}}$};\node at (5,-10) [place] (secondL2_1){$\frac{1}{\var{2}-\lani{2}{1}}$};\node at (5,-12) [place] (secondL2_2){$\frac{1}{\var{2}-\lani{2}{2}}$};\node at (5,-14) [place] (secondL2_3){$\frac{1}{\var{2}-\lani{2}{3}}$};\node at (5,-16) [place] (secondL2_4){$\frac{1}{\var{2}-\lani{2}{4}}$};\node at (10,-2) [place] (third_1){$\prod$};\node at (10,-4) [place] (third_2){$\prod$};\node at (10,-6) [place] (third_3){$\prod$};\node at (10,-8) [place] (third_4){$\prod$};\node at (10,-10) [place] (third_5){$\prod$};\node at (10,-12) [place] (third_6){$\prod$};\node at (10,-14) [place] (third_7){$\prod$};\node at (10,-16) [place] (third_8){$\prod$};\node at (10,-18) [place] (third_9){$\prod$};\node at (10,-20) [place] (third_10){$\prod$};\node at (10,-22) [place] (third_11){$\prod$};\node at (10,-24) [place] (third_12){$\prod$};\node at (10,-26) [place] (third_13){$\prod$};\node at (10,-28) [place] (third_14){$\prod$};\node at (10,-30) [place] (third_15){$\prod$};\node at (10,-32) [place] (third_16){$\prod$};\node at (15,-17) [placeInOut] (output){$\bSigma$};\draw[->] (first_1)--(secondL1_1) node[above,sloped,pos=0.75] { };\draw[->] (first_1)--(secondL1_2) node[above,sloped,pos=0.75] { };\draw[->] (first_1)--(secondL1_3) node[above,sloped,pos=0.75] { };\draw[->] (first_1)--(secondL1_4) node[above,sloped,pos=0.75] { };\draw[->] (first_2)--(secondL2_1) node[above,sloped,pos=0.75] { };\draw[->] (first_2)--(secondL2_2) node[above,sloped,pos=0.75] { };\draw[->] (first_2)--(secondL2_3) node[above,sloped,pos=0.75] { };\draw[->] (first_2)--(secondL2_4) node[above,sloped,pos=0.75] { };\draw[->] (secondL1_1)--(third_1) node[above,sloped,pos=0.25] {};\draw[->] (secondL1_1)--(third_2) node[above,sloped,pos=0.25] {};\draw[->] (secondL1_1)--(third_3) node[above,sloped,pos=0.25] {};\draw[->] (secondL1_1)--(third_4) node[above,sloped,pos=0.25] {};\draw[->] (secondL1_2)--(third_5) node[above,sloped,pos=0.25] {};\draw[->] (secondL1_2)--(third_6) node[above,sloped,pos=0.25] {};\draw[->] (secondL1_2)--(third_7) node[above,sloped,pos=0.25] {};\draw[->] (secondL1_2)--(third_8) node[above,sloped,pos=0.25] {};\draw[->] (secondL1_3)--(third_9) node[above,sloped,pos=0.25] {};\draw[->] (secondL1_3)--(third_10) node[above,sloped,pos=0.25] {};\draw[->] (secondL1_3)--(third_11) node[above,sloped,pos=0.25] {};\draw[->] (secondL1_3)--(third_12) node[above,sloped,pos=0.25] {};\draw[->] (secondL1_4)--(third_13) node[above,sloped,pos=0.25] {};\draw[->] (secondL1_4)--(third_14) node[above,sloped,pos=0.25] {};\draw[->] (secondL1_4)--(third_15) node[above,sloped,pos=0.25] {};\draw[->] (secondL1_4)--(third_16) node[above,sloped,pos=0.25] {};\draw[->] (secondL2_1)--(third_1) node[above,sloped,pos=0.25] {};\draw[->] (secondL2_2)--(third_2) node[above,sloped,pos=0.25] {};\draw[->] (secondL2_3)--(third_3) node[above,sloped,pos=0.25] {};\draw[->] (secondL2_4)--(third_4) node[above,sloped,pos=0.25] {};\draw[->] (secondL2_1)--(third_5) node[above,sloped,pos=0.25] {};\draw[->] (secondL2_2)--(third_6) node[above,sloped,pos=0.25] {};\draw[->] (secondL2_3)--(third_7) node[above,sloped,pos=0.25] {};\draw[->] (secondL2_4)--(third_8) node[above,sloped,pos=0.25] {};\draw[->] (secondL2_1)--(third_9) node[above,sloped,pos=0.25] {};\draw[->] (secondL2_2)--(third_10) node[above,sloped,pos=0.25] {};\draw[->] (secondL2_3)--(third_11) node[above,sloped,pos=0.25] {};\draw[->] (secondL2_4)--(third_12) node[above,sloped,pos=0.25] {};\draw[->] (secondL2_1)--(third_13) node[above,sloped,pos=0.25] {};\draw[->] (secondL2_2)--(third_14) node[above,sloped,pos=0.25] {};\draw[->] (secondL2_3)--(third_15) node[above,sloped,pos=0.25] {};\draw[->] (secondL2_4)--(third_16) node[above,sloped,pos=0.25] {};\draw[->] (third_1)--(output) node[above,sloped,pos=0.25] {-0.2091};\draw[->] (third_2)--(output) node[above,sloped,pos=0.25] {0.32088};\draw[->] (third_3)--(output) node[above,sloped,pos=0.25] {-0.2786};\draw[->] (third_4)--(output) node[above,sloped,pos=0.25] {0};\draw[->] (third_5)--(output) node[above,sloped,pos=0.25] {0.32088};\draw[->] (third_6)--(output) node[above,sloped,pos=0.25] {-0.089343};\draw[->] (third_7)--(output) node[above,sloped,pos=0.25] {0.026364};\draw[->] (third_8)--(output) node[above,sloped,pos=0.25] {0.2297};\draw[->] (third_9)--(output) node[above,sloped,pos=0.25] {-0.2786};\draw[->] (third_10)--(output) node[above,sloped,pos=0.25] {0.026364};\draw[->] (third_11)--(output) node[above,sloped,pos=0.25] {0.028074};\draw[->] (third_12)--(output) node[above,sloped,pos=0.25] {-0.22862};\draw[->] (third_13)--(output) node[above,sloped,pos=0.25] {0};\draw[->] (third_14)--(output) node[above,sloped,pos=0.25] {0.2297};\draw[->] (third_15)--(output) node[above,sloped,pos=0.25] {-0.22862};\draw[->] (third_16)--(output) node[above,sloped,pos=0.25] {0.1309};\end{tikzpicture}} \caption{Equivalent NN representation of the numerator $\bn_{\textrm{lag}}$.}\end{center}\end{figure}

\newpage \subsection{Function \#17 (${\ord=2}$ variables, tensor size: 12.5 \textbf{KB})} $$\frac{\var{1}^4+\var{2}^4+\var{1}^2\var{2}^2+\var{1}\var{2}}{\var{1}^2\var{2}^2-2\var{1}^2-2\var{2}^2+4}$$ \subsubsection{Setup and results overview}\begin{itemize}\item Reference: A/al. 2021 (A.5.12), \cite{Austin:2021}\item Domain: $\mathbb{R}$\item Tensor size: 12.5 \textbf{KB} ($40^{2}$ points)\item Bounds: $ \left(\begin{array}{cc} -1 & 1 \end{array}\right) \times \left(\begin{array}{cc} -1 & 1 \end{array}\right)$ \end{itemize} \begin{table}[H] \centering \begin{tabular}{llllll}
$\#$ & Alg. & Parameters & Dim. & CPU [s] & RMSE \\ 
\hline 
$\mathbf{\#17}$ & A/G/P-V 2025 (A1) & $0.01,3$ & $\mathbf{1 \cdot 10^{02}}$ & $\mathbf{0.0089}$ & $2.7 \cdot 10^{-15}$ \\ 
 & A/G/P-V 2025 (A2) & $1 \cdot 10^{-15},1$ & $1 \cdot 10^{02}$ & $0.23$ & $6.7 \cdot 10^{-13}$ \\ 
 & MDSPACK v1.1.0 & $0.01,1$ & $1 \cdot 10^{02}$ & $0.014$ & $2.7 \cdot 10^{-15}$ \\ 
 & P/P 2025 & $1,0.95,50,0.01,6,12,13$ & $3.2 \cdot 10^{02}$ & $1.4$ & $0.00056$ \\ 
 & C-R/B/G 2023 & $0.001,20$ & $1 \cdot 10^{02}$ & $0.023$ & $2.1 \cdot 10^{-15}$ \\ 
 & B/G 2025 & $1 \cdot 10^{-06},20,3$ & $1 \cdot 10^{02}$ & $0.031$ & $\mathbf{2.5 \cdot 10^{-16}}$ \\ 
 & TensorFlow & $$ & $2.6 \cdot 10^{02}$ & $14$ & $0.025$ \\ 
\hline 
\end{tabular} \caption{Function \#17: best model configuration and performances per methods.} \end{table}\begin{figure}[H] \centering  \includegraphics[width=\textwidth]{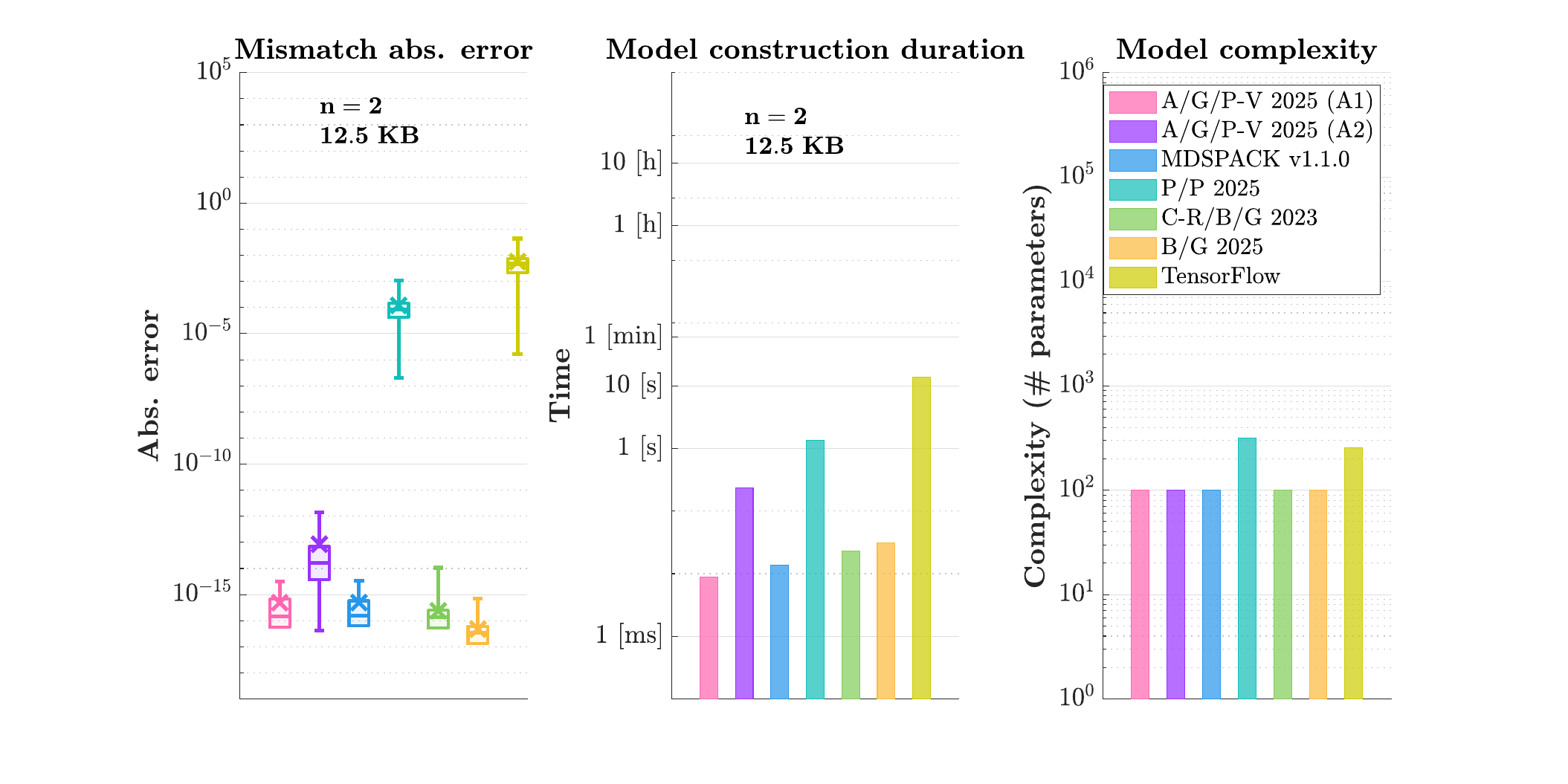} \caption{Function \#17: graphical view of the best model performances.} \end{figure}\begin{figure}[H] \centering  \includegraphics[width=\textwidth]{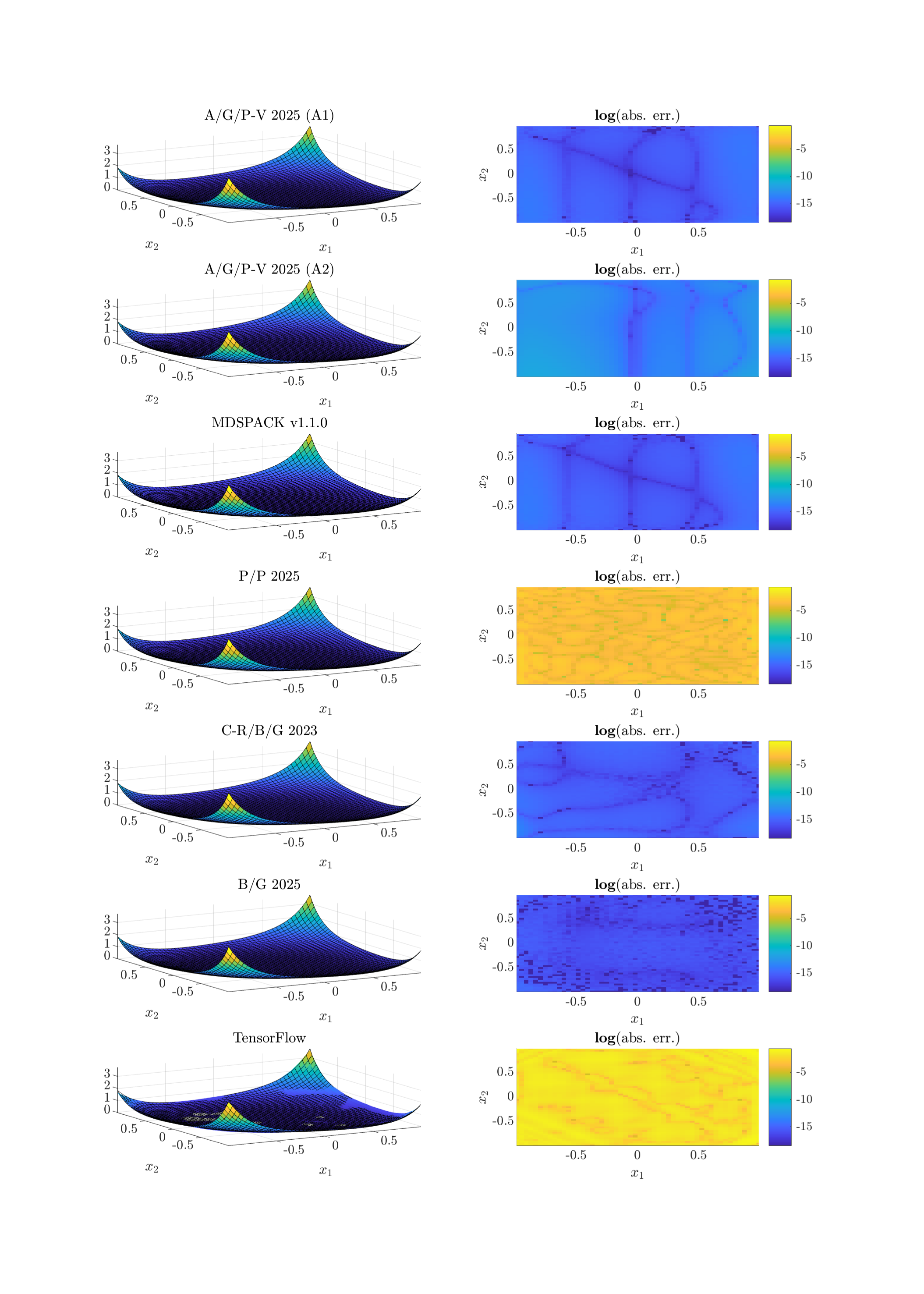} \caption{Function \#17: left side, evaluation of the original (mesh) vs. approximated (coloured surface) and right side, absolute errors (in log-scale).} \end{figure}\subsubsection{mLF detailed informations (M1)} \noindent \textbf{Right interpolation points} ($k_l=\left(\begin{array}{cc} 5 & 5 \end{array}\right)$, where $l=1,\cdots,\ord$):$$ \begin{array}{rcl}\lan{1} &=& \left(\begin{array}{ccccc} -1 & -\frac{11}{19} & -\frac{1}{19} & \frac{9}{19} & 1 \end{array}\right)\\\lan{2} &=& \left(\begin{array}{ccccc} -1 & -\frac{11}{19} & -\frac{1}{19} & \frac{9}{19} & 1 \end{array}\right)\\\end{array} $$\noindent \textbf{Lagrangian weights}: $$\left(\begin{array}{ccc} \bc & \bw & \bc\odot\bw\\ 0.1751 & 4.0 & 0.7003\\ -0.9303 & 1.217 & -1.132\\ 1.488 & 0.5284 & 0.7864\\ -0.8505 & 0.4511 & -0.3837\\ 0.1176 & 2.0 & 0.2353\\ -0.9303 & 1.217 & -1.132\\ 4.944 & 0.2425 & 1.199\\ -7.908 & 0.04323 & -0.3419\\ 4.52 & -0.01229 & -0.05556\\ -0.6252 & 0.5217 & -0.3262\\ 1.488 & 0.5284 & 0.7864\\ -7.908 & 0.04323 & -0.3419\\ 12.65 & 0.0007002 & 0.008857\\ -7.229 & 0.007344 & -0.05309\\ 1.0 & 0.4757 & 0.4757\\ -0.8505 & 0.4511 & -0.3837\\ 4.52 & -0.01229 & -0.05556\\ -7.229 & 0.007344 & -0.05309\\ 4.132 & 0.1191 & 0.492\\ -0.5715 & 0.9847 & -0.5628\\ 0.1176 & 2.0 & 0.2353\\ -0.6252 & 0.5217 & -0.3262\\ 1.0 & 0.4757 & 0.4757\\ -0.5715 & 0.9847 & -0.5628\\ 0.07906 & 4.0 & 0.3162 \end{array}\right)$$\noindent \textbf{Lagrangian form} (basis, numerator and denominator coefficients):$$\left(\begin{array}{ccc}\mathcal{B}_\textrm{lag}(\var{1},\var{2}) & \bN_\textrm{lag} &\bD_\textrm{lag}\end{array}\right) =$$ $$\left(\begin{array}{ccc} \left(\var{1}+1.0\right)\,\left(\var{2}+1.0\right) & 0.7003 & 0.1751\\ \left(\var{1}+1.0\right)\,\left(\var{2}+0.5789\right) & -1.132 & -0.9303\\ \left(\var{1}+1.0\right)\,\left(\var{2}+0.05263\right) & 0.7864 & 1.488\\ \left(\var{1}+1.0\right)\,\left(\var{2}-0.4737\right) & -0.3837 & -0.8505\\ \left(\var{1}+1.0\right)\,\left(\var{2}-1.0\right) & 0.2353 & 0.1176\\ \left(\var{2}+1.0\right)\,\left(\var{1}+0.5789\right) & -1.132 & -0.9303\\ \left(\var{1}+0.5789\right)\,\left(\var{2}+0.5789\right) & 1.199 & 4.944\\ \left(\var{1}+0.5789\right)\,\left(\var{2}+0.05263\right) & -0.3419 & -7.908\\ \left(\var{1}+0.5789\right)\,\left(\var{2}-0.4737\right) & -0.05556 & 4.52\\ \left(\var{2}-1.0\right)\,\left(\var{1}+0.5789\right) & -0.3262 & -0.6252\\ \left(\var{2}+1.0\right)\,\left(\var{1}+0.05263\right) & 0.7864 & 1.488\\ \left(\var{2}+0.5789\right)\,\left(\var{1}+0.05263\right) & -0.3419 & -7.908\\ \left(\var{1}+0.05263\right)\,\left(\var{2}+0.05263\right) & 0.008857 & 12.65\\ \left(\var{1}+0.05263\right)\,\left(\var{2}-0.4737\right) & -0.05309 & -7.229\\ \left(\var{2}-1.0\right)\,\left(\var{1}+0.05263\right) & 0.4757 & 1.0\\ \left(\var{2}+1.0\right)\,\left(\var{1}-0.4737\right) & -0.3837 & -0.8505\\ \left(\var{2}+0.5789\right)\,\left(\var{1}-0.4737\right) & -0.05556 & 4.52\\ \left(\var{2}+0.05263\right)\,\left(\var{1}-0.4737\right) & -0.05309 & -7.229\\ \left(\var{1}-0.4737\right)\,\left(\var{2}-0.4737\right) & 0.492 & 4.132\\ \left(\var{2}-1.0\right)\,\left(\var{1}-0.4737\right) & -0.5628 & -0.5715\\ \left(\var{1}-1.0\right)\,\left(\var{2}+1.0\right) & 0.2353 & 0.1176\\ \left(\var{1}-1.0\right)\,\left(\var{2}+0.5789\right) & -0.3262 & -0.6252\\ \left(\var{1}-1.0\right)\,\left(\var{2}+0.05263\right) & 0.4757 & 1.0\\ \left(\var{1}-1.0\right)\,\left(\var{2}-0.4737\right) & -0.5628 & -0.5715\\ \left(\var{1}-1.0\right)\,\left(\var{2}-1.0\right) & 0.3162 & 0.07906 \end{array}\right).$$\noindent The corresponding function is:$$\begin{array}{rcl}\bG_{\textrm{lag}}(\var{1},\var{2}) &=& \dfrac{\bn_{\textrm{lag}}(\var{1},\var{2})}{\bd_{\textrm{lag}}(\var{1},\var{2})}\\ && \\&=& \dfrac{\sum_{\textrm{row}} \bN_\textrm{lag} \odot\mathcal{B}^{-1}_\textrm{lag}(\var{1},\var{2})}{\sum_{\textrm{row}} \bD_\textrm{lag} \odot\mathcal{B}^{-1}_\textrm{lag}(\var{1},\var{2})}, \end{array}$$\noindent where,\\$\bn_{\textrm{lag}}(\var{1},\var{2}) = 1.063 \cdot 10^{-14}\,{\var{1}}^4\,{\var{2}}^4-4.33 \cdot 10^{-15}\,{\var{1}}^4\,{\var{2}}^3-2.376 \cdot 10^{-14}\,{\var{1}}^4\,{\var{2}}^2+4.481 \cdot 10^{-15}\,{\var{1}}^4\,\var{2}+0.25\,{\var{1}}^4-1.402 \cdot 10^{-14}\,{\var{1}}^3\,{\var{2}}^4-3.66 \cdot 10^{-15}\,{\var{1}}^3\,{\var{2}}^3+1.588 \cdot 10^{-14}\,{\var{1}}^3\,{\var{2}}^2+6.153 \cdot 10^{-15}\,{\var{1}}^3\,\var{2}+4.263 \cdot 10^{-15}\,{\var{1}}^3-9.296 \cdot 10^{-15}\,{\var{1}}^2\,{\var{2}}^4-1.961 \cdot 10^{-15}\,{\var{1}}^2\,{\var{2}}^3+0.25\,{\var{1}}^2\,{\var{2}}^2-5.194 \cdot 10^{-16}\,{\var{1}}^2\,\var{2}-3.146 \cdot 10^{-15}\,{\var{1}}^2+3.902 \cdot 10^{-15}\,\var{1}\,{\var{2}}^4+2.161 \cdot 10^{-15}\,\var{1}\,{\var{2}}^3-5.456 \cdot 10^{-15}\,\var{1}\,{\var{2}}^2+0.25\,\var{1}\,\var{2}-8.249 \cdot 10^{-16}\,\var{1}+0.25\,{\var{2}}^4+1.896 \cdot 10^{-16}\,{\var{2}}^3-2.808 \cdot 10^{-16}\,{\var{2}}^2-1.433 \cdot 10^{-16}\,\var{2}-3.534 \cdot 10^{-17}$ \\~~\\$\bd_{\textrm{lag}}(\var{1},\var{2}) = -1.408 \cdot 10^{-13}\,{\var{1}}^4\,{\var{2}}^4-2.113 \cdot 10^{-14}\,{\var{1}}^4\,{\var{2}}^3+1.148 \cdot 10^{-13}\,{\var{1}}^4\,{\var{2}}^2+1.603 \cdot 10^{-14}\,{\var{1}}^4\,\var{2}+2.601 \cdot 10^{-14}\,{\var{1}}^4+4.001 \cdot 10^{-14}\,{\var{1}}^3\,{\var{2}}^4+5.867 \cdot 10^{-15}\,{\var{1}}^3\,{\var{2}}^3-4.494 \cdot 10^{-14}\,{\var{1}}^3\,{\var{2}}^2-5.769 \cdot 10^{-15}\,{\var{1}}^3\,\var{2}+9.157 \cdot 10^{-15}\,{\var{1}}^3+1.728 \cdot 10^{-13}\,{\var{1}}^2\,{\var{2}}^4+2.366 \cdot 10^{-14}\,{\var{1}}^2\,{\var{2}}^3+0.25\,{\var{1}}^2\,{\var{2}}^2-1.994 \cdot 10^{-14}\,{\var{1}}^2\,\var{2}-0.5\,{\var{1}}^2-4.633 \cdot 10^{-14}\,\var{1}\,{\var{2}}^4-2.391 \cdot 10^{-15}\,\var{1}\,{\var{2}}^3+4.553 \cdot 10^{-14}\,\var{1}\,{\var{2}}^2+3.617 \cdot 10^{-15}\,\var{1}\,\var{2}-1.731 \cdot 10^{-15}\,\var{1}-2.91 \cdot 10^{-15}\,{\var{2}}^4-1.904 \cdot 10^{-16}\,{\var{2}}^3-0.5\,{\var{2}}^2+2.446 \cdot 10^{-16}\,\var{2}+1.0$ \\~~\\\noindent \textbf{Monomial form} (basis, numerator and denominator coefficients - entries $<10^{-12}$ removed):$$\left(\begin{array}{ccc}\mathcal{B}_\textrm{mon}(\var{1},\var{2}) & \bN_\textrm{mon} &\bD_\textrm{mon}\end{array}\right) =$$ $$\left(\begin{array}{ccc} {\var{1}}^4\,{\var{2}}^4 & 0 & 0\\ {\var{1}}^4\,{\var{2}}^3 & 0 & 0\\ {\var{1}}^4\,{\var{2}}^2 & 0 & 0\\ {\var{1}}^4\,\var{2} & 0 & 0\\ {\var{1}}^4 & 0.25 & 0\\ {\var{1}}^3\,{\var{2}}^4 & 0 & 0\\ {\var{1}}^3\,{\var{2}}^3 & 0 & 0\\ {\var{1}}^3\,{\var{2}}^2 & 0 & 0\\ {\var{1}}^3\,\var{2} & 0 & 0\\ {\var{1}}^3 & 0 & 0\\ {\var{1}}^2\,{\var{2}}^4 & 0 & 0\\ {\var{1}}^2\,{\var{2}}^3 & 0 & 0\\ {\var{1}}^2\,{\var{2}}^2 & 0.25 & 0.25\\ {\var{1}}^2\,\var{2} & 0 & 0\\ {\var{1}}^2 & 0 & -0.5\\ \var{1}\,{\var{2}}^4 & 0 & 0\\ \var{1}\,{\var{2}}^3 & 0 & 0\\ \var{1}\,{\var{2}}^2 & 0 & 0\\ \var{1}\,\var{2} & 0.25 & 0\\ \var{1} & 0 & 0\\ {\var{2}}^4 & 0.25 & 0\\ {\var{2}}^3 & 0 & 0\\ {\var{2}}^2 & 0 & -0.5\\ \var{2} & 0 & 0\\ 1.0 & 0 & 1.0 \end{array}\right)$$\noindent The corresponding function is:$$\begin{array}{rcl}\bG_{\textrm{mon}}(\var{1},\var{2}) &=& \dfrac{\bn_{\textrm{mon}}(\var{1},\var{2})}{\bd_{\textrm{mon}}(\var{1},\var{2})}\\ && \\&=& \dfrac{\sum_{\textrm{row}} \bN_\textrm{mon} \odot \mathcal{B}_\textrm{mon}(\var{1},\var{2})}{\sum_{\textrm{row}} \bD_\textrm{mon} \odot\mathcal{B}_\textrm{mon}(\var{1},\var{2})},  \end{array}$$\noindent where,\\$\bn_{\textrm{mon}}(\var{1},\var{2}) = 0.25\,{\var{1}}^4+0.25\,{\var{1}}^2\,{\var{2}}^2+0.25\,\var{1}\,\var{2}+0.25\,{\var{2}}^4$ \\~~\\$\bd_{\textrm{mon}}(\var{1},\var{2}) = 0.25\,{\var{1}}^2\,{\var{2}}^2-0.5\,{\var{1}}^2-0.5\,{\var{2}}^2+1.0$ \\~~\\\noindent \textbf{KST equivalent decoupling pattern} (Barycentric weights $\bc^{\var{l}}$): $$\begin{array}{rclll}\var{2}&: & \left(\begin{array}{ccccc} 1.488 & 1.488 & 1.488 & 1.488 & 1.488\\ -7.908 & -7.908 & -7.908 & -7.908 & -7.908\\ 12.65 & 12.65 & 12.65 & 12.65 & 12.65\\ -7.229 & -7.229 & -7.229 & -7.229 & -7.229\\ 1.0 & 1.0 & 1.0 & 1.0 & 1.0 \end{array}\right)& \textrm{vec}(.) & := \textbf{Bary}(\var{2}) \\\var{1}&: & \left(\begin{array}{c} 0.1176\\ -0.6252\\ 1.0\\ -0.5715\\ 0.07906 \end{array}\right)& \textrm{vec}(.) \otimes \bone_{k_{2}} & := \textbf{Bary}(\var{1}) \\\end{array}$$~\\ Then, with the above notations, one defines the following univariate vector functions:  $$ \left\{ \begin{array}{rcl}\bPhi_{1}(\var{1}) &:=& \textbf{Bary}(\var{1}) \odot \mathbf{Lag}(\var{1}) \\\bPhi_{2}(\var{2}) &:=& \textbf{Bary}(\var{2}) \odot \mathbf{Lag}(\var{2}) \\\end{array} \right. $$\noindent The corresponding function is:$$\begin{array}{rcl}\bG_{\textrm{kst}}(\var{1},\var{2}) &=& \dfrac{\bn_{\textrm{kst}}(\var{1},\var{2})}{\bd_{\textrm{kst}}(\var{1},\var{2})}\\ && \\ &=& \dfrac{\sum_{\text{rows}} \bw \odot \bPhi_{1}(\var{1}) \odot \cdots \odot\bPhi_{2}(\var{2})}{\sum_{\text{rows}} \bPhi_{1}(\var{1}) \odot \cdots \odot\bPhi_{2}(\var{2})} . \end{array}$$~\\ \noindent \textbf{KST-like univariate functions} (equivalent scaled univariate functions $\bphi_{1,\cdots,2}$): $$\left\{\begin{array}{rcrcl}z_{1} &=&\bphi_{1}(\var{1}) &=& \cfrac{\bn_{1}}{\bd_{1}} \\z_{2} &=&\bphi_{2}(\var{2}) &=& \cfrac{\bn_{2}}{\bd_{2}} \\\end{array} \right. .$$\noindent where, \\ \noindent $\bn_{1}=0.5\,{\var{1}}^4+1.735 \cdot 10^{-14}\,{\var{1}}^3+0.5\,{\var{1}}^2+0.5\,\var{1}+0.5$ and \\ \noindent $\bd_{1}=-1.025 \cdot 10^{-14}\,{\var{1}}^4+8.633 \cdot 10^{-15}\,{\var{1}}^3-0.5\,{\var{1}}^2-2.607 \cdot 10^{-15}\,\var{1}+1.0$, \\ \noindent $\bn_{2}=0.5\,{\var{2}}^4+4.332 \cdot 10^{-14}\,{\var{2}}^3+0.5\,{\var{2}}^2-0.5\,\var{2}+0.5$ and \\ \noindent $\bd_{2}=1.826 \cdot 10^{-14}\,{\var{2}}^4+2.347 \cdot 10^{-14}\,{\var{2}}^3-0.5\,{\var{2}}^2-2.347 \cdot 10^{-14}\,\var{2}+1.0$, \\

\newpage \subsection{Function \#18 (${\ord=2}$ variables, tensor size: 12.5 \textbf{KB})} $$\frac{\var{1}^3+\var{2}^3}{\var{1}^2\var{2}^2-2\var{1}^2-2\var{2}^2+4}$$ \subsubsection{Setup and results overview}\begin{itemize}\item Reference: A/al. 2021 (A.5.13), \cite{Austin:2021}\item Domain: $\mathbb{R}$\item Tensor size: 12.5 \textbf{KB} ($40^{2}$ points)\item Bounds: $ \left(\begin{array}{cc} -1 & 1 \end{array}\right) \times \left(\begin{array}{cc} -1 & 1 \end{array}\right)$ \end{itemize} \begin{table}[H] \centering \begin{tabular}{llllll}
$\#$ & Alg. & Parameters & Dim. & CPU [s] & RMSE \\ 
\hline 
$\mathbf{\#18}$ & A/G/P-V 2025 (A1) & $0.1,3$ & $\mathbf{64}$ & $\mathbf{0.0089}$ & $4.4 \cdot 10^{-16}$ \\ 
 & A/G/P-V 2025 (A2) & $1 \cdot 10^{-15},2$ & $64$ & $0.097$ & $1.5 \cdot 10^{-13}$ \\ 
 & MDSPACK v1.1.0 & $0.01,1$ & $64$ & $0.017$ & $5.8 \cdot 10^{-16}$ \\ 
 & P/P 2025 & $1,0.95,50,0.01,6,12,13$ & $3.2 \cdot 10^{02}$ & $1.4$ & $0.00031$ \\ 
 & C-R/B/G 2023 & $0.001,20$ & $80$ & $0.015$ & $3.9 \cdot 10^{-14}$ \\ 
 & B/G 2025 & $1 \cdot 10^{-06},20,3$ & $80$ & $0.02$ & $\mathbf{3.2 \cdot 10^{-16}}$ \\ 
 & TensorFlow & $$ & $2.6 \cdot 10^{02}$ & $14$ & $0.0092$ \\ 
\hline 
\end{tabular} \caption{Function \#18: best model configuration and performances per methods.} \end{table}\begin{figure}[H] \centering  \includegraphics[width=\textwidth]{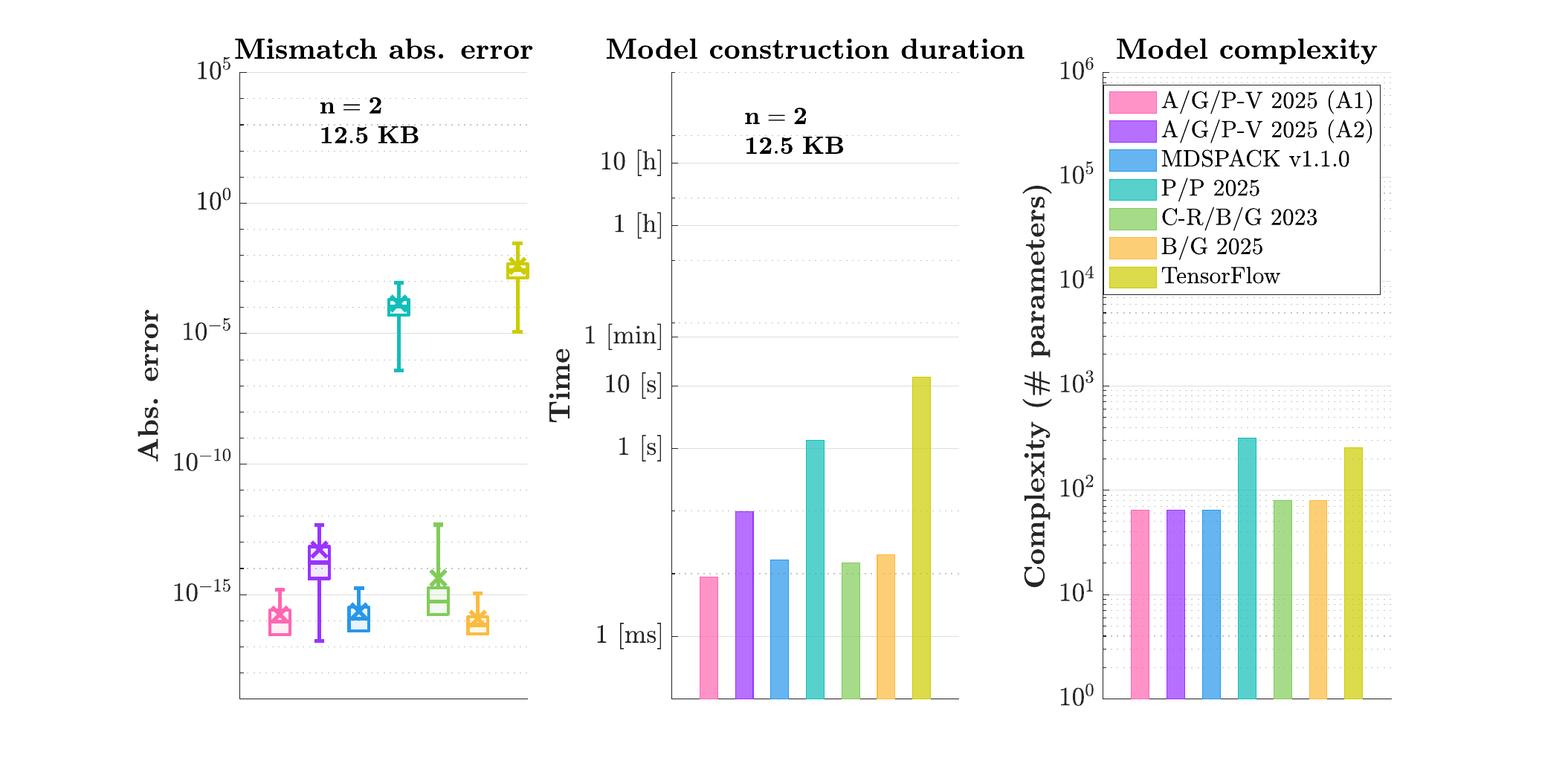} \caption{Function \#18: graphical view of the best model performances.} \end{figure}\begin{figure}[H] \centering  \includegraphics[width=\textwidth]{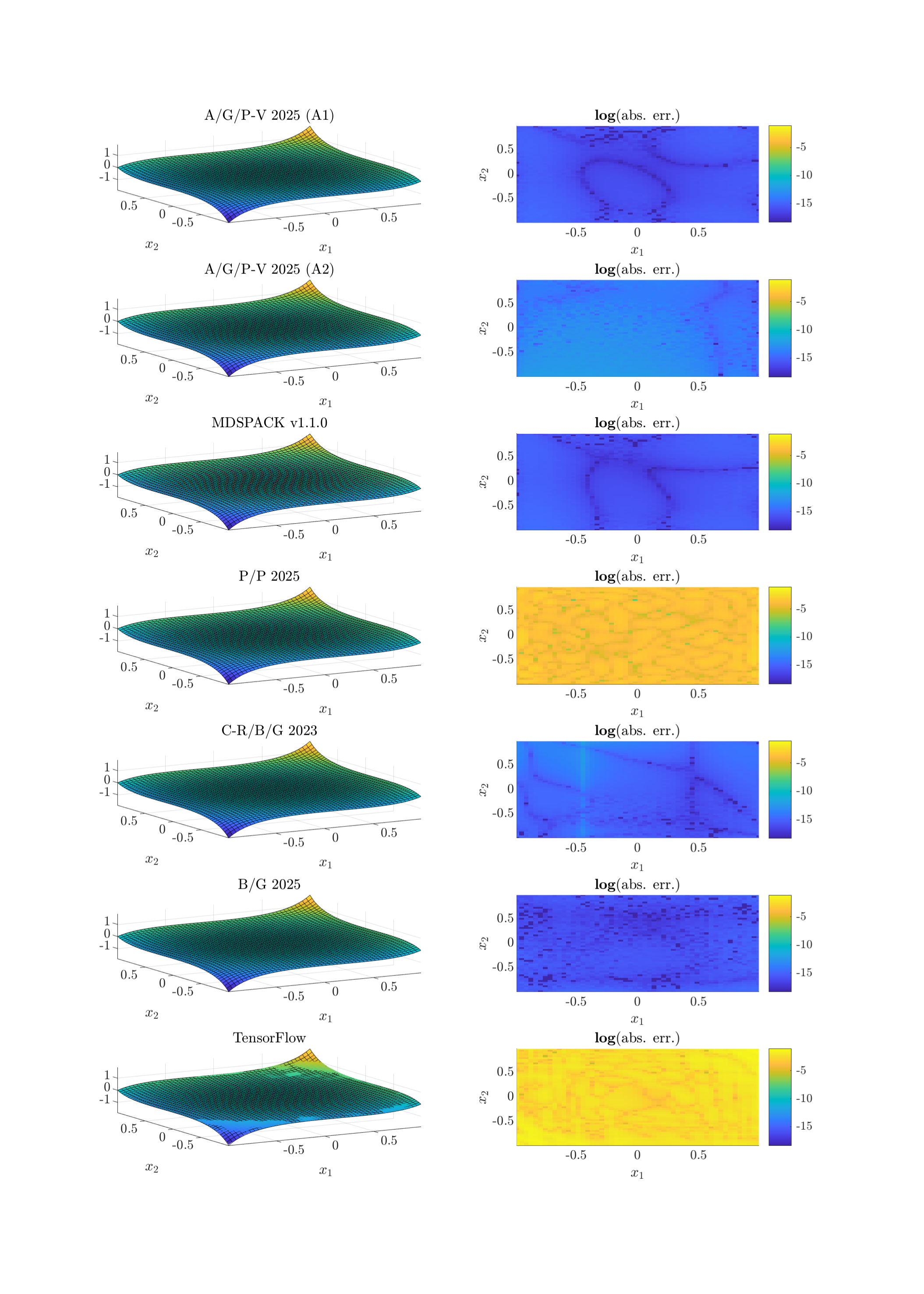} \caption{Function \#18: left side, evaluation of the original (mesh) vs. approximated (coloured surface) and right side, absolute errors (in log-scale).} \end{figure}\subsubsection{mLF detailed informations (M1)} \noindent \textbf{Right interpolation points} ($k_l=\left(\begin{array}{cc} 4 & 4 \end{array}\right)$, where $l=1,\cdots,\ord$):$$ \begin{array}{rcl}\lan{1} &=& \left(\begin{array}{cccc} -1 & -\frac{7}{19} & \frac{5}{19} & 1 \end{array}\right)\\\lan{2} &=& \left(\begin{array}{cccc} -1 & -\frac{7}{19} & \frac{5}{19} & 1 \end{array}\right)\\\end{array} $$\noindent \textbf{Lagrangian weights}: $$\left(\begin{array}{ccc} \bc & \bw & \bc\odot\bw\\ 0.2319 & -2.0 & -0.4639\\ -1.264 & -0.5632 & 0.7119\\ 1.215 & -0.5085 & -0.6181\\ -0.1835 & 0 & 0\\ -1.264 & -0.5632 & 0.7119\\ 6.887 & -0.02878 & -0.1982\\ -6.624 & -0.00883 & 0.05849\\ 1.0 & 0.5096 & 0.5096\\ 1.215 & -0.5085 & -0.6181\\ -6.624 & -0.00883 & 0.05849\\ 6.37 & 0.009778 & 0.06228\\ -0.9617 & 0.5274 & -0.5072\\ -0.1835 & 0 & 0\\ 1.0 & 0.5096 & 0.5096\\ -0.9617 & 0.5274 & -0.5072\\ 0.1452 & 2.0 & 0.2904 \end{array}\right)$$\noindent \textbf{Lagrangian form} (basis, numerator and denominator coefficients):$$\left(\begin{array}{ccc}\mathcal{B}_\textrm{lag}(\var{1},\var{2}) & \bN_\textrm{lag} &\bD_\textrm{lag}\end{array}\right) =$$ $$\left(\begin{array}{ccc} \left(\var{1}+1.0\right)\,\left(\var{2}+1.0\right) & -0.4639 & 0.2319\\ \left(\var{1}+1.0\right)\,\left(\var{2}+0.3684\right) & 0.7119 & -1.264\\ \left(\var{1}+1.0\right)\,\left(\var{2}-0.2632\right) & -0.6181 & 1.215\\ \left(\var{1}+1.0\right)\,\left(\var{2}-1.0\right) & 0 & -0.1835\\ \left(\var{2}+1.0\right)\,\left(\var{1}+0.3684\right) & 0.7119 & -1.264\\ \left(\var{1}+0.3684\right)\,\left(\var{2}+0.3684\right) & -0.1982 & 6.887\\ \left(\var{2}-0.2632\right)\,\left(\var{1}+0.3684\right) & 0.05849 & -6.624\\ \left(\var{2}-1.0\right)\,\left(\var{1}+0.3684\right) & 0.5096 & 1.0\\ \left(\var{2}+1.0\right)\,\left(\var{1}-0.2632\right) & -0.6181 & 1.215\\ \left(\var{1}-0.2632\right)\,\left(\var{2}+0.3684\right) & 0.05849 & -6.624\\ \left(\var{1}-0.2632\right)\,\left(\var{2}-0.2632\right) & 0.06228 & 6.37\\ \left(\var{2}-1.0\right)\,\left(\var{1}-0.2632\right) & -0.5072 & -0.9617\\ \left(\var{1}-1.0\right)\,\left(\var{2}+1.0\right) & 0 & -0.1835\\ \left(\var{1}-1.0\right)\,\left(\var{2}+0.3684\right) & 0.5096 & 1.0\\ \left(\var{1}-1.0\right)\,\left(\var{2}-0.2632\right) & -0.5072 & -0.9617\\ \left(\var{1}-1.0\right)\,\left(\var{2}-1.0\right) & 0.2904 & 0.1452 \end{array}\right).$$\noindent The corresponding function is:$$\begin{array}{rcl}\bG_{\textrm{lag}}(\var{1},\var{2}) &=& \dfrac{\bn_{\textrm{lag}}(\var{1},\var{2})}{\bd_{\textrm{lag}}(\var{1},\var{2})}\\ && \\&=& \dfrac{\sum_{\textrm{row}} \bN_\textrm{lag} \odot\mathcal{B}^{-1}_\textrm{lag}(\var{1},\var{2})}{\sum_{\textrm{row}} \bD_\textrm{lag} \odot\mathcal{B}^{-1}_\textrm{lag}(\var{1},\var{2})}, \end{array}$$\noindent where,\\$\bn_{\textrm{lag}}(\var{1},\var{2}) = -1.873 \cdot 10^{-15}\,{\var{1}}^3\,{\var{2}}^3+1.289 \cdot 10^{-15}\,{\var{1}}^3\,{\var{2}}^2+1.502 \cdot 10^{-15}\,{\var{1}}^3\,\var{2}+0.25\,{\var{1}}^3+1.646 \cdot 10^{-15}\,{\var{1}}^2\,{\var{2}}^3-1.281 \cdot 10^{-15}\,{\var{1}}^2\,{\var{2}}^2-1.24 \cdot 10^{-15}\,{\var{1}}^2\,\var{2}+8.755 \cdot 10^{-16}\,{\var{1}}^2+8.27 \cdot 10^{-16}\,\var{1}\,{\var{2}}^3-5.644 \cdot 10^{-16}\,\var{1}\,{\var{2}}^2-4.559 \cdot 10^{-16}\,\var{1}\,\var{2}+1.933 \cdot 10^{-16}\,\var{1}+0.25\,{\var{2}}^3+5.558 \cdot 10^{-16}\,{\var{2}}^2+1.943 \cdot 10^{-16}\,\var{2}-1.505 \cdot 10^{-16}$ \\~~\\$\bd_{\textrm{lag}}(\var{1},\var{2}) = -1.222 \cdot 10^{-14}\,{\var{1}}^3\,{\var{2}}^3-2.907 \cdot 10^{-15}\,{\var{1}}^3\,{\var{2}}^2+1.122 \cdot 10^{-14}\,{\var{1}}^3\,\var{2}+3.915 \cdot 10^{-15}\,{\var{1}}^3-2.194 \cdot 10^{-15}\,{\var{1}}^2\,{\var{2}}^3+0.25\,{\var{1}}^2\,{\var{2}}^2+1.822 \cdot 10^{-15}\,{\var{1}}^2\,\var{2}-0.5\,{\var{1}}^2+1.077 \cdot 10^{-14}\,\var{1}\,{\var{2}}^3+2.369 \cdot 10^{-15}\,\var{1}\,{\var{2}}^2-9.757 \cdot 10^{-15}\,\var{1}\,\var{2}-3.378 \cdot 10^{-15}\,\var{1}+3.653 \cdot 10^{-15}\,{\var{2}}^3-0.5\,{\var{2}}^2-3.281 \cdot 10^{-15}\,\var{2}+1.0$ \\~~\\\noindent \textbf{Monomial form} (basis, numerator and denominator coefficients - entries $<10^{-12}$ removed):$$\left(\begin{array}{ccc}\mathcal{B}_\textrm{mon}(\var{1},\var{2}) & \bN_\textrm{mon} &\bD_\textrm{mon}\end{array}\right) =$$ $$\left(\begin{array}{ccc} {\var{1}}^3\,{\var{2}}^3 & 0 & 0\\ {\var{1}}^3\,{\var{2}}^2 & 0 & 0\\ {\var{1}}^3\,\var{2} & 0 & 0\\ {\var{1}}^3 & 0.25 & 0\\ {\var{1}}^2\,{\var{2}}^3 & 0 & 0\\ {\var{1}}^2\,{\var{2}}^2 & 0 & 0.25\\ {\var{1}}^2\,\var{2} & 0 & 0\\ {\var{1}}^2 & 0 & -0.5\\ \var{1}\,{\var{2}}^3 & 0 & 0\\ \var{1}\,{\var{2}}^2 & 0 & 0\\ \var{1}\,\var{2} & 0 & 0\\ \var{1} & 0 & 0\\ {\var{2}}^3 & 0.25 & 0\\ {\var{2}}^2 & 0 & -0.5\\ \var{2} & 0 & 0\\ 1.0 & 0 & 1.0 \end{array}\right)$$\noindent The corresponding function is:$$\begin{array}{rcl}\bG_{\textrm{mon}}(\var{1},\var{2}) &=& \dfrac{\bn_{\textrm{mon}}(\var{1},\var{2})}{\bd_{\textrm{mon}}(\var{1},\var{2})}\\ && \\&=& \dfrac{\sum_{\textrm{row}} \bN_\textrm{mon} \odot \mathcal{B}_\textrm{mon}(\var{1},\var{2})}{\sum_{\textrm{row}} \bD_\textrm{mon} \odot\mathcal{B}_\textrm{mon}(\var{1},\var{2})},  \end{array}$$\noindent where,\\$\bn_{\textrm{mon}}(\var{1},\var{2}) = 0.25\,{\var{1}}^3+0.25\,{\var{2}}^3$ \\~~\\$\bd_{\textrm{mon}}(\var{1},\var{2}) = 0.25\,\left({\var{1}}^2-2.0\right)\,\left({\var{2}}^2-2.0\right)$ \\~~\\\noindent \textbf{KST equivalent decoupling pattern} (Barycentric weights $\bc^{\var{l}}$): $$\begin{array}{rclll}\var{2}&: & \left(\begin{array}{cccc} -1.264 & -1.264 & -1.264 & -1.264\\ 6.887 & 6.887 & 6.887 & 6.887\\ -6.624 & -6.624 & -6.624 & -6.624\\ 1.0 & 1.0 & 1.0 & 1.0 \end{array}\right)& \textrm{vec}(.) & := \textbf{Bary}(\var{2}) \\\var{1}&: & \left(\begin{array}{c} -0.1835\\ 1.0\\ -0.9617\\ 0.1452 \end{array}\right)& \textrm{vec}(.) \otimes \bone_{k_{2}} & := \textbf{Bary}(\var{1}) \\\end{array}$$~\\ Then, with the above notations, one defines the following univariate vector functions:  $$ \left\{ \begin{array}{rcl}\bPhi_{1}(\var{1}) &:=& \textbf{Bary}(\var{1}) \odot \mathbf{Lag}(\var{1}) \\\bPhi_{2}(\var{2}) &:=& \textbf{Bary}(\var{2}) \odot \mathbf{Lag}(\var{2}) \\\end{array} \right. $$\noindent The corresponding function is:$$\begin{array}{rcl}\bG_{\textrm{kst}}(\var{1},\var{2}) &=& \dfrac{\bn_{\textrm{kst}}(\var{1},\var{2})}{\bd_{\textrm{kst}}(\var{1},\var{2})}\\ && \\ &=& \dfrac{\sum_{\text{rows}} \bw \odot \bPhi_{1}(\var{1}) \odot \cdots \odot\bPhi_{2}(\var{2})}{\sum_{\text{rows}} \bPhi_{1}(\var{1}) \odot \cdots \odot\bPhi_{2}(\var{2})} . \end{array}$$~\\ \noindent \textbf{KST-like univariate functions} (equivalent scaled univariate functions $\bphi_{1,\cdots,2}$): $$\left\{\begin{array}{rcrcl}z_{1} &=&\bphi_{1}(\var{1}) &=& -\cfrac{1.0\,\left({\var{1}}^3+1.0\right)}{{\var{1}}^2-2.0}\\z_{2} &=&\bphi_{2}(\var{2}) &=& -\cfrac{1.0\,\left({\var{2}}^3-1.0\right)}{{\var{2}}^2-2.0}\\\end{array} \right. .$$\noindent \textbf{Connection with Neural Networks} (equivalent numerator $\bn_{\textrm{lag}}$ representation):\\ \begin{figure}[H]\begin{center} \scalebox{.7}{\begin{tikzpicture}[line width=0.4mm]\tikzstyle{place}=[circle, draw=black, minimum size = 8mm]\tikzstyle{placeInOut}=[circle, draw=orange, minimum size = 8mm]\node at (0,-2) [placeInOut] (first_1){$\var{1}$};\node at (0,-4) [placeInOut] (first_2){$\var{2}$};\node at (5,-2) [place] (secondL1_1){$\frac{1}{\var{1}-\lani{1}{1}}$};\node at (5,-4) [place] (secondL1_2){$\frac{1}{\var{1}-\lani{1}{2}}$};\node at (5,-6) [place] (secondL1_3){$\frac{1}{\var{1}-\lani{1}{3}}$};\node at (5,-8) [place] (secondL1_4){$\frac{1}{\var{1}-\lani{1}{4}}$};\node at (5,-10) [place] (secondL2_1){$\frac{1}{\var{2}-\lani{2}{1}}$};\node at (5,-12) [place] (secondL2_2){$\frac{1}{\var{2}-\lani{2}{2}}$};\node at (5,-14) [place] (secondL2_3){$\frac{1}{\var{2}-\lani{2}{3}}$};\node at (5,-16) [place] (secondL2_4){$\frac{1}{\var{2}-\lani{2}{4}}$};\node at (10,-2) [place] (third_1){$\prod$};\node at (10,-4) [place] (third_2){$\prod$};\node at (10,-6) [place] (third_3){$\prod$};\node at (10,-8) [place] (third_4){$\prod$};\node at (10,-10) [place] (third_5){$\prod$};\node at (10,-12) [place] (third_6){$\prod$};\node at (10,-14) [place] (third_7){$\prod$};\node at (10,-16) [place] (third_8){$\prod$};\node at (10,-18) [place] (third_9){$\prod$};\node at (10,-20) [place] (third_10){$\prod$};\node at (10,-22) [place] (third_11){$\prod$};\node at (10,-24) [place] (third_12){$\prod$};\node at (10,-26) [place] (third_13){$\prod$};\node at (10,-28) [place] (third_14){$\prod$};\node at (10,-30) [place] (third_15){$\prod$};\node at (10,-32) [place] (third_16){$\prod$};\node at (15,-17) [placeInOut] (output){$\bSigma$};\draw[->] (first_1)--(secondL1_1) node[above,sloped,pos=0.75] { };\draw[->] (first_1)--(secondL1_2) node[above,sloped,pos=0.75] { };\draw[->] (first_1)--(secondL1_3) node[above,sloped,pos=0.75] { };\draw[->] (first_1)--(secondL1_4) node[above,sloped,pos=0.75] { };\draw[->] (first_2)--(secondL2_1) node[above,sloped,pos=0.75] { };\draw[->] (first_2)--(secondL2_2) node[above,sloped,pos=0.75] { };\draw[->] (first_2)--(secondL2_3) node[above,sloped,pos=0.75] { };\draw[->] (first_2)--(secondL2_4) node[above,sloped,pos=0.75] { };\draw[->] (secondL1_1)--(third_1) node[above,sloped,pos=0.25] {};\draw[->] (secondL1_1)--(third_2) node[above,sloped,pos=0.25] {};\draw[->] (secondL1_1)--(third_3) node[above,sloped,pos=0.25] {};\draw[->] (secondL1_1)--(third_4) node[above,sloped,pos=0.25] {};\draw[->] (secondL1_2)--(third_5) node[above,sloped,pos=0.25] {};\draw[->] (secondL1_2)--(third_6) node[above,sloped,pos=0.25] {};\draw[->] (secondL1_2)--(third_7) node[above,sloped,pos=0.25] {};\draw[->] (secondL1_2)--(third_8) node[above,sloped,pos=0.25] {};\draw[->] (secondL1_3)--(third_9) node[above,sloped,pos=0.25] {};\draw[->] (secondL1_3)--(third_10) node[above,sloped,pos=0.25] {};\draw[->] (secondL1_3)--(third_11) node[above,sloped,pos=0.25] {};\draw[->] (secondL1_3)--(third_12) node[above,sloped,pos=0.25] {};\draw[->] (secondL1_4)--(third_13) node[above,sloped,pos=0.25] {};\draw[->] (secondL1_4)--(third_14) node[above,sloped,pos=0.25] {};\draw[->] (secondL1_4)--(third_15) node[above,sloped,pos=0.25] {};\draw[->] (secondL1_4)--(third_16) node[above,sloped,pos=0.25] {};\draw[->] (secondL2_1)--(third_1) node[above,sloped,pos=0.25] {};\draw[->] (secondL2_2)--(third_2) node[above,sloped,pos=0.25] {};\draw[->] (secondL2_3)--(third_3) node[above,sloped,pos=0.25] {};\draw[->] (secondL2_4)--(third_4) node[above,sloped,pos=0.25] {};\draw[->] (secondL2_1)--(third_5) node[above,sloped,pos=0.25] {};\draw[->] (secondL2_2)--(third_6) node[above,sloped,pos=0.25] {};\draw[->] (secondL2_3)--(third_7) node[above,sloped,pos=0.25] {};\draw[->] (secondL2_4)--(third_8) node[above,sloped,pos=0.25] {};\draw[->] (secondL2_1)--(third_9) node[above,sloped,pos=0.25] {};\draw[->] (secondL2_2)--(third_10) node[above,sloped,pos=0.25] {};\draw[->] (secondL2_3)--(third_11) node[above,sloped,pos=0.25] {};\draw[->] (secondL2_4)--(third_12) node[above,sloped,pos=0.25] {};\draw[->] (secondL2_1)--(third_13) node[above,sloped,pos=0.25] {};\draw[->] (secondL2_2)--(third_14) node[above,sloped,pos=0.25] {};\draw[->] (secondL2_3)--(third_15) node[above,sloped,pos=0.25] {};\draw[->] (secondL2_4)--(third_16) node[above,sloped,pos=0.25] {};\draw[->] (third_1)--(output) node[above,sloped,pos=0.25] {-0.46386};\draw[->] (third_2)--(output) node[above,sloped,pos=0.25] {0.71186};\draw[->] (third_3)--(output) node[above,sloped,pos=0.25] {-0.61806};\draw[->] (third_4)--(output) node[above,sloped,pos=0.25] {0};\draw[->] (third_5)--(output) node[above,sloped,pos=0.25] {0.71186};\draw[->] (third_6)--(output) node[above,sloped,pos=0.25] {-0.1982};\draw[->] (third_7)--(output) node[above,sloped,pos=0.25] {0.058486};\draw[->] (third_8)--(output) node[above,sloped,pos=0.25] {0.50958};\draw[->] (third_9)--(output) node[above,sloped,pos=0.25] {-0.61806};\draw[->] (third_10)--(output) node[above,sloped,pos=0.25] {0.058486};\draw[->] (third_11)--(output) node[above,sloped,pos=0.25] {0.06228};\draw[->] (third_12)--(output) node[above,sloped,pos=0.25] {-0.50717};\draw[->] (third_13)--(output) node[above,sloped,pos=0.25] {0};\draw[->] (third_14)--(output) node[above,sloped,pos=0.25] {0.50958};\draw[->] (third_15)--(output) node[above,sloped,pos=0.25] {-0.50717};\draw[->] (third_16)--(output) node[above,sloped,pos=0.25] {0.29038};\end{tikzpicture}} \caption{Equivalent NN representation of the numerator $\bn_{\textrm{lag}}$.}\end{center}\end{figure}

\newpage \subsection{Function \#19 (${\ord=2}$ variables, tensor size: 12.5 \textbf{KB})} $$\frac{\var{1}^4+\var{2}^4+\var{1}^2\var{2}^2+\var{1}\var{2}}{\var{1}^3+\var{2}^3+4}$$ \subsubsection{Setup and results overview}\begin{itemize}\item Reference: A/al. 2021 (A.5.14), \cite{Austin:2021}\item Domain: $\mathbb{R}$\item Tensor size: 12.5 \textbf{KB} ($40^{2}$ points)\item Bounds: $ \left(\begin{array}{cc} -1 & 1 \end{array}\right) \times \left(\begin{array}{cc} -1 & 1 \end{array}\right)$ \end{itemize} \begin{table}[H] \centering \begin{tabular}{llllll}
$\#$ & Alg. & Parameters & Dim. & CPU [s] & RMSE \\ 
\hline 
$\mathbf{\#19}$ & A/G/P-V 2025 (A1) & $0.01,2$ & $\mathbf{1 \cdot 10^{02}}$ & $\mathbf{0.0095}$ & $3.1 \cdot 10^{-15}$ \\ 
 & A/G/P-V 2025 (A2) & $1 \cdot 10^{-15},1$ & $1 \cdot 10^{02}$ & $0.15$ & $5 \cdot 10^{-12}$ \\ 
 & MDSPACK v1.1.0 & $0.01,1$ & $1 \cdot 10^{02}$ & $0.013$ & $3 \cdot 10^{-15}$ \\ 
 & P/P 2025 & $1,0.95,50,0.01,6,12,13$ & $3.2 \cdot 10^{02}$ & $1.3$ & $0.00014$ \\ 
 & C-R/B/G 2023 & $0.001,20$ & $1.9 \cdot 10^{02}$ & $0.048$ & $1.3 \cdot 10^{-14}$ \\ 
 & B/G 2025 & $1 \cdot 10^{-06},20,3$ & $1 \cdot 10^{02}$ & $0.095$ & $\mathbf{6.8 \cdot 10^{-16}}$ \\ 
 & TensorFlow & $$ & $2.6 \cdot 10^{02}$ & $14$ & $0.0068$ \\ 
\hline 
\end{tabular} \caption{Function \#19: best model configuration and performances per methods.} \end{table}\begin{figure}[H] \centering  \includegraphics[width=\textwidth]{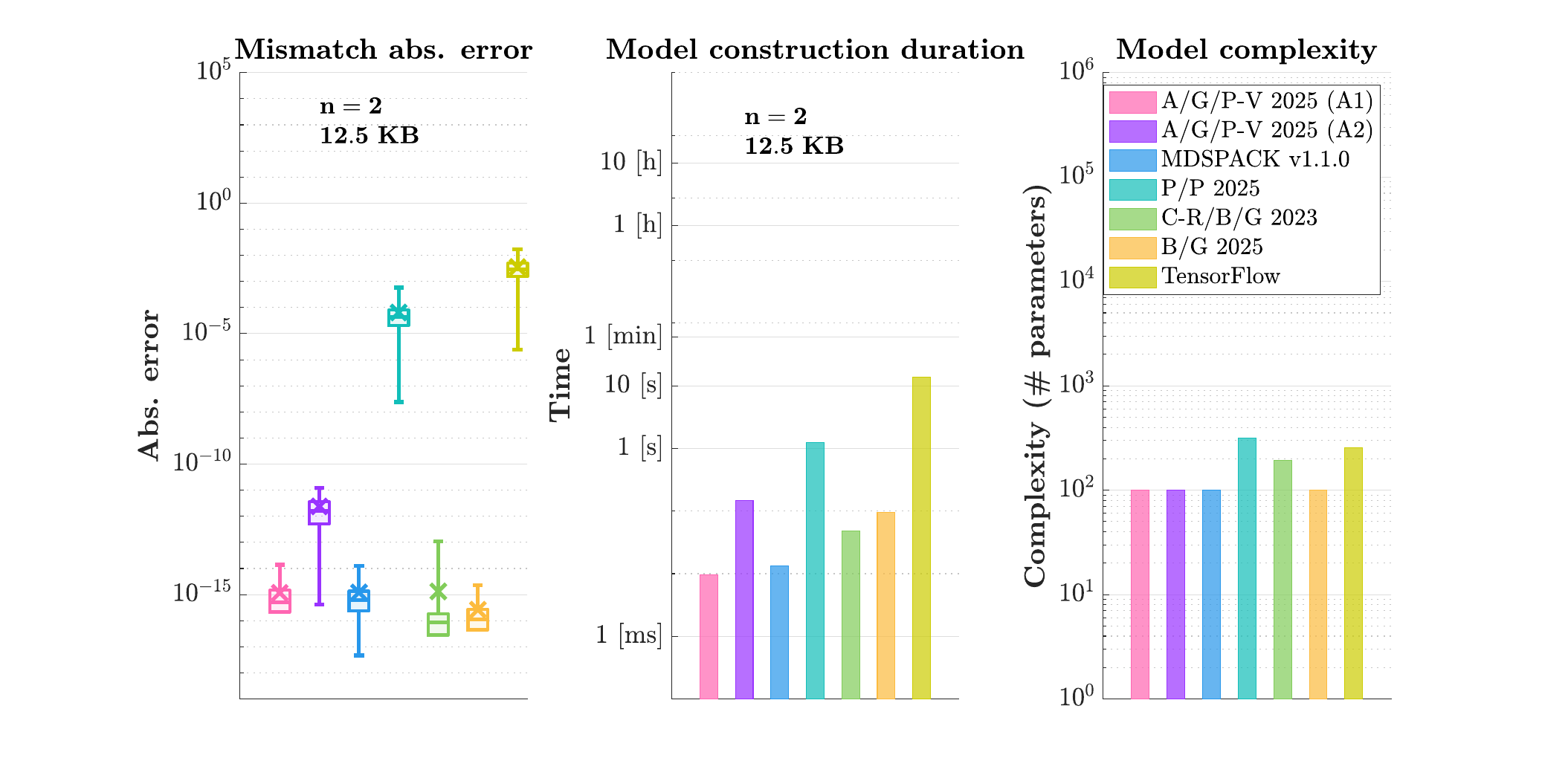} \caption{Function \#19: graphical view of the best model performances.} \end{figure}\begin{figure}[H] \centering  \includegraphics[width=\textwidth]{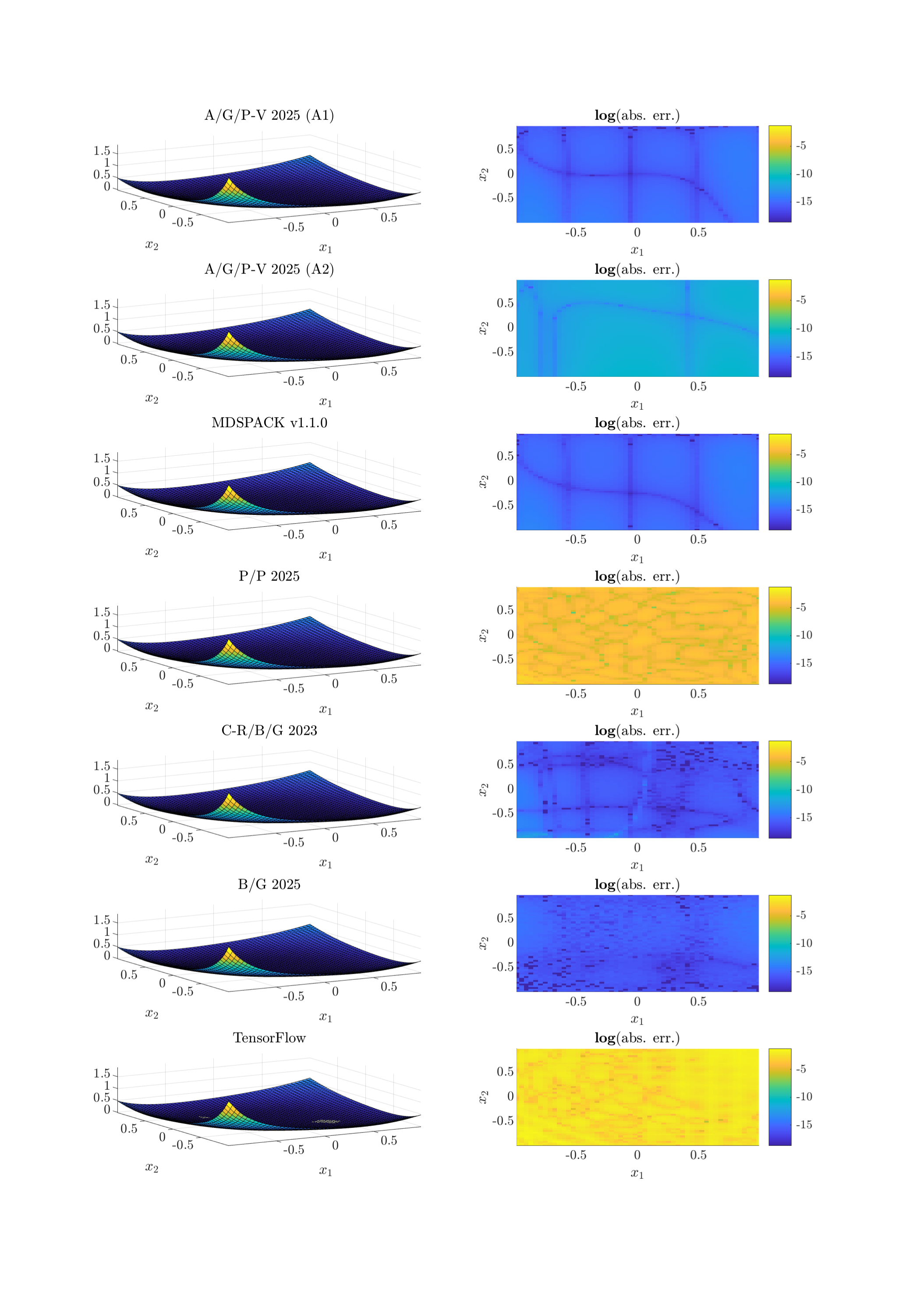} \caption{Function \#19: left side, evaluation of the original (mesh) vs. approximated (coloured surface) and right side, absolute errors (in log-scale).} \end{figure}\subsubsection{mLF detailed informations (M1)} \noindent \textbf{Right interpolation points} ($k_l=\left(\begin{array}{cc} 5 & 5 \end{array}\right)$, where $l=1,\cdots,\ord$):$$ \begin{array}{rcl}\lan{1} &=& \left(\begin{array}{ccccc} -1 & -\frac{11}{19} & -\frac{1}{19} & \frac{9}{19} & 1 \end{array}\right)\\\lan{2} &=& \left(\begin{array}{ccccc} -1 & -\frac{11}{19} & -\frac{1}{19} & \frac{9}{19} & 1 \end{array}\right)\\\end{array} $$\noindent \textbf{Lagrangian weights}: $$\left(\begin{array}{ccc} \bc & \bw & \bc\odot\bw\\ 0.1399 & 2.0 & 0.2797\\ -0.6263 & 0.7222 & -0.4524\\ 0.8928 & 0.3518 & 0.3141\\ -0.5943 & 0.2579 & -0.1533\\ 0.188 & 0.5 & 0.09399\\ -0.6263 & 0.7222 & -0.4524\\ 2.574 & 0.1861 & 0.479\\ -3.616 & 0.03777 & -0.1366\\ 2.389 & -0.009289 & -0.02219\\ -0.7209 & 0.1807 & -0.1303\\ 0.8928 & 0.3518 & 0.3141\\ -3.616 & 0.03777 & -0.1366\\ 5.066 & 0.0006983 & 0.003538\\ -3.344 & 0.006343 & -0.02121\\ 1.0 & 0.19 & 0.19\\ -0.5943 & 0.2579 & -0.1533\\ 2.389 & -0.009289 & -0.02219\\ -3.344 & 0.006343 & -0.02121\\ 2.205 & 0.08912 & 0.1965\\ -0.6565 & 0.3424 & -0.2248\\ 0.188 & 0.5 & 0.09399\\ -0.7209 & 0.1807 & -0.1303\\ 1.0 & 0.19 & 0.19\\ -0.6565 & 0.3424 & -0.2248\\ 0.1895 & 0.6667 & 0.1263 \end{array}\right)$$\noindent \textbf{Lagrangian form} (basis, numerator and denominator coefficients):$$\left(\begin{array}{ccc}\mathcal{B}_\textrm{lag}(\var{1},\var{2}) & \bN_\textrm{lag} &\bD_\textrm{lag}\end{array}\right) =$$ $$\left(\begin{array}{ccc} \left(\var{1}+1.0\right)\,\left(\var{2}+1.0\right) & 0.2797 & 0.1399\\ \left(\var{1}+1.0\right)\,\left(\var{2}+0.5789\right) & -0.4524 & -0.6263\\ \left(\var{1}+1.0\right)\,\left(\var{2}+0.05263\right) & 0.3141 & 0.8928\\ \left(\var{1}+1.0\right)\,\left(\var{2}-0.4737\right) & -0.1533 & -0.5943\\ \left(\var{1}+1.0\right)\,\left(\var{2}-1.0\right) & 0.09399 & 0.188\\ \left(\var{2}+1.0\right)\,\left(\var{1}+0.5789\right) & -0.4524 & -0.6263\\ \left(\var{1}+0.5789\right)\,\left(\var{2}+0.5789\right) & 0.479 & 2.574\\ \left(\var{1}+0.5789\right)\,\left(\var{2}+0.05263\right) & -0.1366 & -3.616\\ \left(\var{1}+0.5789\right)\,\left(\var{2}-0.4737\right) & -0.02219 & 2.389\\ \left(\var{2}-1.0\right)\,\left(\var{1}+0.5789\right) & -0.1303 & -0.7209\\ \left(\var{2}+1.0\right)\,\left(\var{1}+0.05263\right) & 0.3141 & 0.8928\\ \left(\var{2}+0.5789\right)\,\left(\var{1}+0.05263\right) & -0.1366 & -3.616\\ \left(\var{1}+0.05263\right)\,\left(\var{2}+0.05263\right) & 0.003538 & 5.066\\ \left(\var{1}+0.05263\right)\,\left(\var{2}-0.4737\right) & -0.02121 & -3.344\\ \left(\var{2}-1.0\right)\,\left(\var{1}+0.05263\right) & 0.19 & 1.0\\ \left(\var{2}+1.0\right)\,\left(\var{1}-0.4737\right) & -0.1533 & -0.5943\\ \left(\var{2}+0.5789\right)\,\left(\var{1}-0.4737\right) & -0.02219 & 2.389\\ \left(\var{2}+0.05263\right)\,\left(\var{1}-0.4737\right) & -0.02121 & -3.344\\ \left(\var{1}-0.4737\right)\,\left(\var{2}-0.4737\right) & 0.1965 & 2.205\\ \left(\var{2}-1.0\right)\,\left(\var{1}-0.4737\right) & -0.2248 & -0.6565\\ \left(\var{1}-1.0\right)\,\left(\var{2}+1.0\right) & 0.09399 & 0.188\\ \left(\var{1}-1.0\right)\,\left(\var{2}+0.5789\right) & -0.1303 & -0.7209\\ \left(\var{1}-1.0\right)\,\left(\var{2}+0.05263\right) & 0.19 & 1.0\\ \left(\var{1}-1.0\right)\,\left(\var{2}-0.4737\right) & -0.2248 & -0.6565\\ \left(\var{1}-1.0\right)\,\left(\var{2}-1.0\right) & 0.1263 & 0.1895 \end{array}\right).$$\noindent The corresponding function is:$$\begin{array}{rcl}\bG_{\textrm{lag}}(\var{1},\var{2}) &=& \dfrac{\bn_{\textrm{lag}}(\var{1},\var{2})}{\bd_{\textrm{lag}}(\var{1},\var{2})}\\ && \\&=& \dfrac{\sum_{\textrm{row}} \bN_\textrm{lag} \odot\mathcal{B}^{-1}_\textrm{lag}(\var{1},\var{2})}{\sum_{\textrm{row}} \bD_\textrm{lag} \odot\mathcal{B}^{-1}_\textrm{lag}(\var{1},\var{2})}, \end{array}$$\noindent where,\\$\bn_{\textrm{lag}}(\var{1},\var{2}) = 1.797 \cdot 10^{-14}\,{\var{1}}^4\,{\var{2}}^4-5.954 \cdot 10^{-14}\,{\var{1}}^4\,{\var{2}}^3-2.525 \cdot 10^{-14}\,{\var{1}}^4\,{\var{2}}^2+1.319 \cdot 10^{-14}\,{\var{1}}^4\,\var{2}+0.25\,{\var{1}}^4-3.876 \cdot 10^{-15}\,{\var{1}}^3\,{\var{2}}^4-4.595 \cdot 10^{-14}\,{\var{1}}^3\,{\var{2}}^3-5.521 \cdot 10^{-14}\,{\var{1}}^3\,{\var{2}}^2-1.827 \cdot 10^{-15}\,{\var{1}}^3\,\var{2}-9.994 \cdot 10^{-15}\,{\var{1}}^3-6.336 \cdot 10^{-14}\,{\var{1}}^2\,{\var{2}}^4+1.169 \cdot 10^{-14}\,{\var{1}}^2\,{\var{2}}^3+0.25\,{\var{1}}^2\,{\var{2}}^2-3.501 \cdot 10^{-14}\,{\var{1}}^2\,\var{2}+1.215 \cdot 10^{-15}\,{\var{1}}^2-3.83 \cdot 10^{-14}\,\var{1}\,{\var{2}}^4+1.66 \cdot 10^{-14}\,\var{1}\,{\var{2}}^3+5.726 \cdot 10^{-15}\,\var{1}\,{\var{2}}^2+0.25\,\var{1}\,\var{2}-1.232 \cdot 10^{-16}\,\var{1}+0.25\,{\var{2}}^4+1.162 \cdot 10^{-15}\,{\var{2}}^3+5.033 \cdot 10^{-16}\,{\var{2}}^2-1.077 \cdot 10^{-15}\,\var{2}-1.699 \cdot 10^{-17}$ \\~~\\$\bd_{\textrm{lag}}(\var{1},\var{2}) = -5.392 \cdot 10^{-14}\,{\var{1}}^4\,{\var{2}}^4+1.511 \cdot 10^{-14}\,{\var{1}}^4\,{\var{2}}^3-8.307 \cdot 10^{-14}\,{\var{1}}^4\,{\var{2}}^2-5.932 \cdot 10^{-14}\,{\var{1}}^4\,\var{2}+2.48 \cdot 10^{-13}\,{\var{1}}^4+7.801 \cdot 10^{-14}\,{\var{1}}^3\,{\var{2}}^4+4.586 \cdot 10^{-15}\,{\var{1}}^3\,{\var{2}}^3-2.247 \cdot 10^{-13}\,{\var{1}}^3\,{\var{2}}^2-5.023 \cdot 10^{-14}\,{\var{1}}^3\,\var{2}+0.25\,{\var{1}}^3+8.284 \cdot 10^{-14}\,{\var{1}}^2\,{\var{2}}^4-5.202 \cdot 10^{-14}\,{\var{1}}^2\,{\var{2}}^3-4.123 \cdot 10^{-14}\,{\var{1}}^2\,{\var{2}}^2+2.808 \cdot 10^{-15}\,{\var{1}}^2\,\var{2}-3.257 \cdot 10^{-13}\,{\var{1}}^2-4.384 \cdot 10^{-14}\,\var{1}\,{\var{2}}^4-4.84 \cdot 10^{-14}\,\var{1}\,{\var{2}}^3+9.267 \cdot 10^{-14}\,\var{1}\,{\var{2}}^2+2.878 \cdot 10^{-14}\,\var{1}\,\var{2}-2.096 \cdot 10^{-13}\,\var{1}-8.334 \cdot 10^{-15}\,{\var{2}}^4+0.25\,{\var{2}}^3+1.163 \cdot 10^{-14}\,{\var{2}}^2+1.439 \cdot 10^{-15}\,\var{2}+1.0$ \\~~\\\noindent \textbf{Monomial form} (basis, numerator and denominator coefficients - entries $<10^{-12}$ removed):$$\left(\begin{array}{ccc}\mathcal{B}_\textrm{mon}(\var{1},\var{2}) & \bN_\textrm{mon} &\bD_\textrm{mon}\end{array}\right) =$$ $$\left(\begin{array}{ccc} {\var{1}}^4\,{\var{2}}^4 & 0 & 0\\ {\var{1}}^4\,{\var{2}}^3 & 0 & 0\\ {\var{1}}^4\,{\var{2}}^2 & 0 & 0\\ {\var{1}}^4\,\var{2} & 0 & 0\\ {\var{1}}^4 & 0.25 & 0\\ {\var{1}}^3\,{\var{2}}^4 & 0 & 0\\ {\var{1}}^3\,{\var{2}}^3 & 0 & 0\\ {\var{1}}^3\,{\var{2}}^2 & 0 & 0\\ {\var{1}}^3\,\var{2} & 0 & 0\\ {\var{1}}^3 & 0 & 0.25\\ {\var{1}}^2\,{\var{2}}^4 & 0 & 0\\ {\var{1}}^2\,{\var{2}}^3 & 0 & 0\\ {\var{1}}^2\,{\var{2}}^2 & 0.25 & 0\\ {\var{1}}^2\,\var{2} & 0 & 0\\ {\var{1}}^2 & 0 & 0\\ \var{1}\,{\var{2}}^4 & 0 & 0\\ \var{1}\,{\var{2}}^3 & 0 & 0\\ \var{1}\,{\var{2}}^2 & 0 & 0\\ \var{1}\,\var{2} & 0.25 & 0\\ \var{1} & 0 & 0\\ {\var{2}}^4 & 0.25 & 0\\ {\var{2}}^3 & 0 & 0.25\\ {\var{2}}^2 & 0 & 0\\ \var{2} & 0 & 0\\ 1.0 & 0 & 1.0 \end{array}\right)$$\noindent The corresponding function is:$$\begin{array}{rcl}\bG_{\textrm{mon}}(\var{1},\var{2}) &=& \dfrac{\bn_{\textrm{mon}}(\var{1},\var{2})}{\bd_{\textrm{mon}}(\var{1},\var{2})}\\ && \\&=& \dfrac{\sum_{\textrm{row}} \bN_\textrm{mon} \odot \mathcal{B}_\textrm{mon}(\var{1},\var{2})}{\sum_{\textrm{row}} \bD_\textrm{mon} \odot\mathcal{B}_\textrm{mon}(\var{1},\var{2})},  \end{array}$$\noindent where,\\$\bn_{\textrm{mon}}(\var{1},\var{2}) = 0.25\,{\var{1}}^4+0.25\,{\var{1}}^2\,{\var{2}}^2+0.25\,\var{1}\,\var{2}+0.25\,{\var{2}}^4$ \\~~\\$\bd_{\textrm{mon}}(\var{1},\var{2}) = 0.25\,{\var{1}}^3+0.25\,{\var{2}}^3+1.0$ \\~~\\\noindent \textbf{KST equivalent decoupling pattern} (Barycentric weights $\bc^{\var{l}}$): $$\begin{array}{rclll}\var{2}&: & \left(\begin{array}{ccccc} 0.744 & 0.8688 & 0.8928 & 0.9052 & 0.9921\\ -3.332 & -3.57 & -3.616 & -3.639 & -3.805\\ 4.75 & 5.015 & 5.066 & 5.093 & 5.278\\ -3.162 & -3.314 & -3.344 & -3.359 & -3.465\\ 1.0 & 1.0 & 1.0 & 1.0 & 1.0 \end{array}\right)& \textrm{vec}(.) & := \textbf{Bary}(\var{2}) \\\var{1}&: & \left(\begin{array}{c} 0.188\\ -0.7209\\ 1.0\\ -0.6565\\ 0.1895 \end{array}\right)& \textrm{vec}(.) \otimes \bone_{k_{2}} & := \textbf{Bary}(\var{1}) \\\end{array}$$~\\ Then, with the above notations, one defines the following univariate vector functions:  $$ \left\{ \begin{array}{rcl}\bPhi_{1}(\var{1}) &:=& \textbf{Bary}(\var{1}) \odot \mathbf{Lag}(\var{1}) \\\bPhi_{2}(\var{2}) &:=& \textbf{Bary}(\var{2}) \odot \mathbf{Lag}(\var{2}) \\\end{array} \right. $$\noindent The corresponding function is:$$\begin{array}{rcl}\bG_{\textrm{kst}}(\var{1},\var{2}) &=& \dfrac{\bn_{\textrm{kst}}(\var{1},\var{2})}{\bd_{\textrm{kst}}(\var{1},\var{2})}\\ && \\ &=& \dfrac{\sum_{\text{rows}} \bw \odot \bPhi_{1}(\var{1}) \odot \cdots \odot\bPhi_{2}(\var{2})}{\sum_{\text{rows}} \bPhi_{1}(\var{1}) \odot \cdots \odot\bPhi_{2}(\var{2})} . \end{array}$$~\\ \noindent \textbf{KST-like univariate functions} (equivalent scaled univariate functions $\bphi_{1,\cdots,2}$): $$\left\{\begin{array}{rcrcl}z_{1} &=&\bphi_{1}(\var{1}) &=& \cfrac{\bn_{1}}{\bd_{1}} \\z_{2} &=&\bphi_{2}(\var{2}) &=& \cfrac{\bn_{2}}{\bd_{2}} \\\end{array} \right. .$$\noindent where, \\ \noindent $\bn_{1}=0.2\,{\var{1}}^4-9.27 \cdot 10^{-14}\,{\var{1}}^3+0.2\,{\var{1}}^2+0.2\,\var{1}+0.2$ and \\ \noindent $\bd_{1}=5.342 \cdot 10^{-14}\,{\var{1}}^4+0.2\,{\var{1}}^3-2.666 \cdot 10^{-13}\,{\var{1}}^2-1.443 \cdot 10^{-13}\,\var{1}+1.0$, \\ \noindent $\bn_{2}=0.3333\,{\var{2}}^4-2.543 \cdot 10^{-14}\,{\var{2}}^3+0.3333\,{\var{2}}^2-0.3333\,\var{2}+0.3333$ and \\ \noindent $\bd_{2}=-7.577 \cdot 10^{-14}\,{\var{2}}^4+0.3333\,{\var{2}}^3+9.868 \cdot 10^{-14}\,{\var{2}}^2-3.776 \cdot 10^{-14}\,\var{2}+1.0$, \\

\newpage \subsection{Function \#20 (${\ord=3}$ variables, tensor size: 500 \textbf{KB})} $$\texttt{Breit Wigner function}$$ \subsubsection{Setup and results overview}\begin{itemize}\item Reference: A/al. 2021 (A.5.15), \cite{Austin:2021}\item Domain: $\mathbb{R}$\item Tensor size: 500 \textbf{KB} ($40^{3}$ points)\item Bounds: $ \left(\begin{array}{cc} 80 & 100 \end{array}\right) \times \left(\begin{array}{cc} 5 & 10 \end{array}\right) \times \left(\begin{array}{cc} 90 & 93 \end{array}\right)$ \end{itemize} \begin{table}[H] \centering \begin{tabular}{llllll}
$\#$ & Alg. & Parameters & Dim. & CPU [s] & RMSE \\ 
\hline 
$\mathbf{\#20}$ & A/G/P-V 2025 (A1) & $1 \cdot 10^{-06},2$ & $2.4 \cdot 10^{02}$ & $0.029$ & $1.1 \cdot 10^{-06}$ \\ 
 & A/G/P-V 2025 (A2) & $1 \cdot 10^{-15},2$ & $\mathbf{1.4 \cdot 10^{02}}$ & $7.5$ & $3.9 \cdot 10^{-05}$ \\ 
 & MDSPACK v1.1.0 & $1 \cdot 10^{-10},5$ & $2.4 \cdot 10^{02}$ & $\mathbf{0.019}$ & $1.1 \cdot 10^{-06}$ \\ 
 & P/P 2025 & $1,1,50,0.01,6,12,13$ & $3.9 \cdot 10^{02}$ & $16$ & $7.7 \cdot 10^{-05}$ \\ 
 & C-R/B/G 2023 & $1 \cdot 10^{-06},20$ & $1.0 \cdot 10^{03}$ & $5.6$ & $\mathbf{1.5 \cdot 10^{-14}}$ \\ 
 & B/G 2025 & $0.001,20,4$ & $9.6 \cdot 10^{02}$ & $18$ & $1.2 \cdot 10^{-07}$ \\ 
 & TensorFlow & $$ & $3.2 \cdot 10^{02}$ & $2.8 \cdot 10^{02}$ & $0.0033$ \\ 
\hline 
\end{tabular} \caption{Function \#20: best model configuration and performances per methods.} \end{table}\begin{figure}[H] \centering  \includegraphics[width=\textwidth]{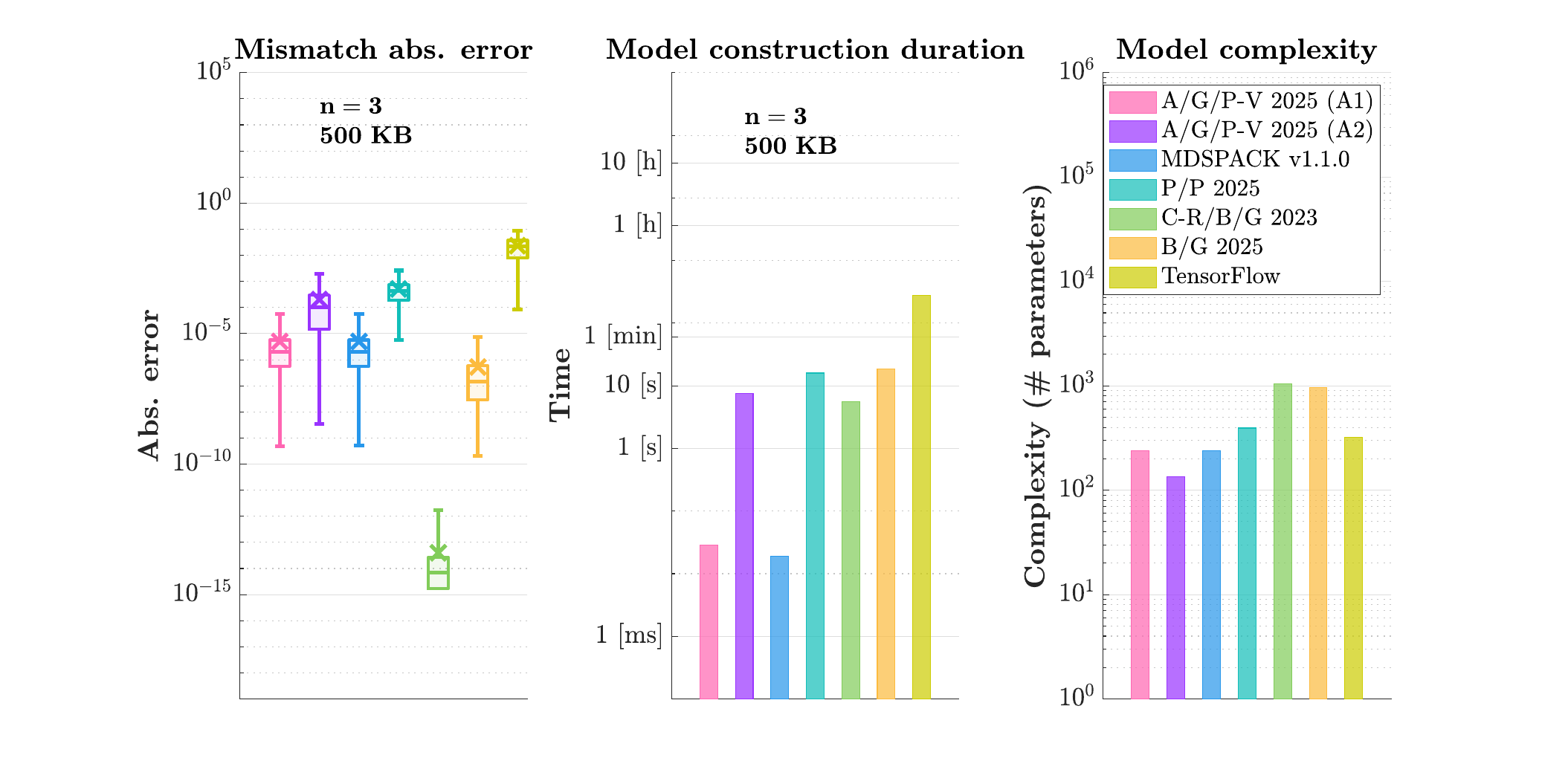} \caption{Function \#20: graphical view of the best model performances.} \end{figure}\begin{figure}[H] \centering  \includegraphics[width=\textwidth]{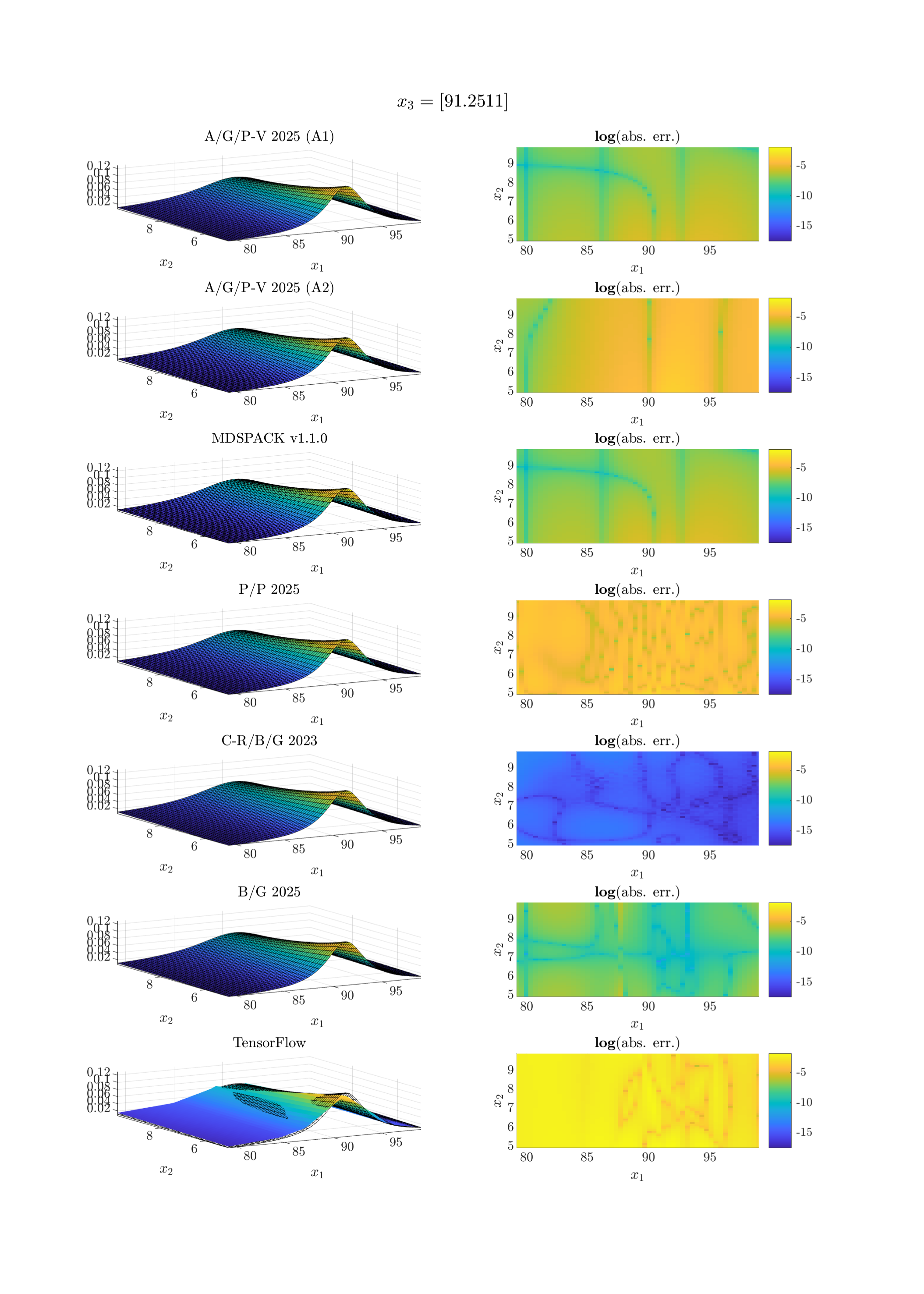} \caption{Function \#20: left side, evaluation of the original (mesh) vs. approximated (coloured surface) and right side, absolute errors (in log-scale).} \end{figure}\subsubsection{mLF detailed informations (M1)} \noindent \textbf{Right interpolation points}: $k_l=\left(\begin{array}{ccc} 4 & 4 & 3 \end{array}\right)$, where $l=1,\cdots,\ord$.$$ \begin{array}{rcl}\lan{1} &\in& \IC^{4} \text{ , linearly spaced between bounds}\\\lan{2} &\in& \IC^{4} \text{ , linearly spaced between bounds}\\\lan{3} &\in& \IC^{3} \text{ , linearly spaced between bounds}\\\end{array} $$\noindent \textbf{$\ord$-D Loewner matrix, barycentric weights and Lagrangian basis}:$$ \begin{array}{rcl}\IL & \in & \IC^{48 \times 48}\\\bc & \in & \IC^{48}\\\bw & \in & \IC^{48}\\\bc\odot \bw & \in & \IC^{48}\\\mathbf{Lag}(\var{1},\var{2},\var{3}) & \in & \IC^{48}\\\end{array} $$

\newpage \subsection{Function \#21 (${\ord=4}$ variables, tensor size: 1.22 \textbf{MB})} $$\frac{\sum_{i=1}^4\mathrm{atan}(x_i)}{\var{1}^2\var{2}^2-\var{1}^2-\var{2}^2+1}$$ \subsubsection{Setup and results overview}\begin{itemize}\item Reference: A/al. 2021 (A.5.16), \cite{Austin:2021}\item Domain: $\mathbb{R}$\item Tensor size: 1.22 \textbf{MB} ($20^{4}$ points)\item Bounds: $ \left(\begin{array}{cc} -\frac{19}{20} & \frac{19}{20} \end{array}\right) \times \left(\begin{array}{cc} -\frac{19}{20} & \frac{19}{20} \end{array}\right) \times \left(\begin{array}{cc} -\frac{19}{20} & \frac{19}{20} \end{array}\right) \times \left(\begin{array}{cc} -\frac{19}{20} & \frac{19}{20} \end{array}\right)$ \end{itemize} \begin{table}[H] \centering \begin{tabular}{llllll}
$\#$ & Alg. & Parameters & Dim. & CPU [s] & RMSE \\ 
\hline 
$\mathbf{\#21}$ & A/G/P-V 2025 (A1) & $1 \cdot 10^{-09},3$ & $2.5 \cdot 10^{04}$ & $0.13$ & $4.1 \cdot 10^{-05}$ \\ 
 & A/G/P-V 2025 (A2) & $1 \cdot 10^{-15},1$ & $1.4 \cdot 10^{04}$ & $88$ & $7.5 \cdot 10^{-05}$ \\ 
 & MDSPACK v1.1.0 & $1 \cdot 10^{-08},4$ & $2.5 \cdot 10^{04}$ & $\mathbf{0.046}$ & $4.1 \cdot 10^{-05}$ \\ 
 & P/P 2025 & $1,1,50,0.01,6,6,13$ & $3.9 \cdot 10^{02}$ & $36$ & $1.8$ \\ 
 & C-R/B/G 2023 & $1 \cdot 10^{-06},20$ & $5.1 \cdot 10^{04}$ & $2.5 \cdot 10^{03}$ & $\mathbf{1.5 \cdot 10^{-05}}$ \\ 
 & B/G 2025 & $1 \cdot 10^{-09},20,4$ & $6.6 \cdot 10^{04}$ & $67$ & $2.9 \cdot 10^{-05}$ \\ 
 & TensorFlow & $$ & $\mathbf{3.8 \cdot 10^{02}}$ & $1.2 \cdot 10^{03}$ & $1.8$ \\ 
\hline 
\end{tabular} \caption{Function \#21: best model configuration and performances per methods.} \end{table}\begin{figure}[H] \centering  \includegraphics[width=\textwidth]{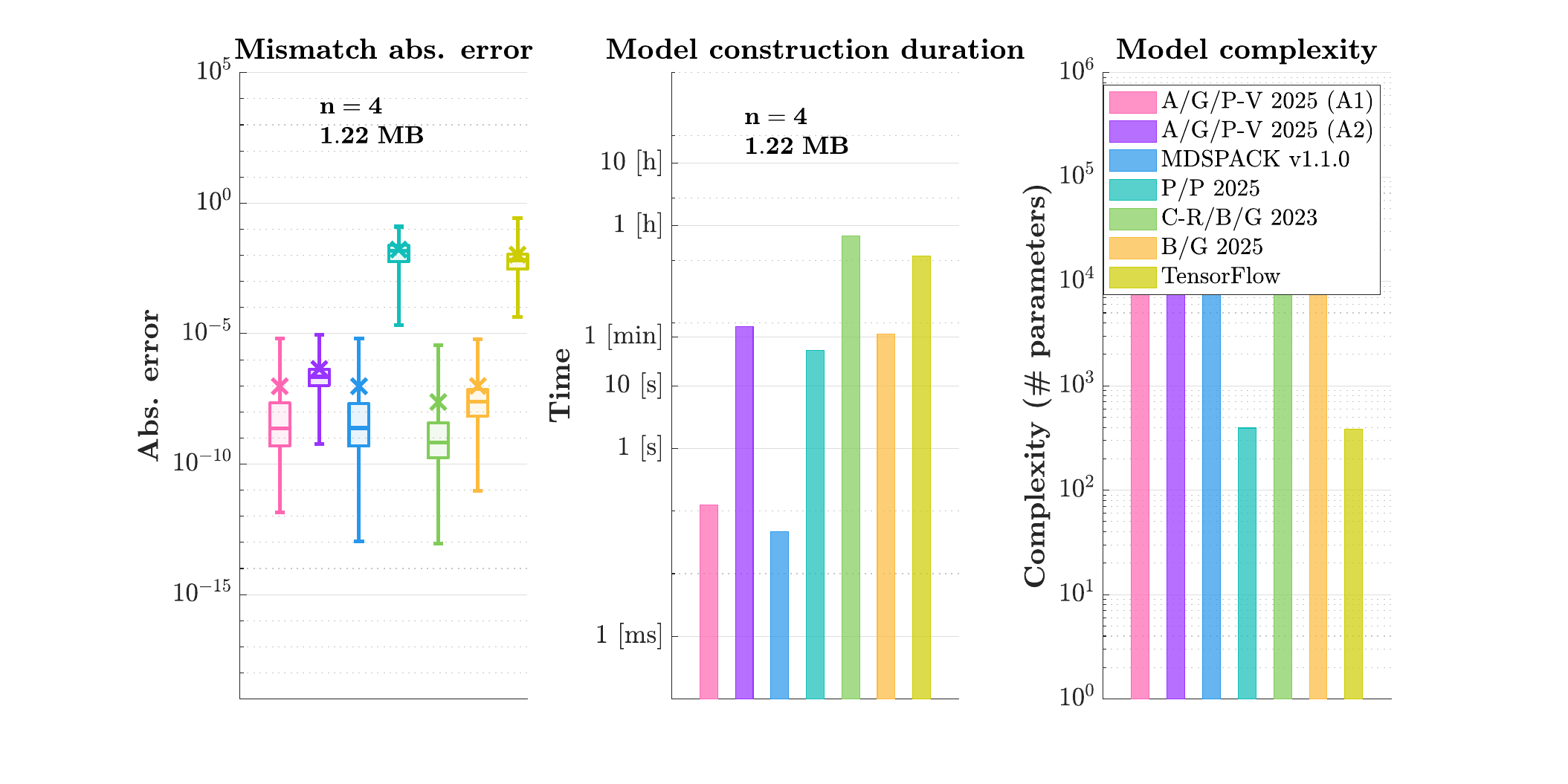} \caption{Function \#21: graphical view of the best model performances.} \end{figure}\begin{figure}[H] \centering  \includegraphics[width=\textwidth]{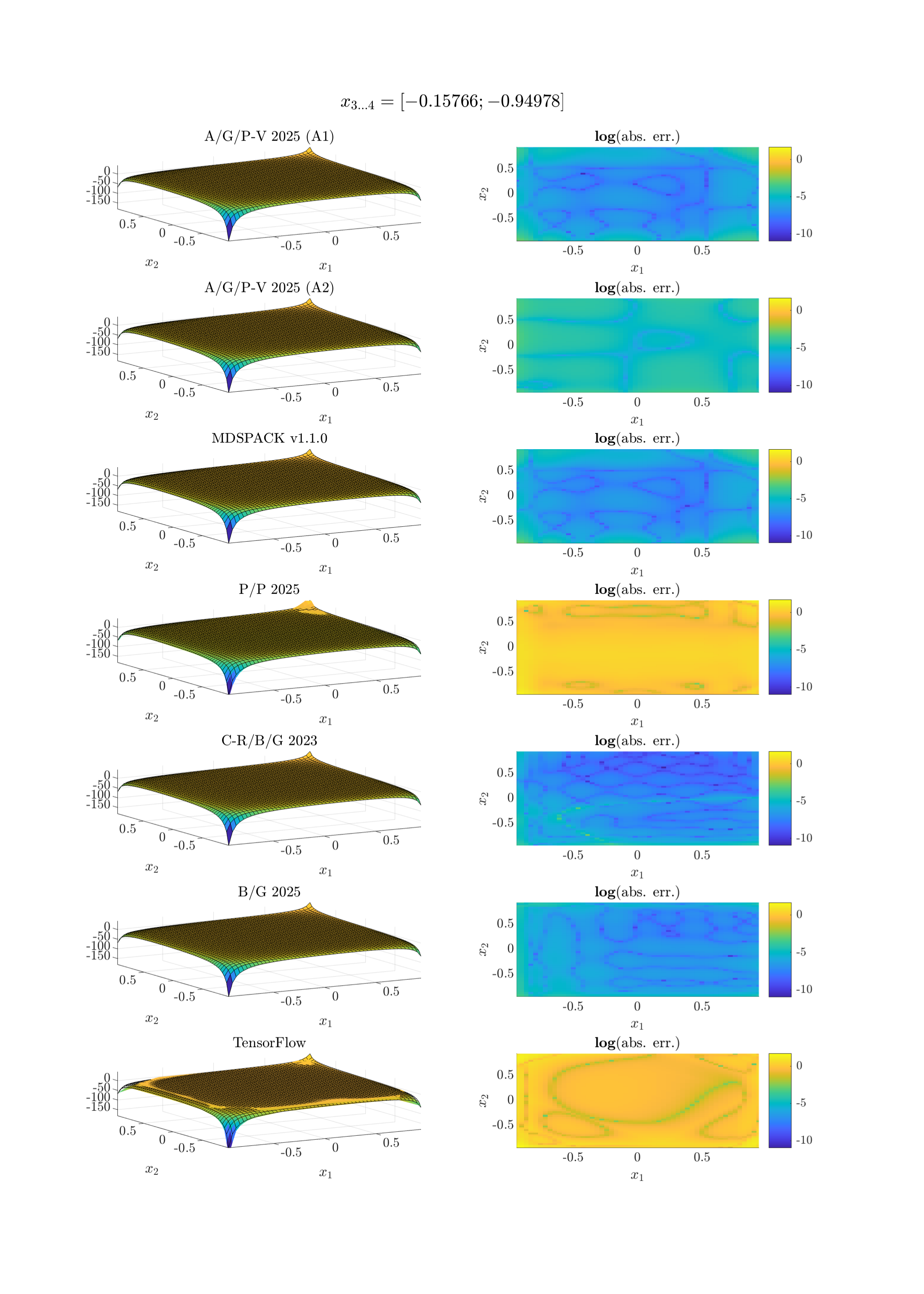} \caption{Function \#21: left side, evaluation of the original (mesh) vs. approximated (coloured surface) and right side, absolute errors (in log-scale).} \end{figure}\subsubsection{mLF detailed informations (M1)} \noindent \textbf{Right interpolation points}: $k_l=\left(\begin{array}{cccc} 8 & 8 & 8 & 8 \end{array}\right)$, where $l=1,\cdots,\ord$.$$ \begin{array}{rcl}\lan{1} &\in& \IC^{8} \text{ , linearly spaced between bounds}\\\lan{2} &\in& \IC^{8} \text{ , linearly spaced between bounds}\\\lan{3} &\in& \IC^{8} \text{ , linearly spaced between bounds}\\\lan{4} &\in& \IC^{8} \text{ , linearly spaced between bounds}\\\end{array} $$\noindent \textbf{$\ord$-D Loewner matrix, barycentric weights and Lagrangian basis}:$$ \begin{array}{rcl}\IL & \in & \IC^{4096 \times 4096}\\\bc & \in & \IC^{4096}\\\bw & \in & \IC^{4096}\\\bc\odot \bw & \in & \IC^{4096}\\\mathbf{Lag}(\var{1},\var{2},\var{3},\var{4}) & \in & \IC^{4096}\\\end{array} $$

\newpage \subsection{Function \#22 (${\ord=4}$ variables, tensor size: 1.22 \textbf{MB})} $$\frac{\mathrm{exp}(\var{1}\var{2}\var{3}\var{4})}{\var{1}^2+\var{2}^2-\var{3}\var{4}+3}$$ \subsubsection{Setup and results overview}\begin{itemize}\item Reference: A/al. 2021 (A.5.17), \cite{Austin:2021}\item Domain: $\mathbb{R}$\item Tensor size: 1.22 \textbf{MB} ($20^{4}$ points)\item Bounds: $ \left(\begin{array}{cc} -1 & 1 \end{array}\right) \times \left(\begin{array}{cc} -1 & 1 \end{array}\right) \times \left(\begin{array}{cc} -1 & 1 \end{array}\right) \times \left(\begin{array}{cc} -1 & 1 \end{array}\right)$ \end{itemize} \begin{table}[H] \centering \begin{tabular}{llllll}
$\#$ & Alg. & Parameters & Dim. & CPU [s] & RMSE \\ 
\hline 
$\mathbf{\#22}$ & A/G/P-V 2025 (A1) & $0.0001,2$ & $1.5 \cdot 10^{03}$ & $0.048$ & $0.0013$ \\ 
 & A/G/P-V 2025 (A2) & $1 \cdot 10^{-15},2$ & $\mathbf{6}$ & $20$ & $0.8$ \\ 
 & MDSPACK v1.1.0 & $1 \cdot 10^{-10},5$ & $1.5 \cdot 10^{03}$ & $\mathbf{0.041}$ & $0.0013$ \\ 
 & P/P 2025 & $1,0.95,50,0.01,4,12,9$ & $2.6 \cdot 10^{02}$ & $27$ & $0.0011$ \\ 
 & C-R/B/G 2023 & $1 \cdot 10^{-09},20$ & $2.4 \cdot 10^{04}$ & $4.3 \cdot 10^{02}$ & $\mathbf{2 \cdot 10^{-13}}$ \\ 
 & B/G 2025 & $1 \cdot 10^{-09},20,3$ & $7 \cdot 10^{04}$ & $63$ & $7.4 \cdot 10^{-07}$ \\ 
 & TensorFlow & $$ & $3.8 \cdot 10^{02}$ & $1.4 \cdot 10^{02}$ & $0.0074$ \\ 
\hline 
\end{tabular} \caption{Function \#22: best model configuration and performances per methods.} \end{table}\begin{figure}[H] \centering  \includegraphics[width=\textwidth]{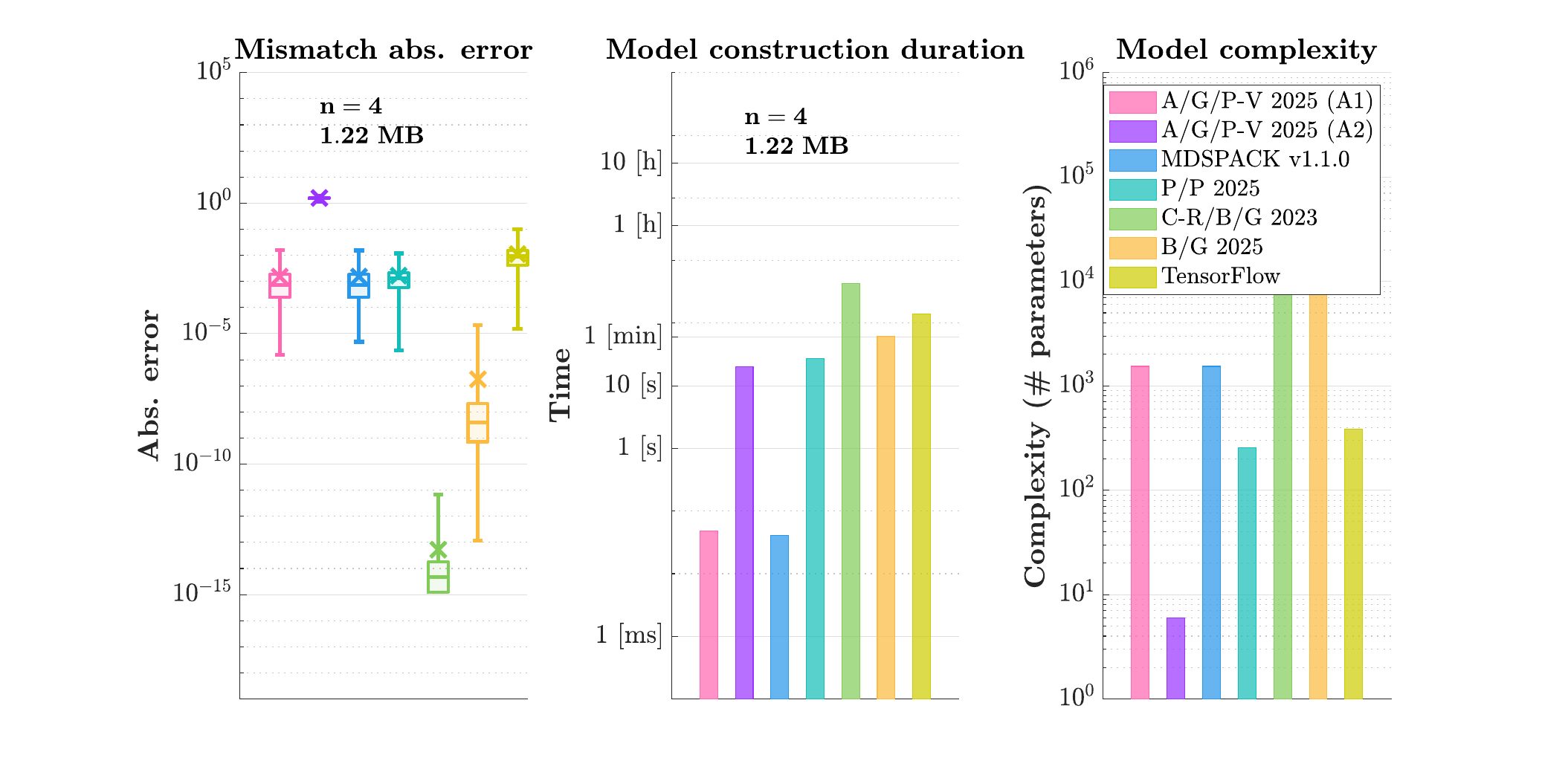} \caption{Function \#22: graphical view of the best model performances.} \end{figure}\begin{figure}[H] \centering  \includegraphics[width=\textwidth]{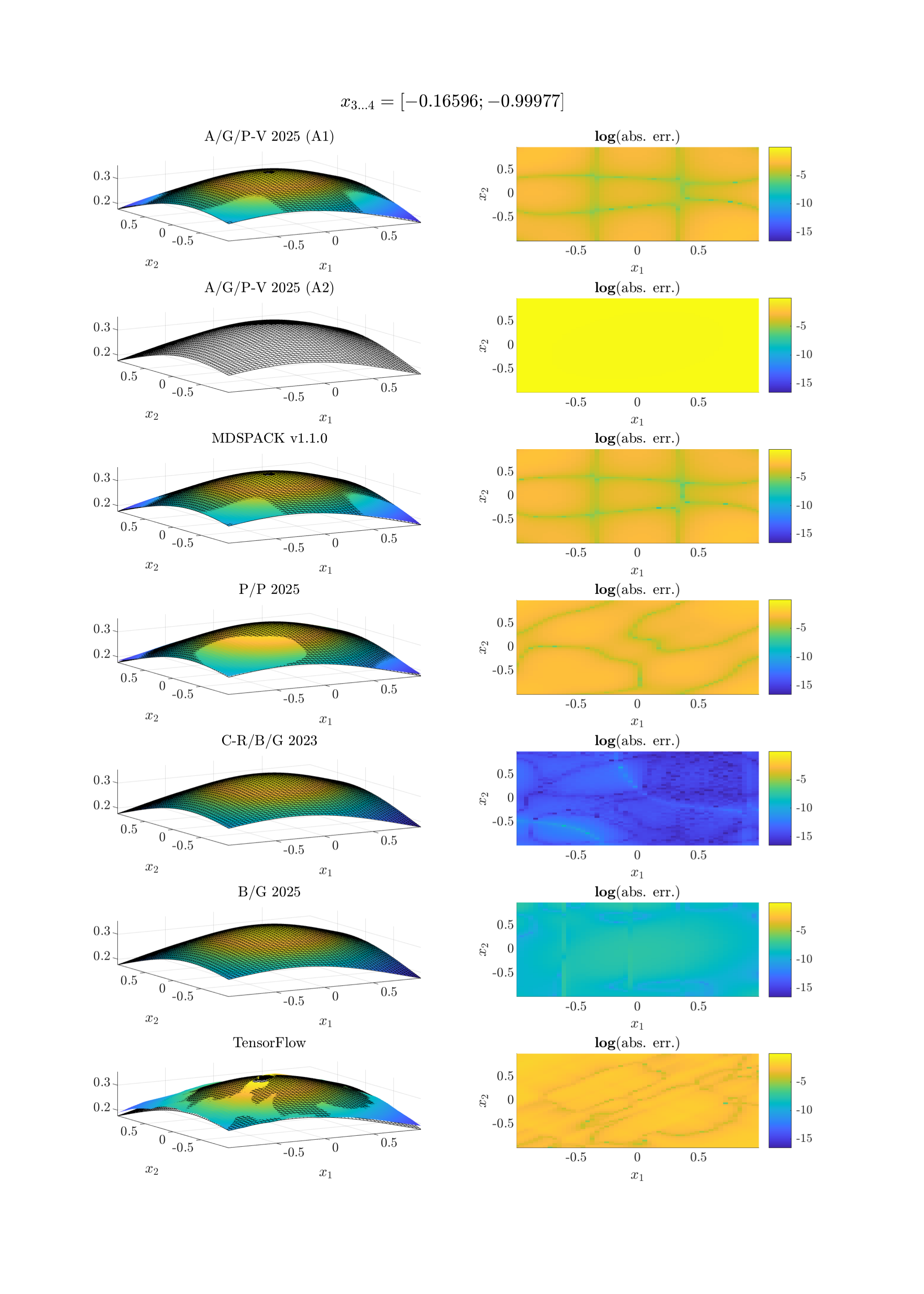} \caption{Function \#22: left side, evaluation of the original (mesh) vs. approximated (coloured surface) and right side, absolute errors (in log-scale).} \end{figure}\subsubsection{mLF detailed informations (M1)} \noindent \textbf{Right interpolation points}: $k_l=\left(\begin{array}{cccc} 4 & 4 & 4 & 4 \end{array}\right)$, where $l=1,\cdots,\ord$.$$ \begin{array}{rcl}\lan{1} &\in& \IC^{4} \text{ , linearly spaced between bounds}\\\lan{2} &\in& \IC^{4} \text{ , linearly spaced between bounds}\\\lan{3} &\in& \IC^{4} \text{ , linearly spaced between bounds}\\\lan{4} &\in& \IC^{4} \text{ , linearly spaced between bounds}\\\end{array} $$\noindent \textbf{$\ord$-D Loewner matrix, barycentric weights and Lagrangian basis}:$$ \begin{array}{rcl}\IL & \in & \IC^{256 \times 256}\\\bc & \in & \IC^{256}\\\bw & \in & \IC^{256}\\\bc\odot \bw & \in & \IC^{256}\\\mathbf{Lag}(\var{1},\var{2},\var{3},\var{4}) & \in & \IC^{256}\\\end{array} $$

\newpage \subsection{Function \#23 (${\ord=4}$ variables, tensor size: 1.79 \textbf{MB})} $$10\prod_{i=1}^4\mathrm{sinc}(x_i)$$ \subsubsection{Setup and results overview}\begin{itemize}\item Reference: A/al. 2021 (A.5.18), \cite{Austin:2021}\item Domain: $\mathbb{R}$\item Tensor size: 1.79 \textbf{MB} ($22^{4}$ points)\item Bounds: $ \left(\begin{array}{cc} \frac{1}{1000000} & 4\,\pi  \end{array}\right) \times \left(\begin{array}{cc} \frac{1}{1000000} & 4\,\pi  \end{array}\right) \times \left(\begin{array}{cc} \frac{1}{1000000} & 4\,\pi  \end{array}\right) \times \left(\begin{array}{cc} \frac{1}{1000000} & 4\,\pi  \end{array}\right)$ \end{itemize} \begin{table}[H] \centering \begin{tabular}{llllll}
$\#$ & Alg. & Parameters & Dim. & CPU [s] & RMSE \\ 
\hline 
$\mathbf{\#23}$ & A/G/P-V 2025 (A1) & $1 \cdot 10^{-09},1$ & $8.8 \cdot 10^{04}$ & $0.33$ & $5.6 \cdot 10^{-08}$ \\ 
 & A/G/P-V 2025 (A2) & $1 \cdot 10^{-15},1$ & $2.5 \cdot 10^{04}$ & $1.9 \cdot 10^{02}$ & $1.6 \cdot 10^{-05}$ \\ 
 & MDSPACK v1.1.0 & $0.01,1$ & $\mathbf{6}$ & $\mathbf{0.062}$ & $0.17$ \\ 
 & P/P 2025 & $1,0.95,50,0.01,6,12,13$ & $4.7 \cdot 10^{02}$ & $82$ & $0.018$ \\ 
 & C-R/B/G 2023 & $0.001,20$ & $NaN$ & $NaN$ & $NaN$ \\ 
 & B/G 2025 & $1 \cdot 10^{-09},20,3$ & $8.8 \cdot 10^{04}$ & $41$ & $\mathbf{9.6 \cdot 10^{-11}}$ \\ 
 & TensorFlow & $$ & $3.8 \cdot 10^{02}$ & $2.2 \cdot 10^{02}$ & $0.067$ \\ 
\hline 
\end{tabular} \caption{Function \#23: best model configuration and performances per methods.} \end{table}\begin{figure}[H] \centering  \includegraphics[width=\textwidth]{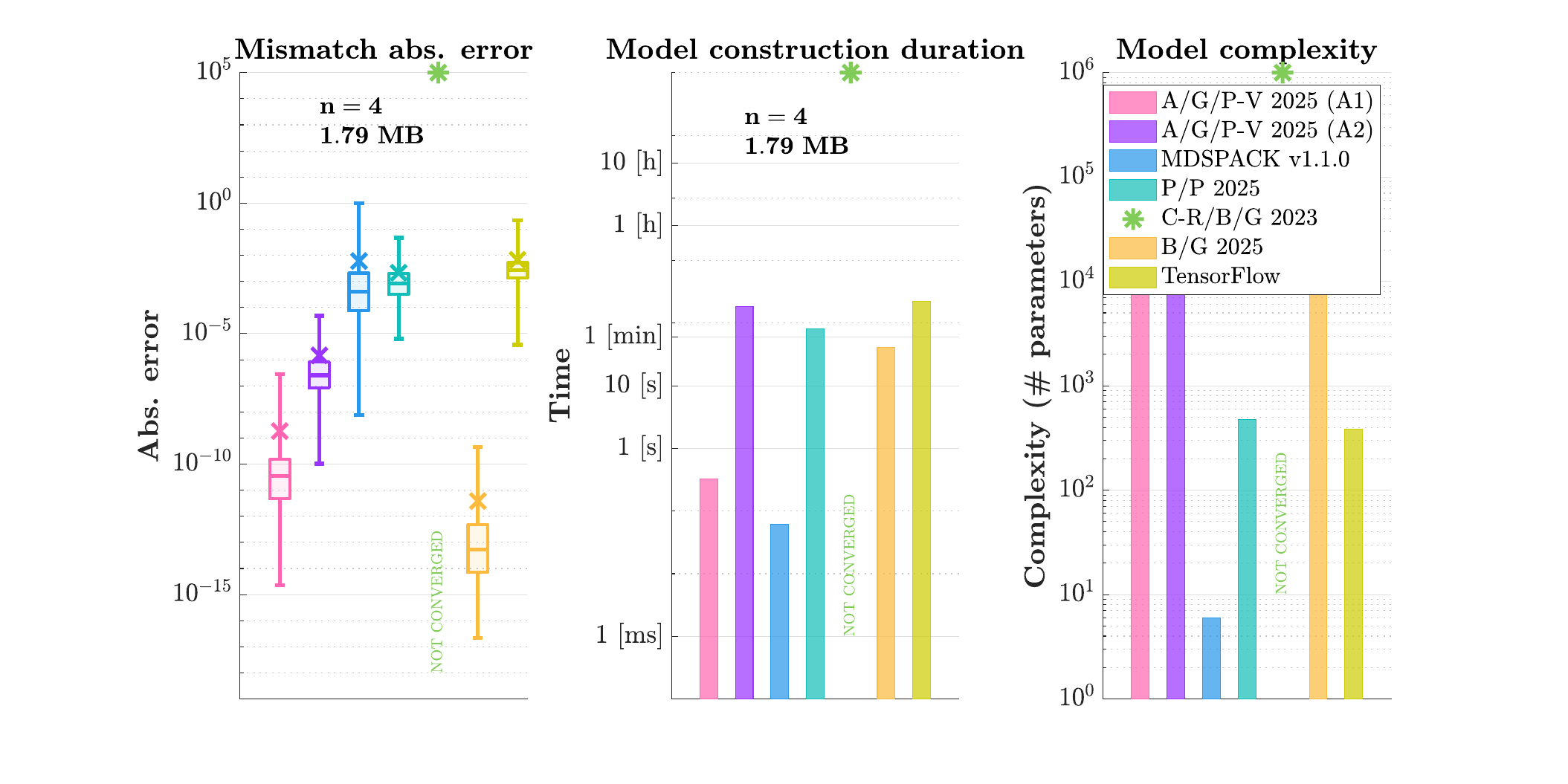} \caption{Function \#23: graphical view of the best model performances.} \end{figure}\begin{figure}[H] \centering  \includegraphics[width=\textwidth]{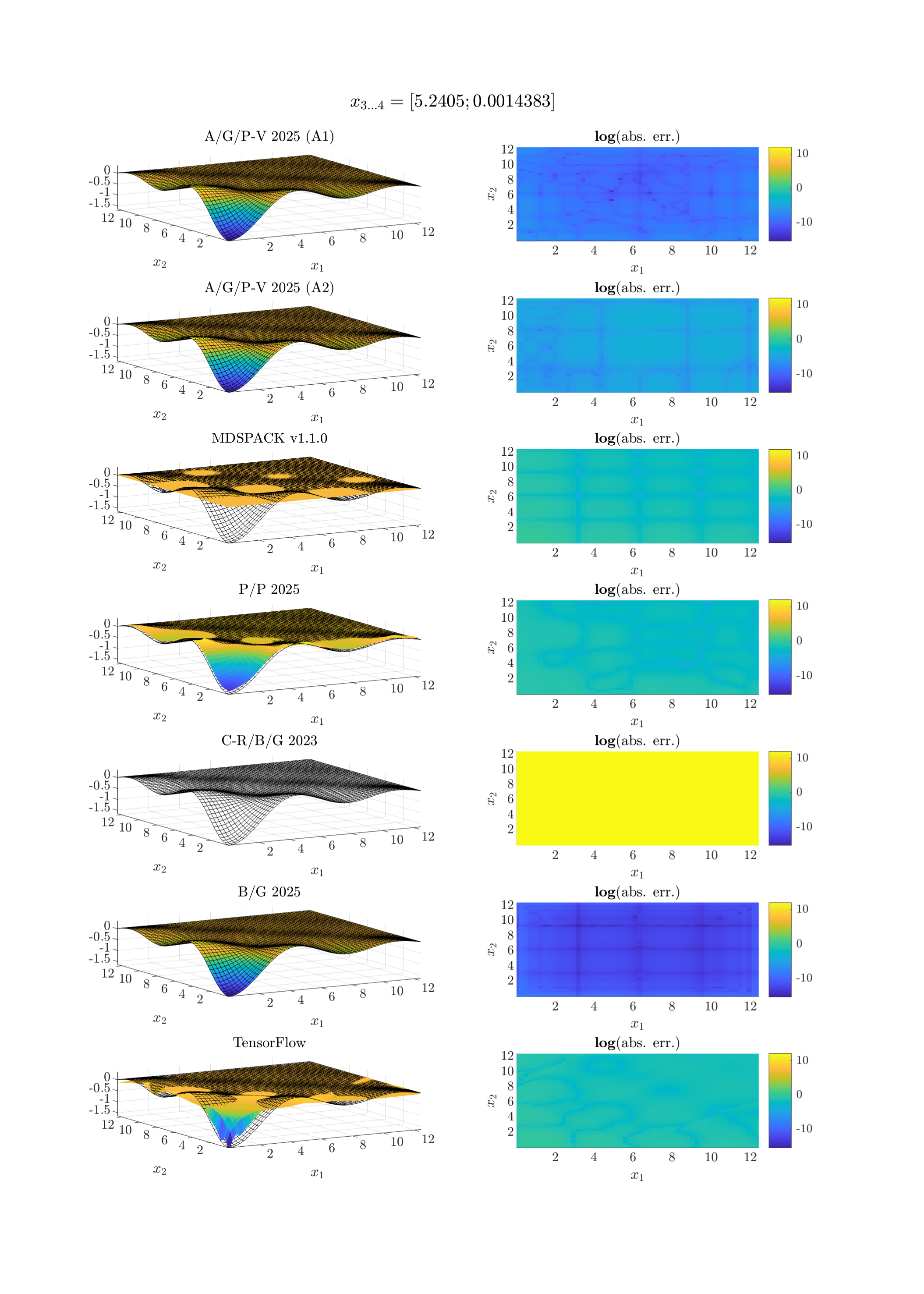} \caption{Function \#23: left side, evaluation of the original (mesh) vs. approximated (coloured surface) and right side, absolute errors (in log-scale).} \end{figure}\subsubsection{mLF detailed informations (M1)} \noindent \textbf{Right interpolation points}: $k_l=\left(\begin{array}{cccc} 11 & 11 & 11 & 11 \end{array}\right)$, where $l=1,\cdots,\ord$.$$ \begin{array}{rcl}\lan{1} &\in& \IC^{11} \text{ , linearly spaced between bounds}\\\lan{2} &\in& \IC^{11} \text{ , linearly spaced between bounds}\\\lan{3} &\in& \IC^{11} \text{ , linearly spaced between bounds}\\\lan{4} &\in& \IC^{11} \text{ , linearly spaced between bounds}\\\end{array} $$\noindent \textbf{$\ord$-D Loewner matrix, barycentric weights and Lagrangian basis}:$$ \begin{array}{rcl}\IL & \in & \IC^{14641 \times 14641}\\\bc & \in & \IC^{14641}\\\bw & \in & \IC^{14641}\\\bc\odot \bw & \in & \IC^{14641}\\\mathbf{Lag}(\var{1},\var{2},\var{3},\var{4}) & \in & \IC^{14641}\\\end{array} $$

\newpage \subsection{Function \#24 (${\ord=2}$ variables, tensor size: 13.8 \textbf{KB})} $$10\mathrm{sinc}(\var{1})\mathrm{sinc}(\var{2})$$ \subsubsection{Setup and results overview}\begin{itemize}\item Reference: A/al. 2021 (A.5.19), \cite{Austin:2021}\item Domain: $\mathbb{R}$\item Tensor size: 13.8 \textbf{KB} ($42^{2}$ points)\item Bounds: $ \left(\begin{array}{cc} \frac{1}{1000000} & 4\,\pi  \end{array}\right) \times \left(\begin{array}{cc} \frac{1}{1000000} & 4\,\pi  \end{array}\right)$ \end{itemize} \begin{table}[H] \centering \begin{tabular}{llllll}
$\#$ & Alg. & Parameters & Dim. & CPU [s] & RMSE \\ 
\hline 
$\mathbf{\#24}$ & A/G/P-V 2025 (A1) & $1 \cdot 10^{-09},1$ & $4.8 \cdot 10^{02}$ & $0.012$ & $6.9 \cdot 10^{-08}$ \\ 
 & A/G/P-V 2025 (A2) & $1 \cdot 10^{-15},3$ & $3.2 \cdot 10^{02}$ & $0.37$ & $0.0016$ \\ 
 & MDSPACK v1.1.0 & $1 \cdot 10^{-14},7$ & $4.3 \cdot 10^{02}$ & $\mathbf{0.0035}$ & $1.4$ \\ 
 & P/P 2025 & $1,0.95,50,0.01,6,12,13$ & $3.2 \cdot 10^{02}$ & $1.2$ & $0.0027$ \\ 
 & C-R/B/G 2023 & $1 \cdot 10^{-09},20$ & $4.8 \cdot 10^{02}$ & $0.12$ & $2 \cdot 10^{-10}$ \\ 
 & B/G 2025 & $1 \cdot 10^{-09},20,3$ & $5.3 \cdot 10^{02}$ & $0.25$ & $\mathbf{1.9 \cdot 10^{-10}}$ \\ 
 & TensorFlow & $$ & $\mathbf{2.6 \cdot 10^{02}}$ & $16$ & $0.24$ \\ 
\hline 
\end{tabular} \caption{Function \#24: best model configuration and performances per methods.} \end{table}\begin{figure}[H] \centering  \includegraphics[width=\textwidth]{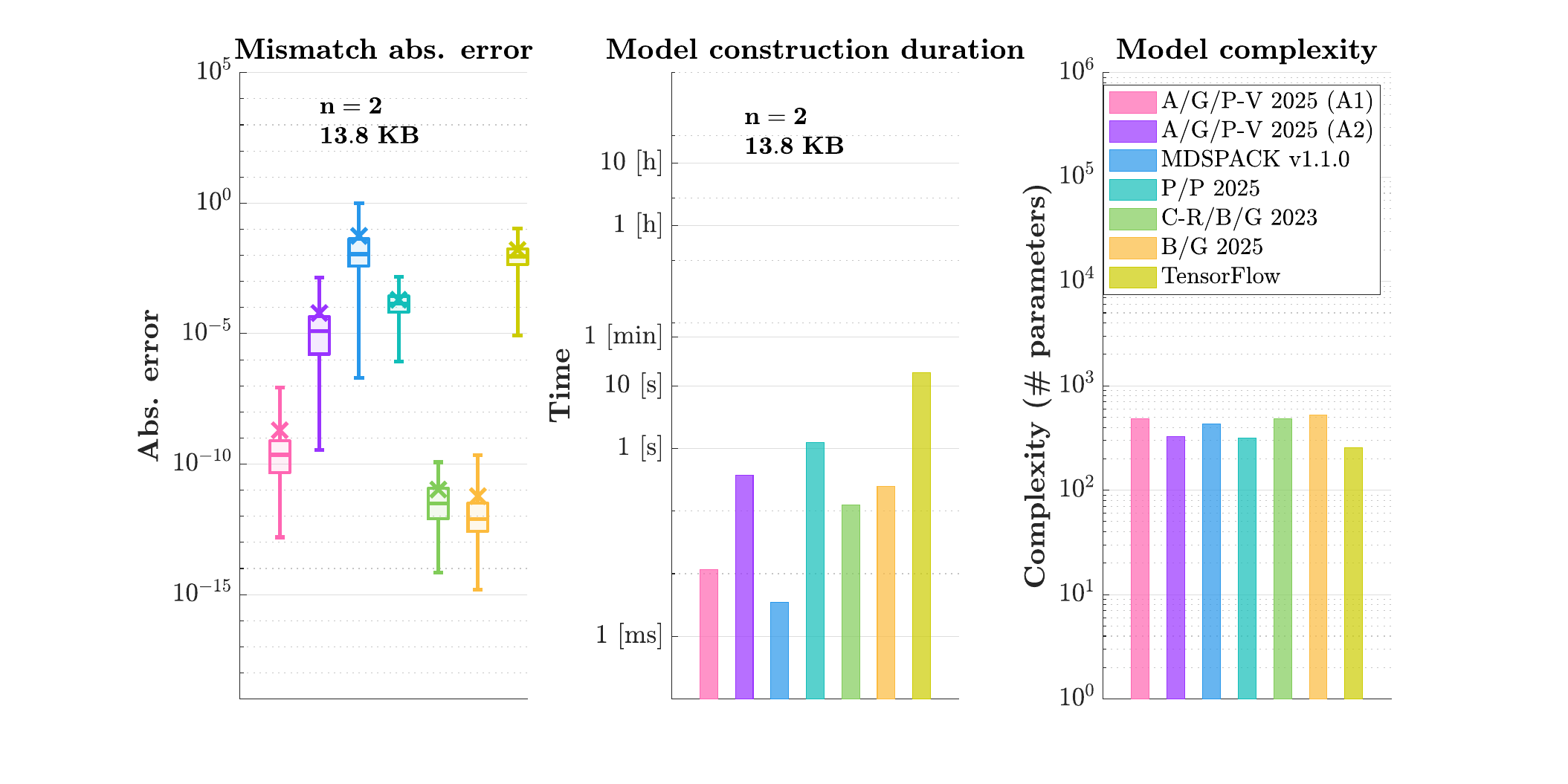} \caption{Function \#24: graphical view of the best model performances.} \end{figure}\begin{figure}[H] \centering  \includegraphics[width=\textwidth]{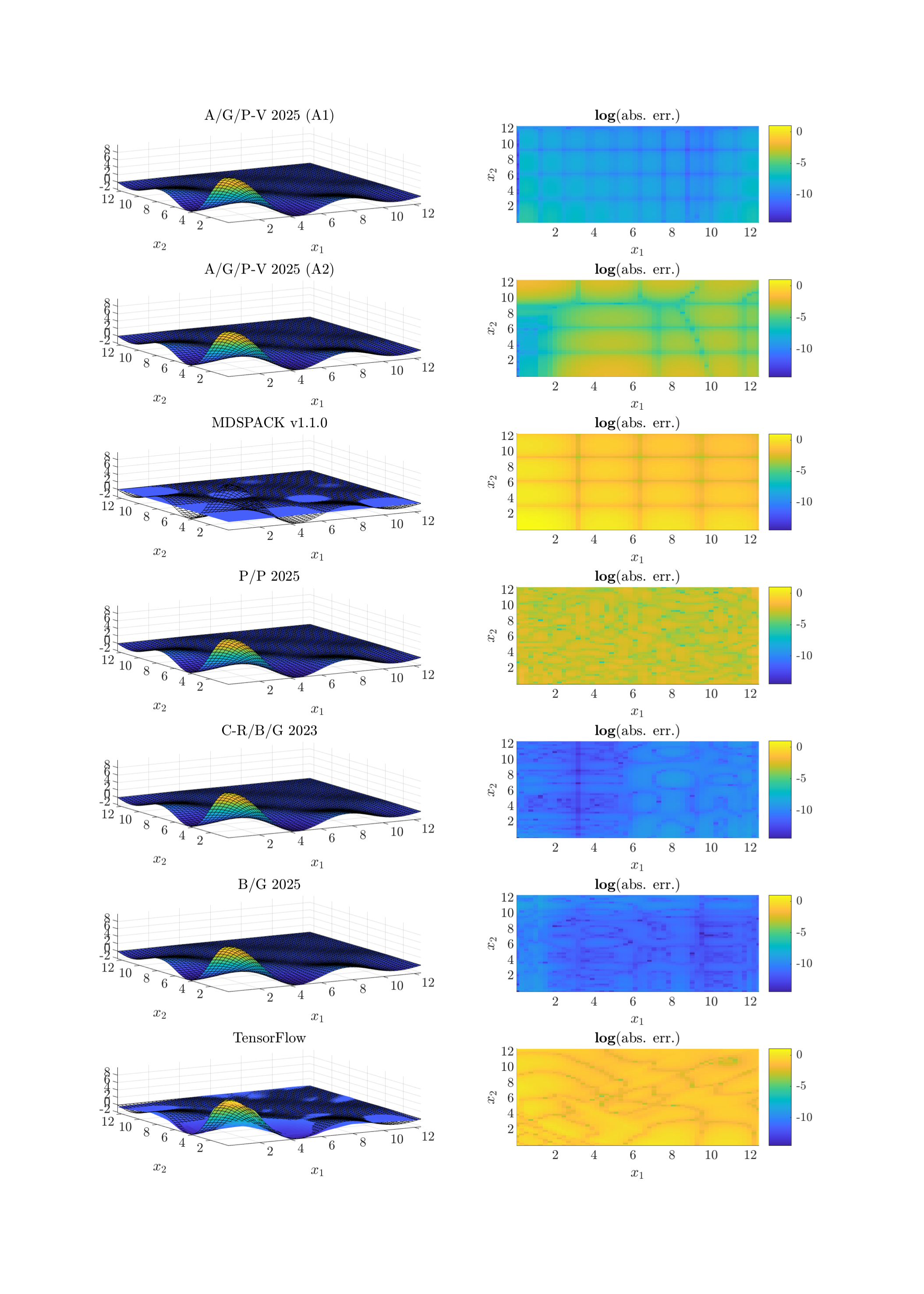} \caption{Function \#24: left side, evaluation of the original (mesh) vs. approximated (coloured surface) and right side, absolute errors (in log-scale).} \end{figure}\subsubsection{mLF detailed informations (M1)} \noindent \textbf{Right interpolation points}: $k_l=\left(\begin{array}{cc} 11 & 11 \end{array}\right)$, where $l=1,\cdots,\ord$.$$ \begin{array}{rcl}\lan{1} &\in& \IC^{11} \text{ , linearly spaced between bounds}\\\lan{2} &\in& \IC^{11} \text{ , linearly spaced between bounds}\\\end{array} $$\noindent \textbf{$\ord$-D Loewner matrix, barycentric weights and Lagrangian basis}:$$ \begin{array}{rcl}\IL & \in & \IC^{121 \times 121}\\\bc & \in & \IC^{121}\\\bw & \in & \IC^{121}\\\bc\odot \bw & \in & \IC^{121}\\\mathbf{Lag}(\var{1},\var{2}) & \in & \IC^{121}\\\end{array} $$

\newpage \subsection{Function \#25 (${\ord=2}$ variables, tensor size: 12.5 \textbf{KB})} $$\var{1}^2+\var{2}^2+\var{1}\var{2}-\var{2}+1$$ \subsubsection{Setup and results overview}\begin{itemize}\item Reference: A/al. 2021 (A.5.20), \cite{Austin:2021}\item Domain: $\mathbb{R}$\item Tensor size: 12.5 \textbf{KB} ($40^{2}$ points)\item Bounds: $ \left(\begin{array}{cc} -1 & 1 \end{array}\right) \times \left(\begin{array}{cc} -1 & 1 \end{array}\right)$ \end{itemize} \begin{table}[H] \centering \begin{tabular}{llllll}
$\#$ & Alg. & Parameters & Dim. & CPU [s] & RMSE \\ 
\hline 
$\mathbf{\#25}$ & A/G/P-V 2025 (A1) & $0.5,1$ & $\mathbf{36}$ & $0.036$ & $\mathbf{1 \cdot 10^{-15}}$ \\ 
 & A/G/P-V 2025 (A2) & $1 \cdot 10^{-15},3$ & $36$ & $0.067$ & $2 \cdot 10^{-15}$ \\ 
 & MDSPACK v1.1.0 & $0.01,1$ & $36$ & $0.022$ & $1.1 \cdot 10^{-15}$ \\ 
 & P/P 2025 & $1,0.95,50,0.01,10,4,21$ & $5.1 \cdot 10^{02}$ & $1.4$ & $3.9 \cdot 10^{-06}$ \\ 
 & C-R/B/G 2023 & $0.001,20$ & $72$ & $\mathbf{0.015}$ & $1.4 \cdot 10^{-13}$ \\ 
 & B/G 2025 & $1 \cdot 10^{-06},20,3$ & $72$ & $0.015$ & $7.6 \cdot 10^{-15}$ \\ 
 & TensorFlow & $$ & $2.6 \cdot 10^{02}$ & $14$ & $0.0071$ \\ 
\hline 
\end{tabular} \caption{Function \#25: best model configuration and performances per methods.} \end{table}\begin{figure}[H] \centering  \includegraphics[width=\textwidth]{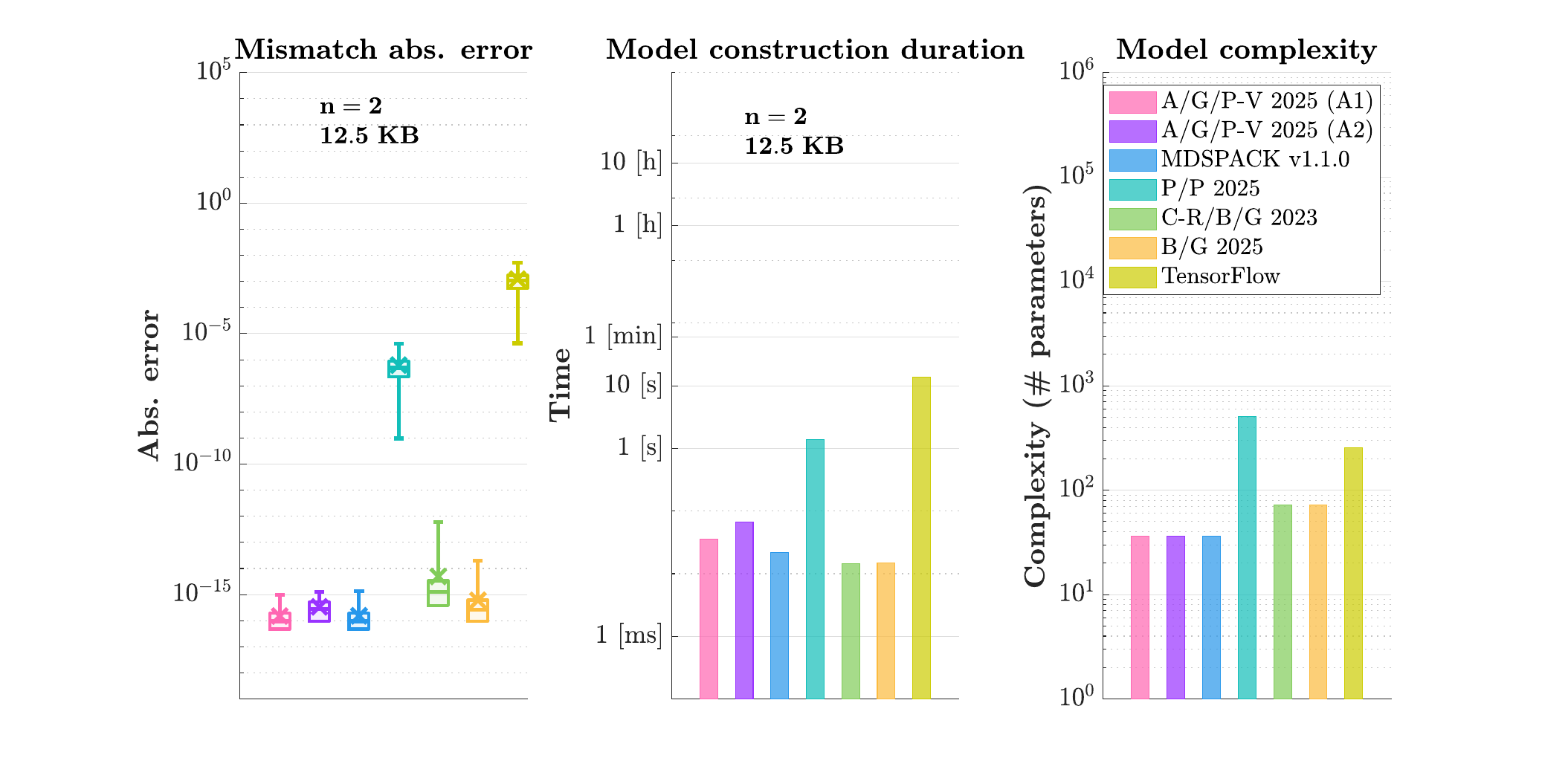} \caption{Function \#25: graphical view of the best model performances.} \end{figure}\begin{figure}[H] \centering  \includegraphics[width=\textwidth]{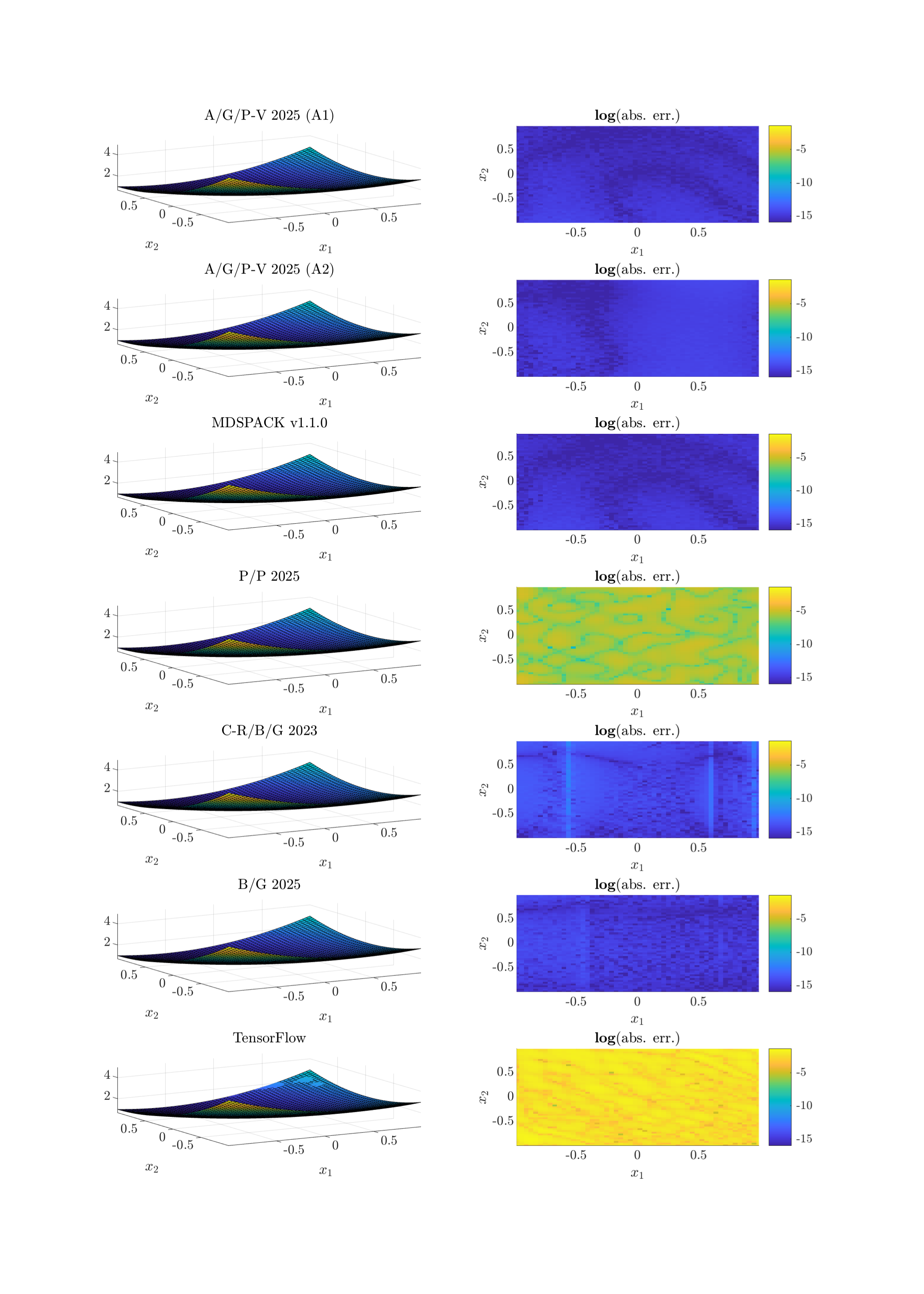} \caption{Function \#25: left side, evaluation of the original (mesh) vs. approximated (coloured surface) and right side, absolute errors (in log-scale).} \end{figure}\subsubsection{mLF detailed informations (M1)} \noindent \textbf{Right interpolation points} ($k_l=\left(\begin{array}{cc} 3 & 3 \end{array}\right)$, where $l=1,\cdots,\ord$):$$ \begin{array}{rcl}\lan{1} &=& \left(\begin{array}{ccc} -1 & -\frac{1}{19} & 1 \end{array}\right)\\\lan{2} &=& \left(\begin{array}{ccc} -1 & -\frac{1}{19} & 1 \end{array}\right)\\\end{array} $$\noindent \textbf{Lagrangian weights}: $$\left(\begin{array}{ccc} \bc & \bw & \bc\odot\bw\\ -0.5848 & 5.0 & -2.924\\ 1.111 & 2.108 & 2.342\\ -0.5263 & 1.0 & -0.5263\\ 1.111 & 3.055 & 3.395\\ -2.111 & 1.061 & -2.24\\ 1.0 & 0.9501 & 0.9501\\ -0.5263 & 3.0 & -1.579\\ 1.0 & 2.003 & 2.003\\ -0.4737 & 3.0 & -1.421 \end{array}\right)$$\noindent \textbf{Lagrangian form} (basis, numerator and denominator coefficients):$$\left(\begin{array}{ccc}\mathcal{B}_\textrm{lag}(\var{1},\var{2}) & \bN_\textrm{lag} &\bD_\textrm{lag}\end{array}\right) =$$ $$\left(\begin{array}{ccc} \left(\var{1}+1.0\right)\,\left(\var{2}+1.0\right) & -2.924 & -0.5848\\ \left(\var{1}+1.0\right)\,\left(\var{2}+0.05263\right) & 2.342 & 1.111\\ \left(\var{1}+1.0\right)\,\left(\var{2}-1.0\right) & -0.5263 & -0.5263\\ \left(\var{2}+1.0\right)\,\left(\var{1}+0.05263\right) & 3.395 & 1.111\\ \left(\var{1}+0.05263\right)\,\left(\var{2}+0.05263\right) & -2.24 & -2.111\\ \left(\var{2}-1.0\right)\,\left(\var{1}+0.05263\right) & 0.9501 & 1.0\\ \left(\var{1}-1.0\right)\,\left(\var{2}+1.0\right) & -1.579 & -0.5263\\ \left(\var{1}-1.0\right)\,\left(\var{2}+0.05263\right) & 2.003 & 1.0\\ \left(\var{1}-1.0\right)\,\left(\var{2}-1.0\right) & -1.421 & -0.4737 \end{array}\right).$$\noindent The corresponding function is:$$\begin{array}{rcl}\bG_{\textrm{lag}}(\var{1},\var{2}) &=& \dfrac{\bn_{\textrm{lag}}(\var{1},\var{2})}{\bd_{\textrm{lag}}(\var{1},\var{2})}\\ && \\&=& \dfrac{\sum_{\textrm{row}} \bN_\textrm{lag} \odot\mathcal{B}^{-1}_\textrm{lag}(\var{1},\var{2})}{\sum_{\textrm{row}} \bD_\textrm{lag} \odot\mathcal{B}^{-1}_\textrm{lag}(\var{1},\var{2})}, \end{array}$$\noindent where,\\$\bn_{\textrm{lag}}(\var{1},\var{2}) = {\var{1}}^2+\var{1}\,\var{2}+{\var{2}}^2-1.0\,\var{2}+1.0$ \\~~\\$\bd_{\textrm{lag}}(\var{1},\var{2}) = 1.0$ \\~~\\\noindent \textbf{Monomial form} (basis, numerator and denominator coefficients - entries $<10^{-12}$ removed):$$\left(\begin{array}{ccc}\mathcal{B}_\textrm{mon}(\var{1},\var{2}) & \bN_\textrm{mon} &\bD_\textrm{mon}\end{array}\right) =$$ $$\left(\begin{array}{ccc} {\var{1}}^2\,{\var{2}}^2 & 0 & 0\\ {\var{1}}^2\,\var{2} & 0 & 0\\ {\var{1}}^2 & -1.0 & 0\\ \var{1}\,{\var{2}}^2 & 0 & 0\\ \var{1}\,\var{2} & -1.0 & 0\\ \var{1} & 0 & 0\\ {\var{2}}^2 & -1.0 & 0\\ \var{2} & 1.0 & 0\\ 1.0 & -1.0 & -1.0 \end{array}\right)$$\noindent The corresponding function is:$$\begin{array}{rcl}\bG_{\textrm{mon}}(\var{1},\var{2}) &=& \dfrac{\bn_{\textrm{mon}}(\var{1},\var{2})}{\bd_{\textrm{mon}}(\var{1},\var{2})}\\ && \\&=& \dfrac{\sum_{\textrm{row}} \bN_\textrm{mon} \odot \mathcal{B}_\textrm{mon}(\var{1},\var{2})}{\sum_{\textrm{row}} \bD_\textrm{mon} \odot\mathcal{B}_\textrm{mon}(\var{1},\var{2})},  \end{array}$$\noindent where,\\$\bn_{\textrm{mon}}(\var{1},\var{2}) = {\var{1}}^2+\var{1}\,\var{2}+{\var{2}}^2-1.0\,\var{2}+1.0$ \\~~\\$\bd_{\textrm{mon}}(\var{1},\var{2}) = 1.0$ \\~~\\\noindent \textbf{KST equivalent decoupling pattern} (Barycentric weights $\bc^{\var{l}}$): $$\begin{array}{rclll}\var{2}&: & \left(\begin{array}{ccc} 1.111 & 1.111 & 1.111\\ -2.111 & -2.111 & -2.111\\ 1.0 & 1.0 & 1.0 \end{array}\right)& \textrm{vec}(.) & := \textbf{Bary}(\var{2}) \\\var{1}&: & \left(\begin{array}{c} -0.5263\\ 1.0\\ -0.4737 \end{array}\right)& \textrm{vec}(.) \otimes \bone_{k_{2}} & := \textbf{Bary}(\var{1}) \\\end{array}$$~\\ Then, with the above notations, one defines the following univariate vector functions:  $$ \left\{ \begin{array}{rcl}\bPhi_{1}(\var{1}) &:=& \textbf{Bary}(\var{1}) \odot \mathbf{Lag}(\var{1}) \\\bPhi_{2}(\var{2}) &:=& \textbf{Bary}(\var{2}) \odot \mathbf{Lag}(\var{2}) \\\end{array} \right. $$\noindent The corresponding function is:$$\begin{array}{rcl}\bG_{\textrm{kst}}(\var{1},\var{2}) &=& \dfrac{\bn_{\textrm{kst}}(\var{1},\var{2})}{\bd_{\textrm{kst}}(\var{1},\var{2})}\\ && \\ &=& \dfrac{\sum_{\text{rows}} \bw \odot \bPhi_{1}(\var{1}) \odot \cdots \odot\bPhi_{2}(\var{2})}{\sum_{\text{rows}} \bPhi_{1}(\var{1}) \odot \cdots \odot\bPhi_{2}(\var{2})} . \end{array}$$~\\ \noindent \textbf{KST-like univariate functions} (equivalent scaled univariate functions $\bphi_{1,\cdots,2}$): $$\left\{\begin{array}{rcrcl}z_{1} &=&\bphi_{1}(\var{1}) &=& {\var{1}}^2+\var{1}+1.0\\z_{2} &=&\bphi_{2}(\var{2}) &=& {\var{2}}^2-2.0\,\var{2}+2.0\\\end{array} \right. .$$\noindent \textbf{Connection with Neural Networks} (equivalent numerator $\bn_{\textrm{lag}}$ representation):\\ \begin{figure}[H]\begin{center} \scalebox{.7}{\begin{tikzpicture}[line width=0.4mm]\tikzstyle{place}=[circle, draw=black, minimum size = 8mm]\tikzstyle{placeInOut}=[circle, draw=orange, minimum size = 8mm]\node at (0,-2) [placeInOut] (first_1){$\var{1}$};\node at (0,-4) [placeInOut] (first_2){$\var{2}$};\node at (5,-2) [place] (secondL1_1){$\frac{1}{\var{1}-\lani{1}{1}}$};\node at (5,-4) [place] (secondL1_2){$\frac{1}{\var{1}-\lani{1}{2}}$};\node at (5,-6) [place] (secondL1_3){$\frac{1}{\var{1}-\lani{1}{3}}$};\node at (5,-8) [place] (secondL2_1){$\frac{1}{\var{2}-\lani{2}{1}}$};\node at (5,-10) [place] (secondL2_2){$\frac{1}{\var{2}-\lani{2}{2}}$};\node at (5,-12) [place] (secondL2_3){$\frac{1}{\var{2}-\lani{2}{3}}$};\node at (10,-2) [place] (third_1){$\prod$};\node at (10,-4) [place] (third_2){$\prod$};\node at (10,-6) [place] (third_3){$\prod$};\node at (10,-8) [place] (third_4){$\prod$};\node at (10,-10) [place] (third_5){$\prod$};\node at (10,-12) [place] (third_6){$\prod$};\node at (10,-14) [place] (third_7){$\prod$};\node at (10,-16) [place] (third_8){$\prod$};\node at (10,-18) [place] (third_9){$\prod$};\node at (15,-10) [placeInOut] (output){$\bSigma$};\draw[->] (first_1)--(secondL1_1) node[above,sloped,pos=0.75] { };\draw[->] (first_1)--(secondL1_2) node[above,sloped,pos=0.75] { };\draw[->] (first_1)--(secondL1_3) node[above,sloped,pos=0.75] { };\draw[->] (first_2)--(secondL2_1) node[above,sloped,pos=0.75] { };\draw[->] (first_2)--(secondL2_2) node[above,sloped,pos=0.75] { };\draw[->] (first_2)--(secondL2_3) node[above,sloped,pos=0.75] { };\draw[->] (secondL1_1)--(third_1) node[above,sloped,pos=0.25] {};\draw[->] (secondL1_1)--(third_2) node[above,sloped,pos=0.25] {};\draw[->] (secondL1_1)--(third_3) node[above,sloped,pos=0.25] {};\draw[->] (secondL1_2)--(third_4) node[above,sloped,pos=0.25] {};\draw[->] (secondL1_2)--(third_5) node[above,sloped,pos=0.25] {};\draw[->] (secondL1_2)--(third_6) node[above,sloped,pos=0.25] {};\draw[->] (secondL1_3)--(third_7) node[above,sloped,pos=0.25] {};\draw[->] (secondL1_3)--(third_8) node[above,sloped,pos=0.25] {};\draw[->] (secondL1_3)--(third_9) node[above,sloped,pos=0.25] {};\draw[->] (secondL2_1)--(third_1) node[above,sloped,pos=0.25] {};\draw[->] (secondL2_2)--(third_2) node[above,sloped,pos=0.25] {};\draw[->] (secondL2_3)--(third_3) node[above,sloped,pos=0.25] {};\draw[->] (secondL2_1)--(third_4) node[above,sloped,pos=0.25] {};\draw[->] (secondL2_2)--(third_5) node[above,sloped,pos=0.25] {};\draw[->] (secondL2_3)--(third_6) node[above,sloped,pos=0.25] {};\draw[->] (secondL2_1)--(third_7) node[above,sloped,pos=0.25] {};\draw[->] (secondL2_2)--(third_8) node[above,sloped,pos=0.25] {};\draw[->] (secondL2_3)--(third_9) node[above,sloped,pos=0.25] {};\draw[->] (third_1)--(output) node[above,sloped,pos=0.25] {-2.924};\draw[->] (third_2)--(output) node[above,sloped,pos=0.25] {2.3423};\draw[->] (third_3)--(output) node[above,sloped,pos=0.25] {-0.52632};\draw[->] (third_4)--(output) node[above,sloped,pos=0.25] {3.3949};\draw[->] (third_5)--(output) node[above,sloped,pos=0.25] {-2.2398};\draw[->] (third_6)--(output) node[above,sloped,pos=0.25] {0.95014};\draw[->] (third_7)--(output) node[above,sloped,pos=0.25] {-1.5789};\draw[->] (third_8)--(output) node[above,sloped,pos=0.25] {2.0028};\draw[->] (third_9)--(output) node[above,sloped,pos=0.25] {-1.4211};\end{tikzpicture}} \caption{Equivalent NN representation of the numerator $\bn_{\textrm{lag}}$.}\end{center}\end{figure}

\newpage \subsection{Function \#26 (${\ord=3}$ variables, tensor size: 1.65 \textbf{MB})} $$\frac{\var{1}+\var{2}+\var{3}}{6+\cos(\var{1})+\cos(\var{2})+\cos(\var{3})}$$ \subsubsection{Setup and results overview}\begin{itemize}\item Reference: B/G 2025, \cite{Balicki:2025}\item Domain: $\mathbb{R}$\item Tensor size: 1.65 \textbf{MB} ($60^{3}$ points)\item Bounds: $ \left(\begin{array}{cc} -10 & 10 \end{array}\right) \times \left(\begin{array}{cc} -10 & 10 \end{array}\right) \times \left(\begin{array}{cc} -10 & 10 \end{array}\right)$ \end{itemize} \begin{table}[H] \centering \begin{tabular}{llllll}
$\#$ & Alg. & Parameters & Dim. & CPU [s] & RMSE \\ 
\hline 
$\mathbf{\#26}$ & A/G/P-V 2025 (A1) & $0.001,1$ & $5 \cdot 10^{03}$ & $0.079$ & $0.075$ \\ 
 & A/G/P-V 2025 (A2) & $1 \cdot 10^{-15},1$ & $\mathbf{5}$ & $18$ & $9.6$ \\ 
 & MDSPACK v1.1.0 & $0.0001,2$ & $5 \cdot 10^{03}$ & $\mathbf{0.053}$ & $0.075$ \\ 
 & P/P 2025 & $1,0.95,50,0.01,10,6,21$ & $7.6 \cdot 10^{02}$ & $90$ & $0.0032$ \\ 
 & C-R/B/G 2023 & $1 \cdot 10^{-09},20$ & $2 \cdot 10^{04}$ & $2.2 \cdot 10^{03}$ & $\mathbf{1.6 \cdot 10^{-10}}$ \\ 
 & B/G 2025 & $1 \cdot 10^{-09},20,4$ & $1.9 \cdot 10^{04}$ & $80$ & $1.8 \cdot 10^{-08}$ \\ 
 & TensorFlow & $$ & $3.2 \cdot 10^{02}$ & $1.3 \cdot 10^{03}$ & $0.49$ \\ 
\hline 
\end{tabular} \caption{Function \#26: best model configuration and performances per methods.} \end{table}\begin{figure}[H] \centering  \includegraphics[width=\textwidth]{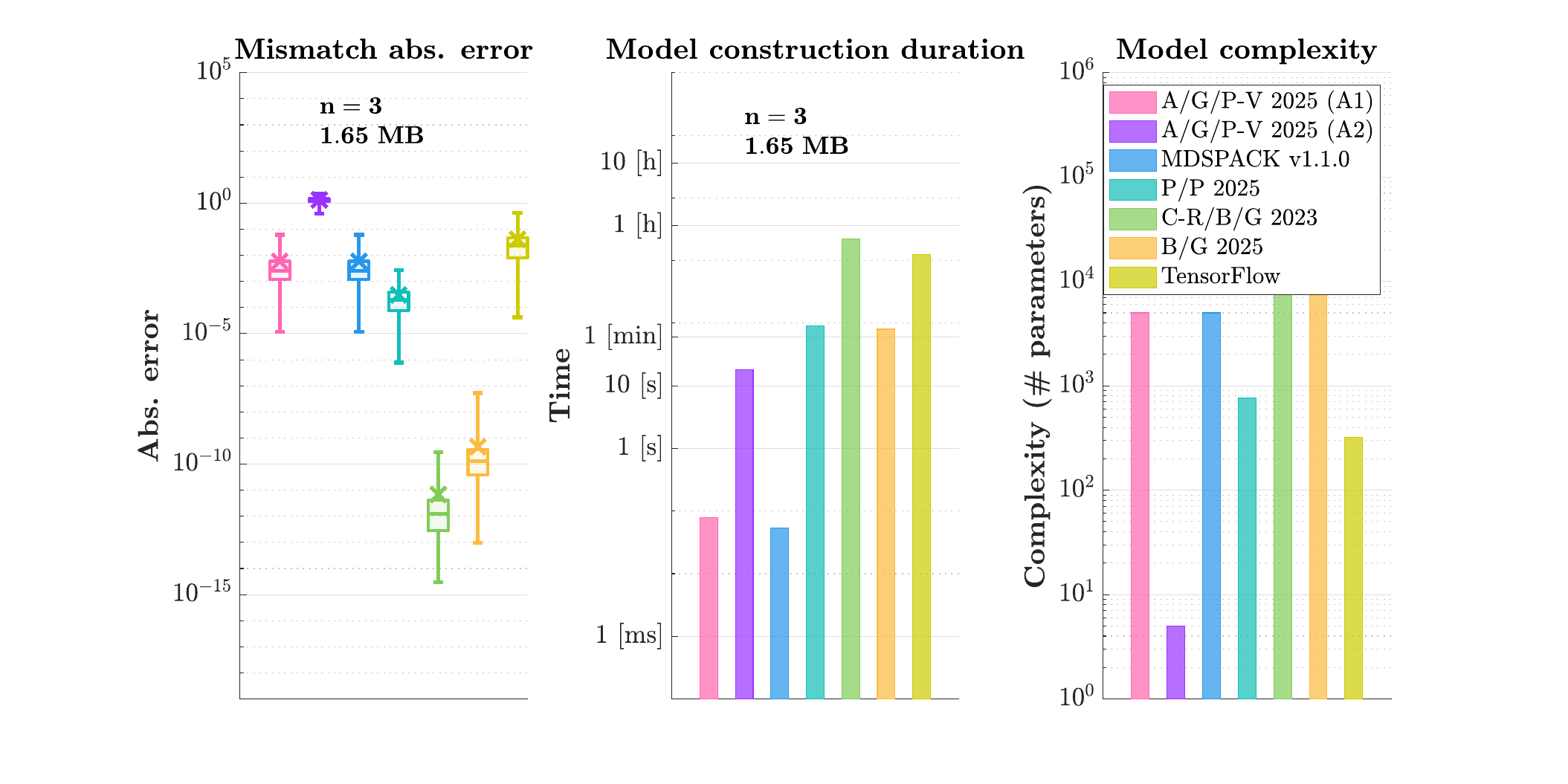} \caption{Function \#26: graphical view of the best model performances.} \end{figure}\begin{figure}[H] \centering  \includegraphics[width=\textwidth]{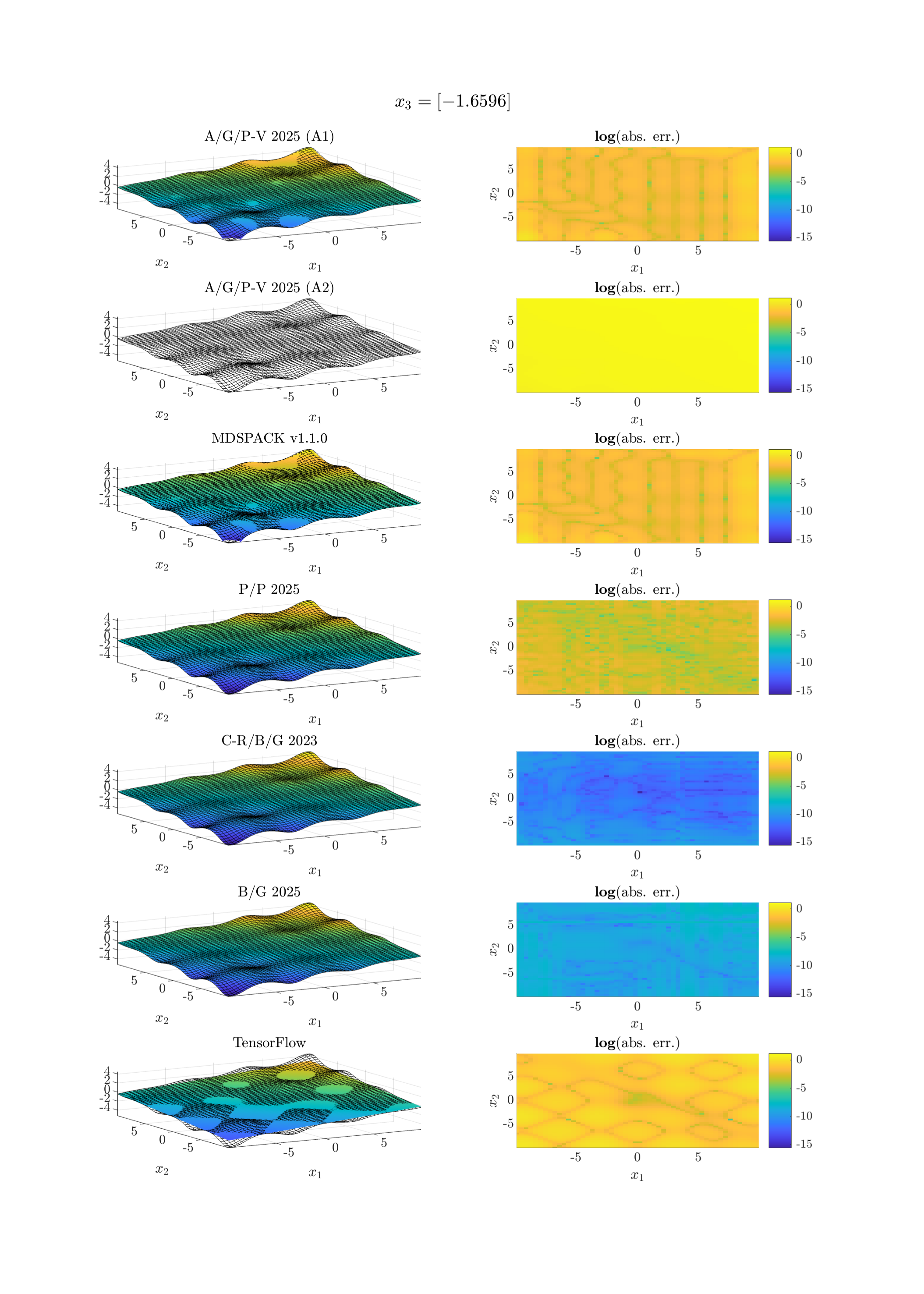} \caption{Function \#26: left side, evaluation of the original (mesh) vs. approximated (coloured surface) and right side, absolute errors (in log-scale).} \end{figure}\subsubsection{mLF detailed informations (M1)} \noindent \textbf{Right interpolation points}: $k_l=\left(\begin{array}{ccc} 10 & 10 & 10 \end{array}\right)$, where $l=1,\cdots,\ord$.$$ \begin{array}{rcl}\lan{1} &\in& \IC^{10} \text{ , linearly spaced between bounds}\\\lan{2} &\in& \IC^{10} \text{ , linearly spaced between bounds}\\\lan{3} &\in& \IC^{10} \text{ , linearly spaced between bounds}\\\end{array} $$\noindent \textbf{$\ord$-D Loewner matrix, barycentric weights and Lagrangian basis}:$$ \begin{array}{rcl}\IL & \in & \IC^{1000 \times 1000}\\\bc & \in & \IC^{1000}\\\bw & \in & \IC^{1000}\\\bc\odot \bw & \in & \IC^{1000}\\\mathbf{Lag}(\var{1},\var{2},\var{3}) & \in & \IC^{1000}\\\end{array} $$

\newpage \subsection{Function \#27 (${\ord=5}$ variables, tensor size: 90.6 \textbf{MB})} $$\frac{\var{1}+\var{2}+\var{3}+\var{4}+\var{5}}{10+\cos(\var{1})+\cos(\var{2})+\cos(\var{3})+\cos(\var{4})+\cos(\var{5})}$$ \subsubsection{Setup and results overview}\begin{itemize}\item Reference: B/G 2025, \cite{Balicki:2025}\item Domain: $\mathbb{R}$\item Tensor size: 90.6 \textbf{MB} ($26^{5}$ points)\item Bounds: $ \left(\begin{array}{cc} -4 & 4 \end{array}\right) \times \left(\begin{array}{cc} -4 & 4 \end{array}\right) \times \left(\begin{array}{cc} -4 & 4 \end{array}\right) \times \left(\begin{array}{cc} -4 & 4 \end{array}\right) \times \left(\begin{array}{cc} -4 & 4 \end{array}\right)$ \end{itemize} \begin{table}[H] \centering \begin{tabular}{llllll}
$\#$ & Alg. & Parameters & Dim. & CPU [s] & RMSE \\ 
\hline 
$\mathbf{\#27}$ & A/G/P-V 2025 (A1) & $1 \cdot 10^{-09},3$ & $4.1 \cdot 10^{05}$ & $3.9$ & $0.003$ \\ 
 & A/G/P-V 2025 (A2) & $1 \cdot 10^{-15},1$ & $\mathbf{2.2 \cdot 10^{02}}$ & $9.6 \cdot 10^{03}$ & $0.44$ \\ 
 & MDSPACK v1.1.0 & $1 \cdot 10^{-10},5$ & $4.1 \cdot 10^{05}$ & $\mathbf{3.6}$ & $\mathbf{0.003}$ \\ 
 & P/P 2025 & $NaN$ & $NaN$ & $NaN$ & $NaN$ \\ 
 & C-R/B/G 2023 & $NaN$ & $NaN$ & $NaN$ & $NaN$ \\ 
 & B/G 2025 & $NaN$ & $NaN$ & $NaN$ & $NaN$ \\ 
 & TensorFlow & $NaN$ & $NaN$ & $NaN$ & $NaN$ \\ 
\hline 
\end{tabular} \caption{Function \#27: best model configuration and performances per methods.} \end{table}\begin{figure}[H] \centering  \includegraphics[width=\textwidth]{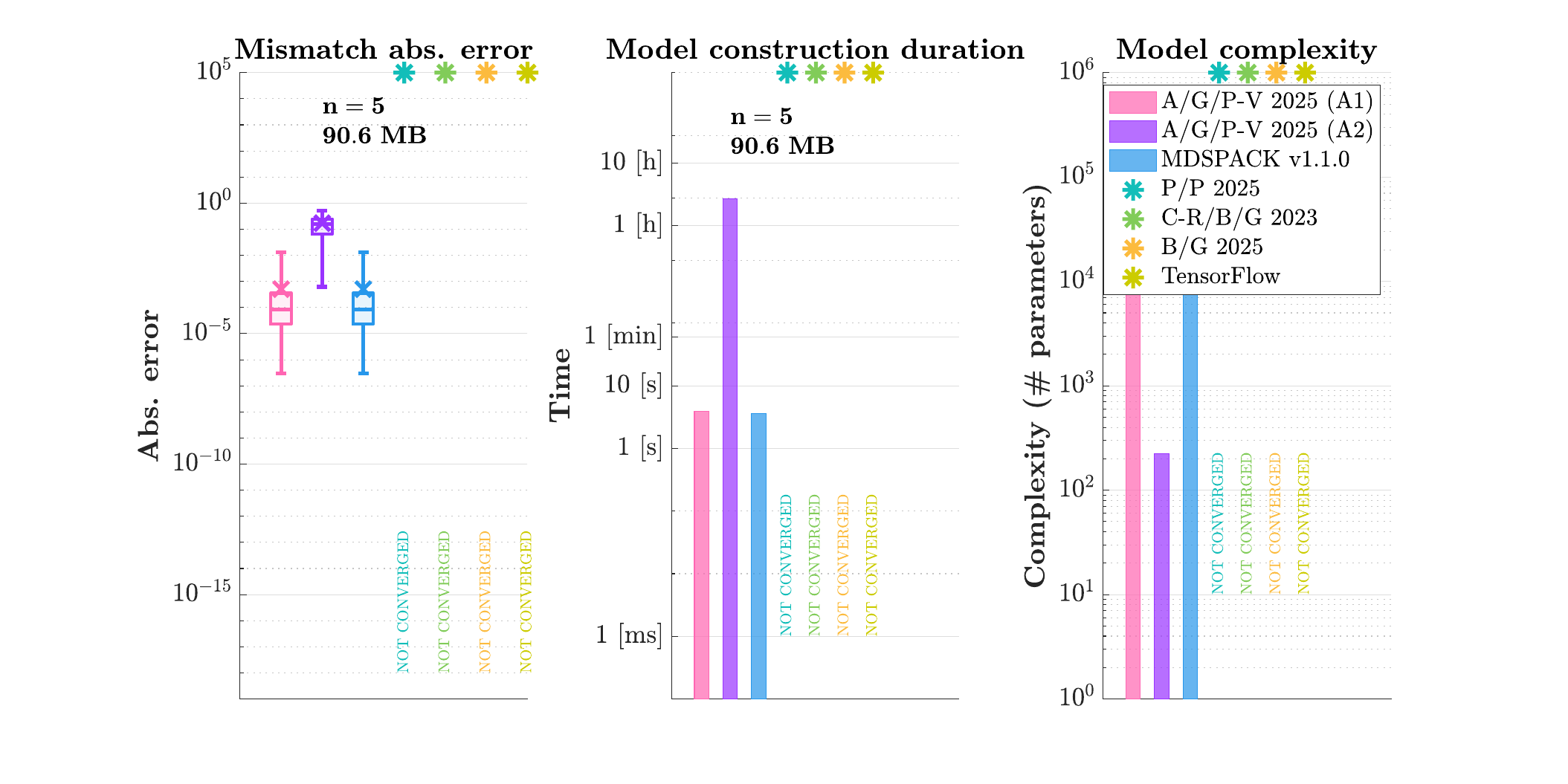} \caption{Function \#27: graphical view of the best model performances.} \end{figure}\begin{figure}[H] \centering  \includegraphics[width=\textwidth]{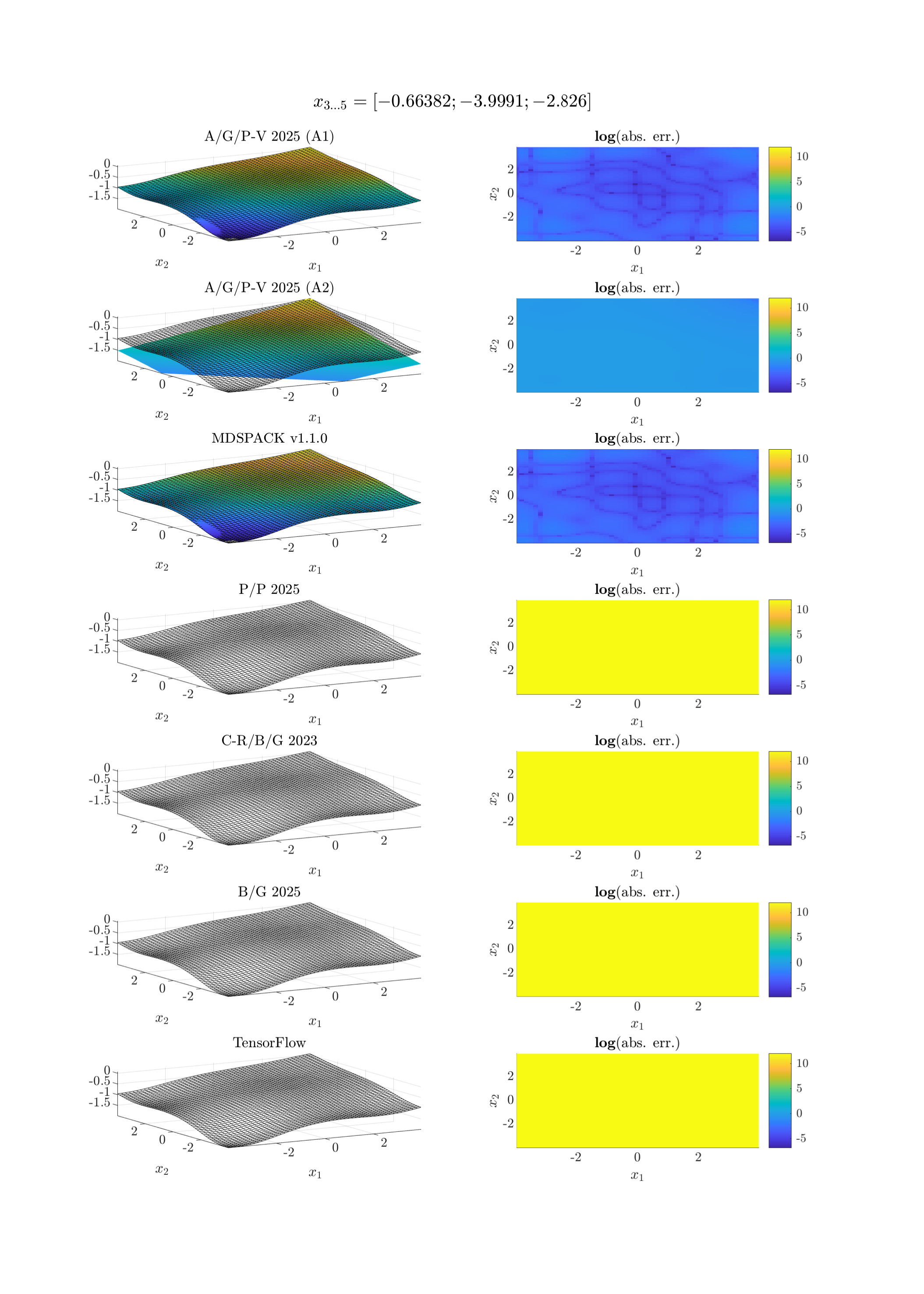} \caption{Function \#27: left side, evaluation of the original (mesh) vs. approximated (coloured surface) and right side, absolute errors (in log-scale).} \end{figure}\subsubsection{mLF detailed informations (M1)} \noindent \textbf{Right interpolation points}: $k_l=\left(\begin{array}{ccccc} 9 & 9 & 9 & 9 & 9 \end{array}\right)$, where $l=1,\cdots,\ord$.$$ \begin{array}{rcl}\lan{1} &\in& \IC^{9} \text{ , linearly spaced between bounds}\\\lan{2} &\in& \IC^{9} \text{ , linearly spaced between bounds}\\\lan{3} &\in& \IC^{9} \text{ , linearly spaced between bounds}\\\lan{4} &\in& \IC^{9} \text{ , linearly spaced between bounds}\\\lan{5} &\in& \IC^{9} \text{ , linearly spaced between bounds}\\\end{array} $$\noindent \textbf{$\ord$-D Loewner matrix, barycentric weights and Lagrangian basis}:$$ \begin{array}{rcl}\IL & \in & \IC^{59049 \times 59049}\\\bc & \in & \IC^{59049}\\\bw & \in & \IC^{59049}\\\bc\odot \bw & \in & \IC^{59049}\\\mathbf{Lag}(\var{1},\var{2},\var{3},\var{4},\var{5}) & \in & \IC^{59049}\\\end{array} $$

\newpage \subsection{Function \#28 (${\ord=2}$ variables, tensor size: 30 \textbf{KB})} $$\left(\frac{\var{1}}{\var{1}+1}\right)^4 (1+\mathrm{exp}(-\var{2}^2)) \left(1+\var{2} \cos(\var{2}) \mathrm{exp}\frac{(-\var{1}\var{2})}{\var{1}+1}\right)$$ \subsubsection{Setup and results overview}\begin{itemize}\item Reference: J/al. 2024 (Toy function), [none]\item Domain: $\mathbb{R}$\item Tensor size: 30 \textbf{KB} ($62^{2}$ points)\item Bounds: $ \left(\begin{array}{cc} \frac{1}{10000000000} & 10 \end{array}\right) \times \left(\begin{array}{cc} \frac{1}{10000000000} & 10 \end{array}\right)$ \end{itemize} \begin{table}[H] \centering \begin{tabular}{llllll}
$\#$ & Alg. & Parameters & Dim. & CPU [s] & RMSE \\ 
\hline 
$\mathbf{\#28}$ & A/G/P-V 2025 (A1) & $1 \cdot 10^{-10},3$ & $5.1 \cdot 10^{02}$ & $0.021$ & $0.2$ \\ 
 & A/G/P-V 2025 (A2) & $1 \cdot 10^{-15},3$ & $3.2 \cdot 10^{02}$ & $0.73$ & $0.0045$ \\ 
 & MDSPACK v1.1.0 & $1 \cdot 10^{-06},3$ & $\mathbf{1.6 \cdot 10^{02}}$ & $\mathbf{0.019}$ & $0.016$ \\ 
 & P/P 2025 & $1,0.95,50,0.01,10,12,21$ & $6.8 \cdot 10^{02}$ & $3.3$ & $0.00024$ \\ 
 & C-R/B/G 2023 & $1 \cdot 10^{-09},20$ & $9.9 \cdot 10^{02}$ & $0.53$ & $\mathbf{4.2 \cdot 10^{-08}}$ \\ 
 & B/G 2025 & $1 \cdot 10^{-09},20,4$ & $1.1 \cdot 10^{03}$ & $2.5$ & $8.9 \cdot 10^{-07}$ \\ 
 & TensorFlow & $$ & $2.6 \cdot 10^{02}$ & $26$ & $0.15$ \\ 
\hline 
\end{tabular} \caption{Function \#28: best model configuration and performances per methods.} \end{table}\begin{figure}[H] \centering  \includegraphics[width=\textwidth]{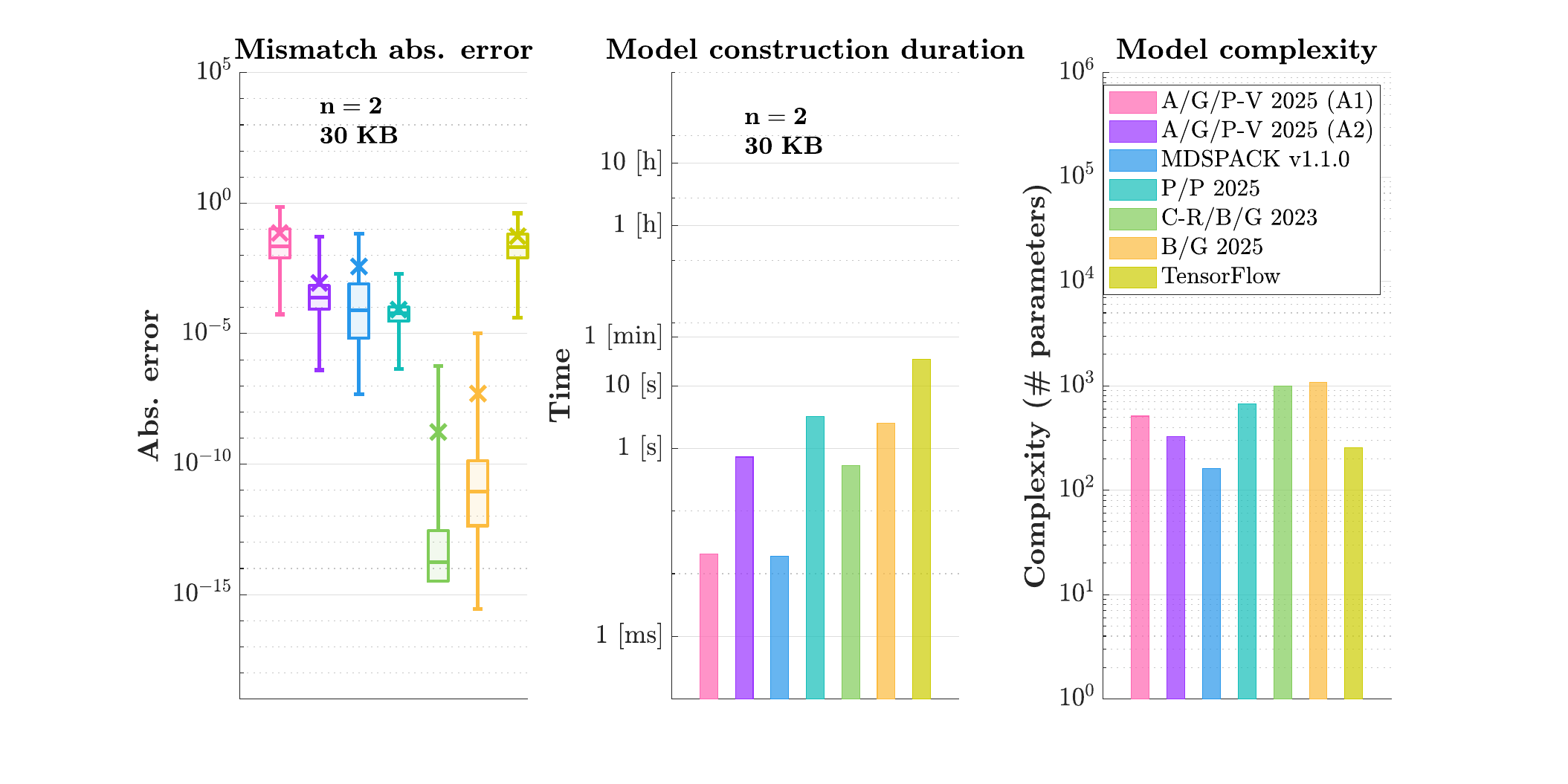} \caption{Function \#28: graphical view of the best model performances.} \end{figure}\begin{figure}[H] \centering  \includegraphics[width=\textwidth]{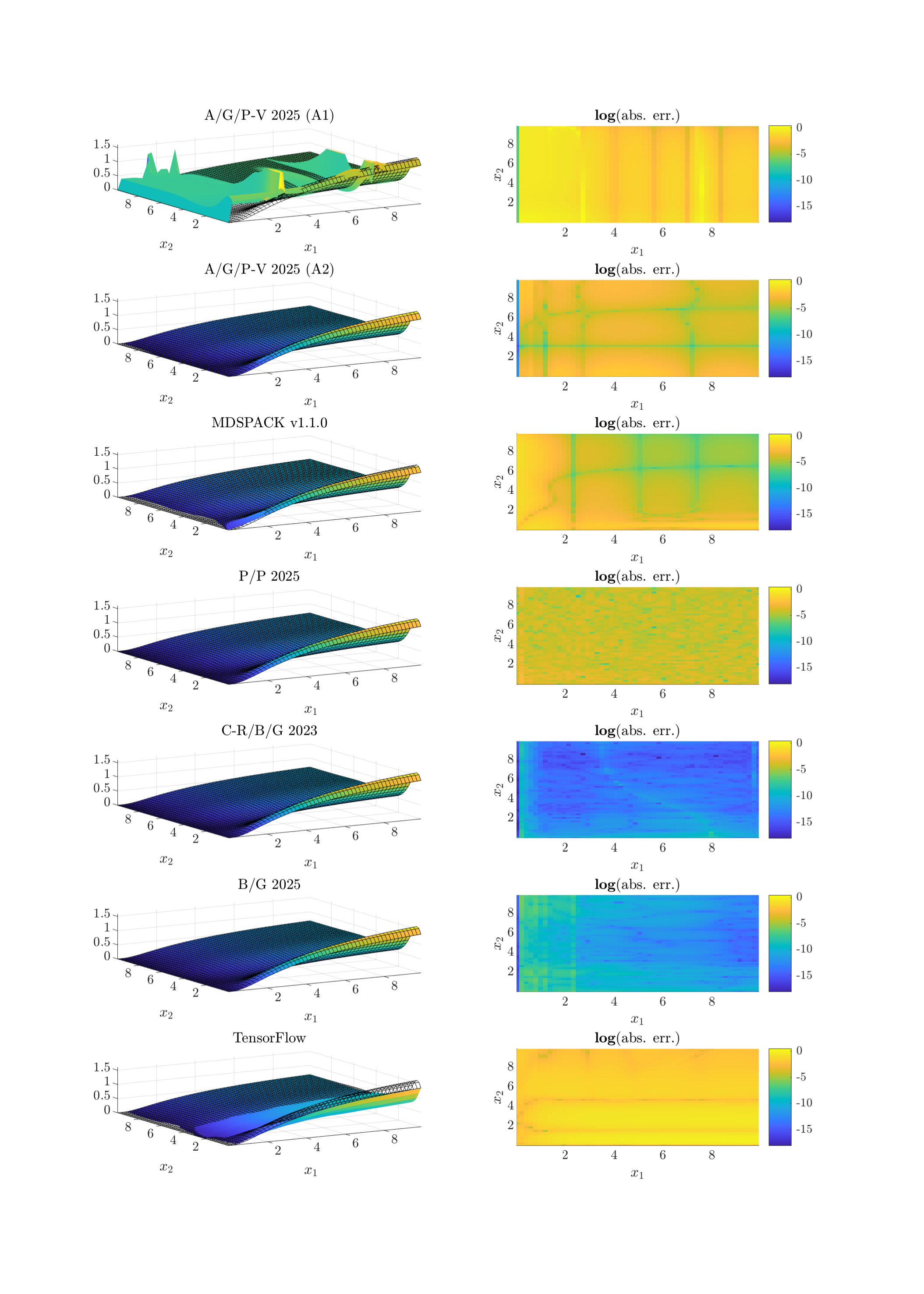} \caption{Function \#28: left side, evaluation of the original (mesh) vs. approximated (coloured surface) and right side, absolute errors (in log-scale).} \end{figure}\subsubsection{mLF detailed informations (M1)} \noindent \textbf{Right interpolation points}: $k_l=\left(\begin{array}{cc} 8 & 16 \end{array}\right)$, where $l=1,\cdots,\ord$.$$ \begin{array}{rcl}\lan{1} &\in& \IC^{8} \text{ , linearly spaced between bounds}\\\lan{2} &\in& \IC^{16} \text{ , linearly spaced between bounds}\\\end{array} $$\noindent \textbf{$\ord$-D Loewner matrix, barycentric weights and Lagrangian basis}:$$ \begin{array}{rcl}\IL & \in & \IC^{128 \times 128}\\\bc & \in & \IC^{128}\\\bw & \in & \IC^{128}\\\bc\odot \bw & \in & \IC^{128}\\\mathbf{Lag}(\var{1},\var{2}) & \in & \IC^{128}\\\end{array} $$

\newpage \subsection{Function \#29 (${\ord=2}$ variables, tensor size: 12.5 \textbf{KB})} $$\min(10|\var{1}|,1)\mathrm{sign}(\var{1}) + \frac{\var{1}\var{2}^3}{10}$$ \subsubsection{Setup and results overview}\begin{itemize}\item Reference: Personal communication, [none]\item Domain: $\mathbb{R}$\item Tensor size: 12.5 \textbf{KB} ($40^{2}$ points)\item Bounds: $ \left(\begin{array}{cc} -1 & 1 \end{array}\right) \times \left(\begin{array}{cc} -1 & 1 \end{array}\right)$ \end{itemize} \begin{table}[H] \centering \begin{tabular}{llllll}
$\#$ & Alg. & Parameters & Dim. & CPU [s] & RMSE \\ 
\hline 
$\mathbf{\#29}$ & A/G/P-V 2025 (A1) & $0.01,2$ & $\mathbf{1.1 \cdot 10^{02}}$ & $0.01$ & $0.038$ \\ 
 & A/G/P-V 2025 (A2) & $1 \cdot 10^{-15},3$ & $3.2 \cdot 10^{02}$ & $0.45$ & $0.2$ \\ 
 & MDSPACK v1.1.0 & $0.0001,2$ & $1.1 \cdot 10^{02}$ & $\mathbf{0.01}$ & $0.038$ \\ 
 & P/P 2025 & $4,1,50,0.01,4,6,9$ & $1.3 \cdot 10^{02}$ & $0.46$ & $\mathbf{0.017}$ \\ 
 & C-R/B/G 2023 & $0.001,20$ & $6.6 \cdot 10^{02}$ & $0.19$ & $1.4$ \\ 
 & B/G 2025 & $0.001,20,3$ & $5.7 \cdot 10^{02}$ & $0.12$ & $0.49$ \\ 
 & TensorFlow & $NaN$ & $NaN$ & $NaN$ & $NaN$ \\ 
\hline 
\end{tabular} \caption{Function \#29: best model configuration and performances per methods.} \end{table}\begin{figure}[H] \centering  \includegraphics[width=\textwidth]{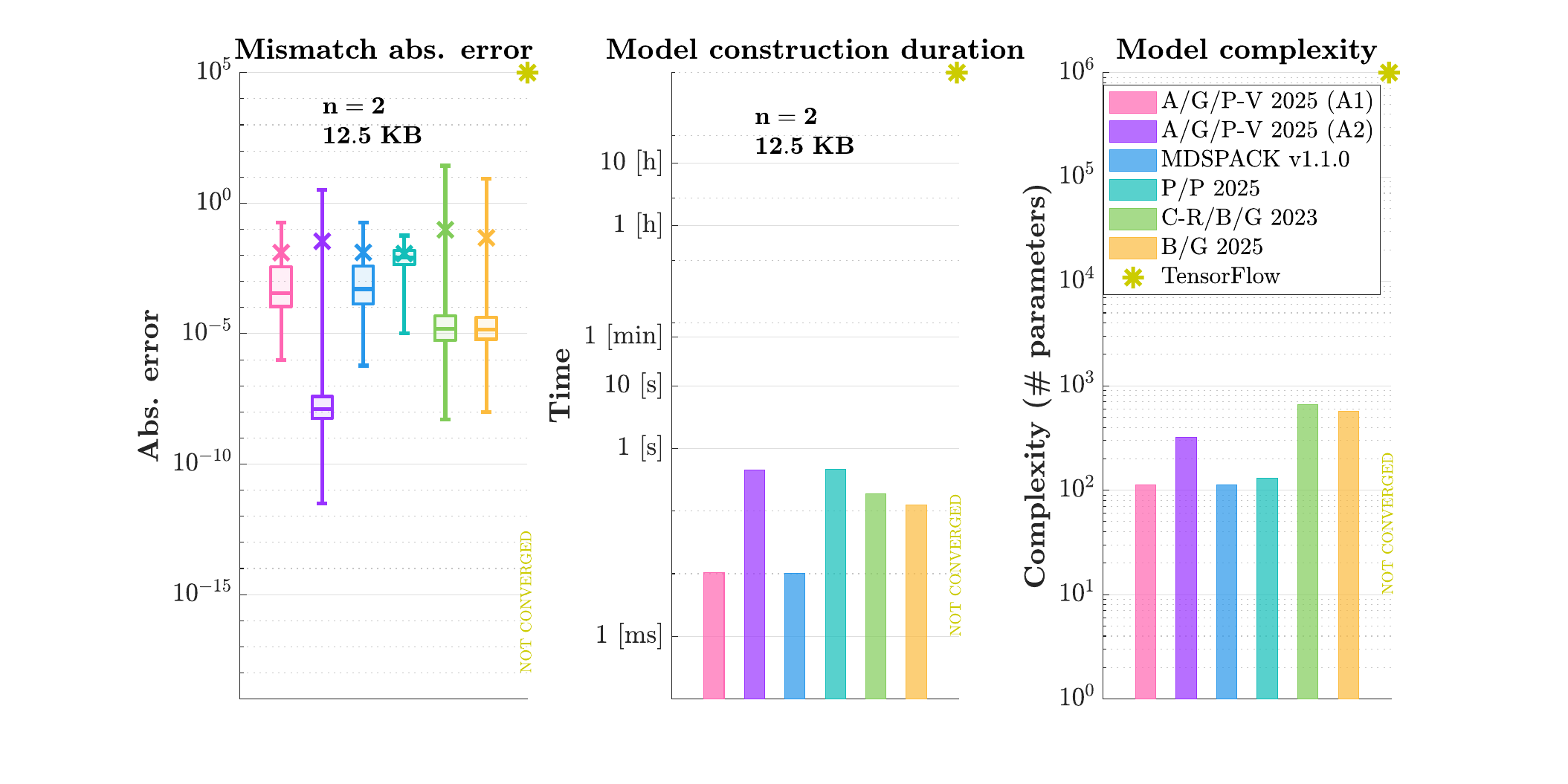} \caption{Function \#29: graphical view of the best model performances.} \end{figure}\begin{figure}[H] \centering  \includegraphics[width=\textwidth]{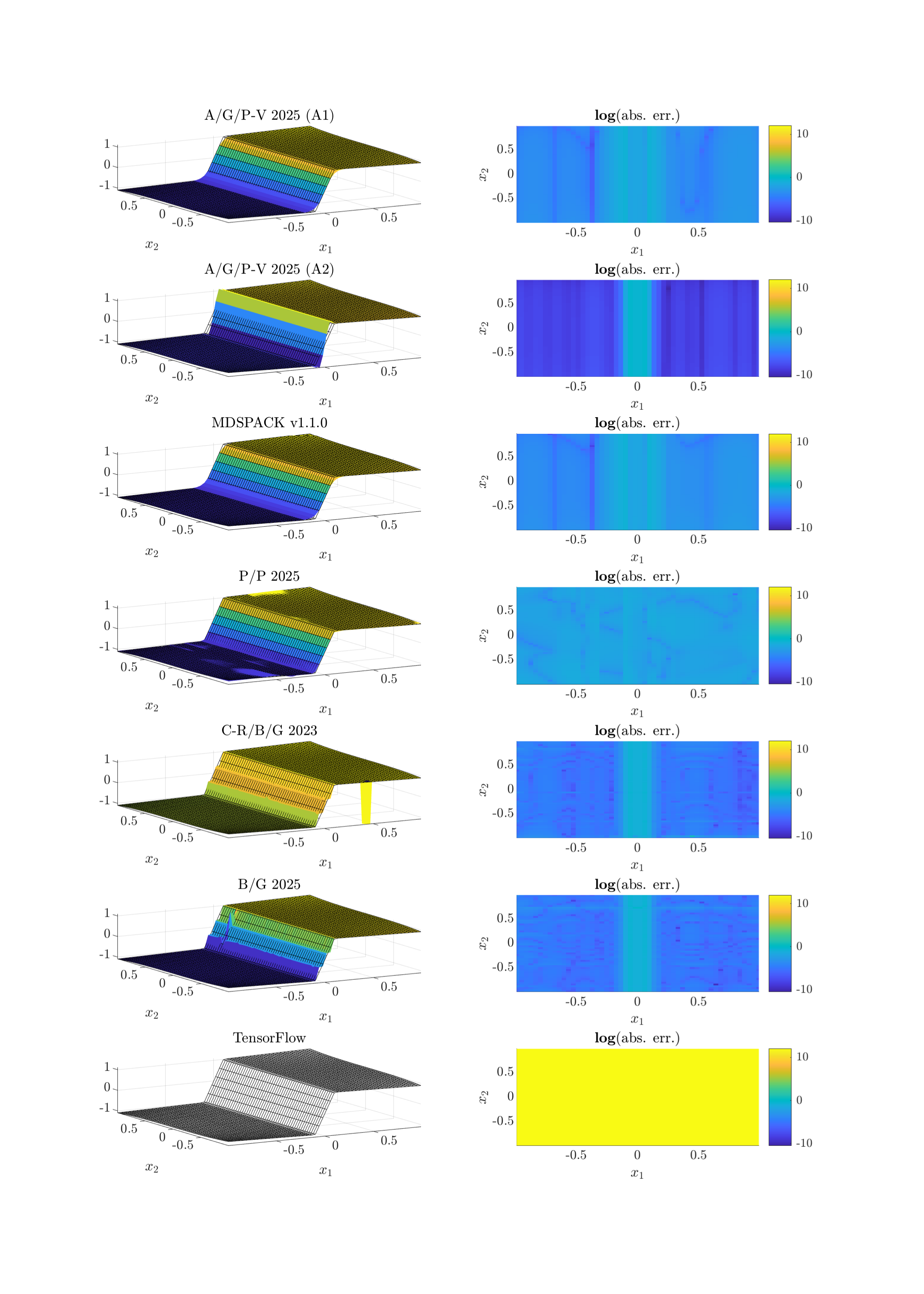} \caption{Function \#29: left side, evaluation of the original (mesh) vs. approximated (coloured surface) and right side, absolute errors (in log-scale).} \end{figure}\subsubsection{mLF detailed informations (M1)} \noindent \textbf{Right interpolation points} ($k_l=\left(\begin{array}{cc} 7 & 4 \end{array}\right)$, where $l=1,\cdots,\ord$):$$ \begin{array}{rcl}\lan{1} &=& \left(\begin{array}{ccccccc} -1 & -\frac{13}{19} & -\frac{7}{19} & -\frac{1}{19} & \frac{5}{19} & \frac{11}{19} & 1 \end{array}\right)\\\lan{2} &=& \left(\begin{array}{cccc} -1 & -\frac{7}{19} & \frac{5}{19} & 1 \end{array}\right)\\\end{array} $$\noindent \textbf{Lagrangian weights}: $$\left(\begin{array}{ccc} \bc & \bw & \bc\odot\bw\\ 0.689 & -0.9 & -0.6201\\ -2.014 & -0.995 & 2.004\\ 1.87 & -1.002 & -1.874\\ -0.5452 & -1.1 & 0.5997\\ -1.264 & -0.9316 & 1.177\\ 3.694 & -0.9966 & -3.682\\ -3.431 & -1.001 & 3.435\\ 1.0 & -1.068 & -1.068\\ 0.6523 & -0.9632 & -0.6283\\ -1.907 & -0.9982 & 1.903\\ 1.771 & -1.001 & -1.772\\ -0.5161 & -1.037 & 0.5351\\ -0.1037 & -0.5211 & 0.05401\\ 0.303 & -0.5261 & -0.1594\\ -0.2813 & -0.5264 & 0.1481\\ 0.08201 & -0.5316 & -0.04359\\ 0.3185 & 0.9737 & 0.3102\\ -0.9311 & 0.9987 & -0.9299\\ 0.8646 & 1.0 & 0.865\\ -0.252 & 1.026 & -0.2587\\ -0.7253 & 0.9421 & -0.6833\\ 2.12 & 0.9971 & 2.114\\ -1.969 & 1.001 & -1.971\\ 0.5739 & 1.058 & 0.6071\\ 0.4635 & 0.9 & 0.4171\\ -1.355 & 0.995 & -1.348\\ 1.258 & 1.002 & 1.26\\ -0.3667 & 1.1 & -0.4034 \end{array}\right)$$\noindent \textbf{Lagrangian form} (basis, numerator and denominator coefficients):$$\left(\begin{array}{ccc}\mathcal{B}_\textrm{lag}(\var{1},\var{2}) & \bN_\textrm{lag} &\bD_\textrm{lag}\end{array}\right) =$$ $$\left(\begin{array}{ccc} \left(\var{1}+1.0\right)\,\left(\var{2}+1.0\right) & -0.6201 & 0.689\\ \left(\var{1}+1.0\right)\,\left(\var{2}+0.3684\right) & 2.004 & -2.014\\ \left(\var{1}+1.0\right)\,\left(\var{2}-0.2632\right) & -1.874 & 1.87\\ \left(\var{1}+1.0\right)\,\left(\var{2}-1.0\right) & 0.5997 & -0.5452\\ \left(\var{2}+1.0\right)\,\left(\var{1}+0.6842\right) & 1.177 & -1.264\\ \left(\var{1}+0.6842\right)\,\left(\var{2}+0.3684\right) & -3.682 & 3.694\\ \left(\var{2}-0.2632\right)\,\left(\var{1}+0.6842\right) & 3.435 & -3.431\\ \left(\var{2}-1.0\right)\,\left(\var{1}+0.6842\right) & -1.068 & 1.0\\ \left(\var{2}+1.0\right)\,\left(\var{1}+0.3684\right) & -0.6283 & 0.6523\\ \left(\var{1}+0.3684\right)\,\left(\var{2}+0.3684\right) & 1.903 & -1.907\\ \left(\var{2}-0.2632\right)\,\left(\var{1}+0.3684\right) & -1.772 & 1.771\\ \left(\var{2}-1.0\right)\,\left(\var{1}+0.3684\right) & 0.5351 & -0.5161\\ \left(\var{2}+1.0\right)\,\left(\var{1}+0.05263\right) & 0.05401 & -0.1037\\ \left(\var{2}+0.3684\right)\,\left(\var{1}+0.05263\right) & -0.1594 & 0.303\\ \left(\var{2}-0.2632\right)\,\left(\var{1}+0.05263\right) & 0.1481 & -0.2813\\ \left(\var{2}-1.0\right)\,\left(\var{1}+0.05263\right) & -0.04359 & 0.08201\\ \left(\var{2}+1.0\right)\,\left(\var{1}-0.2632\right) & 0.3102 & 0.3185\\ \left(\var{1}-0.2632\right)\,\left(\var{2}+0.3684\right) & -0.9299 & -0.9311\\ \left(\var{1}-0.2632\right)\,\left(\var{2}-0.2632\right) & 0.865 & 0.8646\\ \left(\var{2}-1.0\right)\,\left(\var{1}-0.2632\right) & -0.2587 & -0.252\\ \left(\var{2}+1.0\right)\,\left(\var{1}-0.5789\right) & -0.6833 & -0.7253\\ \left(\var{1}-0.5789\right)\,\left(\var{2}+0.3684\right) & 2.114 & 2.12\\ \left(\var{2}-0.2632\right)\,\left(\var{1}-0.5789\right) & -1.971 & -1.969\\ \left(\var{2}-1.0\right)\,\left(\var{1}-0.5789\right) & 0.6071 & 0.5739\\ \left(\var{1}-1.0\right)\,\left(\var{2}+1.0\right) & 0.4171 & 0.4635\\ \left(\var{1}-1.0\right)\,\left(\var{2}+0.3684\right) & -1.348 & -1.355\\ \left(\var{1}-1.0\right)\,\left(\var{2}-0.2632\right) & 1.26 & 1.258\\ \left(\var{1}-1.0\right)\,\left(\var{2}-1.0\right) & -0.4034 & -0.3667 \end{array}\right).$$\noindent The corresponding function is:$$\begin{array}{rcl}\bG_{\textrm{lag}}(\var{1},\var{2}) &=& \dfrac{\bn_{\textrm{lag}}(\var{1},\var{2})}{\bd_{\textrm{lag}}(\var{1},\var{2})}\\ && \\&=& \dfrac{\sum_{\textrm{row}} \bN_\textrm{lag} \odot\mathcal{B}^{-1}_\textrm{lag}(\var{1},\var{2})}{\sum_{\textrm{row}} \bD_\textrm{lag} \odot\mathcal{B}^{-1}_\textrm{lag}(\var{1},\var{2})}, \end{array}$$\noindent where,\\$\bn_{\textrm{lag}}(\var{1},\var{2}) = 1.898\,{\var{1}}^6\,{\var{2}}^3-4.99 \cdot 10^{-12}\,{\var{1}}^6\,{\var{2}}^2-3.382 \cdot 10^{-11}\,{\var{1}}^6\,\var{2}+9.386\,{\var{1}}^6+12.41\,{\var{1}}^5\,{\var{2}}^3+3.944 \cdot 10^{-12}\,{\var{1}}^5\,{\var{2}}^2-2.844 \cdot 10^{-11}\,{\var{1}}^5\,\var{2}+55.08\,{\var{1}}^5+0.6682\,{\var{1}}^4\,{\var{2}}^3+5.471 \cdot 10^{-12}\,{\var{1}}^4\,{\var{2}}^2+3.668 \cdot 10^{-11}\,{\var{1}}^4\,\var{2}+15.8\,{\var{1}}^4+4.392\,{\var{1}}^3\,{\var{2}}^3-2.102 \cdot 10^{-12}\,{\var{1}}^3\,{\var{2}}^2+3.218 \cdot 10^{-11}\,{\var{1}}^3\,\var{2}+103.1\,{\var{1}}^3-0.08665\,{\var{1}}^2\,{\var{2}}^3-2.198 \cdot 10^{-13}\,{\var{1}}^2\,{\var{2}}^2-2.554 \cdot 10^{-12}\,{\var{1}}^2\,\var{2}-0.376\,{\var{1}}^2+0.1303\,\var{1}\,{\var{2}}^3+7.95 \cdot 10^{-14}\,\var{1}\,{\var{2}}^2-2.367 \cdot 10^{-12}\,\var{1}\,\var{2}+11.1\,\var{1}+0.001715\,{\var{2}}^3+2.372 \cdot 10^{-14}\,{\var{2}}^2-1.099 \cdot 10^{-13}\,\var{2}-0.001715$ \\~~\\$\bd_{\textrm{lag}}(\var{1},\var{2}) = 5.204 \cdot 10^{-12}\,{\var{1}}^6\,{\var{2}}^3+5.974 \cdot 10^{-12}\,{\var{1}}^6\,{\var{2}}^2-2.268 \cdot 10^{-11}\,{\var{1}}^6\,\var{2}+8.484\,{\var{1}}^6+1.639 \cdot 10^{-11}\,{\var{1}}^5\,{\var{2}}^3-1.575 \cdot 10^{-12}\,{\var{1}}^5\,{\var{2}}^2-2.896 \cdot 10^{-11}\,{\var{1}}^5\,\var{2}+21.21\,{\var{1}}^5-2.658 \cdot 10^{-12}\,{\var{1}}^4\,{\var{2}}^3-4.615 \cdot 10^{-12}\,{\var{1}}^4\,{\var{2}}^2+2.094 \cdot 10^{-11}\,{\var{1}}^4\,\var{2}+111.6\,{\var{1}}^4-1.847 \cdot 10^{-11}\,{\var{1}}^3\,{\var{2}}^3+1.8 \cdot 10^{-12}\,{\var{1}}^3\,{\var{2}}^2+3.253 \cdot 10^{-11}\,{\var{1}}^3\,\var{2}+3.794\,{\var{1}}^3-3.354 \cdot 10^{-12}\,{\var{1}}^2\,{\var{2}}^3+4.967 \cdot 10^{-13}\,{\var{1}}^2\,{\var{2}}^2+3.343 \cdot 10^{-12}\,{\var{1}}^2\,\var{2}+48.22\,{\var{1}}^2+1.988 \cdot 10^{-12}\,\var{1}\,{\var{2}}^3+5.146 \cdot 10^{-14}\,\var{1}\,{\var{2}}^2-3.375 \cdot 10^{-12}\,\var{1}\,\var{2}-0.1938\,\var{1}+1.667 \cdot 10^{-13}\,{\var{2}}^3-3.565 \cdot 10^{-14}\,{\var{2}}^2-1.883 \cdot 10^{-13}\,\var{2}+1.0$ \\~~\\\noindent \textbf{Monomial form} (basis, numerator and denominator coefficients - entries $<10^{-12}$ removed):$$\left(\begin{array}{ccc}\mathcal{B}_\textrm{mon}(\var{1},\var{2}) & \bN_\textrm{mon} &\bD_\textrm{mon}\end{array}\right) =$$ $$\left(\begin{array}{ccc} {\var{1}}^6\,{\var{2}}^3 & -0.017 & 0\\ {\var{1}}^6\,{\var{2}}^2 & 0 & 0\\ {\var{1}}^6\,\var{2} & 0 & 0\\ {\var{1}}^6 & -0.0841 & -0.07601\\ {\var{1}}^5\,{\var{2}}^3 & -0.1112 & 0\\ {\var{1}}^5\,{\var{2}}^2 & 0 & 0\\ {\var{1}}^5\,\var{2} & 0 & 0\\ {\var{1}}^5 & -0.4935 & -0.19\\ {\var{1}}^4\,{\var{2}}^3 & -0.005987 & 0\\ {\var{1}}^4\,{\var{2}}^2 & 0 & 0\\ {\var{1}}^4\,\var{2} & 0 & 0\\ {\var{1}}^4 & -0.1416 & -1.0\\ {\var{1}}^3\,{\var{2}}^3 & -0.03935 & 0\\ {\var{1}}^3\,{\var{2}}^2 & 0 & 0\\ {\var{1}}^3\,\var{2} & 0 & 0\\ {\var{1}}^3 & -0.924 & -0.03399\\ {\var{1}}^2\,{\var{2}}^3 & 0.0007764 & 0\\ {\var{1}}^2\,{\var{2}}^2 & 0 & 0\\ {\var{1}}^2\,\var{2} & 0 & 0\\ {\var{1}}^2 & 0.003368 & -0.432\\ \var{1}\,{\var{2}}^3 & -0.001167 & 0\\ \var{1}\,{\var{2}}^2 & 0 & 0\\ \var{1}\,\var{2} & 0 & 0\\ \var{1} & -0.0995 & 0.001736\\ {\var{2}}^3 & -1.536 \cdot 10^{-5} & 0\\ {\var{2}}^2 & 0 & 0\\ \var{2} & 0 & 0\\ 1.0 & 1.536 \cdot 10^{-5} & -0.00896 \end{array}\right)$$\noindent The corresponding function is:$$\begin{array}{rcl}\bG_{\textrm{mon}}(\var{1},\var{2}) &=& \dfrac{\bn_{\textrm{mon}}(\var{1},\var{2})}{\bd_{\textrm{mon}}(\var{1},\var{2})}\\ && \\&=& \dfrac{\sum_{\textrm{row}} \bN_\textrm{mon} \odot \mathcal{B}_\textrm{mon}(\var{1},\var{2})}{\sum_{\textrm{row}} \bD_\textrm{mon} \odot\mathcal{B}_\textrm{mon}(\var{1},\var{2})},  \end{array}$$\noindent where,\\$\bn_{\textrm{mon}}(\var{1},\var{2}) = 1.898\,{\var{1}}^6\,{\var{2}}^3+9.386\,{\var{1}}^6+12.41\,{\var{1}}^5\,{\var{2}}^3+55.08\,{\var{1}}^5+0.6682\,{\var{1}}^4\,{\var{2}}^3+15.8\,{\var{1}}^4+4.392\,{\var{1}}^3\,{\var{2}}^3+103.1\,{\var{1}}^3-0.08665\,{\var{1}}^2\,{\var{2}}^3-0.376\,{\var{1}}^2+0.1303\,\var{1}\,{\var{2}}^3+11.1\,\var{1}+0.001715\,{\var{2}}^3-0.001715$ \\~~\\$\bd_{\textrm{mon}}(\var{1},\var{2}) = 8.484\,{\var{1}}^6+21.21\,{\var{1}}^5+111.6\,{\var{1}}^4+3.794\,{\var{1}}^3+48.22\,{\var{1}}^2-0.1938\,\var{1}+1.0$ \\~~\\\noindent \textbf{KST equivalent decoupling pattern} (Barycentric weights $\bc^{\var{l}}$): $$\begin{array}{rclll}\var{2}&: & \left(\begin{array}{ccccccc} -1.264 & -1.264 & -1.264 & -1.264 & -1.264 & -1.264 & -1.264\\ 3.694 & 3.694 & 3.694 & 3.694 & 3.694 & 3.694 & 3.694\\ -3.431 & -3.431 & -3.431 & -3.431 & -3.431 & -3.431 & -3.431\\ 1.0 & 1.0 & 1.0 & 1.0 & 1.0 & 1.0 & 1.0 \end{array}\right)& \textrm{vec}(.) & := \textbf{Bary}(\var{2}) \\\var{1}&: & \left(\begin{array}{c} -0.5452\\ 1.0\\ -0.5161\\ 0.08201\\ -0.252\\ 0.5739\\ -0.3667 \end{array}\right)& \textrm{vec}(.) \otimes \bone_{k_{2}} & := \textbf{Bary}(\var{1}) \\\end{array}$$~\\ Then, with the above notations, one defines the following univariate vector functions:  $$ \left\{ \begin{array}{rcl}\bPhi_{1}(\var{1}) &:=& \textbf{Bary}(\var{1}) \odot \mathbf{Lag}(\var{1}) \\\bPhi_{2}(\var{2}) &:=& \textbf{Bary}(\var{2}) \odot \mathbf{Lag}(\var{2}) \\\end{array} \right. $$\noindent The corresponding function is:$$\begin{array}{rcl}\bG_{\textrm{kst}}(\var{1},\var{2}) &=& \dfrac{\bn_{\textrm{kst}}(\var{1},\var{2})}{\bd_{\textrm{kst}}(\var{1},\var{2})}\\ && \\ &=& \dfrac{\sum_{\text{rows}} \bw \odot \bPhi_{1}(\var{1}) \odot \cdots \odot\bPhi_{2}(\var{2})}{\sum_{\text{rows}} \bPhi_{1}(\var{1}) \odot \cdots \odot\bPhi_{2}(\var{2})} . \end{array}$$~\\ \noindent \textbf{KST-like univariate functions} (equivalent scaled univariate functions $\bphi_{1,\cdots,2}$): $$\left\{\begin{array}{rcrcl}z_{1} &=&\bphi_{1}(\var{1}) &=& \cfrac{\bn_{1}}{\bd_{1}} \\z_{2} &=&\bphi_{2}(\var{2}) &=& -0.1\,{\var{2}}^3-1.0\\\end{array} \right. .$$\noindent where, \\ \noindent $\bn_{1}=11.28\,{\var{1}}^6+67.49\,{\var{1}}^5+16.47\,{\var{1}}^4+107.5\,{\var{1}}^3-0.4626\,{\var{1}}^2+11.23\,\var{1}-1.001 \cdot 10^{-15}$ and \\ \noindent $\bd_{1}=8.484\,{\var{1}}^6+21.21\,{\var{1}}^5+111.6\,{\var{1}}^4+3.794\,{\var{1}}^3+48.22\,{\var{1}}^2-0.1938\,\var{1}+1.0$, \\

\newpage \subsection{Function \#30 (${\ord=8}$ variables, tensor size: 128 \textbf{MB})} $$\texttt{Borehole function}$$ \subsubsection{Setup and results overview}\begin{itemize}\item Reference: Borehole function, \cite{Surjanovic}\item Domain: $\mathbb{R}$\item Tensor size: 128 \textbf{MB} ($8^{8}$ points)\item Bounds: $ \left(\begin{array}{cc} \frac{1}{20} & \frac{3}{20} \end{array}\right) \times \left(\begin{array}{cc} 100 & 50000 \end{array}\right) \times \left(\begin{array}{cc} 63070 & 115600 \end{array}\right) \times \left(\begin{array}{cc} 990 & 1110 \end{array}\right) \times \left(\begin{array}{cc} \frac{631}{10} & 116 \end{array}\right) \times \left(\begin{array}{cc} 700 & 820 \end{array}\right) \times \left(\begin{array}{cc} 1120 & 1680 \end{array}\right) \times \left(\begin{array}{cc} 9855 & 12045 \end{array}\right)$ \end{itemize} \begin{table}[H] \centering \begin{tabular}{llllll}
$\#$ & Alg. & Parameters & Dim. & CPU [s] & RMSE \\ 
\hline 
$\mathbf{\#30}$ & A/G/P-V 2025 (A1) & $1 \cdot 10^{-09},1$ & $\mathbf{1 \cdot 10^{04}}$ & $18$ & $0.0025$ \\ 
 & A/G/P-V 2025 (A2) & $1 \cdot 10^{-15},1$ & $NaN$ & $NaN$ & $NaN$ \\ 
 & MDSPACK v1.1.0 & $1 \cdot 10^{-10},5$ & $1 \cdot 10^{04}$ & $\mathbf{18}$ & $\mathbf{0.0025}$ \\ 
 & P/P 2025 & $NaN$ & $NaN$ & $NaN$ & $NaN$ \\ 
 & C-R/B/G 2023 & $NaN$ & $NaN$ & $NaN$ & $NaN$ \\ 
 & B/G 2025 & $NaN$ & $NaN$ & $NaN$ & $NaN$ \\ 
 & TensorFlow & $NaN$ & $NaN$ & $NaN$ & $NaN$ \\ 
\hline 
\end{tabular} \caption{Function \#30: best model configuration and performances per methods.} \end{table}\begin{figure}[H] \centering  \includegraphics[width=\textwidth]{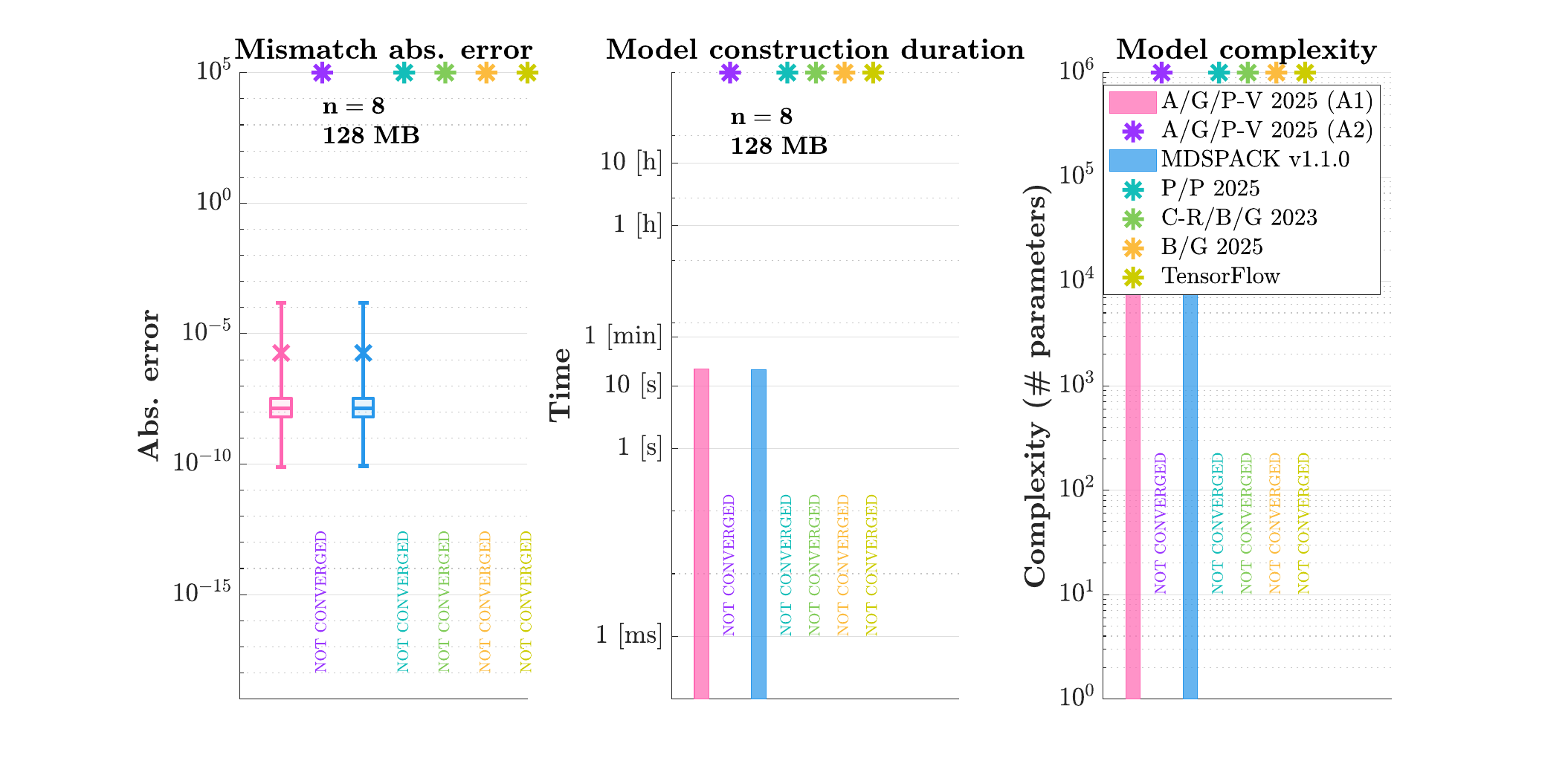} \caption{Function \#30: graphical view of the best model performances.} \end{figure}\begin{figure}[H] \centering  \includegraphics[width=\textwidth]{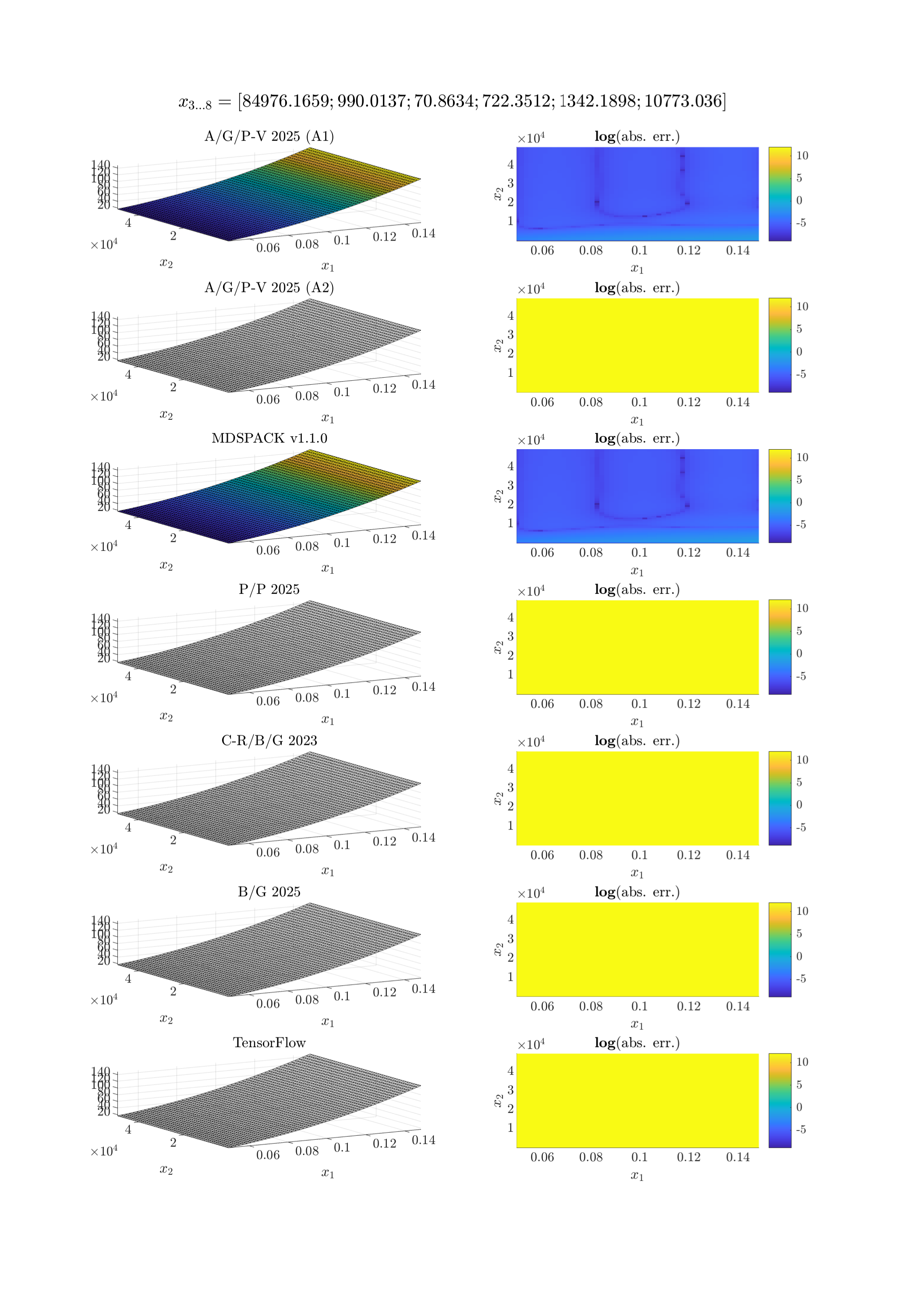} \caption{Function \#30: left side, evaluation of the original (mesh) vs. approximated (coloured surface) and right side, absolute errors (in log-scale).} \end{figure}\subsubsection{mLF detailed informations (M1)} \noindent \textbf{Right interpolation points}: $k_l=\left(\begin{array}{cccccccc} 4 & 4 & 2 & 2 & 2 & 2 & 2 & 2 \end{array}\right)$, where $l=1,\cdots,\ord$.$$ \begin{array}{rcl}\lan{1} &\in& \IC^{4} \text{ , linearly spaced between bounds}\\\lan{2} &\in& \IC^{4} \text{ , linearly spaced between bounds}\\\lan{3} &\in& \IC^{2} \text{ , linearly spaced between bounds}\\\lan{4} &\in& \IC^{2} \text{ , linearly spaced between bounds}\\\lan{5} &\in& \IC^{2} \text{ , linearly spaced between bounds}\\\lan{6} &\in& \IC^{2} \text{ , linearly spaced between bounds}\\\lan{7} &\in& \IC^{2} \text{ , linearly spaced between bounds}\\\lan{8} &\in& \IC^{2} \text{ , linearly spaced between bounds}\\\end{array} $$\noindent \textbf{$\ord$-D Loewner matrix, barycentric weights and Lagrangian basis}:$$ \begin{array}{rcl}\IL & \in & \IC^{1024 \times 1024}\\\bc & \in & \IC^{1024}\\\bw & \in & \IC^{1024}\\\bc\odot \bw & \in & \IC^{1024}\\\mathbf{Lag}(\var{1},\var{2},\var{3},\var{4},\var{5},\var{6},\var{7},\var{8}) & \in & \IC^{1024}\\\end{array} $$

\newpage \subsection{Function \#31 (${\ord=6}$ variables, tensor size: 128 \textbf{MB})} $$\var{1}^2 \var{2}^3 \var{3} \var{4} - \var{5}^2 + \var{6}$$ \subsubsection{Setup and results overview}\begin{itemize}\item Reference: Personal communication, [none]\item Domain: $\mathbb{R}$\item Tensor size: 128 \textbf{MB} ($16^{6}$ points)\item Bounds: $ \left(\begin{array}{cc} -2 & 2 \end{array}\right) \times \left(\begin{array}{cc} -2 & 2 \end{array}\right) \times \left(\begin{array}{cc} -2 & 2 \end{array}\right) \times \left(\begin{array}{cc} -2 & 2 \end{array}\right) \times \left(\begin{array}{cc} -2 & 2 \end{array}\right) \times \left(\begin{array}{cc} -2 & 2 \end{array}\right)$ \end{itemize} \begin{table}[H] \centering \begin{tabular}{llllll}
$\#$ & Alg. & Parameters & Dim. & CPU [s] & RMSE \\ 
\hline 
$\mathbf{\#31}$ & A/G/P-V 2025 (A1) & $0.1,3$ & $\mathbf{2.3 \cdot 10^{03}}$ & $\mathbf{7.8}$ & $\mathbf{1.3 \cdot 10^{-14}}$ \\ 
 & A/G/P-V 2025 (A2) & $1 \cdot 10^{-15},2$ & $2.3 \cdot 10^{03}$ & $1.6 \cdot 10^{03}$ & $7.5 \cdot 10^{-13}$ \\ 
 & MDSPACK v1.1.0 & $0.01,1$ & $2.3 \cdot 10^{03}$ & $8.7$ & $1.8 \cdot 10^{-14}$ \\ 
 & P/P 2025 & $NaN$ & $NaN$ & $NaN$ & $NaN$ \\ 
 & C-R/B/G 2023 & $NaN$ & $NaN$ & $NaN$ & $NaN$ \\ 
 & B/G 2025 & $NaN$ & $NaN$ & $NaN$ & $NaN$ \\ 
 & TensorFlow & $NaN$ & $NaN$ & $NaN$ & $NaN$ \\ 
\hline 
\end{tabular} \caption{Function \#31: best model configuration and performances per methods.} \end{table}\begin{figure}[H] \centering  \includegraphics[width=\textwidth]{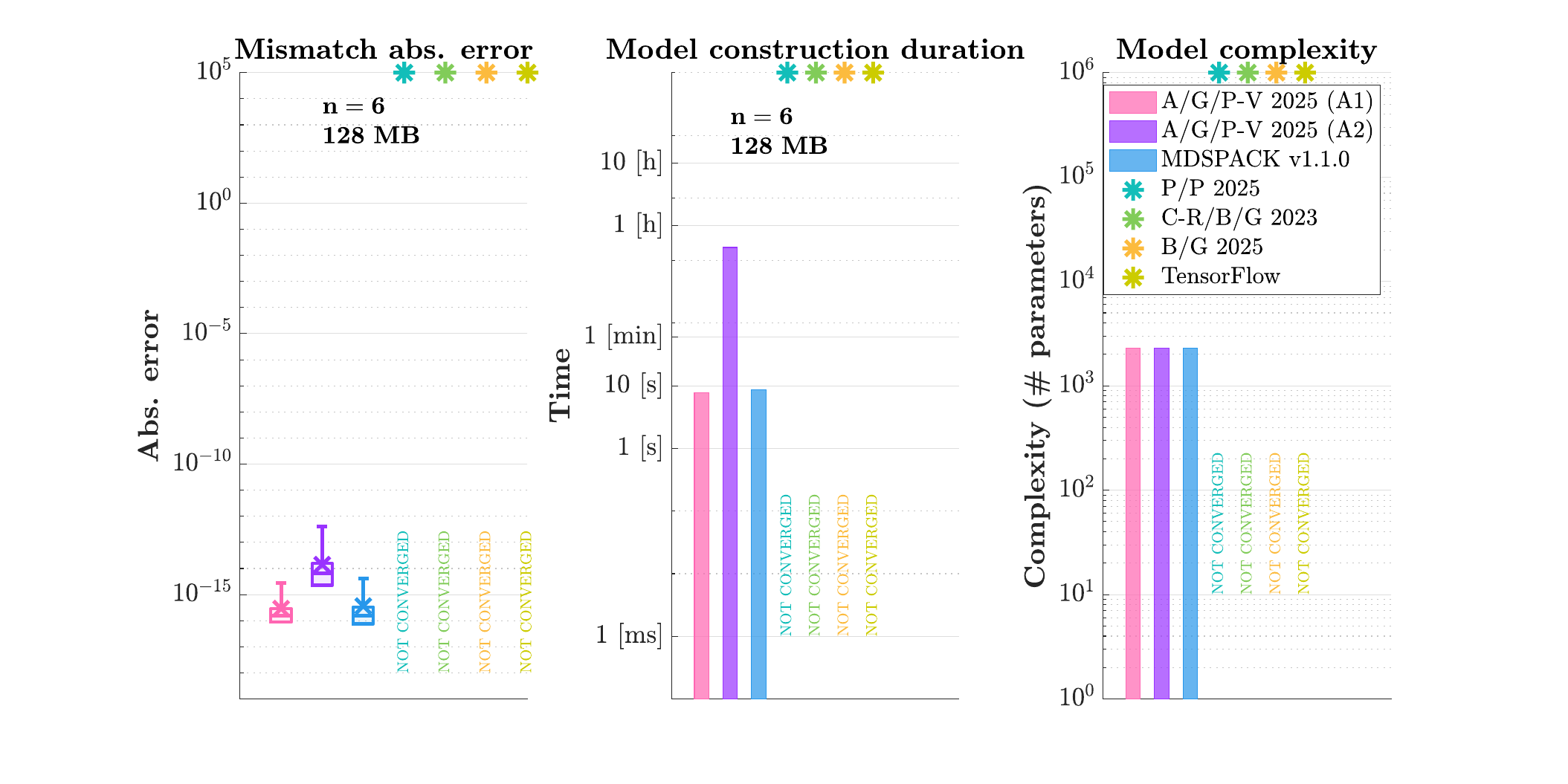} \caption{Function \#31: graphical view of the best model performances.} \end{figure}\begin{figure}[H] \centering  \includegraphics[width=\textwidth]{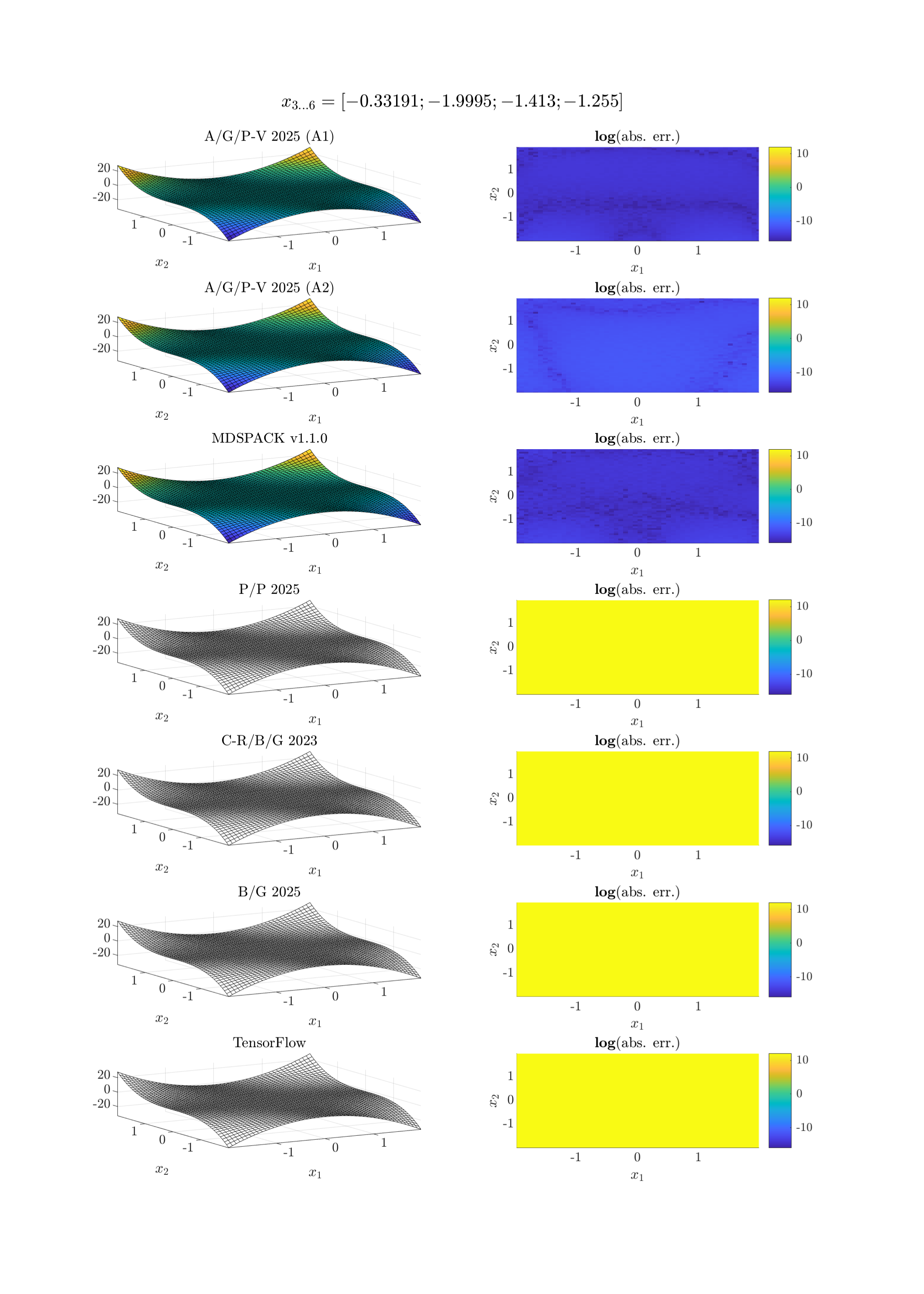} \caption{Function \#31: left side, evaluation of the original (mesh) vs. approximated (coloured surface) and right side, absolute errors (in log-scale).} \end{figure}\subsubsection{mLF detailed informations (M1)} \noindent \textbf{Right interpolation points}: $k_l=\left(\begin{array}{cccccc} 3 & 4 & 2 & 2 & 3 & 2 \end{array}\right)$, where $l=1,\cdots,\ord$.$$ \begin{array}{rcl}\lan{1} &\in& \IC^{3} \text{ , linearly spaced between bounds}\\\lan{2} &\in& \IC^{4} \text{ , linearly spaced between bounds}\\\lan{3} &\in& \IC^{2} \text{ , linearly spaced between bounds}\\\lan{4} &\in& \IC^{2} \text{ , linearly spaced between bounds}\\\lan{5} &\in& \IC^{3} \text{ , linearly spaced between bounds}\\\lan{6} &\in& \IC^{2} \text{ , linearly spaced between bounds}\\\end{array} $$\noindent \textbf{$\ord$-D Loewner matrix, barycentric weights and Lagrangian basis}:$$ \begin{array}{rcl}\IL & \in & \IC^{288 \times 288}\\\bc & \in & \IC^{288}\\\bw & \in & \IC^{288}\\\bc\odot \bw & \in & \IC^{288}\\\mathbf{Lag}(\var{1},\var{2},\var{3},\var{4},\var{5},\var{6}) & \in & \IC^{288}\\\end{array} $$

\newpage \subsection{Function \#32 (${\ord=2}$ variables, tensor size: 12.5 \textbf{KB})} $$\mathrm{atan}(\var{1}) + \var{2}^3$$ \subsubsection{Setup and results overview}\begin{itemize}\item Reference: Personal communication, [none]\item Domain: $\mathbb{R}$\item Tensor size: 12.5 \textbf{KB} ($40^{2}$ points)\item Bounds: $ \left(\begin{array}{cc} -2 & 2 \end{array}\right) \times \left(\begin{array}{cc} -2 & 2 \end{array}\right)$ \end{itemize} \begin{table}[H] \centering \begin{tabular}{llllll}
$\#$ & Alg. & Parameters & Dim. & CPU [s] & RMSE \\ 
\hline 
$\mathbf{\#32}$ & A/G/P-V 2025 (A1) & $1 \cdot 10^{-14},3$ & $\mathbf{2.6 \cdot 10^{02}}$ & $0.012$ & $\mathbf{8.3 \cdot 10^{-14}}$ \\ 
 & A/G/P-V 2025 (A2) & $1 \cdot 10^{-15},1$ & $3.2 \cdot 10^{02}$ & $0.52$ & $8.6 \cdot 10^{-13}$ \\ 
 & MDSPACK v1.1.0 & $1 \cdot 10^{-14},7$ & $2.6 \cdot 10^{02}$ & $\mathbf{0.0089}$ & $1.1 \cdot 10^{-13}$ \\ 
 & P/P 2025 & $1,1,50,0.01,10,4,21$ & $5.1 \cdot 10^{02}$ & $1.9$ & $0.00016$ \\ 
 & C-R/B/G 2023 & $1 \cdot 10^{-06},20$ & $6.1 \cdot 10^{02}$ & $0.18$ & $4.7 \cdot 10^{-07}$ \\ 
 & B/G 2025 & $1 \cdot 10^{-09},20,4$ & $7.2 \cdot 10^{02}$ & $0.31$ & $1.1 \cdot 10^{-08}$ \\ 
 & TensorFlow & $$ & $2.6 \cdot 10^{02}$ & $14$ & $0.028$ \\ 
\hline 
\end{tabular} \caption{Function \#32: best model configuration and performances per methods.} \end{table}\begin{figure}[H] \centering  \includegraphics[width=\textwidth]{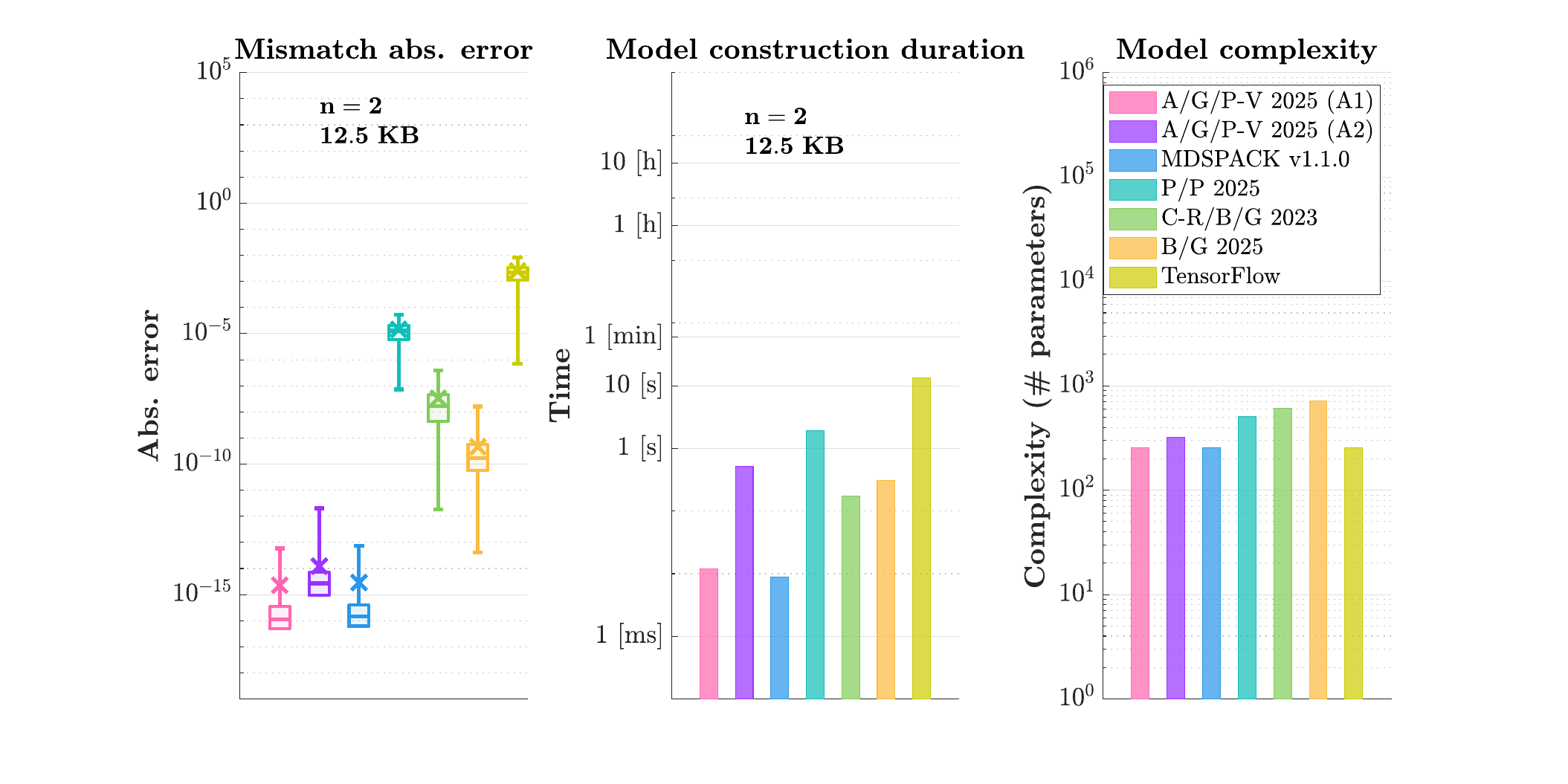} \caption{Function \#32: graphical view of the best model performances.} \end{figure}\begin{figure}[H] \centering  \includegraphics[width=\textwidth]{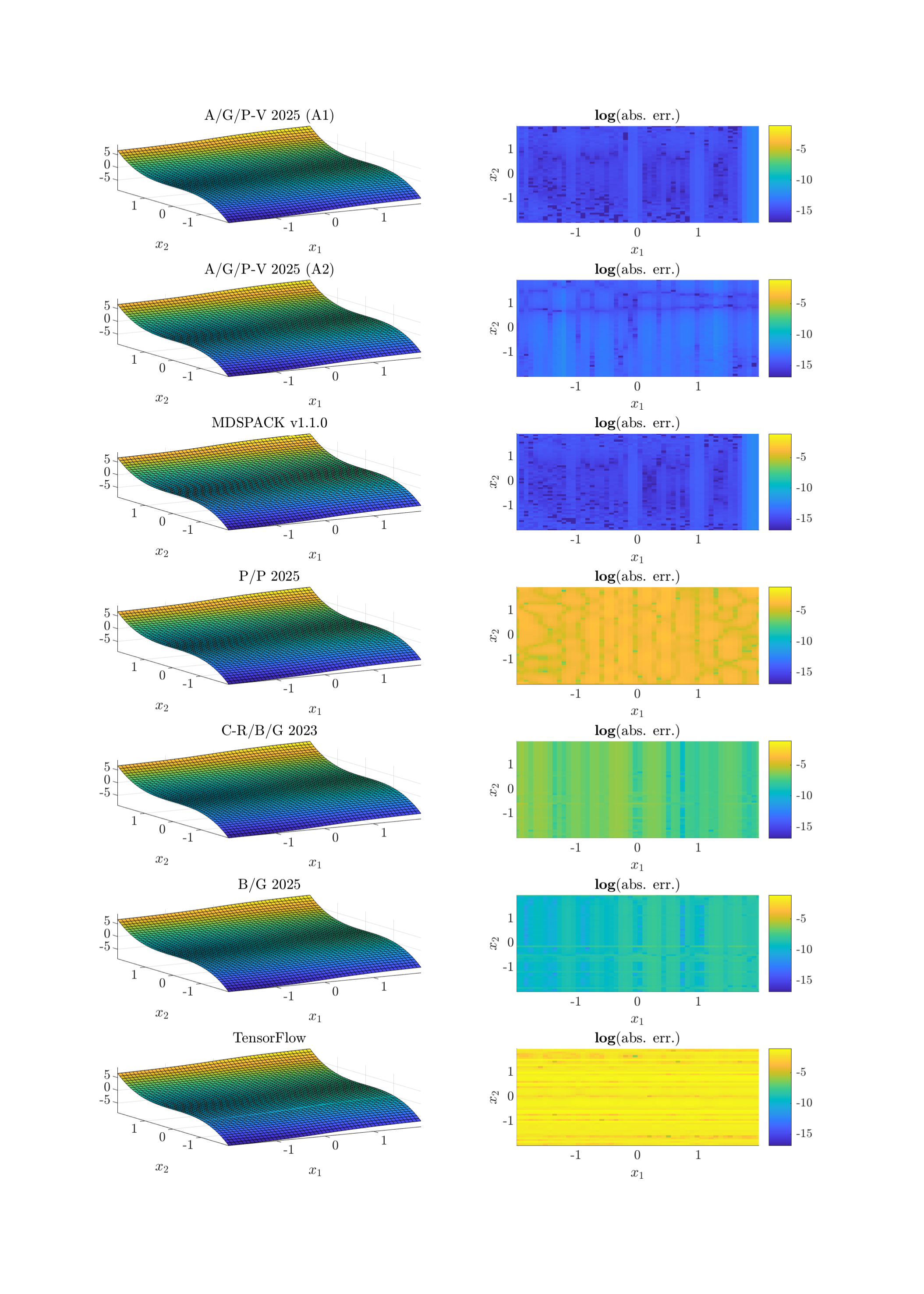} \caption{Function \#32: left side, evaluation of the original (mesh) vs. approximated (coloured surface) and right side, absolute errors (in log-scale).} \end{figure}\subsubsection{mLF detailed informations (M1)} \noindent \textbf{Right interpolation points}: $k_l=\left(\begin{array}{cc} 16 & 4 \end{array}\right)$, where $l=1,\cdots,\ord$.$$ \begin{array}{rcl}\lan{1} &\in& \IC^{16} \text{ , linearly spaced between bounds}\\\lan{2} &\in& \IC^{4} \text{ , linearly spaced between bounds}\\\end{array} $$\noindent \textbf{$\ord$-D Loewner matrix, barycentric weights and Lagrangian basis}:$$ \begin{array}{rcl}\IL & \in & \IC^{64 \times 64}\\\bc & \in & \IC^{64}\\\bw & \in & \IC^{64}\\\bc\odot \bw & \in & \IC^{64}\\\mathbf{Lag}(\var{1},\var{2}) & \in & \IC^{64}\\\end{array} $$

\newpage \subsection{Function \#33 (${\ord=2}$ variables, tensor size: 28.1 \textbf{KB})} $$\frac{\var{1}+\var{2}}{\cos(\var{1})^2+\cos(\var{2}) + 3}$$ \subsubsection{Setup and results overview}\begin{itemize}\item Reference: Personal communication, [none]\item Domain: $\mathbb{R}$\item Tensor size: 28.1 \textbf{KB} ($60^{2}$ points)\item Bounds: $ \left(\begin{array}{cc} -10 & 10 \end{array}\right) \times \left(\begin{array}{cc} -10 & 10 \end{array}\right)$ \end{itemize} \begin{table}[H] \centering \begin{tabular}{llllll}
$\#$ & Alg. & Parameters & Dim. & CPU [s] & RMSE \\ 
\hline 
$\mathbf{\#33}$ & A/G/P-V 2025 (A1) & $1 \cdot 10^{-10},2$ & $1.3 \cdot 10^{03}$ & $0.023$ & $0.0011$ \\ 
 & A/G/P-V 2025 (A2) & $1 \cdot 10^{-15},1$ & $\mathbf{4}$ & $0.29$ & $9$ \\ 
 & MDSPACK v1.1.0 & $1 \cdot 10^{-12},6$ & $1.3 \cdot 10^{03}$ & $\mathbf{0.018}$ & $0.00083$ \\ 
 & P/P 2025 & $1,0.95,50,0.01,10,6,21$ & $5.5 \cdot 10^{02}$ & $2.1$ & $0.054$ \\ 
 & C-R/B/G 2023 & $0.001,20$ & $1.2 \cdot 10^{03}$ & $0.49$ & $0.00025$ \\ 
 & B/G 2025 & $1 \cdot 10^{-06},20,2$ & $1.3 \cdot 10^{03}$ & $0.88$ & $\mathbf{6.7 \cdot 10^{-06}}$ \\ 
 & TensorFlow & $$ & $2.6 \cdot 10^{02}$ & $25$ & $0.32$ \\ 
\hline 
\end{tabular} \caption{Function \#33: best model configuration and performances per methods.} \end{table}\begin{figure}[H] \centering  \includegraphics[width=\textwidth]{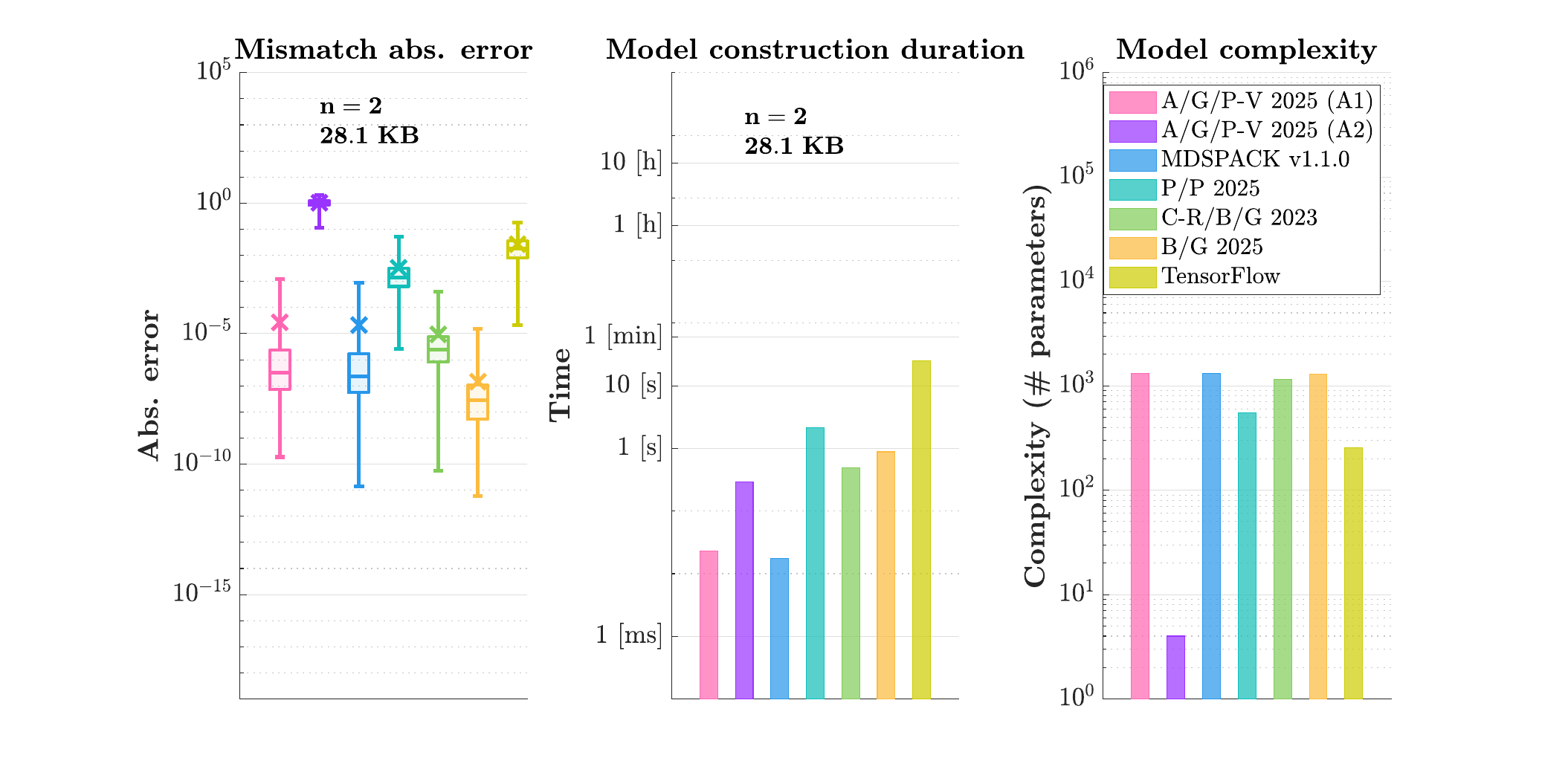} \caption{Function \#33: graphical view of the best model performances.} \end{figure}\begin{figure}[H] \centering  \includegraphics[width=\textwidth]{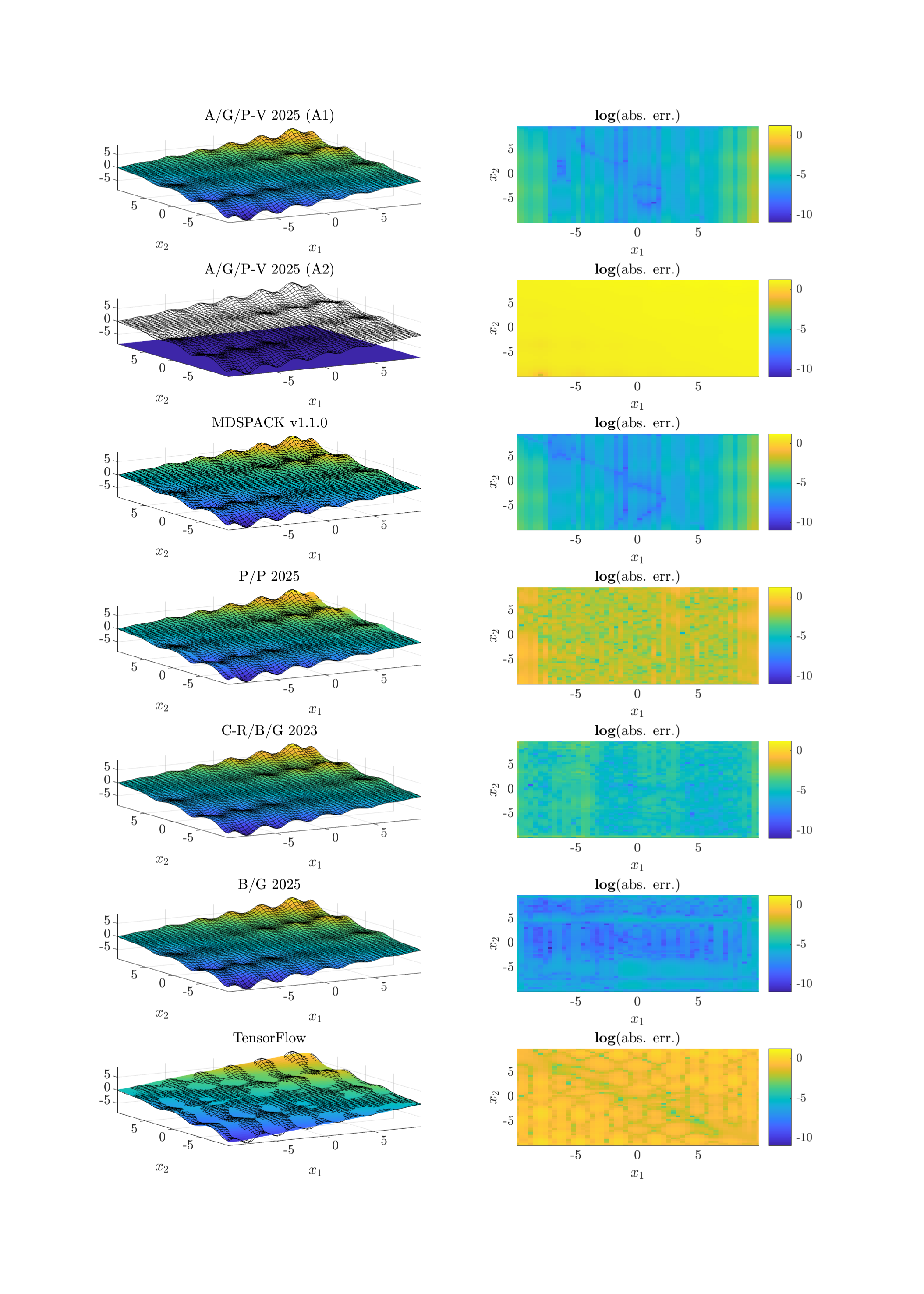} \caption{Function \#33: left side, evaluation of the original (mesh) vs. approximated (coloured surface) and right side, absolute errors (in log-scale).} \end{figure}\subsubsection{mLF detailed informations (M1)} \noindent \textbf{Right interpolation points}: $k_l=\left(\begin{array}{cc} 22 & 15 \end{array}\right)$, where $l=1,\cdots,\ord$.$$ \begin{array}{rcl}\lan{1} &\in& \IC^{22} \text{ , linearly spaced between bounds}\\\lan{2} &\in& \IC^{15} \text{ , linearly spaced between bounds}\\\end{array} $$\noindent \textbf{$\ord$-D Loewner matrix, barycentric weights and Lagrangian basis}:$$ \begin{array}{rcl}\IL & \in & \IC^{330 \times 330}\\\bc & \in & \IC^{330}\\\bw & \in & \IC^{330}\\\bc\odot \bw & \in & \IC^{330}\\\mathbf{Lag}(\var{1},\var{2}) & \in & \IC^{330}\\\end{array} $$

\newpage \subsection{Function \#34 (${\ord=2}$ variables, tensor size: 1.22 \textbf{MB})} $$\texttt{Riemann $\zeta$ function (real part)}$$ \subsubsection{Setup and results overview}\begin{itemize}\item Reference: Riemann $\zeta$ function (real part), [none]\item Domain: $\mathbb{R}$\item Tensor size: 1.22 \textbf{MB} ($400^{2}$ points)\item Bounds: $ \left(\begin{array}{cc} \frac{9}{20} & \frac{11}{20} \end{array}\right) \times \left(\begin{array}{cc} 1 & 50 \end{array}\right)$ \end{itemize} \begin{table}[H] \centering \begin{tabular}{llllll}
$\#$ & Alg. & Parameters & Dim. & CPU [s] & RMSE \\ 
\hline 
$\mathbf{\#34}$ & A/G/P-V 2025 (A1) & $1 \cdot 10^{-10},3$ & $2.3 \cdot 10^{03}$ & $0.82$ & $4.6 \cdot 10^{-05}$ \\ 
 & A/G/P-V 2025 (A2) & $1 \cdot 10^{-15},1$ & $\mathbf{4}$ & $10$ & $3.1$ \\ 
 & MDSPACK v1.1.0 & $1 \cdot 10^{-08},4$ & $1.8 \cdot 10^{03}$ & $\mathbf{0.77}$ & $\mathbf{2.4 \cdot 10^{-05}}$ \\ 
 & P/P 2025 & $1,1,50,0.01,10,12,21$ & $6.8 \cdot 10^{02}$ & $77$ & $0.068$ \\ 
 & C-R/B/G 2023 & $0.001,20$ & $1.2 \cdot 10^{03}$ & $84$ & $1.2$ \\ 
 & B/G 2025 & $1 \cdot 10^{-09},20,4$ & $1.4 \cdot 10^{03}$ & $17$ & $2$ \\ 
 & TensorFlow & $NaN$ & $NaN$ & $NaN$ & $NaN$ \\ 
\hline 
\end{tabular} \caption{Function \#34: best model configuration and performances per methods.} \end{table}\begin{figure}[H] \centering  \includegraphics[width=\textwidth]{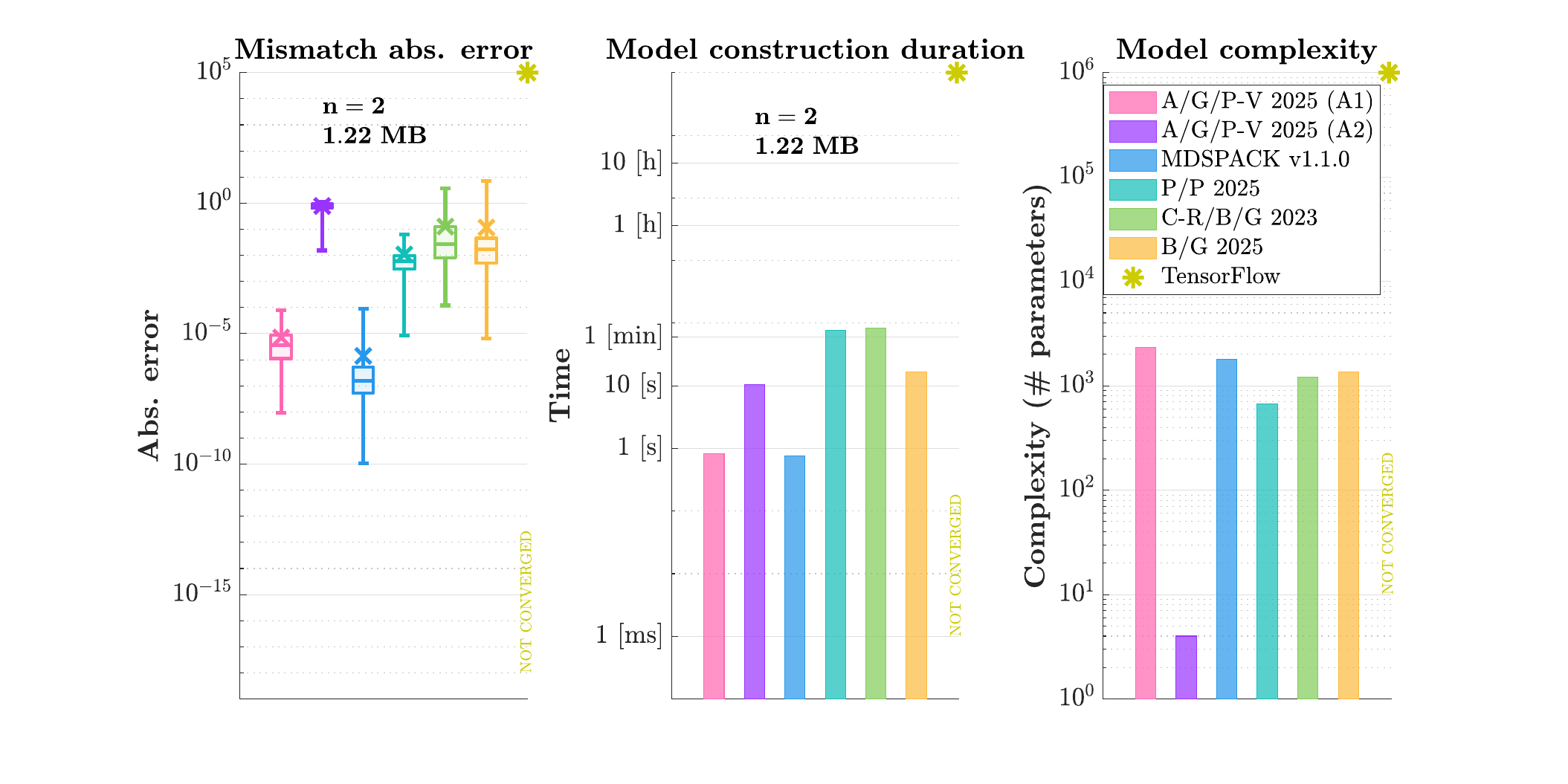} \caption{Function \#34: graphical view of the best model performances.} \end{figure}\begin{figure}[H] \centering  \includegraphics[width=\textwidth]{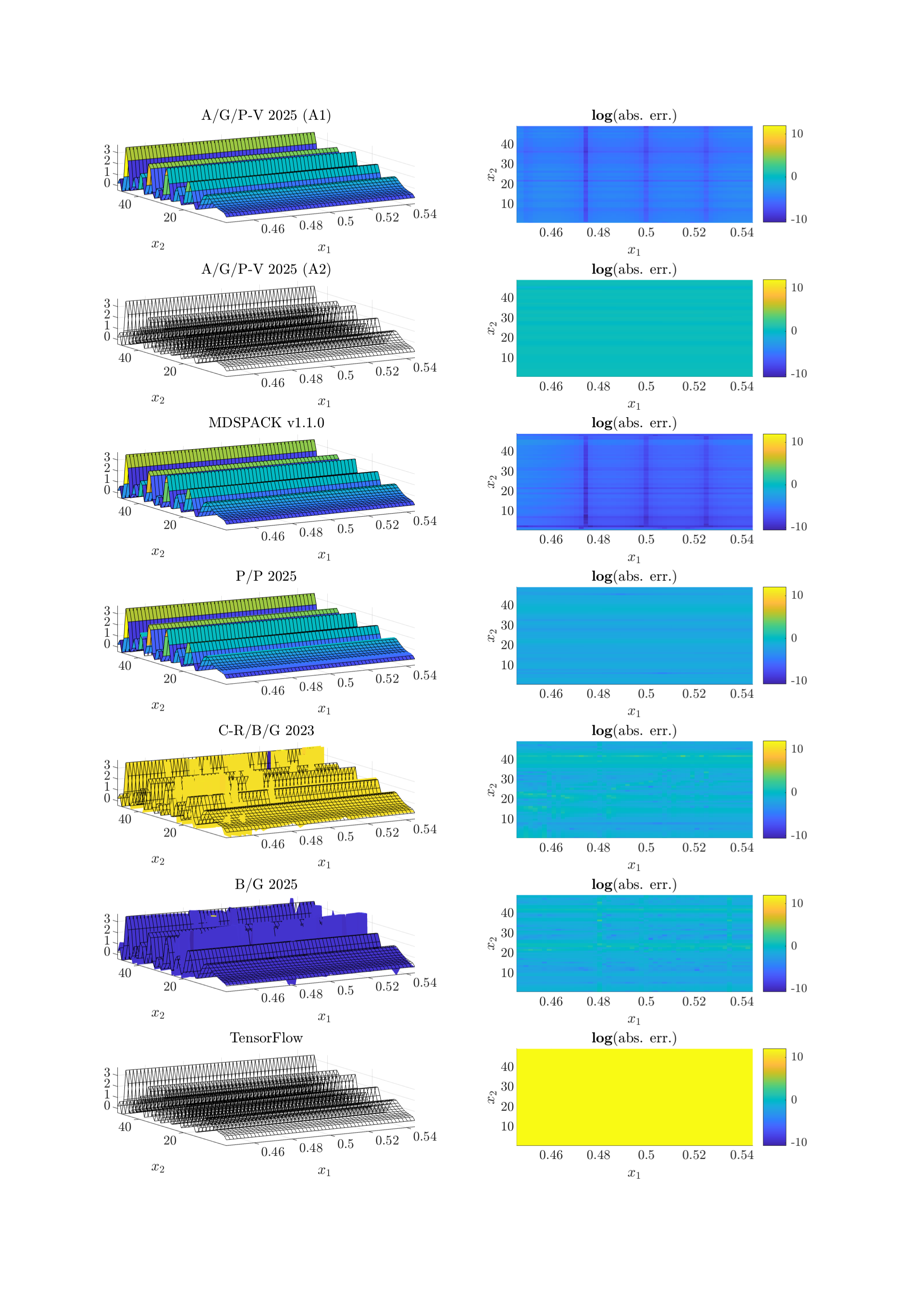} \caption{Function \#34: left side, evaluation of the original (mesh) vs. approximated (coloured surface) and right side, absolute errors (in log-scale).} \end{figure}\subsubsection{mLF detailed informations (M1)} \noindent \textbf{Right interpolation points}: $k_l=\left(\begin{array}{cc} 5 & 116 \end{array}\right)$, where $l=1,\cdots,\ord$.$$ \begin{array}{rcl}\lan{1} &\in& \IC^{5} \text{ , linearly spaced between bounds}\\\lan{2} &\in& \IC^{116} \text{ , linearly spaced between bounds}\\\end{array} $$\noindent \textbf{$\ord$-D Loewner matrix, barycentric weights and Lagrangian basis}:$$ \begin{array}{rcl}\IL & \in & \IC^{580 \times 580}\\\bc & \in & \IC^{580}\\\bw & \in & \IC^{580}\\\bc\odot \bw & \in & \IC^{580}\\\mathbf{Lag}(\var{1},\var{2}) & \in & \IC^{580}\\\end{array} $$

\newpage \subsection{Function \#35 (${\ord=2}$ variables, tensor size: 1.22 \textbf{MB})} $$\texttt{Riemann $\zeta$ function (imaginary part)}$$ \subsubsection{Setup and results overview}\begin{itemize}\item Reference: Riemann $\zeta$ function (imaginary part), [none]\item Domain: $\mathbb{R}$\item Tensor size: 1.22 \textbf{MB} ($400^{2}$ points)\item Bounds: $ \left(\begin{array}{cc} \frac{9}{20} & \frac{11}{20} \end{array}\right) \times \left(\begin{array}{cc} 1 & 50 \end{array}\right)$ \end{itemize} \begin{table}[H] \centering \begin{tabular}{llllll}
$\#$ & Alg. & Parameters & Dim. & CPU [s] & RMSE \\ 
\hline 
$\mathbf{\#35}$ & A/G/P-V 2025 (A1) & $1 \cdot 10^{-09},3$ & $1.8 \cdot 10^{03}$ & $0.82$ & $7.3 \cdot 10^{-05}$ \\ 
 & A/G/P-V 2025 (A2) & $1 \cdot 10^{-15},2$ & $\mathbf{1.4 \cdot 10^{02}}$ & $19$ & $0.95$ \\ 
 & MDSPACK v1.1.0 & $1 \cdot 10^{-12},6$ & $1.8 \cdot 10^{03}$ & $\mathbf{0.76}$ & $\mathbf{7.2 \cdot 10^{-05}}$ \\ 
 & P/P 2025 & $1,0.95,50,0.01,10,12,21$ & $6.8 \cdot 10^{02}$ & $78$ & $0.02$ \\ 
 & C-R/B/G 2023 & $0.001,20$ & $1.5 \cdot 10^{03}$ & $83$ & $0.79$ \\ 
 & B/G 2025 & $1 \cdot 10^{-06},20,4$ & $1.4 \cdot 10^{03}$ & $19$ & $0.94$ \\ 
 & TensorFlow & $NaN$ & $NaN$ & $NaN$ & $NaN$ \\ 
\hline 
\end{tabular} \caption{Function \#35: best model configuration and performances per methods.} \end{table}\begin{figure}[H] \centering  \includegraphics[width=\textwidth]{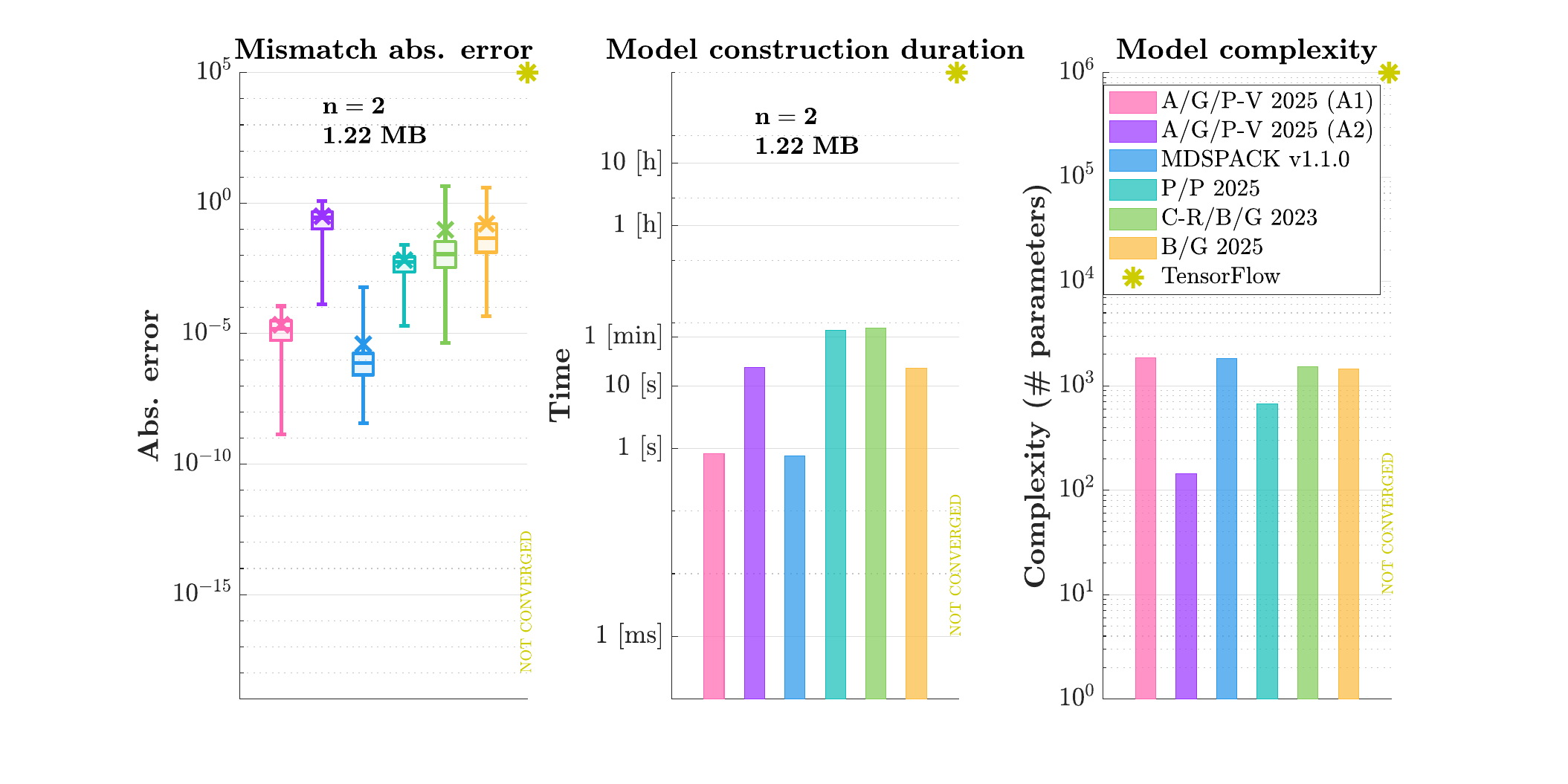} \caption{Function \#35: graphical view of the best model performances.} \end{figure}\begin{figure}[H] \centering  \includegraphics[width=\textwidth]{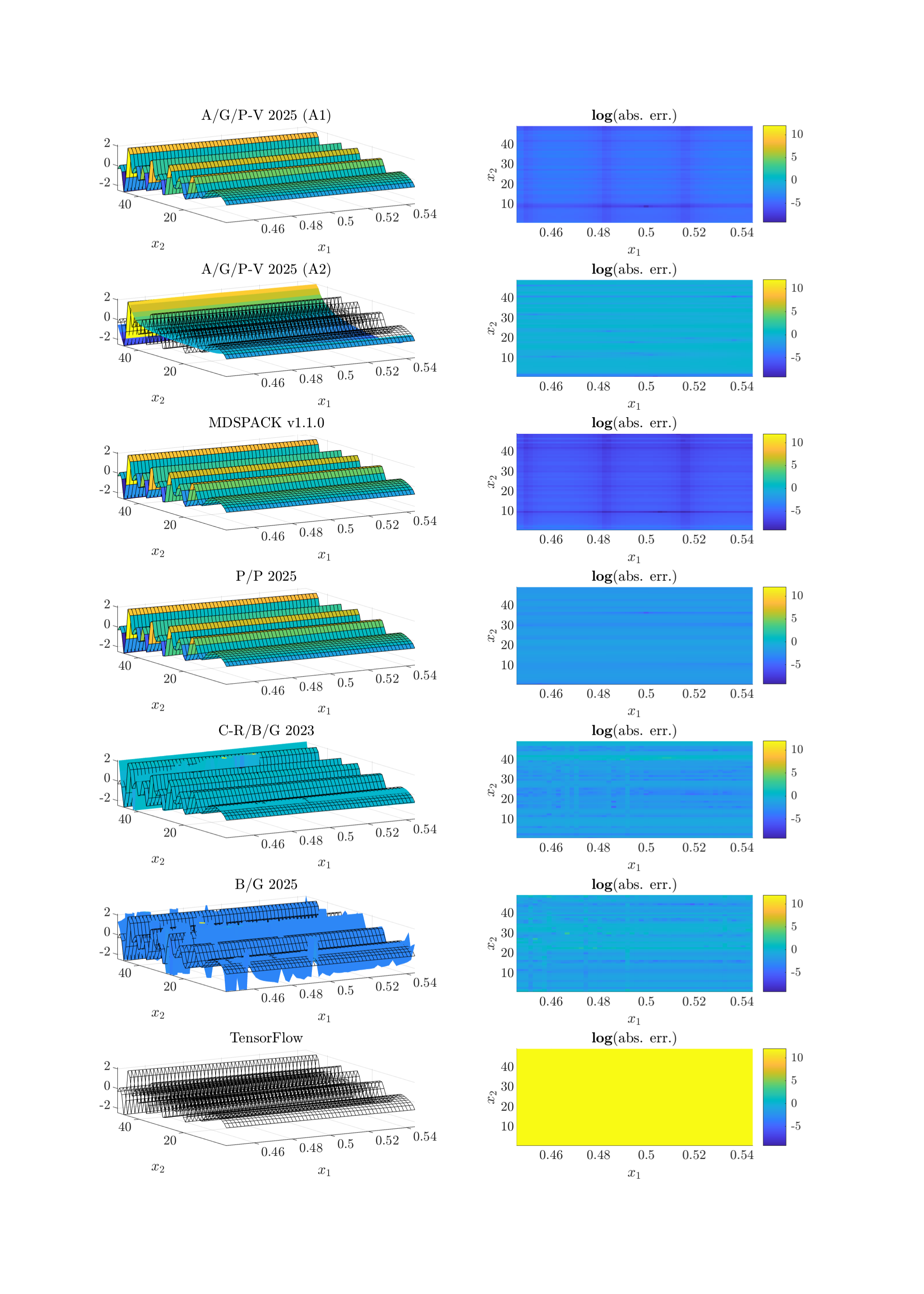} \caption{Function \#35: left side, evaluation of the original (mesh) vs. approximated (coloured surface) and right side, absolute errors (in log-scale).} \end{figure}\subsubsection{mLF detailed informations (M1)} \noindent \textbf{Right interpolation points}: $k_l=\left(\begin{array}{cc} 4 & 115 \end{array}\right)$, where $l=1,\cdots,\ord$.$$ \begin{array}{rcl}\lan{1} &\in& \IC^{4} \text{ , linearly spaced between bounds}\\\lan{2} &\in& \IC^{115} \text{ , linearly spaced between bounds}\\\end{array} $$\noindent \textbf{$\ord$-D Loewner matrix, barycentric weights and Lagrangian basis}:$$ \begin{array}{rcl}\IL & \in & \IC^{460 \times 460}\\\bc & \in & \IC^{460}\\\bw & \in & \IC^{460}\\\bc\odot \bw & \in & \IC^{460}\\\mathbf{Lag}(\var{1},\var{2}) & \in & \IC^{460}\\\end{array} $$

\newpage \subsection{Function \#36 (${\ord=3}$ variables, tensor size: 62.5 \textbf{KB})} $$\frac{\var{2}}{3+1/3 \var{2}\var{1}-\var{3}^2}$$ \subsubsection{Setup and results overview}\begin{itemize}\item Reference: Personal communication, [none]\item Domain: $\mathbb{R}$\item Tensor size: 62.5 \textbf{KB} ($20^{3}$ points)\item Bounds: $ \left(\begin{array}{cc} \frac{1}{10} & 1 \end{array}\right) \times \left(\begin{array}{cc} \frac{1}{10} & 1 \end{array}\right) \times \left(\begin{array}{cc} \frac{1}{10} & 1 \end{array}\right)$ \end{itemize} \begin{table}[H] \centering \begin{tabular}{llllll}
$\#$ & Alg. & Parameters & Dim. & CPU [s] & RMSE \\ 
\hline 
$\mathbf{\#36}$ & A/G/P-V 2025 (A1) & $0.01,2$ & $\mathbf{60}$ & $0.0085$ & $1.4 \cdot 10^{-15}$ \\ 
 & A/G/P-V 2025 (A2) & $1 \cdot 10^{-15},2$ & $60$ & $0.33$ & $\mathbf{5.2 \cdot 10^{-16}}$ \\ 
 & MDSPACK v1.1.0 & $0.0001,2$ & $60$ & $\mathbf{0.0077}$ & $1.4 \cdot 10^{-15}$ \\ 
 & P/P 2025 & $1,1,50,0.01,4,12,9$ & $2.2 \cdot 10^{02}$ & $2.3$ & $6.6 \cdot 10^{-06}$ \\ 
 & C-R/B/G 2023 & $0.001,20$ & $6.4 \cdot 10^{02}$ & $0.26$ & $3.7 \cdot 10^{-14}$ \\ 
 & B/G 2025 & $1 \cdot 10^{-09},20,4$ & $4.2 \cdot 10^{03}$ & $3.6$ & $2.4 \cdot 10^{-10}$ \\ 
 & TensorFlow & $$ & $3.2 \cdot 10^{02}$ & $40$ & $0.0054$ \\ 
\hline 
\end{tabular} \caption{Function \#36: best model configuration and performances per methods.} \end{table}\begin{figure}[H] \centering  \includegraphics[width=\textwidth]{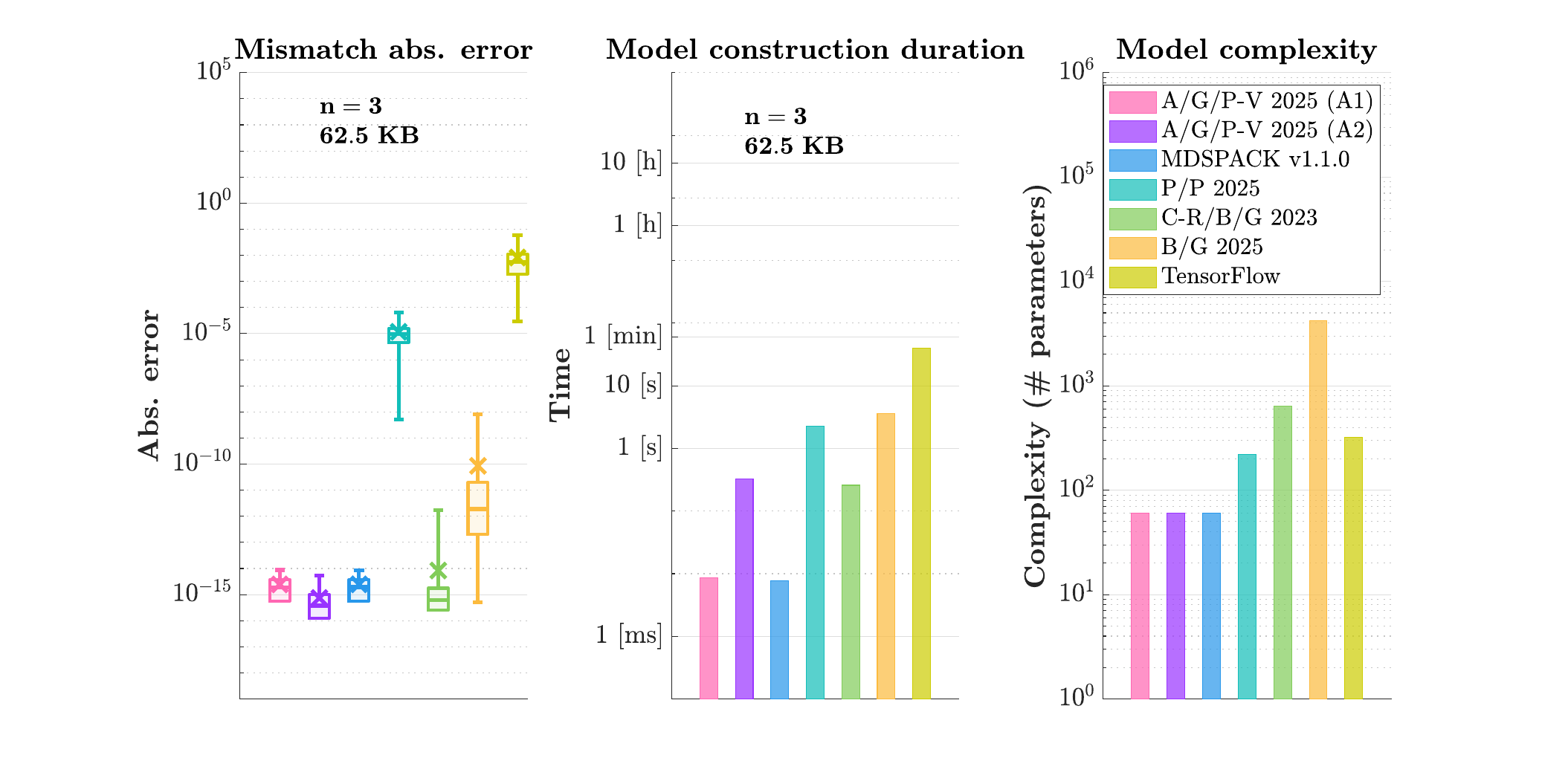} \caption{Function \#36: graphical view of the best model performances.} \end{figure}\begin{figure}[H] \centering  \includegraphics[width=\textwidth]{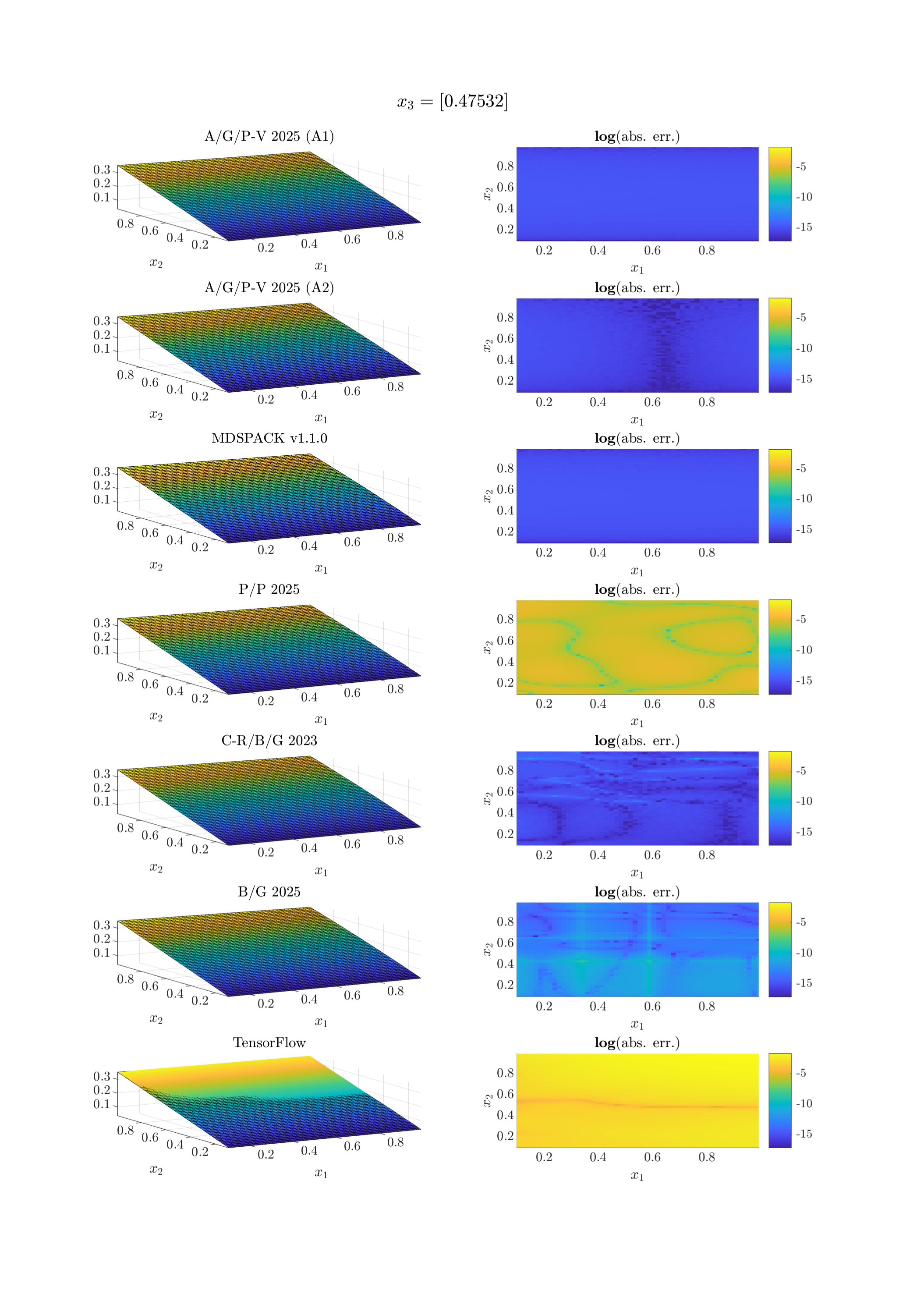} \caption{Function \#36: left side, evaluation of the original (mesh) vs. approximated (coloured surface) and right side, absolute errors (in log-scale).} \end{figure}\subsubsection{mLF detailed informations (M1)} \noindent \textbf{Right interpolation points} ($k_l=\left(\begin{array}{ccc} 2 & 2 & 3 \end{array}\right)$, where $l=1,\cdots,\ord$):$$ \begin{array}{rcl}\lan{1} &=& \left(\begin{array}{cc} \frac{1}{10} & 1 \end{array}\right)\\\lan{2} &=& \left(\begin{array}{cc} \frac{1}{10} & 1 \end{array}\right)\\\lan{3} &=& \left(\begin{array}{ccc} \frac{1}{10} & \frac{1}{2} & 1 \end{array}\right)\\\end{array} $$\noindent \textbf{Lagrangian weights}: $$\left(\begin{array}{ccc} \bc & \bw & \bc\odot\bw\\ 1.604 & 0.03341 & 0.05357\\ -2.655 & 0.03632 & -0.09643\\ 0.8586 & 0.04992 & 0.04286\\ -1.62 & 0.3308 & -0.5357\\ 2.684 & 0.3593 & 0.9643\\ -0.8714 & 0.4918 & -0.4286\\ -1.62 & 0.03308 & -0.05357\\ 2.684 & 0.03593 & 0.09643\\ -0.8714 & 0.04918 & -0.04286\\ 1.78 & 0.3009 & 0.5357\\ -2.973 & 0.3243 & -0.9643\\ 1.0 & 0.4286 & 0.4286 \end{array}\right)$$\noindent \textbf{Lagrangian form} (basis, numerator and denominator coefficients):$$\left(\begin{array}{ccc}\mathcal{B}_\textrm{lag}(\var{1},\var{2},\var{3}) & \bN_\textrm{lag} &\bD_\textrm{lag}\end{array}\right) =$$ $$\left(\begin{array}{ccc} \left(\var{1}-0.1\right)\,\left(\var{2}-0.1\right)\,\left(\var{3}-0.1\right) & 0.05357 & 1.604\\ \left(\var{3}-0.5\right)\,\left(\var{1}-0.1\right)\,\left(\var{2}-0.1\right) & -0.09643 & -2.655\\ \left(\var{3}-1.0\right)\,\left(\var{1}-0.1\right)\,\left(\var{2}-0.1\right) & 0.04286 & 0.8586\\ \left(\var{2}-1.0\right)\,\left(\var{1}-0.1\right)\,\left(\var{3}-0.1\right) & -0.5357 & -1.62\\ \left(\var{2}-1.0\right)\,\left(\var{3}-0.5\right)\,\left(\var{1}-0.1\right) & 0.9643 & 2.684\\ \left(\var{2}-1.0\right)\,\left(\var{3}-1.0\right)\,\left(\var{1}-0.1\right) & -0.4286 & -0.8714\\ \left(\var{1}-1.0\right)\,\left(\var{2}-0.1\right)\,\left(\var{3}-0.1\right) & -0.05357 & -1.62\\ \left(\var{1}-1.0\right)\,\left(\var{3}-0.5\right)\,\left(\var{2}-0.1\right) & 0.09643 & 2.684\\ \left(\var{1}-1.0\right)\,\left(\var{3}-1.0\right)\,\left(\var{2}-0.1\right) & -0.04286 & -0.8714\\ \left(\var{1}-1.0\right)\,\left(\var{2}-1.0\right)\,\left(\var{3}-0.1\right) & 0.5357 & 1.78\\ \left(\var{1}-1.0\right)\,\left(\var{2}-1.0\right)\,\left(\var{3}-0.5\right) & -0.9643 & -2.973\\ \left(\var{1}-1.0\right)\,\left(\var{2}-1.0\right)\,\left(\var{3}-1.0\right) & 0.4286 & 1.0 \end{array}\right).$$\noindent The corresponding function is:$$\begin{array}{rcl}\bG_{\textrm{lag}}(\var{1},\var{2},\var{3}) &=& \dfrac{\bn_{\textrm{lag}}(\var{1},\var{2},\var{3})}{\bd_{\textrm{lag}}(\var{1},\var{2},\var{3})}\\ && \\&=& \dfrac{\sum_{\textrm{row}} \bN_\textrm{lag} \odot\mathcal{B}^{-1}_\textrm{lag}(\var{1},\var{2},\var{3})}{\sum_{\textrm{row}} \bD_\textrm{lag} \odot\mathcal{B}^{-1}_\textrm{lag}(\var{1},\var{2},\var{3})}, \end{array}$$\noindent where,\\$\bn_{\textrm{lag}}(\var{1},\var{2},\var{3}) = 0.3333\,\var{2}-1.655 \cdot 10^{-15}\,\var{1}-9.543 \cdot 10^{-15}\,\var{3}-1.01 \cdot 10^{-15}\,\var{1}\,\var{2}+4.965 \cdot 10^{-15}\,\var{1}\,\var{3}+4.275 \cdot 10^{-14}\,\var{2}\,\var{3}-3.31 \cdot 10^{-15}\,\var{1}\,{\var{3}}^2-2.85 \cdot 10^{-14}\,\var{2}\,{\var{3}}^2+6.362 \cdot 10^{-15}\,{\var{3}}^2-2.02 \cdot 10^{-15}\,\var{1}\,\var{2}\,{\var{3}}^2+3.03 \cdot 10^{-15}\,\var{1}\,\var{2}\,\var{3}+3.181 \cdot 10^{-15}$ \\~~\\$\bd_{\textrm{lag}}(\var{1},\var{2},\var{3}) = 0.1111\,\var{1}\,\var{2}-9.514 \cdot 10^{-14}\,\var{2}-1.868 \cdot 10^{-13}\,\var{3}-5.721 \cdot 10^{-14}\,\var{1}+1.716 \cdot 10^{-13}\,\var{1}\,\var{3}+2.854 \cdot 10^{-13}\,\var{2}\,\var{3}-1.144 \cdot 10^{-13}\,\var{1}\,{\var{3}}^2-1.903 \cdot 10^{-13}\,\var{2}\,{\var{3}}^2-0.3333\,{\var{3}}^2+8.944 \cdot 10^{-14}\,\var{1}\,\var{2}\,{\var{3}}^2-1.342 \cdot 10^{-13}\,\var{1}\,\var{2}\,\var{3}+1.0$ \\~~\\\noindent \textbf{Monomial form} (basis, numerator and denominator coefficients - entries $<10^{-12}$ removed):$$\left(\begin{array}{ccc}\mathcal{B}_\textrm{mon}(\var{1},\var{2},\var{3}) & \bN_\textrm{mon} &\bD_\textrm{mon}\end{array}\right) =$$ $$\left(\begin{array}{ccc} \var{1}\,\var{2}\,{\var{3}}^2 & 0 & 0\\ \var{1}\,\var{2}\,\var{3} & 0 & 0\\ \var{1}\,\var{2} & 0 & 0.1111\\ \var{1}\,{\var{3}}^2 & 0 & 0\\ \var{1}\,\var{3} & 0 & 0\\ \var{1} & 0 & 0\\ \var{2}\,{\var{3}}^2 & 0 & 0\\ \var{2}\,\var{3} & 0 & 0\\ \var{2} & 0.3333 & 0\\ {\var{3}}^2 & 0 & -0.3333\\ \var{3} & 0 & 0\\ 1.0 & 0 & 1.0 \end{array}\right)$$\noindent The corresponding function is:$$\begin{array}{rcl}\bG_{\textrm{mon}}(\var{1},\var{2},\var{3}) &=& \dfrac{\bn_{\textrm{mon}}(\var{1},\var{2},\var{3})}{\bd_{\textrm{mon}}(\var{1},\var{2},\var{3})}\\ && \\&=& \dfrac{\sum_{\textrm{row}} \bN_\textrm{mon} \odot \mathcal{B}_\textrm{mon}(\var{1},\var{2},\var{3})}{\sum_{\textrm{row}} \bD_\textrm{mon} \odot\mathcal{B}_\textrm{mon}(\var{1},\var{2},\var{3})},  \end{array}$$\noindent where,\\$\bn_{\textrm{mon}}(\var{1},\var{2},\var{3}) = 0.3333\,\var{2}$ \\~~\\$\bd_{\textrm{mon}}(\var{1},\var{2},\var{3}) = -0.3333\,{\var{3}}^2+0.1111\,\var{1}\,\var{2}+1.0$ \\~~\\\noindent \textbf{KST equivalent decoupling pattern} (Barycentric weights $\bc^{\var{l}}$): $$\begin{array}{rclll}\var{3}&: & \left(\begin{array}{cccc} 1.868 & 1.859 & 1.859 & 1.78\\ -3.092 & -3.08 & -3.08 & -2.973\\ 1.0 & 1.0 & 1.0 & 1.0 \end{array}\right)& \textrm{vec}(.) & := \textbf{Bary}(\var{3}) \\\var{2}&: & \left(\begin{array}{cc} -0.9852 & -0.8714\\ 1.0 & 1.0 \end{array}\right)& \textrm{vec}(.) \otimes \bone_{k_{3}} & := \textbf{Bary}(\var{2}) \\\var{1}&: & \left(\begin{array}{c} -0.8714\\ 1.0 \end{array}\right)& \textrm{vec}(.) \otimes \bone_{k_{3}k_{2}} & := \textbf{Bary}(\var{1}) \\\end{array}$$~\\ Then, with the above notations, one defines the following univariate vector functions:  $$ \left\{ \begin{array}{rcl}\bPhi_{1}(\var{1}) &:=& \textbf{Bary}(\var{1}) \odot \mathbf{Lag}(\var{1}) \\\bPhi_{2}(\var{2}) &:=& \textbf{Bary}(\var{2}) \odot \mathbf{Lag}(\var{2}) \\\bPhi_{3}(\var{3}) &:=& \textbf{Bary}(\var{3}) \odot \mathbf{Lag}(\var{3}) \\\end{array} \right. $$\noindent The corresponding function is:$$\begin{array}{rcl}\bG_{\textrm{kst}}(\var{1},\var{2},\var{3}) &=& \dfrac{\bn_{\textrm{kst}}(\var{1},\var{2},\var{3})}{\bd_{\textrm{kst}}(\var{1},\var{2},\var{3})}\\ && \\ &=& \dfrac{\sum_{\text{rows}} \bw \odot \bPhi_{1}(\var{1}) \odot \cdots \odot\bPhi_{3}(\var{3})}{\sum_{\text{rows}} \bPhi_{1}(\var{1}) \odot \cdots \odot\bPhi_{3}(\var{3})} . \end{array}$$~\\ \noindent \textbf{KST-like univariate functions} (equivalent scaled univariate functions $\bphi_{1,\cdots,3}$): $$\left\{\begin{array}{rcrcl}z_{1} &=&\bphi_{1}(\var{1}) &=& \cfrac{3.0}{\var{1}+6.0}\\z_{2} &=&\bphi_{2}(\var{2}) &=& \cfrac{30.0\,\var{2}}{\var{2}+60.0}\\z_{3} &=&\bphi_{3}(\var{3}) &=& \cfrac{\bn_{3}}{\bd_{3}} \\\end{array} \right. .$$\noindent where, \\ \noindent $\bn_{3}=3.37 \cdot 10^{-15}\,{\var{3}}^2-5.056 \cdot 10^{-15}\,\var{3}+0.0333$ and \\ \noindent $\bd_{3}=-0.333\,{\var{3}}^2-1.513 \cdot 10^{-13}\,\var{3}+1.0$, \\ \noindent \textbf{Connection with Neural Networks} (equivalent numerator $\bn_{\textrm{lag}}$ representation):\\ \begin{figure}[H]\begin{center} \scalebox{.7}{\begin{tikzpicture}[line width=0.4mm]\tikzstyle{place}=[circle, draw=black, minimum size = 8mm]\tikzstyle{placeInOut}=[circle, draw=orange, minimum size = 8mm]\node at (0,-2) [placeInOut] (first_1){$\var{1}$};\node at (0,-4) [placeInOut] (first_2){$\var{2}$};\node at (0,-6) [placeInOut] (first_3){$\var{3}$};\node at (5,-2) [place] (secondL1_1){$\frac{1}{\var{1}-\lani{1}{1}}$};\node at (5,-4) [place] (secondL1_2){$\frac{1}{\var{1}-\lani{1}{2}}$};\node at (5,-6) [place] (secondL2_1){$\frac{1}{\var{2}-\lani{2}{1}}$};\node at (5,-8) [place] (secondL2_2){$\frac{1}{\var{2}-\lani{2}{2}}$};\node at (5,-10) [place] (secondL3_1){$\frac{1}{\var{3}-\lani{3}{1}}$};\node at (5,-12) [place] (secondL3_2){$\frac{1}{\var{3}-\lani{3}{2}}$};\node at (5,-14) [place] (secondL3_3){$\frac{1}{\var{3}-\lani{3}{3}}$};\node at (10,-2) [place] (third_1){$\prod$};\node at (10,-4) [place] (third_2){$\prod$};\node at (10,-6) [place] (third_3){$\prod$};\node at (10,-8) [place] (third_4){$\prod$};\node at (10,-10) [place] (third_5){$\prod$};\node at (10,-12) [place] (third_6){$\prod$};\node at (10,-14) [place] (third_7){$\prod$};\node at (10,-16) [place] (third_8){$\prod$};\node at (10,-18) [place] (third_9){$\prod$};\node at (10,-20) [place] (third_10){$\prod$};\node at (10,-22) [place] (third_11){$\prod$};\node at (10,-24) [place] (third_12){$\prod$};\node at (15,-13) [placeInOut] (output){$\bSigma$};\draw[->] (first_1)--(secondL1_1) node[above,sloped,pos=0.75] { };\draw[->] (first_1)--(secondL1_2) node[above,sloped,pos=0.75] { };\draw[->] (first_2)--(secondL2_1) node[above,sloped,pos=0.75] { };\draw[->] (first_2)--(secondL2_2) node[above,sloped,pos=0.75] { };\draw[->] (first_3)--(secondL3_1) node[above,sloped,pos=0.75] { };\draw[->] (first_3)--(secondL3_2) node[above,sloped,pos=0.75] { };\draw[->] (first_3)--(secondL3_3) node[above,sloped,pos=0.75] { };\draw[->] (secondL1_1)--(third_1) node[above,sloped,pos=0.25] {};\draw[->] (secondL1_1)--(third_2) node[above,sloped,pos=0.25] {};\draw[->] (secondL1_1)--(third_3) node[above,sloped,pos=0.25] {};\draw[->] (secondL1_1)--(third_4) node[above,sloped,pos=0.25] {};\draw[->] (secondL1_1)--(third_5) node[above,sloped,pos=0.25] {};\draw[->] (secondL1_1)--(third_6) node[above,sloped,pos=0.25] {};\draw[->] (secondL1_2)--(third_7) node[above,sloped,pos=0.25] {};\draw[->] (secondL1_2)--(third_8) node[above,sloped,pos=0.25] {};\draw[->] (secondL1_2)--(third_9) node[above,sloped,pos=0.25] {};\draw[->] (secondL1_2)--(third_10) node[above,sloped,pos=0.25] {};\draw[->] (secondL1_2)--(third_11) node[above,sloped,pos=0.25] {};\draw[->] (secondL1_2)--(third_12) node[above,sloped,pos=0.25] {};\draw[->] (secondL2_1)--(third_1) node[above,sloped,pos=0.25] {};\draw[->] (secondL2_1)--(third_2) node[above,sloped,pos=0.25] {};\draw[->] (secondL2_1)--(third_3) node[above,sloped,pos=0.25] {};\draw[->] (secondL2_2)--(third_4) node[above,sloped,pos=0.25] {};\draw[->] (secondL2_2)--(third_5) node[above,sloped,pos=0.25] {};\draw[->] (secondL2_2)--(third_6) node[above,sloped,pos=0.25] {};\draw[->] (secondL2_1)--(third_7) node[above,sloped,pos=0.25] {};\draw[->] (secondL2_1)--(third_8) node[above,sloped,pos=0.25] {};\draw[->] (secondL2_1)--(third_9) node[above,sloped,pos=0.25] {};\draw[->] (secondL2_2)--(third_10) node[above,sloped,pos=0.25] {};\draw[->] (secondL2_2)--(third_11) node[above,sloped,pos=0.25] {};\draw[->] (secondL2_2)--(third_12) node[above,sloped,pos=0.25] {};\draw[->] (secondL3_1)--(third_1) node[above,sloped,pos=0.25] {};\draw[->] (secondL3_2)--(third_2) node[above,sloped,pos=0.25] {};\draw[->] (secondL3_3)--(third_3) node[above,sloped,pos=0.25] {};\draw[->] (secondL3_1)--(third_4) node[above,sloped,pos=0.25] {};\draw[->] (secondL3_2)--(third_5) node[above,sloped,pos=0.25] {};\draw[->] (secondL3_3)--(third_6) node[above,sloped,pos=0.25] {};\draw[->] (secondL3_1)--(third_7) node[above,sloped,pos=0.25] {};\draw[->] (secondL3_2)--(third_8) node[above,sloped,pos=0.25] {};\draw[->] (secondL3_3)--(third_9) node[above,sloped,pos=0.25] {};\draw[->] (secondL3_1)--(third_10) node[above,sloped,pos=0.25] {};\draw[->] (secondL3_2)--(third_11) node[above,sloped,pos=0.25] {};\draw[->] (secondL3_3)--(third_12) node[above,sloped,pos=0.25] {};\draw[->] (third_1)--(output) node[above,sloped,pos=0.25] {0.053571};\draw[->] (third_2)--(output) node[above,sloped,pos=0.25] {-0.096429};\draw[->] (third_3)--(output) node[above,sloped,pos=0.25] {0.042857};\draw[->] (third_4)--(output) node[above,sloped,pos=0.25] {-0.53571};\draw[->] (third_5)--(output) node[above,sloped,pos=0.25] {0.96429};\draw[->] (third_6)--(output) node[above,sloped,pos=0.25] {-0.42857};\draw[->] (third_7)--(output) node[above,sloped,pos=0.25] {-0.053571};\draw[->] (third_8)--(output) node[above,sloped,pos=0.25] {0.096429};\draw[->] (third_9)--(output) node[above,sloped,pos=0.25] {-0.042857};\draw[->] (third_10)--(output) node[above,sloped,pos=0.25] {0.53571};\draw[->] (third_11)--(output) node[above,sloped,pos=0.25] {-0.96429};\draw[->] (third_12)--(output) node[above,sloped,pos=0.25] {0.42857};\end{tikzpicture}} \caption{Equivalent NN representation of the numerator $\bn_{\textrm{lag}}$.}\end{center}\end{figure}

\newpage \subsection{Function \#37 (${\ord=4}$ variables, tensor size: 1.22 \textbf{MB})} $$\var{1}\var{4}^3+\sin(2\var{2})\var{3}$$ \subsubsection{Setup and results overview}\begin{itemize}\item Reference: Personal communication, [none]\item Domain: $\mathbb{R}$\item Tensor size: 1.22 \textbf{MB} ($20^{4}$ points)\item Bounds: $ \left(\begin{array}{cc} \frac{1}{1000} & 1 \end{array}\right) \times \left(\begin{array}{cc} \frac{1}{1000} & 1 \end{array}\right) \times \left(\begin{array}{cc} \frac{1}{1000} & 1 \end{array}\right) \times \left(\begin{array}{cc} \frac{1}{1000} & 1 \end{array}\right)$ \end{itemize} \begin{table}[H] \centering \begin{tabular}{llllll}
$\#$ & Alg. & Parameters & Dim. & CPU [s] & RMSE \\ 
\hline 
$\mathbf{\#37}$ & A/G/P-V 2025 (A1) & $1 \cdot 10^{-09},1$ & $5.8 \cdot 10^{02}$ & $0.041$ & $4.9 \cdot 10^{-10}$ \\ 
 & A/G/P-V 2025 (A2) & $1 \cdot 10^{-15},1$ & $5.8 \cdot 10^{02}$ & $28$ & $2.2 \cdot 10^{-09}$ \\ 
 & MDSPACK v1.1.0 & $1 \cdot 10^{-08},4$ & $5.8 \cdot 10^{02}$ & $\mathbf{0.039}$ & $4.9 \cdot 10^{-10}$ \\ 
 & P/P 2025 & $1,1,50,0.01,4,6,9$ & $\mathbf{2 \cdot 10^{02}}$ & $24$ & $6.9 \cdot 10^{-06}$ \\ 
 & C-R/B/G 2023 & $1 \cdot 10^{-06},20$ & $4.6 \cdot 10^{04}$ & $2 \cdot 10^{03}$ & $2.9 \cdot 10^{-11}$ \\ 
 & B/G 2025 & $1 \cdot 10^{-09},20,4$ & $9.4 \cdot 10^{03}$ & $14$ & $\mathbf{1.6 \cdot 10^{-13}}$ \\ 
 & TensorFlow & $$ & $3.8 \cdot 10^{02}$ & $1.4 \cdot 10^{02}$ & $0.003$ \\ 
\hline 
\end{tabular} \caption{Function \#37: best model configuration and performances per methods.} \end{table}\begin{figure}[H] \centering  \includegraphics[width=\textwidth]{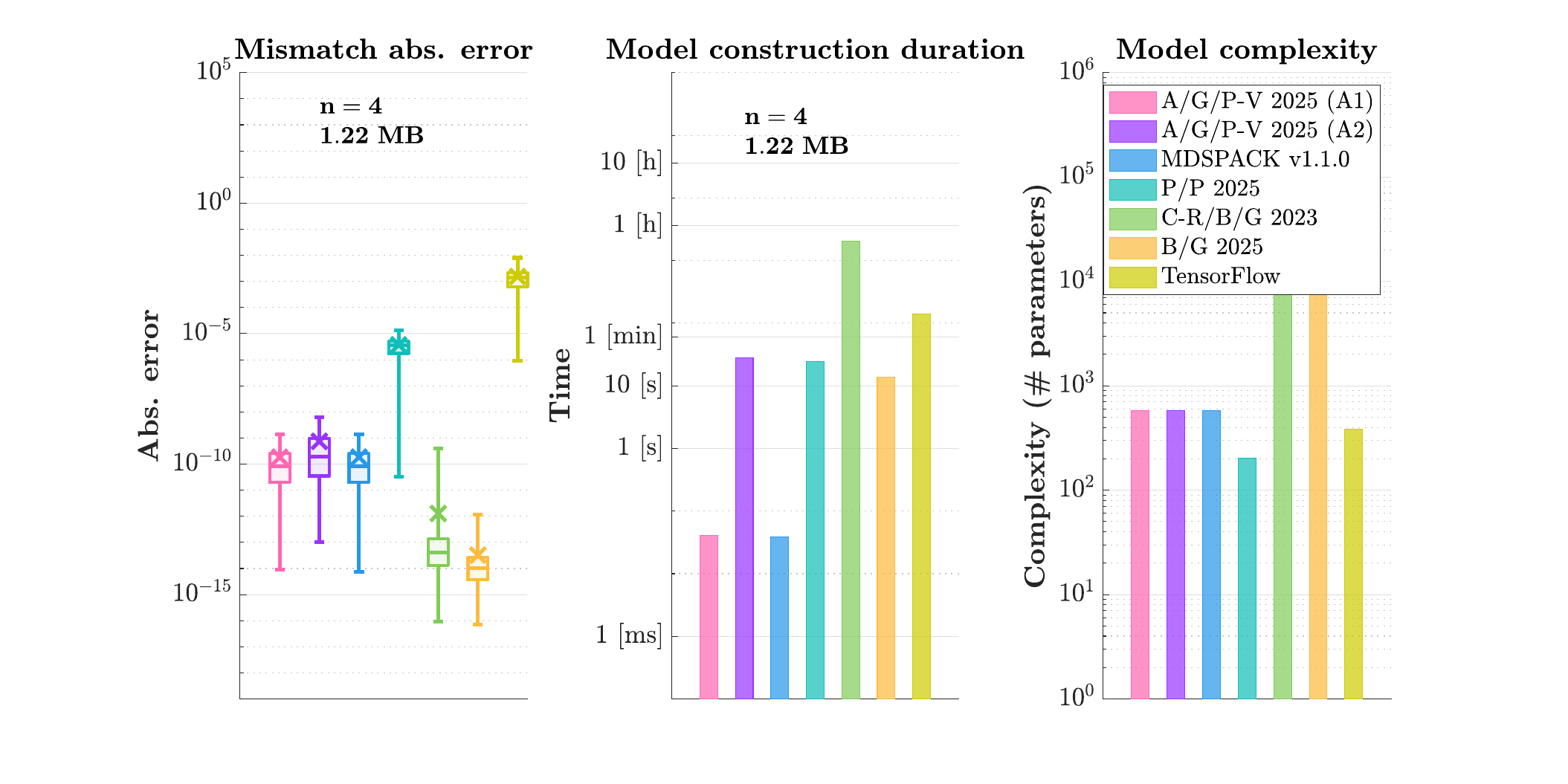} \caption{Function \#37: graphical view of the best model performances.} \end{figure}\begin{figure}[H] \centering  \includegraphics[width=\textwidth]{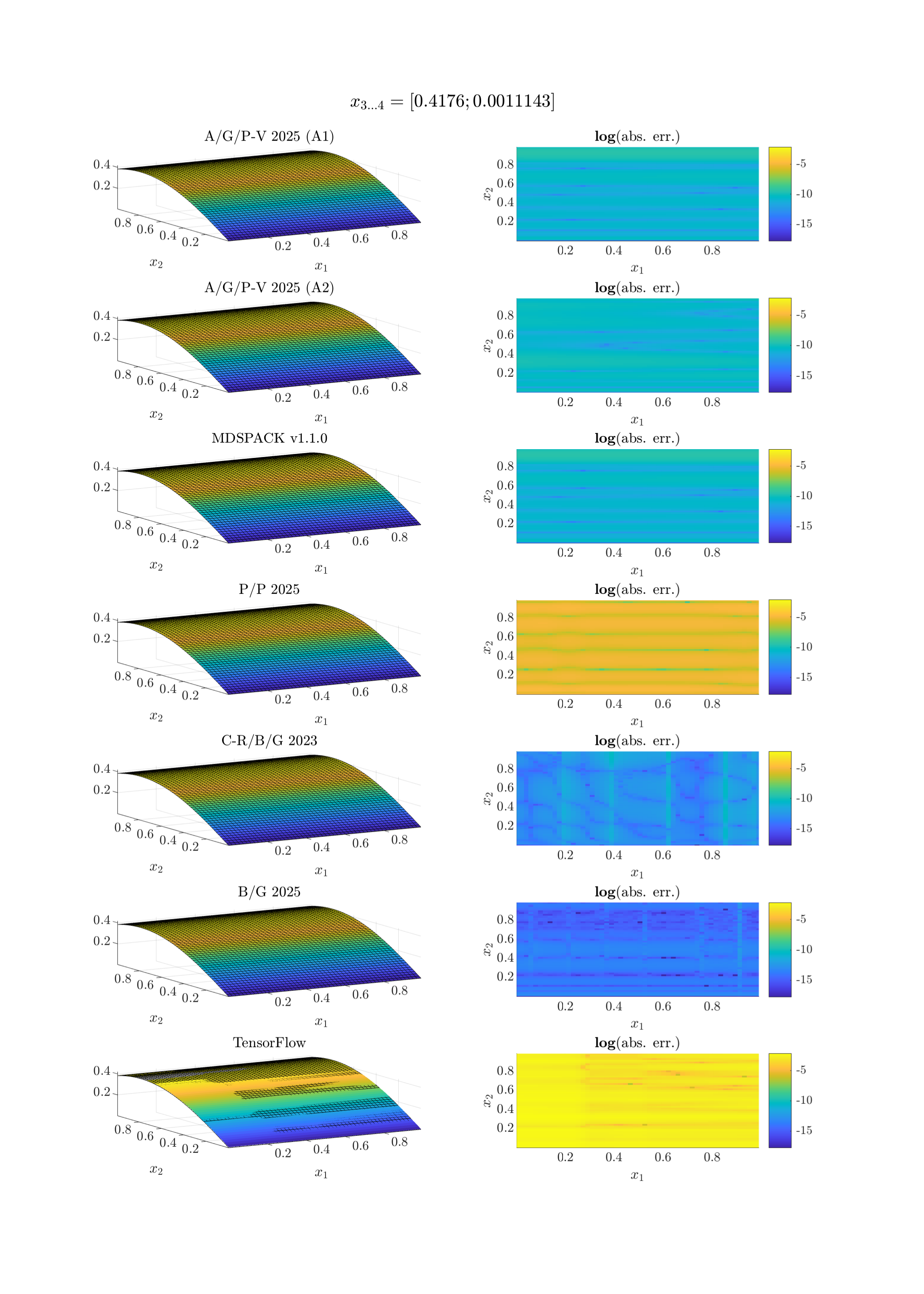} \caption{Function \#37: left side, evaluation of the original (mesh) vs. approximated (coloured surface) and right side, absolute errors (in log-scale).} \end{figure}\subsubsection{mLF detailed informations (M1)} \noindent \textbf{Right interpolation points}: $k_l=\left(\begin{array}{cccc} 2 & 6 & 2 & 4 \end{array}\right)$, where $l=1,\cdots,\ord$.$$ \begin{array}{rcl}\lan{1} &\in& \IC^{2} \text{ , linearly spaced between bounds}\\\lan{2} &\in& \IC^{6} \text{ , linearly spaced between bounds}\\\lan{3} &\in& \IC^{2} \text{ , linearly spaced between bounds}\\\lan{4} &\in& \IC^{4} \text{ , linearly spaced between bounds}\\\end{array} $$\noindent \textbf{$\ord$-D Loewner matrix, barycentric weights and Lagrangian basis}:$$ \begin{array}{rcl}\IL & \in & \IC^{96 \times 96}\\\bc & \in & \IC^{96}\\\bw & \in & \IC^{96}\\\bc\odot \bw & \in & \IC^{96}\\\mathbf{Lag}(\var{1},\var{2},\var{3},\var{4}) & \in & \IC^{96}\\\end{array} $$

\newpage \subsection{Function \#38 (${\ord=3}$ variables, tensor size: 1.65 \textbf{MB})} $$\frac{\var{1}^9 \var{2}^7 + \var{1}^3 + 5 \var{3}^2}{5 \var{1}^4 + 4 \var{1}^2 + \var{3}\var{2}^3 + 1}$$ \subsubsection{Setup and results overview}\begin{itemize}\item Reference: A.C. Antoulas presentation, [none]\item Domain: $\mathbb{R}$\item Tensor size: 1.65 \textbf{MB} ($60^{3}$ points)\item Bounds: $ \left(\begin{array}{cc} -\frac{11}{10} & \frac{11}{10} \end{array}\right) \times \left(\begin{array}{cc} -\frac{11}{10} & \frac{11}{10} \end{array}\right) \times \left(\begin{array}{cc} -\frac{11}{10} & \frac{11}{10} \end{array}\right)$ \end{itemize} \begin{table}[H] \centering \begin{tabular}{llllll}
$\#$ & Alg. & Parameters & Dim. & CPU [s] & RMSE \\ 
\hline 
$\mathbf{\#38}$ & A/G/P-V 2025 (A1) & $1 \cdot 10^{-06},3$ & $1.2 \cdot 10^{03}$ & $0.08$ & $0.0023$ \\ 
 & A/G/P-V 2025 (A2) & $1 \cdot 10^{-15},2$ & $3.8 \cdot 10^{02}$ & $31$ & $13$ \\ 
 & MDSPACK v1.1.0 & $1 \cdot 10^{-12},6$ & $1.2 \cdot 10^{03}$ & $\mathbf{0.05}$ & $0.0023$ \\ 
 & P/P 2025 & $1,1,50,0.01,10,6,21$ & $7.6 \cdot 10^{02}$ & $90$ & $1.3$ \\ 
 & C-R/B/G 2023 & $1 \cdot 10^{-09},20$ & $7.9 \cdot 10^{03}$ & $2.2 \cdot 10^{02}$ & $\mathbf{4.6 \cdot 10^{-09}}$ \\ 
 & B/G 2025 & $1 \cdot 10^{-06},20,3$ & $1.4 \cdot 10^{04}$ & $1 \cdot 10^{02}$ & $1.8$ \\ 
 & TensorFlow & $$ & $\mathbf{3.2 \cdot 10^{02}}$ & $9.1 \cdot 10^{02}$ & $2.1$ \\ 
\hline 
\end{tabular} \caption{Function \#38: best model configuration and performances per methods.} \end{table}\begin{figure}[H] \centering  \includegraphics[width=\textwidth]{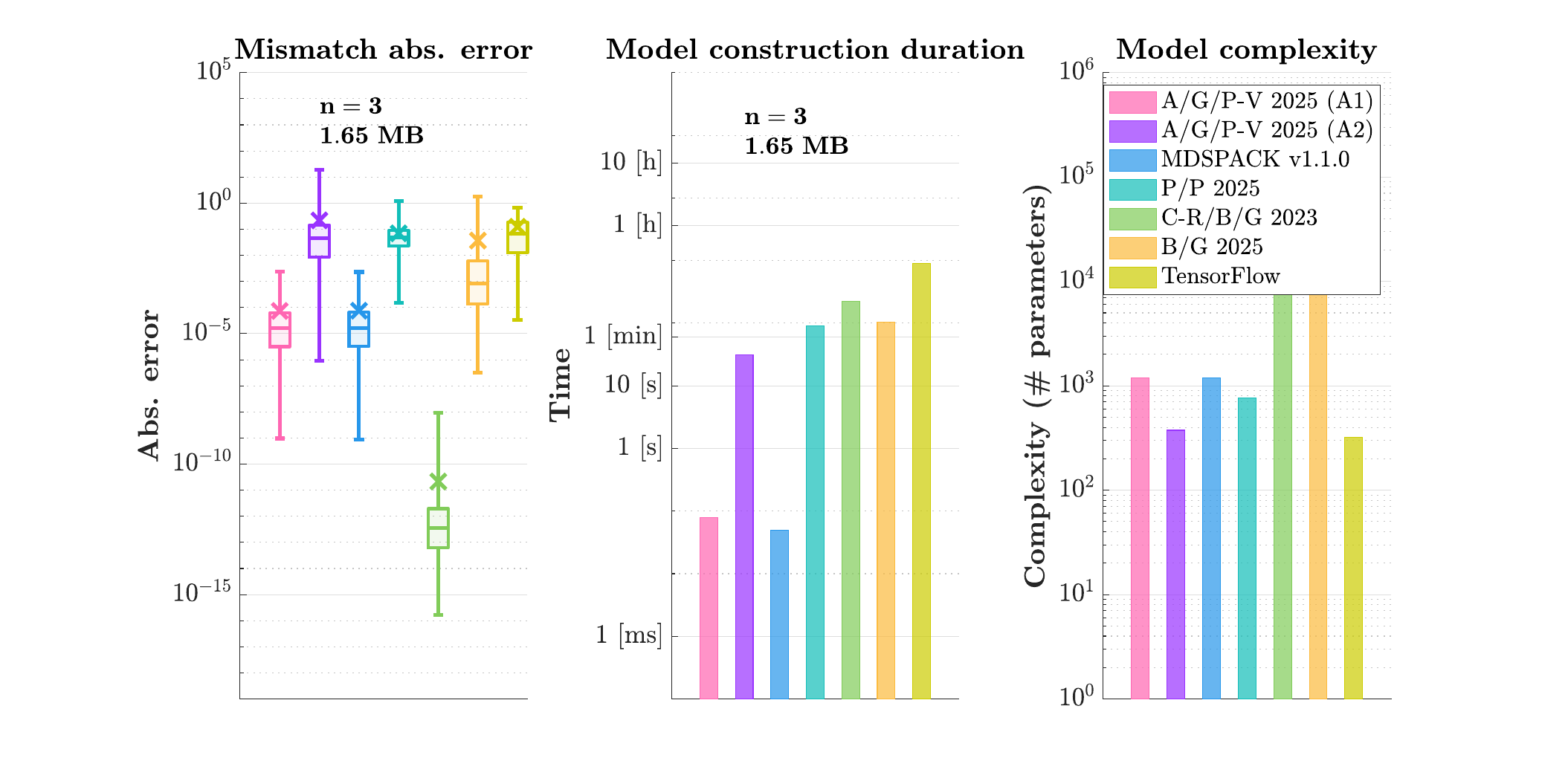} \caption{Function \#38: graphical view of the best model performances.} \end{figure}\begin{figure}[H] \centering  \includegraphics[width=\textwidth]{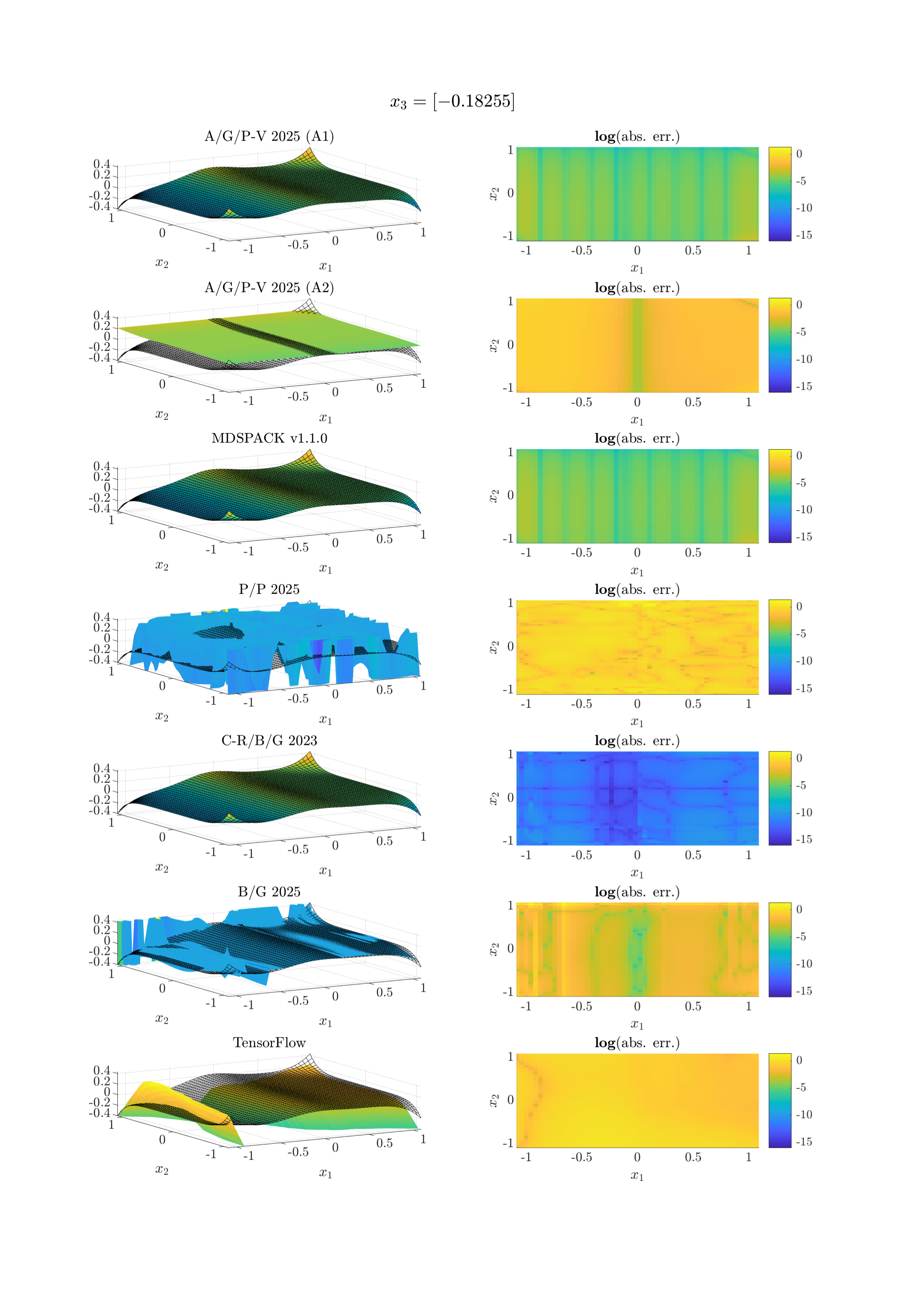} \caption{Function \#38: left side, evaluation of the original (mesh) vs. approximated (coloured surface) and right side, absolute errors (in log-scale).} \end{figure}\subsubsection{mLF detailed informations (M1)} \noindent \textbf{Right interpolation points}: $k_l=\left(\begin{array}{ccc} 10 & 8 & 3 \end{array}\right)$, where $l=1,\cdots,\ord$.$$ \begin{array}{rcl}\lan{1} &\in& \IC^{10} \text{ , linearly spaced between bounds}\\\lan{2} &\in& \IC^{8} \text{ , linearly spaced between bounds}\\\lan{3} &\in& \IC^{3} \text{ , linearly spaced between bounds}\\\end{array} $$\noindent \textbf{$\ord$-D Loewner matrix, barycentric weights and Lagrangian basis}:$$ \begin{array}{rcl}\IL & \in & \IC^{240 \times 240}\\\bc & \in & \IC^{240}\\\bw & \in & \IC^{240}\\\bc\odot \bw & \in & \IC^{240}\\\mathbf{Lag}(\var{1},\var{2},\var{3}) & \in & \IC^{240}\\\end{array} $$

\newpage \subsection{Function \#39 (${\ord=3}$ variables, tensor size: 500 \textbf{KB})} $$\frac{\var{3}+\var{1}^4}{\var{1}^3+\var{2}^2+1}$$ \subsubsection{Setup and results overview}\begin{itemize}\item Reference: Personal communication, [none]\item Domain: $\mathbb{R}$\item Tensor size: 500 \textbf{KB} ($40^{3}$ points)\item Bounds: $ \left(\begin{array}{cc} \frac{1}{10} & 10 \end{array}\right) \times \left(\begin{array}{cc} \frac{1}{10} & 10 \end{array}\right) \times \left(\begin{array}{cc} \frac{1}{10} & 10 \end{array}\right)$ \end{itemize} \begin{table}[H] \centering \begin{tabular}{llllll}
$\#$ & Alg. & Parameters & Dim. & CPU [s] & RMSE \\ 
\hline 
$\mathbf{\#39}$ & A/G/P-V 2025 (A1) & $0.0001,2$ & $\mathbf{1.5 \cdot 10^{02}}$ & $0.026$ & $9.2 \cdot 10^{-12}$ \\ 
 & A/G/P-V 2025 (A2) & $1 \cdot 10^{-15},3$ & $1.5 \cdot 10^{02}$ & $4.2$ & $6.3 \cdot 10^{-11}$ \\ 
 & MDSPACK v1.1.0 & $0.0001,2$ & $1.5 \cdot 10^{02}$ & $\mathbf{0.021}$ & $4.6 \cdot 10^{-12}$ \\ 
 & P/P 2025 & $1,1,50,0.01,6,12,13$ & $3.9 \cdot 10^{02}$ & $17$ & $0.0041$ \\ 
 & C-R/B/G 2023 & $0.001,20$ & $6 \cdot 10^{02}$ & $5.4$ & $\mathbf{3.6 \cdot 10^{-14}}$ \\ 
 & B/G 2025 & $0.001,20,3$ & $4 \cdot 10^{02}$ & $4$ & $5.4 \cdot 10^{-14}$ \\ 
 & TensorFlow & $$ & $3.2 \cdot 10^{02}$ & $2.8 \cdot 10^{02}$ & $1.2$ \\ 
\hline 
\end{tabular} \caption{Function \#39: best model configuration and performances per methods.} \end{table}\begin{figure}[H] \centering  \includegraphics[width=\textwidth]{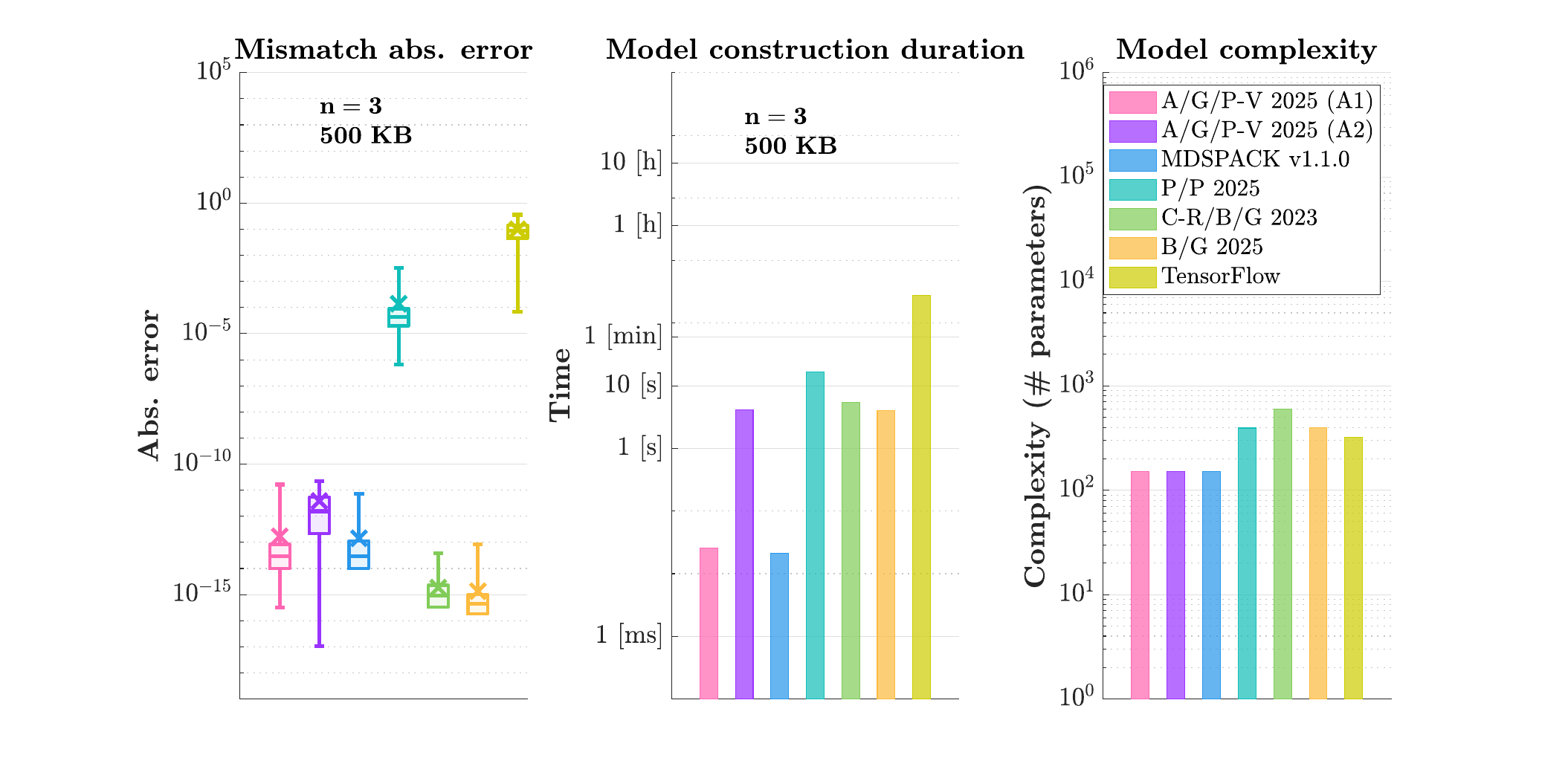} \caption{Function \#39: graphical view of the best model performances.} \end{figure}\begin{figure}[H] \centering  \includegraphics[width=\textwidth]{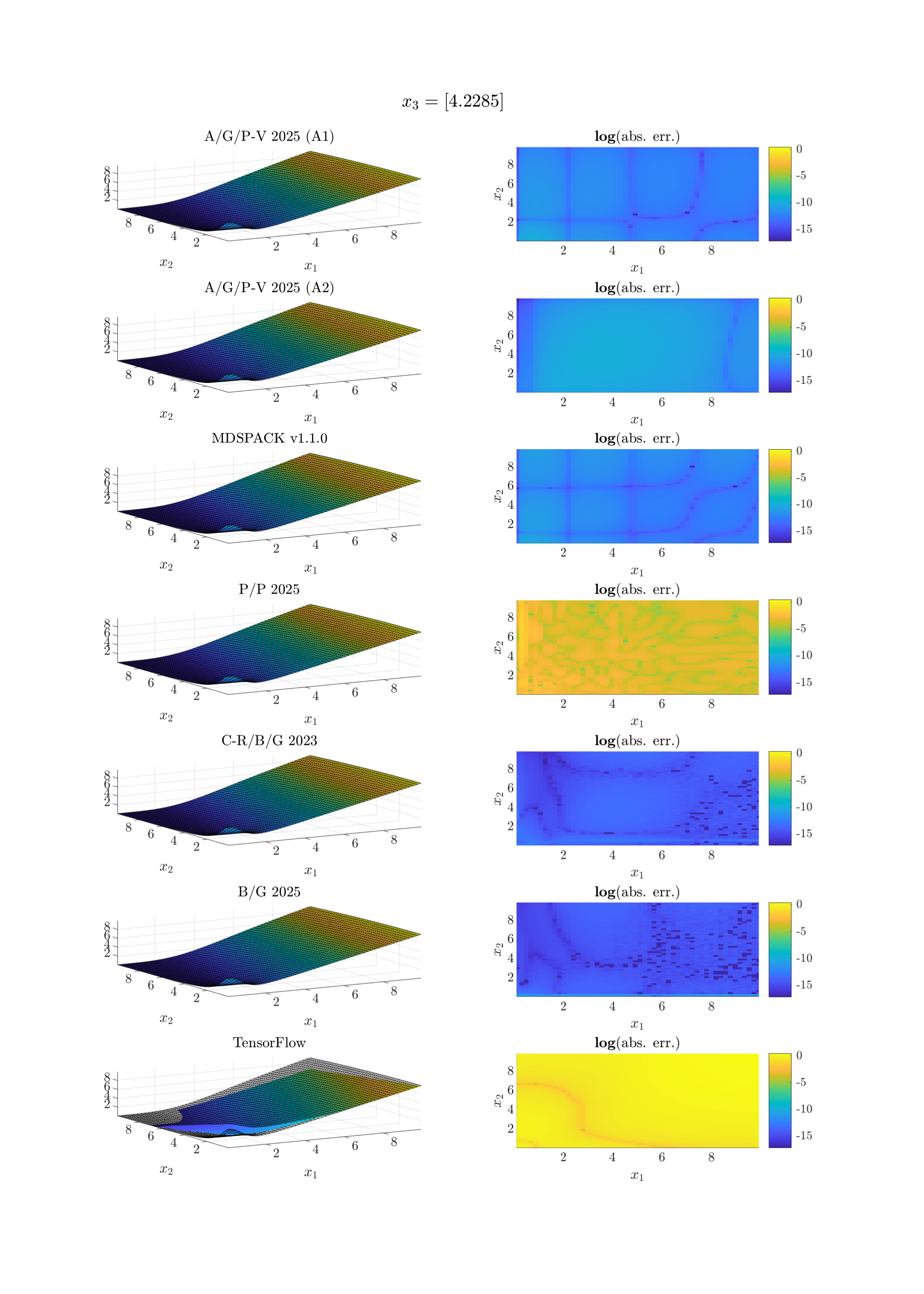} \caption{Function \#39: left side, evaluation of the original (mesh) vs. approximated (coloured surface) and right side, absolute errors (in log-scale).} \end{figure}\subsubsection{mLF detailed informations (M1)} \noindent \textbf{Right interpolation points} ($k_l=\left(\begin{array}{ccc} 5 & 3 & 2 \end{array}\right)$, where $l=1,\cdots,\ord$):$$ \begin{array}{rcl}\lan{1} &=& \left(\begin{array}{ccccc} \frac{1}{10} & \frac{83}{38} & \frac{91}{19} & \frac{281}{38} & 10 \end{array}\right)\\\lan{2} &=& \left(\begin{array}{ccc} \frac{1}{10} & \frac{91}{19} & 10 \end{array}\right)\\\lan{3} &=& \left(\begin{array}{cc} \frac{1}{10} & 10 \end{array}\right)\\\end{array} $$\noindent \textbf{Lagrangian weights}: $$\left(\begin{array}{ccc} \bc & \bw & \bc\odot\bw\\ 0.0008124 & 0.09901 & 8.044 \cdot 10^{-5}\\ -0.0008124 & 9.891 & -0.008036\\ -0.03655 & 0.004181 & -0.0001528\\ 0.03655 & 0.4177 & 0.01527\\ 0.07305 & 0.0009911 & 7.24 \cdot 10^{-5}\\ -0.07305 & 0.09901 & -0.007232\\ -0.02932 & 2.0 & -0.05864\\ 0.02932 & 2.866 & 0.08403\\ 0.1675 & 0.6653 & 0.1114\\ -0.1675 & 0.9535 & -0.1597\\ -0.2572 & 0.2052 & -0.05277\\ 0.2572 & 0.294 & 0.07563\\ 0.3792 & 4.747 & 1.8\\ -0.3792 & 4.836 & -1.834\\ -0.8695 & 3.933 & -3.42\\ 0.8695 & 4.007 & 3.484\\ 0.6491 & 2.496 & 1.62\\ -0.6491 & 2.543 & -1.65\\ -0.8913 & 7.377 & -6.574\\ 0.8913 & 7.401 & 6.596\\ 1.789 & 6.982 & 12.49\\ -1.789 & 7.005 & -12.53\\ -1.0 & 5.917 & -5.917\\ 1.0 & 5.937 & 5.937\\ 0.5406 & 9.99 & 5.4\\ -0.5406 & 10.0 & -5.406\\ -1.051 & 9.766 & -10.26\\ 1.051 & 9.776 & 10.27\\ 0.5351 & 9.083 & 4.86\\ -0.5351 & 9.092 & -4.865 \end{array}\right)$$\noindent \textbf{Lagrangian form} (basis, numerator and denominator coefficients):$$\left(\begin{array}{ccc}\mathcal{B}_\textrm{lag}(\var{1},\var{2},\var{3}) & \bN_\textrm{lag} &\bD_\textrm{lag}\end{array}\right) =$$ $$\left(\begin{array}{ccc} \left(\var{1}-0.1\right)\,\left(\var{2}-0.1\right)\,\left(\var{3}-0.1\right) & 8.044 \cdot 10^{-5} & 0.0008124\\ \left(\var{3}-10.0\right)\,\left(\var{1}-0.1\right)\,\left(\var{2}-0.1\right) & -0.008036 & -0.0008124\\ \left(\var{2}-4.789\right)\,\left(\var{1}-0.1\right)\,\left(\var{3}-0.1\right) & -0.0001528 & -0.03655\\ \left(\var{2}-4.789\right)\,\left(\var{3}-10.0\right)\,\left(\var{1}-0.1\right) & 0.01527 & 0.03655\\ \left(\var{2}-10.0\right)\,\left(\var{1}-0.1\right)\,\left(\var{3}-0.1\right) & 7.24 \cdot 10^{-5} & 0.07305\\ \left(\var{2}-10.0\right)\,\left(\var{3}-10.0\right)\,\left(\var{1}-0.1\right) & -0.007232 & -0.07305\\ \left(\var{1}-2.184\right)\,\left(\var{2}-0.1\right)\,\left(\var{3}-0.1\right) & -0.05864 & -0.02932\\ \left(\var{1}-2.184\right)\,\left(\var{3}-10.0\right)\,\left(\var{2}-0.1\right) & 0.08403 & 0.02932\\ \left(\var{1}-2.184\right)\,\left(\var{2}-4.789\right)\,\left(\var{3}-0.1\right) & 0.1114 & 0.1675\\ \left(\var{1}-2.184\right)\,\left(\var{2}-4.789\right)\,\left(\var{3}-10.0\right) & -0.1597 & -0.1675\\ \left(\var{1}-2.184\right)\,\left(\var{2}-10.0\right)\,\left(\var{3}-0.1\right) & -0.05277 & -0.2572\\ \left(\var{1}-2.184\right)\,\left(\var{2}-10.0\right)\,\left(\var{3}-10.0\right) & 0.07563 & 0.2572\\ \left(\var{1}-4.789\right)\,\left(\var{2}-0.1\right)\,\left(\var{3}-0.1\right) & 1.8 & 0.3792\\ \left(\var{1}-4.789\right)\,\left(\var{3}-10.0\right)\,\left(\var{2}-0.1\right) & -1.834 & -0.3792\\ \left(\var{1}-4.789\right)\,\left(\var{2}-4.789\right)\,\left(\var{3}-0.1\right) & -3.42 & -0.8695\\ \left(\var{1}-4.789\right)\,\left(\var{2}-4.789\right)\,\left(\var{3}-10.0\right) & 3.484 & 0.8695\\ \left(\var{1}-4.789\right)\,\left(\var{2}-10.0\right)\,\left(\var{3}-0.1\right) & 1.62 & 0.6491\\ \left(\var{1}-4.789\right)\,\left(\var{2}-10.0\right)\,\left(\var{3}-10.0\right) & -1.65 & -0.6491\\ \left(\var{2}-0.1\right)\,\left(\var{3}-0.1\right)\,\left(\var{1}-7.395\right) & -6.574 & -0.8913\\ \left(\var{3}-10.0\right)\,\left(\var{2}-0.1\right)\,\left(\var{1}-7.395\right) & 6.596 & 0.8913\\ \left(\var{2}-4.789\right)\,\left(\var{3}-0.1\right)\,\left(\var{1}-7.395\right) & 12.49 & 1.789\\ \left(\var{2}-4.789\right)\,\left(\var{3}-10.0\right)\,\left(\var{1}-7.395\right) & -12.53 & -1.789\\ \left(\var{2}-10.0\right)\,\left(\var{3}-0.1\right)\,\left(\var{1}-7.395\right) & -5.917 & -1.0\\ \left(\var{2}-10.0\right)\,\left(\var{3}-10.0\right)\,\left(\var{1}-7.395\right) & 5.937 & 1.0\\ \left(\var{1}-10.0\right)\,\left(\var{2}-0.1\right)\,\left(\var{3}-0.1\right) & 5.4 & 0.5406\\ \left(\var{1}-10.0\right)\,\left(\var{3}-10.0\right)\,\left(\var{2}-0.1\right) & -5.406 & -0.5406\\ \left(\var{2}-4.789\right)\,\left(\var{1}-10.0\right)\,\left(\var{3}-0.1\right) & -10.26 & -1.051\\ \left(\var{2}-4.789\right)\,\left(\var{1}-10.0\right)\,\left(\var{3}-10.0\right) & 10.27 & 1.051\\ \left(\var{1}-10.0\right)\,\left(\var{2}-10.0\right)\,\left(\var{3}-0.1\right) & 4.86 & 0.5351\\ \left(\var{1}-10.0\right)\,\left(\var{2}-10.0\right)\,\left(\var{3}-10.0\right) & -4.865 & -0.5351 \end{array}\right).$$\noindent The corresponding function is:$$\begin{array}{rcl}\bG_{\textrm{lag}}(\var{1},\var{2},\var{3}) &=& \dfrac{\bn_{\textrm{lag}}(\var{1},\var{2},\var{3})}{\bd_{\textrm{lag}}(\var{1},\var{2},\var{3})}\\ && \\&=& \dfrac{\sum_{\textrm{row}} \bN_\textrm{lag} \odot\mathcal{B}^{-1}_\textrm{lag}(\var{1},\var{2},\var{3})}{\sum_{\textrm{row}} \bD_\textrm{lag} \odot\mathcal{B}^{-1}_\textrm{lag}(\var{1},\var{2},\var{3})}, \end{array}$$\noindent where,\\$\bn_{\textrm{lag}}(\var{1},\var{2},\var{3}) = 1.775 \cdot 10^{-9}\,\var{1}+8.593 \cdot 10^{-11}\,\var{2}+1.0\,\var{3}-5.544 \cdot 10^{-11}\,{\var{1}}^2\,{\var{2}}^2+1.244 \cdot 10^{-11}\,{\var{1}}^3\,{\var{2}}^2-8.345 \cdot 10^{-13}\,{\var{1}}^4\,{\var{2}}^2-9.28 \cdot 10^{-10}\,\var{1}\,\var{2}-1.822 \cdot 10^{-10}\,\var{1}\,\var{3}-8.676 \cdot 10^{-12}\,\var{2}\,\var{3}+7.352 \cdot 10^{-11}\,\var{1}\,{\var{2}}^2+7.028 \cdot 10^{-10}\,{\var{1}}^2\,\var{2}+1.421 \cdot 10^{-10}\,{\var{1}}^2\,\var{3}-1.592 \cdot 10^{-10}\,{\var{1}}^3\,\var{2}-3.376 \cdot 10^{-11}\,{\var{1}}^3\,\var{3}+6.842 \cdot 10^{-13}\,{\var{2}}^2\,\var{3}+1.088 \cdot 10^{-11}\,{\var{1}}^4\,\var{2}+2.485 \cdot 10^{-12}\,{\var{1}}^4\,\var{3}-1.375 \cdot 10^{-9}\,{\var{1}}^2+3.215 \cdot 10^{-10}\,{\var{1}}^3-6.81 \cdot 10^{-12}\,{\var{2}}^2+1.0\,{\var{1}}^4-7.387 \cdot 10^{-12}\,\var{1}\,{\var{2}}^2\,\var{3}-7.097 \cdot 10^{-11}\,{\var{1}}^2\,\var{2}\,\var{3}+1.609 \cdot 10^{-11}\,{\var{1}}^3\,\var{2}\,\var{3}-1.099 \cdot 10^{-12}\,{\var{1}}^4\,\var{2}\,\var{3}+5.571 \cdot 10^{-12}\,{\var{1}}^2\,{\var{2}}^2\,\var{3}-1.251 \cdot 10^{-12}\,{\var{1}}^3\,{\var{2}}^2\,\var{3}+8.397 \cdot 10^{-14}\,{\var{1}}^4\,{\var{2}}^2\,\var{3}+9.37 \cdot 10^{-11}\,\var{1}\,\var{2}\,\var{3}-1.641 \cdot 10^{-10}$ \\~~\\$\bd_{\textrm{lag}}(\var{1},\var{2},\var{3}) = 1.497 \cdot 10^{-10}\,\var{1}+8.93 \cdot 10^{-12}\,\var{2}+1.381 \cdot 10^{-12}\,\var{3}-6.046 \cdot 10^{-12}\,{\var{1}}^2\,{\var{2}}^2+1.343 \cdot 10^{-12}\,{\var{1}}^3\,{\var{2}}^2-8.942 \cdot 10^{-14}\,{\var{1}}^4\,{\var{2}}^2-9.663 \cdot 10^{-11}\,\var{1}\,\var{2}-1.497 \cdot 10^{-11}\,\var{1}\,\var{3}-8.998 \cdot 10^{-13}\,\var{2}\,\var{3}+8.784 \cdot 10^{-12}\,\var{1}\,{\var{2}}^2+7.327 \cdot 10^{-11}\,{\var{1}}^2\,\var{2}+1.179 \cdot 10^{-11}\,{\var{1}}^2\,\var{3}-1.662 \cdot 10^{-11}\,{\var{1}}^3\,\var{2}-2.873 \cdot 10^{-12}\,{\var{1}}^3\,\var{3}+7.183 \cdot 10^{-14}\,{\var{2}}^2\,\var{3}+1.131 \cdot 10^{-12}\,{\var{1}}^4\,\var{2}+2.191 \cdot 10^{-13}\,{\var{1}}^4\,\var{3}-1.202 \cdot 10^{-10}\,{\var{1}}^2+1.0\,{\var{1}}^3+1.0\,{\var{2}}^2-2.123 \cdot 10^{-12}\,{\var{1}}^4-7.788 \cdot 10^{-13}\,\var{1}\,{\var{2}}^2\,\var{3}-7.385 \cdot 10^{-12}\,{\var{1}}^2\,\var{2}\,\var{3}+1.675 \cdot 10^{-12}\,{\var{1}}^3\,\var{2}\,\var{3}-1.141 \cdot 10^{-13}\,{\var{1}}^4\,\var{2}\,\var{3}+5.906 \cdot 10^{-13}\,{\var{1}}^2\,{\var{2}}^2\,\var{3}-1.335 \cdot 10^{-13}\,{\var{1}}^3\,{\var{2}}^2\,\var{3}+8.941 \cdot 10^{-15}\,{\var{1}}^4\,{\var{2}}^2\,\var{3}+9.739 \cdot 10^{-12}\,\var{1}\,\var{2}\,\var{3}+1.0$ \\~~\\\noindent \textbf{Monomial form} (basis, numerator and denominator coefficients - entries $<10^{-12}$ removed):$$\left(\begin{array}{ccc}\mathcal{B}_\textrm{mon}(\var{1},\var{2},\var{3}) & \bN_\textrm{mon} &\bD_\textrm{mon}\end{array}\right) =$$ $$\left(\begin{array}{ccc} {\var{1}}^4\,{\var{2}}^2\,\var{3} & 0 & 0\\ {\var{1}}^4\,{\var{2}}^2 & 0 & 0\\ {\var{1}}^4\,\var{2}\,\var{3} & 1.099 \cdot 10^{-12} & 0\\ {\var{1}}^4\,\var{2} & -1.088 \cdot 10^{-11} & -1.128 \cdot 10^{-12}\\ {\var{1}}^4\,\var{3} & -2.485 \cdot 10^{-12} & 0\\ {\var{1}}^4 & -1.0 & 2.123 \cdot 10^{-12}\\ {\var{1}}^3\,{\var{2}}^2\,\var{3} & 1.251 \cdot 10^{-12} & 0\\ {\var{1}}^3\,{\var{2}}^2 & -1.244 \cdot 10^{-11} & -1.334 \cdot 10^{-12}\\ {\var{1}}^3\,\var{2}\,\var{3} & -1.609 \cdot 10^{-11} & -1.671 \cdot 10^{-12}\\ {\var{1}}^3\,\var{2} & 1.592 \cdot 10^{-10} & 1.657 \cdot 10^{-11}\\ {\var{1}}^3\,\var{3} & 3.377 \cdot 10^{-11} & 2.873 \cdot 10^{-12}\\ {\var{1}}^3 & -3.215 \cdot 10^{-10} & -1.0\\ {\var{1}}^2\,{\var{2}}^2\,\var{3} & -5.571 \cdot 10^{-12} & 0\\ {\var{1}}^2\,{\var{2}}^2 & 5.543 \cdot 10^{-11} & 6.001 \cdot 10^{-12}\\ {\var{1}}^2\,\var{2}\,\var{3} & 7.097 \cdot 10^{-11} & 7.363 \cdot 10^{-12}\\ {\var{1}}^2\,\var{2} & -7.028 \cdot 10^{-10} & -7.305 \cdot 10^{-11}\\ {\var{1}}^2\,\var{3} & -1.421 \cdot 10^{-10} & -1.179 \cdot 10^{-11}\\ {\var{1}}^2 & 1.375 \cdot 10^{-9} & 1.202 \cdot 10^{-10}\\ \var{1}\,{\var{2}}^2\,\var{3} & 7.388 \cdot 10^{-12} & 0\\ \var{1}\,{\var{2}}^2 & -7.353 \cdot 10^{-11} & -8.724 \cdot 10^{-12}\\ \var{1}\,\var{2}\,\var{3} & -9.37 \cdot 10^{-11} & -9.708 \cdot 10^{-12}\\ \var{1}\,\var{2} & 9.28 \cdot 10^{-10} & 9.632 \cdot 10^{-11}\\ \var{1}\,\var{3} & 1.822 \cdot 10^{-10} & 1.497 \cdot 10^{-11}\\ \var{1} & -1.775 \cdot 10^{-9} & -1.497 \cdot 10^{-10}\\ {\var{2}}^2\,\var{3} & 0 & 0\\ {\var{2}}^2 & 6.81 \cdot 10^{-12} & -1.0\\ \var{2}\,\var{3} & 8.676 \cdot 10^{-12} & 0\\ \var{2} & -8.594 \cdot 10^{-11} & -8.901 \cdot 10^{-12}\\ \var{3} & -1.0 & -1.381 \cdot 10^{-12}\\ 1.0 & 1.641 \cdot 10^{-10} & -1.0 \end{array}\right)$$\noindent The corresponding function is:$$\begin{array}{rcl}\bG_{\textrm{mon}}(\var{1},\var{2},\var{3}) &=& \dfrac{\bn_{\textrm{mon}}(\var{1},\var{2},\var{3})}{\bd_{\textrm{mon}}(\var{1},\var{2},\var{3})}\\ && \\&=& \dfrac{\sum_{\textrm{row}} \bN_\textrm{mon} \odot \mathcal{B}_\textrm{mon}(\var{1},\var{2},\var{3})}{\sum_{\textrm{row}} \bD_\textrm{mon} \odot\mathcal{B}_\textrm{mon}(\var{1},\var{2},\var{3})},  \end{array}$$\noindent where,\\$\bn_{\textrm{mon}}(\var{1},\var{2},\var{3}) = 1.775 \cdot 10^{-9}\,\var{1}+8.594 \cdot 10^{-11}\,\var{2}+1.0\,\var{3}-5.543 \cdot 10^{-11}\,{\var{1}}^2\,{\var{2}}^2+1.244 \cdot 10^{-11}\,{\var{1}}^3\,{\var{2}}^2-9.28 \cdot 10^{-10}\,\var{1}\,\var{2}-1.822 \cdot 10^{-10}\,\var{1}\,\var{3}-8.676 \cdot 10^{-12}\,\var{2}\,\var{3}+7.353 \cdot 10^{-11}\,\var{1}\,{\var{2}}^2+7.028 \cdot 10^{-10}\,{\var{1}}^2\,\var{2}+1.421 \cdot 10^{-10}\,{\var{1}}^2\,\var{3}-1.592 \cdot 10^{-10}\,{\var{1}}^3\,\var{2}-3.377 \cdot 10^{-11}\,{\var{1}}^3\,\var{3}+1.088 \cdot 10^{-11}\,{\var{1}}^4\,\var{2}+2.485 \cdot 10^{-12}\,{\var{1}}^4\,\var{3}-1.375 \cdot 10^{-9}\,{\var{1}}^2+3.215 \cdot 10^{-10}\,{\var{1}}^3-6.81 \cdot 10^{-12}\,{\var{2}}^2+1.0\,{\var{1}}^4-7.388 \cdot 10^{-12}\,\var{1}\,{\var{2}}^2\,\var{3}-7.097 \cdot 10^{-11}\,{\var{1}}^2\,\var{2}\,\var{3}+1.609 \cdot 10^{-11}\,{\var{1}}^3\,\var{2}\,\var{3}-1.099 \cdot 10^{-12}\,{\var{1}}^4\,\var{2}\,\var{3}+5.571 \cdot 10^{-12}\,{\var{1}}^2\,{\var{2}}^2\,\var{3}-1.251 \cdot 10^{-12}\,{\var{1}}^3\,{\var{2}}^2\,\var{3}+9.37 \cdot 10^{-11}\,\var{1}\,\var{2}\,\var{3}-1.641 \cdot 10^{-10}$ \\~~\\$\bd_{\textrm{mon}}(\var{1},\var{2},\var{3}) = 1.497 \cdot 10^{-10}\,\var{1}+8.901 \cdot 10^{-12}\,\var{2}+1.381 \cdot 10^{-12}\,\var{3}-6.001 \cdot 10^{-12}\,{\var{1}}^2\,{\var{2}}^2+1.334 \cdot 10^{-12}\,{\var{1}}^3\,{\var{2}}^2-9.632 \cdot 10^{-11}\,\var{1}\,\var{2}-1.497 \cdot 10^{-11}\,\var{1}\,\var{3}+8.724 \cdot 10^{-12}\,\var{1}\,{\var{2}}^2+7.305 \cdot 10^{-11}\,{\var{1}}^2\,\var{2}+1.179 \cdot 10^{-11}\,{\var{1}}^2\,\var{3}-1.657 \cdot 10^{-11}\,{\var{1}}^3\,\var{2}-2.873 \cdot 10^{-12}\,{\var{1}}^3\,\var{3}+1.128 \cdot 10^{-12}\,{\var{1}}^4\,\var{2}-1.202 \cdot 10^{-10}\,{\var{1}}^2+1.0\,{\var{1}}^3+1.0\,{\var{2}}^2-2.123 \cdot 10^{-12}\,{\var{1}}^4-7.363 \cdot 10^{-12}\,{\var{1}}^2\,\var{2}\,\var{3}+1.671 \cdot 10^{-12}\,{\var{1}}^3\,\var{2}\,\var{3}+9.708 \cdot 10^{-12}\,\var{1}\,\var{2}\,\var{3}+1.0$ \\~~\\\begin{landscape} \noindent \textbf{KST equivalent decoupling pattern} (Barycentric weights $\bc^{\var{l}}$): $$\begin{array}{rclll}\var{3}&: & \left(\begin{array}{ccccccccccccccc} -1.0 & -1.0 & -1.0 & -1.0 & -1.0 & -1.0 & -1.0 & -1.0 & -1.0 & -1.0 & -1.0 & -1.0 & -1.0 & -1.0 & -1.0\\ 1.0 & 1.0 & 1.0 & 1.0 & 1.0 & 1.0 & 1.0 & 1.0 & 1.0 & 1.0 & 1.0 & 1.0 & 1.0 & 1.0 & 1.0 \end{array}\right)& \textrm{vec}(.) & := \textbf{Bary}(\var{3}) \\\var{2}&: & \left(\begin{array}{ccccc} 0.01112 & 0.114 & 0.5842 & 0.8913 & 1.01\\ -0.5004 & -0.651 & -1.34 & -1.789 & -1.963\\ 1.0 & 1.0 & 1.0 & 1.0 & 1.0 \end{array}\right)& \textrm{vec}(.) \otimes \bone_{k_{3}} & := \textbf{Bary}(\var{2}) \\\var{1}&: & \left(\begin{array}{c} -0.07305\\ 0.2572\\ -0.6491\\ 1.0\\ -0.5351 \end{array}\right)& \textrm{vec}(.) \otimes \bone_{k_{3}k_{2}} & := \textbf{Bary}(\var{1}) \\\end{array}$$\end{landscape} ~\\ Then, with the above notations, one defines the following univariate vector functions:  $$ \left\{ \begin{array}{rcl}\bPhi_{1}(\var{1}) &:=& \textbf{Bary}(\var{1}) \odot \mathbf{Lag}(\var{1}) \\\bPhi_{2}(\var{2}) &:=& \textbf{Bary}(\var{2}) \odot \mathbf{Lag}(\var{2}) \\\bPhi_{3}(\var{3}) &:=& \textbf{Bary}(\var{3}) \odot \mathbf{Lag}(\var{3}) \\\end{array} \right. $$\noindent The corresponding function is:$$\begin{array}{rcl}\bG_{\textrm{kst}}(\var{1},\var{2},\var{3}) &=& \dfrac{\bn_{\textrm{kst}}(\var{1},\var{2},\var{3})}{\bd_{\textrm{kst}}(\var{1},\var{2},\var{3})}\\ && \\ &=& \dfrac{\sum_{\text{rows}} \bw \odot \bPhi_{1}(\var{1}) \odot \cdots \odot\bPhi_{3}(\var{3})}{\sum_{\text{rows}} \bPhi_{1}(\var{1}) \odot \cdots \odot\bPhi_{3}(\var{3})} . \end{array}$$~\\ \noindent \textbf{KST-like univariate functions} (equivalent scaled univariate functions $\bphi_{1,\cdots,3}$): $$\left\{\begin{array}{rcrcl}z_{1} &=&\bphi_{1}(\var{1}) &=& \cfrac{\bn_{1}}{\bd_{1}} \\z_{2} &=&\bphi_{2}(\var{2}) &=& \cfrac{\bn_{2}}{\bd_{2}} \\z_{3} &=&\bphi_{3}(\var{3}) &=& 0.9891\,\var{3}+9.891 \cdot 10^{-5}\\\end{array} \right. .$$\noindent where, \\ \noindent $\bn_{1}=0.009901\,{\var{1}}^4-6.609 \cdot 10^{-14}\,{\var{1}}^3+4.588 \cdot 10^{-14}\,{\var{1}}^2+8.034 \cdot 10^{-14}\,\var{1}+0.09901$ and \\ \noindent $\bd_{1}=-3.588 \cdot 10^{-16}\,{\var{1}}^4+0.009901\,{\var{1}}^3-2.188 \cdot 10^{-13}\,{\var{1}}^2+1.061 \cdot 10^{-12}\,\var{1}+1.0$, \\ \noindent $\bn_{2}=-6.84 \cdot 10^{-16}\,{\var{2}}^2+1.033 \cdot 10^{-14}\,\var{2}+9.99$ and \\ \noindent $\bd_{2}=0.999\,{\var{2}}^2+1.077 \cdot 10^{-15}\,\var{2}+1.0$, \\

\newpage \subsection{Function \#40 (${\ord=4}$ variables, tensor size: 19.5 \textbf{MB})} $$\frac{\var{3}\var{1}}{\var{1}^2+\var{2}+\var{3}^2+1}+\var{4}^3$$ \subsubsection{Setup and results overview}\begin{itemize}\item Reference: Personal communication, [none]\item Domain: $\mathbb{R}$\item Tensor size: 19.5 \textbf{MB} ($40^{4}$ points)\item Bounds: $ \left(\begin{array}{cc} 1 & 4 \end{array}\right) \times \left(\begin{array}{cc} 1 & 4 \end{array}\right) \times \left(\begin{array}{cc} 1 & 4 \end{array}\right) \times \left(\begin{array}{cc} 1 & 4 \end{array}\right)$ \end{itemize} \begin{table}[H] \centering \begin{tabular}{llllll}
$\#$ & Alg. & Parameters & Dim. & CPU [s] & RMSE \\ 
\hline 
$\mathbf{\#40}$ & A/G/P-V 2025 (A1) & $0.0001,1$ & $4.3 \cdot 10^{02}$ & $\mathbf{0.5}$ & $\mathbf{8.9 \cdot 10^{-14}}$ \\ 
 & A/G/P-V 2025 (A2) & $1 \cdot 10^{-15},1$ & $4.3 \cdot 10^{02}$ & $1.7 \cdot 10^{02}$ & $1.6 \cdot 10^{-12}$ \\ 
 & MDSPACK v1.1.0 & $0.01,1$ & $4.3 \cdot 10^{02}$ & $0.54$ & $9.1 \cdot 10^{-14}$ \\ 
 & P/P 2025 & $1,0.95,50,0.01,4,12,9$ & $\mathbf{2.6 \cdot 10^{02}}$ & $1.8 \cdot 10^{03}$ & $0.00011$ \\ 
 & C-R/B/G 2023 & $NaN$ & $NaN$ & $NaN$ & $NaN$ \\ 
 & B/G 2025 & $1 \cdot 10^{-06},20,4$ & $8.6 \cdot 10^{02}$ & $77$ & $4.9 \cdot 10^{-10}$ \\ 
 & TensorFlow & $$ & $3.8 \cdot 10^{02}$ & $1.4 \cdot 10^{02}$ & $0.11$ \\ 
\hline 
\end{tabular} \caption{Function \#40: best model configuration and performances per methods.} \end{table}\begin{figure}[H] \centering  \includegraphics[width=\textwidth]{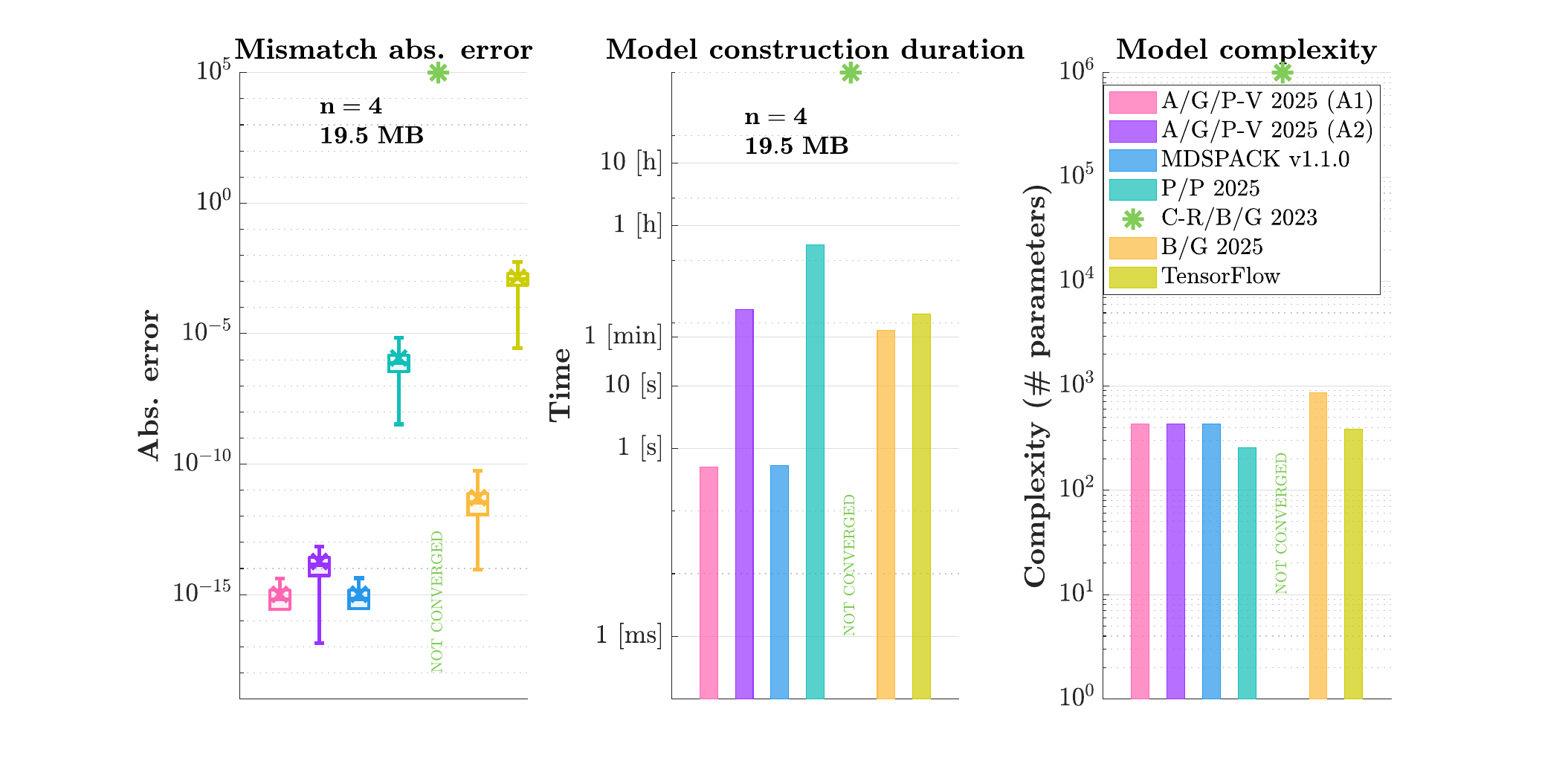} \caption{Function \#40: graphical view of the best model performances.} \end{figure}\begin{figure}[H] \centering  \includegraphics[width=\textwidth]{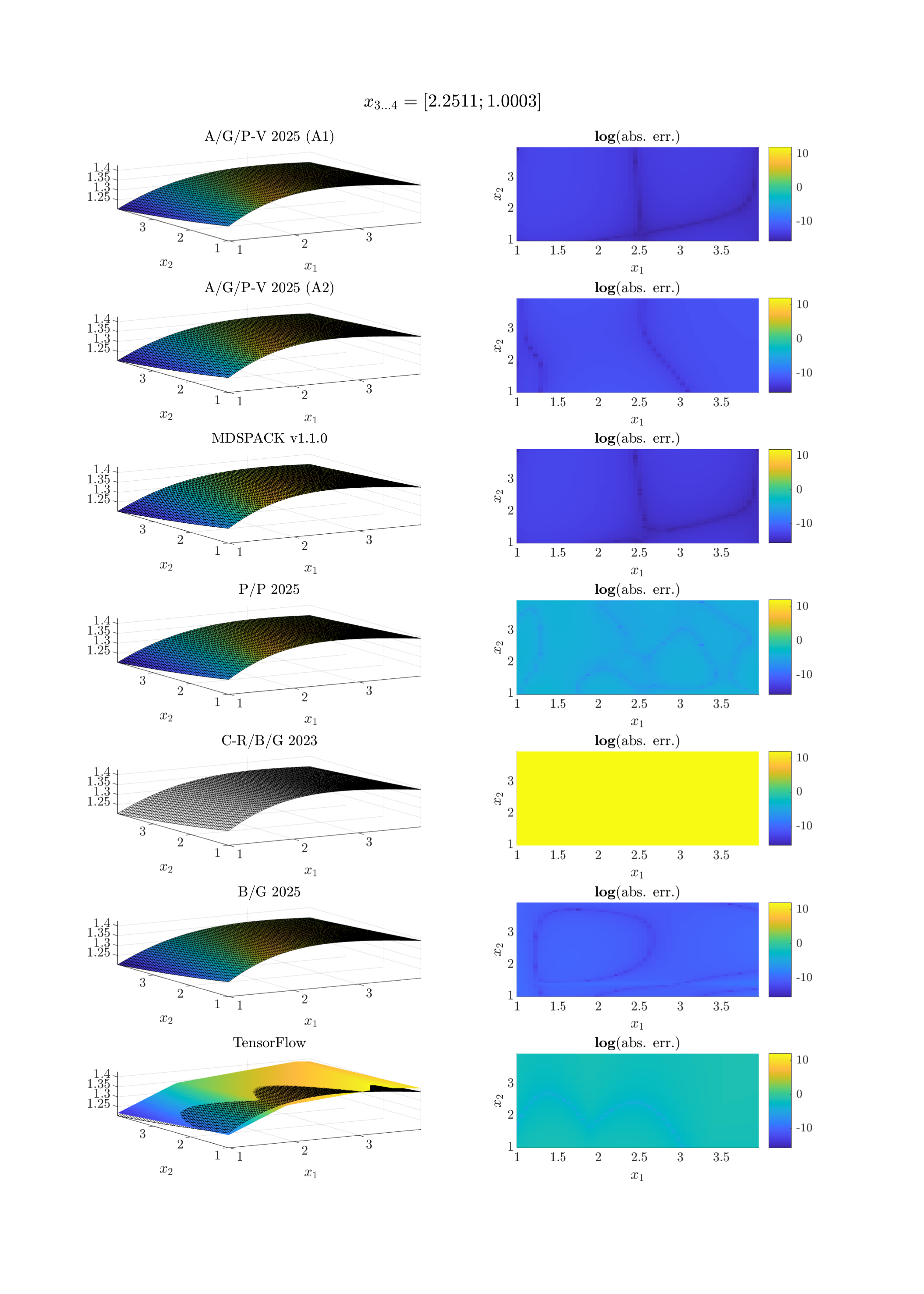} \caption{Function \#40: left side, evaluation of the original (mesh) vs. approximated (coloured surface) and right side, absolute errors (in log-scale).} \end{figure}\subsubsection{mLF detailed informations (M1)} \noindent \textbf{Right interpolation points}: $k_l=\left(\begin{array}{cccc} 3 & 2 & 3 & 4 \end{array}\right)$, where $l=1,\cdots,\ord$.$$ \begin{array}{rcl}\lan{1} &\in& \IC^{3} \text{ , linearly spaced between bounds}\\\lan{2} &\in& \IC^{2} \text{ , linearly spaced between bounds}\\\lan{3} &\in& \IC^{3} \text{ , linearly spaced between bounds}\\\lan{4} &\in& \IC^{4} \text{ , linearly spaced between bounds}\\\end{array} $$\noindent \textbf{$\ord$-D Loewner matrix, barycentric weights and Lagrangian basis}:$$ \begin{array}{rcl}\IL & \in & \IC^{72 \times 72}\\\bc & \in & \IC^{72}\\\bw & \in & \IC^{72}\\\bc\odot \bw & \in & \IC^{72}\\\mathbf{Lag}(\var{1},\var{2},\var{3},\var{4}) & \in & \IC^{72}\\\end{array} $$

\newpage \subsection{Function \#41 (${\ord=5}$ variables, tensor size: 781 \textbf{KB})} $$\frac{\var{5}^3\var{3}\var{1}+\var{3}^2}{\var{1}^3+\var{2}\var{3}+\var{4}}$$ \subsubsection{Setup and results overview}\begin{itemize}\item Reference: Personal communication, [none]\item Domain: $\mathbb{R}$\item Tensor size: 781 \textbf{KB} ($10^{5}$ points)\item Bounds: $ \left(\begin{array}{cc} \frac{1}{10} & 1 \end{array}\right) \times \left(\begin{array}{cc} \frac{1}{10} & 1 \end{array}\right) \times \left(\begin{array}{cc} \frac{1}{10} & 1 \end{array}\right) \times \left(\begin{array}{cc} \frac{1}{10} & 1 \end{array}\right) \times \left(\begin{array}{cc} \frac{1}{10} & 1 \end{array}\right)$ \end{itemize} \begin{table}[H] \centering \begin{tabular}{llllll}
$\#$ & Alg. & Parameters & Dim. & CPU [s] & RMSE \\ 
\hline 
$\mathbf{\#41}$ & A/G/P-V 2025 (A1) & $0.001,2$ & $1.3 \cdot 10^{03}$ & $0.051$ & $5.3 \cdot 10^{-14}$ \\ 
 & A/G/P-V 2025 (A2) & $1 \cdot 10^{-15},1$ & $1.3 \cdot 10^{03}$ & $7.9$ & $\mathbf{4.5 \cdot 10^{-14}}$ \\ 
 & MDSPACK v1.1.0 & $0.0001,2$ & $1.3 \cdot 10^{03}$ & $\mathbf{0.037}$ & $5.4 \cdot 10^{-14}$ \\ 
 & P/P 2025 & $1,0.95,50,0.01,4,12,9$ & $\mathbf{2.9 \cdot 10^{02}}$ & $17$ & $0.002$ \\ 
 & C-R/B/G 2023 & $0.001,20$ & $1.4 \cdot 10^{04}$ & $54$ & $5.5 \cdot 10^{-13}$ \\ 
 & B/G 2025 & $0.001,20,4$ & $9.7 \cdot 10^{04}$ & $65$ & $1.4 \cdot 10^{-05}$ \\ 
 & TensorFlow & $NaN$ & $NaN$ & $NaN$ & $NaN$ \\ 
\hline 
\end{tabular} \caption{Function \#41: best model configuration and performances per methods.} \end{table}\begin{figure}[H] \centering  \includegraphics[width=\textwidth]{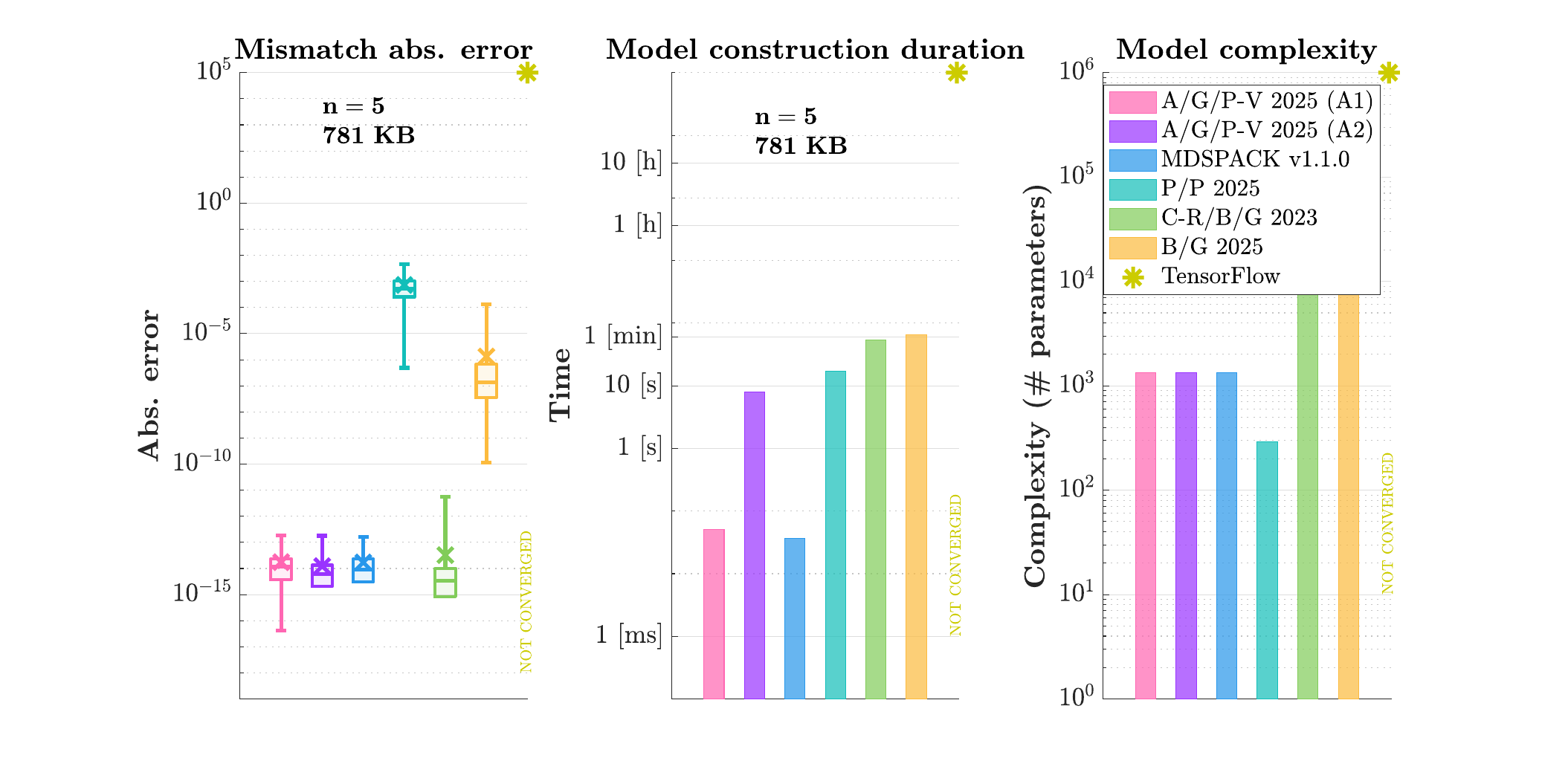} \caption{Function \#41: graphical view of the best model performances.} \end{figure}\begin{figure}[H] \centering  \includegraphics[width=\textwidth]{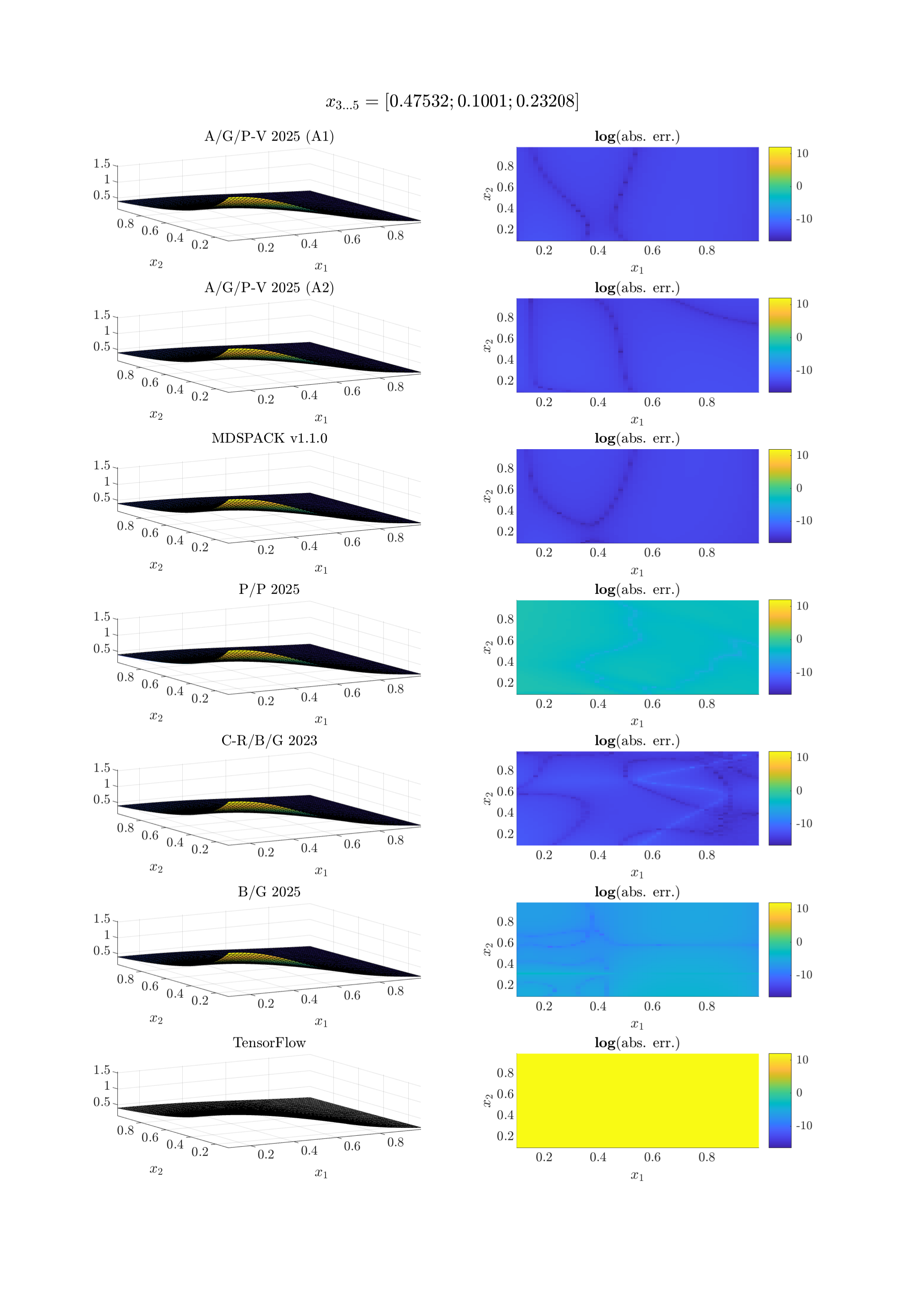} \caption{Function \#41: left side, evaluation of the original (mesh) vs. approximated (coloured surface) and right side, absolute errors (in log-scale).} \end{figure}\subsubsection{mLF detailed informations (M1)} \noindent \textbf{Right interpolation points}: $k_l=\left(\begin{array}{ccccc} 4 & 2 & 3 & 2 & 4 \end{array}\right)$, where $l=1,\cdots,\ord$.$$ \begin{array}{rcl}\lan{1} &\in& \IC^{4} \text{ , linearly spaced between bounds}\\\lan{2} &\in& \IC^{2} \text{ , linearly spaced between bounds}\\\lan{3} &\in& \IC^{3} \text{ , linearly spaced between bounds}\\\lan{4} &\in& \IC^{2} \text{ , linearly spaced between bounds}\\\lan{5} &\in& \IC^{4} \text{ , linearly spaced between bounds}\\\end{array} $$\noindent \textbf{$\ord$-D Loewner matrix, barycentric weights and Lagrangian basis}:$$ \begin{array}{rcl}\IL & \in & \IC^{192 \times 192}\\\bc & \in & \IC^{192}\\\bw & \in & \IC^{192}\\\bc\odot \bw & \in & \IC^{192}\\\mathbf{Lag}(\var{1},\var{2},\var{3},\var{4},\var{5}) & \in & \IC^{192}\\\end{array} $$

\newpage \subsection{Function \#42 (${\ord=6}$ variables, tensor size: 7.63 \textbf{MB})} $$\frac{\var{1}+\var{3}-\sqrt{2}\var{6}^2}{\var{1}^4+\var{2}\var{3}+\var{4}^3+\var{5}^2+\var{6}}$$ \subsubsection{Setup and results overview}\begin{itemize}\item Reference: Personal communication, [none]\item Domain: $\mathbb{R}$\item Tensor size: 7.63 \textbf{MB} ($10^{6}$ points)\item Bounds: $ \left(\begin{array}{cc} \frac{1}{10} & 1 \end{array}\right) \times \left(\begin{array}{cc} \frac{1}{10} & 1 \end{array}\right) \times \left(\begin{array}{cc} \frac{1}{10} & 1 \end{array}\right) \times \left(\begin{array}{cc} \frac{1}{10} & 1 \end{array}\right) \times \left(\begin{array}{cc} \frac{1}{10} & 1 \end{array}\right) \times \left(\begin{array}{cc} \frac{1}{10} & 1 \end{array}\right)$ \end{itemize} \begin{table}[H] \centering \begin{tabular}{llllll}
$\#$ & Alg. & Parameters & Dim. & CPU [s] & RMSE \\ 
\hline 
$\mathbf{\#42}$ & A/G/P-V 2025 (A1) & $0.0001,2$ & $5.8 \cdot 10^{03}$ & $\mathbf{0.38}$ & $\mathbf{4 \cdot 10^{-14}}$ \\ 
 & A/G/P-V 2025 (A2) & $1 \cdot 10^{-15},1$ & $5.8 \cdot 10^{03}$ & $1.4 \cdot 10^{02}$ & $4.1 \cdot 10^{-13}$ \\ 
 & MDSPACK v1.1.0 & $1 \cdot 10^{-06},3$ & $5.8 \cdot 10^{03}$ & $0.38$ & $4 \cdot 10^{-14}$ \\ 
 & P/P 2025 & $1,0.95,50,0.01,4,12,9$ & $\mathbf{3.3 \cdot 10^{02}}$ & $3.9 \cdot 10^{02}$ & $0.0032$ \\ 
 & C-R/B/G 2023 & $NaN$ & $NaN$ & $NaN$ & $NaN$ \\ 
 & B/G 2025 & $1 \cdot 10^{-06},20,4$ & $1.3 \cdot 10^{06}$ & $1.2 \cdot 10^{03}$ & $0.00032$ \\ 
 & TensorFlow & $NaN$ & $NaN$ & $NaN$ & $NaN$ \\ 
\hline 
\end{tabular} \caption{Function \#42: best model configuration and performances per methods.} \end{table}\begin{figure}[H] \centering  \includegraphics[width=\textwidth]{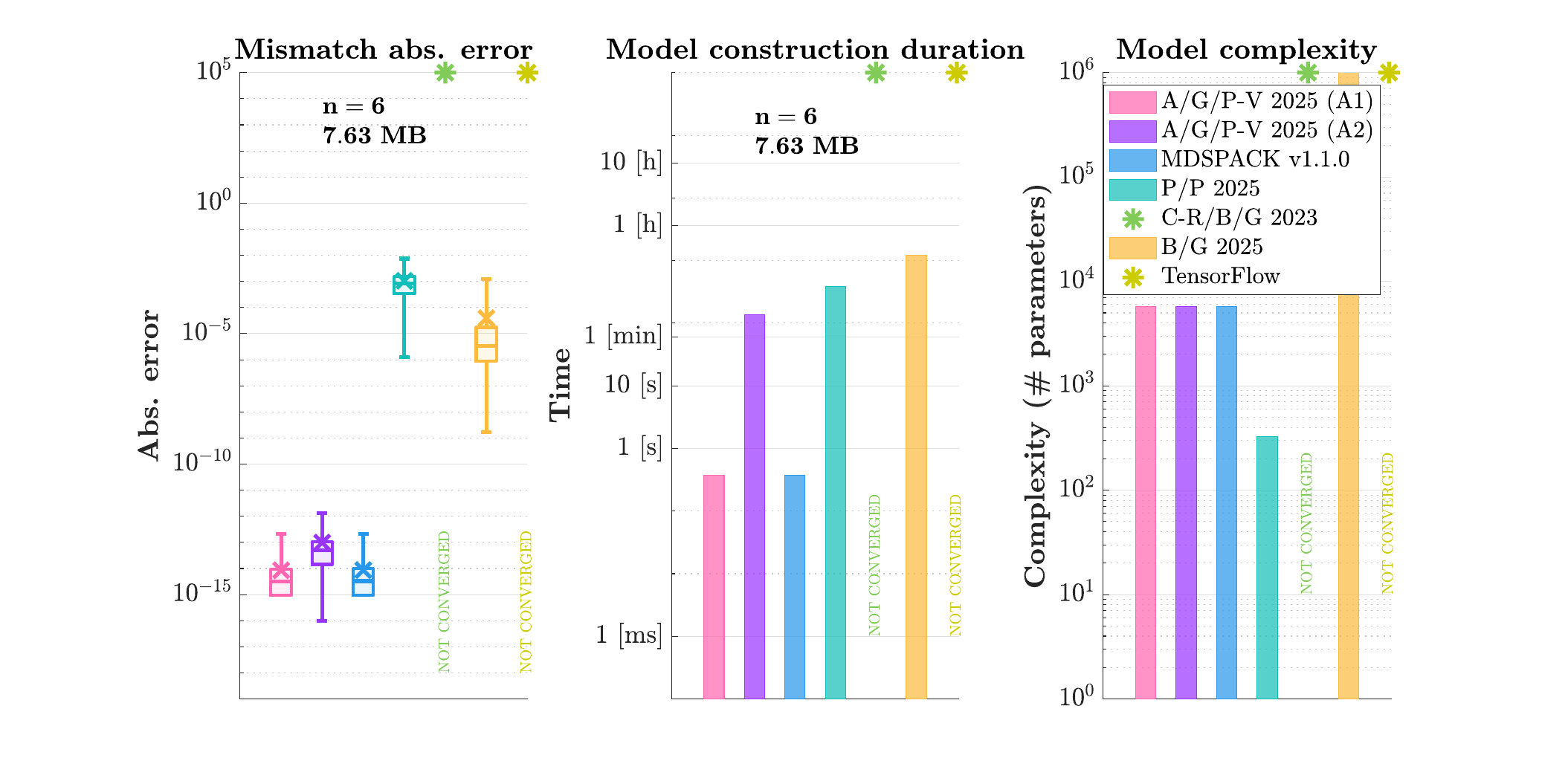} \caption{Function \#42: graphical view of the best model performances.} \end{figure}\begin{figure}[H] \centering  \includegraphics[width=\textwidth]{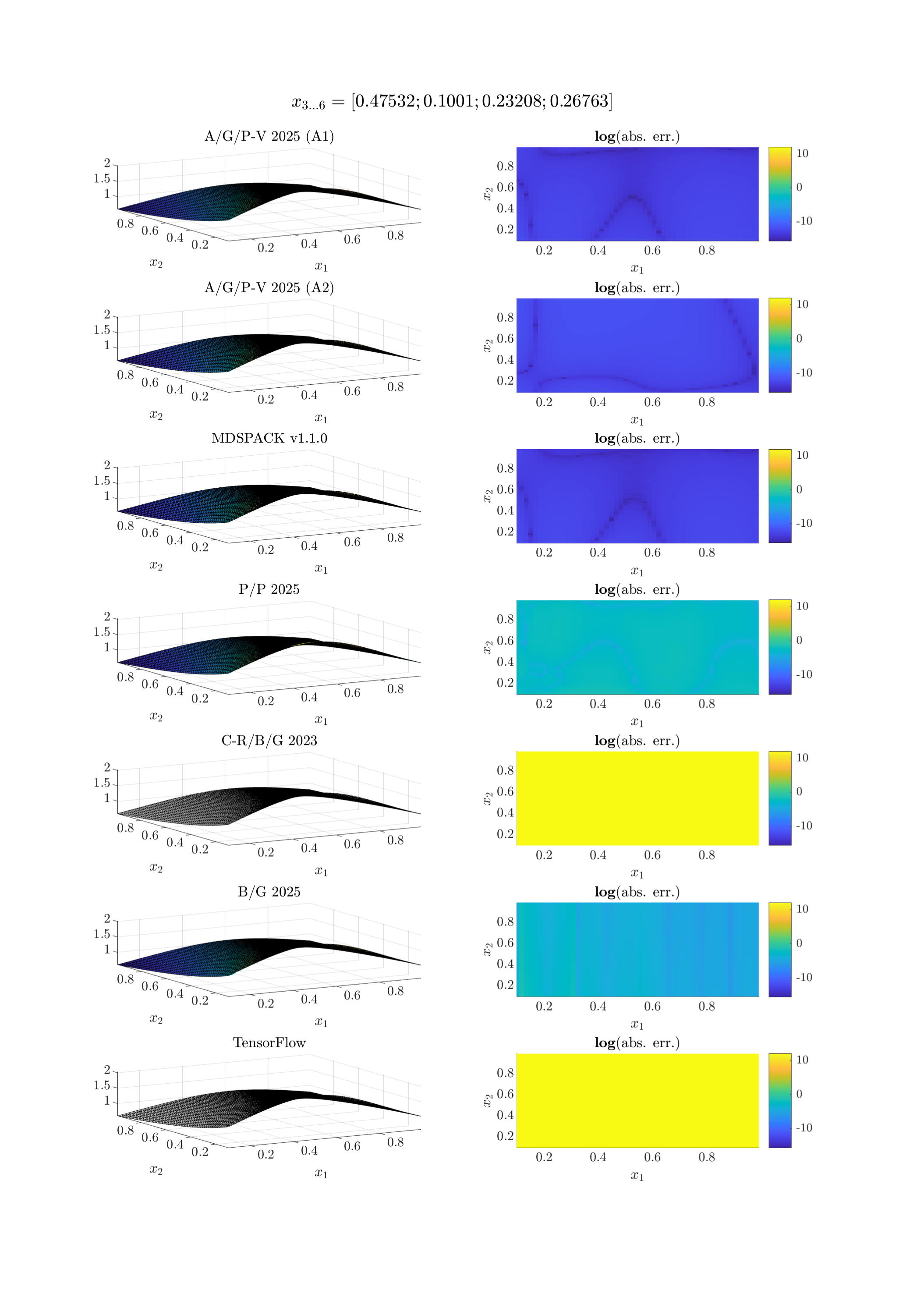} \caption{Function \#42: left side, evaluation of the original (mesh) vs. approximated (coloured surface) and right side, absolute errors (in log-scale).} \end{figure}\subsubsection{mLF detailed informations (M1)} \noindent \textbf{Right interpolation points}: $k_l=\left(\begin{array}{cccccc} 5 & 2 & 2 & 4 & 3 & 3 \end{array}\right)$, where $l=1,\cdots,\ord$.$$ \begin{array}{rcl}\lan{1} &\in& \IC^{5} \text{ , linearly spaced between bounds}\\\lan{2} &\in& \IC^{2} \text{ , linearly spaced between bounds}\\\lan{3} &\in& \IC^{2} \text{ , linearly spaced between bounds}\\\lan{4} &\in& \IC^{4} \text{ , linearly spaced between bounds}\\\lan{5} &\in& \IC^{3} \text{ , linearly spaced between bounds}\\\lan{6} &\in& \IC^{3} \text{ , linearly spaced between bounds}\\\end{array} $$\noindent \textbf{$\ord$-D Loewner matrix, barycentric weights and Lagrangian basis}:$$ \begin{array}{rcl}\IL & \in & \IC^{720 \times 720}\\\bc & \in & \IC^{720}\\\bw & \in & \IC^{720}\\\bc\odot \bw & \in & \IC^{720}\\\mathbf{Lag}(\var{1},\var{2},\var{3},\var{4},\var{5},\var{6}) & \in & \IC^{720}\\\end{array} $$

\newpage \subsection{Function \#43 (${\ord=7}$ variables, tensor size: 76.3 \textbf{MB})} $$\frac{\var{3}\var{2}^3+1}{\var{1}^4+\var{2}^2\var{3}+\var{4}^2+\var{5}+\var{6}^3+\var{7}}$$ \subsubsection{Setup and results overview}\begin{itemize}\item Reference: Personal communication, [none]\item Domain: $\mathbb{R}$\item Tensor size: 76.3 \textbf{MB} ($10^{7}$ points)\item Bounds: $ \left(\begin{array}{cc} 1 & 10 \end{array}\right) \times \left(\begin{array}{cc} 1 & 10 \end{array}\right) \times \left(\begin{array}{cc} 1 & 10 \end{array}\right) \times \left(\begin{array}{cc} 1 & 10 \end{array}\right) \times \left(\begin{array}{cc} 1 & 10 \end{array}\right) \times \left(\begin{array}{cc} 1 & 10 \end{array}\right) \times \left(\begin{array}{cc} 1 & 10 \end{array}\right)$ \end{itemize} \begin{table}[H] \centering \begin{tabular}{llllll}
$\#$ & Alg. & Parameters & Dim. & CPU [s] & RMSE \\ 
\hline 
$\mathbf{\#43}$ & A/G/P-V 2025 (A1) & $0.0001,1$ & $\mathbf{1.7 \cdot 10^{04}}$ & $\mathbf{4.5}$ & $\mathbf{1 \cdot 10^{-12}}$ \\ 
 & A/G/P-V 2025 (A2) & $1 \cdot 10^{-15},1$ & $1.7 \cdot 10^{04}$ & $2.3 \cdot 10^{03}$ & $4.9 \cdot 10^{-11}$ \\ 
 & MDSPACK v1.1.0 & $1 \cdot 10^{-08},4$ & $1.7 \cdot 10^{04}$ & $5.4$ & $1 \cdot 10^{-12}$ \\ 
 & P/P 2025 & $NaN$ & $NaN$ & $NaN$ & $NaN$ \\ 
 & C-R/B/G 2023 & $NaN$ & $NaN$ & $NaN$ & $NaN$ \\ 
 & B/G 2025 & $NaN$ & $NaN$ & $NaN$ & $NaN$ \\ 
 & TensorFlow & $NaN$ & $NaN$ & $NaN$ & $NaN$ \\ 
\hline 
\end{tabular} \caption{Function \#43: best model configuration and performances per methods.} \end{table}\begin{figure}[H] \centering  \includegraphics[width=\textwidth]{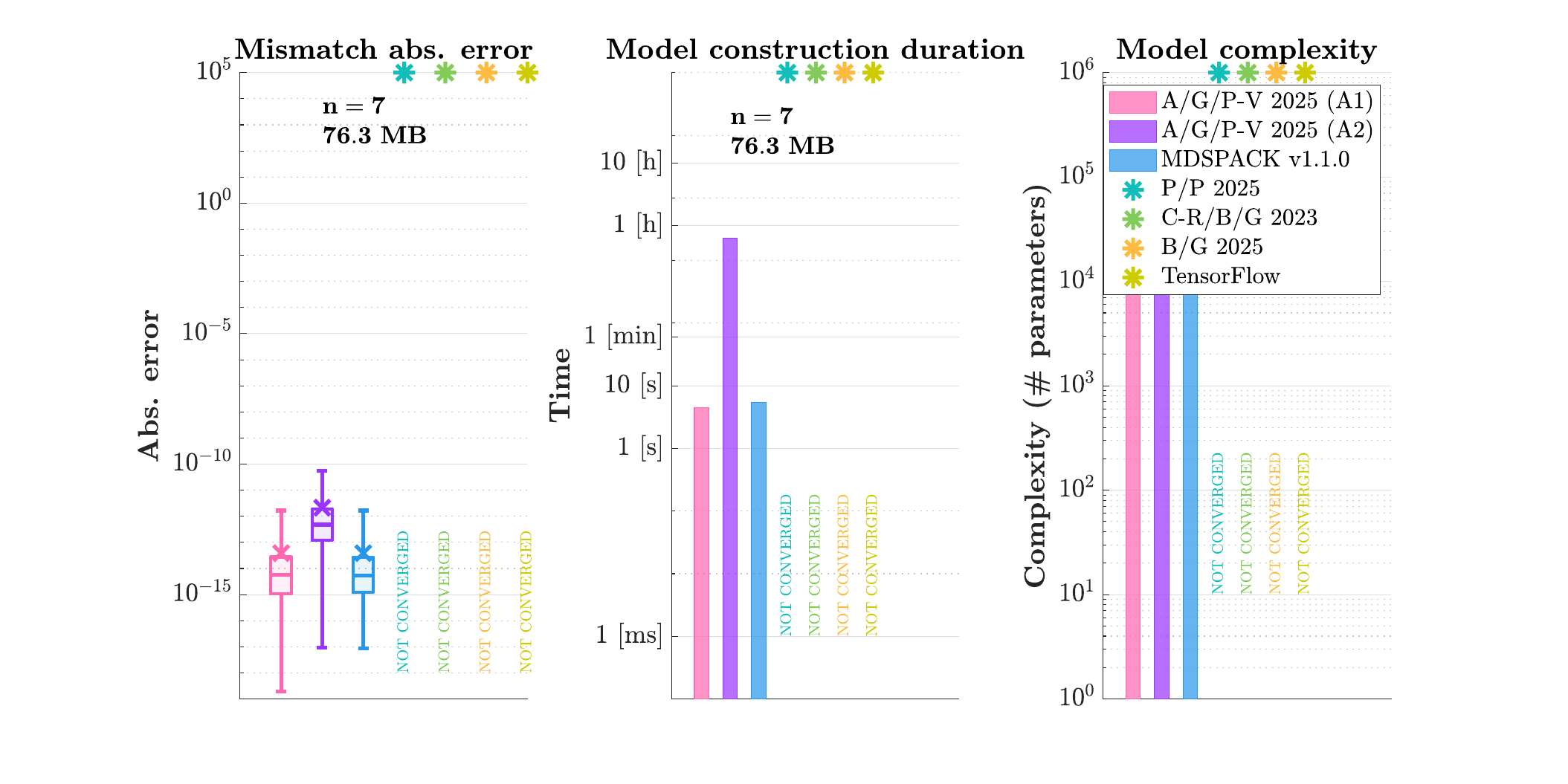} \caption{Function \#43: graphical view of the best model performances.} \end{figure}\begin{figure}[H] \centering  \includegraphics[width=\textwidth]{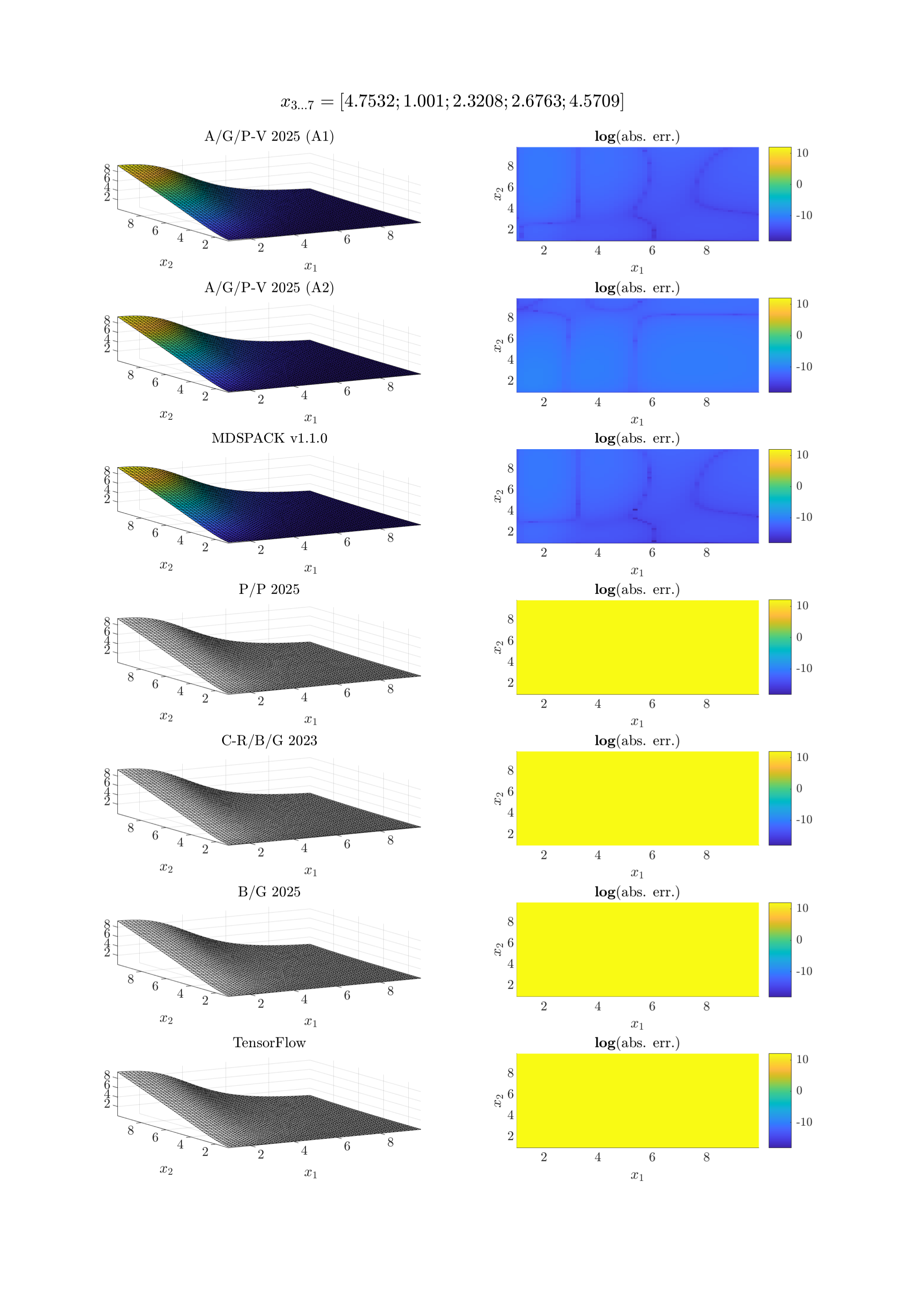} \caption{Function \#43: left side, evaluation of the original (mesh) vs. approximated (coloured surface) and right side, absolute errors (in log-scale).} \end{figure}\subsubsection{mLF detailed informations (M1)} \noindent \textbf{Right interpolation points}: $k_l=\left(\begin{array}{ccccccc} 5 & 4 & 2 & 3 & 2 & 4 & 2 \end{array}\right)$, where $l=1,\cdots,\ord$.$$ \begin{array}{rcl}\lan{1} &\in& \IC^{5} \text{ , linearly spaced between bounds}\\\lan{2} &\in& \IC^{4} \text{ , linearly spaced between bounds}\\\lan{3} &\in& \IC^{2} \text{ , linearly spaced between bounds}\\\lan{4} &\in& \IC^{3} \text{ , linearly spaced between bounds}\\\lan{5} &\in& \IC^{2} \text{ , linearly spaced between bounds}\\\lan{6} &\in& \IC^{4} \text{ , linearly spaced between bounds}\\\lan{7} &\in& \IC^{2} \text{ , linearly spaced between bounds}\\\end{array} $$\noindent \textbf{$\ord$-D Loewner matrix, barycentric weights and Lagrangian basis}:$$ \begin{array}{rcl}\IL & \in & \IC^{1920 \times 1920}\\\bc & \in & \IC^{1920}\\\bw & \in & \IC^{1920}\\\bc\odot \bw & \in & \IC^{1920}\\\mathbf{Lag}(\var{1},\var{2},\var{3},\var{4},\var{5},\var{6},\var{7}) & \in & \IC^{1920}\\\end{array} $$

\newpage \subsection{Function \#44 (${\ord=8}$ variables, tensor size: 763 \textbf{MB})} $$\frac{1}{\var{1}^4+\var{2}^2\var{3}+\var{4}^2+\var{5}+\var{6}+\var{7}+\var{8}}$$ \subsubsection{Setup and results overview}\begin{itemize}\item Reference: Personal communication, [none]\item Domain: $\mathbb{R}$\item Tensor size: 763 \textbf{MB} ($10^{8}$ points)\item Bounds: $ \left(\begin{array}{cc} \frac{1}{10} & 20 \end{array}\right) \times \left(\begin{array}{cc} \frac{1}{10} & 20 \end{array}\right) \times \left(\begin{array}{cc} \frac{1}{10} & 20 \end{array}\right) \times \left(\begin{array}{cc} \frac{1}{10} & 20 \end{array}\right) \times \left(\begin{array}{cc} \frac{1}{10} & 20 \end{array}\right) \times \left(\begin{array}{cc} \frac{1}{10} & 20 \end{array}\right) \times \left(\begin{array}{cc} \frac{1}{10} & 20 \end{array}\right) \times \left(\begin{array}{cc} \frac{1}{10} & 20 \end{array}\right)$ \end{itemize} \begin{table}[H] \centering \begin{tabular}{llllll}
$\#$ & Alg. & Parameters & Dim. & CPU [s] & RMSE \\ 
\hline 
$\mathbf{\#44}$ & A/G/P-V 2025 (A1) & $0.0001,2$ & $\mathbf{1.4 \cdot 10^{04}}$ & $\mathbf{1.7 \cdot 10^{02}}$ & $\mathbf{2.3 \cdot 10^{-12}}$ \\ 
 & A/G/P-V 2025 (A2) & $1 \cdot 10^{-15},1$ & $NaN$ & $NaN$ & $NaN$ \\ 
 & MDSPACK v1.1.0 & $1 \cdot 10^{-14},7$ & $1.4 \cdot 10^{04}$ & $1.9 \cdot 10^{02}$ & $2.3 \cdot 10^{-12}$ \\ 
 & P/P 2025 & $NaN$ & $NaN$ & $NaN$ & $NaN$ \\ 
 & C-R/B/G 2023 & $NaN$ & $NaN$ & $NaN$ & $NaN$ \\ 
 & B/G 2025 & $NaN$ & $NaN$ & $NaN$ & $NaN$ \\ 
 & TensorFlow & $NaN$ & $NaN$ & $NaN$ & $NaN$ \\ 
\hline 
\end{tabular} \caption{Function \#44: best model configuration and performances per methods.} \end{table}\begin{figure}[H] \centering  \includegraphics[width=\textwidth]{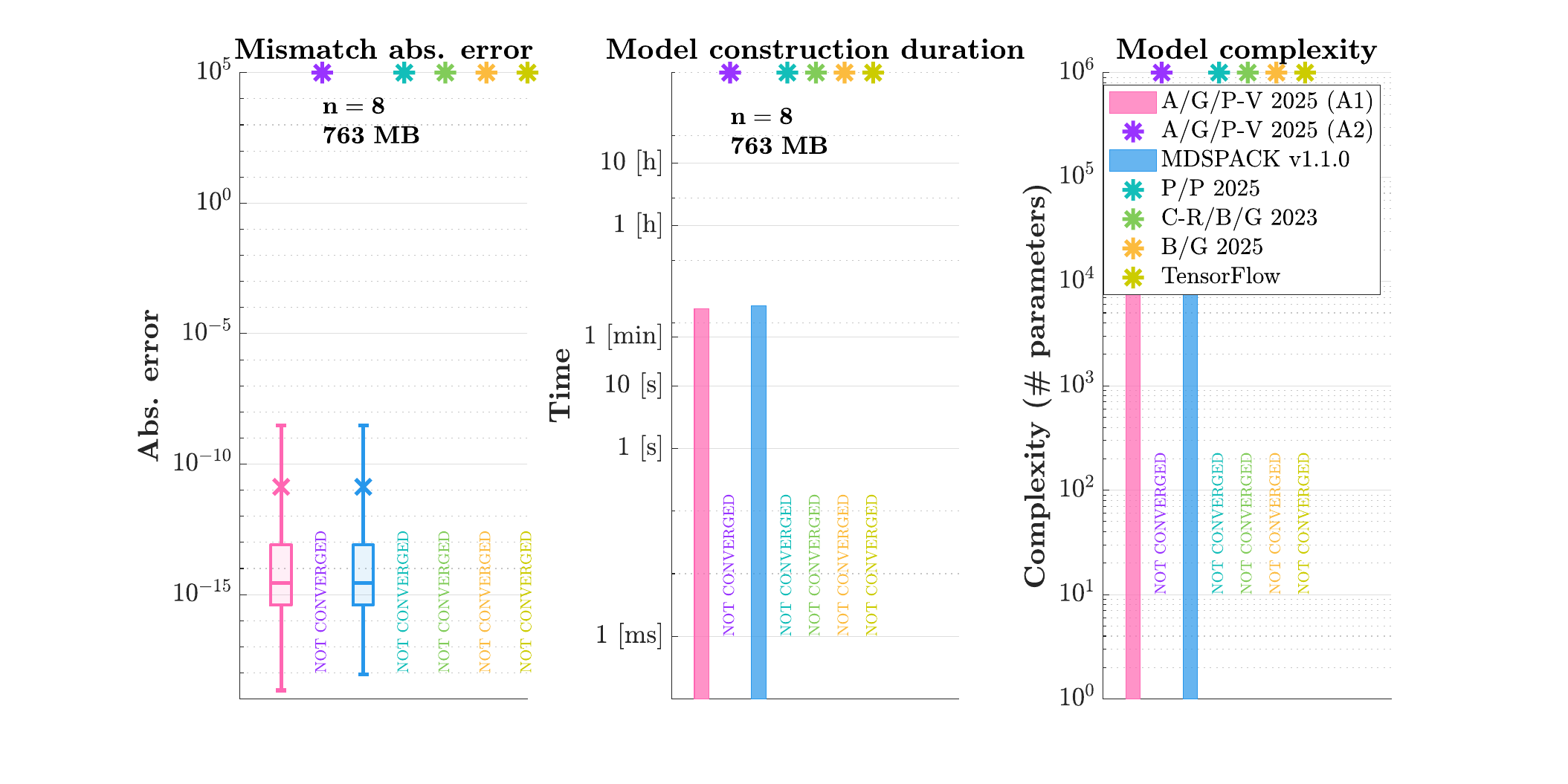} \caption{Function \#44: graphical view of the best model performances.} \end{figure}\begin{figure}[H] \centering  \includegraphics[width=\textwidth]{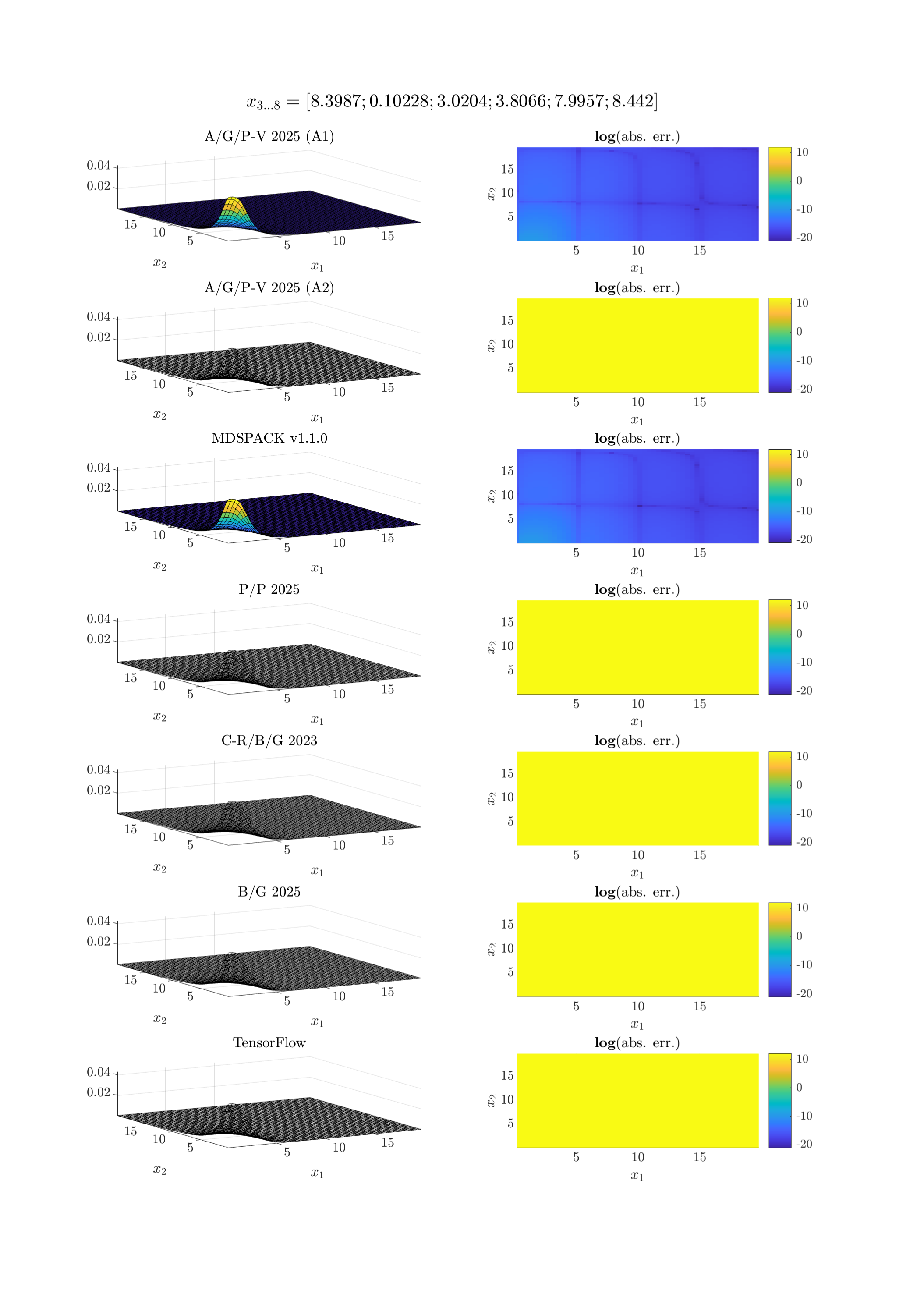} \caption{Function \#44: left side, evaluation of the original (mesh) vs. approximated (coloured surface) and right side, absolute errors (in log-scale).} \end{figure}\subsubsection{mLF detailed informations (M1)} \noindent \textbf{Right interpolation points}: $k_l=\left(\begin{array}{cccccccc} 5 & 3 & 2 & 3 & 2 & 2 & 2 & 2 \end{array}\right)$, where $l=1,\cdots,\ord$.$$ \begin{array}{rcl}\lan{1} &\in& \IC^{5} \text{ , linearly spaced between bounds}\\\lan{2} &\in& \IC^{3} \text{ , linearly spaced between bounds}\\\lan{3} &\in& \IC^{2} \text{ , linearly spaced between bounds}\\\lan{4} &\in& \IC^{3} \text{ , linearly spaced between bounds}\\\lan{5} &\in& \IC^{2} \text{ , linearly spaced between bounds}\\\lan{6} &\in& \IC^{2} \text{ , linearly spaced between bounds}\\\lan{7} &\in& \IC^{2} \text{ , linearly spaced between bounds}\\\lan{8} &\in& \IC^{2} \text{ , linearly spaced between bounds}\\\end{array} $$\noindent \textbf{$\ord$-D Loewner matrix, barycentric weights and Lagrangian basis}:$$ \begin{array}{rcl}\IL & \in & \IC^{1440 \times 1440}\\\bc & \in & \IC^{1440}\\\bw & \in & \IC^{1440}\\\bc\odot \bw & \in & \IC^{1440}\\\mathbf{Lag}(\var{1},\var{2},\var{3},\var{4},\var{5},\var{6},\var{7},\var{8}) & \in & \IC^{1440}\\\end{array} $$

\newpage \subsection{Function \#45 (${\ord=9}$ variables, tensor size: 76.9 \textbf{MB})} $$\frac{1}{\var{1}^2+\var{2}^2\var{3}+\var{4}^2+\var{5}+\var{6}+\var{7}+\var{8}+\var{9}}$$ \subsubsection{Setup and results overview}\begin{itemize}\item Reference: Personal communication, [none]\item Domain: $\mathbb{R}$\item Tensor size: 76.9 \textbf{MB} ($6^{9}$ points)\item Bounds: $ \left(\begin{array}{cc} 1 & 5 \end{array}\right) \times \left(\begin{array}{cc} 1 & 5 \end{array}\right) \times \left(\begin{array}{cc} 1 & 5 \end{array}\right) \times \left(\begin{array}{cc} 1 & 5 \end{array}\right) \times \left(\begin{array}{cc} 1 & 5 \end{array}\right) \times \left(\begin{array}{cc} 1 & 5 \end{array}\right) \times \left(\begin{array}{cc} 1 & 5 \end{array}\right) \times \left(\begin{array}{cc} 1 & 5 \end{array}\right) \times \left(\begin{array}{cc} 1 & 5 \end{array}\right)$ \end{itemize} \begin{table}[H] \centering \begin{tabular}{llllll}
$\#$ & Alg. & Parameters & Dim. & CPU [s] & RMSE \\ 
\hline 
$\mathbf{\#45}$ & A/G/P-V 2025 (A1) & $0.01,1$ & $\mathbf{1.9 \cdot 10^{04}}$ & $\mathbf{8.1}$ & $\mathbf{3.4 \cdot 10^{-17}}$ \\ 
 & A/G/P-V 2025 (A2) & $1 \cdot 10^{-15},1$ & $NaN$ & $NaN$ & $NaN$ \\ 
 & MDSPACK v1.1.0 & $1 \cdot 10^{-06},3$ & $1.9 \cdot 10^{04}$ & $9.6$ & $4.9 \cdot 10^{-17}$ \\ 
 & P/P 2025 & $NaN$ & $NaN$ & $NaN$ & $NaN$ \\ 
 & C-R/B/G 2023 & $NaN$ & $NaN$ & $NaN$ & $NaN$ \\ 
 & B/G 2025 & $NaN$ & $NaN$ & $NaN$ & $NaN$ \\ 
 & TensorFlow & $NaN$ & $NaN$ & $NaN$ & $NaN$ \\ 
\hline 
\end{tabular} \caption{Function \#45: best model configuration and performances per methods.} \end{table}\begin{figure}[H] \centering  \includegraphics[width=\textwidth]{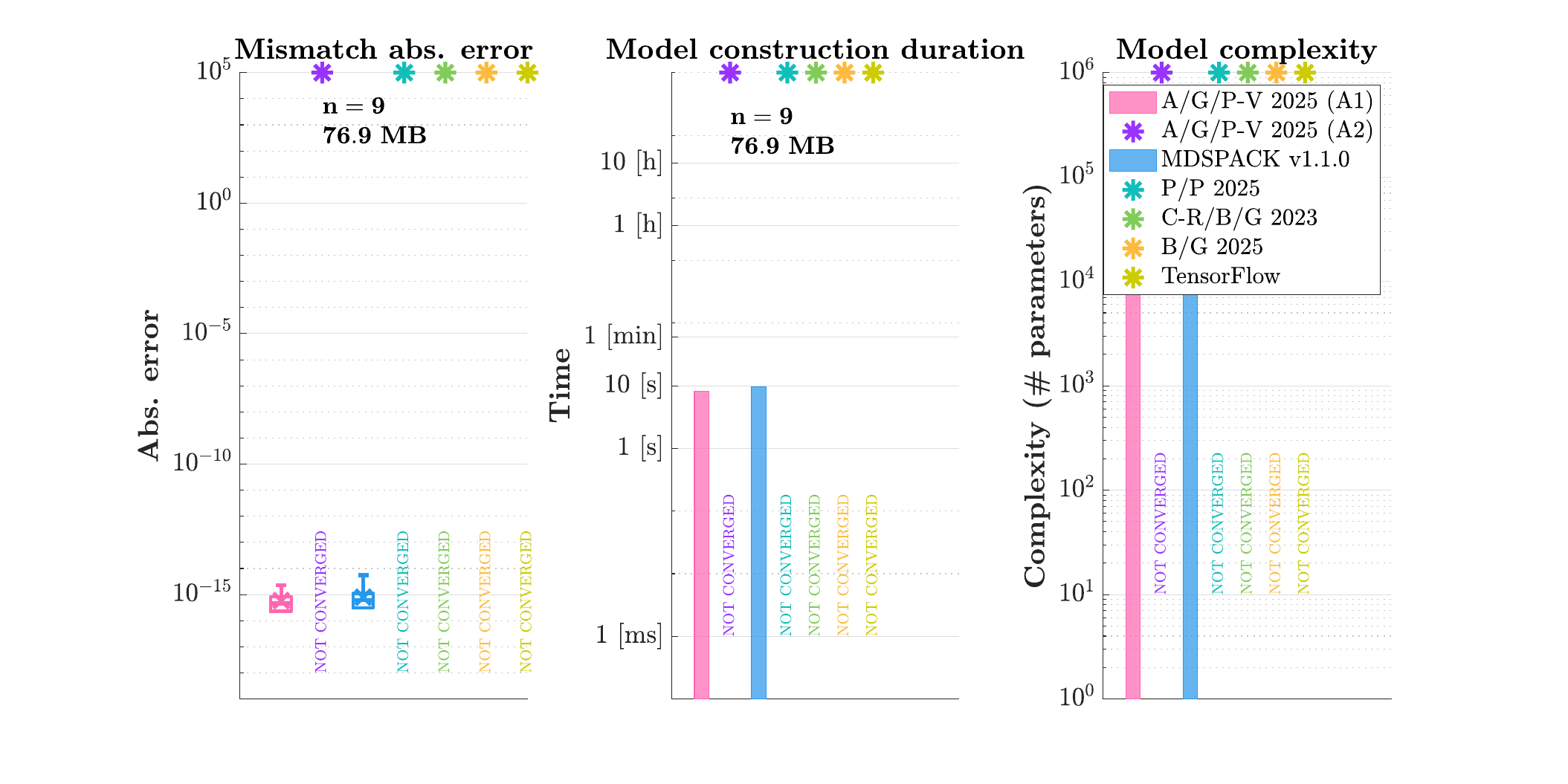} \caption{Function \#45: graphical view of the best model performances.} \end{figure}\begin{figure}[H] \centering  \includegraphics[width=\textwidth]{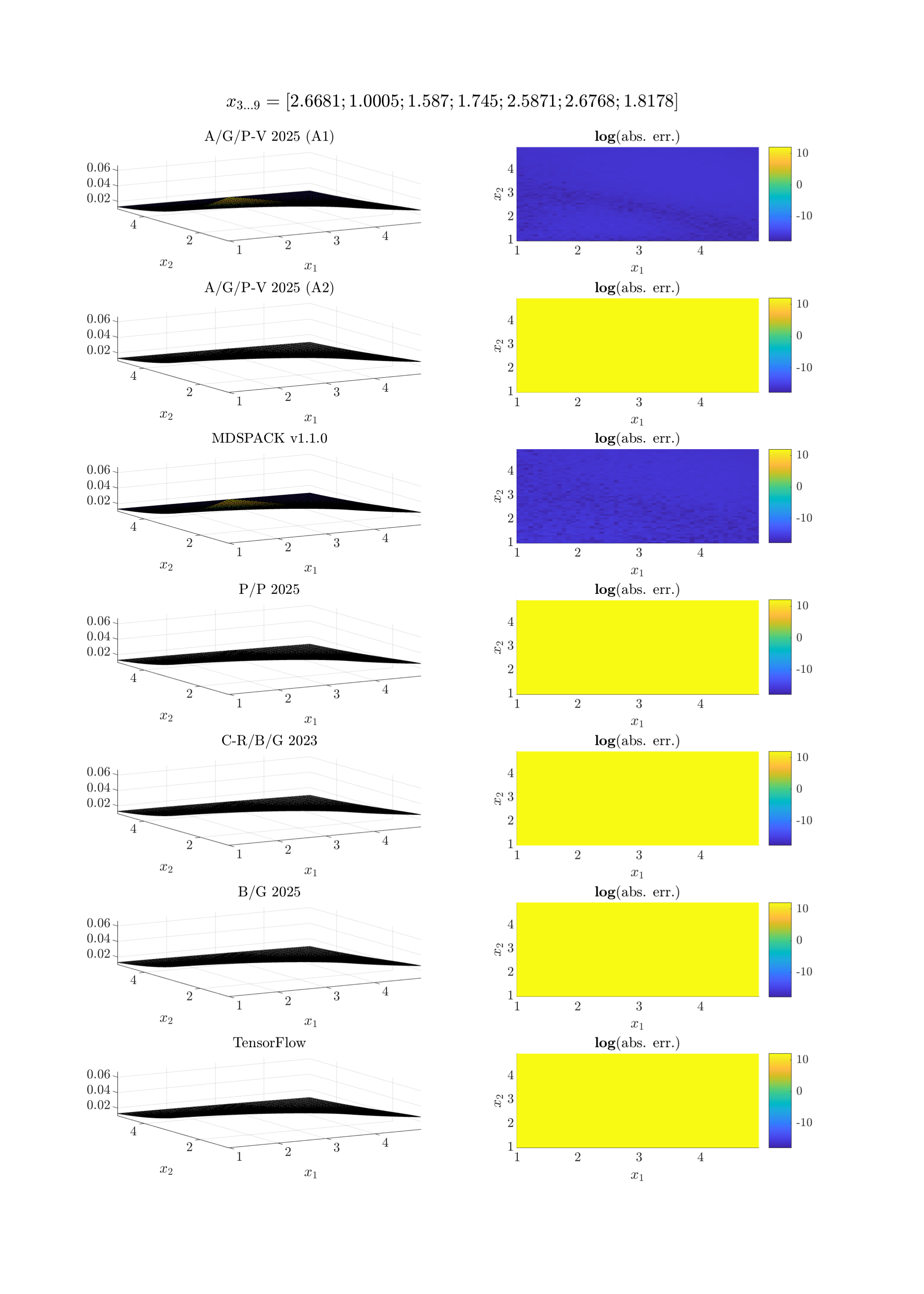} \caption{Function \#45: left side, evaluation of the original (mesh) vs. approximated (coloured surface) and right side, absolute errors (in log-scale).} \end{figure}\subsubsection{mLF detailed informations (M1)} \noindent \textbf{Right interpolation points}: $k_l=\left(\begin{array}{ccccccccc} 3 & 3 & 2 & 3 & 2 & 2 & 2 & 2 & 2 \end{array}\right)$, where $l=1,\cdots,\ord$.$$ \begin{array}{rcl}\lan{1} &\in& \IC^{3} \text{ , linearly spaced between bounds}\\\lan{2} &\in& \IC^{3} \text{ , linearly spaced between bounds}\\\lan{3} &\in& \IC^{2} \text{ , linearly spaced between bounds}\\\lan{4} &\in& \IC^{3} \text{ , linearly spaced between bounds}\\\lan{5} &\in& \IC^{2} \text{ , linearly spaced between bounds}\\\lan{6} &\in& \IC^{2} \text{ , linearly spaced between bounds}\\\lan{7} &\in& \IC^{2} \text{ , linearly spaced between bounds}\\\lan{8} &\in& \IC^{2} \text{ , linearly spaced between bounds}\\\lan{9} &\in& \IC^{2} \text{ , linearly spaced between bounds}\\\end{array} $$\noindent \textbf{$\ord$-D Loewner matrix, barycentric weights and Lagrangian basis}:$$ \begin{array}{rcl}\IL & \in & \IC^{1728 \times 1728}\\\bc & \in & \IC^{1728}\\\bw & \in & \IC^{1728}\\\bc\odot \bw & \in & \IC^{1728}\\\mathbf{Lag}(\var{1},\var{2},\var{3},\var{4},\var{5},\var{6},\var{7},\var{8},\var{9}) & \in & \IC^{1728}\\\end{array} $$

\newpage \subsection{Function \#46 (${\ord=10}$ variables, tensor size: 461 \textbf{MB})} $$\frac{1}{\var{1}+\var{1}^2\var{2}\var{3}+\var{4}+\var{5}+\var{6}+\var{7}\var{8}+\var{9}^2+\var{10}}$$ \subsubsection{Setup and results overview}\begin{itemize}\item Reference: Personal communication, [none]\item Domain: $\mathbb{R}$\item Tensor size: 461 \textbf{MB} ($6^{10}$ points)\item Bounds: $ \left(\begin{array}{cc} 1 & 5 \end{array}\right) \times \left(\begin{array}{cc} 1 & 5 \end{array}\right) \times \left(\begin{array}{cc} 1 & 5 \end{array}\right) \times \left(\begin{array}{cc} 1 & 5 \end{array}\right) \times \left(\begin{array}{cc} 1 & 5 \end{array}\right) \times \left(\begin{array}{cc} 1 & 5 \end{array}\right) \times \left(\begin{array}{cc} 1 & 5 \end{array}\right) \times \left(\begin{array}{cc} 1 & 5 \end{array}\right) \times \left(\begin{array}{cc} 1 & 5 \end{array}\right) \times \left(\begin{array}{cc} 1 & 5 \end{array}\right)$ \end{itemize} \begin{table}[H] \centering \begin{tabular}{llllll}
$\#$ & Alg. & Parameters & Dim. & CPU [s] & RMSE \\ 
\hline 
$\mathbf{\#46}$ & A/G/P-V 2025 (A1) & $0.01,3$ & $\mathbf{2.8 \cdot 10^{04}}$ & $1.6 \cdot 10^{02}$ & $\mathbf{5.7 \cdot 10^{-17}}$ \\ 
 & A/G/P-V 2025 (A2) & $1 \cdot 10^{-15},1$ & $NaN$ & $NaN$ & $NaN$ \\ 
 & MDSPACK v1.1.0 & $1 \cdot 10^{-06},3$ & $2.8 \cdot 10^{04}$ & $\mathbf{1.4 \cdot 10^{02}}$ & $7.2 \cdot 10^{-17}$ \\ 
 & P/P 2025 & $NaN$ & $NaN$ & $NaN$ & $NaN$ \\ 
 & C-R/B/G 2023 & $NaN$ & $NaN$ & $NaN$ & $NaN$ \\ 
 & B/G 2025 & $NaN$ & $NaN$ & $NaN$ & $NaN$ \\ 
 & TensorFlow & $NaN$ & $NaN$ & $NaN$ & $NaN$ \\ 
\hline 
\end{tabular} \caption{Function \#46: best model configuration and performances per methods.} \end{table}\begin{figure}[H] \centering  \includegraphics[width=\textwidth]{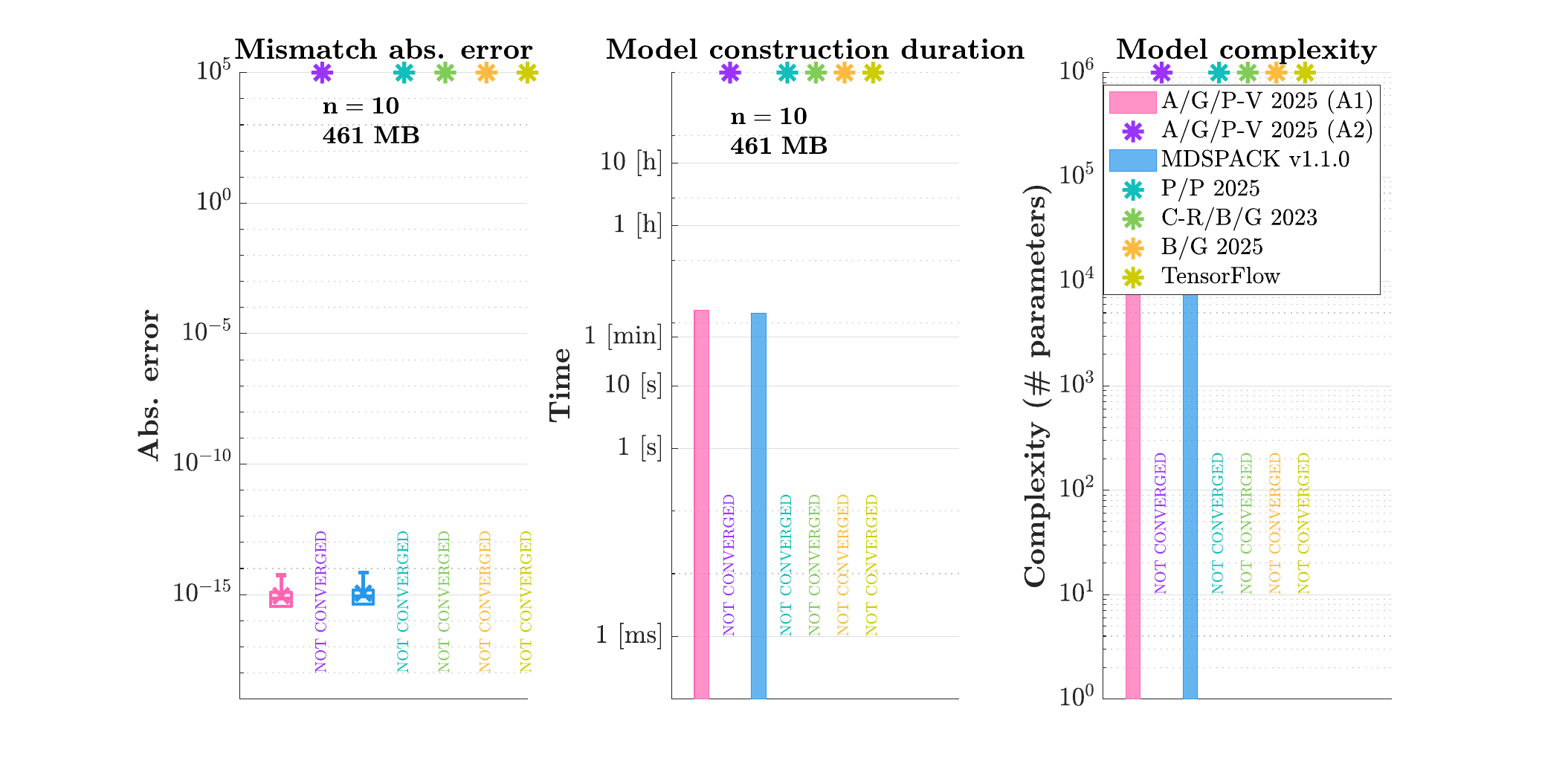} \caption{Function \#46: graphical view of the best model performances.} \end{figure}\begin{figure}[H] \centering  \includegraphics[width=\textwidth]{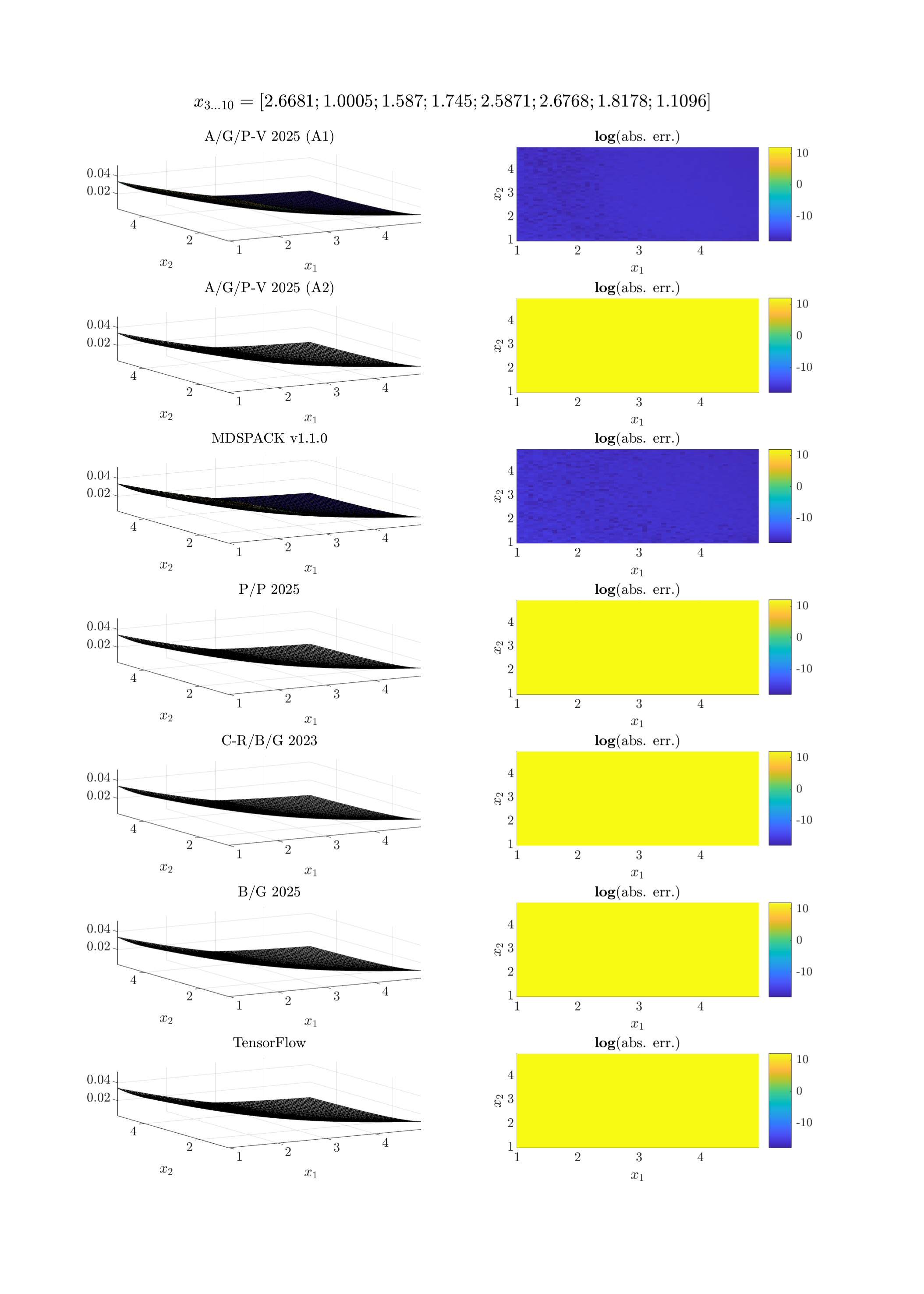} \caption{Function \#46: left side, evaluation of the original (mesh) vs. approximated (coloured surface) and right side, absolute errors (in log-scale).} \end{figure}\subsubsection{mLF detailed informations (M1)} \noindent \textbf{Right interpolation points}: $k_l=\left(\begin{array}{cccccccccc} 3 & 2 & 2 & 2 & 2 & 2 & 2 & 2 & 3 & 2 \end{array}\right)$, where $l=1,\cdots,\ord$.$$ \begin{array}{rcl}\lan{1} &\in& \IC^{3} \text{ , linearly spaced between bounds}\\\lan{2} &\in& \IC^{2} \text{ , linearly spaced between bounds}\\\lan{3} &\in& \IC^{2} \text{ , linearly spaced between bounds}\\\lan{4} &\in& \IC^{2} \text{ , linearly spaced between bounds}\\\lan{5} &\in& \IC^{2} \text{ , linearly spaced between bounds}\\\lan{6} &\in& \IC^{2} \text{ , linearly spaced between bounds}\\\lan{7} &\in& \IC^{2} \text{ , linearly spaced between bounds}\\\lan{8} &\in& \IC^{2} \text{ , linearly spaced between bounds}\\\lan{9} &\in& \IC^{3} \text{ , linearly spaced between bounds}\\\lan{10} &\in& \IC^{2} \text{ , linearly spaced between bounds}\\\end{array} $$\noindent \textbf{$\ord$-D Loewner matrix, barycentric weights and Lagrangian basis}:$$ \begin{array}{rcl}\IL & \in & \IC^{2304 \times 2304}\\\bc & \in & \IC^{2304}\\\bw & \in & \IC^{2304}\\\bc\odot \bw & \in & \IC^{2304}\\\mathbf{Lag}(\var{1},\var{2},\var{3},\var{4},\var{5},\var{6},\var{7},\var{8},\var{9},\var{10}) & \in & \IC^{2304}\\\end{array} $$

\newpage \subsection{Function \#47 (${\ord=5}$ variables, tensor size: 1.9 \textbf{MB})} $$\begin{array}{c}(1 + 2\var{1})(-2 + \var{2})(-\var{3})(3 + \var{4})(2- 3\var{5}) \\ + (-1 + \var{1})(2\var{2})(1 + 3\var{3})(-\var{4})(1 -\var{5})\end{array}$$ \subsubsection{Setup and results overview}\begin{itemize}\item Reference: G/al. 2025 (Ex 3.1), \cite{GHK:2025}\item Domain: $\mathbb{R}$\item Tensor size: 1.9 \textbf{MB} ($12^{5}$ points)\item Bounds: $ \left(\begin{array}{cc} -2 & 2 \end{array}\right) \times \left(\begin{array}{cc} -2 & 2 \end{array}\right) \times \left(\begin{array}{cc} -2 & 2 \end{array}\right) \times \left(\begin{array}{cc} -2 & 2 \end{array}\right) \times \left(\begin{array}{cc} -2 & 2 \end{array}\right)$ \end{itemize} \begin{table}[H] \centering \begin{tabular}{llllll}
$\#$ & Alg. & Parameters & Dim. & CPU [s] & RMSE \\ 
\hline 
$\mathbf{\#47}$ & A/G/P-V 2025 (A1) & $0.5,2$ & $\mathbf{2.2 \cdot 10^{02}}$ & $\mathbf{0.069}$ & $\mathbf{1 \cdot 10^{-13}}$ \\ 
 & A/G/P-V 2025 (A2) & $1 \cdot 10^{-15},1$ & $NaN$ & $NaN$ & $NaN$ \\ 
 & MDSPACK v1.1.0 & $0.01,1$ & $2.2 \cdot 10^{02}$ & $0.11$ & $1.1 \cdot 10^{-13}$ \\ 
 & P/P 2025 & $1,1,50,0.01,6,12,13$ & $5.5 \cdot 10^{02}$ & $91$ & $37$ \\ 
 & C-R/B/G 2023 & $0.001,20$ & $2.2 \cdot 10^{02}$ & $1.3 \cdot 10^{02}$ & $2.4 \cdot 10^{-12}$ \\ 
 & B/G 2025 & $1 \cdot 10^{-06},20,2$ & $2.2 \cdot 10^{02}$ & $1.8$ & $4.9 \cdot 10^{-13}$ \\ 
 & TensorFlow & $NaN$ & $NaN$ & $NaN$ & $NaN$ \\ 
\hline 
\end{tabular} \caption{Function \#47: best model configuration and performances per methods.} \end{table}\begin{figure}[H] \centering  \includegraphics[width=\textwidth]{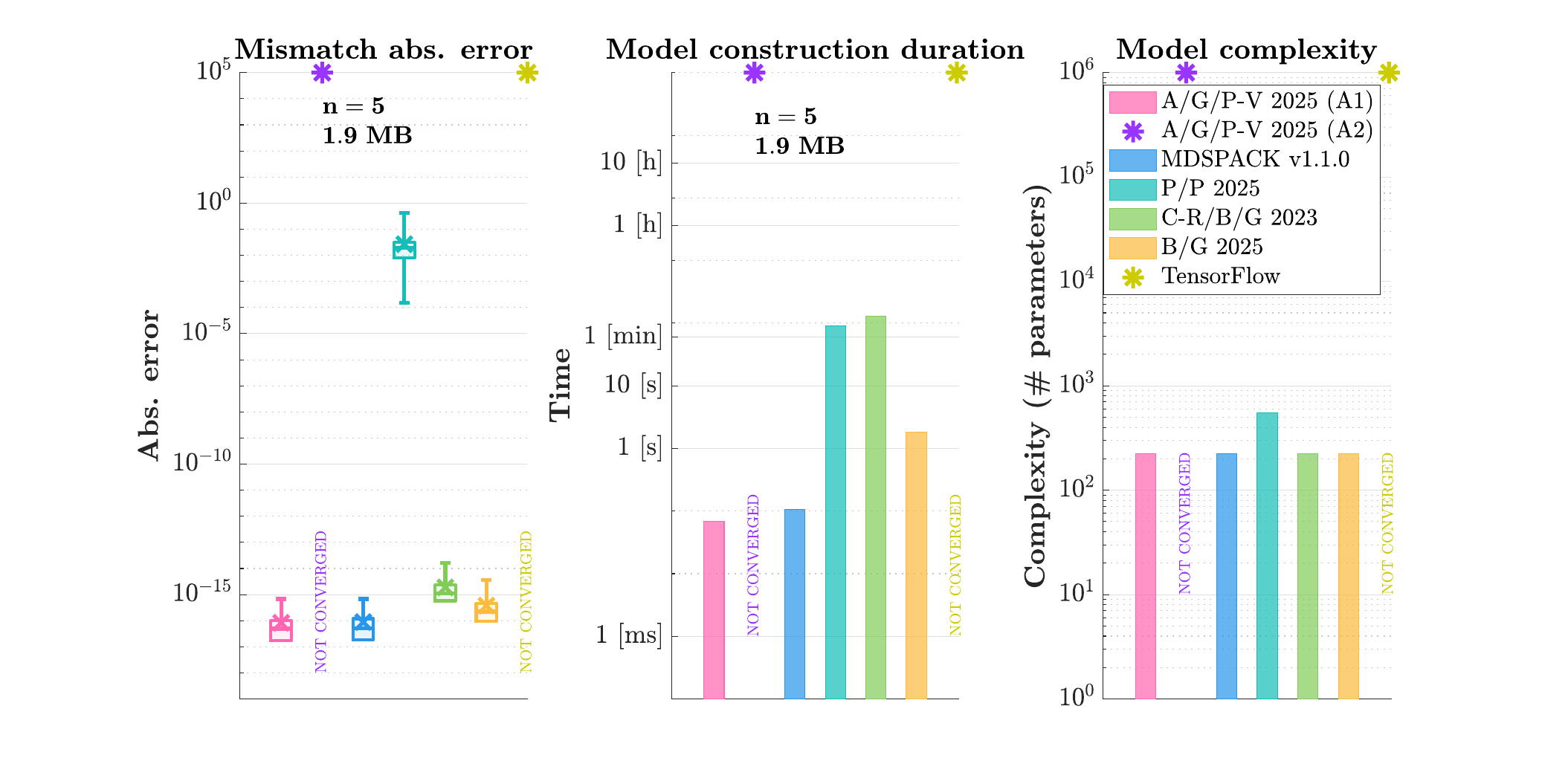} \caption{Function \#47: graphical view of the best model performances.} \end{figure}\begin{figure}[H] \centering  \includegraphics[width=\textwidth]{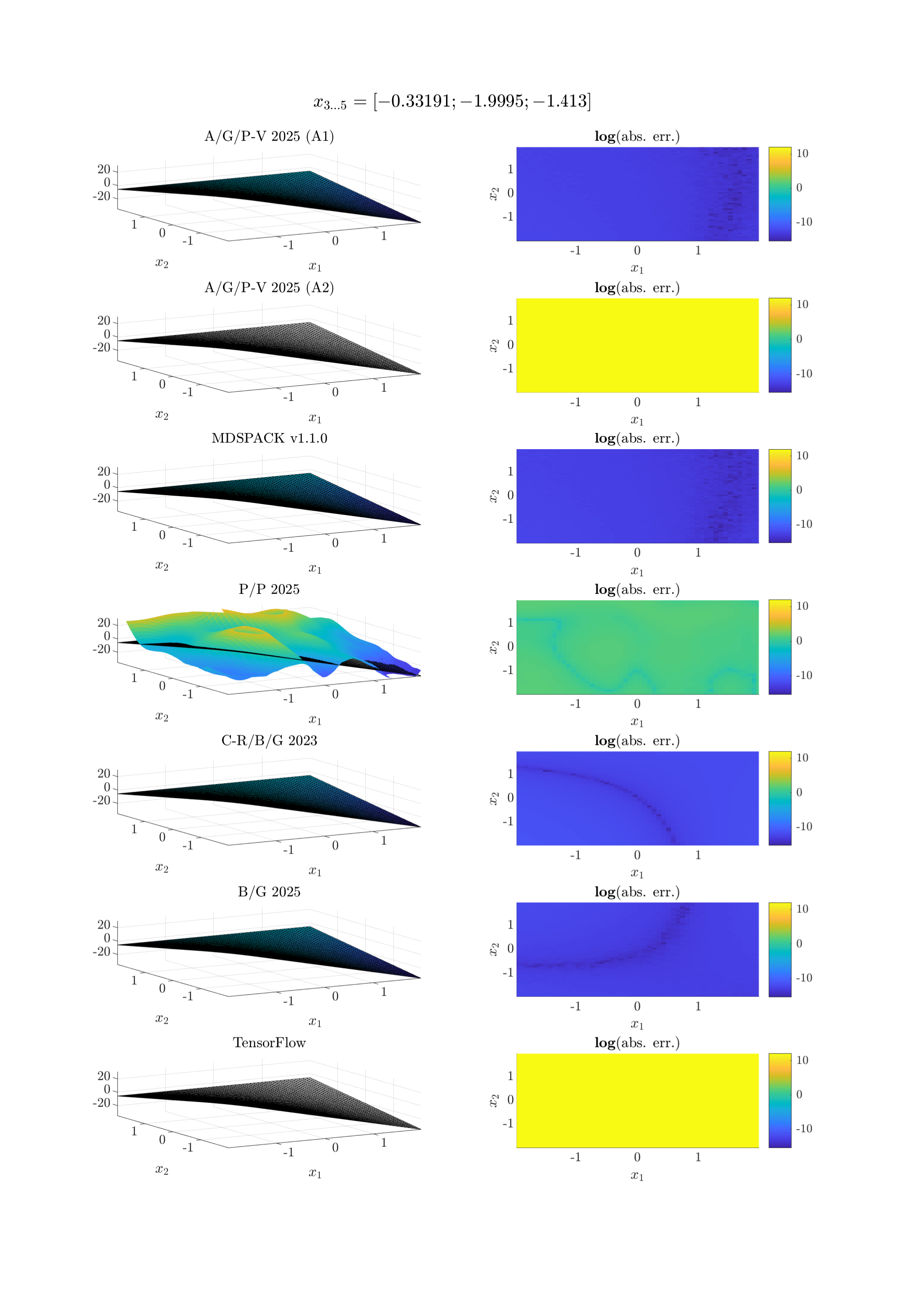} \caption{Function \#47: left side, evaluation of the original (mesh) vs. approximated (coloured surface) and right side, absolute errors (in log-scale).} \end{figure}\subsubsection{mLF detailed informations (M1)} \noindent \textbf{Right interpolation points}: $k_l=\left(\begin{array}{ccccc} 2 & 2 & 2 & 2 & 2 \end{array}\right)$, where $l=1,\cdots,\ord$.$$ \begin{array}{rcl}\lan{1} &\in& \IC^{2} \text{ , linearly spaced between bounds}\\\lan{2} &\in& \IC^{2} \text{ , linearly spaced between bounds}\\\lan{3} &\in& \IC^{2} \text{ , linearly spaced between bounds}\\\lan{4} &\in& \IC^{2} \text{ , linearly spaced between bounds}\\\lan{5} &\in& \IC^{2} \text{ , linearly spaced between bounds}\\\end{array} $$\noindent \textbf{$\ord$-D Loewner matrix, barycentric weights and Lagrangian basis}:$$ \begin{array}{rcl}\IL & \in & \IC^{32 \times 32}\\\bc & \in & \IC^{32}\\\bw & \in & \IC^{32}\\\bc\odot \bw & \in & \IC^{32}\\\mathbf{Lag}(\var{1},\var{2},\var{3},\var{4},\var{5}) & \in & \IC^{32}\\\end{array} $$

\newpage \subsection{Function \#48 (${\ord=3}$ variables, tensor size: 13.5 \textbf{KB})} $$\var{1}\var{2}+\var{1}\var{3}+\var{2}\var{3}$$ \subsubsection{Setup and results overview}\begin{itemize}\item Reference: G. P\'olya and G.Szeg\"o, \cite{Polya:1925}\item Domain: $\mathbb{R}$\item Tensor size: 13.5 \textbf{KB} ($12^{3}$ points)\item Bounds: $ \left(\begin{array}{cc} -\frac{1}{2} & 1 \end{array}\right) \times \left(\begin{array}{cc} -\frac{1}{2} & 1 \end{array}\right) \times \left(\begin{array}{cc} -\frac{1}{2} & 1 \end{array}\right)$ \end{itemize} \begin{table}[H] \centering \begin{tabular}{llllll}
$\#$ & Alg. & Parameters & Dim. & CPU [s] & RMSE \\ 
\hline 
$\mathbf{\#48}$ & A/G/P-V 2025 (A1) & $0.5,1$ & $\mathbf{40}$ & $0.051$ & $1.4 \cdot 10^{-16}$ \\ 
 & A/G/P-V 2025 (A2) & $1 \cdot 10^{-15},1$ & $40$ & $0.11$ & $1.2 \cdot 10^{-15}$ \\ 
 & MDSPACK v1.1.0 & $0.01,1$ & $40$ & $\mathbf{0.016}$ & $\mathbf{1.3 \cdot 10^{-16}}$ \\ 
 & P/P 2025 & $1,1,50,0.01,6,4,13$ & $2.9 \cdot 10^{02}$ & $0.49$ & $2.5 \cdot 10^{-15}$ \\ 
 & C-R/B/G 2023 & $0.001,20$ & $40$ & $0.029$ & $9 \cdot 10^{-16}$ \\ 
 & B/G 2025 & $1 \cdot 10^{-09},20,2$ & $40$ & $0.022$ & $1.7 \cdot 10^{-16}$ \\ 
 & TensorFlow & $NaN$ & $NaN$ & $NaN$ & $NaN$ \\ 
\hline 
\end{tabular} \caption{Function \#48: best model configuration and performances per methods.} \end{table}\begin{figure}[H] \centering  \includegraphics[width=\textwidth]{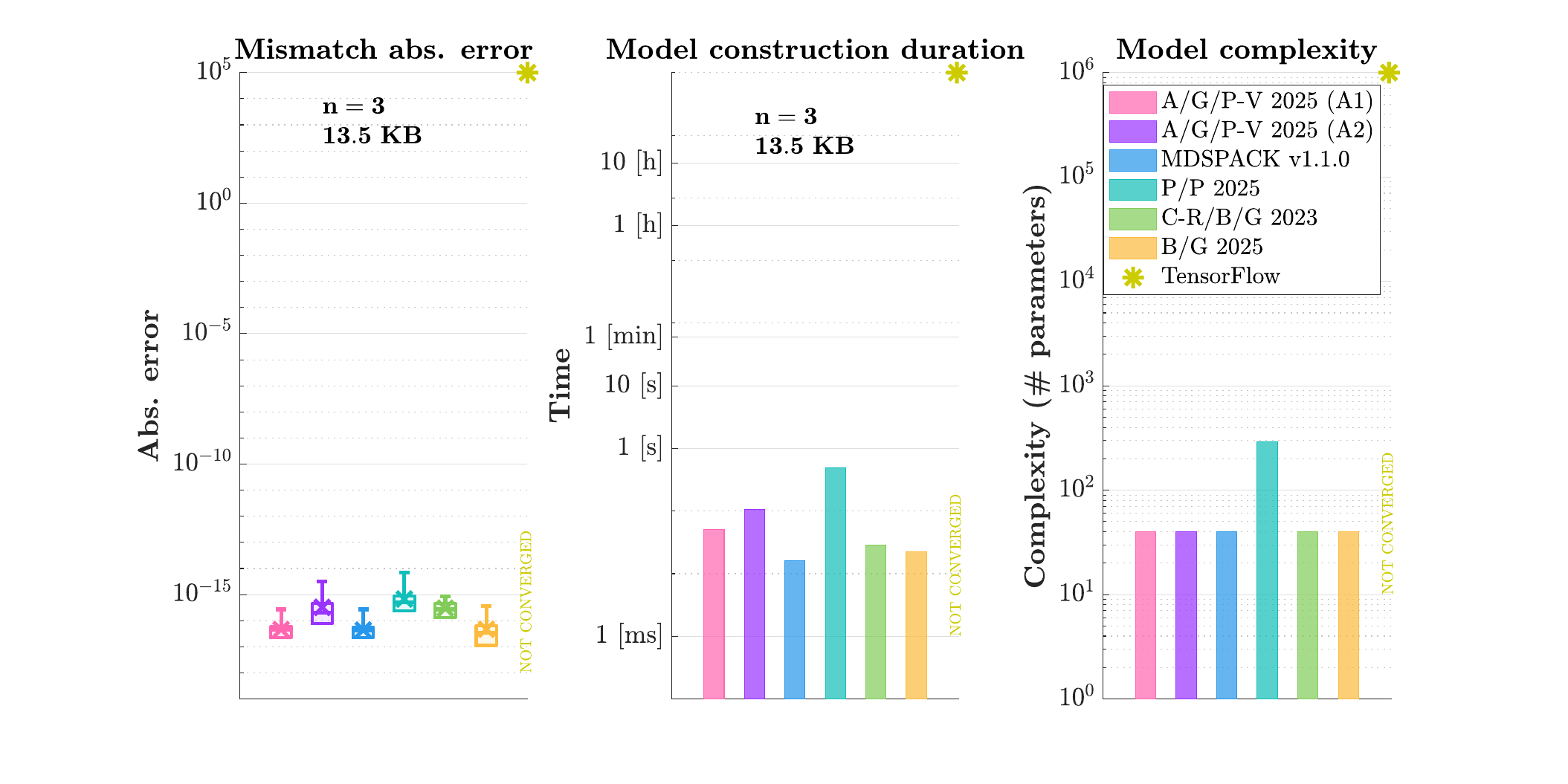} \caption{Function \#48: graphical view of the best model performances.} \end{figure}\begin{figure}[H] \centering  \includegraphics[width=\textwidth]{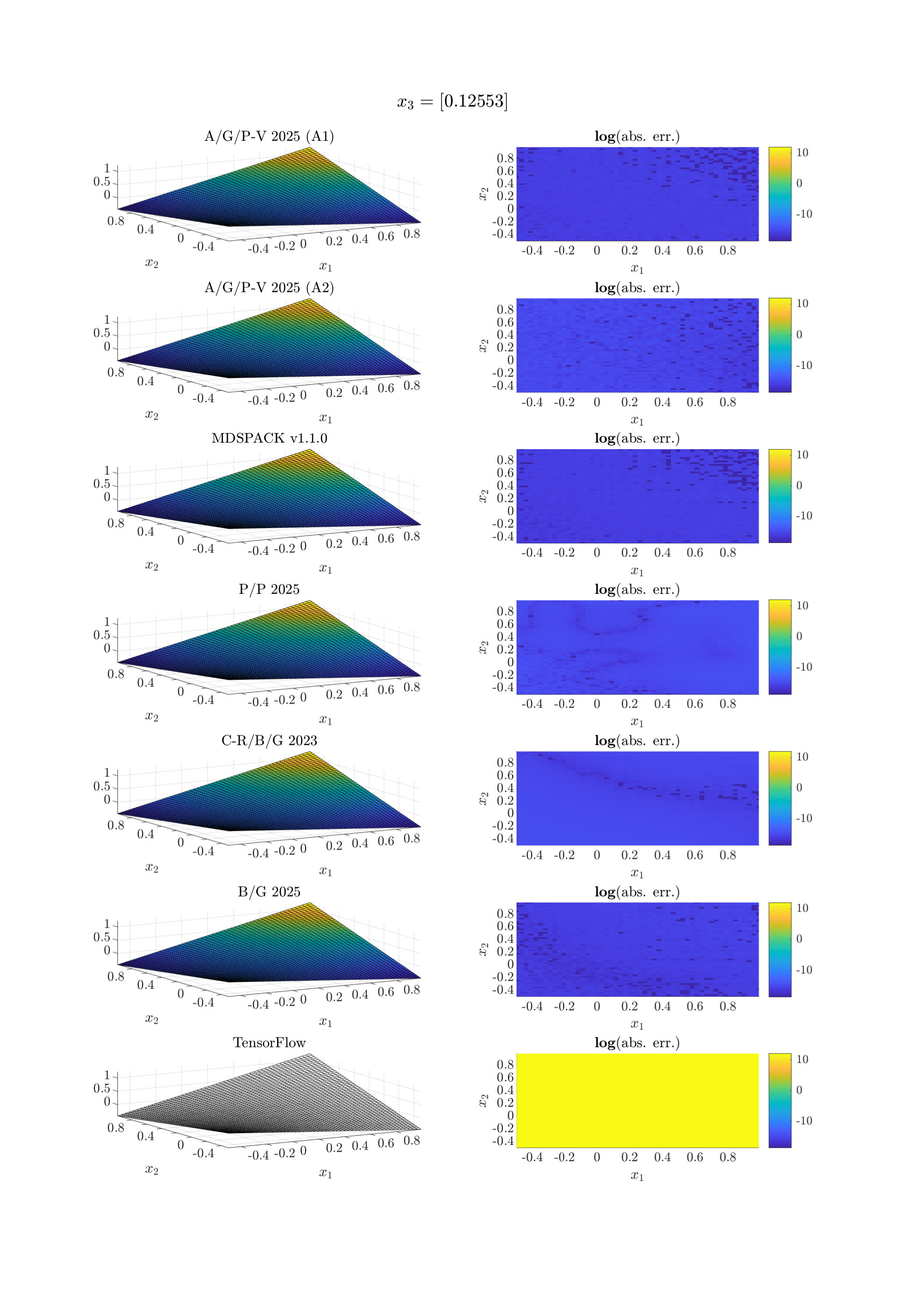} \caption{Function \#48: left side, evaluation of the original (mesh) vs. approximated (coloured surface) and right side, absolute errors (in log-scale).} \end{figure}\subsubsection{mLF detailed informations (M1)} \noindent \textbf{Right interpolation points} ($k_l=\left(\begin{array}{ccc} 2 & 2 & 2 \end{array}\right)$, where $l=1,\cdots,\ord$):$$ \begin{array}{rcl}\lan{1} &=& \left(\begin{array}{cc} -\frac{1}{2} & 1 \end{array}\right)\\\lan{2} &=& \left(\begin{array}{cc} -\frac{1}{2} & 1 \end{array}\right)\\\lan{3} &=& \left(\begin{array}{cc} -\frac{1}{2} & 1 \end{array}\right)\\\end{array} $$\noindent \textbf{Lagrangian weights}: $$\left(\begin{array}{ccc} \bc & \bw & \bc\odot\bw\\ -1.0 & 0.75 & -0.75\\ 1.0 & -0.75 & -0.75\\ 1.0 & -0.75 & -0.75\\ -1.0 & 0 & 0\\ 1.0 & -0.75 & -0.75\\ -1.0 & 0 & 0\\ -1.0 & 0 & 0\\ 1.0 & 3.0 & 3.0 \end{array}\right)$$\noindent \textbf{Lagrangian form} (basis, numerator and denominator coefficients):$$\left(\begin{array}{ccc}\mathcal{B}_\textrm{lag}(\var{1},\var{2},\var{3}) & \bN_\textrm{lag} &\bD_\textrm{lag}\end{array}\right) =$$ $$\left(\begin{array}{ccc} \left(\var{1}+0.5\right)\,\left(\var{2}+0.5\right)\,\left(\var{3}+0.5\right) & -0.75 & -1.0\\ \left(\var{1}+0.5\right)\,\left(\var{2}+0.5\right)\,\left(\var{3}-1.0\right) & -0.75 & 1.0\\ \left(\var{1}+0.5\right)\,\left(\var{2}-1.0\right)\,\left(\var{3}+0.5\right) & -0.75 & 1.0\\ \left(\var{1}+0.5\right)\,\left(\var{2}-1.0\right)\,\left(\var{3}-1.0\right) & 0 & -1.0\\ \left(\var{1}-1.0\right)\,\left(\var{2}+0.5\right)\,\left(\var{3}+0.5\right) & -0.75 & 1.0\\ \left(\var{1}-1.0\right)\,\left(\var{2}+0.5\right)\,\left(\var{3}-1.0\right) & 0 & -1.0\\ \left(\var{1}-1.0\right)\,\left(\var{2}-1.0\right)\,\left(\var{3}+0.5\right) & 0 & -1.0\\ \left(\var{1}-1.0\right)\,\left(\var{2}-1.0\right)\,\left(\var{3}-1.0\right) & 3.0 & 1.0 \end{array}\right).$$\noindent The corresponding function is:$$\begin{array}{rcl}\bG_{\textrm{lag}}(\var{1},\var{2},\var{3}) &=& \dfrac{\bn_{\textrm{lag}}(\var{1},\var{2},\var{3})}{\bd_{\textrm{lag}}(\var{1},\var{2},\var{3})}\\ && \\&=& \dfrac{\sum_{\textrm{row}} \bN_\textrm{lag} \odot\mathcal{B}^{-1}_\textrm{lag}(\var{1},\var{2},\var{3})}{\sum_{\textrm{row}} \bD_\textrm{lag} \odot\mathcal{B}^{-1}_\textrm{lag}(\var{1},\var{2},\var{3})}, \end{array}$$\noindent where,\\$\bn_{\textrm{lag}}(\var{1},\var{2},\var{3}) = \var{1}\,\var{2}+\var{1}\,\var{3}+\var{2}\,\var{3}$ \\~~\\$\bd_{\textrm{lag}}(\var{1},\var{2},\var{3}) = 1.0$ \\~~\\\noindent \textbf{Monomial form} (basis, numerator and denominator coefficients - entries $<10^{-12}$ removed):$$\left(\begin{array}{ccc}\mathcal{B}_\textrm{mon}(\var{1},\var{2},\var{3}) & \bN_\textrm{mon} &\bD_\textrm{mon}\end{array}\right) =$$ $$\left(\begin{array}{ccc} \var{1}\,\var{2}\,\var{3} & 0 & 0\\ \var{1}\,\var{2} & 1.0 & 0\\ \var{1}\,\var{3} & 1.0 & 0\\ \var{1} & 0 & 0\\ \var{2}\,\var{3} & 1.0 & 0\\ \var{2} & 0 & 0\\ \var{3} & 0 & 0\\ 1.0 & 0 & 1.0 \end{array}\right)$$\noindent The corresponding function is:$$\begin{array}{rcl}\bG_{\textrm{mon}}(\var{1},\var{2},\var{3}) &=& \dfrac{\bn_{\textrm{mon}}(\var{1},\var{2},\var{3})}{\bd_{\textrm{mon}}(\var{1},\var{2},\var{3})}\\ && \\&=& \dfrac{\sum_{\textrm{row}} \bN_\textrm{mon} \odot \mathcal{B}_\textrm{mon}(\var{1},\var{2},\var{3})}{\sum_{\textrm{row}} \bD_\textrm{mon} \odot\mathcal{B}_\textrm{mon}(\var{1},\var{2},\var{3})},  \end{array}$$\noindent where,\\$\bn_{\textrm{mon}}(\var{1},\var{2},\var{3}) = \var{1}\,\var{2}+\var{1}\,\var{3}+\var{2}\,\var{3}$ \\~~\\$\bd_{\textrm{mon}}(\var{1},\var{2},\var{3}) = 1.0$ \\~~\\\noindent \textbf{KST equivalent decoupling pattern} (Barycentric weights $\bc^{\var{l}}$): $$\begin{array}{rclll}\var{3}&: & \left(\begin{array}{cccc} -1.0 & -1.0 & -1.0 & -1.0\\ 1.0 & 1.0 & 1.0 & 1.0 \end{array}\right)& \textrm{vec}(.) & := \textbf{Bary}(\var{3}) \\\var{2}&: & \left(\begin{array}{cc} -1.0 & -1.0\\ 1.0 & 1.0 \end{array}\right)& \textrm{vec}(.) \otimes \bone_{k_{3}} & := \textbf{Bary}(\var{2}) \\\var{1}&: & \left(\begin{array}{c} -1.0\\ 1.0 \end{array}\right)& \textrm{vec}(.) \otimes \bone_{k_{3}k_{2}} & := \textbf{Bary}(\var{1}) \\\end{array}$$~\\ Then, with the above notations, one defines the following univariate vector functions:  $$ \left\{ \begin{array}{rcl}\bPhi_{1}(\var{1}) &:=& \textbf{Bary}(\var{1}) \odot \mathbf{Lag}(\var{1}) \\\bPhi_{2}(\var{2}) &:=& \textbf{Bary}(\var{2}) \odot \mathbf{Lag}(\var{2}) \\\bPhi_{3}(\var{3}) &:=& \textbf{Bary}(\var{3}) \odot \mathbf{Lag}(\var{3}) \\\end{array} \right. $$\noindent The corresponding function is:$$\begin{array}{rcl}\bG_{\textrm{kst}}(\var{1},\var{2},\var{3}) &=& \dfrac{\bn_{\textrm{kst}}(\var{1},\var{2},\var{3})}{\bd_{\textrm{kst}}(\var{1},\var{2},\var{3})}\\ && \\ &=& \dfrac{\sum_{\text{rows}} \bw \odot \bPhi_{1}(\var{1}) \odot \cdots \odot\bPhi_{3}(\var{3})}{\sum_{\text{rows}} \bPhi_{1}(\var{1}) \odot \cdots \odot\bPhi_{3}(\var{3})} . \end{array}$$~\\ \noindent \textbf{KST-like univariate functions} (equivalent scaled univariate functions $\bphi_{1,\cdots,3}$): $$\left\{\begin{array}{rcrcl}z_{1} &=&\bphi_{1}(\var{1}) &=& 2.0\,\var{1}+1.0\\z_{2} &=&\bphi_{2}(\var{2}) &=& 0.5\,\var{2}-0.5\\z_{3} &=&\bphi_{3}(\var{3}) &=& 0.25-1.0\,\var{3}\\\end{array} \right. .$$\noindent \textbf{Connection with Neural Networks} (equivalent numerator $\bn_{\textrm{lag}}$ representation):\\ \begin{figure}[H]\begin{center} \scalebox{.7}{\begin{tikzpicture}[line width=0.4mm]\tikzstyle{place}=[circle, draw=black, minimum size = 8mm]\tikzstyle{placeInOut}=[circle, draw=orange, minimum size = 8mm]\node at (0,-2) [placeInOut] (first_1){$\var{1}$};\node at (0,-4) [placeInOut] (first_2){$\var{2}$};\node at (0,-6) [placeInOut] (first_3){$\var{3}$};\node at (5,-2) [place] (secondL1_1){$\frac{1}{\var{1}-\lani{1}{1}}$};\node at (5,-4) [place] (secondL1_2){$\frac{1}{\var{1}-\lani{1}{2}}$};\node at (5,-6) [place] (secondL2_1){$\frac{1}{\var{2}-\lani{2}{1}}$};\node at (5,-8) [place] (secondL2_2){$\frac{1}{\var{2}-\lani{2}{2}}$};\node at (5,-10) [place] (secondL3_1){$\frac{1}{\var{3}-\lani{3}{1}}$};\node at (5,-12) [place] (secondL3_2){$\frac{1}{\var{3}-\lani{3}{2}}$};\node at (10,-2) [place] (third_1){$\prod$};\node at (10,-4) [place] (third_2){$\prod$};\node at (10,-6) [place] (third_3){$\prod$};\node at (10,-8) [place] (third_4){$\prod$};\node at (10,-10) [place] (third_5){$\prod$};\node at (10,-12) [place] (third_6){$\prod$};\node at (10,-14) [place] (third_7){$\prod$};\node at (10,-16) [place] (third_8){$\prod$};\node at (15,-9) [placeInOut] (output){$\bSigma$};\draw[->] (first_1)--(secondL1_1) node[above,sloped,pos=0.75] { };\draw[->] (first_1)--(secondL1_2) node[above,sloped,pos=0.75] { };\draw[->] (first_2)--(secondL2_1) node[above,sloped,pos=0.75] { };\draw[->] (first_2)--(secondL2_2) node[above,sloped,pos=0.75] { };\draw[->] (first_3)--(secondL3_1) node[above,sloped,pos=0.75] { };\draw[->] (first_3)--(secondL3_2) node[above,sloped,pos=0.75] { };\draw[->] (secondL1_1)--(third_1) node[above,sloped,pos=0.25] {};\draw[->] (secondL1_1)--(third_2) node[above,sloped,pos=0.25] {};\draw[->] (secondL1_1)--(third_3) node[above,sloped,pos=0.25] {};\draw[->] (secondL1_1)--(third_4) node[above,sloped,pos=0.25] {};\draw[->] (secondL1_2)--(third_5) node[above,sloped,pos=0.25] {};\draw[->] (secondL1_2)--(third_6) node[above,sloped,pos=0.25] {};\draw[->] (secondL1_2)--(third_7) node[above,sloped,pos=0.25] {};\draw[->] (secondL1_2)--(third_8) node[above,sloped,pos=0.25] {};\draw[->] (secondL2_1)--(third_1) node[above,sloped,pos=0.25] {};\draw[->] (secondL2_1)--(third_2) node[above,sloped,pos=0.25] {};\draw[->] (secondL2_2)--(third_3) node[above,sloped,pos=0.25] {};\draw[->] (secondL2_2)--(third_4) node[above,sloped,pos=0.25] {};\draw[->] (secondL2_1)--(third_5) node[above,sloped,pos=0.25] {};\draw[->] (secondL2_1)--(third_6) node[above,sloped,pos=0.25] {};\draw[->] (secondL2_2)--(third_7) node[above,sloped,pos=0.25] {};\draw[->] (secondL2_2)--(third_8) node[above,sloped,pos=0.25] {};\draw[->] (secondL3_1)--(third_1) node[above,sloped,pos=0.25] {};\draw[->] (secondL3_2)--(third_2) node[above,sloped,pos=0.25] {};\draw[->] (secondL3_1)--(third_3) node[above,sloped,pos=0.25] {};\draw[->] (secondL3_2)--(third_4) node[above,sloped,pos=0.25] {};\draw[->] (secondL3_1)--(third_5) node[above,sloped,pos=0.25] {};\draw[->] (secondL3_2)--(third_6) node[above,sloped,pos=0.25] {};\draw[->] (secondL3_1)--(third_7) node[above,sloped,pos=0.25] {};\draw[->] (secondL3_2)--(third_8) node[above,sloped,pos=0.25] {};\draw[->] (third_1)--(output) node[above,sloped,pos=0.25] {-0.75};\draw[->] (third_2)--(output) node[above,sloped,pos=0.25] {-0.75};\draw[->] (third_3)--(output) node[above,sloped,pos=0.25] {-0.75};\draw[->] (third_4)--(output) node[above,sloped,pos=0.25] {0};\draw[->] (third_5)--(output) node[above,sloped,pos=0.25] {-0.75};\draw[->] (third_6)--(output) node[above,sloped,pos=0.25] {0};\draw[->] (third_7)--(output) node[above,sloped,pos=0.25] {0};\draw[->] (third_8)--(output) node[above,sloped,pos=0.25] {3};\end{tikzpicture}} \caption{Equivalent NN representation of the numerator $\bn_{\textrm{lag}}$.}\end{center}\end{figure}

\newpage \subsection{Function \#49 (${\ord=2}$ variables, tensor size: 50 \textbf{KB})} $$\texttt{Hankel function $H_0$ (real part)}$$ \subsubsection{Setup and results overview}\begin{itemize}\item Reference: Hankel function, [none]\item Domain: $\mathbb{R}$\item Tensor size: 50 \textbf{KB} ($80^{2}$ points)\item Bounds: $ \left(\begin{array}{cc} 1 & 10 \end{array}\right) \times \left(\begin{array}{cc} \frac{1}{10} & 1 \end{array}\right)$ \end{itemize} \begin{table}[H] \centering \begin{tabular}{llllll}
$\#$ & Alg. & Parameters & Dim. & CPU [s] & RMSE \\ 
\hline 
$\mathbf{\#49}$ & A/G/P-V 2025 (A1) & $0.0001,3$ & $1.3 \cdot 10^{02}$ & $0.038$ & $0.016$ \\ 
 & A/G/P-V 2025 (A2) & $1 \cdot 10^{-15},1$ & $\mathbf{1.2 \cdot 10^{02}}$ & $1.7$ & $0.11$ \\ 
 & MDSPACK v1.1.0 & $1 \cdot 10^{-06},3$ & $1.3 \cdot 10^{02}$ & $\mathbf{0.031}$ & $0.017$ \\ 
 & P/P 2025 & $1,1,50,0.01,10,4,21$ & $5.1 \cdot 10^{02}$ & $3$ & $4.6 \cdot 10^{-05}$ \\ 
 & C-R/B/G 2023 & $1 \cdot 10^{-09},20$ & $7.2 \cdot 10^{02}$ & $0.45$ & $\mathbf{2.7 \cdot 10^{-13}}$ \\ 
 & B/G 2025 & $1 \cdot 10^{-06},20,4$ & $5.3 \cdot 10^{02}$ & $7.6$ & $1 \cdot 10^{-08}$ \\ 
 & TensorFlow & $NaN$ & $NaN$ & $NaN$ & $NaN$ \\ 
\hline 
\end{tabular} \caption{Function \#49: best model configuration and performances per methods.} \end{table}\begin{figure}[H] \centering  \includegraphics[width=\textwidth]{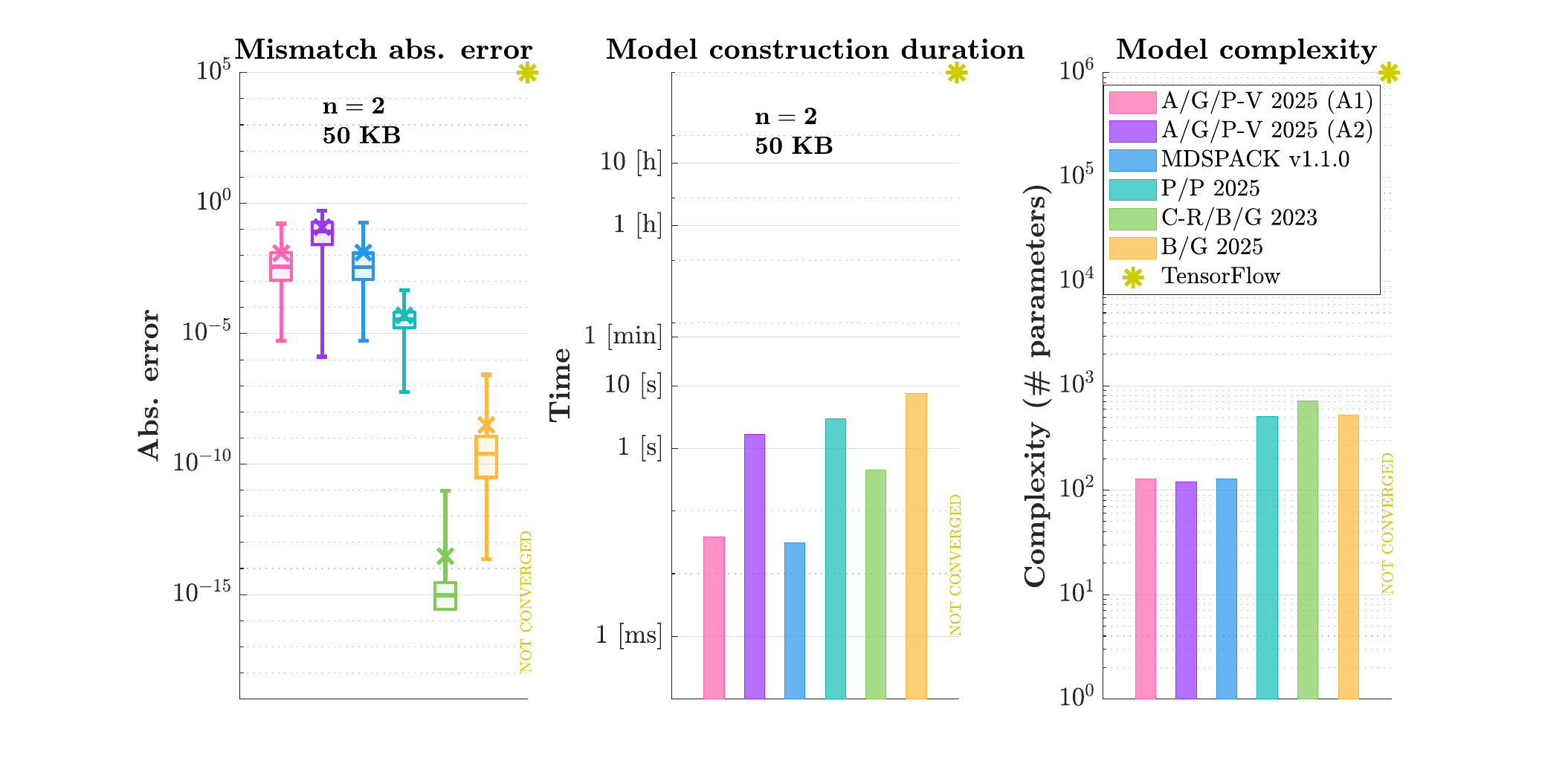} \caption{Function \#49: graphical view of the best model performances.} \end{figure}\begin{figure}[H] \centering  \includegraphics[width=\textwidth]{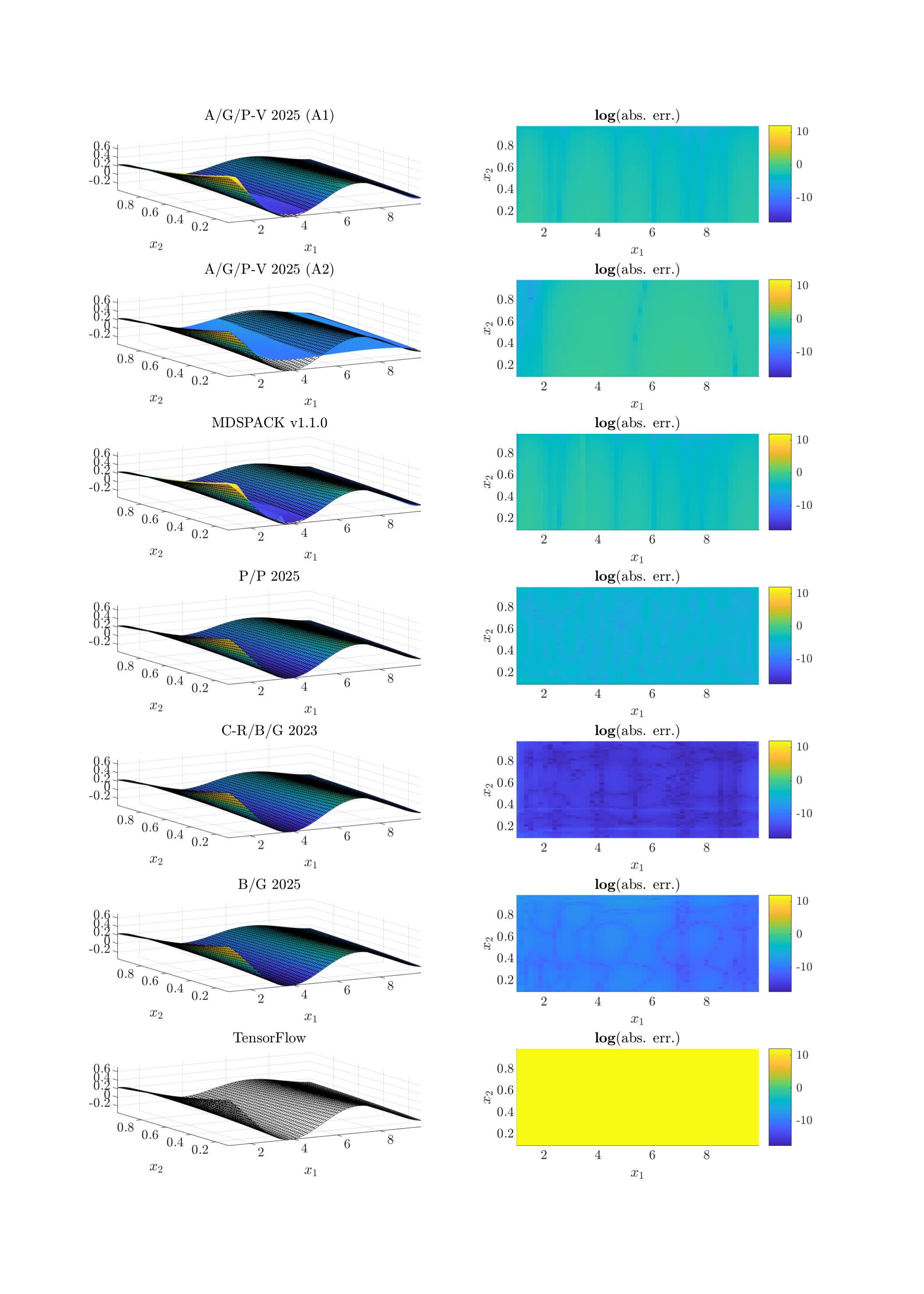} \caption{Function \#49: left side, evaluation of the original (mesh) vs. approximated (coloured surface) and right side, absolute errors (in log-scale).} \end{figure}\subsubsection{mLF detailed informations (M1)} \noindent \textbf{Right interpolation points}: $k_l=\left(\begin{array}{cc} 8 & 4 \end{array}\right)$, where $l=1,\cdots,\ord$.$$ \begin{array}{rcl}\lan{1} &\in& \IC^{8} \text{ , linearly spaced between bounds}\\\lan{2} &\in& \IC^{4} \text{ , linearly spaced between bounds}\\\end{array} $$\noindent \textbf{$\ord$-D Loewner matrix, barycentric weights and Lagrangian basis}:$$ \begin{array}{rcl}\IL & \in & \IC^{32 \times 32}\\\bc & \in & \IC^{32}\\\bw & \in & \IC^{32}\\\bc\odot \bw & \in & \IC^{32}\\\mathbf{Lag}(\var{1},\var{2}) & \in & \IC^{32}\\\end{array} $$

\newpage \subsection{Function \#50 (${\ord=2}$ variables, tensor size: 50 \textbf{KB})} $$\texttt{Hankel function $H_0$ (imaginary part)}$$ \subsubsection{Setup and results overview}\begin{itemize}\item Reference: Hankel function, [none]\item Domain: $\mathbb{R}$\item Tensor size: 50 \textbf{KB} ($80^{2}$ points)\item Bounds: $ \left(\begin{array}{cc} 1 & 10 \end{array}\right) \times \left(\begin{array}{cc} \frac{1}{10} & 1 \end{array}\right)$ \end{itemize} \begin{table}[H] \centering \begin{tabular}{llllll}
$\#$ & Alg. & Parameters & Dim. & CPU [s] & RMSE \\ 
\hline 
$\mathbf{\#50}$ & A/G/P-V 2025 (A1) & $0.0001,3$ & $1.1 \cdot 10^{02}$ & $0.037$ & $0.0098$ \\ 
 & A/G/P-V 2025 (A2) & $1 \cdot 10^{-15},1$ & $\mathbf{4}$ & $1.4$ & $0.45$ \\ 
 & MDSPACK v1.1.0 & $1 \cdot 10^{-08},4$ & $1.1 \cdot 10^{02}$ & $\mathbf{0.031}$ & $0.0098$ \\ 
 & P/P 2025 & $1,1,50,0.01,6,6,13$ & $2.4 \cdot 10^{02}$ & $2.1$ & $4 \cdot 10^{-05}$ \\ 
 & C-R/B/G 2023 & $1 \cdot 10^{-06},20$ & $5.4 \cdot 10^{02}$ & $0.48$ & $1 \cdot 10^{-06}$ \\ 
 & B/G 2025 & $1 \cdot 10^{-09},20,2$ & $7.3 \cdot 10^{02}$ & $1.8$ & $\mathbf{7.5 \cdot 10^{-08}}$ \\ 
 & TensorFlow & $NaN$ & $NaN$ & $NaN$ & $NaN$ \\ 
\hline 
\end{tabular} \caption{Function \#50: best model configuration and performances per methods.} \end{table}\begin{figure}[H] \centering  \includegraphics[width=\textwidth]{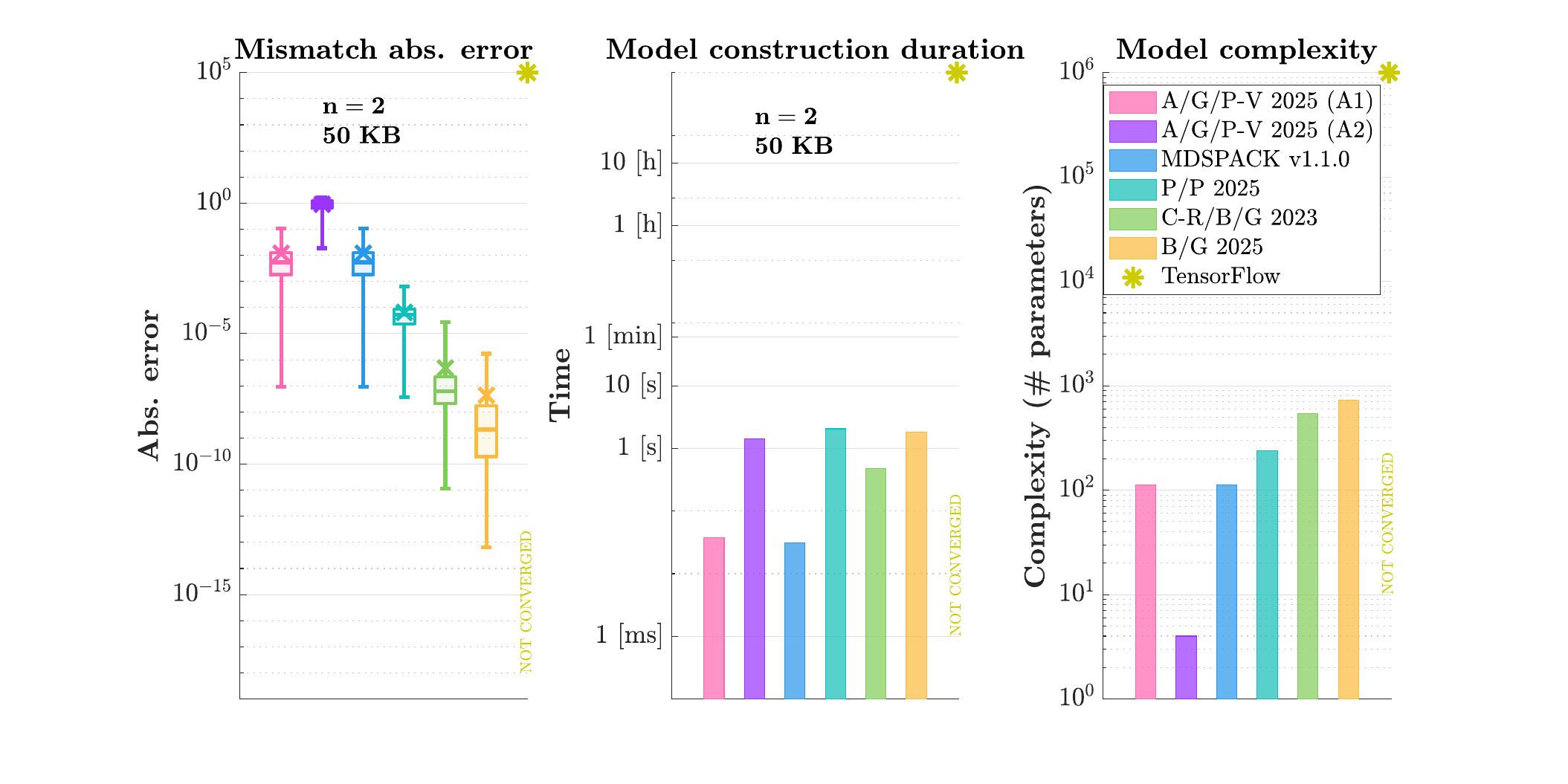} \caption{Function \#50: graphical view of the best model performances.} \end{figure}\begin{figure}[H] \centering  \includegraphics[width=\textwidth]{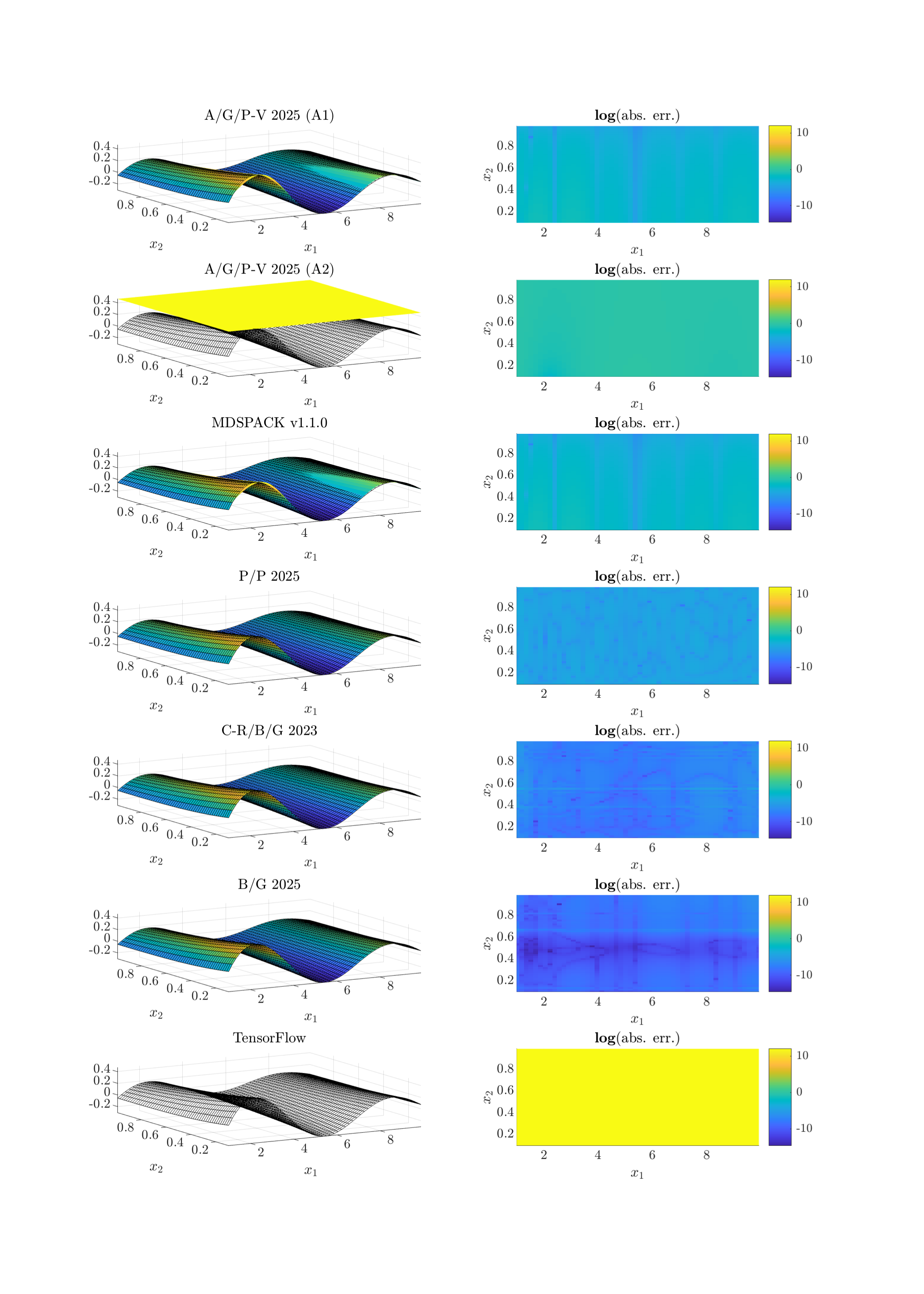} \caption{Function \#50: left side, evaluation of the original (mesh) vs. approximated (coloured surface) and right side, absolute errors (in log-scale).} \end{figure}\subsubsection{mLF detailed informations (M1)} \noindent \textbf{Right interpolation points} ($k_l=\left(\begin{array}{cc} 7 & 4 \end{array}\right)$, where $l=1,\cdots,\ord$):$$ \begin{array}{rcl}\lan{1} &=& \left(\begin{array}{ccccccc} 1 & \frac{31}{13} & 4 & \frac{70}{13} & 7 & \frac{109}{13} & 10 \end{array}\right)\\\lan{2} &=& \left(\begin{array}{cccc} \frac{1}{10} & \frac{2}{5} & \frac{7}{10} & 1 \end{array}\right)\\\end{array} $$\noindent \textbf{Lagrangian weights}: $$\left(\begin{array}{ccc} \bc & \bw & \bc\odot\bw\\ 0.01018 & 0.04853 & 0.0004939\\ -0.04269 & -0.02484 & 0.00106\\ 0.06216 & -0.04978 & -0.003094\\ -0.03036 & -0.05106 & 0.00155\\ -0.1355 & 0.4619 & -0.06259\\ 0.4679 & 0.3369 & 0.1576\\ -0.5414 & 0.2435 & -0.1318\\ 0.2096 & 0.1746 & 0.0366\\ 0.392 & -0.01089 & -0.00427\\ -1.402 & 0.001665 & -0.002334\\ 1.67 & 0.008203 & 0.0137\\ -0.6616 & 0.01097 & -0.00726\\ -0.6529 & -0.3071 & 0.2005\\ 2.265 & -0.2259 & -0.5118\\ -2.612 & -0.1659 & 0.4331\\ 1.0 & -0.1215 & -0.1215\\ 0.5474 & -0.0254 & -0.0139\\ -1.976 & -0.02305 & 0.04554\\ 2.367 & -0.02015 & -0.04769\\ -0.9409 & -0.01715 & 0.01613\\ -0.3498 & 0.2397 & -0.08384\\ 1.223 & 0.1764 & 0.2157\\ -1.419 & 0.1298 & -0.1841\\ 0.5463 & 0.09537 & 0.0521\\ 0.09286 & 0.05147 & 0.00478\\ -0.3064 & 0.04056 & -0.01243\\ 0.337 & 0.03181 & 0.01072\\ -0.1235 & 0.02485 & -0.003069 \end{array}\right)$$\noindent \textbf{Lagrangian form} (basis, numerator and denominator coefficients):$$\left(\begin{array}{ccc}\mathcal{B}_\textrm{lag}(\var{1},\var{2}) & \bN_\textrm{lag} &\bD_\textrm{lag}\end{array}\right) =$$ $$\left(\begin{array}{ccc} \left(\var{1}-1.0\right)\,\left(\var{2}-0.1\right) & 0.0004939 & 0.01018\\ \left(\var{1}-1.0\right)\,\left(\var{2}-0.4\right) & 0.00106 & -0.04269\\ \left(\var{1}-1.0\right)\,\left(\var{2}-0.7\right) & -0.003094 & 0.06216\\ \left(\var{1}-1.0\right)\,\left(\var{2}-1.0\right) & 0.00155 & -0.03036\\ \left(\var{1}-2.385\right)\,\left(\var{2}-0.1\right) & -0.06259 & -0.1355\\ \left(\var{1}-2.385\right)\,\left(\var{2}-0.4\right) & 0.1576 & 0.4679\\ \left(\var{1}-2.385\right)\,\left(\var{2}-0.7\right) & -0.1318 & -0.5414\\ \left(\var{2}-1.0\right)\,\left(\var{1}-2.385\right) & 0.0366 & 0.2096\\ \left(\var{1}-4.0\right)\,\left(\var{2}-0.1\right) & -0.00427 & 0.392\\ \left(\var{1}-4.0\right)\,\left(\var{2}-0.4\right) & -0.002334 & -1.402\\ \left(\var{1}-4.0\right)\,\left(\var{2}-0.7\right) & 0.0137 & 1.67\\ \left(\var{2}-1.0\right)\,\left(\var{1}-4.0\right) & -0.00726 & -0.6616\\ \left(\var{1}-5.385\right)\,\left(\var{2}-0.1\right) & 0.2005 & -0.6529\\ \left(\var{1}-5.385\right)\,\left(\var{2}-0.4\right) & -0.5118 & 2.265\\ \left(\var{1}-5.385\right)\,\left(\var{2}-0.7\right) & 0.4331 & -2.612\\ \left(\var{2}-1.0\right)\,\left(\var{1}-5.385\right) & -0.1215 & 1.0\\ \left(\var{2}-0.1\right)\,\left(\var{1}-7.0\right) & -0.0139 & 0.5474\\ \left(\var{2}-0.4\right)\,\left(\var{1}-7.0\right) & 0.04554 & -1.976\\ \left(\var{2}-0.7\right)\,\left(\var{1}-7.0\right) & -0.04769 & 2.367\\ \left(\var{2}-1.0\right)\,\left(\var{1}-7.0\right) & 0.01613 & -0.9409\\ \left(\var{1}-8.385\right)\,\left(\var{2}-0.1\right) & -0.08384 & -0.3498\\ \left(\var{1}-8.385\right)\,\left(\var{2}-0.4\right) & 0.2157 & 1.223\\ \left(\var{1}-8.385\right)\,\left(\var{2}-0.7\right) & -0.1841 & -1.419\\ \left(\var{2}-1.0\right)\,\left(\var{1}-8.385\right) & 0.0521 & 0.5463\\ \left(\var{1}-10.0\right)\,\left(\var{2}-0.1\right) & 0.00478 & 0.09286\\ \left(\var{1}-10.0\right)\,\left(\var{2}-0.4\right) & -0.01243 & -0.3064\\ \left(\var{1}-10.0\right)\,\left(\var{2}-0.7\right) & 0.01072 & 0.337\\ \left(\var{2}-1.0\right)\,\left(\var{1}-10.0\right) & -0.003069 & -0.1235 \end{array}\right).$$\noindent The corresponding function is:$$\begin{array}{rcl}\bG_{\textrm{lag}}(\var{1},\var{2}) &=& \dfrac{\bn_{\textrm{lag}}(\var{1},\var{2})}{\bd_{\textrm{lag}}(\var{1},\var{2})}\\ && \\&=& \dfrac{\sum_{\textrm{row}} \bN_\textrm{lag} \odot\mathcal{B}^{-1}_\textrm{lag}(\var{1},\var{2})}{\sum_{\textrm{row}} \bD_\textrm{lag} \odot\mathcal{B}^{-1}_\textrm{lag}(\var{1},\var{2})}, \end{array}$$\noindent where,\\$\bn_{\textrm{lag}}(\var{1},\var{2}) = -2.962 \cdot 10^{-7}\,{\var{1}}^6\,{\var{2}}^3-2.327 \cdot 10^{-6}\,{\var{1}}^6\,{\var{2}}^2+2.579 \cdot 10^{-5}\,{\var{1}}^6\,\var{2}-5.661 \cdot 10^{-5}\,{\var{1}}^6+1.347 \cdot 10^{-5}\,{\var{1}}^5\,{\var{2}}^3+5.447 \cdot 10^{-5}\,{\var{1}}^5\,{\var{2}}^2-0.0008236\,{\var{1}}^5\,\var{2}+0.001932\,{\var{1}}^5-0.0002432\,{\var{1}}^4\,{\var{2}}^3-0.0001936\,{\var{1}}^4\,{\var{2}}^2+0.009288\,{\var{1}}^4\,\var{2}-0.02394\,{\var{1}}^4+0.002159\,{\var{1}}^3\,{\var{2}}^3-0.003791\,{\var{1}}^3\,{\var{2}}^2-0.04201\,{\var{1}}^3\,\var{2}+0.1286\,{\var{1}}^3-0.009479\,{\var{1}}^2\,{\var{2}}^3+0.03398\,{\var{1}}^2\,{\var{2}}^2+0.051\,{\var{1}}^2\,\var{2}-0.2701\,{\var{1}}^2+0.01811\,\var{1}\,{\var{2}}^3-0.08409\,\var{1}\,{\var{2}}^2+0.06805\,\var{1}\,\var{2}+0.1074\,\var{1}-0.009959\,{\var{2}}^3+0.04434\,{\var{2}}^2-0.05386\,\var{2}+0.04834$ \\~~\\$\bd_{\textrm{lag}}(\var{1},\var{2}) = -2.162 \cdot 10^{-5}\,{\var{1}}^6\,{\var{2}}^3-5.912 \cdot 10^{-5}\,{\var{1}}^6\,{\var{2}}^2-5.154 \cdot 10^{-5}\,{\var{1}}^6\,\var{2}+0.0001316\,{\var{1}}^6+0.0007204\,{\var{1}}^5\,{\var{2}}^3+0.00196\,{\var{1}}^5\,{\var{2}}^2+0.001726\,{\var{1}}^5\,\var{2}-0.004425\,{\var{1}}^5-0.009428\,{\var{1}}^4\,{\var{2}}^3-0.02537\,{\var{1}}^4\,{\var{2}}^2-0.02239\,{\var{1}}^4\,\var{2}+0.05783\,{\var{1}}^4+0.06157\,{\var{1}}^3\,{\var{2}}^3+0.1622\,{\var{1}}^3\,{\var{2}}^2+0.1412\,{\var{1}}^3\,\var{2}-0.3737\,{\var{1}}^3-0.2103\,{\var{1}}^2\,{\var{2}}^3-0.5333\,{\var{1}}^2\,{\var{2}}^2-0.4469\,{\var{1}}^2\,\var{2}+1.243\,{\var{1}}^2+0.3552\,\var{1}\,{\var{2}}^3+0.8409\,\var{1}\,{\var{2}}^2+0.63\,\var{1}\,\var{2}-2.011\,\var{1}-0.24\,{\var{2}}^3-0.5247\,{\var{2}}^2-0.3844\,\var{2}+1.0$ \\~~\\\noindent \textbf{Monomial form} (basis, numerator and denominator coefficients - entries $<10^{-12}$ removed):$$\left(\begin{array}{ccc}\mathcal{B}_\textrm{mon}(\var{1},\var{2}) & \bN_\textrm{mon} &\bD_\textrm{mon}\end{array}\right) =$$ $$\left(\begin{array}{ccc} {\var{1}}^6\,{\var{2}}^3 & -1.473 \cdot 10^{-7} & -1.075 \cdot 10^{-5}\\ {\var{1}}^6\,{\var{2}}^2 & -1.157 \cdot 10^{-6} & -2.939 \cdot 10^{-5}\\ {\var{1}}^6\,\var{2} & 1.282 \cdot 10^{-5} & -2.563 \cdot 10^{-5}\\ {\var{1}}^6 & -2.815 \cdot 10^{-5} & 6.543 \cdot 10^{-5}\\ {\var{1}}^5\,{\var{2}}^3 & 6.696 \cdot 10^{-6} & 0.0003582\\ {\var{1}}^5\,{\var{2}}^2 & 2.708 \cdot 10^{-5} & 0.0009745\\ {\var{1}}^5\,\var{2} & -0.0004095 & 0.0008582\\ {\var{1}}^5 & 0.0009608 & -0.0022\\ {\var{1}}^4\,{\var{2}}^3 & -0.0001209 & -0.004687\\ {\var{1}}^4\,{\var{2}}^2 & -9.623 \cdot 10^{-5} & -0.01261\\ {\var{1}}^4\,\var{2} & 0.004618 & -0.01113\\ {\var{1}}^4 & -0.0119 & 0.02875\\ {\var{1}}^3\,{\var{2}}^3 & 0.001073 & 0.03061\\ {\var{1}}^3\,{\var{2}}^2 & -0.001885 & 0.08062\\ {\var{1}}^3\,\var{2} & -0.02089 & 0.07021\\ {\var{1}}^3 & 0.06392 & -0.1858\\ {\var{1}}^2\,{\var{2}}^3 & -0.004713 & -0.1045\\ {\var{1}}^2\,{\var{2}}^2 & 0.01689 & -0.2651\\ {\var{1}}^2\,\var{2} & 0.02535 & -0.2222\\ {\var{1}}^2 & -0.1343 & 0.6182\\ \var{1}\,{\var{2}}^3 & 0.009004 & 0.1766\\ \var{1}\,{\var{2}}^2 & -0.04181 & 0.4181\\ \var{1}\,\var{2} & 0.03383 & 0.3132\\ \var{1} & 0.05341 & -1.0\\ {\var{2}}^3 & -0.004951 & -0.1193\\ {\var{2}}^2 & 0.02205 & -0.2609\\ \var{2} & -0.02678 & -0.1911\\ 1.0 & 0.02404 & 0.4972 \end{array}\right)$$\noindent The corresponding function is:$$\begin{array}{rcl}\bG_{\textrm{mon}}(\var{1},\var{2}) &=& \dfrac{\bn_{\textrm{mon}}(\var{1},\var{2})}{\bd_{\textrm{mon}}(\var{1},\var{2})}\\ && \\&=& \dfrac{\sum_{\textrm{row}} \bN_\textrm{mon} \odot \mathcal{B}_\textrm{mon}(\var{1},\var{2})}{\sum_{\textrm{row}} \bD_\textrm{mon} \odot\mathcal{B}_\textrm{mon}(\var{1},\var{2})},  \end{array}$$\noindent where,\\$\bn_{\textrm{mon}}(\var{1},\var{2}) = -2.962 \cdot 10^{-7}\,{\var{1}}^6\,{\var{2}}^3-2.327 \cdot 10^{-6}\,{\var{1}}^6\,{\var{2}}^2+2.579 \cdot 10^{-5}\,{\var{1}}^6\,\var{2}-5.661 \cdot 10^{-5}\,{\var{1}}^6+1.347 \cdot 10^{-5}\,{\var{1}}^5\,{\var{2}}^3+5.447 \cdot 10^{-5}\,{\var{1}}^5\,{\var{2}}^2-0.0008236\,{\var{1}}^5\,\var{2}+0.001932\,{\var{1}}^5-0.0002432\,{\var{1}}^4\,{\var{2}}^3-0.0001936\,{\var{1}}^4\,{\var{2}}^2+0.009288\,{\var{1}}^4\,\var{2}-0.02394\,{\var{1}}^4+0.002159\,{\var{1}}^3\,{\var{2}}^3-0.003791\,{\var{1}}^3\,{\var{2}}^2-0.04201\,{\var{1}}^3\,\var{2}+0.1286\,{\var{1}}^3-0.009479\,{\var{1}}^2\,{\var{2}}^3+0.03398\,{\var{1}}^2\,{\var{2}}^2+0.051\,{\var{1}}^2\,\var{2}-0.2701\,{\var{1}}^2+0.01811\,\var{1}\,{\var{2}}^3-0.08409\,\var{1}\,{\var{2}}^2+0.06805\,\var{1}\,\var{2}+0.1074\,\var{1}-0.009959\,{\var{2}}^3+0.04434\,{\var{2}}^2-0.05386\,\var{2}+0.04834$ \\~~\\$\bd_{\textrm{mon}}(\var{1},\var{2}) = -2.162 \cdot 10^{-5}\,{\var{1}}^6\,{\var{2}}^3-5.912 \cdot 10^{-5}\,{\var{1}}^6\,{\var{2}}^2-5.154 \cdot 10^{-5}\,{\var{1}}^6\,\var{2}+0.0001316\,{\var{1}}^6+0.0007204\,{\var{1}}^5\,{\var{2}}^3+0.00196\,{\var{1}}^5\,{\var{2}}^2+0.001726\,{\var{1}}^5\,\var{2}-0.004425\,{\var{1}}^5-0.009428\,{\var{1}}^4\,{\var{2}}^3-0.02537\,{\var{1}}^4\,{\var{2}}^2-0.02239\,{\var{1}}^4\,\var{2}+0.05783\,{\var{1}}^4+0.06157\,{\var{1}}^3\,{\var{2}}^3+0.1622\,{\var{1}}^3\,{\var{2}}^2+0.1412\,{\var{1}}^3\,\var{2}-0.3737\,{\var{1}}^3-0.2103\,{\var{1}}^2\,{\var{2}}^3-0.5333\,{\var{1}}^2\,{\var{2}}^2-0.4469\,{\var{1}}^2\,\var{2}+1.243\,{\var{1}}^2+0.3552\,\var{1}\,{\var{2}}^3+0.8409\,\var{1}\,{\var{2}}^2+0.63\,\var{1}\,\var{2}-2.011\,\var{1}-0.24\,{\var{2}}^3-0.5247\,{\var{2}}^2-0.3844\,\var{2}+1.0$ \\~~\\\noindent \textbf{KST equivalent decoupling pattern} (Barycentric weights $\bc^{\var{l}}$): $$\begin{array}{rclll}\var{2}&: & \left(\begin{array}{ccccccc} -0.3352 & -0.6466 & -0.5925 & -0.6529 & -0.5818 & -0.6404 & -0.7518\\ 1.406 & 2.233 & 2.12 & 2.265 & 2.1 & 2.238 & 2.481\\ -2.047 & -2.583 & -2.525 & -2.612 & -2.516 & -2.597 & -2.729\\ 1.0 & 1.0 & 1.0 & 1.0 & 1.0 & 1.0 & 1.0 \end{array}\right)& \textrm{vec}(.) & := \textbf{Bary}(\var{2}) \\\var{1}&: & \left(\begin{array}{c} -0.03036\\ 0.2096\\ -0.6616\\ 1.0\\ -0.9409\\ 0.5463\\ -0.1235 \end{array}\right)& \textrm{vec}(.) \otimes \bone_{k_{2}} & := \textbf{Bary}(\var{1}) \\\end{array}$$~\\ Then, with the above notations, one defines the following univariate vector functions:  $$ \left\{ \begin{array}{rcl}\bPhi_{1}(\var{1}) &:=& \textbf{Bary}(\var{1}) \odot \mathbf{Lag}(\var{1}) \\\bPhi_{2}(\var{2}) &:=& \textbf{Bary}(\var{2}) \odot \mathbf{Lag}(\var{2}) \\\end{array} \right. $$\noindent The corresponding function is:$$\begin{array}{rcl}\bG_{\textrm{kst}}(\var{1},\var{2}) &=& \dfrac{\bn_{\textrm{kst}}(\var{1},\var{2})}{\bd_{\textrm{kst}}(\var{1},\var{2})}\\ && \\ &=& \dfrac{\sum_{\text{rows}} \bw \odot \bPhi_{1}(\var{1}) \odot \cdots \odot\bPhi_{2}(\var{2})}{\sum_{\text{rows}} \bPhi_{1}(\var{1}) \odot \cdots \odot\bPhi_{2}(\var{2})} . \end{array}$$~\\ \noindent \textbf{KST-like univariate functions} (equivalent scaled univariate functions $\bphi_{1,\cdots,2}$): $$\left\{\begin{array}{rcrcl}z_{1} &=&\bphi_{1}(\var{1}) &=& \cfrac{\bn_{1}}{\bd_{1}} \\z_{2} &=&\bphi_{2}(\var{2}) &=& \cfrac{\bn_{2}}{\bd_{2}} \\\end{array} \right. .$$\noindent where, \\ \noindent $\bn_{1}=0.0002243\,{\var{1}}^6-0.007892\,{\var{1}}^5+0.1012\,{\var{1}}^4-0.5696\,{\var{1}}^3+1.305\,{\var{1}}^2-0.7344\,\var{1}-0.1936$ and \\ \noindent $\bd_{1}=4.465 \cdot 10^{-6}\,{\var{1}}^6+0.0001251\,{\var{1}}^5-0.004381\,{\var{1}}^4+0.05892\,{\var{1}}^3-0.3551\,{\var{1}}^2+1.243\,\var{1}+1.0$, \\ \noindent $\bn_{2}=-0.006807\,{\var{2}}^3+0.11\,{\var{2}}^2-0.3591\,\var{2}+0.08826$ and \\ \noindent $\bd_{2}=0.4793\,{\var{2}}^3+0.889\,{\var{2}}^2+0.9167\,\var{2}+1.0$, \\

%

%%%%%%%%%%%%%%%%
\newpage
\section{Discussions and conclusions}
\label{sec:conclusion}
The following comments would benefit  further investigations in the future:
\begin{itemize}
\item When data are obtained from a \textbf{polynomial or rational function} $\bH$, M1, M2 and M3 are by far the most efficient methods since they are fast, accurate, and recover the exact complexity (without over-fitting). In addition, we demonstrate that whatever the tensor size, the solution is perfectly recovered with model construction time largely acceptable on a standard computer.
\item When data are obtained from a \textbf{non-rational function} $\bH$, the adaptive interpolation point selection  scheme seems a good candidate, and M2 and M5, M6  reveal to be very efficient methods. M2 benefits from the recursive barycentric values construction and is much faster and scalable when complexity increases ($\ord$, tensor size, etc.). %However, M2's current implementation suffers from issues in some configurations.
\item The size (dimension and size on the disk) of the original tensor is the main limit to M4, M5, M6 and M7. Indeed, the computation time is largely dictated by the dimension $\ord$ of the tensor and of its associated size $N_1\times N_2 \times \cdots \times N_\ord$. This stresses the importance of M1, M2 and M3 benefiting from variable decoupling.
\item Algorithmic strategies are under investigation to automatize as much as possible the parameter tuning, the order estimation and interpolation point selection. To make the user experience smoother and propose a solution that robustly solves the tensor approximation problem.
%\item Try to compare the methods with the very same dimension (i.e. complexity)
%\item M1, M2 and M3 never fail
\end{itemize}

\noindent \textbf{General conclusions.} In this report we presented a survey and benchmark methodology to evaluate multivariate model construction on the basis of tensors. In addition, we provided a digest of the variable decoupling feature presented in \cite{Antoulas:2025}, including some new results not presented in the original work. %We believe we paved a first approach to evaluate and we now provide some preliminary conclusions and discussions. 

More in detail, we reported on different methods allowing to construct a surrogate approximate model directly from tensors. Each method optimizes a specific model structure (rational in barycentric basis, MLP and KAN with different splines). Still, each approach and code share the very same input: \textbf{a $\ord$-D tensor}. We believe  that a complete comparison over a large set of tensors, constructed from different functions, is provided, and that the metrics used are appropriate for evaluating the efficiency of the methods. Among the \CAS \, cases considered in this study, most of the approaches successfully reached an appropriate and accurate surrogate. We believe that the approach proposed in \cite{Antoulas:2025} is  a viable candidate to deal with very complex real-life tensors. This method shows a very fast computation time, large flexibility, few tuning parameters while still providing very accurate approximating functions, easily interpretable, scalable to very-large tensors. Improvements of M1, M2 and M3 to meet the practical expectations for non-expert users and reach its full potential will be sought.

Lastly, we want to briefly comment on third the party methods (namely M4, M5, M6 and M7) used in this report: (i) we warmly thank authors for making their code available; (ii) we would like to stress their easy accessibility and usage, together with their sufficiently detailed documentation. %we report on the fact that their use was actually quite simple with documentation sufficiently detailed.%; (iii) we repeat that we may have badly / non-optimally use their code and apologize if so and remain open for modifications and comments.

\textbf{Software availability:} the (research oriented) \texttt{MATLAB} code used to generate the figures and illustrations corresponding to the numerical results presented in this paper is available at \url{https://github.com/cpoussot/mLF} and \url{https://github.com/cpoussot/benchmark_tensor}. %Compiled and industrial-oriented code is available on request at \url{https://mordigitalsystems.fr/}.

%The report aims at being updated. Please share comments, remarks or suggestions to:\\  \texttt{charles.poussot-vassal@onera.fr}.

%%%%%%%%%%%%%%%%
\bibliographystyle{plain}
\bibliography{references}

\end{document}